\documentclass[preprint]{elsarticle}

\usepackage{lineno,hyperref}
\modulolinenumbers[5]

\journal{Journal of Computational Physics}

\usepackage{amsmath}
\usepackage{amssymb}
\usepackage{graphicx}
\graphicspath{{./pics/}}
\usepackage[usenames, dvipsnames]{color}
\usepackage{amsfonts}
\usepackage{amscd}
\usepackage{mathrsfs}
\usepackage{bm}
\usepackage{enumerate}
\usepackage{algorithmic, algorithm}
\newcommand{\R}{\mathbb{R}}
\newcommand{\M}{\mathcal{M}}
\newcommand{\N}{\mathcal{N}}
\newcommand{\CP}{\mathcal{P}}
\newcommand{\bx}{\bm{x}}
\newcommand{\by}{\bm{y}}

\newtheorem{thm}{Theorem}

\begin{document}

\begin{frontmatter}

\title{Scientific Data Interpolation with Low Dimensional Manifold Model\tnoteref{tnote:USgov}}
\tnotetext[tnote:USgov]{The submitted manuscript has been authored by a contractor of the U.S. Government under Contract No. DE-AC05-00OR22725. Accordingly, the U.S. Government retains a non-exclusive, royalty-free license to publish or reproduce the published form of this contribution, or allow others to do so, for U.S. Government purposes.}

\author[UCLA]{Wei Zhu\corref{cor_wei}\fnref{fn_wei}}
\cortext[cor_wei]{Corresponding author}
\ead{weizhu731@math.ucla.edu}

\author[UCLA]{Bao Wang}
\ead{wangbao@math.ucla.edu}

\author[addr:cam-ornl]{Richard Barnard}
\ead{barnardrc@ornl.gov}

\author[addr:cam-ornl,addr:math-utk]{Cory D. Hauck\fnref{fn_cory}}
\ead{hauckc@ornl.gov}

\author[UCLA_physics]{Frank Jenko}
\ead{jenko@physics.ucla.edu}

\author[UCLA]{Stanley Osher\fnref{fn_stan}}
\ead{sjo@math.ucla.edu}

\address[UCLA]{Department of Mathematics, University of California Los Angeles, Los Angeles, CA 90095, USA}
\address[addr:cam-ornl]{Computational and Applied Mathematics Group, Oak Ridge National Laboratory, Oak Ridge, TN 37831, USA}
\address[addr:math-utk]{Department of Mathematics, University of Tennessee, Knoxville, TN 37996-1320, USA}
\address[UCLA_physics]{Department of Physics and Astronomy, University of California, Los Angeles, CA 90095, USA}

\fntext[fn_wei]{This work is supported by ONR Grant N00014- 14-1-0444.}

\fntext[fn_cory]{This material is based, in part, upon work supported by the U.S. Department of Energy, Office of Science, Office of Advanced Scientific Computing and by the Laboratory Directed Research and Development Program of Oak Ridge National Laboratory (ORNL), managed by UT-Battelle, LLC for the U. S. Department of Energy under Contract No. De-AC05-00OR22725.}

\fntext[fn_stan]{This material is based, in part, upon work supported by the U.S. Department of Energy, Office of Science and by National Science Foundation, under Grant Numbers DOE-SC0013838 and DMS-1554564.}

\begin{abstract}
We propose to apply a low dimensional manifold model to scientific data interpolation from regular and irregular samplings with a significant amount of missing information. The low dimensionality of the patch manifold for general scientific data sets has been used as a regularizer in a variational formulation. The problem is solved via  alternating  minimization with respect to the manifold and the data set, and the Laplace-Beltrami operator in the Euler-Lagrange equation is discretized using the weighted graph Laplacian. Various scientific data sets from different fields of study are used to illustrate the performance of the proposed algorithm on data compression and interpolation from both regular and irregular samplings.
\end{abstract}


\end{frontmatter}

\linenumbers

\section{Introduction}

Interpolation and reconstruction of scientific data sets from sparse sampling is of great interest to many researchers from various communities. In many situations,  data are only partially sampled due to logistic, economic, or computational constraints: limited number of sensors in seismic data or hyperspectral data acquisition, low-dose radiographs in medical imaging, coarse-grid solutions of partial differential equations due to computational complexity, etc. Moreover, sometimes one may also intentionally sample partial information of the scientific data set as a straightforward data compression technique. As a result, it has become an important topic to reconstruct the original data set from regular or irregular samplings.

There are typically two ways to approach this problem. The first one is to use the underlying physics to infer the missing data \cite{curvelet_seismic,Ronen1987,Bagaini1996,Stolt2002,Sergey2003}. The drawback is that such techniques are usually problem-specific and not generally applicable to similar problems in other fields of study. Signal and data processing techniques, on the other hand, usually do not require too much prior information of the governing physics. These models intend to fill in the missing information by the properties manifested by the sampled data themselves, while implicitly enforcing common structures from physical intuition in the regularization.

Many signal processing approaches to data interpolation have been studied in the context of image inpainting and seismic data interpolation. Popular interpolation models have been proposed through total variation \cite{ROF,chan_inpaint}, wavelets \cite{mallat2008wavelet,wavelet_inpaint}, and curvelets \cite{curvelets,mac_inpaint,curvelet_inpaint,curvelet_seismic}. After the introduction of the nonlocal mean by Buades et al. in \cite{nlmean}, patch-based techniques exploiting similarity and redundancy of local patches have been extensively studied for inpainting and reconstruction \cite{osher_nonlocal,peyre2008,ebi}. This also leads to a wide  variety of sparse-signal models which assume that patches can be sparsely represented by atoms in a prefixed or learned dictionary \cite{mallat2008wavelet,sparse_inpaint}. Patch-based Bayesian models have also been proposed in image and data interpolation \cite{ple,bpfa}. However, as reported in \cite{ple}, some of the algorithms can only be applied to the interpolation of randomly selected samples, and fail to achieve satisfactory results for uniform grid interpolation. Moreover, most of the methods perform poorly when a significant amount of information ($\ge 95\%$) is missing.

Recently, a low dimensional manifold model (LDMM) has been proposed for general image processing problems \cite{ldmm}. In particular, it achieved state-of-the-art results for image interpolation problems with a significant number of missing pixels. The main idea behind LDMM is that the patch manifold (to be explained in Section \ref{sec:ldmm}) of a real-world 2D image has a much lower intrinsic dimension than that of the ambient space. Based on this observation, the authors used the dimension of the patch manifold as a regularizer in the variational formulation, and the optimization problem is solved using alternating minimization with respect to the image and the manifold. The key step in the algorithm, which involves solving a Laplace-Beltrami equation over an unstructured point cloud sampling the patch manifold, is solved via either the point integral method \cite{pim} or the weighted graph Laplacian \cite{wgl}.

In this work, we apply LDMM to the interpolation of 2D and 3D scientific data sets from either regular or irregular samplings, and demonstrate its superiority when compared to other methods. Moreover, we also compare the performance of LDMM as a sampling-based data compression technique to other standard compression methods. Unlike the other compression methods, sampling-based methods do not require access to the full data set. Although the results of sampling-based algorithms are generally inferior to standard compression methods, they have the advantage of easy implementation in the compression step, and they are also faster in the reconstruction step if only the reconstruction of a small portion of the data set is required. A useful by-product of this comparison is that the standard compression methods are implicitly compared against one another on a set of physically meaningful test cases that can be used for future benchmarks.

The rest of the paper is organized as follows. Section \ref{sec:ldmm} reviews the low dimensional manifold model and justifies its application to scientific data interpolation through a dimension analysis. Section \ref{sec:implementation} outlines the detailed numerical implementaion of LDMM with weighted graph Laplacian which was missing in \cite{wgl}. A comparison of the numerical results on various scientific data interpolation and compression is reported in Section \ref{sec:results}. Finally, we draw our conclusion in Section \ref{sec:conclusion}.

\section{Low Dimensional Manifold Model}
\label{sec:ldmm}
Low dimensional manifold model (LDMM) is a recently proposed mathematical image processing technique which performs particularly well on natural image inpainting \cite{ldmm, ldmm_wgl}. The main observation is that the intrinsic dimension of the patch manifold of a natural image is much smaller than that of the ambient Euclidean space. Therefore it is intuitive to use the dimension of the patch manifold as a regularizer to recover the degraded image. We argue that the same property holds true for scientific data sets. Throughout the entire paper, we present our analysis and algorithm for 3D scientific data sets. The formulation for 2D and higher dimensional data sets follows in a natural way.

\subsection{Patch Manifold and Dimension Analysis}
Consider a 3D datacube $f \in \R^{m\times n \times r}$. For any voxel $\bx \in \bar{\Omega} = \left\{1,2, \ldots, m\right\}\times \left\{1,2, \ldots, n\right\}\times \left\{1,2, \ldots, r\right\}\footnote{The notation $\Omega$ is reserved for the sampled subset of $\bar{\Omega}$.}$, the patch $\CP f(\bx)$ is defined as a vector storing the data values in a 3D cube of size $s_1\times s_2 \times s_3$, with  $\bx$ being the first voxel of the 3D cube in the lexicographic order, i.e. $\bx$ is in one particular corner of the cube\footnote{One can also choose $\bx$ to be the center of the cube, and the result will be similar. The reason is that the reconstruction is performed on patches instead of on voxels. This will be clear in Section \ref{sec:implementation}.}. The patch set $\CP(f)$ of $f$ is the collection of all patches:
\begin{align*}
  \CP(f)  = \left\{ \CP{f}(\bx): \bx \in \bar{\Omega} \right\} \subset \R^d, \quad d = s_1\times s_2\times s_3.
\end{align*}
We assume that the patch set $\CP(f)$, which is a point cloud in $\R^d$, samples an underlying structure $\M$, which is refered to as the patch manifold of $f$. Rigorously speaking, $\M$ is not a smooth manifold, but instead is a collection of manifolds, $\left(\N_l\right)_{l=1}^L$, with different dimensions corresponding to various patterns in the data set, $\M = \cup_{l=1}^L\N_l$. For any $\bm{p} \in \M$, we use the notation $\M(\bm{p})$ to denote the smooth manifold $\N_l$ to which $\bm{p}$ belongs, and $\dim(\M(\bm{p}))$ is the dimension of  $\M(\bm{p})$.

An important assumption is that for scientific data sets, the intrinsic dimension of the patch manifold $\M$ is often much smaller than the dimension of the embedding space $\R^d$. For example, if $f$ is locally smooth at $\bx$ corresponding to smoothly variant region of the data set, then $\CP f(\bx)$ can be approximated by a linear function via Taylor expansion:
\begin{align*}
  \CP f(\bx)(\by) \approx f(\bx) + \nabla f(\bx) \cdot (\by-\bx).
\end{align*}
Therefore, $\M$ can be approximated by a 4D manifold locally at $\CP f(\bx)$. If $f$ is a piecewise smooth function with a sharp interface corresponding to a shock wave, then the patches can be parameterized by the location and orientation of the shock, as well as the gradient and voxel value information in the two regions. This implies that $\M$ is locally close to an 11D manifold. If $f = a(\bx)\cos\left(\theta(\bx)\right)$ models oscillatory structures, then Taylor expansion with respect to $a$ and $\theta$ implies that $\M$ can be locally approximated by a smooth manifold of dimension $8$.

When dealing with 3D data sets of size $256\times 256\times 32$ in our numerical tests, we typically choose patches of size  $6\times 6 \times 4$. This implies that the dimension $d$ of the ambient space is $144$. The dimension analysis above justifies the claim that the patch manifold $\M$ is a low dimensional manifold.

\subsection{Variational Formulation}
Based on the discussion in the previous section, we use the dimension of the patch manifold $\M$ as a regularizer in the following variational formula:
\begin{align}
  \label{eq:manifold-model}
  \min_{f\in \mathbb{R}^{m\times n \times r},\atop \M\subset \mathbb{R}^d} \quad\int_\M\dim(\mathcal{M}(\bm{p}))d\bm{p},\quad \quad \text{subject to:}\quad  
b=\Phi_{\Omega} f,\quad \CP(f)\subset \M,
\end{align}
where
\begin{align*}
\int_\M\dim(\mathcal{M}(\bm{p}))d\bm{p}=\sum_{l=1}^L\int_{\N_l}\dim(\N_l)d\mu_{\N_l}(\bm{p})=\sum_{l=1}^L|\N_l|\dim(\N_l),
\end{align*}
$\mu_{\N_l}$ is the surface measure on $\N_l$, $\Phi_\Omega$ is the sampling operator on the subset $\Omega \subset \bar{\Omega}$, and $b$ is the partially observed data. It is worth mentioning that $\int_\M\dim(\mathcal{M}(\bm{p}))d\bm{p}$ can be thought of as the $L^1$ norm of the local dimension of the manifold $\M$. It has been shown in \cite{ldmm} that the dimension $\dim(\N)$ of any smooth manifold $\N$ can be calculated by the following simple formula:

\begin{thm} 
\label{thm:dim}
Let $\N$ be a smooth submanifold isometrically embedded in $\R^d$. For any $\bm{p}=(p_1,p_2,\cdots,p_d)\in \N$,  
  \begin{align*}
    \dim(\N)=\sum_{i=1}^d\left|\nabla_\N\alpha_i(\bm{p})\right|^2,
  \end{align*}
where $\alpha_i(\bm{p})=p_i$ is the coordinate function, and $\nabla_\N$ is the gradient operator on the manifold $\N$. More specifically, $\nabla_\N \alpha_i = \sum_{s,t=1}^{k}g^{st}\partial_t\alpha_i\partial_s$, where $k$ is the intrinsic dimension of $\N$, and $g^{st}$ is the inverse of the metric tensor.
\end{thm}

The interested reader can refer to \cite{lee2012} for manifold calculus and \cite{ldmm} for the proof. As a result of Theorem \ref{thm:dim}, (\ref{eq:manifold-model}) can be reformulated as:
\begin{align}
  \label{eq:manifold-model-alpha}
    \min_{f\in \mathbb{R}^{m\times n \times r},\atop \M\subset \mathbb{R}^d} \quad\sum_{i=1}^{d}\|\nabla_\M \alpha_i\|_{L^2(\M)}^2, \quad \text{subject to:}\quad  
b=\Phi_{\Omega} f,\quad  \CP(f)\subset \M,
\end{align}
where
\begin{align}\nonumber
\sum_{i=1}^d\|\nabla_\M \alpha_i\|_{L^2(\M)}^2&=\sum_{l=1}^L\sum_{i=1}^d\|\nabla_{\N_l} \alpha_i\|_{L^2(\N_l)}^2=\sum_{l=1}^L\sum_{i=1}^d\int_{\N_l}\left|\nabla_{\N_l}\alpha_i(\bm{p})\right|^2d\mu_{\N_l}(\bm{p})\\\label{eq:l1-dim}
 &= \sum_{l=1}^L|\N_{l}|\dim(\N_l) = \int_\M\dim(\M(\bm{p}))d\bm{p}
\end{align}

The variational problem (\ref{eq:manifold-model-alpha}) can be solved by alternating minimization with respect to $\M$ and $f$. More specifically, given $\M^k$ and $f^k$ at step $k$ satisfying $\CP(f^k) \subset \M^k$:
\begin{itemize}
\item With fixed $\M^k$, update the data $f^{k+1}$ by solving:
\begin{align}
  \label{eq:update-f}
  &\min_{f\in \mathbb{R}^{m\times n\times r}} 
\quad \sum_{i=1}^d\|\nabla_{\M^k} \alpha_i^f\|_{L^2(\M^{k})}^2,\\ \nonumber
&\text{subject to:}\quad \alpha_i^f(\CP(f^k)(\bx))=\CP_i f(\bx),\quad \bx\in \overline{\Omega},\quad i=1,\cdots,d,\\ \nonumber
&\hspace{3.9cm} f(\bx)=b(\bx),\hspace{.3cm} \bx\in \Omega,
\end{align}
where $\CP_i f(\bx)$ is the $i$-th element in the patch at the voxel $\bx$.
\item Update the manifold $\M^{k+1}$ by setting:
  \begin{align*}
    \M^{k+1} = \alpha^{f^{k+1}}(\M^{k})
  \end{align*}
\end{itemize}
If $f^k$ converges to a solution $f^*$, then $\alpha^{f^*}=\alpha$, the identity map, so that $\M^{k}$ converges to a manifold $\M^*$. $f^*$ is then the LDMM approximation of the unknown data.
%Notice that when the iteration above converges, $\bm{\alpha}$ will be close to the identity function which maps the manifold $\M$ onto itself.

The remaining question is how to solve (\ref{eq:update-f}). In \cite{ldmm}, the authors transformed the Euler-Lagrange equation of (\ref{eq:update-f}) into an integral equation, which was solved by the point integral method \cite{pim}. This procedure avoids discretizing the manifold gradient operator $\nabla_\M$, and is shown to perform very well on image inpainting. However, the point integral method involves solving $d$ linear equations on the patch domain per iteration, which makes the numerical procedure very computationally expensive. In \cite{ldmm_wgl}, the authors presented an alternative solution procedure by using the weighted graph Laplacian (WGL) \cite{wgl} to discretize $\nabla_\M$ directly. This speeds up the numerical computation significantly because only one linear equation is to be solved every iteration. We hereby briefly introduce for completeness the intuition and implementation of WGL.

\subsection{Weighted Graph Laplacian}
\label{sec:wgl}
The weighted graph Laplacian (WGL) was recently proposed in \cite{wgl} to smoothly interpolate functions on a point cloud. Let $C = \left\{\bm{c}_1, \bm{c}_2, \ldots, \bm{c}_n\right\}$ be a set of points in $\R^d$, and let $g$ be a function defined on a subset $S = \left\{\bm{s}_1, \bm{s}_2, \ldots, \bm{s}_n\right\}\subset C$. The goal is to extend $g$ to $C$ by finding a smooth function $u$ on $\M$ that agrees with $g$ when restricted to $S$.

The widely used harmonic extension model \cite{zhu2003semi,chung1997spectral} seeks to solve the interpolation problem by minimizing the following energy:
\begin{align}
  \mathcal{J}(u) = \|\nabla_\M u\|_{L^2(\M)}^2, \quad \text{subject to:} \quad u(\bm{p}) = g(\bm{p}) \quad \text{on } S.
\end{align}

A common way to discretize the manifold gradient $\nabla_\M u$ is to use the non-local approximation:
\begin{align*}
  \nabla_\M u (\bm{p})(\bm{q}) \approx \sqrt{w(\bm{p},\bm{q})}\left(u(\bm{p})-u(\bm{q})\right),
\end{align*}
where $w$ is a positive weight function, e.g. $w(\bm{p},\bm{q})=\exp\left(-\frac{\|\bm{p}-\bm{q}\|^2}{\sigma^2}\right)$. With this approximation
  \begin{align}
    \label{eq:gl}
    \mathcal{J}(u) \approx \sum_{\bm{p},\bm{q} \in P} w(\bm{p},\bm{q})\left(u(\bm{p})-u(\bm{q})\right)^2.
  \end{align}
Such discretization leads to the well-known graph Laplacian method \cite{zhu2003semi,Buhler2009laplacian,bertozzi2012diffuse}.

A closer look into the energy $\mathcal{J}$ in (\ref{eq:gl}) reveals that the model will fail to achieve satisfactory results when the sample rate $|S|/|C|$ is very low. More specifically, after rewriting (\ref{eq:gl}) in the following form:
\begin{align}
  \label{eq:gl-split}
  \mathcal{J}(u) = \sum_{\bm{p} \in S}\sum_{\bm{q} \in C} w(\bm{p},\bm{q})\left(u(\bm{p})-u(\bm{q})\right)^2+\sum_{\bm{p} \in C\setminus S}\sum_{\bm{q} \in C} w(\bm{p},\bm{q})\left(u(\bm{p})-u(\bm{q})\right)^2,
\end{align}
one can see that the first term in (\ref{eq:gl-split}) is much smaller than the second term when $|S| \ll |C|$. As a result, the minimizing procedure will prioritize the second term, and therefore sacrifice the continuity of $u$ on the sampled set $S$. An easy remedy for this scenario is to add a large weight $\mu = |C|/|S|$ in front of the first term in (\ref{eq:gl-split}) to balance the two terms:
\begin{align}
  \label{eq:wgl}
  \mathcal{J}_{\text{WGL}}(u) = \mu\sum_{\bm{p} \in S}\sum_{\bm{q} \in C} w(\bm{p},\bm{q})\left(u(\bm{p})-u(\bm{q})\right)^2+\sum_{\bm{p} \in C\setminus S}\sum_{\bm{q} \in C} w(\bm{p},\bm{q})\left(u(\bm{p})-u(\bm{q})\right)^2.
\end{align}
It is readily checked that $\mathcal{J}_{\text{WGL}}$ generalizes the graph Laplacian $\mathcal{J}$ in the sense that $\mathcal{J}_{\text{WGL}}=\mathcal{J}$ when $|S|=|C|$. The generalized energy functional $\mathcal{J}_{\text{WGL}}$ is called  the weighted graph Laplacian.

We point out that such intuition can be made precise by deriving (\ref{eq:wgl}) through the point integral method. The interested reader can refer to \cite{wgl} for the details.

\section{Numerical Implementation}
\label{sec:implementation}

In this section, we provide a detailed explaination of the numerical implementation of LDMM. Using the terminology introduced in Section \ref{sec:wgl}, the functions to be interpolated in (\ref{eq:update-f}) are $\alpha_i$, the point cloud $C$  is $\CP(f^k)$, and the sampled set for $\alpha_i$ is $S_i = \left\{ \CP f^k(\bx) : \CP_if^k(\bx) \text{ is sampled} \right\}$. Based on the discussion in Section \ref{sec:wgl}, (\ref{eq:update-f}) can be discretized into the following problem:
\begin{align}
  \label{eq:discretize-1}
  \min_{f\in \mathbb{R}^{m\times n \times r}} 
\quad &\sum_{i=1}^d\left(\sum_{\bx\in \overline{\Omega}\backslash \Omega_i}\sum_{\by\in \overline{\Omega}} \overline{w}(\bx,\by)(\CP_if(\bx)-\CP_if(\by))^2
\right.\\
&\hspace{2cm}\left.+\mu\sum_{\bx\in \Omega_i}\sum_{\by\in \overline{\Omega}} \overline{w}(\bx,\by)(\CP_if(\bx)-\CP_if(\by))^2\right),\nonumber
\end{align}
\begin{align*}
\text{Subject to:}\quad \quad f(\bx)=b(\bx),\quad  \bx\in \Omega\subset \overline{\Omega},
\end{align*}
where $\mu = \frac{|\bar{\Omega}|}{|\Omega|}$, $\Omega_i = \left\{ \bm{x} \in \bar{\Omega}: \CP_if^k(\bx) \text{ is sampled}\right\}$, the values $\bar{w}(\bx,\by) = w(\CP f(\bx),\CP f(\by))$ form the elements of a matrix $\bm{\bar{W}}$, and $w$ is a symmetric sparse weight function computed from the point cloud $\CP f^k$. More specifically,
\begin{align}
\label{eq:weight}
  w(\bm{p},\bm{q}) = \exp \left(-\frac{\|\bm{p}-\bm{q}\|^2}{\sigma(\bm{p})\sigma(\bm{q})}\right),
\end{align}
where $\sigma(\bm{p})$ is the normalizing factor. In the numerical experiments, the weight $w$ has been truncated to $20$ nearest neighbors using the space-partitioning data structure $k$-d tree \cite{Friedman1977kdtree}. We employ a randomized and approximate version of the algorithm \cite{Muja2009ANN,Zhu2017HSI} implemented in the open source VLFeat package\footnote{http://ww.vlfeat.org} \cite{vedaldi08vlfeat}. The normalizing factor is chosen as the distance between $\bm{p}$ and its $10$th nearest neighbor. 

In order to derive the Euler-Lagrange equation of (\ref{eq:discretize-1}), we define $\mathcal{P}_i$ as the translation operator that maps $f$ into the shifted data set $\CP_i f = \left(\CP_if(\bx)\right)_{ \bx \in \bar{\Omega}}$, where $\CP_if(\bx)$  is the $i$-th element in the patch at the voxel $\bx$ defined in (\ref{eq:update-f}), and a periodic padding is used when patches exceed the domain of the 3D data set. With such padding, the adjoint operator $\CP_i^*$ of $\CP_i$ is equal to its inverse $\CP_i^{-1}$. It is readily checked by standard variational techniques that the Euler-Lagrange equation of (\ref{eq:discretize-1}) is:
\begin{align}
  \label{eq:discretized-el}
  \left\{
  \begin{aligned}
    \left[\sum_{i=1}^d\mathcal{P}_i^*(h_i)+(\mu-1) \sum_{i=1}^d\mathcal{P}_i^*(g_i)\right](\bx) = 0, \quad & \bx\in\bar{\Omega}\setminus \Omega\\
   f(\bx)=b(\bx), \quad&\bx\in \Omega
  \end{aligned}\right.
\end{align}
where
\begin{align}
  \label{eq:def_hg}
  \left\{
  \begin{aligned}
      h_i(\bx)&=\sum_{\by\in \bar{\Omega}}2\bar{w}(\bx,\by)(\mathcal{P}_if(\bx)-\mathcal{P}_if(\by))\\
      g_i(\bx)& =\sum_{\by\in \Omega_i}\bar{w}(\bx,\by)(\mathcal{P}_if(\bx)-\mathcal{P}_if(\by)).
    \end{aligned}\right.
\end{align}

We use the notation $\bx_{\widehat{j}}$ to denote the $j$-th element after $\bx$ in the  patch. It is easy to verify that $\CP_if(\bx) = f(\bx_{\widehat{i-1}})$, and  $\CP_i^*f(\bx) = \CP_i^{-1}f(\bx) = f(\bx_{\widehat{1-i}})$.

Using such notation, we have:
  \begin{align*}
    \mathcal{P}_i^*h_i(\bx) &= h_i(\bx_{\widehat{1-i}})=\sum_{\by\in \bar{\Omega}}2\bar{w}(\bx_{\widehat{1-i}},\by)\left(\mathcal{P}_if(\bx_{\widehat{1-i}})-\mathcal{P}_if(\by)\right)\\
    &=\sum_{\by\in \bar{\Omega}}2\bar{w}(\bx_{\widehat{1-i}},\by)\left(f(\bx)-f(\by_{\widehat{i-1}})\right)\\
    &=\sum_{\by\in \bar{\Omega}}2\bar{w}(\bx_{\widehat{1-i}},\by_{\widehat{1-i}})\left(f(\bx)-f(\by)\right).
  \end{align*}
Therefore
\begin{align}
\label{eq:ph}
  \sum_{i=1}^d\mathcal{P}_i^*(h_i)(\bx) = \sum_{i=1}^d\sum_{\by\in \bar{\Omega}}2\bar{w}(\bx_{\widehat{1-i}},\by_{\widehat{1-i}})\left(f(\bx)-f(\by)\right).
\end{align}
Similarly,
\begin{align}
\label{eq:pg}
  \sum_{i=1}^d\mathcal{P}_i^*(g_i)(\bx) = \sum_{i=1}^d\sum_{\by\in \Omega}\bar{w}(\bx_{\widehat{1-i}},\by_{\widehat{1-i}})\left(f(\bx)-f(\by)\right).
\end{align}

When we substitute (\ref{eq:ph}) and (\ref{eq:pg}) into (\ref{eq:discretized-el}), the Euler-Lagrange equation becomes:
\begin{align*}
  \left\{
  \begin{aligned}
    &\sum_{\by\in\bar{\Omega}}\left(\sum_{i=1}^d2\bar{w}(\bx_{\widehat{1-i}},\by_{\widehat{1-i}})\right)\left(f(\bx)-f(\by)\right)\\
    &+(\mu-1)\sum_{\by\in\Omega}\left(\sum_{i=1}^d\bar{w}(\bx_{\widehat{1-i}},\by_{\widehat{1-i}})\right)\left(f(\bx)-f(\by)\right)=0, &&\bx\in\bar{\Omega}\setminus\Omega\\
   &\hspace{5.5cm}f(\bx)=b(\bx), &&\bx\in \Omega.
  \end{aligned}\right.
\end{align*}
Let $\tilde{w}(\bx,\by)=\sum_{i=1}^d\bar{w}(\bx_{\widehat{1-i}},\by_{\widehat{1-i}})$, i.e. $\bm{\tilde{W}}$ is assembled from translated versions of the original matrix $\bm{\bar{W}}$, then
\begin{align}
  \label{eq:final_el}
  \left\{
  \begin{aligned}
  2\sum_{\by\in\bar{\Omega}}\tilde{w}(\bx,\by)\left( f(\bx)-f(\by)\right)+  {red}(\mu-1) \sum_{\by\in\Omega}\tilde{w}(\bx,\by)\left( f(\bx)-f(\by)\right) = 0, \quad &\bx \in \bar{\Omega}\setminus\Omega\\
  f(\bx)=b(\bx), \quad&\bx\in \Omega.
  \end{aligned}\right.
\end{align}

\begin{figure}[H]
  \centering
  \includegraphics[width=12cm]{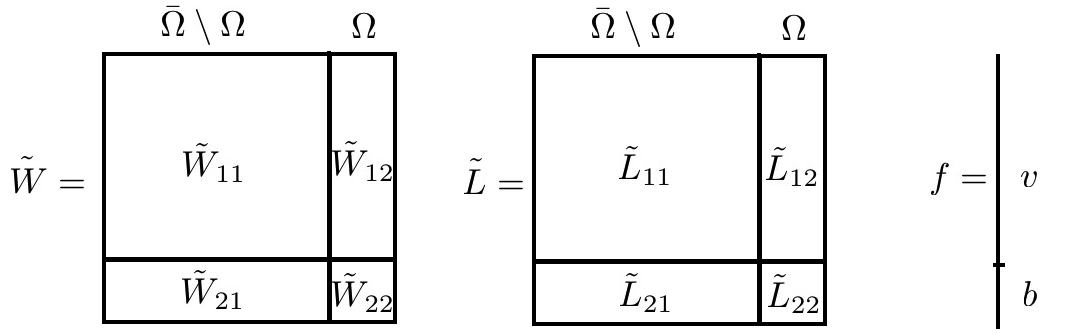}  
  \caption{A visual illustration of the matrix/vector definitions. The matrices $\bm{\tilde{W}}$, $\tilde{\bm{L}}$ and $\bm{f}$ are partitioned into sampled ($\Omega$) and unsampled ($\bar{\Omega} \setminus \Omega$) blocks. For example, $\bm{\tilde{W}}_{12}$ is the matrix corresponding to the weights between unsampled and sampled points.}
  \label{fig:matrix}
\end{figure}

Define the graph Laplacian matrix $\tilde{\bm{L}}$ associated with the new weight matrix $\tilde{\bm{W}}$ as $\tilde{\bm{L}} = \tilde{\bm{D}} - \tilde{\bm{W}}$, where $\tilde{\bm{D}}$ is the diagonal matrix with diagonal entries $\tilde{\bm{D}}(\bx,\bx) = \sum_{\by \in \bar{\Omega}}\tilde{w}(\bx,\by)$. It is easy to check that (\ref{eq:final_el}) can be written in the matrix form:
\begin{align}
\label{eq:final-linear-sys}
  \left(2\bm{\tilde{L}}_{11}+(\mu-1)\bm{\Delta}\right)\bm{v} = (\mu+1)\bm{\tilde{W}}_{12}\bm{b}
\end{align}
where $\bm{\tilde{W}}_{ij}$ and $\bm{\tilde{L}}_{ij}$ are submatrices corresponding to unsampled ($i,j =1$) or sampled ($i,j = 2$) parts of $\bm{\tilde{W}}$ and $\bm{\tilde{L}}$, $\bm{v}$ and $\bm{b}$ correspond to unsampled and sampled parts of $\bm{f}$, and $\bm{\Delta}$ is the diagonal matrix with its diagonal entries equaling the sums of the rows of $\bm{\tilde{W}}_{12}$. See Figure \ref{fig:matrix} for a visual illustration of the definitions of the matrices.

The final LDMM algorithm for 3D scientific data reconstruction from partial sampling is shown in Algorithm \ref{alg:inpaint}. As a remark, we point out that in our current Matlab and C++ implementation, the most time consuming part of the algorithm is step 3, the assembling of the weight matrices, which involves permutations of  sparse weight matrices. We reduce this cost with a parallelization implementation in the matrix assembly step.

\begin{algorithm}
\floatname{algorithm}{Algorithm}
\caption{LDMM for 3D scientific data reconstruction from partial sampling}
\label{alg:inpaint}
\begin{algorithmic}%[1]
\REQUIRE A subsampled data $f|_\Omega = b$.
\ENSURE  Reconstructed data $f$.
\vspace{2mm}
\STATE Initial guess $f^0$.
\vspace{1mm}
\WHILE {not converge}
\vspace{1mm}
\STATE 1. Compute the patch set $\CP(f^k)$ from the current iterate $f^k$.
\vspace{1mm}
\STATE 2. Compute the weight function 
$$\overline{w}(\bx,\by)=w(\CP f^k(\bx),\CP f^k(\by)),\quad \bx,\by\in \overline{\Omega}.$$ 
\STATE 3. Assemble the new weight function
$$\tilde{w}(\bx,\by)=\sum_{i=1}^d\bar{w}(\bx_{\widehat{1-i}},\by_{\widehat{1-i}})$$
%\vspace{1mm}
\STATE 4. Update the data $f^{k+1}$ by solving for variable $v$ in equation \eqref{eq:final-linear-sys}.
\STATE 5. $k\leftarrow k+1$.
\ENDWHILE
\vspace{1mm}
\STATE $f=f^k$.
\end{algorithmic}
\end{algorithm}

\section{Numerical Results}
\label{sec:results}

In this section, we present the numerical results of LDMM on various 2D and 3D scientific data interpolation from either regular or irregular samplings. The performance of LDMM is compared to that of the exemplar-based interpolation (EBI) \cite{ebi} and the piecewise linear estimator (PLE) \cite{ple} in the case of random sampling interpolation. As pointed out in \cite{ple}, PLE fails to work on regular sampling interpolation without a proper initialization (bicubic interpolation in their case). We also noticed in our experiment that the result of EBI on regular sampling interpolation is inferior to that of the simple cubic spline interpolation. Therefore, in the case of regular sampling interpolation,  we instead compare the results of LDMM to the standard methods including cubic spline interpolation, discrete Fourier transform (DFT), discrete cosine transform (DCT), and wavelet transform. Moreover, we also examine the effectiveness of LDMM as a data compression technique and compare it to other standard compression methods including  DFT, DCT, wavelet transform, and tensor decomposition. As for the tensor decomposition methods, we use the singular value decomposition (SVD) for 2D data sets, and the Tucker decomposition \cite{TTB_Software,TTB_Dense} for 3D data sets. The Tucker decomposition is a form of higher-order SVD, which decomposes a tensor into a core tensor multiplied by a matrix along each mode.

\subsection{Description of the Testing Data sets and Parameter Setup}
\begin{figure}[H]
  \centering
  \begin{tabular}{ccc}
    \includegraphics[width=3.5cm]{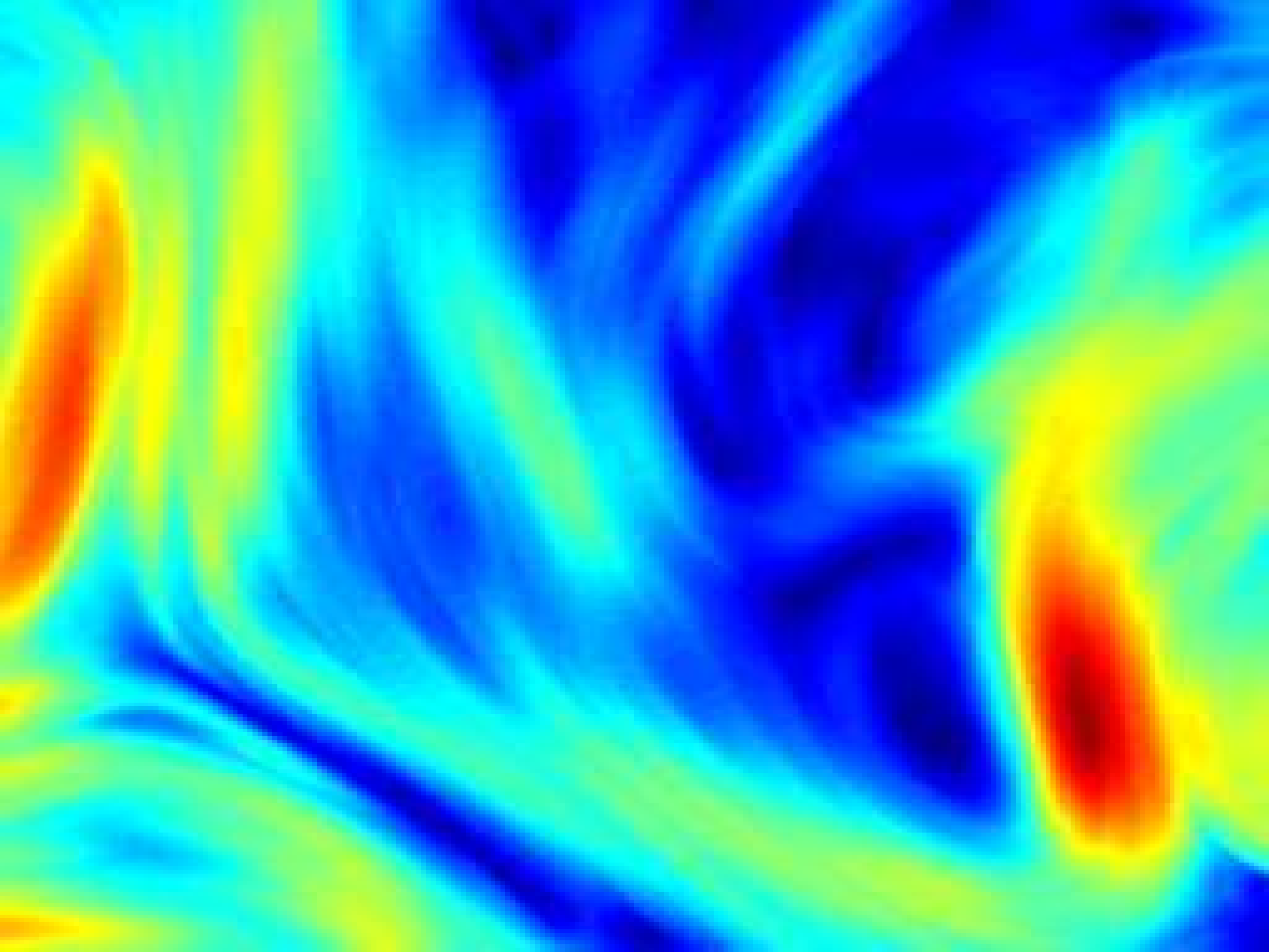}&
    \includegraphics[width=3.5cm]{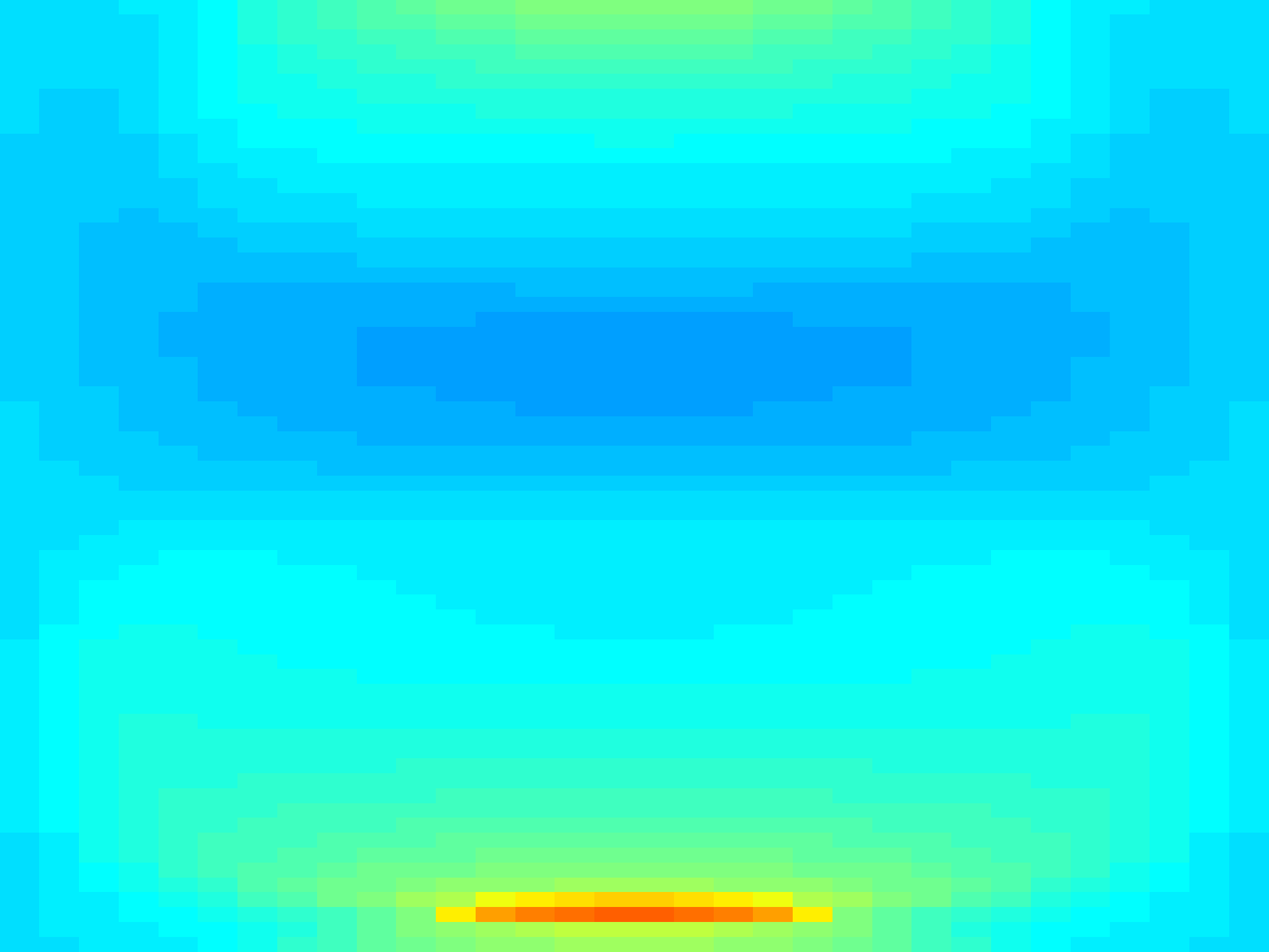}&
    \includegraphics[width=3.5cm]{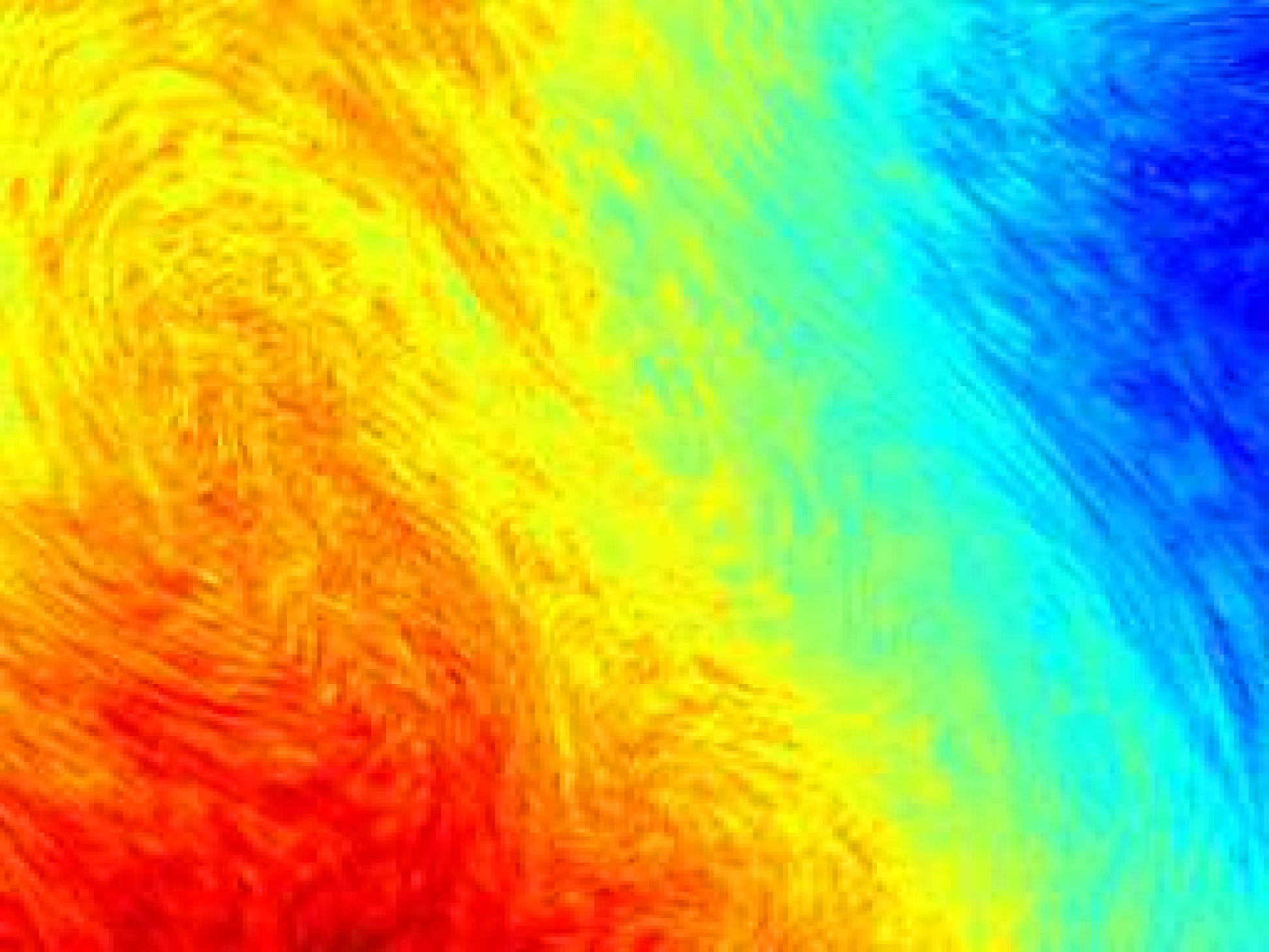}\\
    \includegraphics[width=3.5cm]{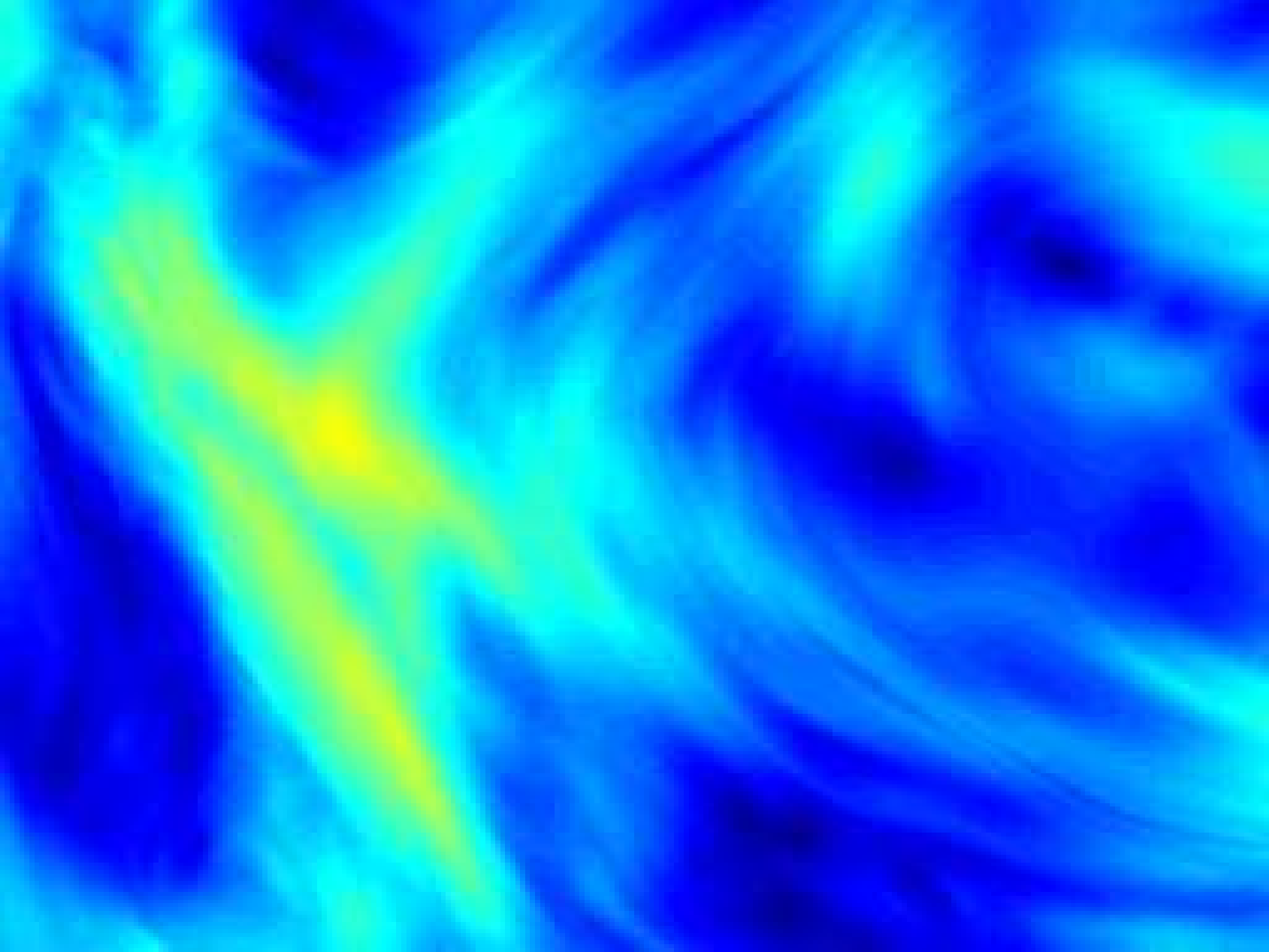}&
    \includegraphics[width=3.5cm]{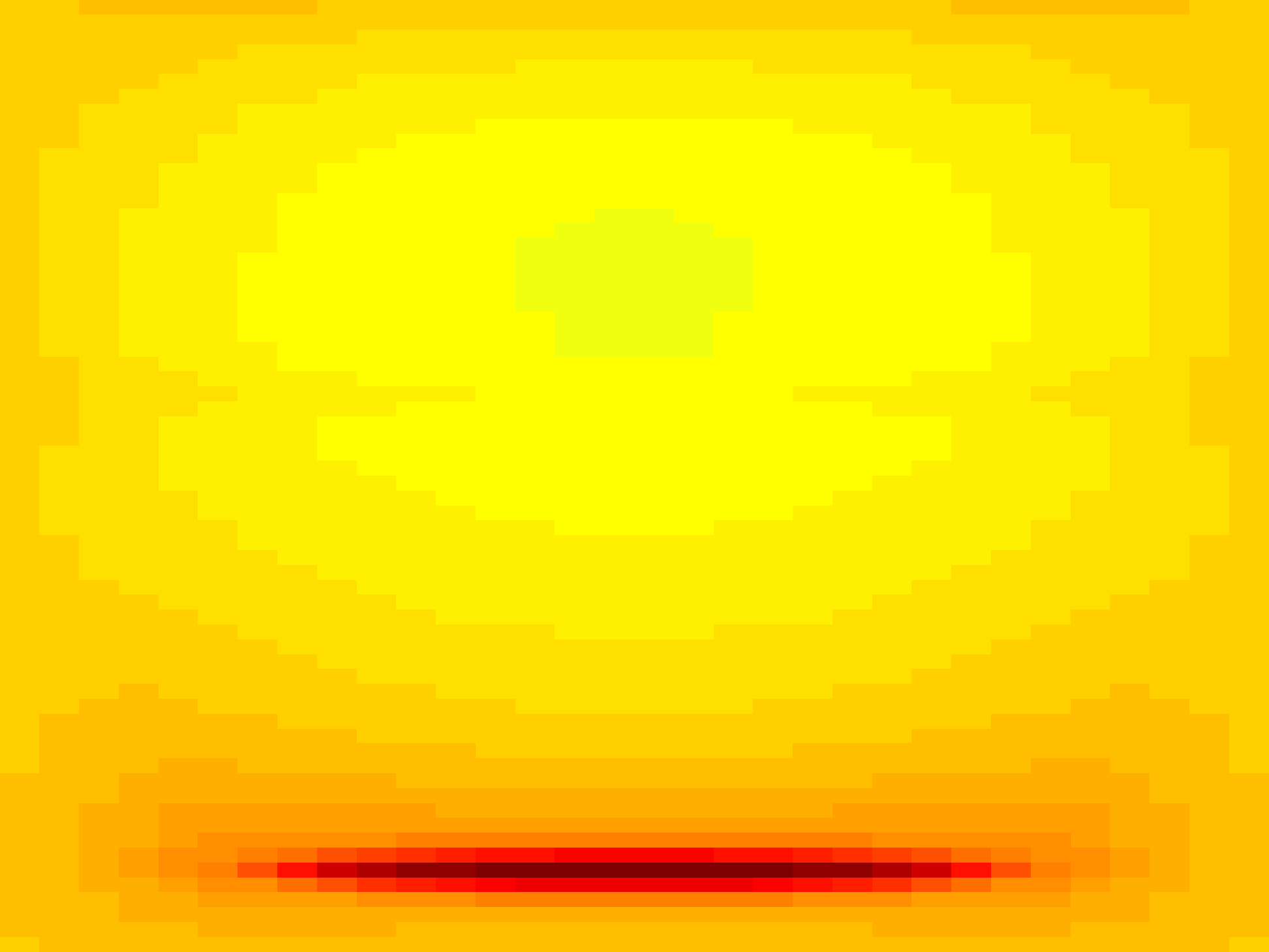}&
    \includegraphics[width=3.5cm]{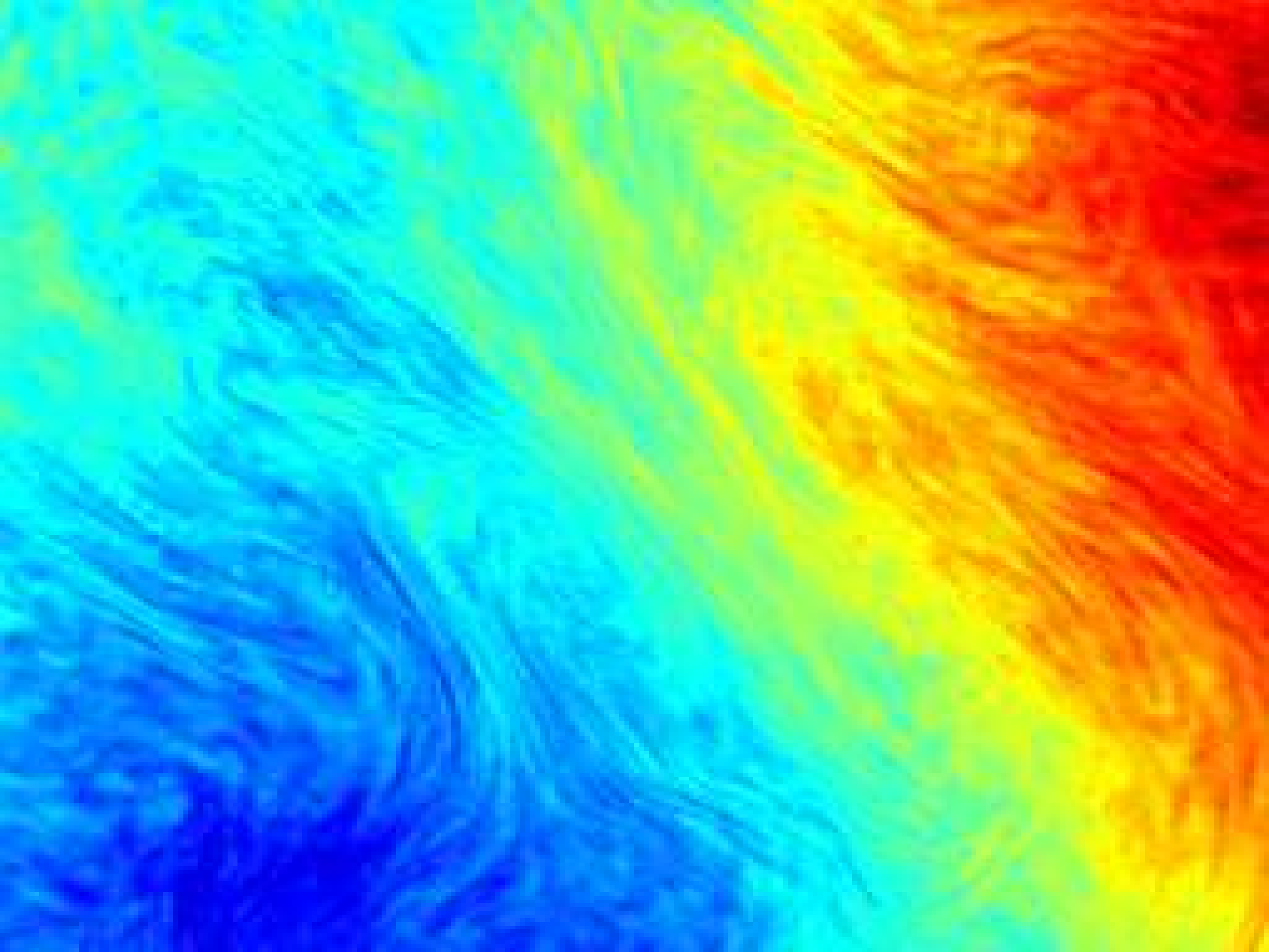}\\
    (a) & (b) & (c)
  \end{tabular}
  \caption{Visual illustrations of the 3D data sets. The two figures in column (a) are 2D spatial cross sections of the 3D plasma (magnetic field) data set at different $z$ coordinates. The figures in column (b) are 2D cross sections of the 3D lattice data set correspongding to angular flux at $x = 0.24$ and $x = 1.18$. The figures in column (c) are 2D spatial cross sections of the 3D plasma (distribution function) data set.}
  \label{fig:datasets_3d}
\end{figure}

\begin{figure}[H]
  \centering
  \begin{tabular}{ccc}
    \includegraphics[width=3.5cm]{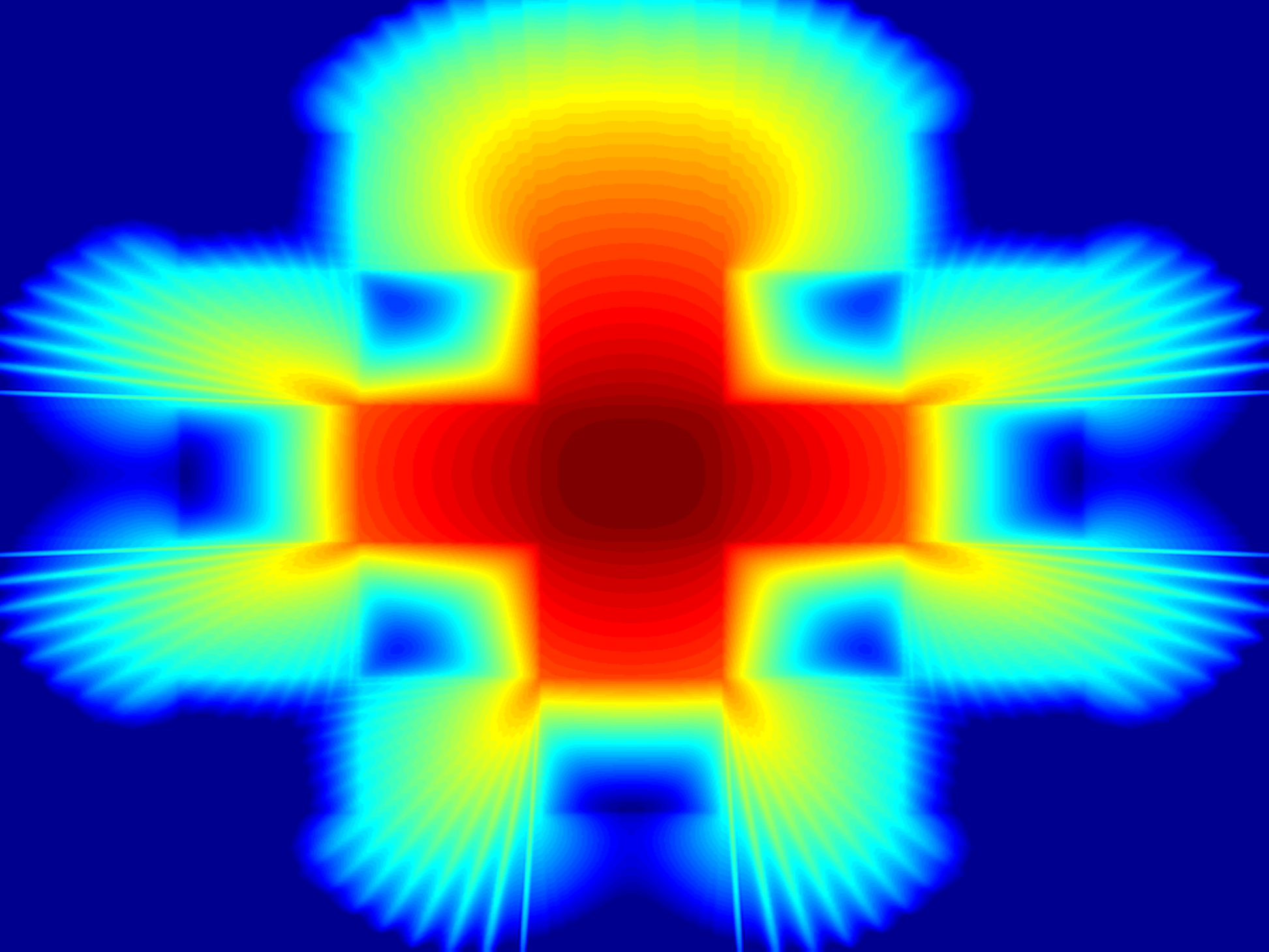}&
    \includegraphics[width=3.5cm]{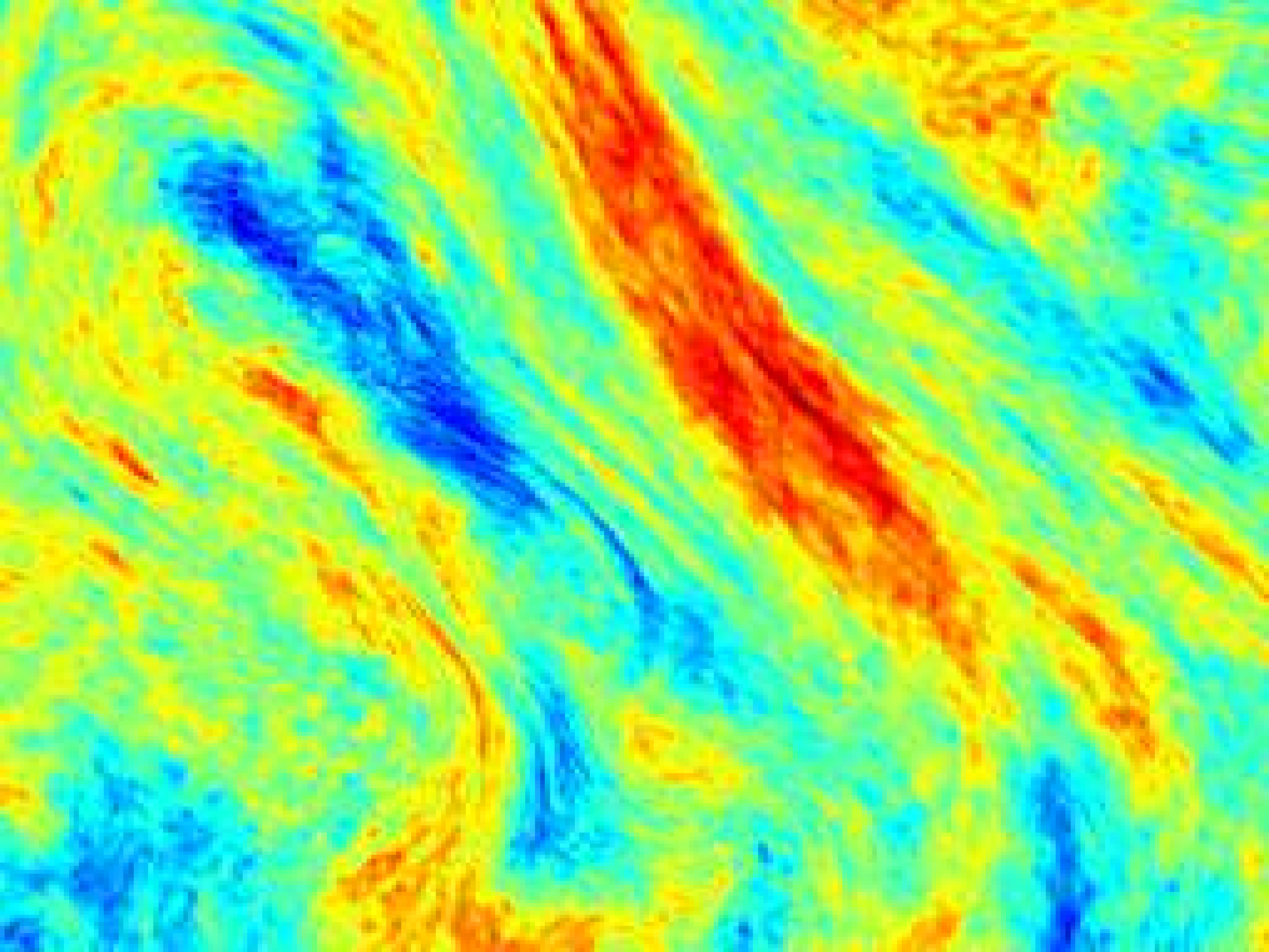}&
    \includegraphics[width=3.5cm]{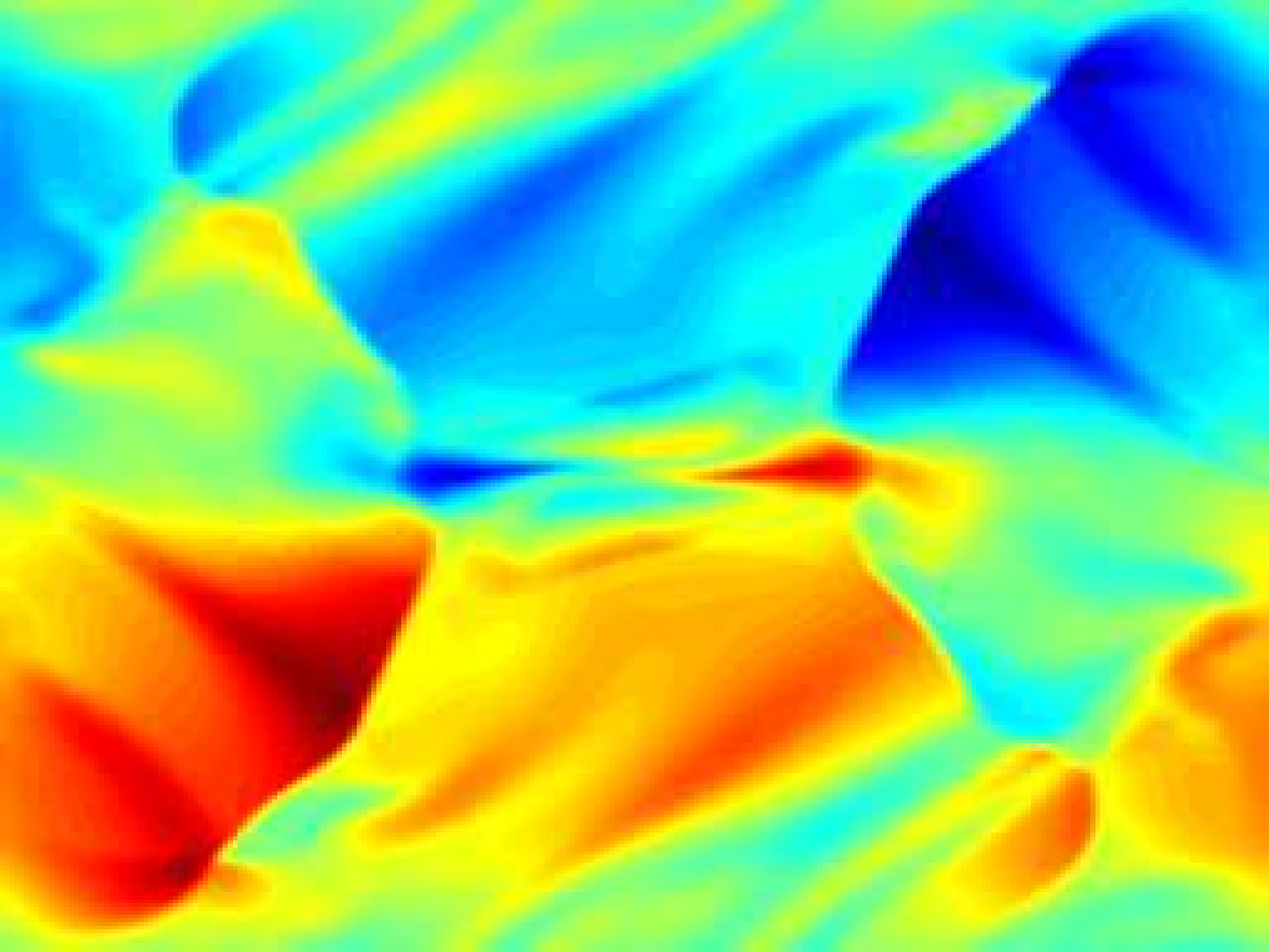}\\
    (a) & (b) & (c)
  \end{tabular}
  \caption{Visual illustrations of the 2D data sets. (a) 2D lattice. (b) 2D plasma (distribution function). (c) 2D vortex.}
  \label{fig:datasets_2d}
\end{figure}

\begin{figure}[H]
  \centering
  \includegraphics[width=3.5cm]{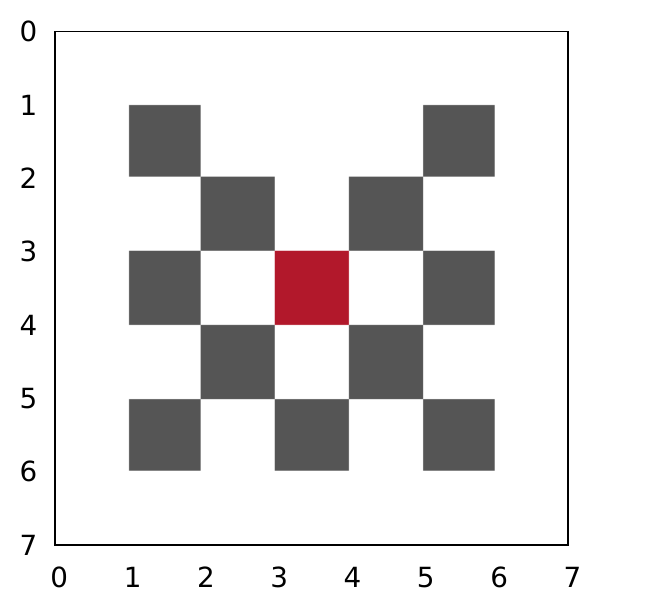}
  \caption{The schematic of the lattice benchmark problem. A central particle source region (red) is surrounded by a checkerboard array of highly absorbing material (gray) embedded within a lightly scattering material (white).}
  \label{fig:checkerboard}
\end{figure}

The algorithms are tested on six scientific data sets, three of which are three-dimensional. See Figure \ref{fig:datasets_3d} and Figure \ref{fig:datasets_2d} for visual illustrations of the data sets.
\begin{itemize}
  \item \textbf{3D plasma (magnetic field)}: The data set is taken from a gyrokinetic simulation of Alfv{\'e}nic turbulence in 5D phase space (3D real space plus 2D velocity space, with the fast gyroangle dependence removed) \cite{told15}, carried out with the GENE code \cite{jenko00}. It represents a snapshot of the magnitude of magnetic field fluctuations in real space during the statistically quasi-stationary state of fully developed turbulence. In this simulation, the focus is on the dissipation range of this weakly collisional turbulent plasma which cannot be described adequately by magnetohydrodynamics (MHD). Gyrokinetics offers an efficient description of the very tail of the MHD cascade. The size of this data is $256 \times 256 \times 32$.
  \item \textbf{3D/2D lattice}:
The lattice benchmark problem, originally due to Brunner \cite{Brunner2002,Brunner2005}, is a two-dimensional cartoon of a nuclear reactor assembly that has become a common test problem of angular discretization methods for kinetic equations of radiation transport \cite{hauck2013collision, McClarren:2010p2380,McClarren:2010p2381,Schaefer:2009p2385}.   

A schematic of the problem is shown in Figure \ref{fig:checkerboard}.  It involves a particle source surrounded by a checkerboard array of highly absorbing material (gray) embedded within a lightly scattering material (white).  Particles are emitted into the domain through a central source region (red).   

The simulated quantity is a distribution function that depends on five independent variables:  two spatial, two angular, plus time.   The data used here was generated using the algorithm described in \cite{crockatt2017arbitrary} which combines a third-order space-time discretization (discontinuous Galerkin in space and integral deferred correction in time) and an angular discretization based on a tensor product collocation scheme.   

We consider for this problem two quantities of interest.  The first (\textbf{2D lattice}) is the angular average of the distribution function at a fixed time; this is a two-dimensional data set of size $896 \times 896$.  The second is the distribution function at a fixed time and fixed vertical location along the line $y=4.5$.  This is a three dimensional data set of size $188 \times 64 \times 32$.  Both sets of data are given in log scale.
  \item \textbf{3D/2D plasma (distribution function)}:
This data set is again taken from a gyrokinetic simulation of Alfv{\'e}nic turbulence in 5D phase space (3D real space plus 2D velocity space, with the fast gyroangle dependence removed) described in \cite{told15}. The 3D data set describes the distribution function for the ion species as a function of the two spatial coordinates perpendicular to the background magnetic field and of the velocity parallel to this guide field at a given value of perpendicular velocity and time.  Meanwhile, the 2D data set describes a snapshot of the same distribution function for the ion species as a function of the two perpendicular spatial coordinates integrated over velocity space. The sizes of the 3D and 2D data sets are $256\times 256\times 32$ and $256\times 256$ respectively.
  \item \textbf{2D vortex}: This data set comes from a numerical solution of the Orszag-Tang vortex system \cite{orszag1979small}, which provides a model of complex flow with many features of magnetohydrodynamics systems. Starting from a smooth state, the system evolves into turbulance, generating complex interactions between different shock waves. The data set used in this paper is the numerical solution at time $t = 2$ of the density component obtained with the third order Chebyshev polynomial approximate Osher-Solomon scheme \cite{Castro2016347} on a $256 \times 256$ uniform mesh.
\end{itemize}


For irregular sampling interpolation, the algorithms are tested to reconstruct the original data sets from $5\%$ and $10\%$ random subsamples. For the regular aliased sampling, the original 2D data sets are decimated by a factor of $4$ in both directions; for 3D data sets, we consider two types of sampling procedures: downsampling by a factor of $2$ in all directions, or by a factor of $4$ in only the first two dimensions.

For all the data sets listed above, the weight matrices in LDMM are truncated to $20$ nearest neighbors, and the normalizing factor $\sigma(\bm{p})$ in (\ref{eq:weight}) is chosen as the distance between $\bx$ and its $10$th nearest neighbor. The patch sizes chosen for different data sets are listed in Table \ref{tab:patch_size}. The reason why the 2D plasma (distribution function) data set uses a much larger patch size, $16\times 16$ instead of $6 \times 6$, is that the structures in this data set are much more complicated than the other data sets. This complexity implies a much  higher intrinsic dimension of the patch manifold. Therefore a larger patch size is chosen so that the manifold dimension can be still smaller than that of the embedding space. Notice also that $6\times 6\times 1$ patch size is chosen for the 3D plasma (magnetic field) data set. This is because of the low resolution of the data set in the third dimension. However, $6\times 6\times 4$ patches are chosen in the $2\times 2\times 2$ regular down sampling. This is because we want to avoid patches that do not contain any sampled voxels.

\begin{table}[H]
  \centering
  \begin{tabular}{||c|c c c c c||} 
    \hline
     & 5\% & 10\% & $4\times 4$ & $4\times 4 \times 1$ & $2\times 2\times 2$\\
    \hline
    2D lattice & $6\times 6$ & $6\times 6$ & $6\times 6$ & N/A & N/A\\
    \hline
    2D plasma (D) & $16 \times 16$ & $16 \times 16$ & $16 \times 16$ & N/A & N/A\\
    \hline
    2D vortex & $6\times 6$ & $6\times 6$ & $6\times 6$ & N/A & N/A\\
    \hline
    3D plasma (M)& $6\times 6 \times 1$ & $6\times 6 \times 1$ & N/A & $6\times 6 \times 1$ & $6\times 6 \times 4$\\
    \hline
    3D lattice & $4\times 4\times 4$ & $4\times 4\times 4$ & N/A & $4\times 4\times 4$&$4\times 4\times 4$\\
    \hline
    3D plasma (D) & $6\times 6 \times 4$ & $6\times 6 \times 4$ & N/A & $6\times 6 \times 4$ & $6\times 6 \times 4$\\
    \hline
 \end{tabular}
  \caption{Patch sizes for different datasets. The first row of the table indicates the different types of irregular and regular downsampling procedures. 3D/2D plasma (D) stands for 3D/2D plasma (distribution function), and 3D plasma (M) stands for 3D plasma (magetic field).}
  \label{tab:patch_size}
\end{table}

The quality of the reconstruction $\hat{f}$ of the original data $f\in \R^{m\times n \times r}$ ($r=1$ for 2D data sets) is evaluated in the following three norms:
\begin{align}
  \label{eq:l1}
  &\|e\|_1 = \frac{1}{mnr}\sum_{i,j,k} |e_{i,j,k}/R|,\\
  \label{eq:l2}
  &\|e\|_2 = \left(\frac{1}{mnr}\sum_{i,j,k}|e_{i,j,k}/R|^2\right)^{\frac{1}{2}},\\
  \label{eq:linfty}
  &\|e\|_\infty = \max_{i,j,k} |e_{i,j,k}/R|,
\end{align}
where $e = f-\hat{f}$ is the error of the reconstruction, $R = \max_{i,j,k}\hat{f}_{i,j,k} - \min_{i,j,k}\hat{f}_{i,j,k}$ is the numerical range of the data set. Moreover, the peak signal-to-noise ratio (PSNR), which is related to (\ref{eq:l2}), is also given to measure the performance of the algorithms:
\begin{align}
  \label{eq:psnr}
  PSNR = 10 \log_{10}\left(\frac{1}{\|f-\hat{f}\|_2^2}\right).
\end{align}

\subsection{Interpolation with Random Sampling}

\begin{figure}[H]
  \centering
  \begin{tabular}{cccc}
    Original& EBI (36.32dB)& PLE (40.01dB)  & LDMM (\textbf{42.55dB})\\
    \includegraphics[width = 2.5cm]{antonio_2d_original}&
    \includegraphics[width = 2.5cm]{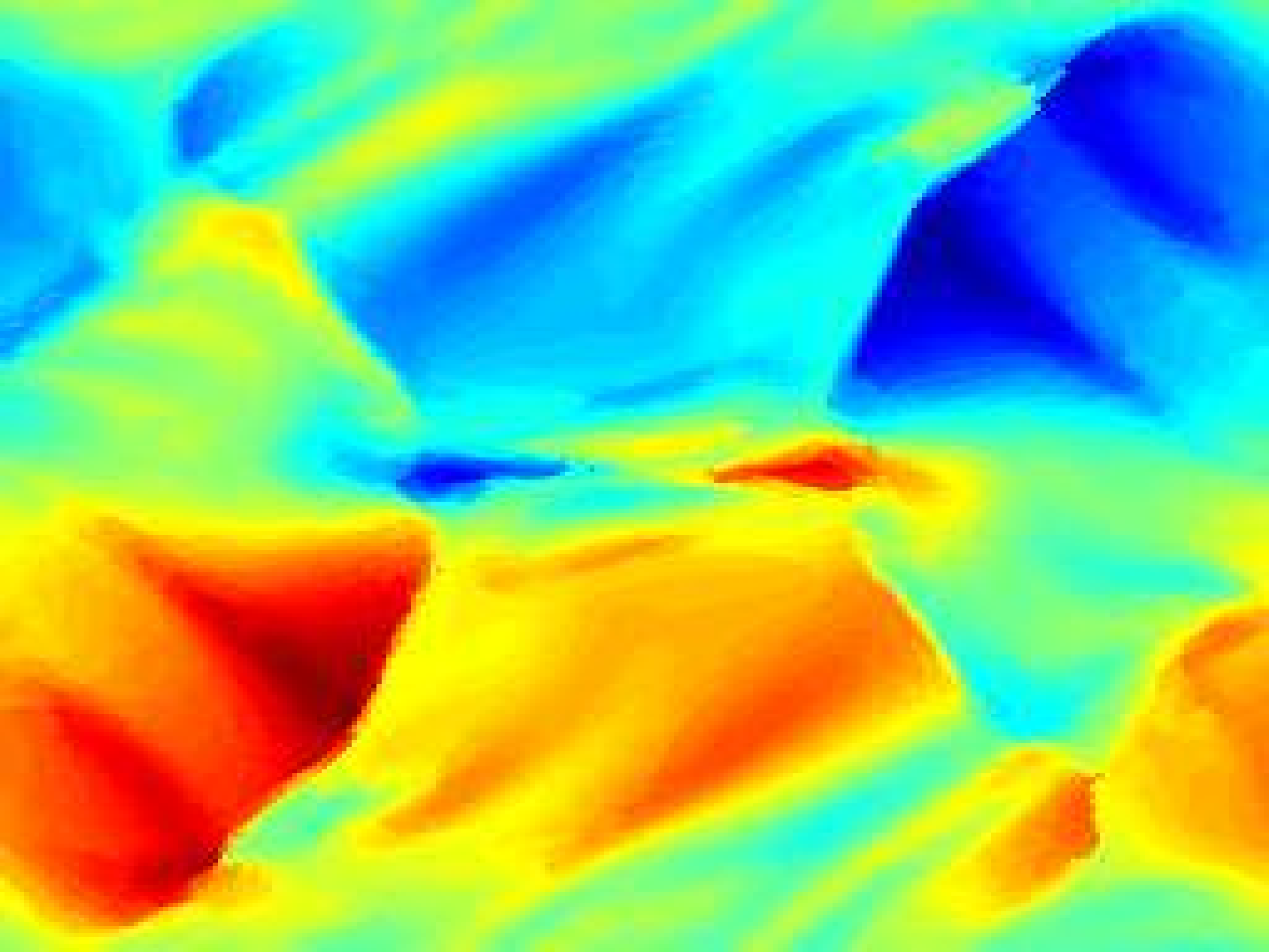}&
    \includegraphics[width = 2.5cm]{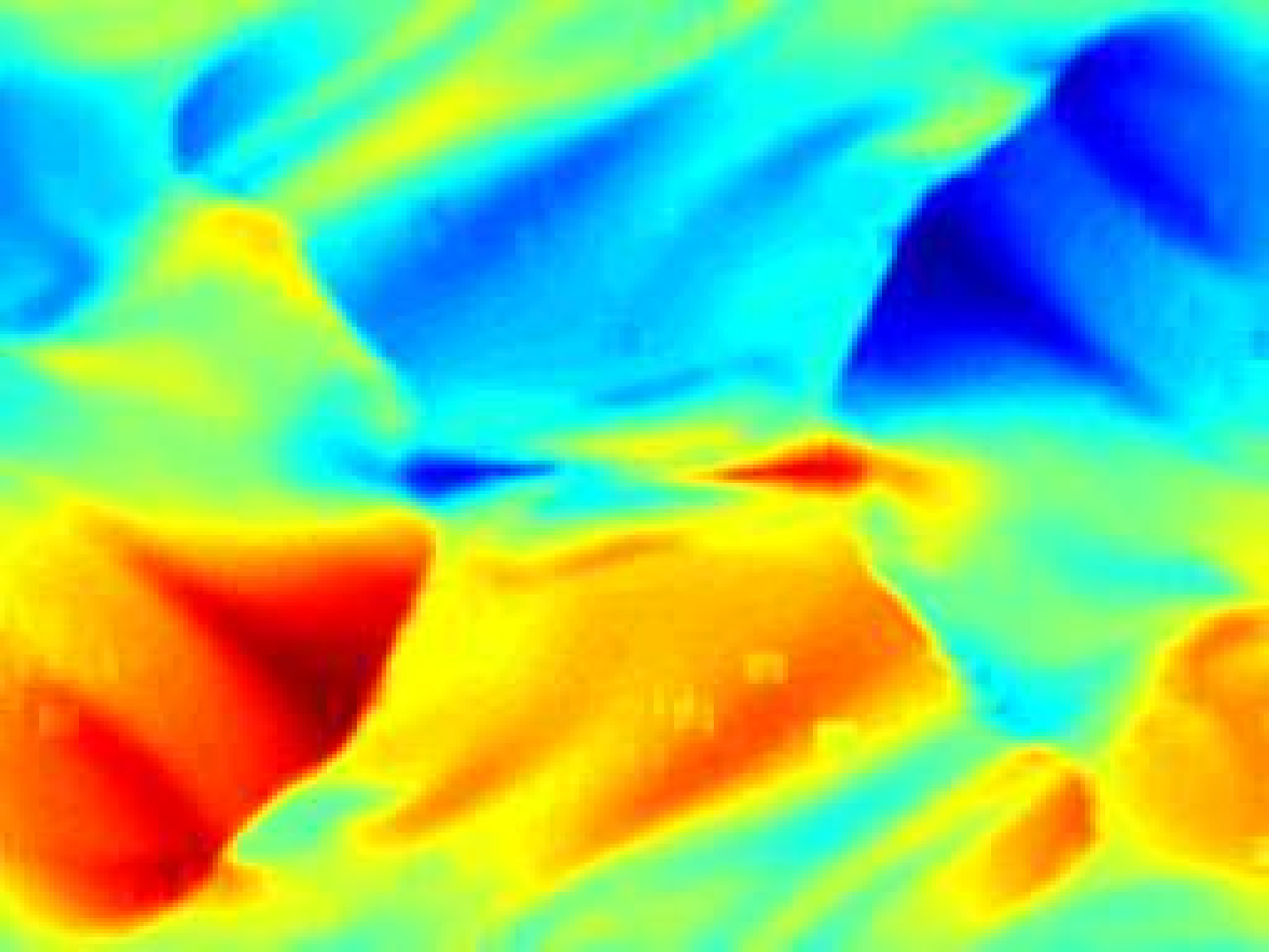}&
    \includegraphics[width = 2.5cm]{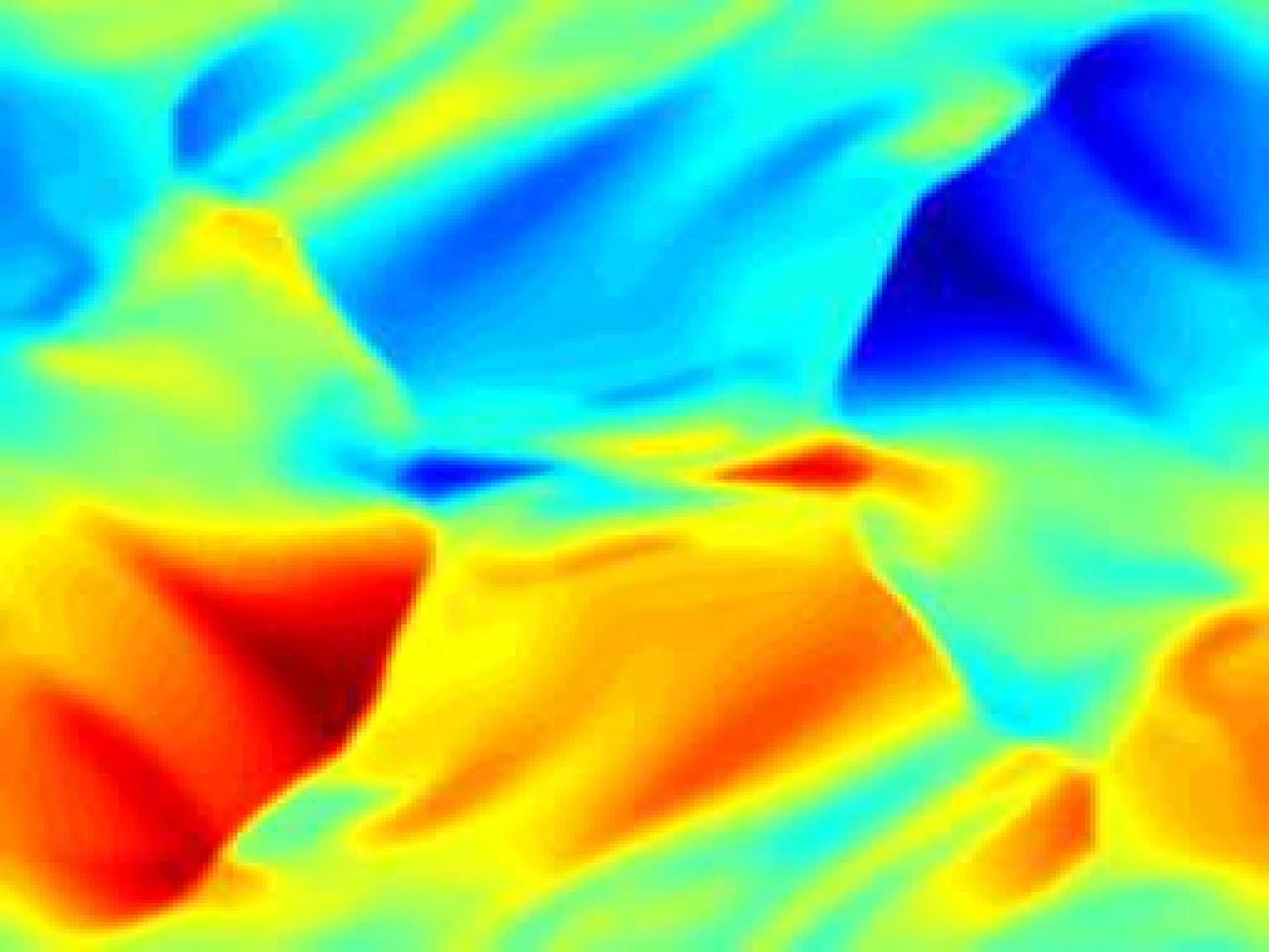}\\
    Subsample & Error & Error & Error\\    \vspace{.5cm}
    \includegraphics[width = 2.5cm]{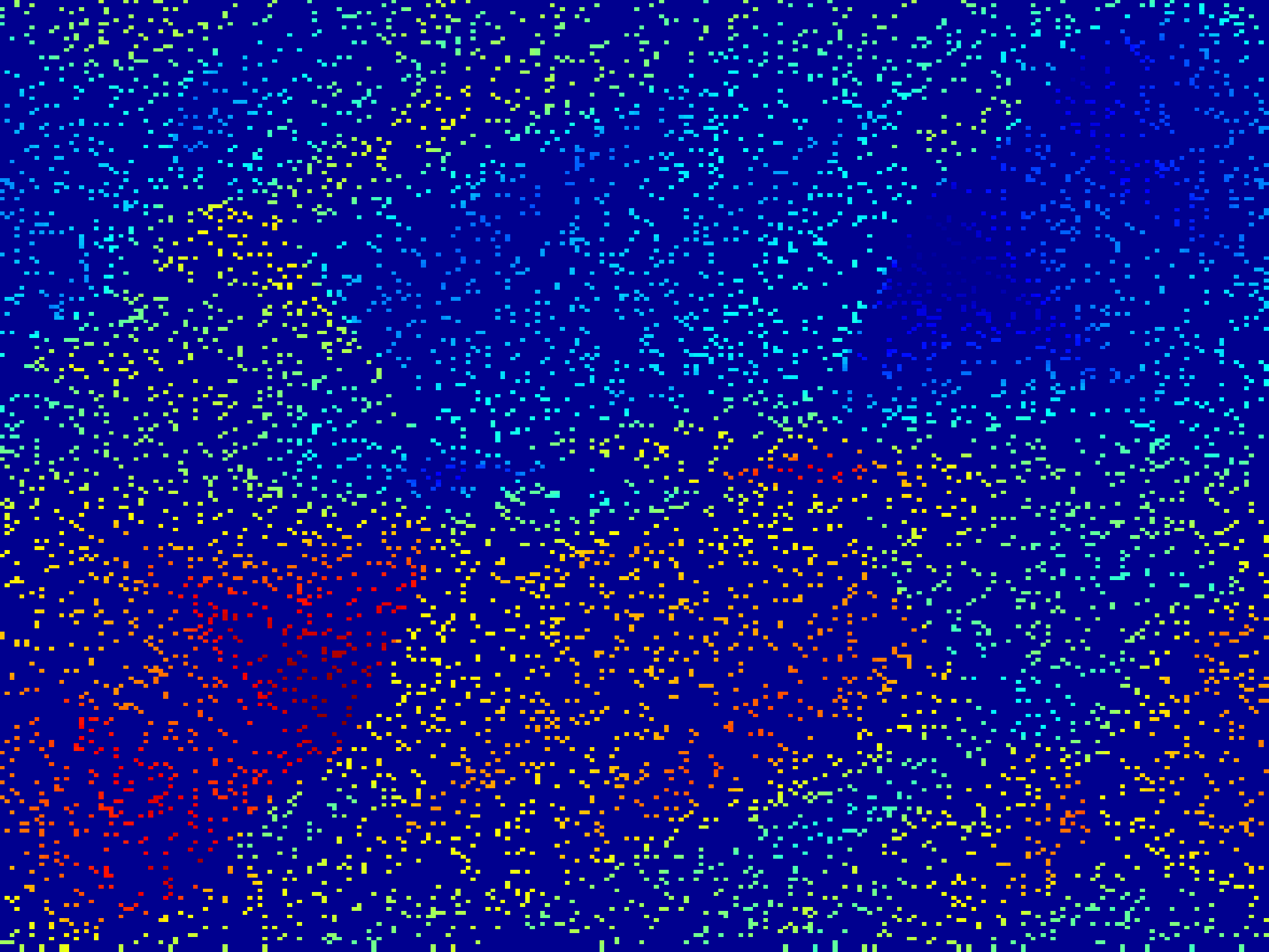}&    
    \includegraphics[width = 2.5cm]{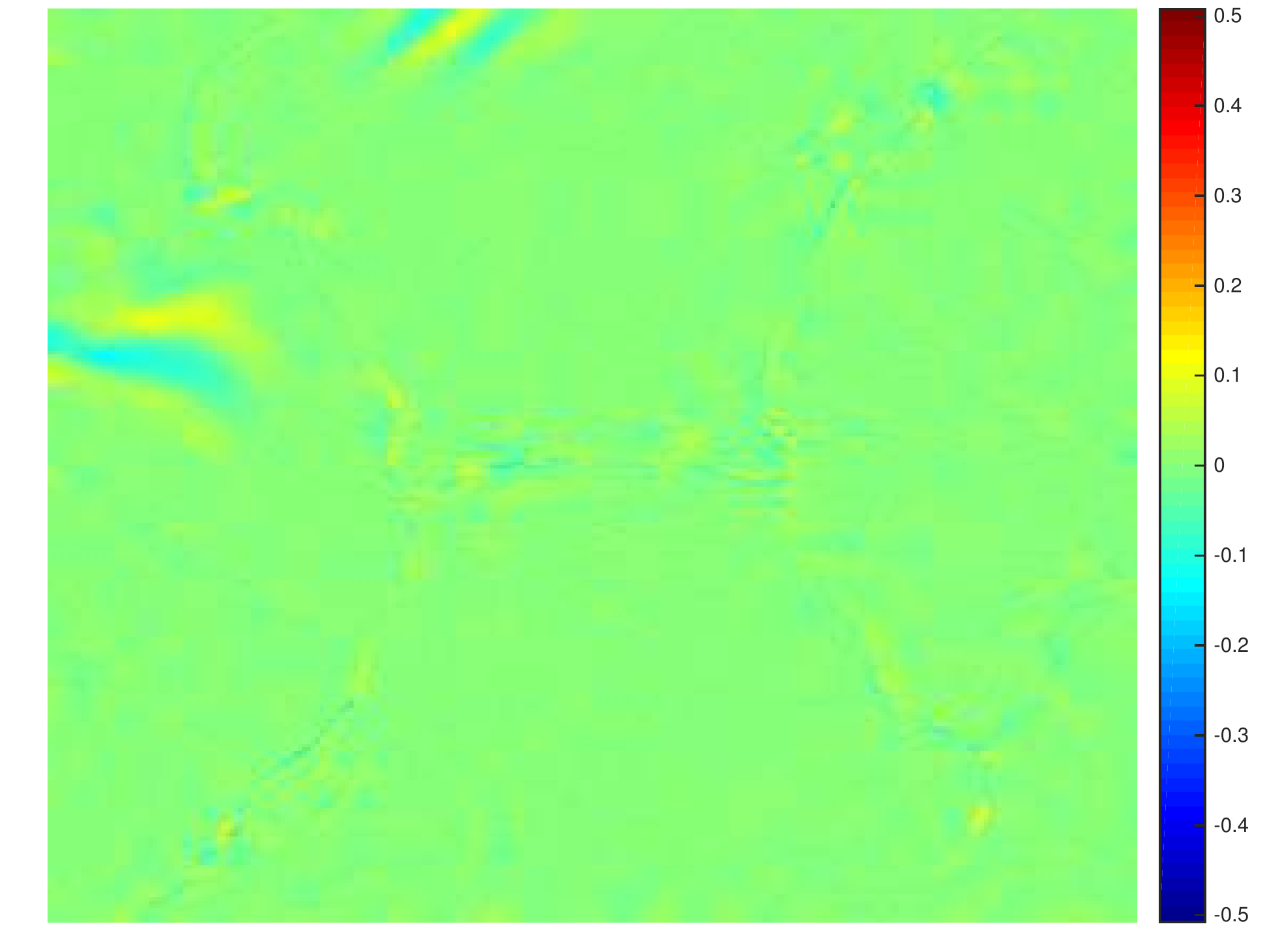}&    
    \includegraphics[width = 2.5cm]{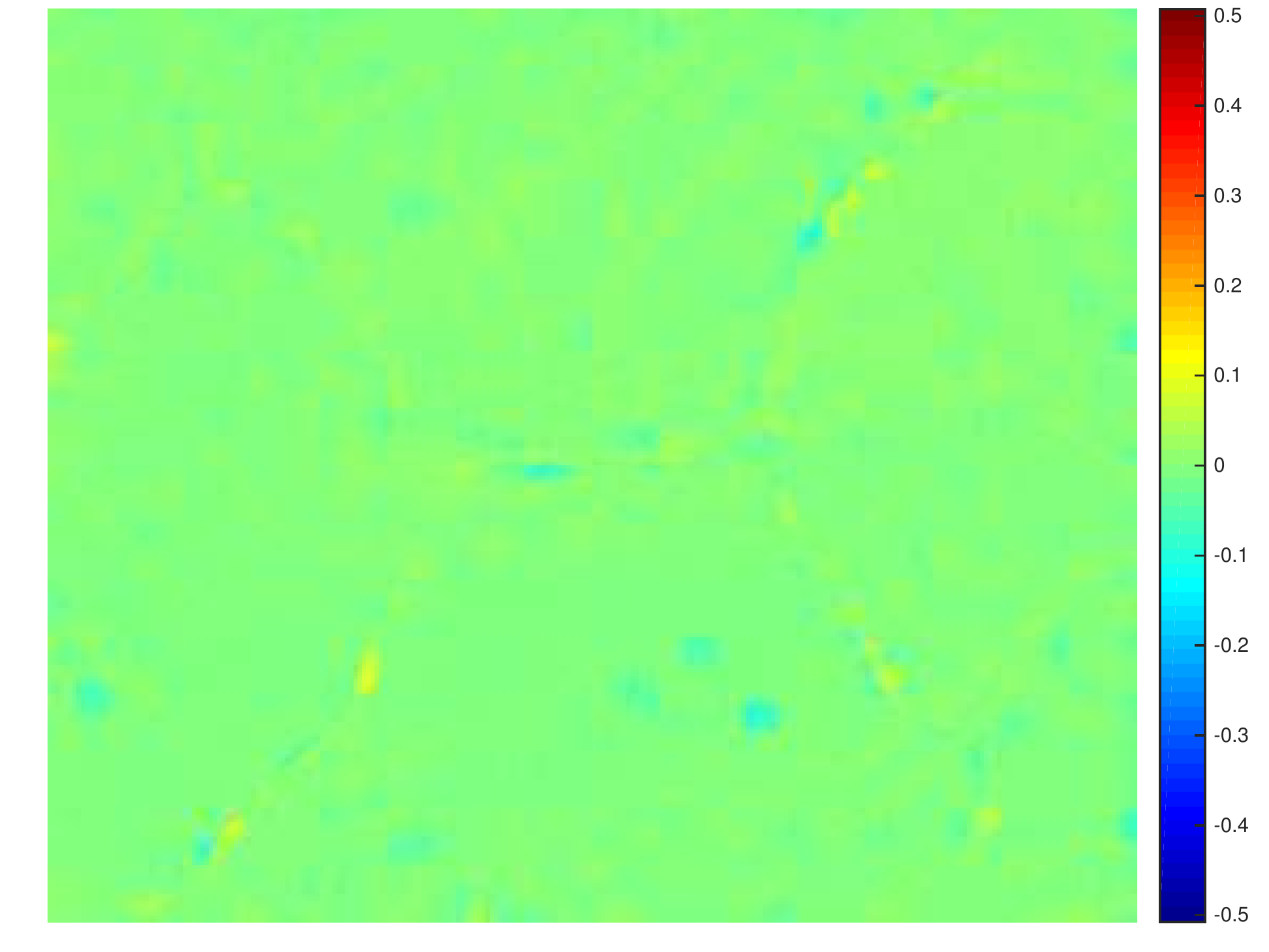}&
    \includegraphics[width = 2.5cm]{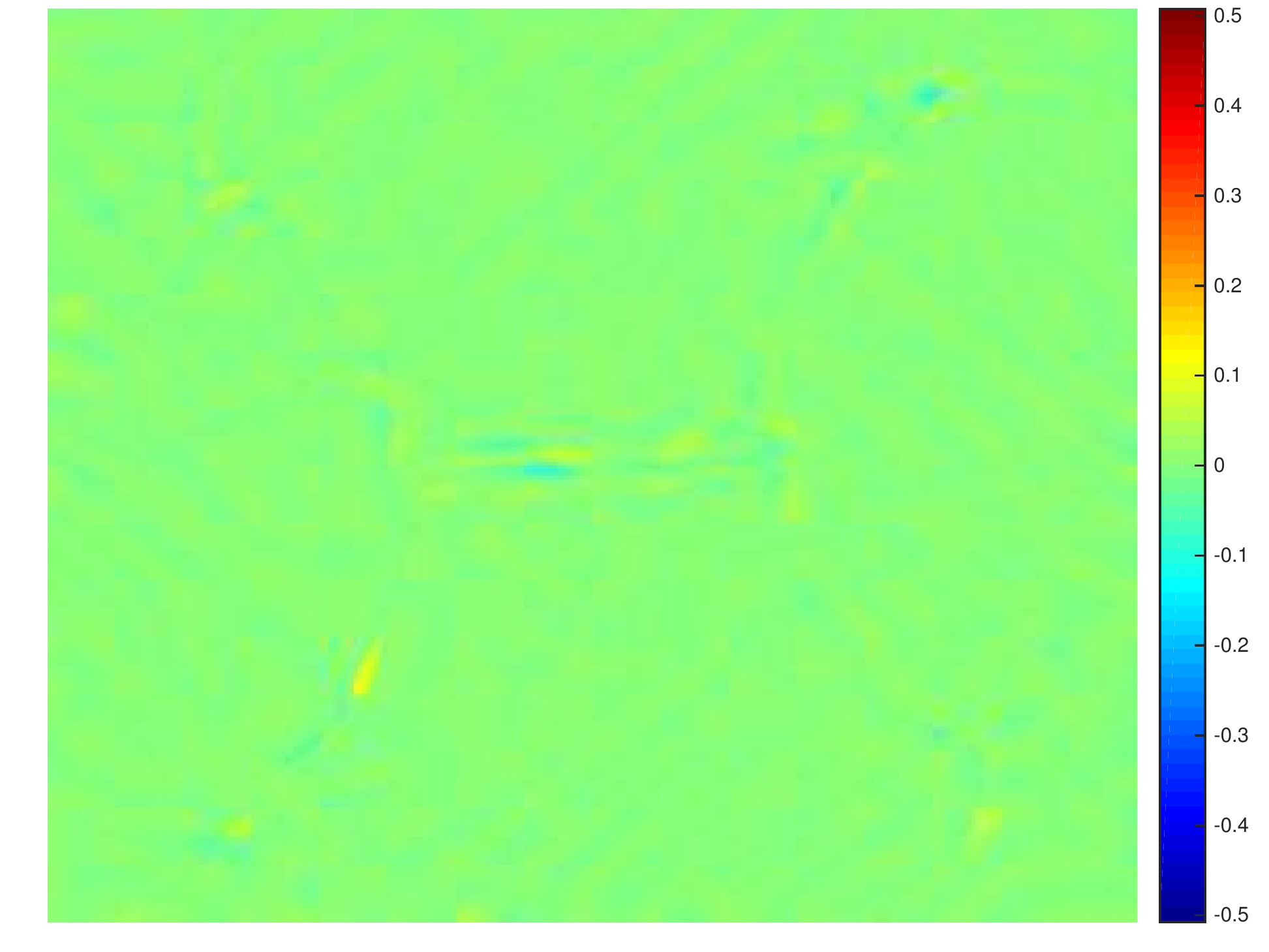}\\
    Original& EBI (26.77dB)& PLE (28.48dB)  & LDMM (\textbf{29.56dB})\\
    \includegraphics[width = 2.5cm]{shock_2d_original}&
    \includegraphics[width = 2.5cm]{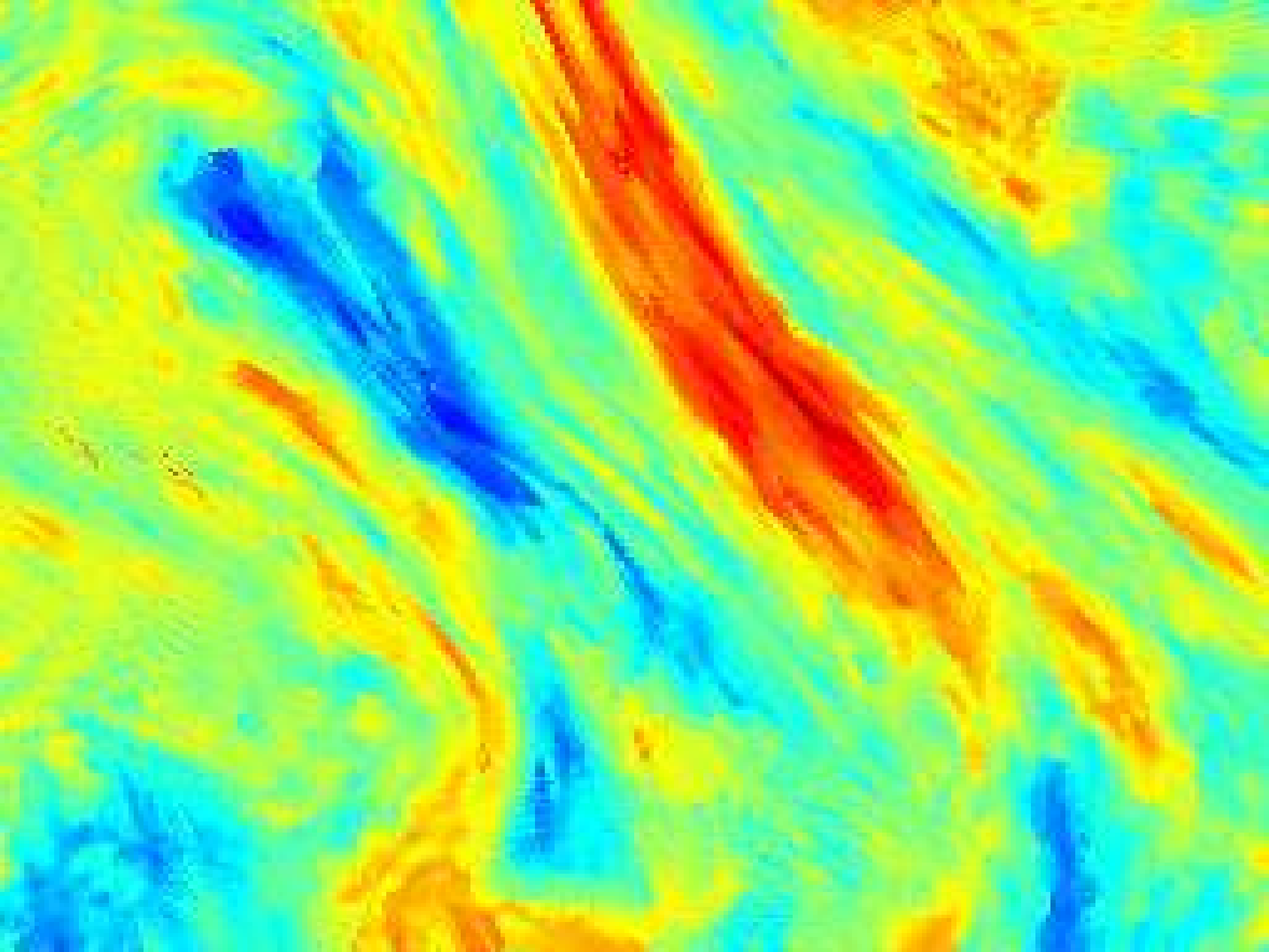}&
    \includegraphics[width = 2.5cm]{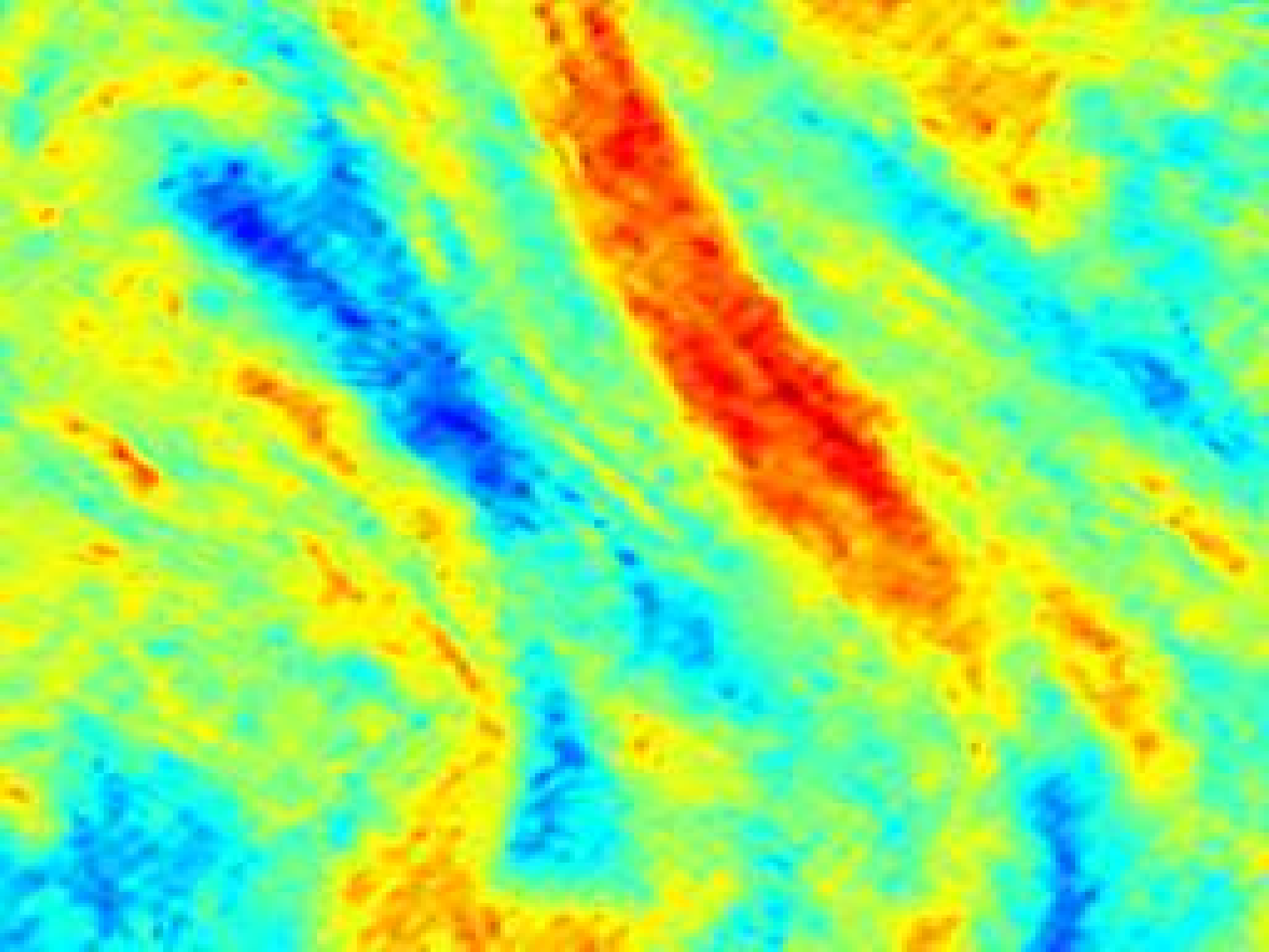}&
    \includegraphics[width = 2.5cm]{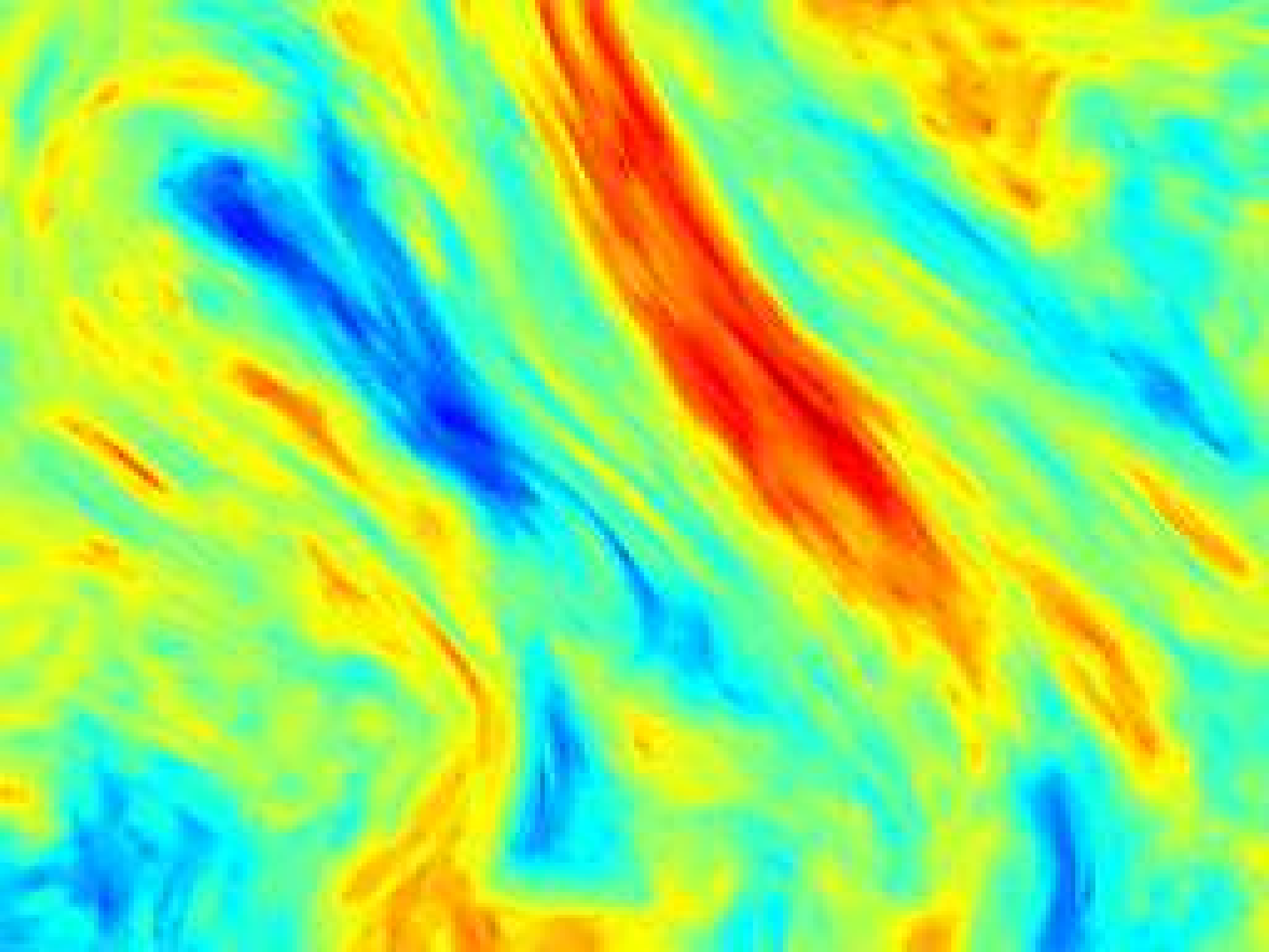}\\
    Subsample & Error & Error & Error\\    \vspace{.5cm}
    \includegraphics[width = 2.5cm]{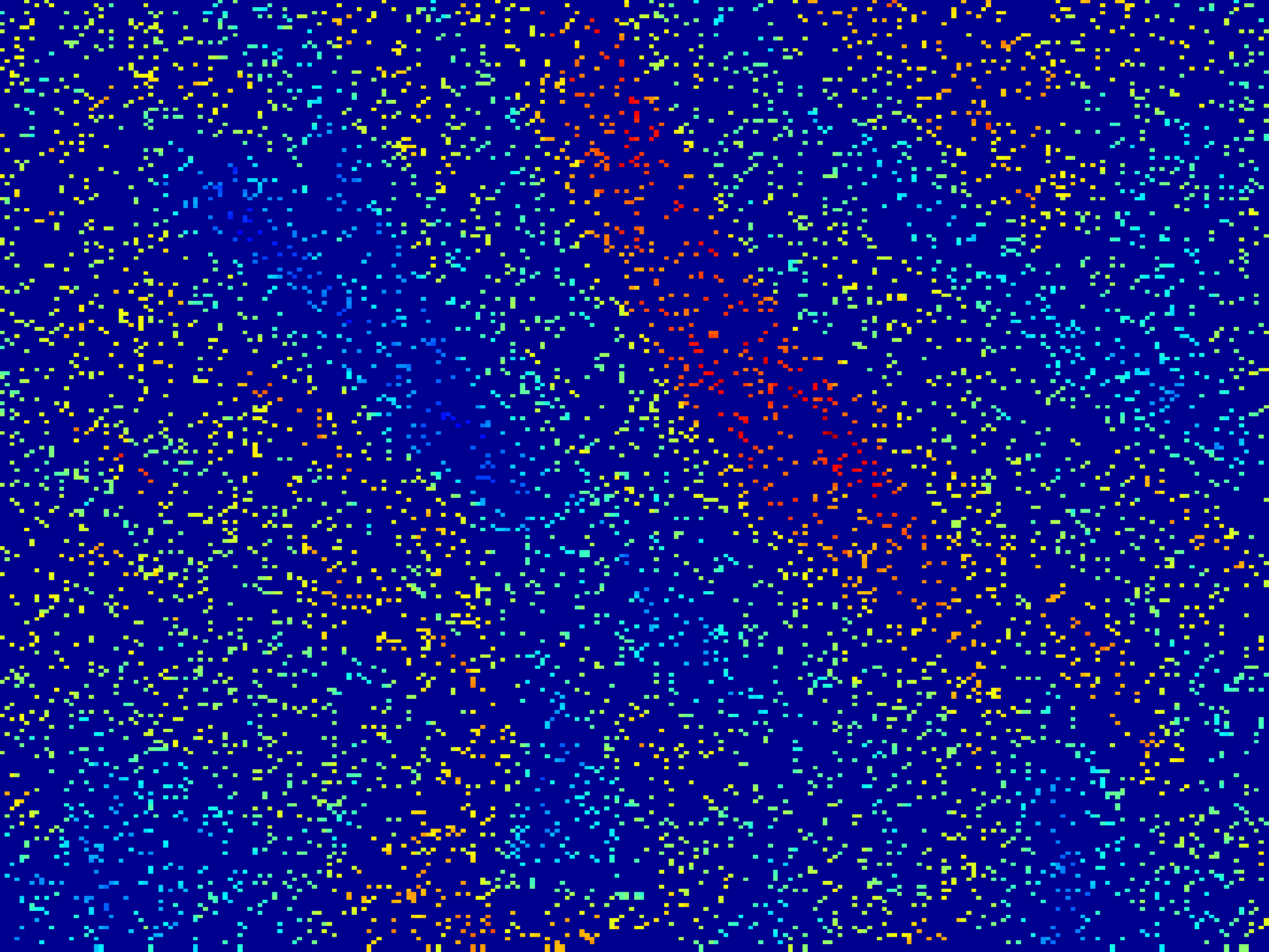}&    
    \includegraphics[width = 2.5cm]{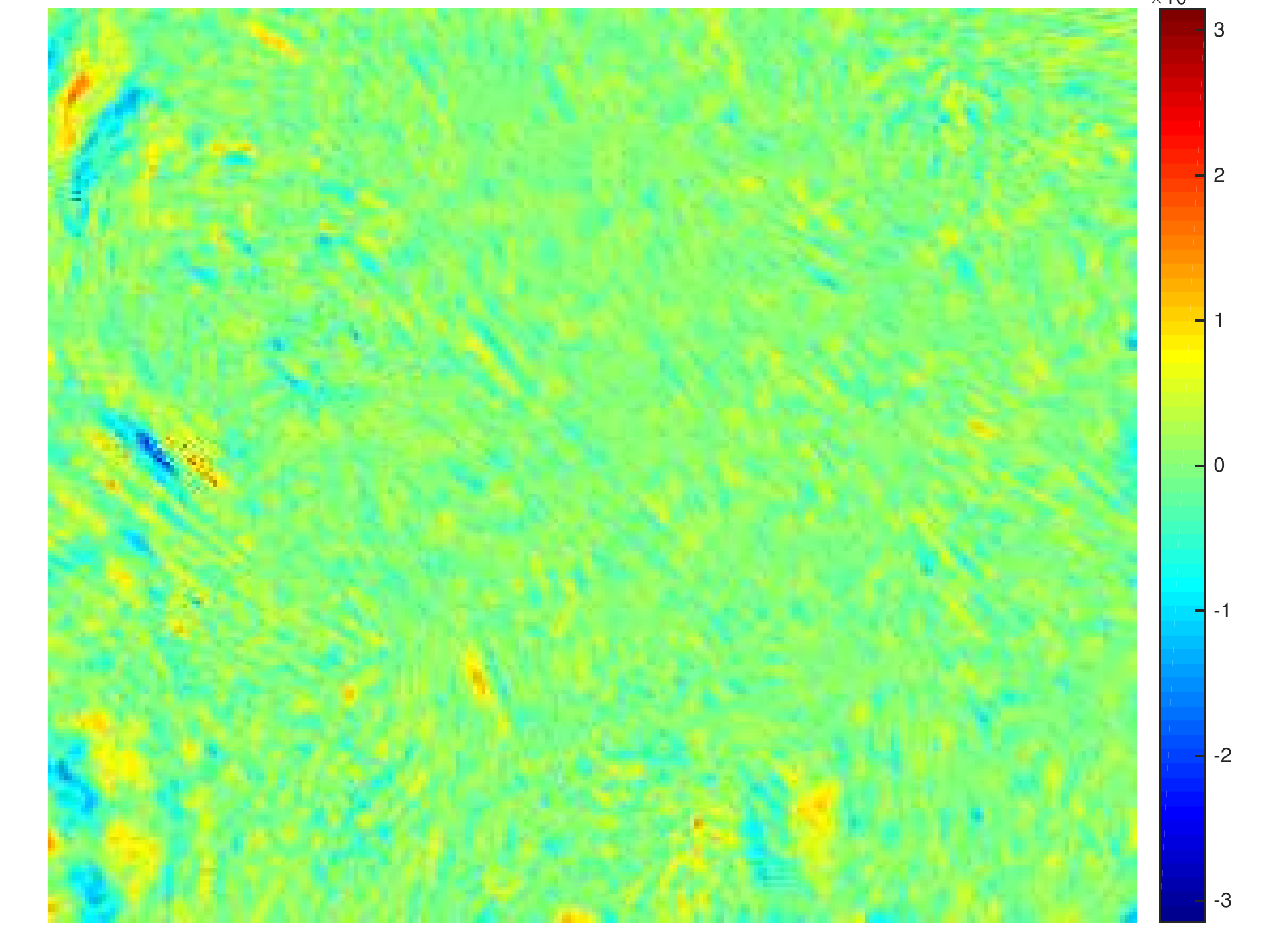}&    
    \includegraphics[width = 2.5cm]{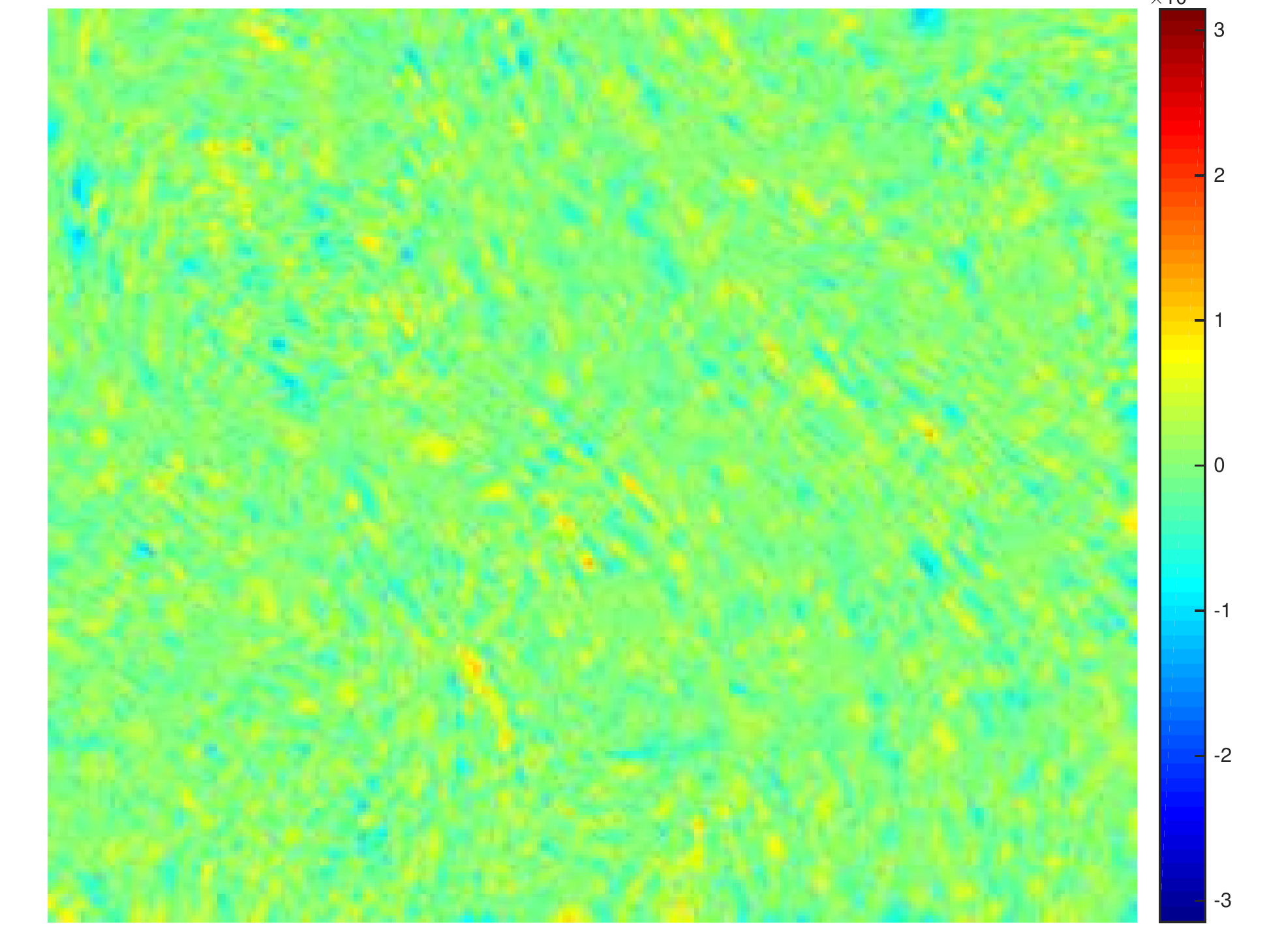}&
    \includegraphics[width = 2.5cm]{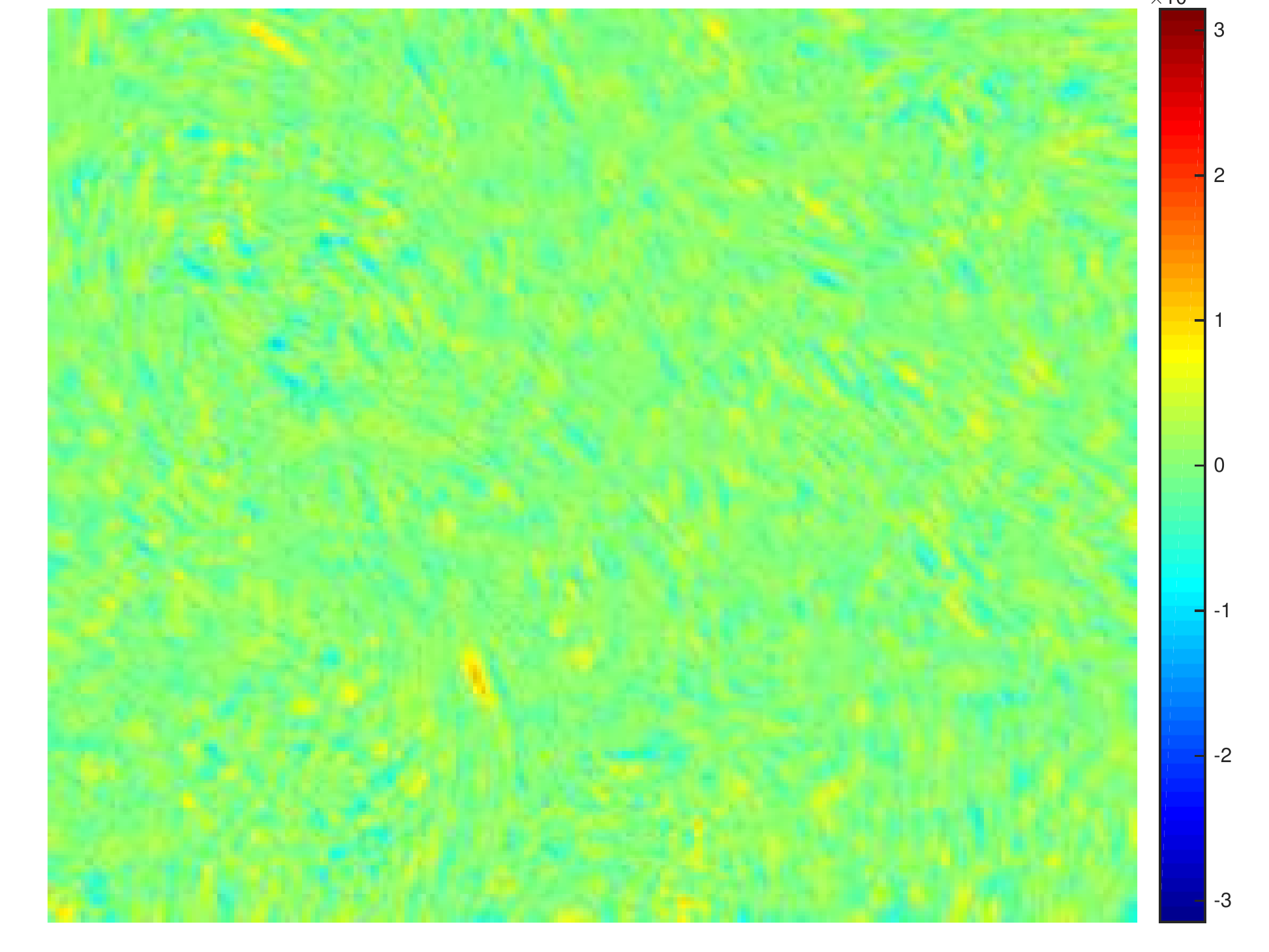}\\
    Original& EBI (43.62dB)& PLE (42.32dB)  & LDMM (\textbf{47.98dB})\\
    \includegraphics[width = 2.5cm]{latticebig_2d_original}&
    \includegraphics[width = 2.5cm]{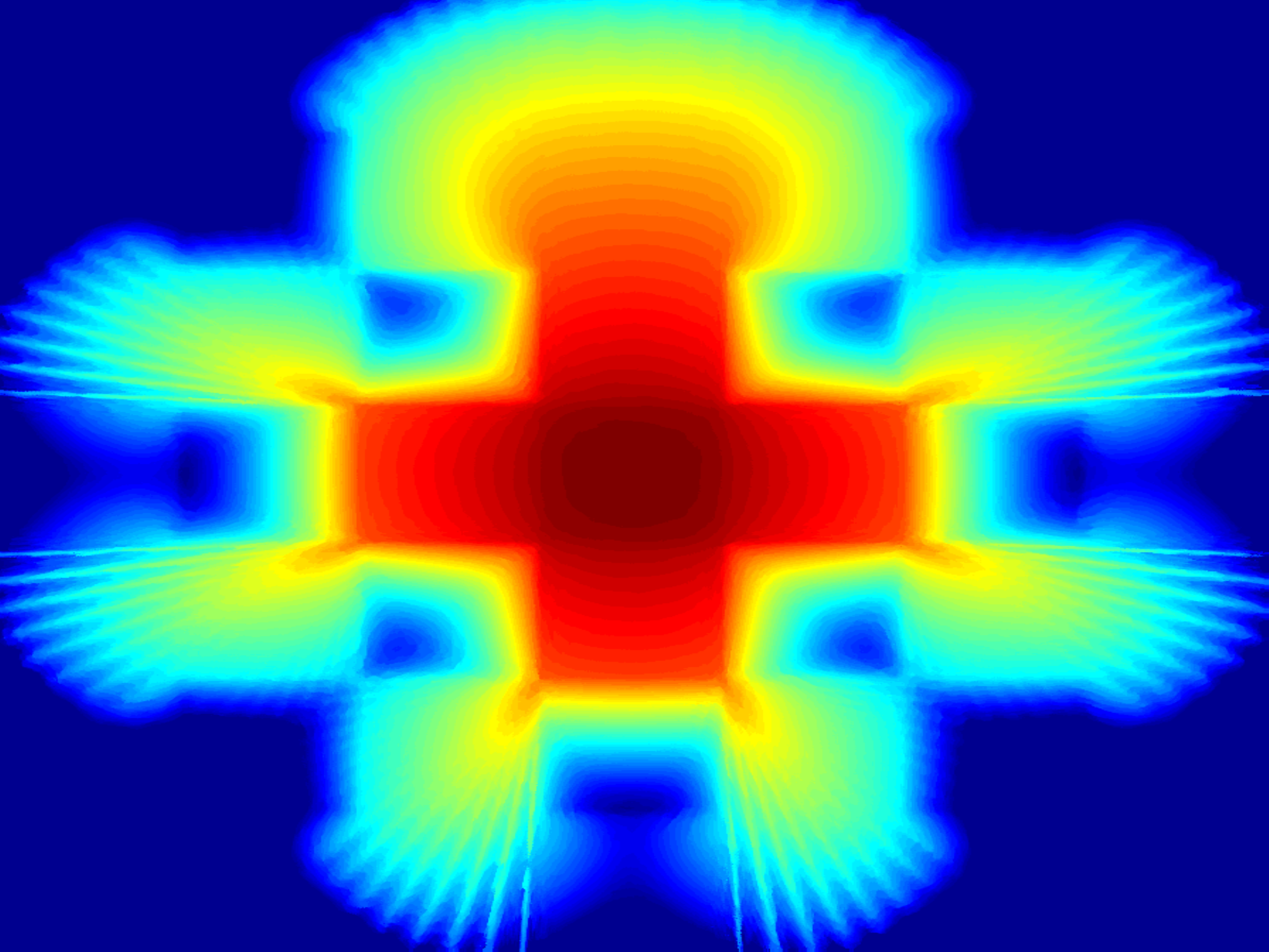}&
    \includegraphics[width = 2.5cm]{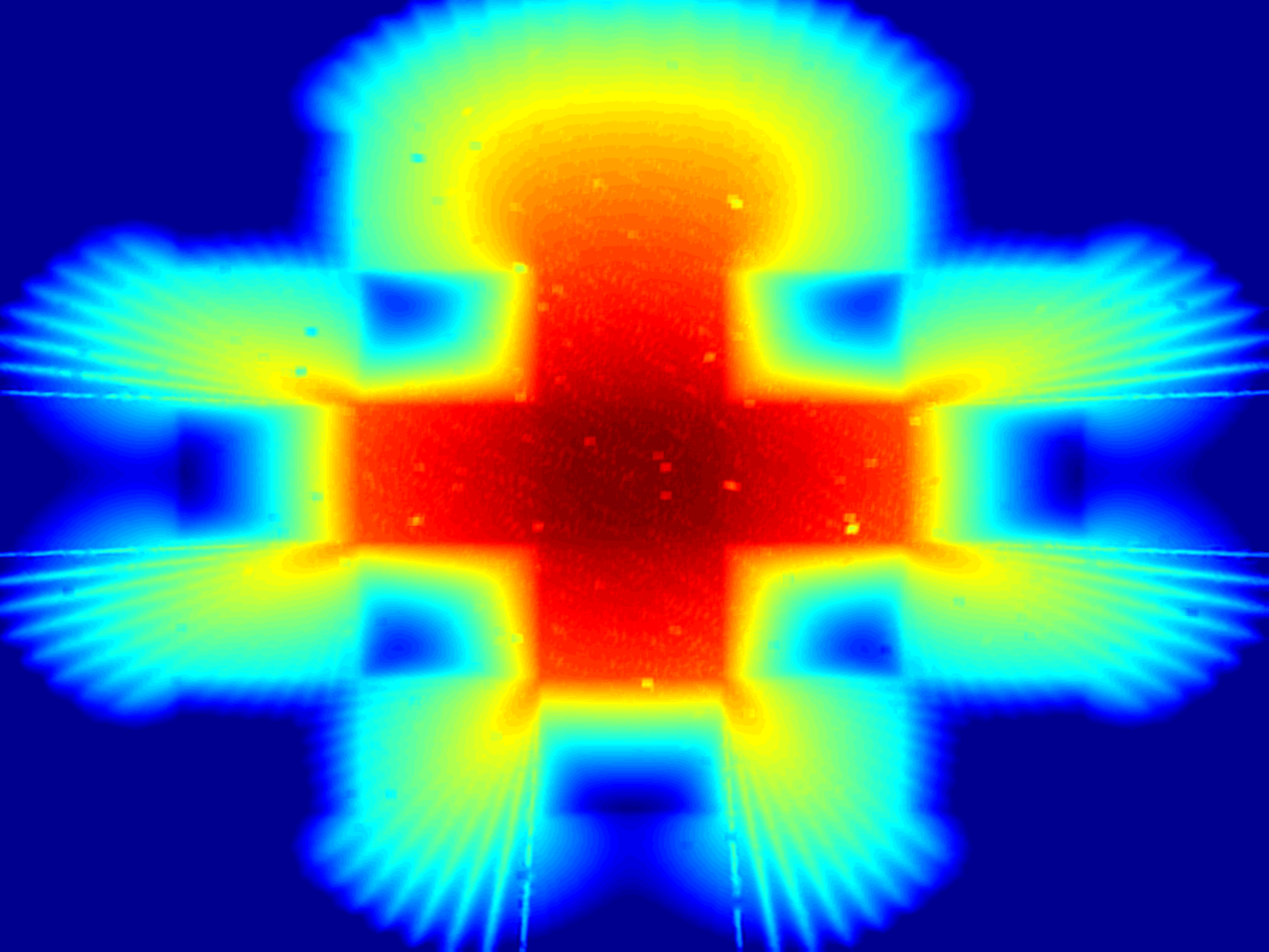}&
    \includegraphics[width = 2.5cm]{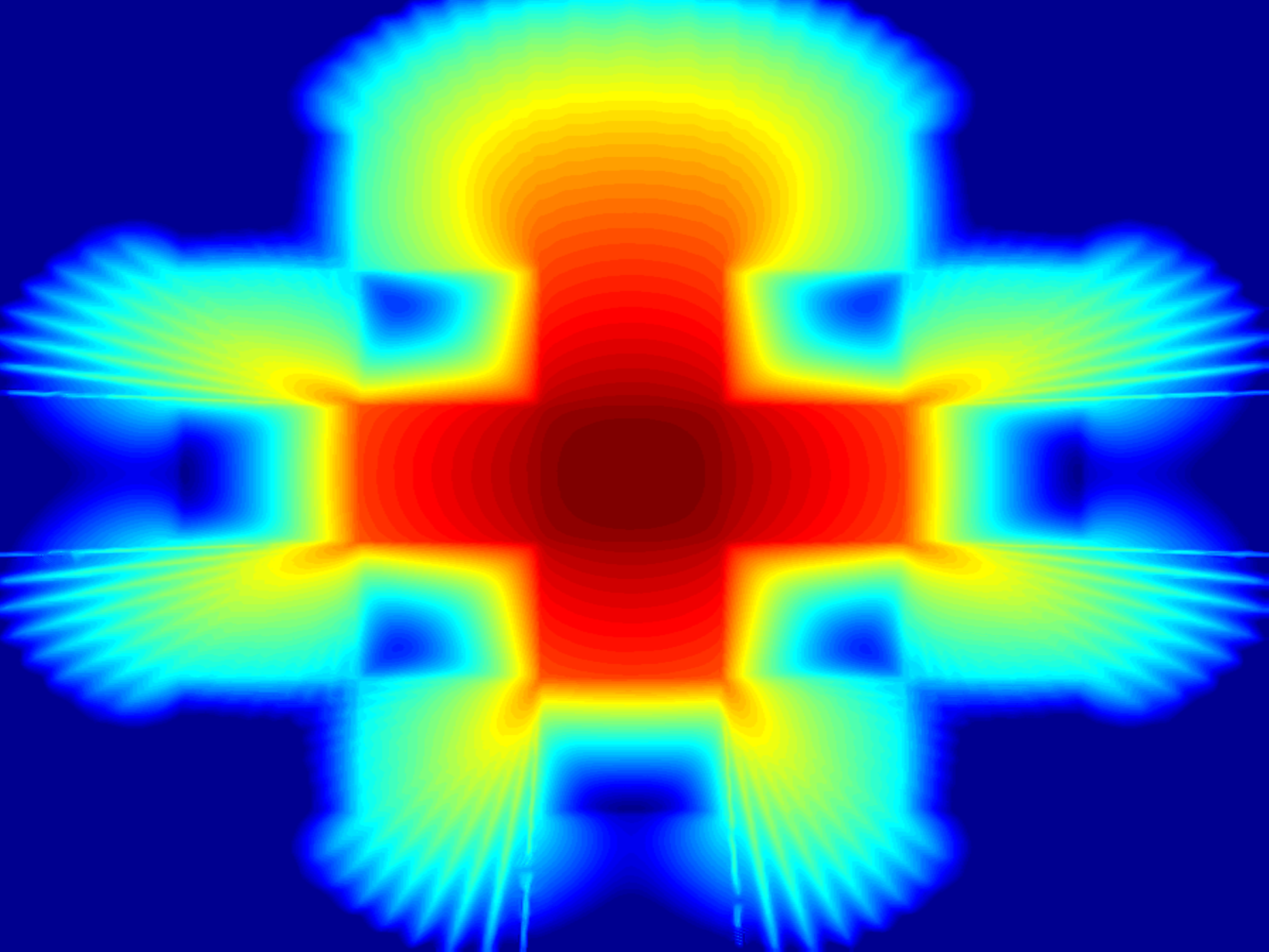}\\
    Subsample & Error & Error & Error\\
    \includegraphics[width = 2.5cm]{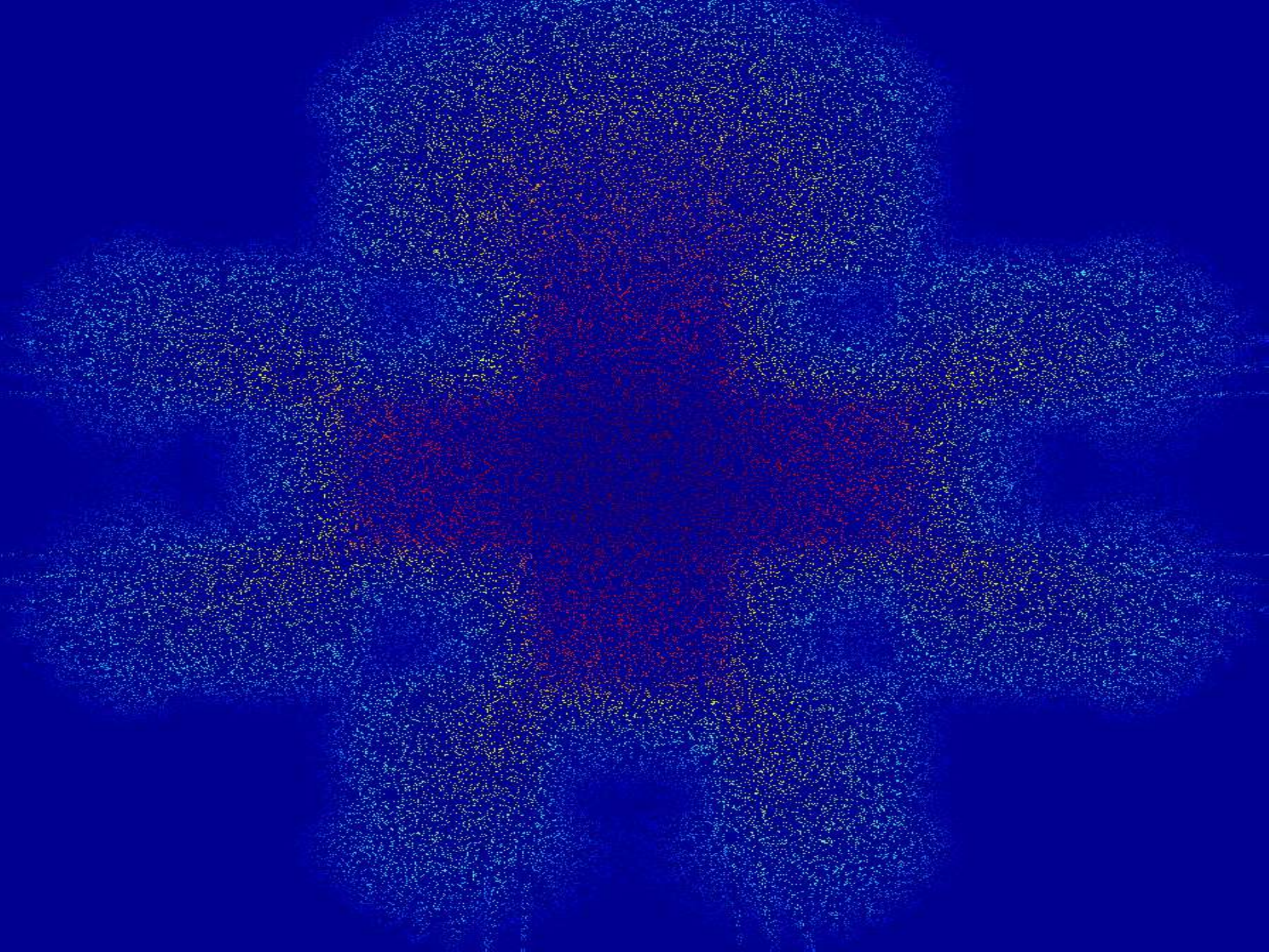}&    
    \includegraphics[width = 2.5cm]{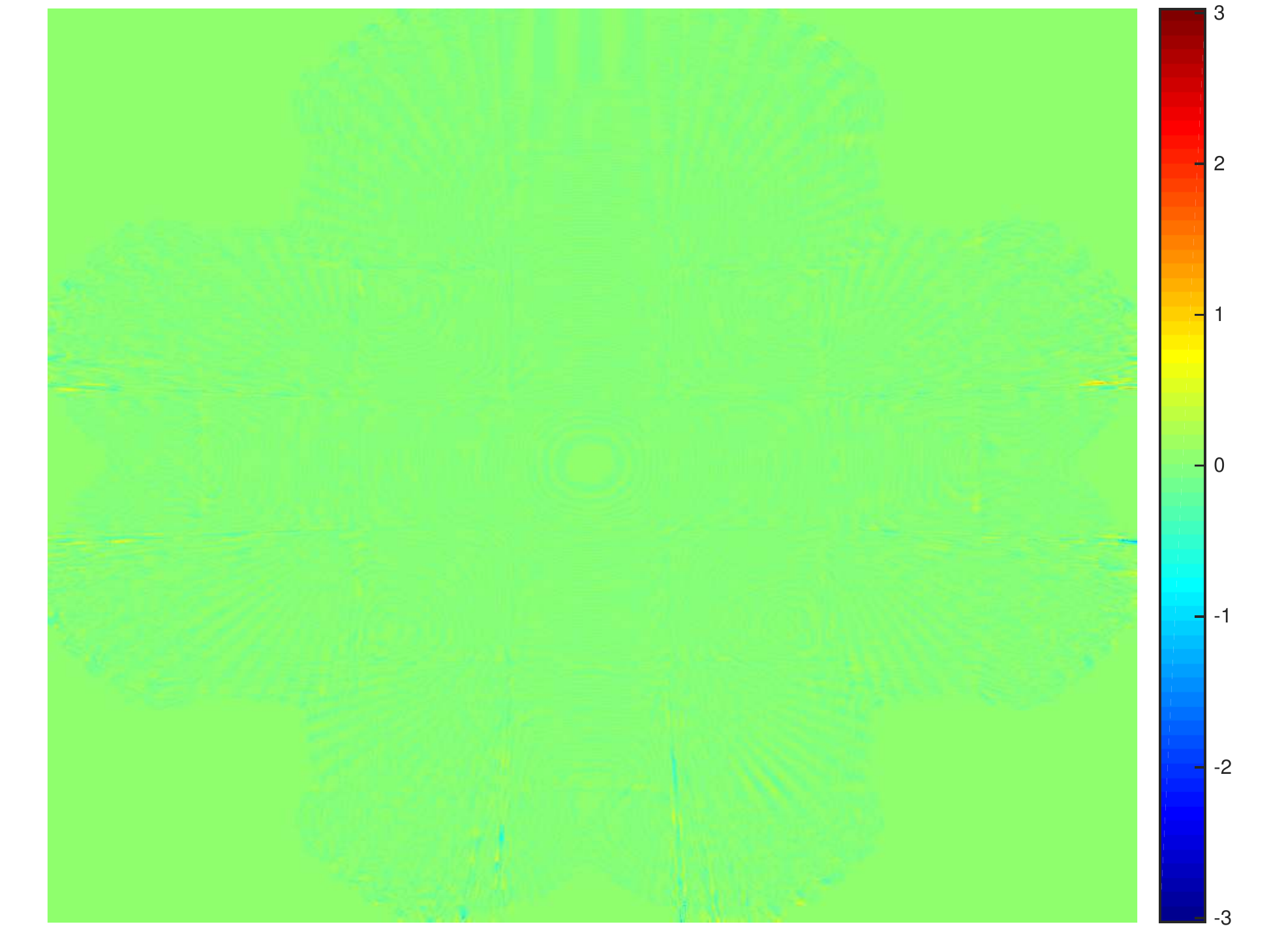}&    
    \includegraphics[width = 2.5cm]{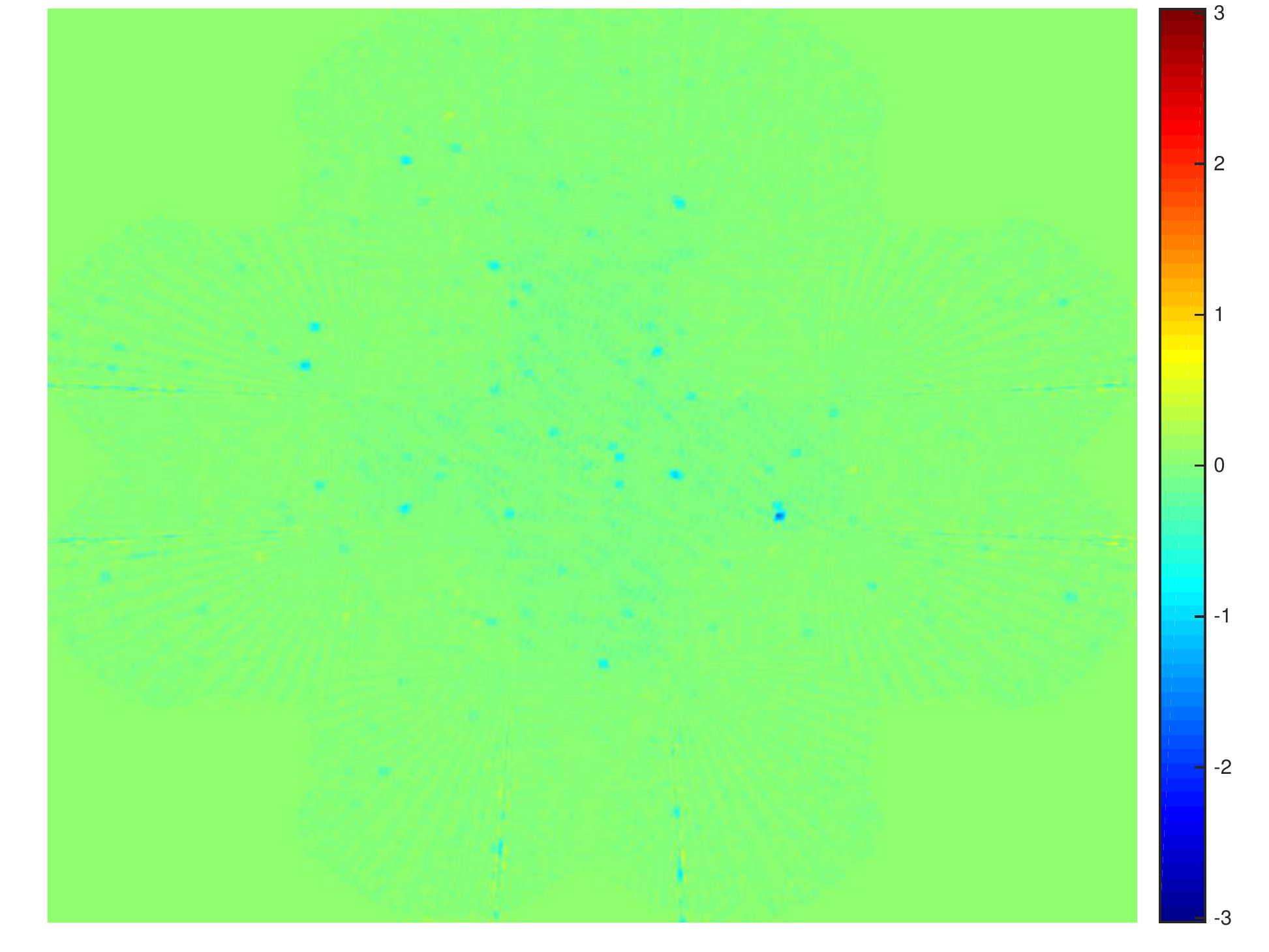}&
    \includegraphics[width = 2.5cm]{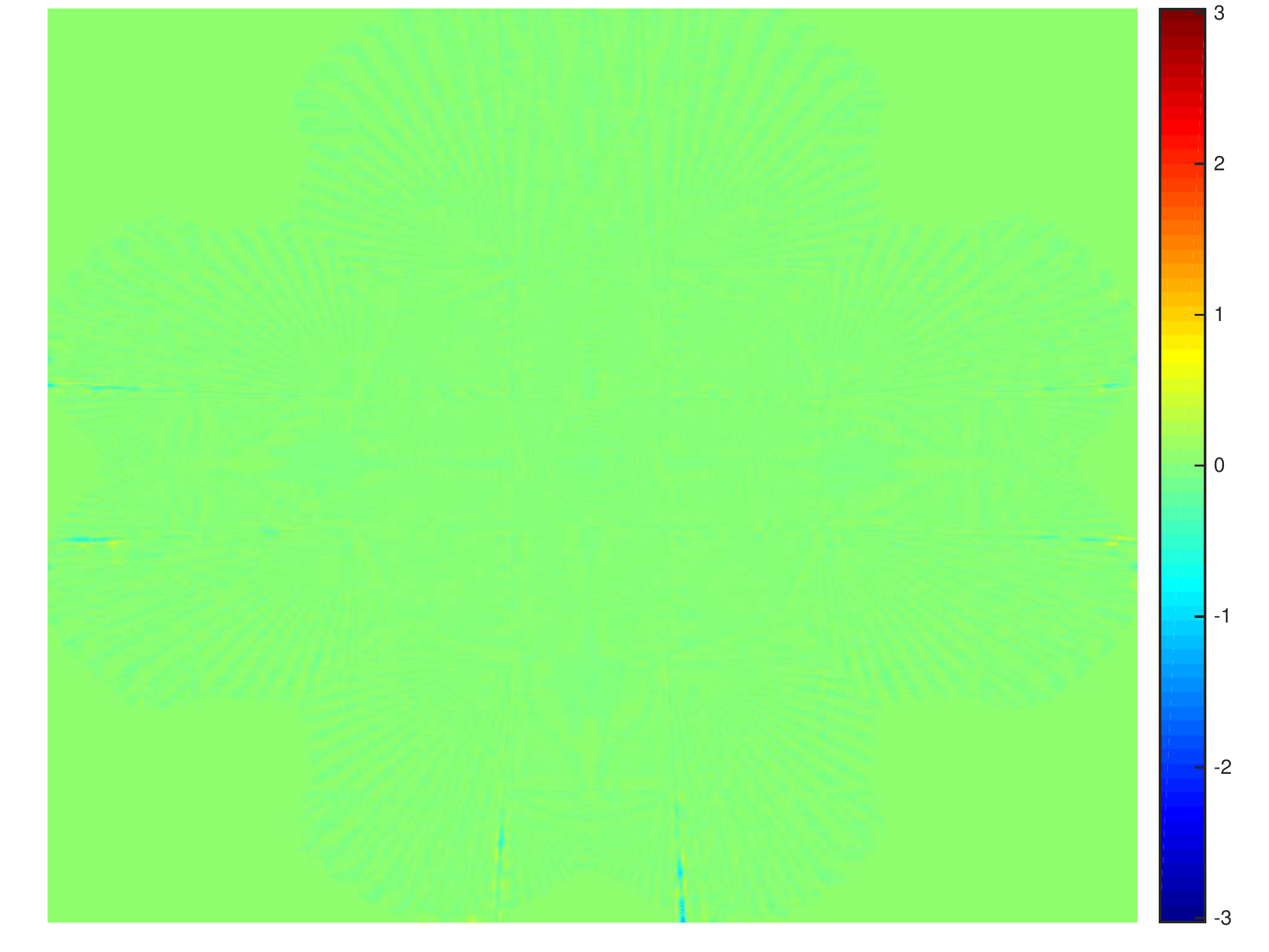}
  \end{tabular}
  \caption{Interpolation of 2D scientific data sets from $10\%$ random sampling. The figures in the first column are the original and subsampled data. The figures in the other three columns are the results and errors of the competing algorithms.}
  \label{fig:result_random_2d_10p}
\end{figure}

\begin{figure}[H]
  \centering
  \begin{tabular}{cccc}
    Original& EBI (31.36dB)& PLE (24.84dB)  & LDMM (\textbf{39.09dB})\\
    \includegraphics[width = 2.5cm]{antonio_2d_original}&
    \includegraphics[width = 2.5cm]{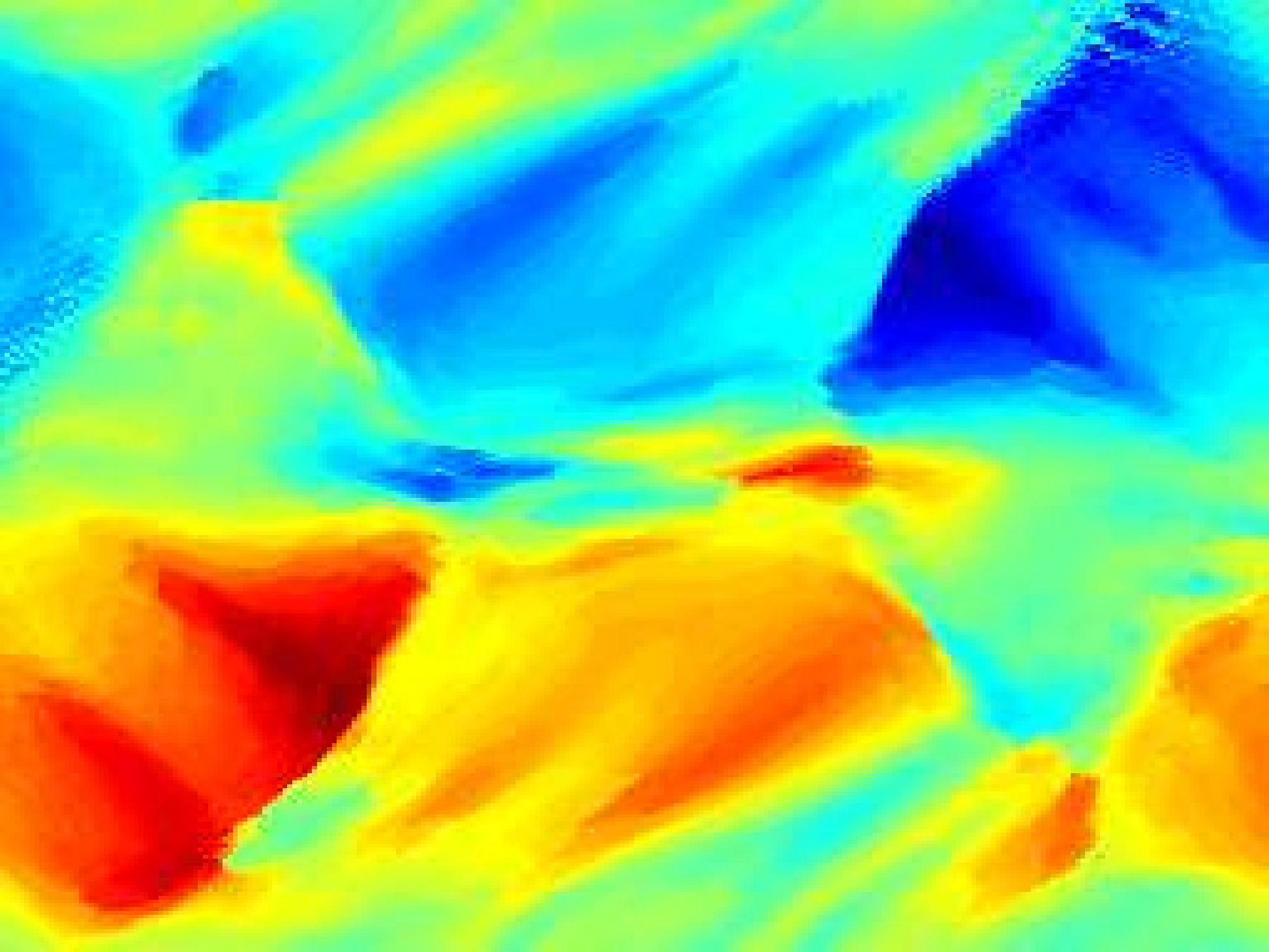}&
    \includegraphics[width = 2.5cm]{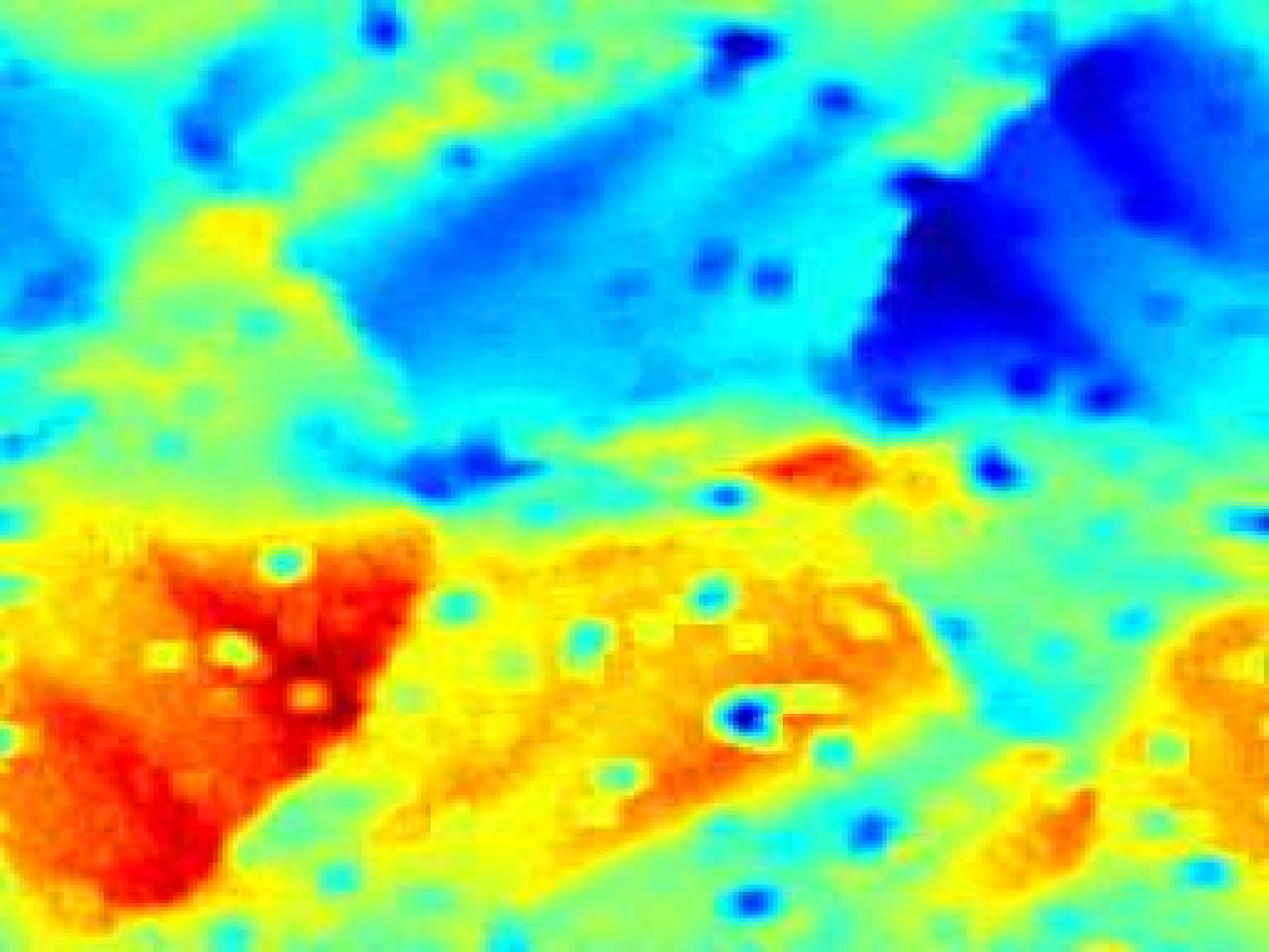}&
    \includegraphics[width = 2.5cm]{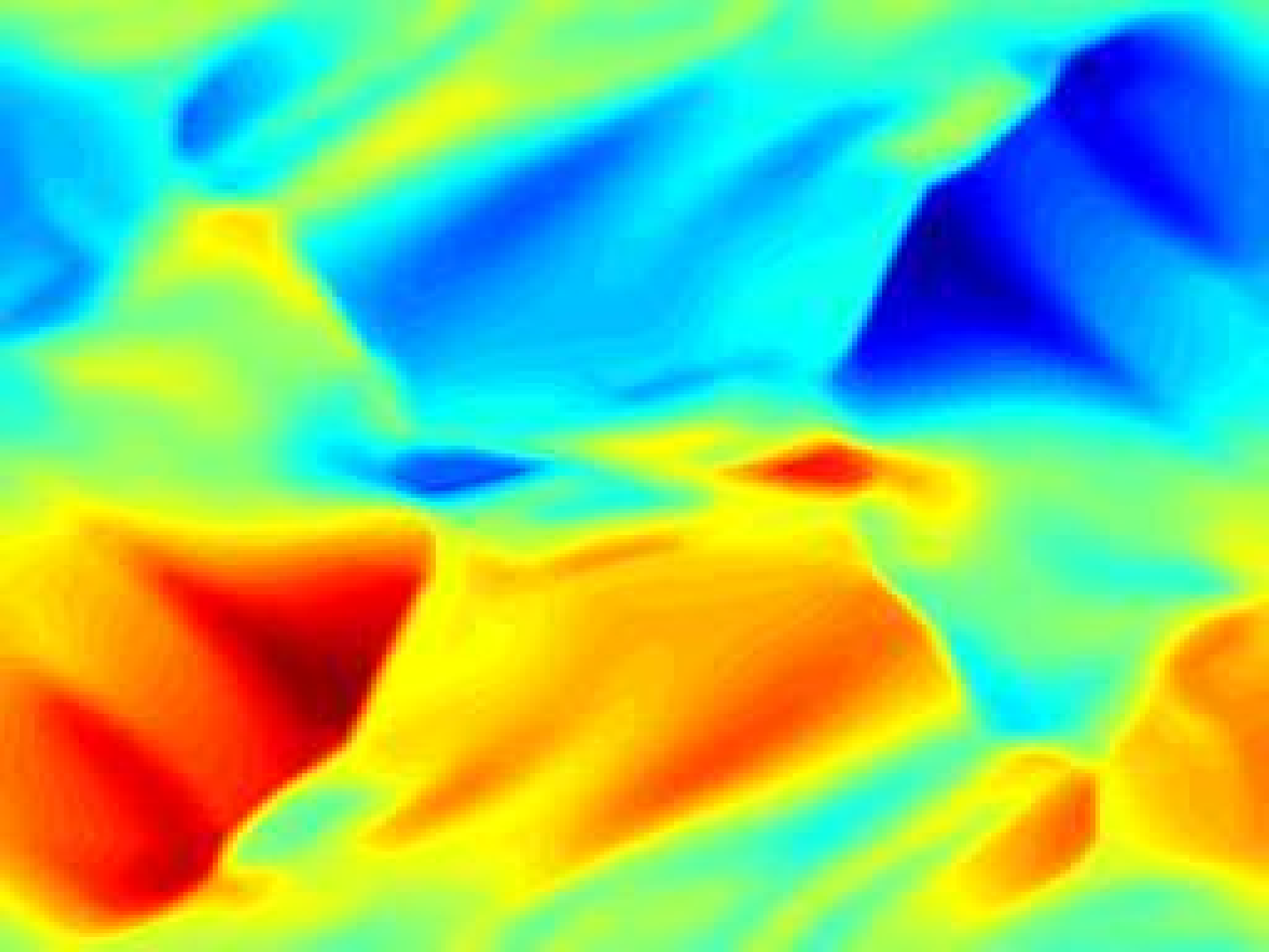}\\
    Subsample & Error & Error & Error\\     \vspace{.5cm}
    \includegraphics[width = 2.5cm]{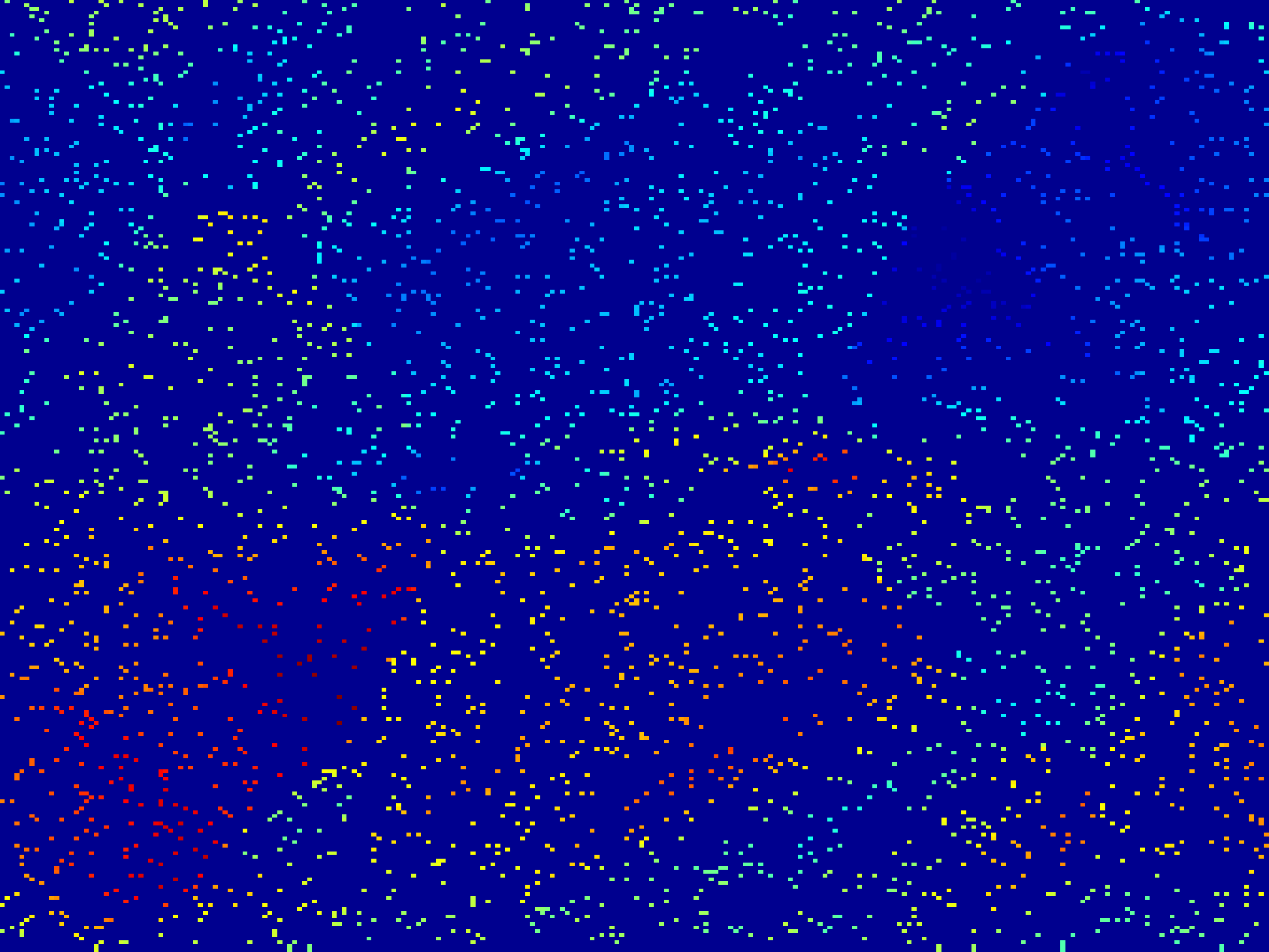}&    
    \includegraphics[width = 2.5cm]{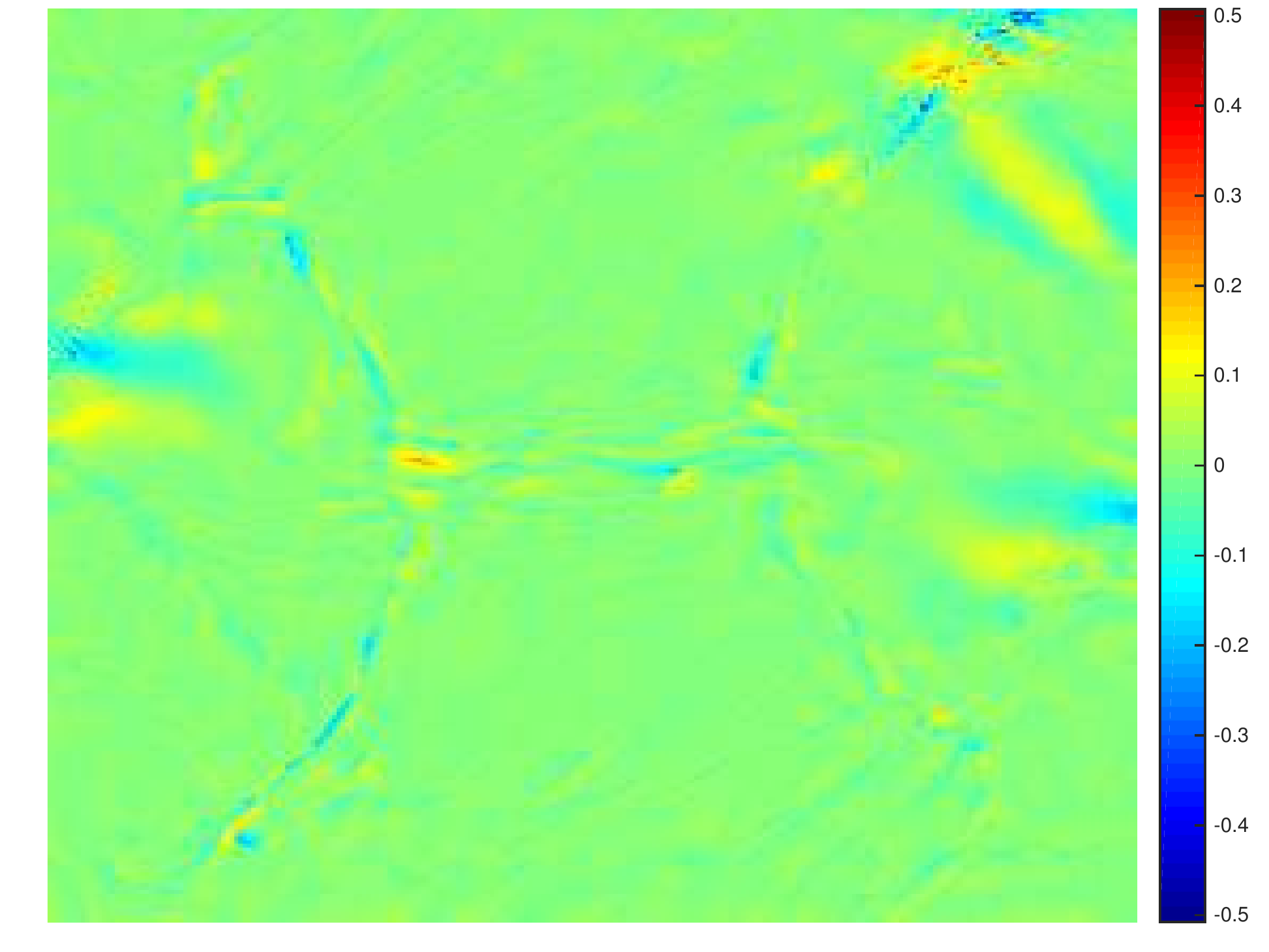}&    
    \includegraphics[width = 2.5cm]{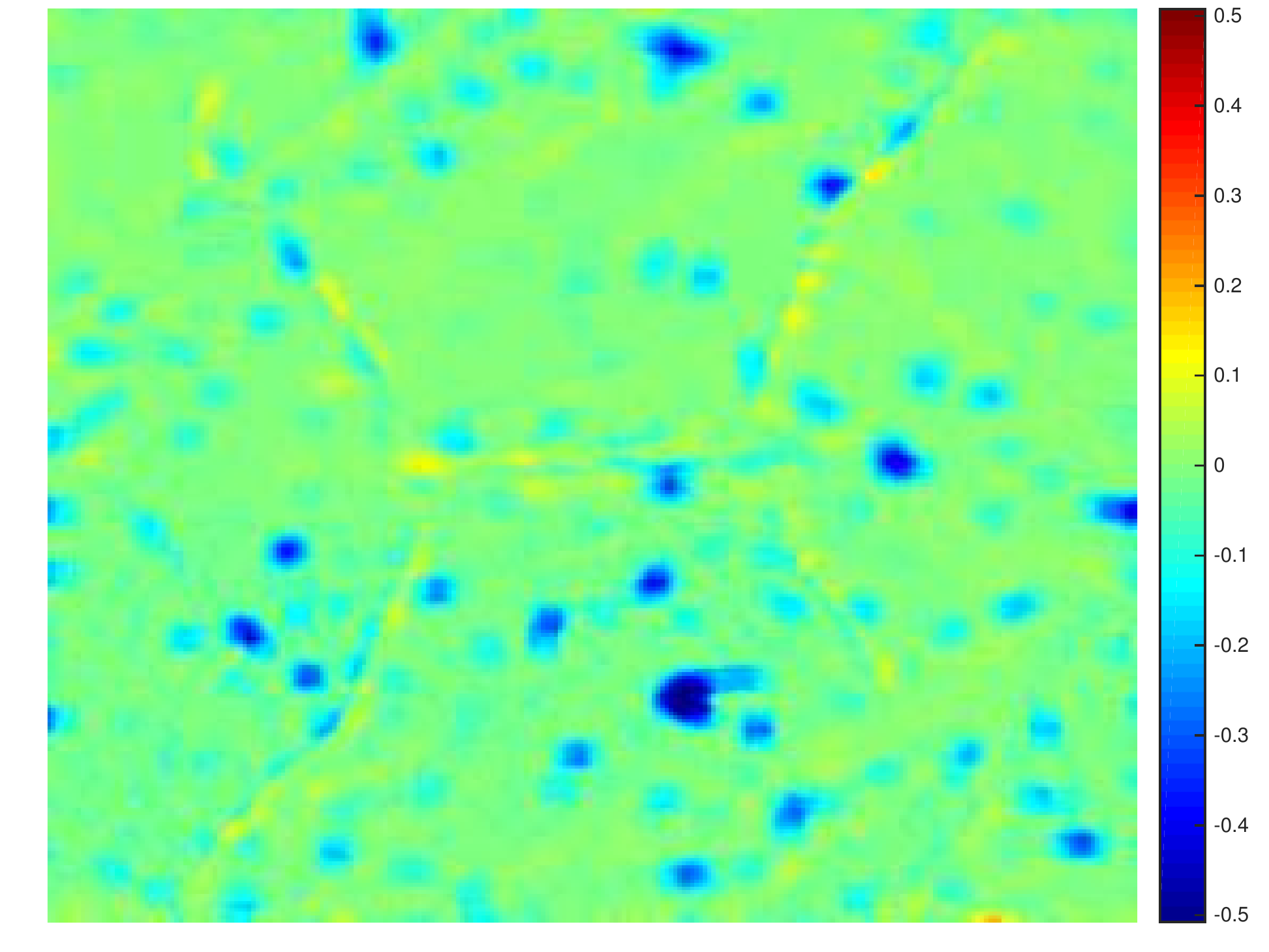}&
    \includegraphics[width = 2.5cm]{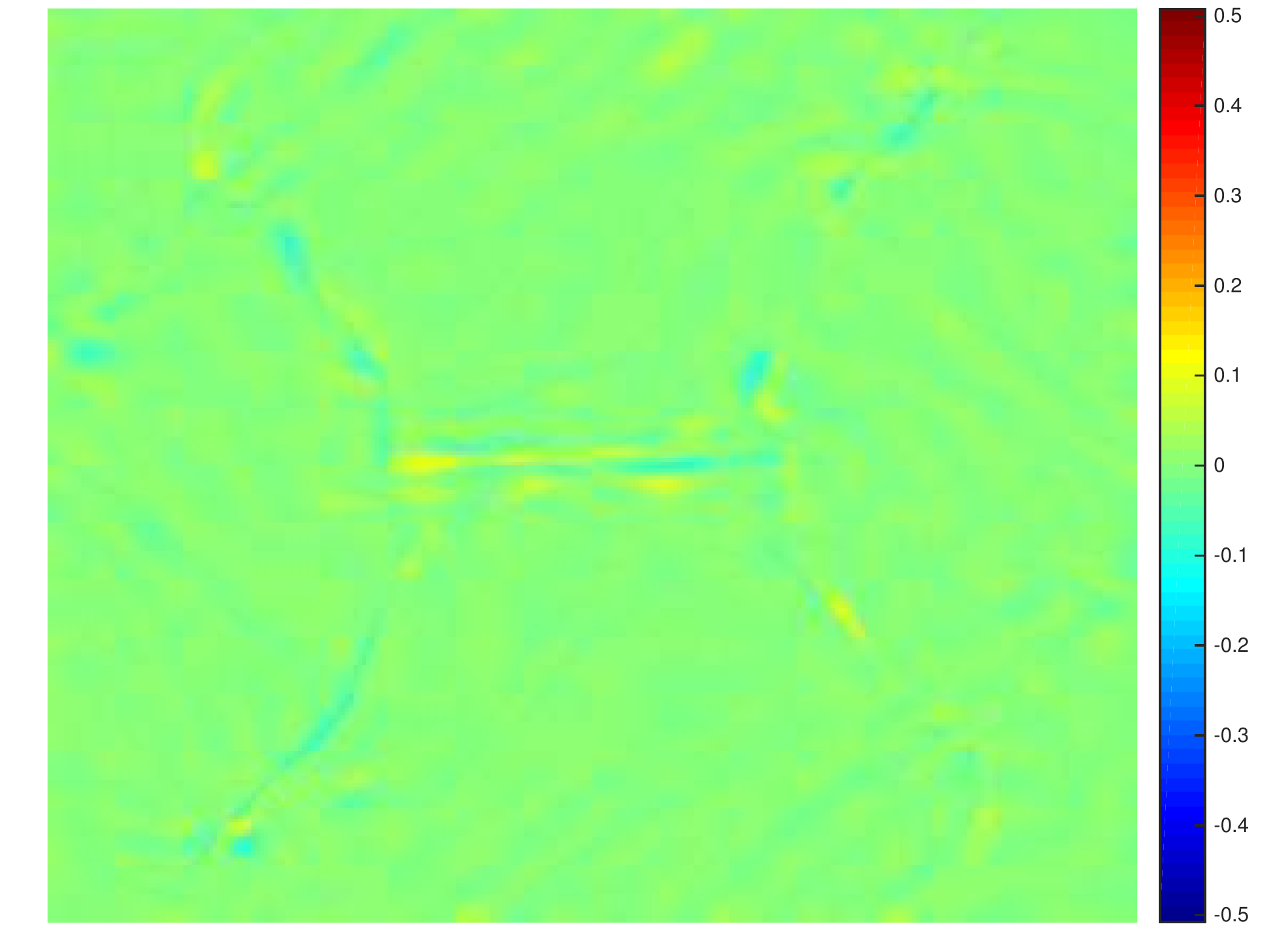}\\
    Original& EBI (25.65dB)& PLE (21.88dB)  & LDMM (\textbf{27.93dB})\\
    \includegraphics[width = 2.5cm]{shock_2d_original}&
    \includegraphics[width = 2.5cm]{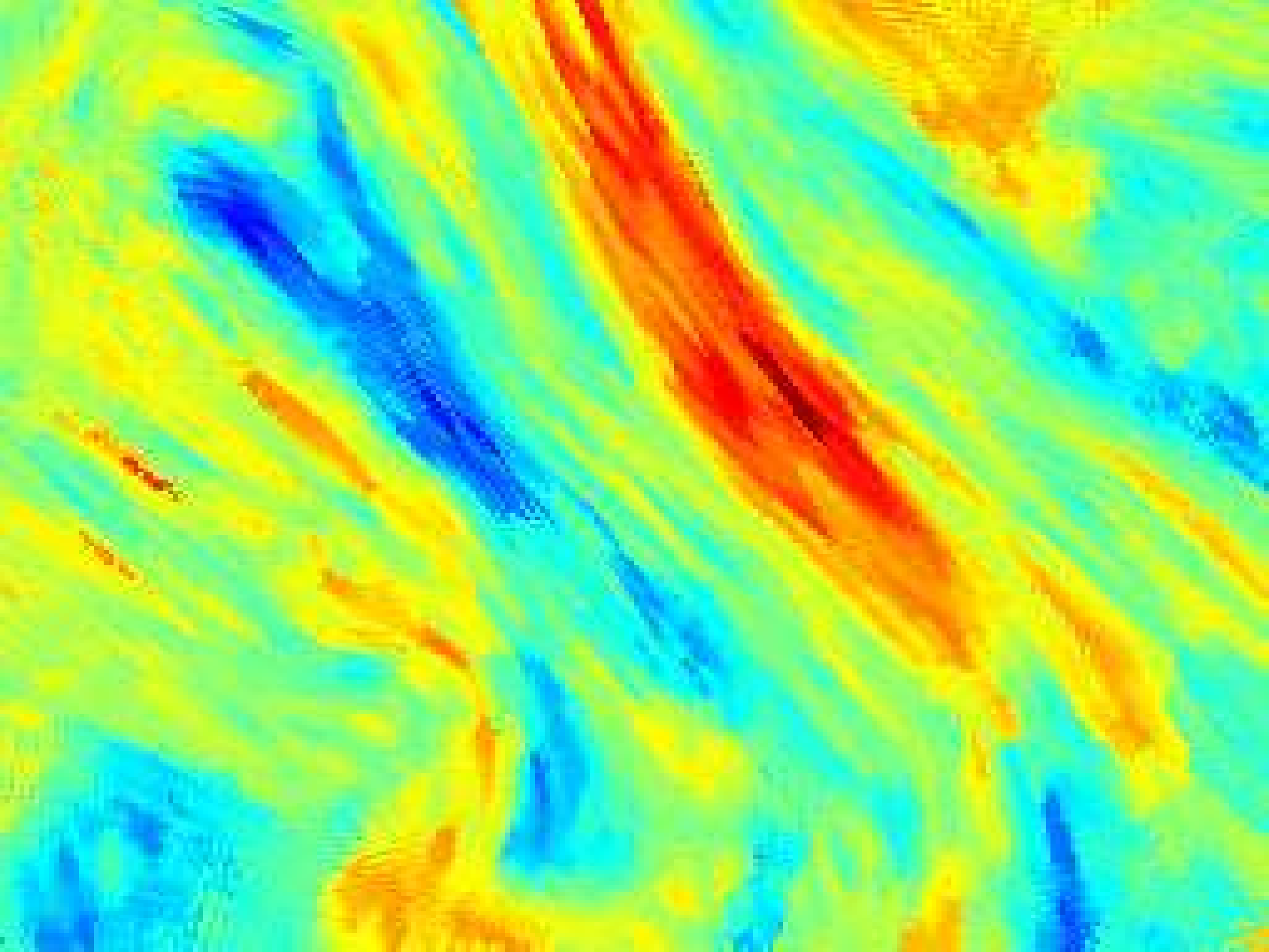}&
    \includegraphics[width = 2.5cm]{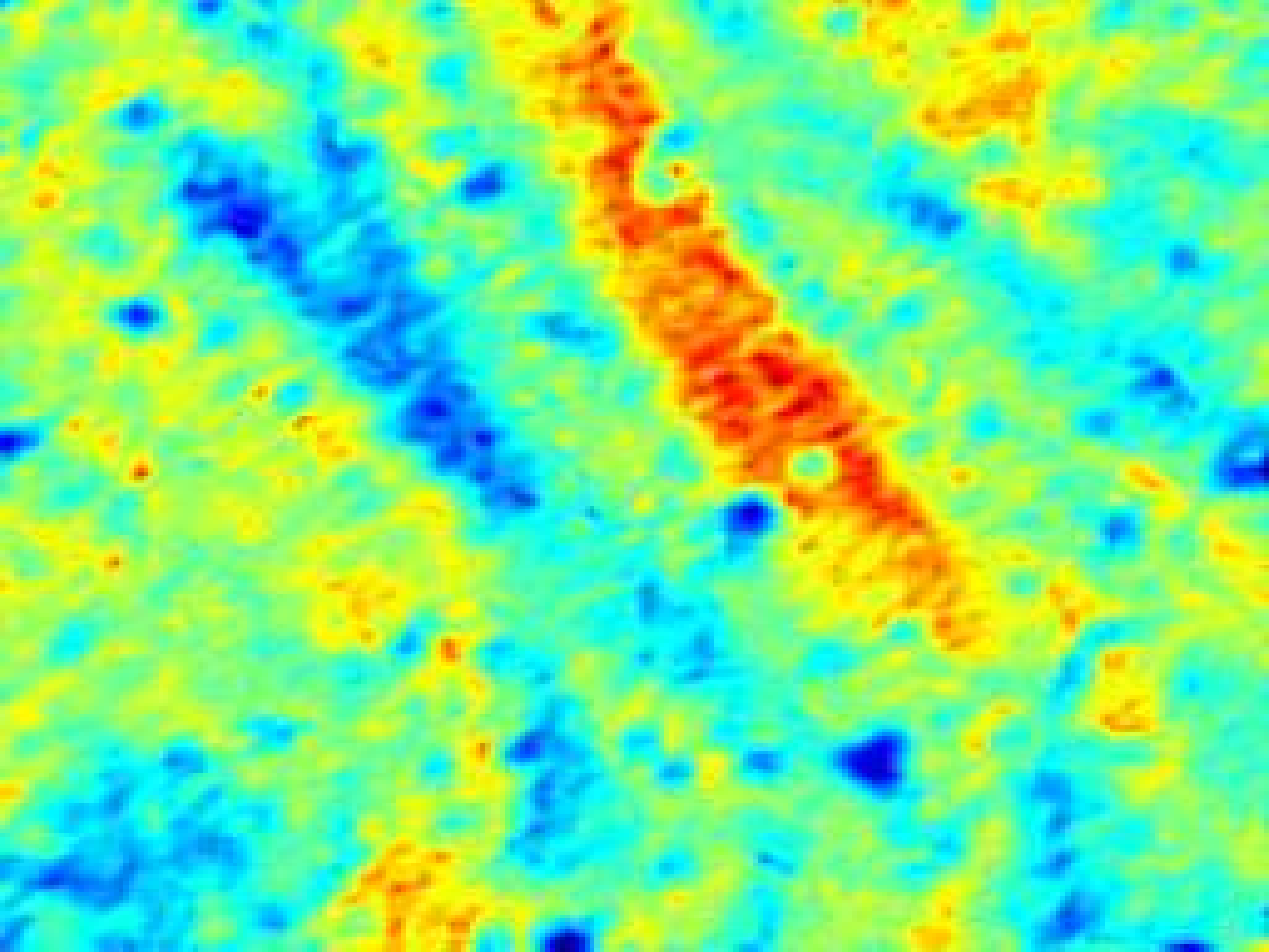}&
    \includegraphics[width = 2.5cm]{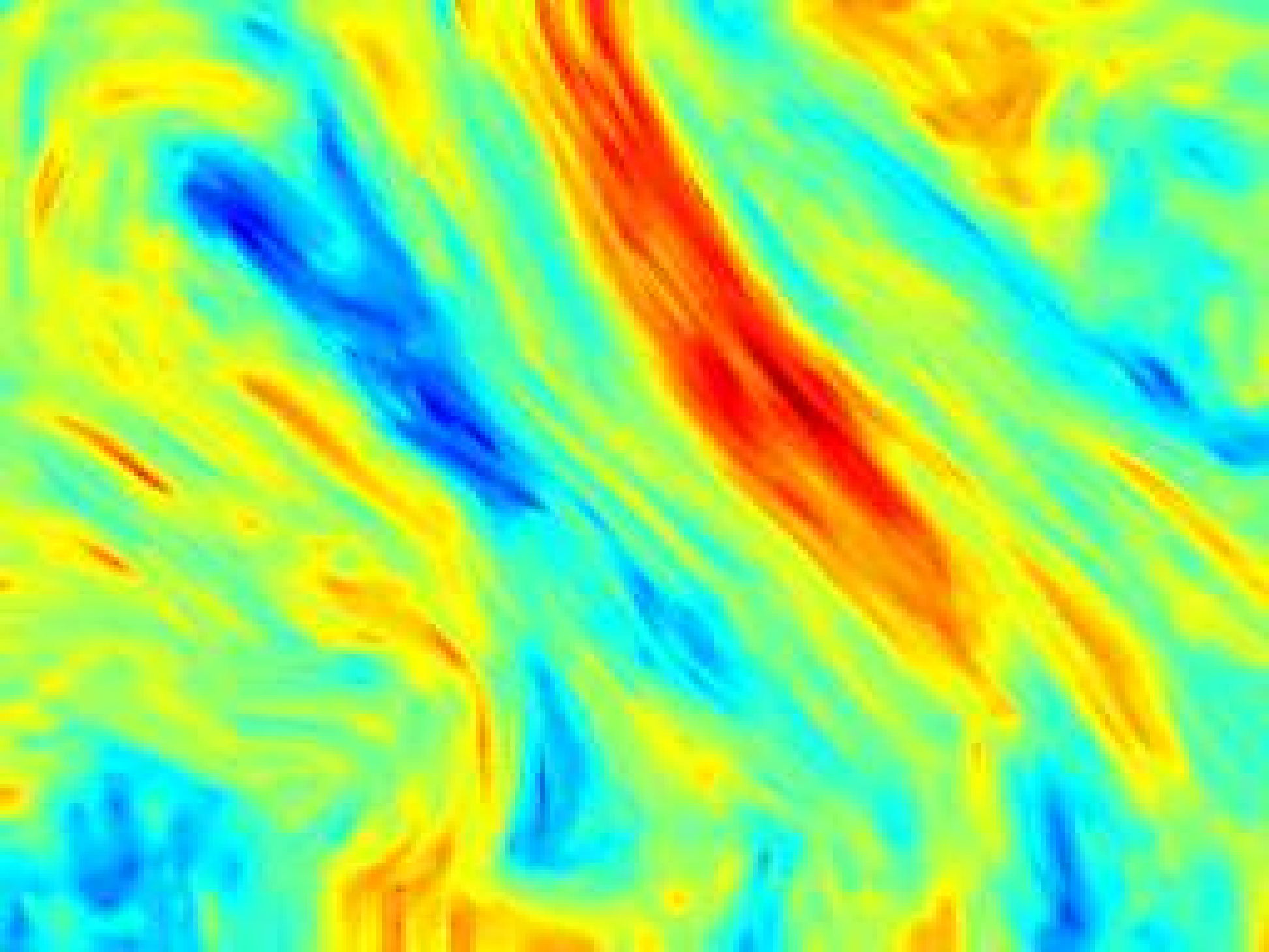}\\
    Subsample & Error & Error & Error\\    \vspace{.5cm}
    \includegraphics[width = 2.5cm]{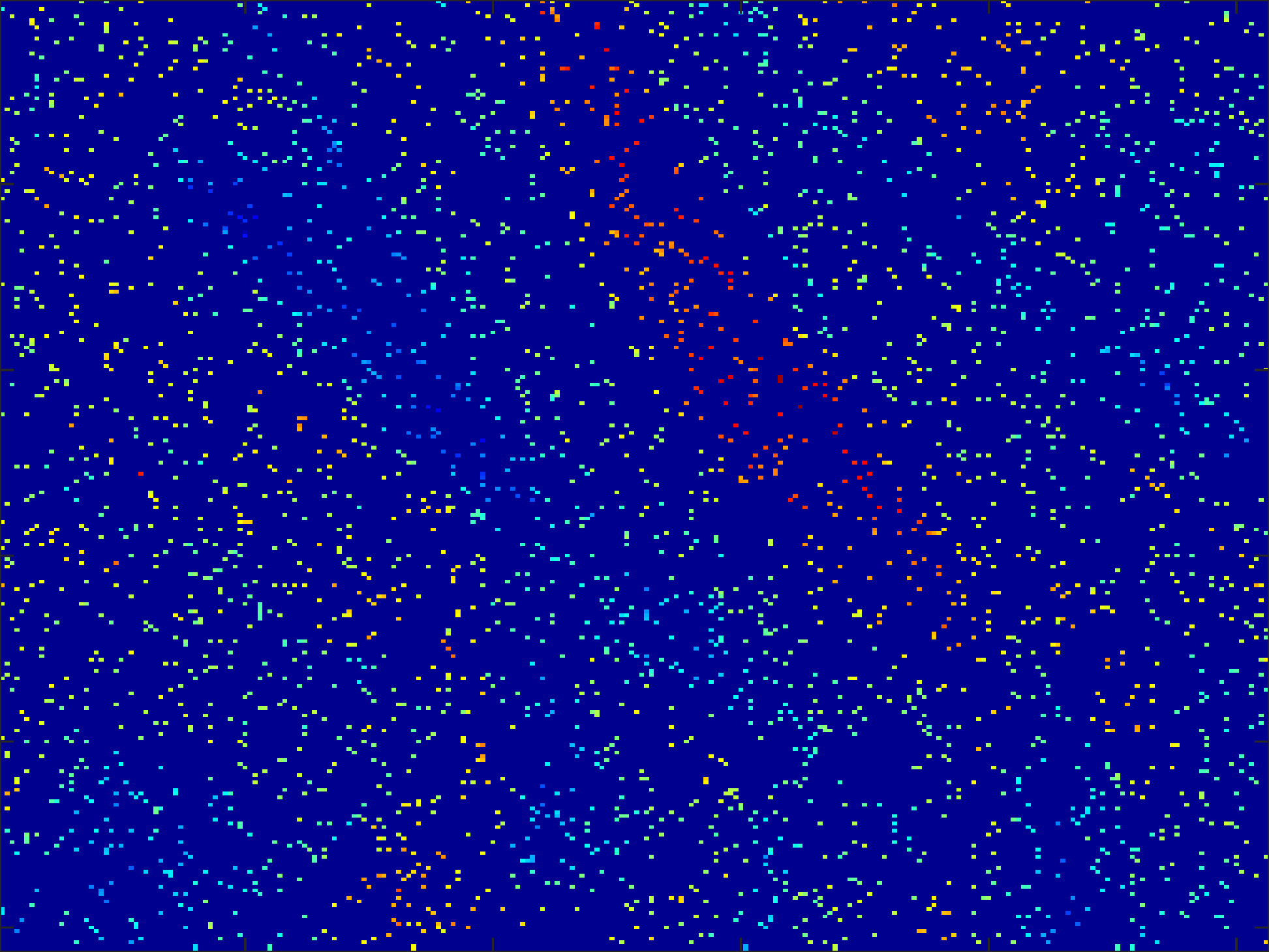}&    
    \includegraphics[width = 2.5cm]{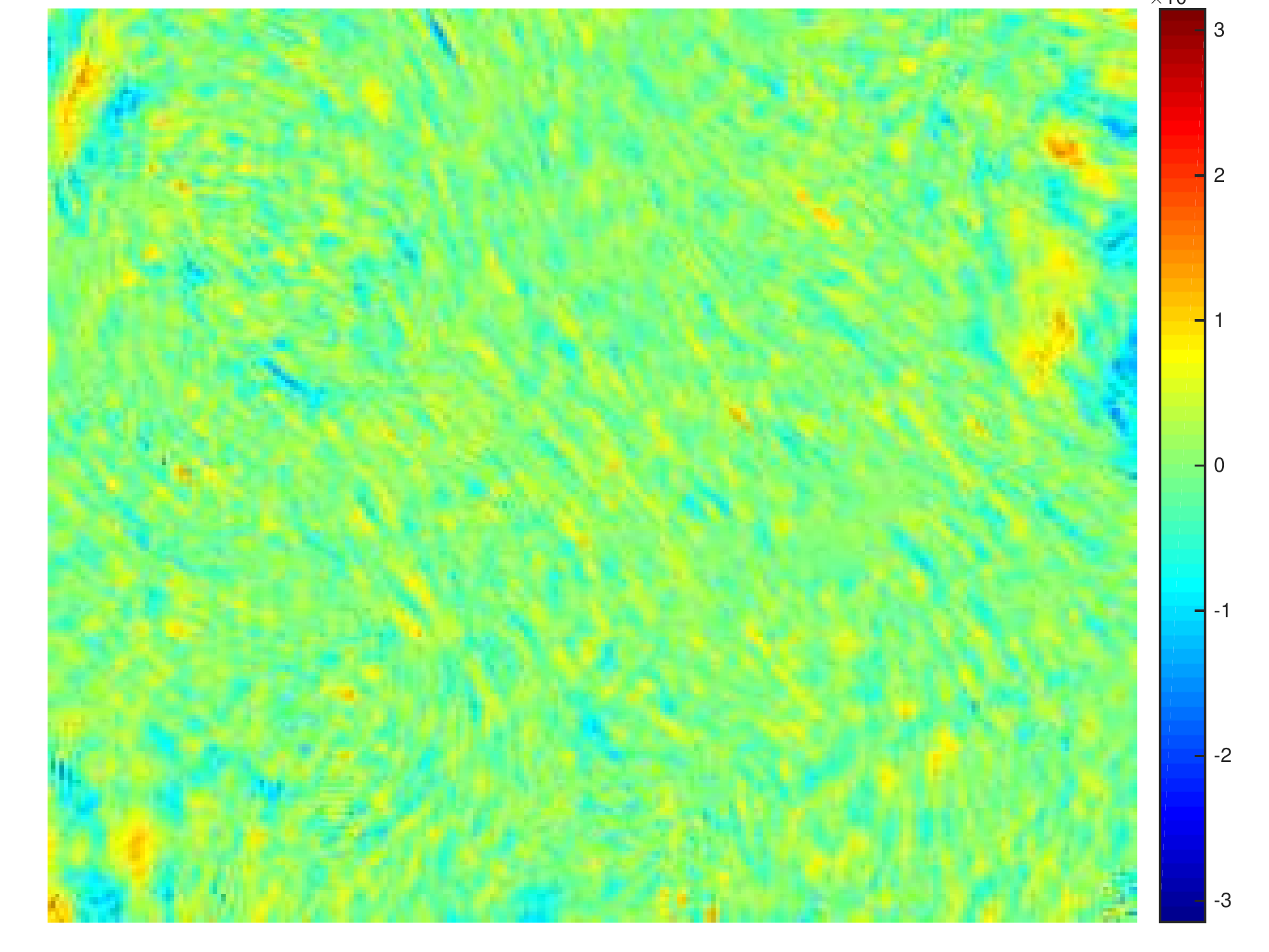}&    
    \includegraphics[width = 2.5cm]{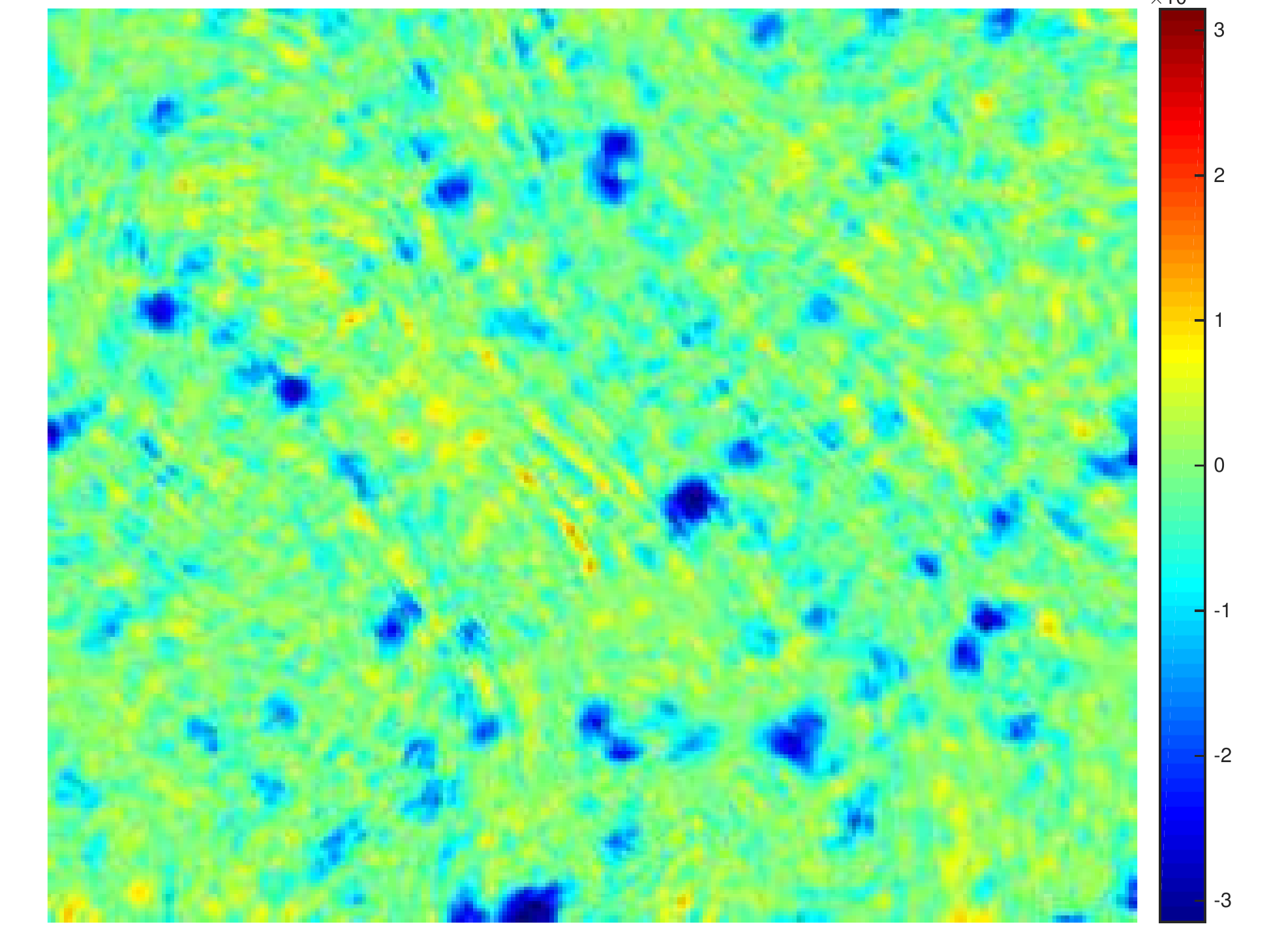}&
    \includegraphics[width = 2.5cm]{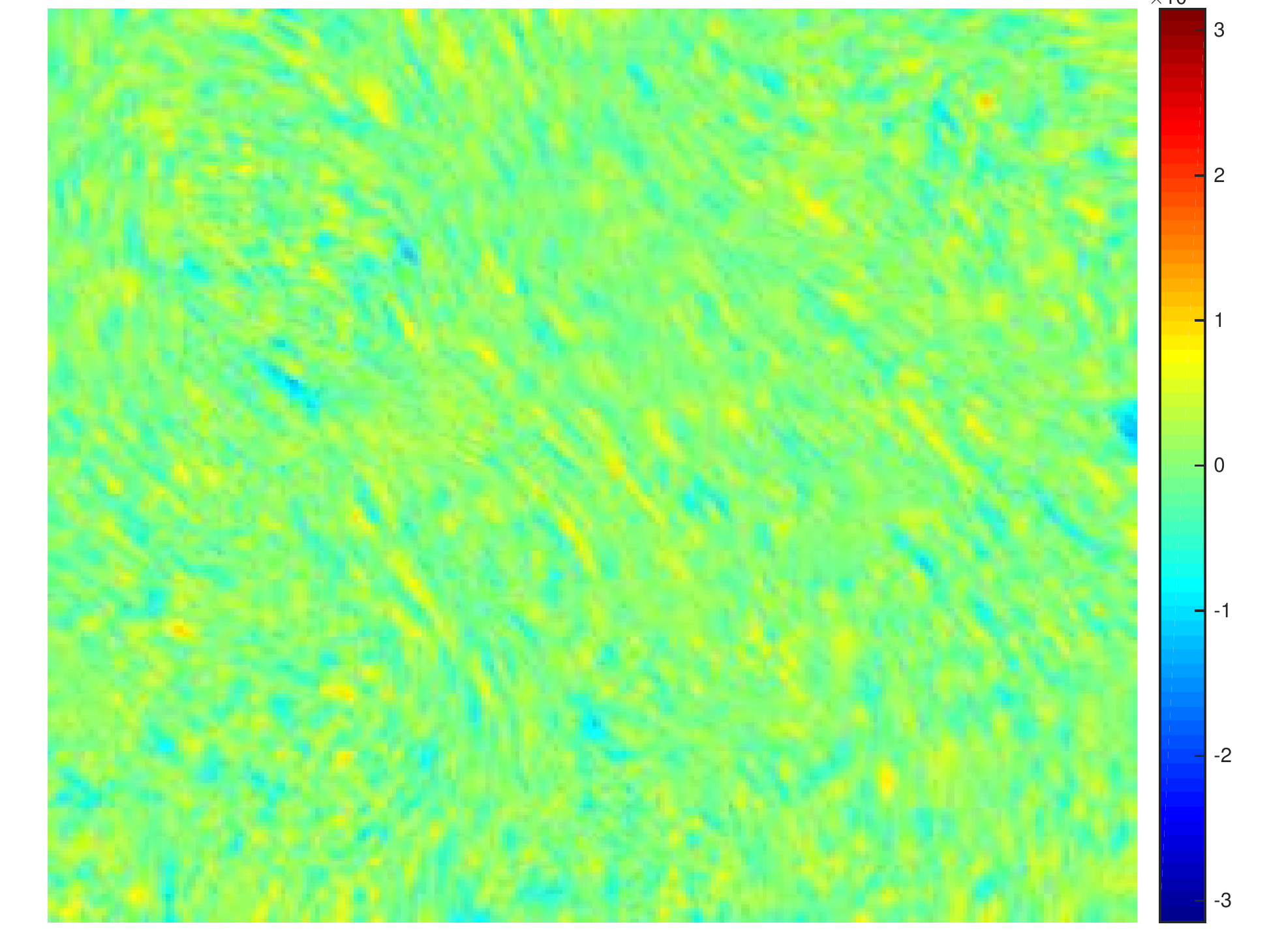}\\
    Original& EBI (38.16dB)& PLE (27.08dB)  & LDMM (\textbf{44.15dB})\\
    \includegraphics[width = 2.5cm]{latticebig_2d_original}&
    \includegraphics[width = 2.5cm]{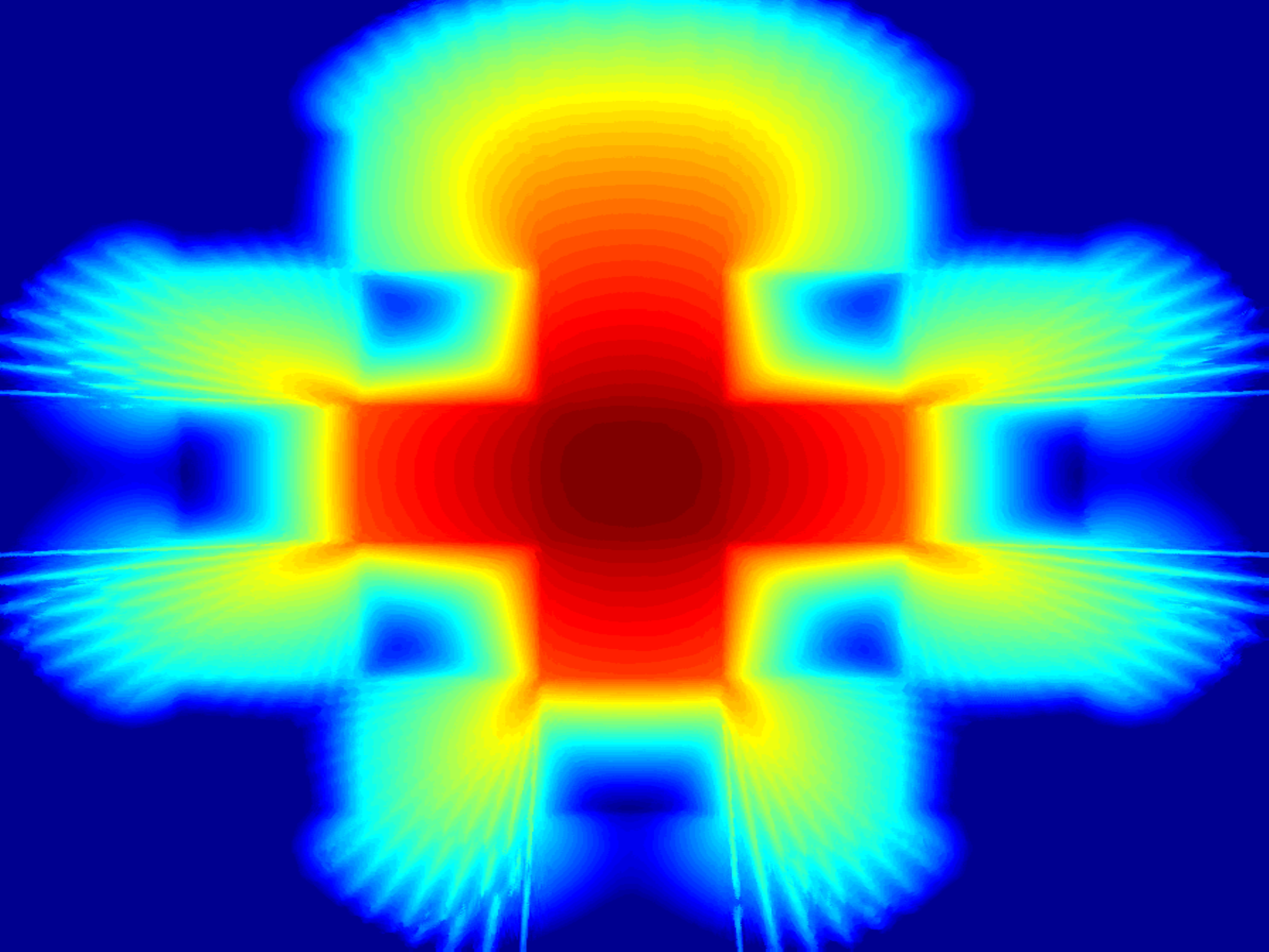}&
    \includegraphics[width = 2.5cm]{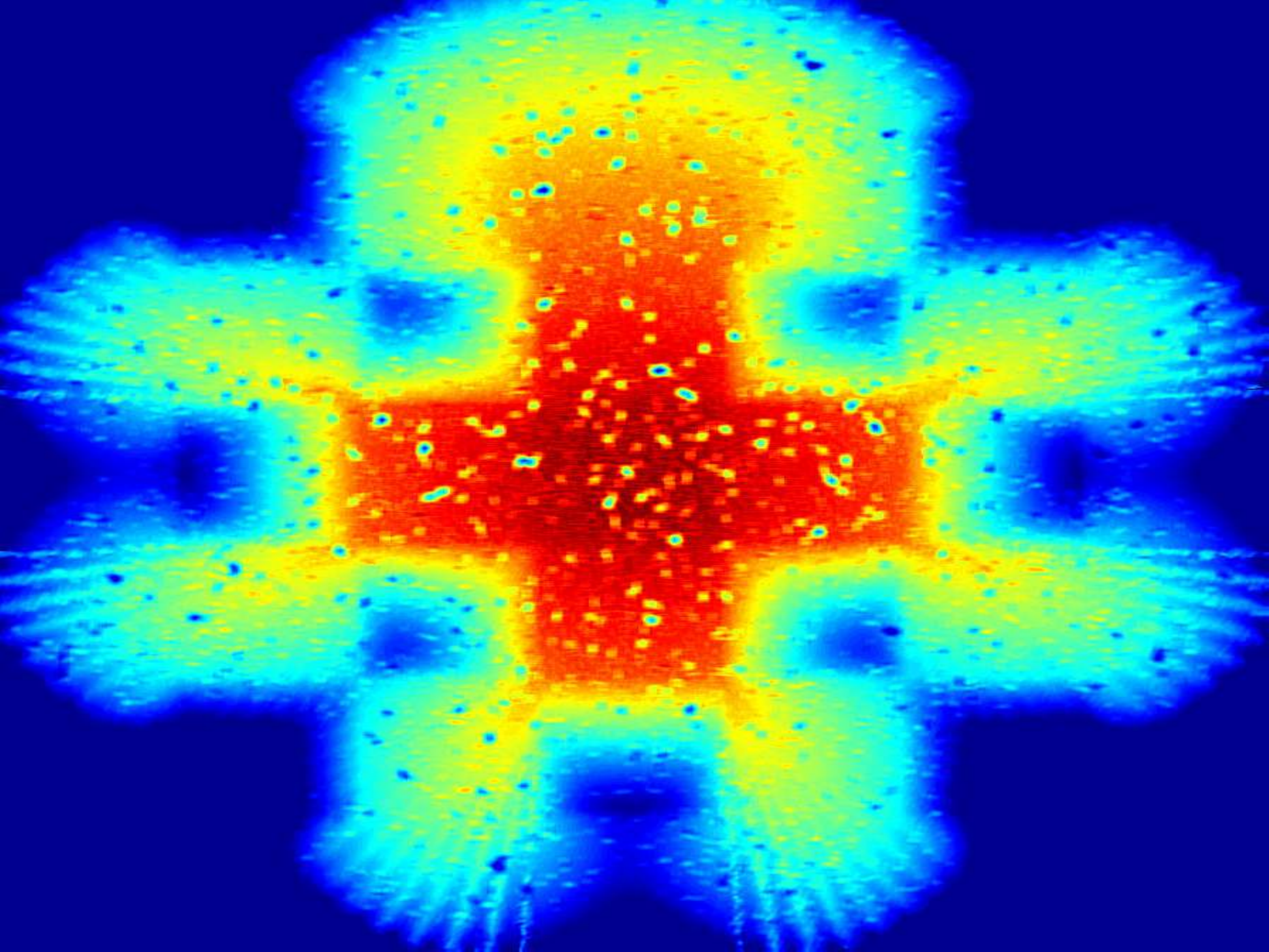}&
    \includegraphics[width = 2.5cm]{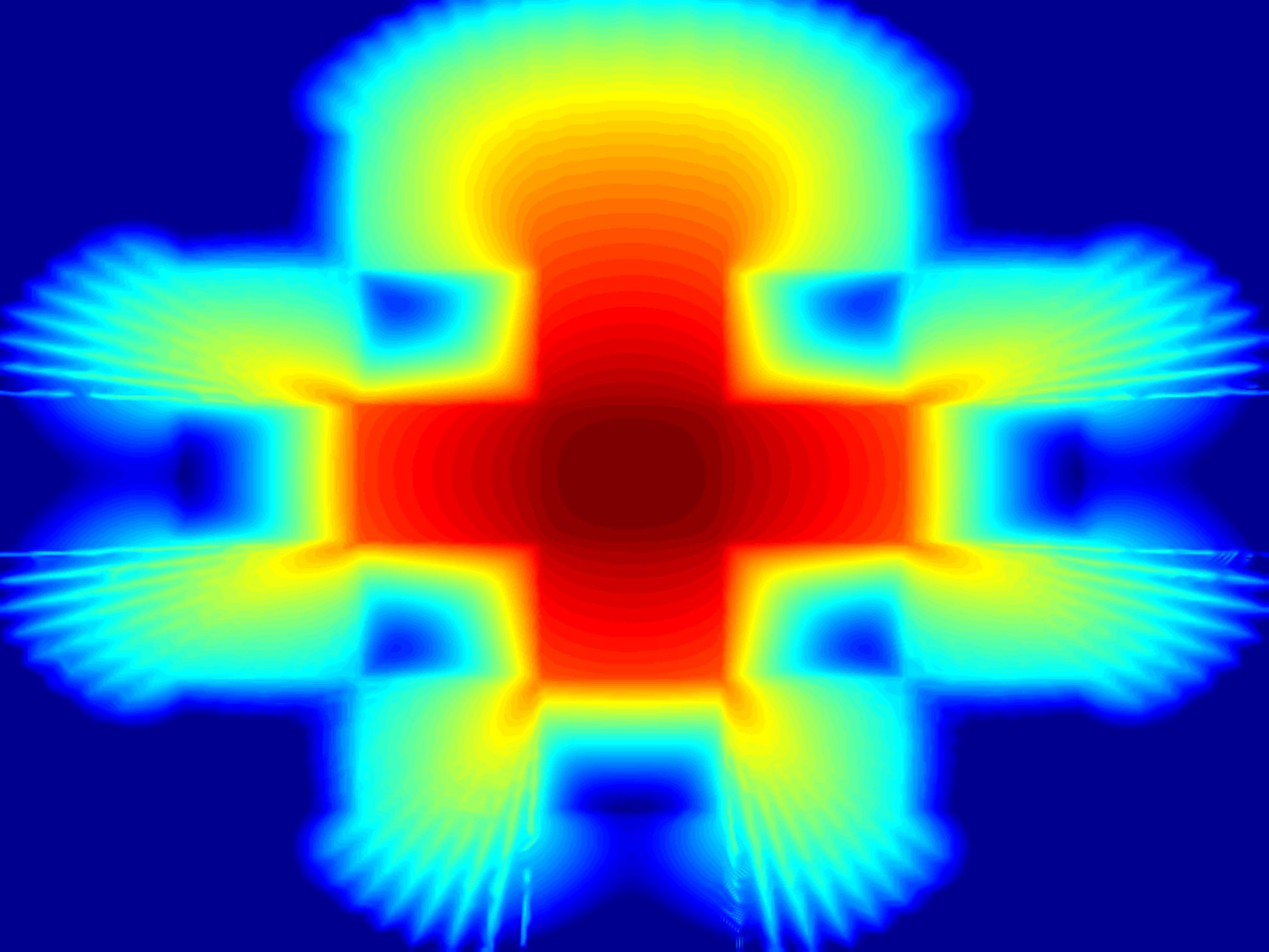}\\
    Subsample & Error & Error & Error\\
    \includegraphics[width = 2.5cm]{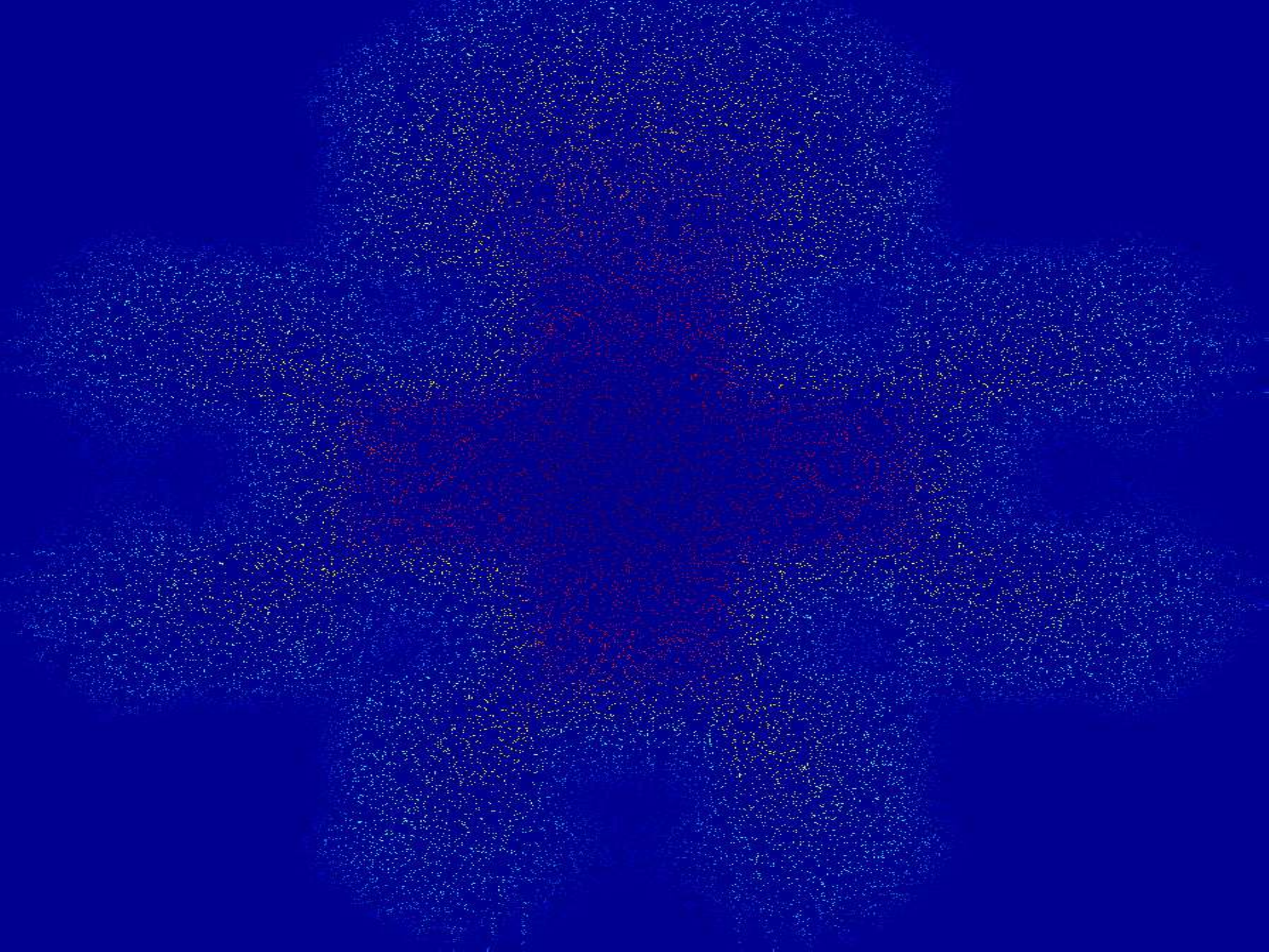}&    
    \includegraphics[width = 2.5cm]{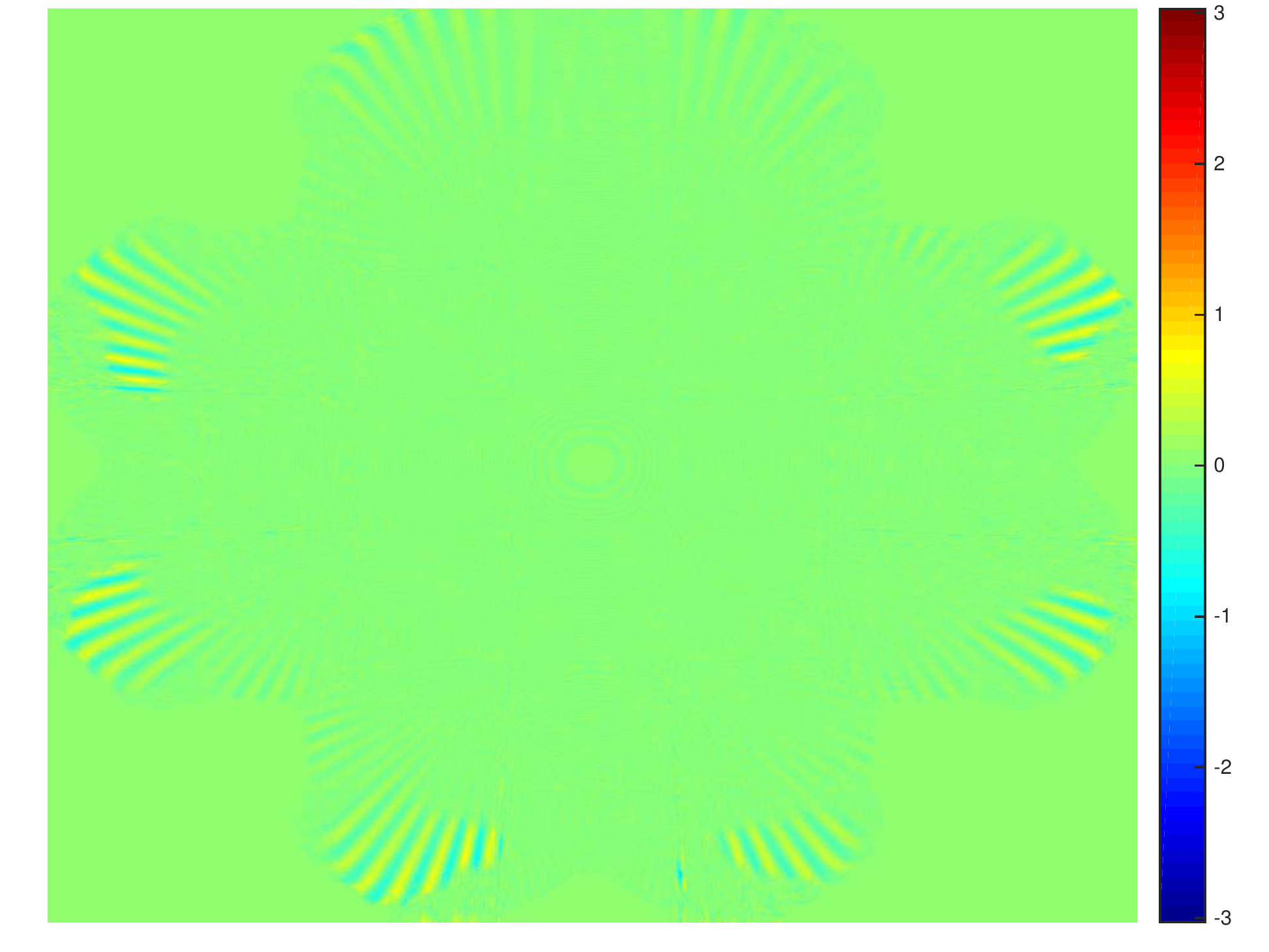}&    
    \includegraphics[width = 2.5cm]{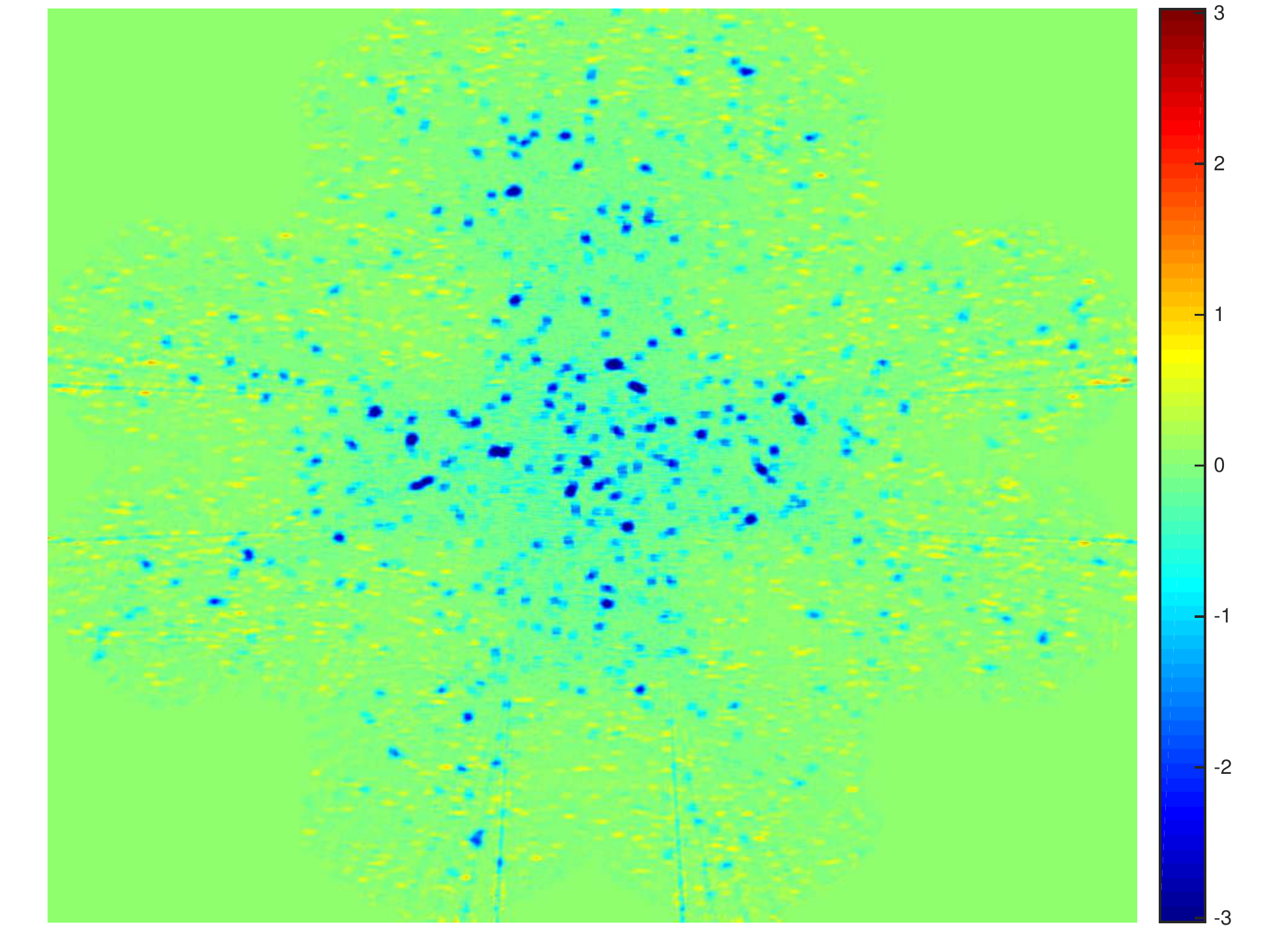}&
    \includegraphics[width = 2.5cm]{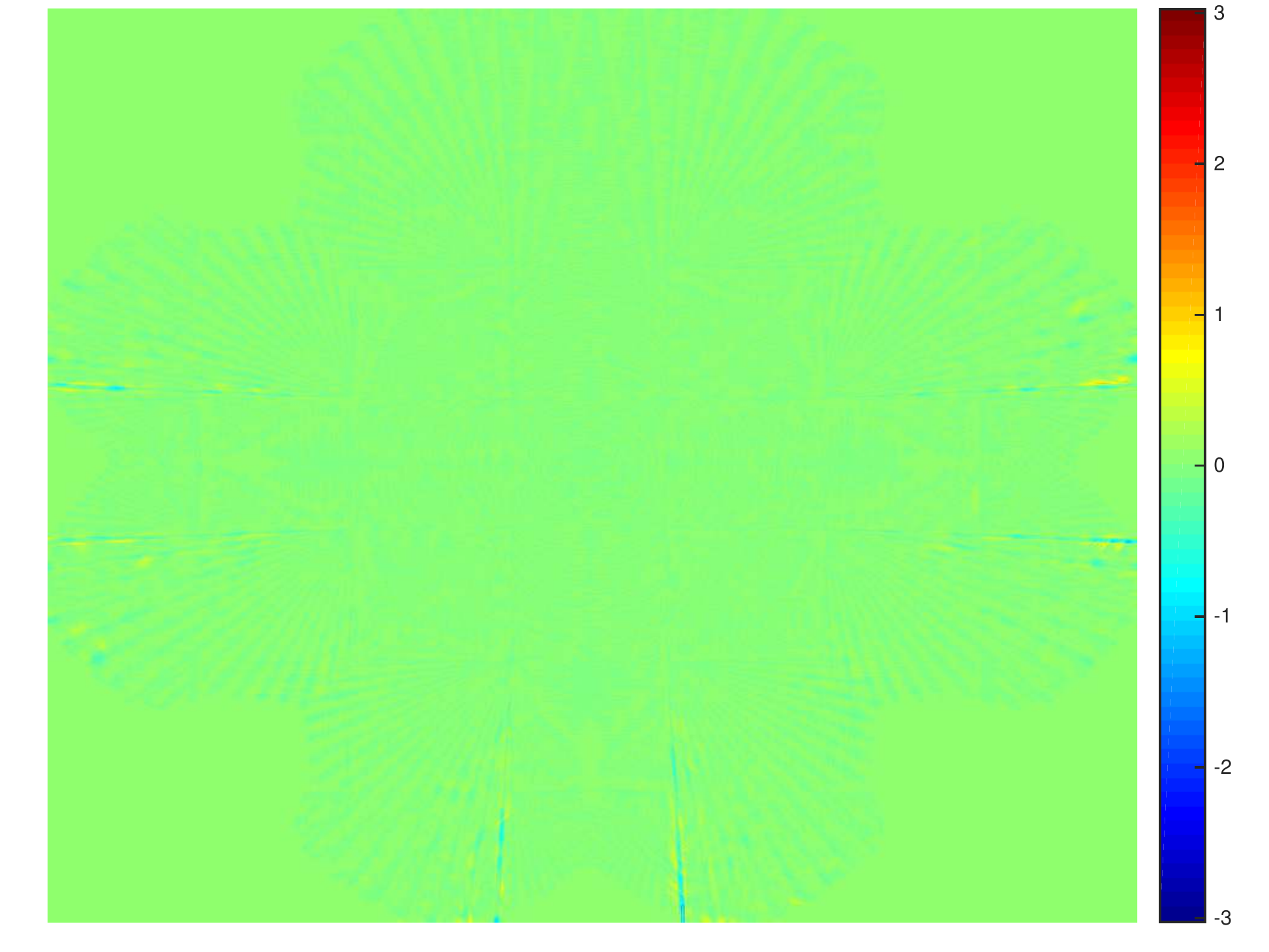}
  \end{tabular}
  \caption{Interpolation of 2D scientific data sets from $5\%$ random sampling. The figures in the first column are the original and subsampled data. The figures in the other three columns are the results and errors of the competing algorithms.}
  \label{fig:result_random_2d_5p}
\end{figure}

\begin{figure}[H]
  \centering
  \begin{tabular}{cccc}
    Original& EBI (37.88dB)& PLE (37.96dB)  & LDMM (\textbf{44.18dB})\\
    \includegraphics[width = 2.5cm]{plasma_3d_original_band_19}&
    \includegraphics[width = 2.5cm]{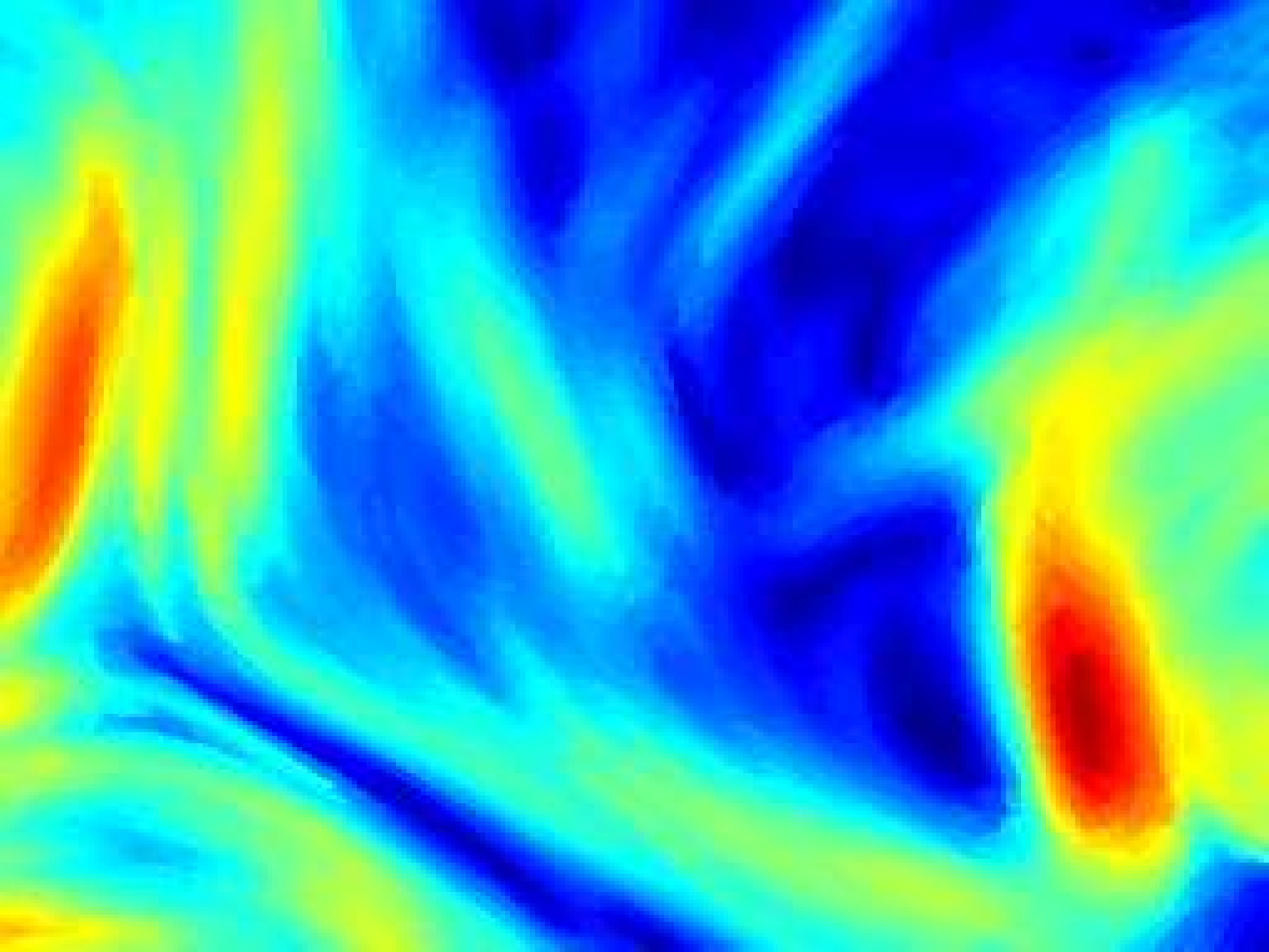}&
    \includegraphics[width = 2.5cm]{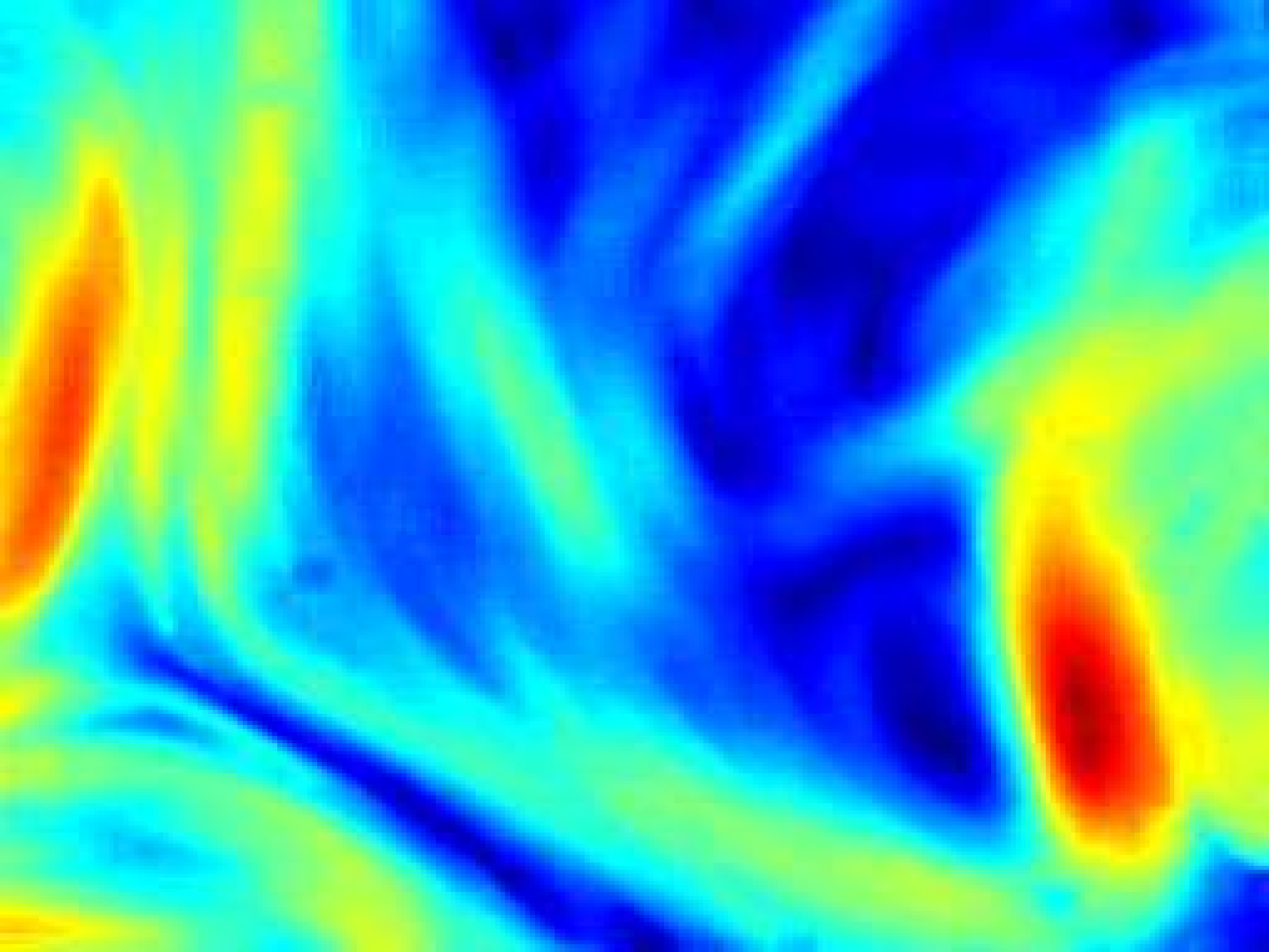}&
    \includegraphics[width = 2.5cm]{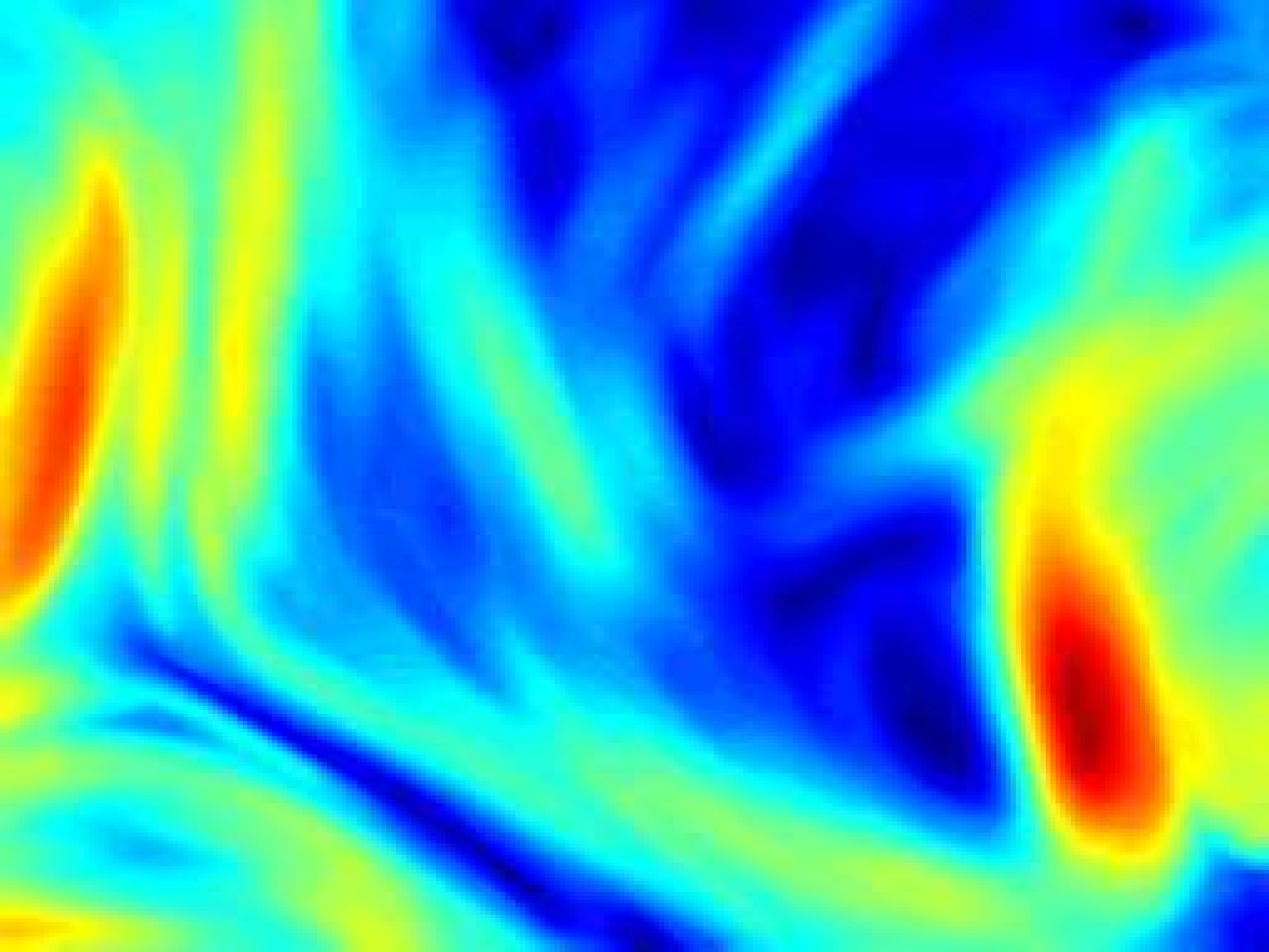}\\
    Subsample & Error & Error & Error\\    
    \includegraphics[width = 2.5cm]{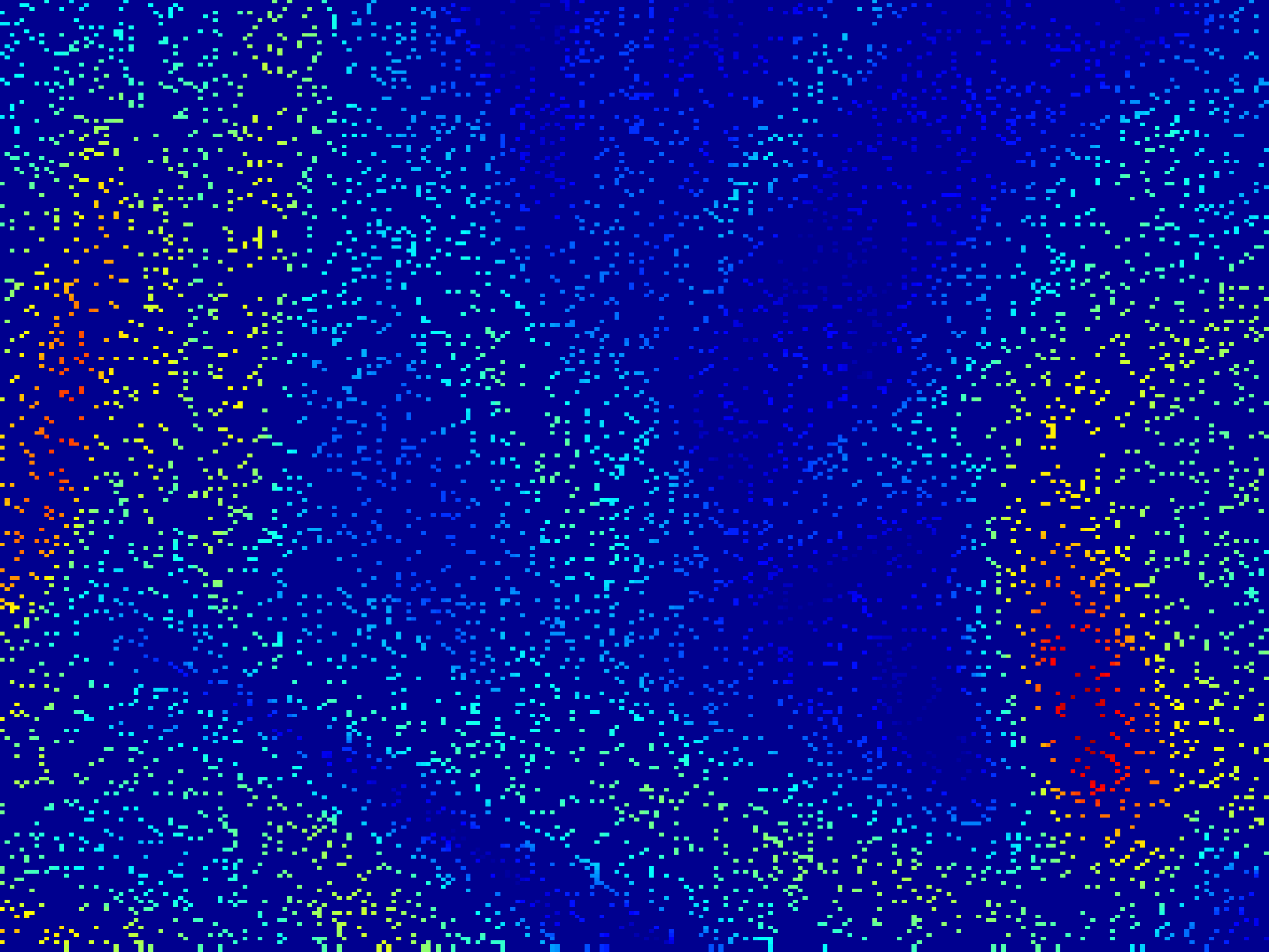}&    
    \includegraphics[width = 2.5cm]{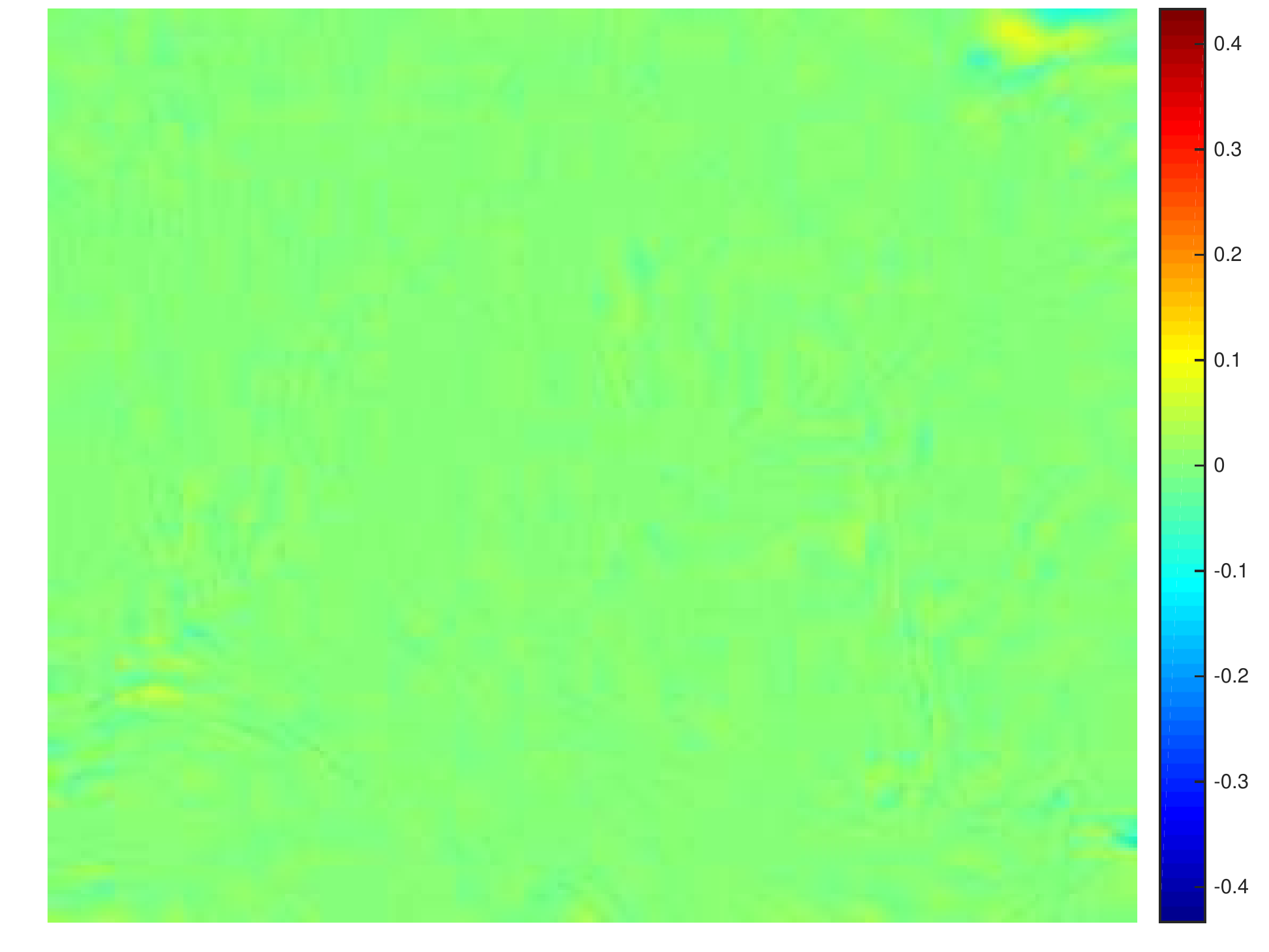}&    
    \includegraphics[width = 2.5cm]{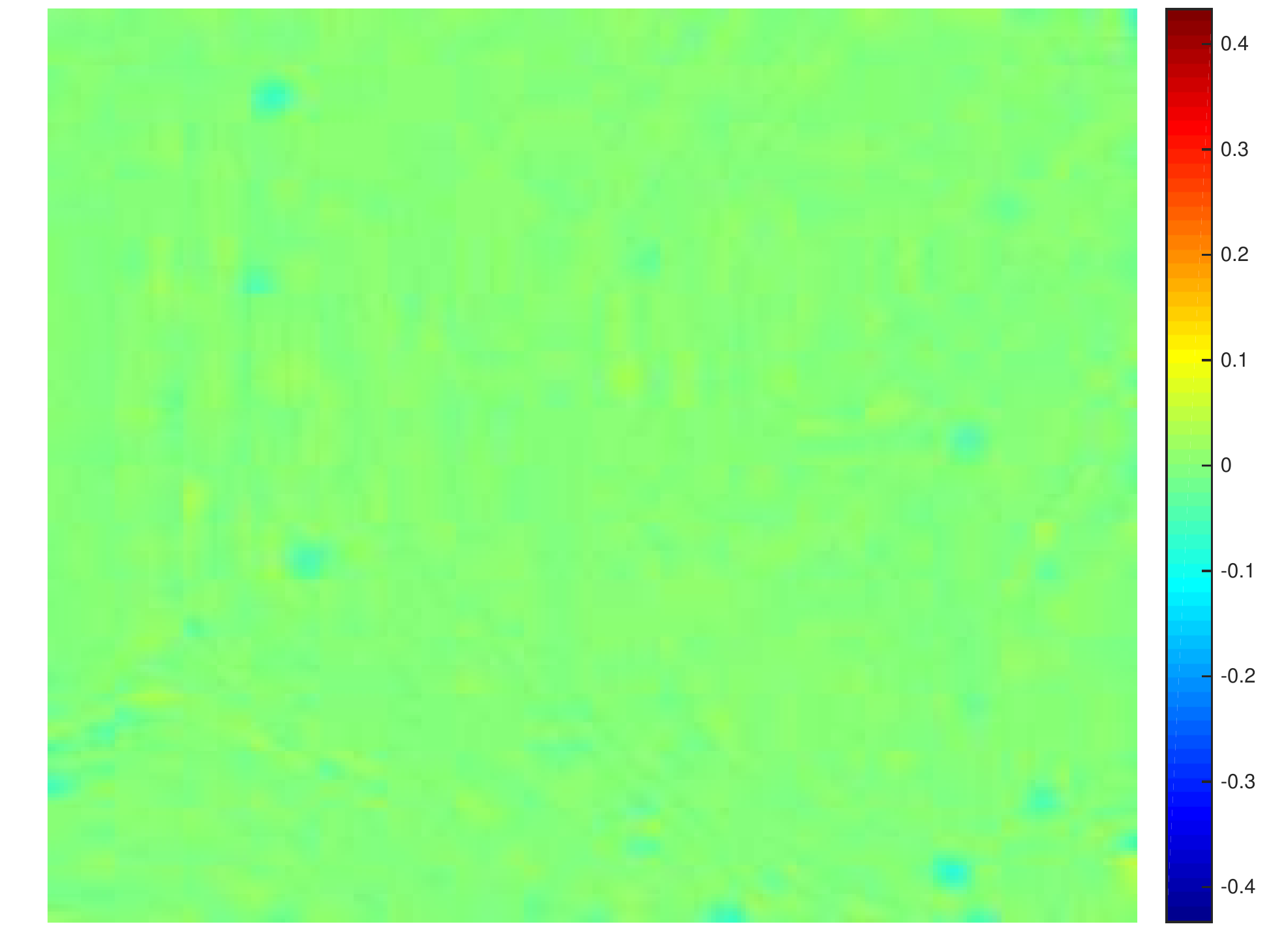}&
    \includegraphics[width = 2.5cm]{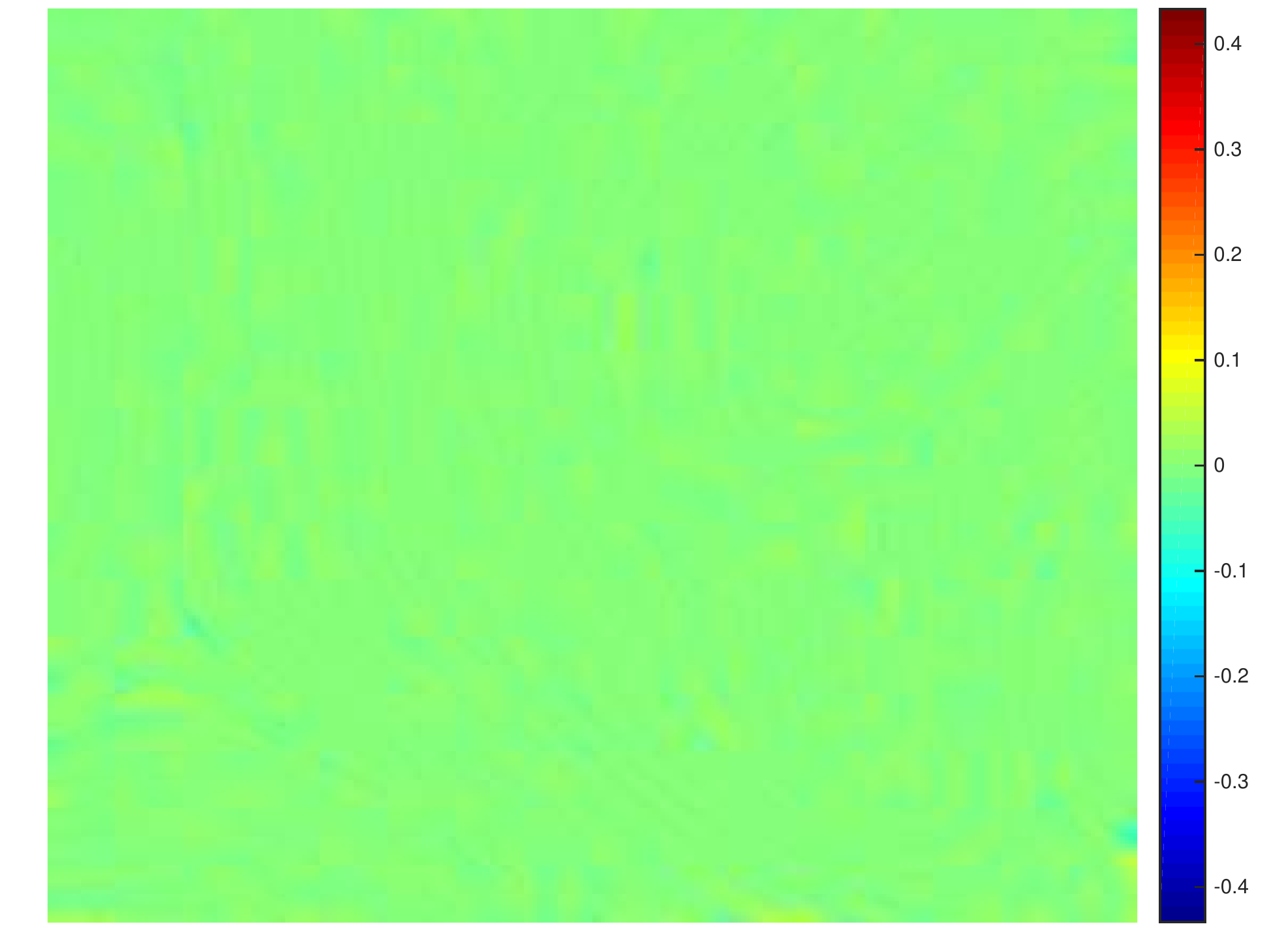}\\
    Original& EBI (37.88dB)& PLE (37.96dB)  & LDMM (\textbf{44.18dB})\\
    \includegraphics[width = 2.5cm]{plasma_3d_original_band_29}&
    \includegraphics[width = 2.5cm]{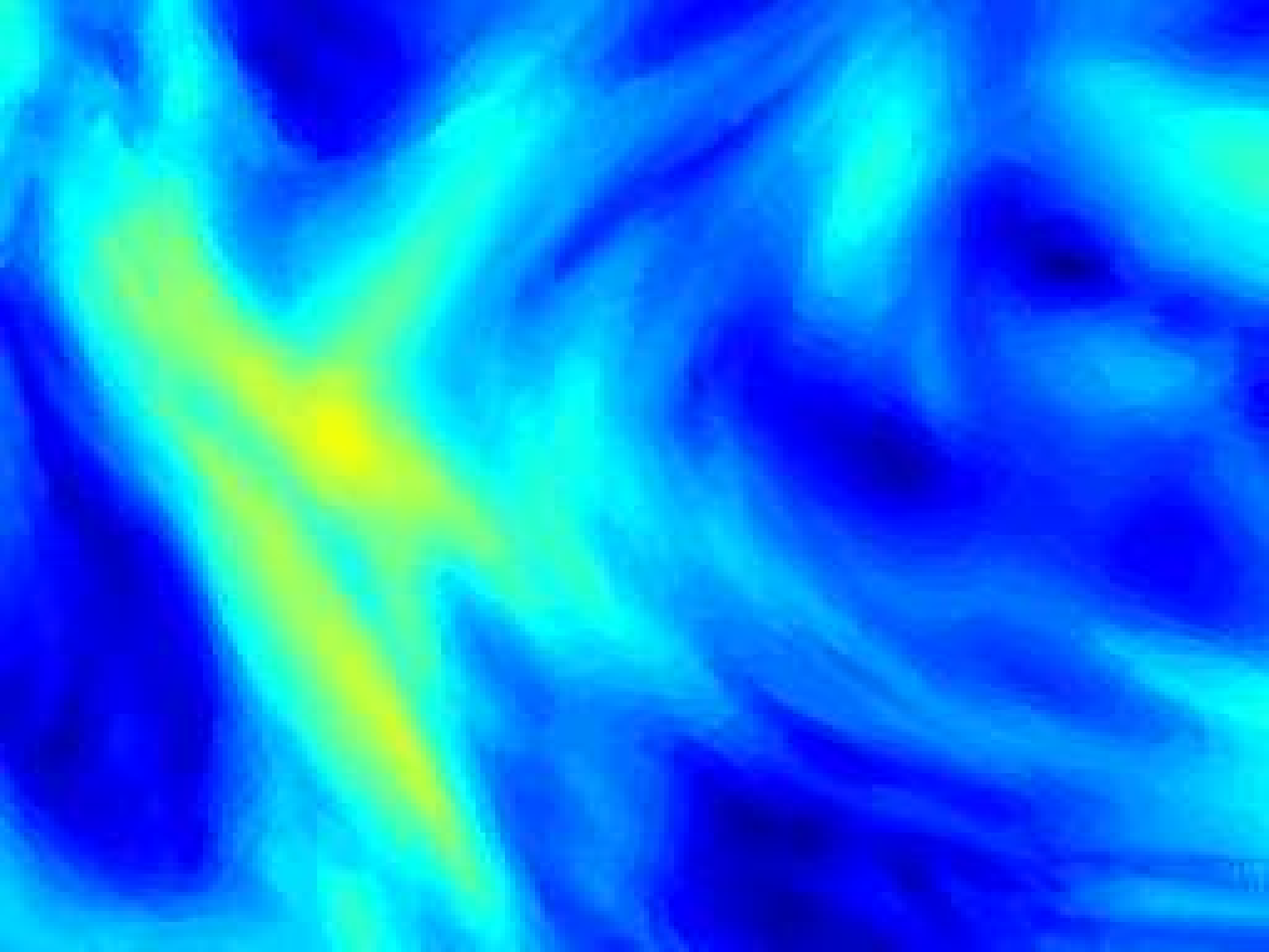}&
    \includegraphics[width = 2.5cm]{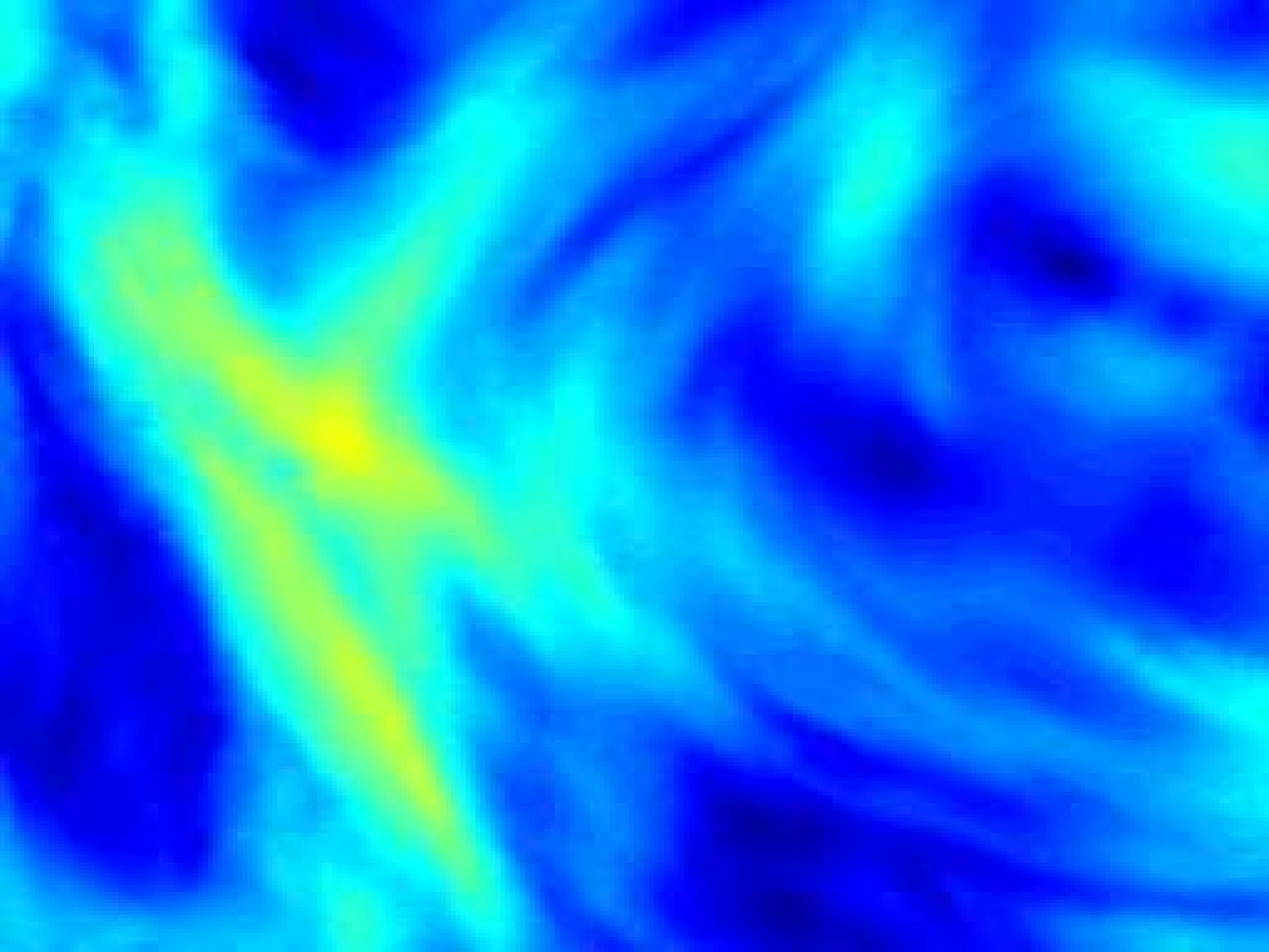}&
    \includegraphics[width = 2.5cm]{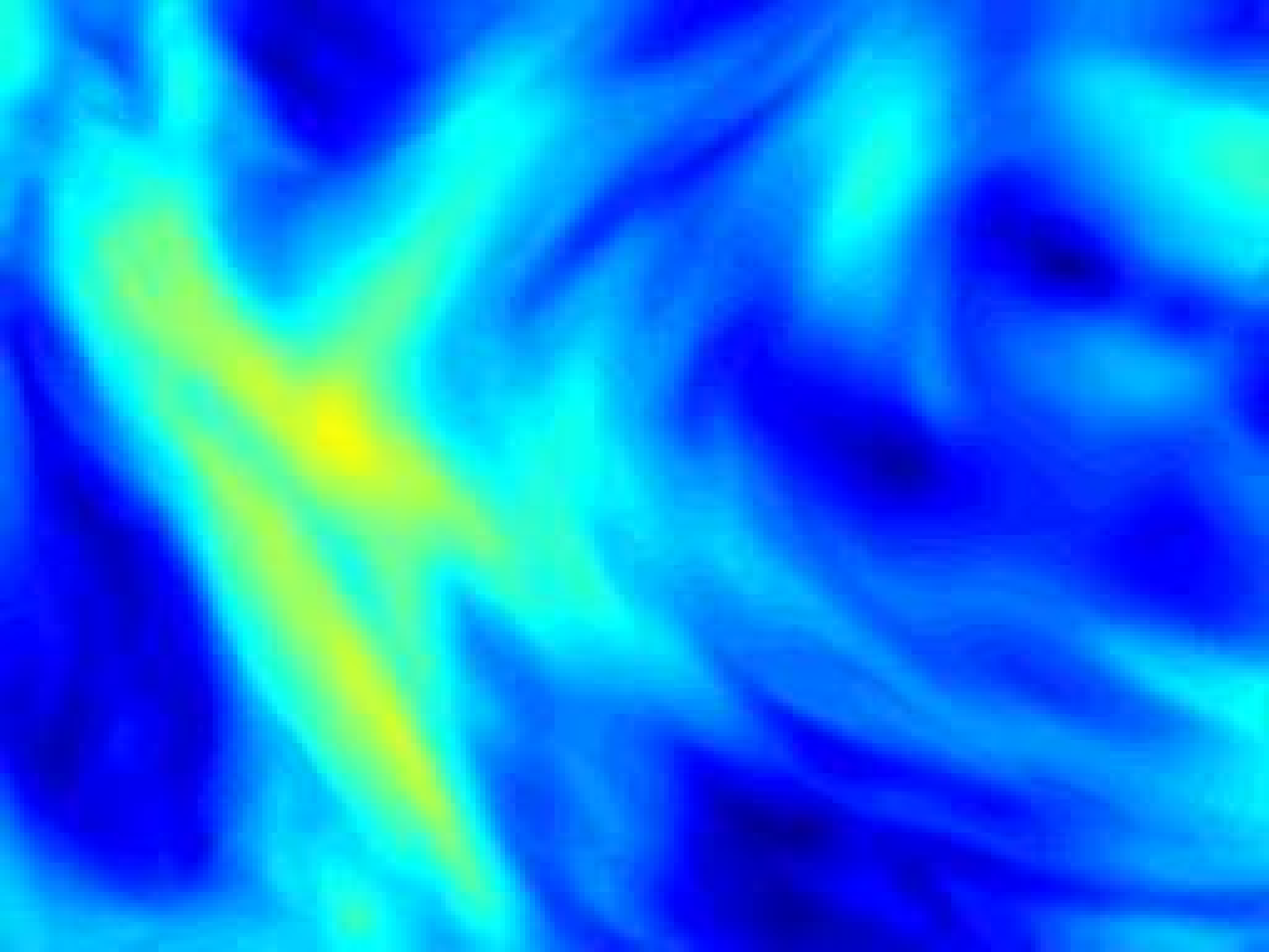}\\
    Subsample & Error & Error & Error\\
    \includegraphics[width = 2.5cm]{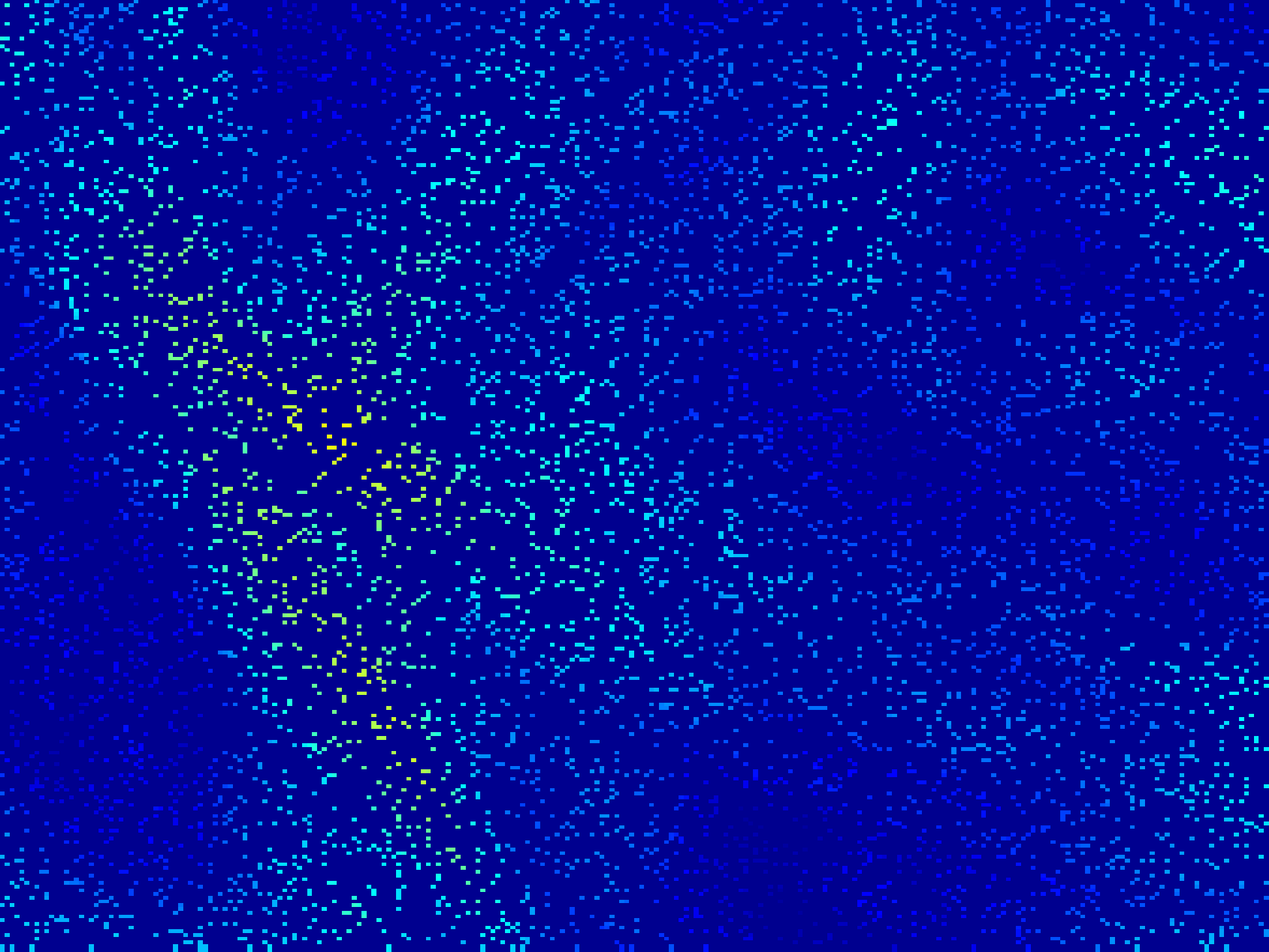}&    
    \includegraphics[width = 2.5cm]{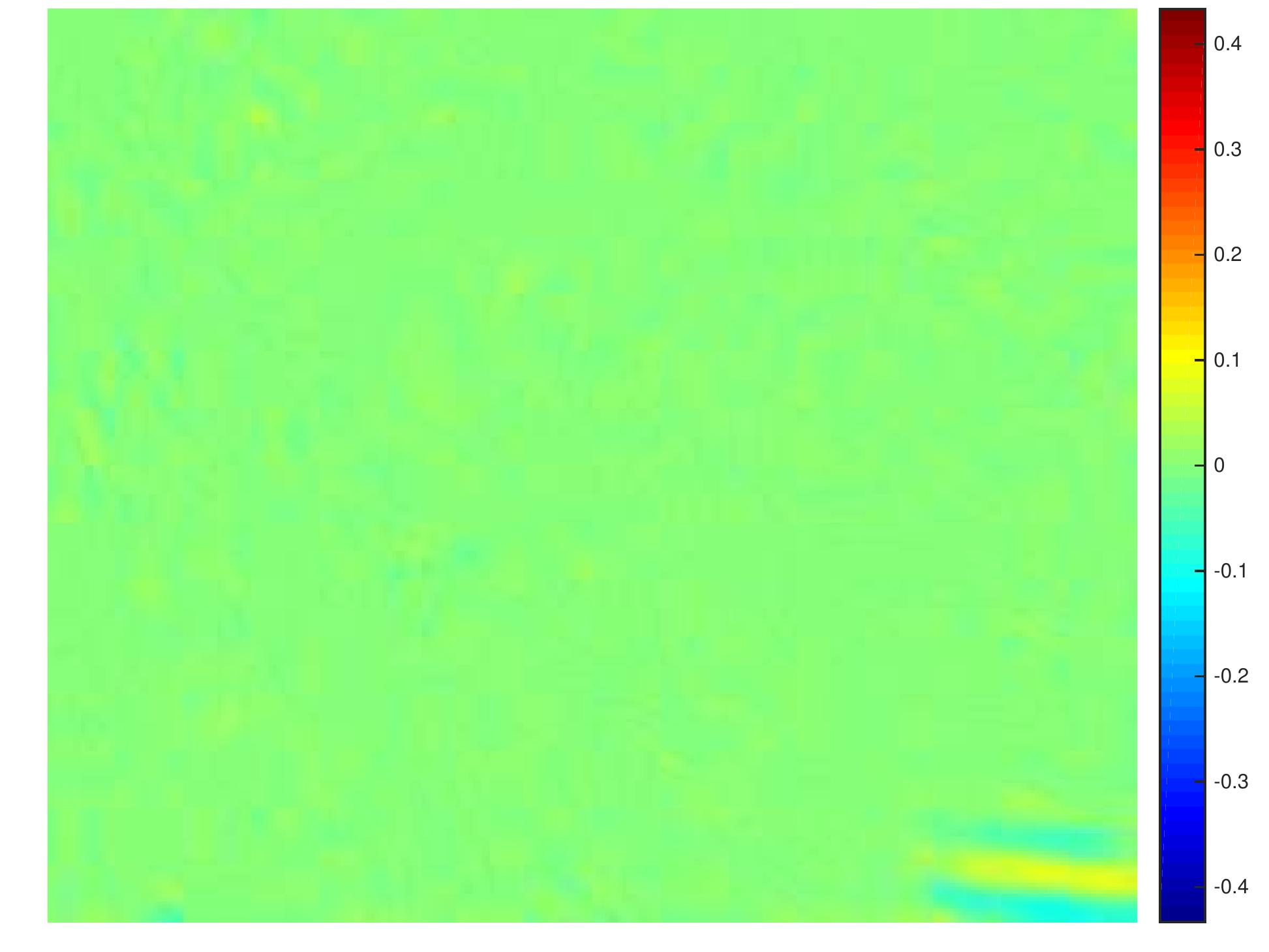}&    
    \includegraphics[width = 2.5cm]{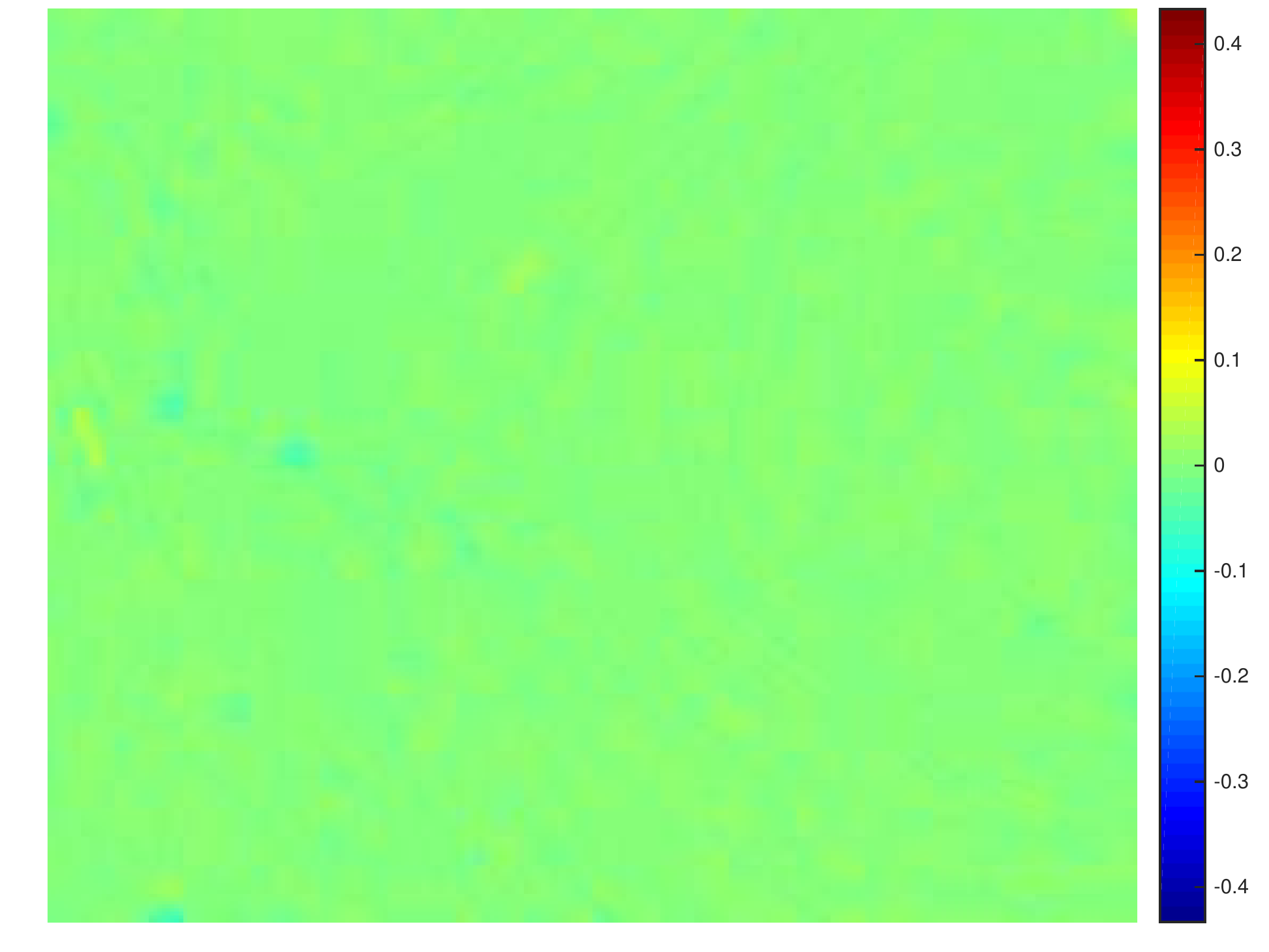}&
    \includegraphics[width = 2.5cm]{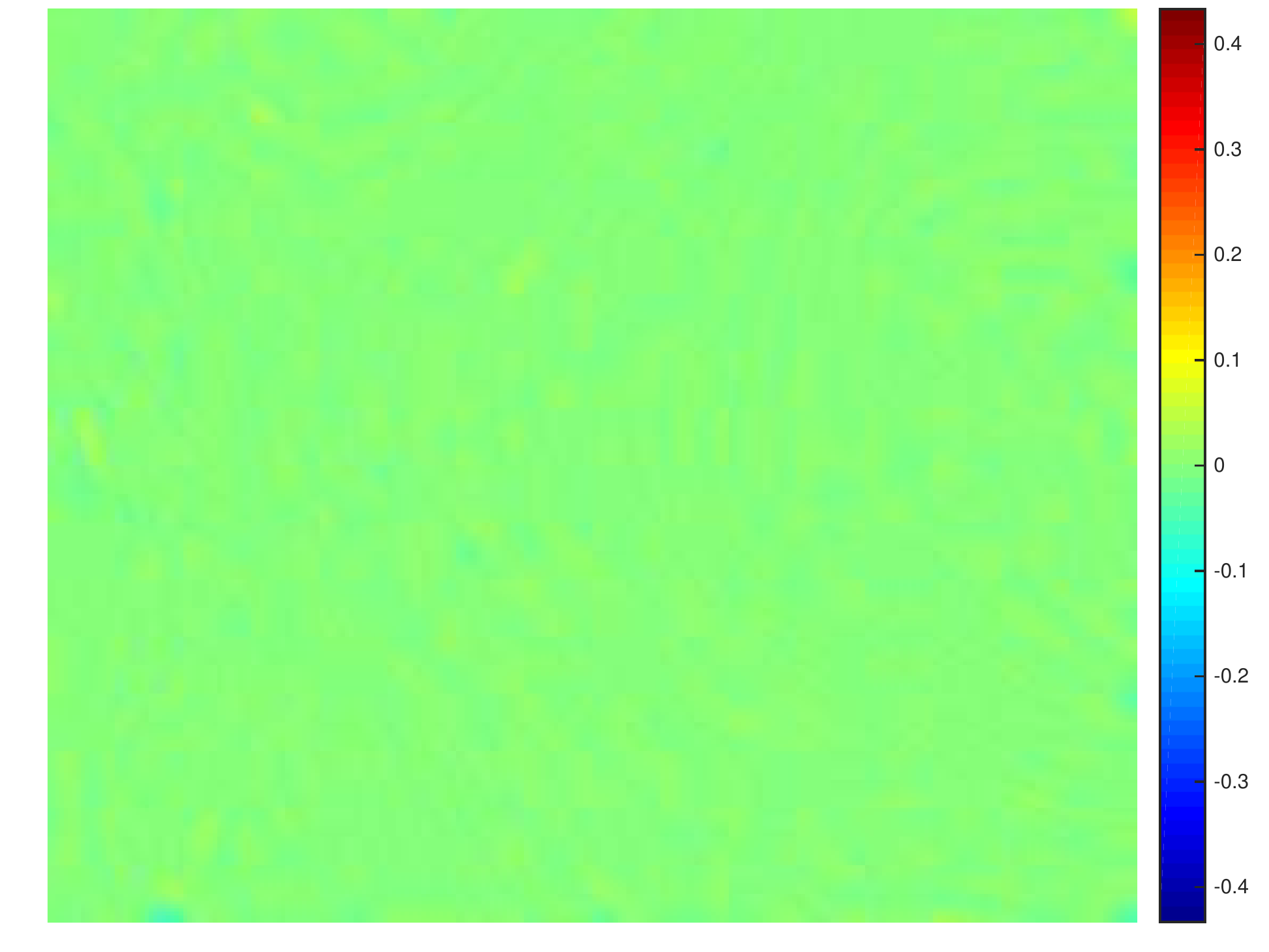}\\
  \end{tabular}
  \caption{Interpolation of the 3D plasma (magnetic field) data set from $10\%$ random sampling. The figures in the first column are two spatial cross sections of the original and subsampled data. The figures in the other three columns are the results and errors of the competing algorithms.}
  \label{fig:result_random_plasma_3d_10p}
\end{figure}

\begin{figure}[H]
  \centering
  \begin{tabular}{cccc}
    Original& EBI (33.93dB)& PLE (25.80dB)  & LDMM (\textbf{40.07dB})\\
    \includegraphics[width = 2.5cm]{plasma_3d_original_band_19}&
    \includegraphics[width = 2.5cm]{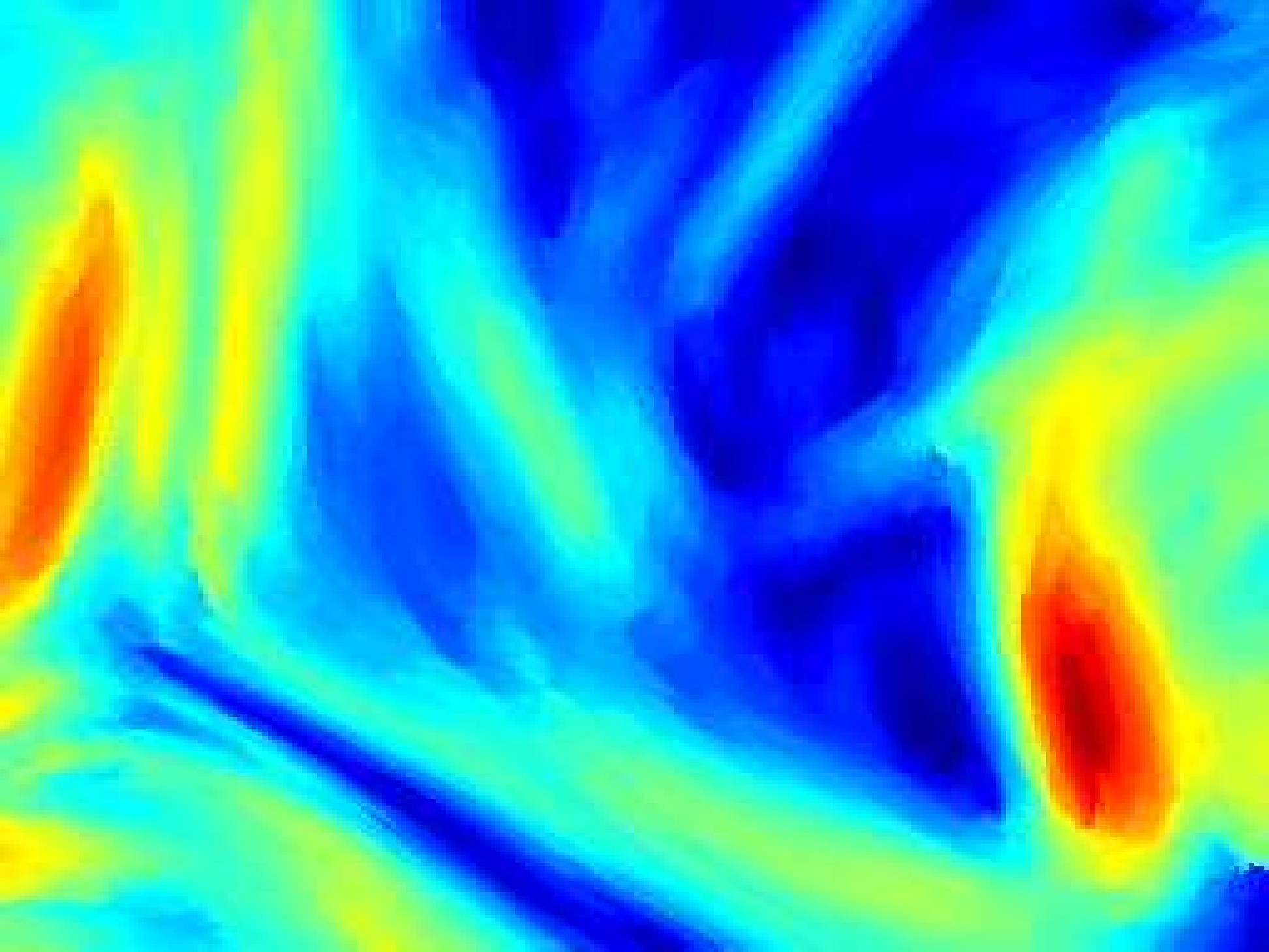}&
    \includegraphics[width = 2.5cm]{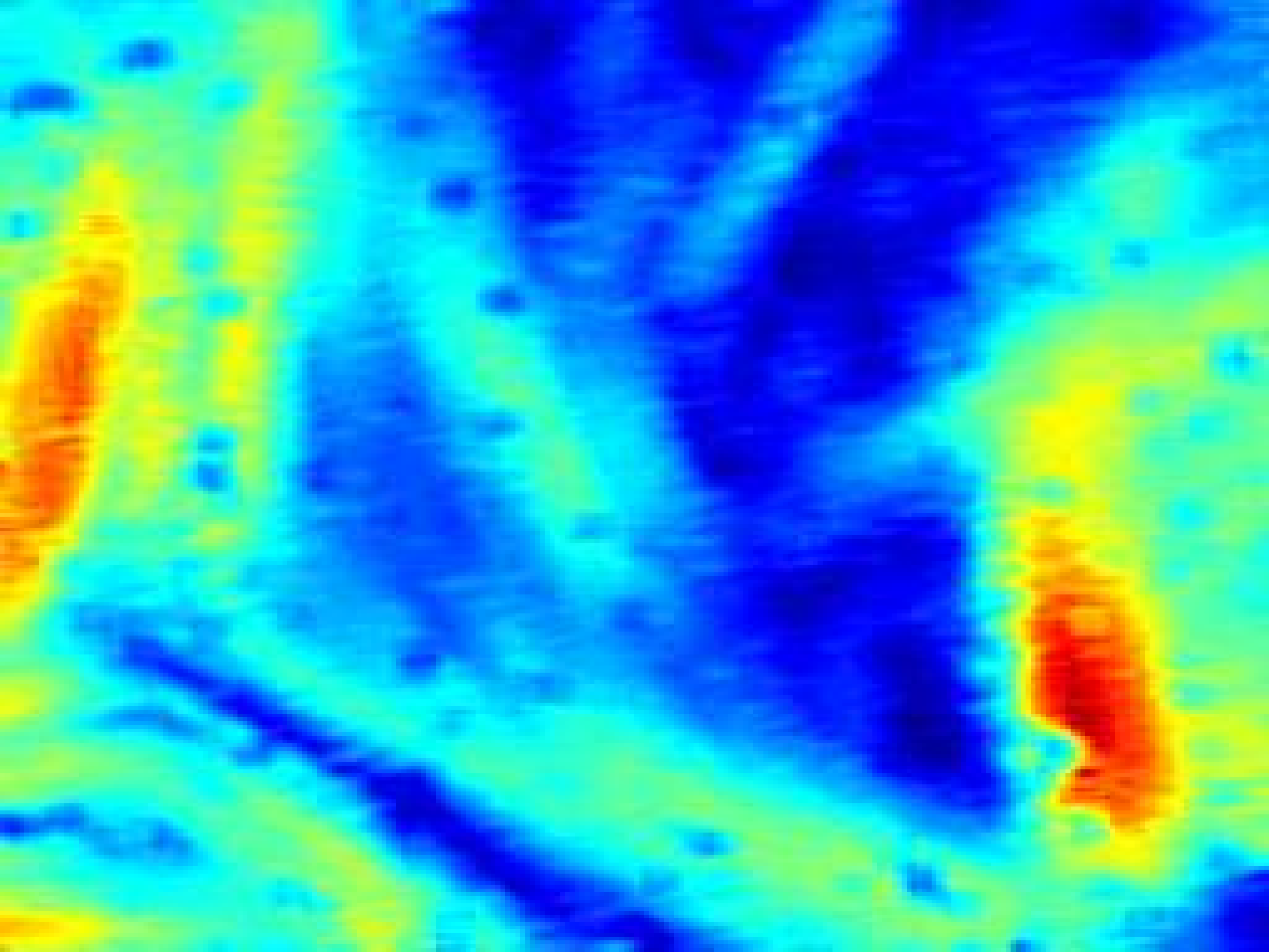}&
    \includegraphics[width = 2.5cm]{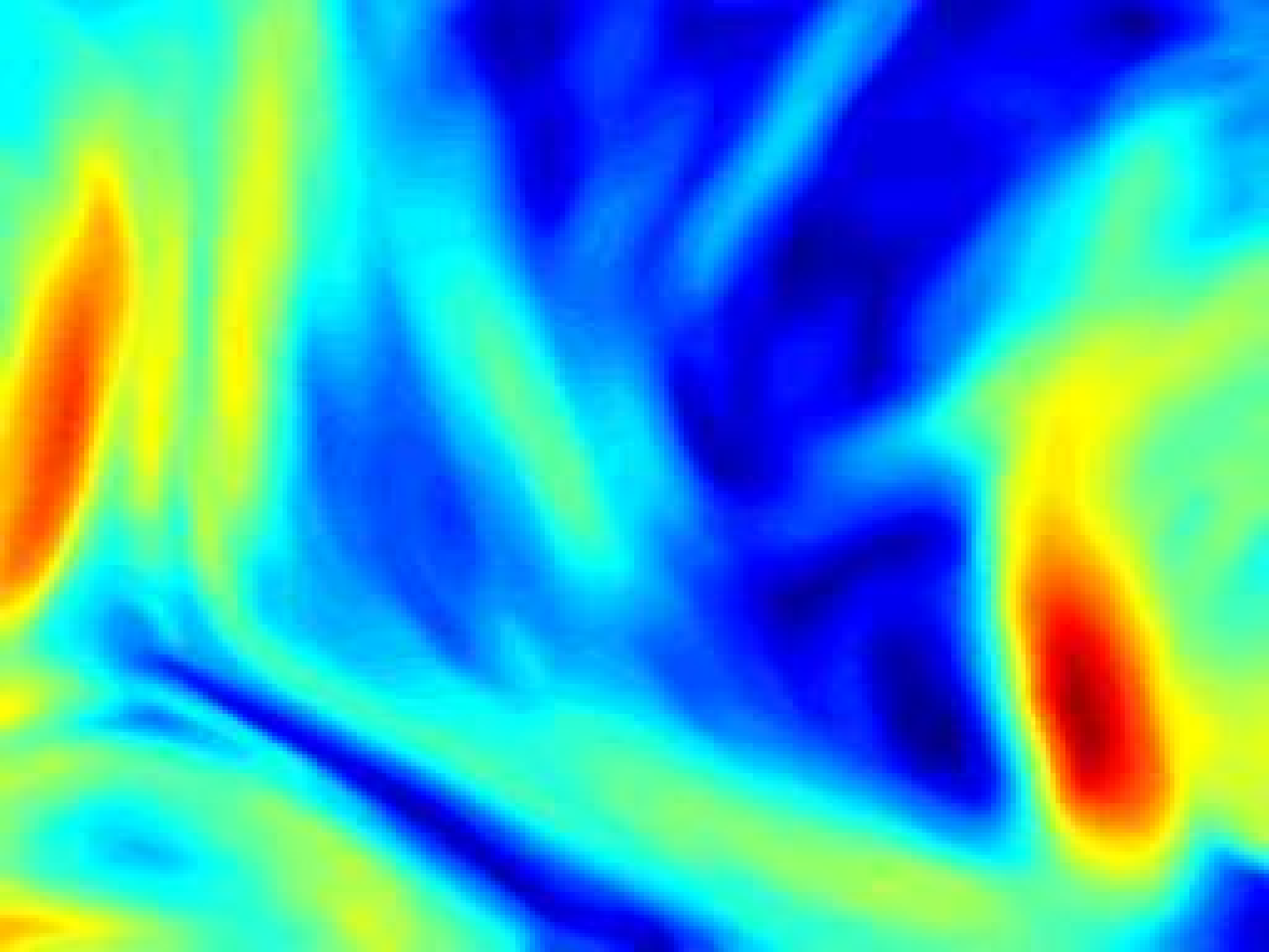}\\
    Subsample & Error & Error & Error\\
    \includegraphics[width = 2.5cm]{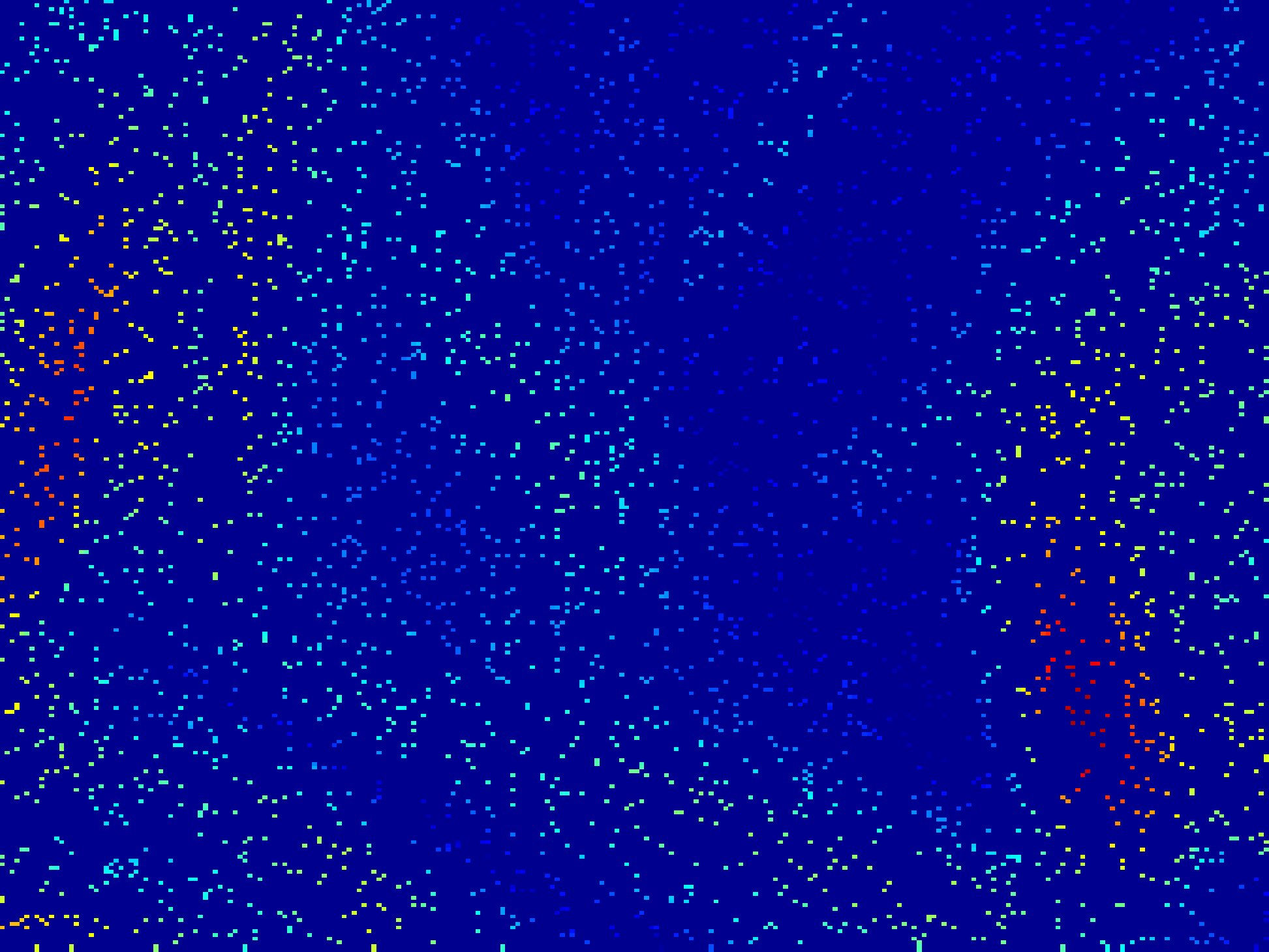}&    
    \includegraphics[width = 2.5cm]{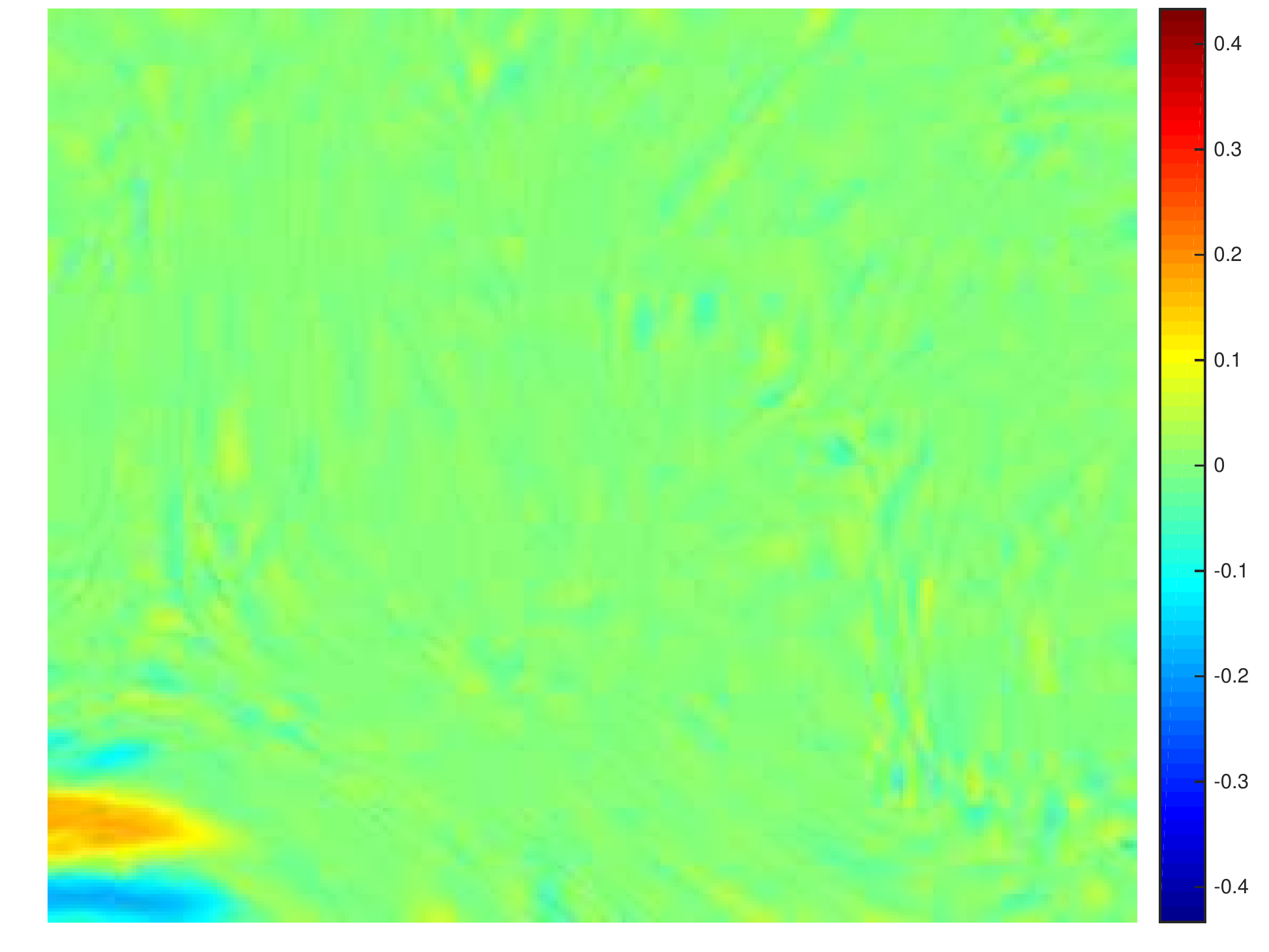}&    
    \includegraphics[width = 2.5cm]{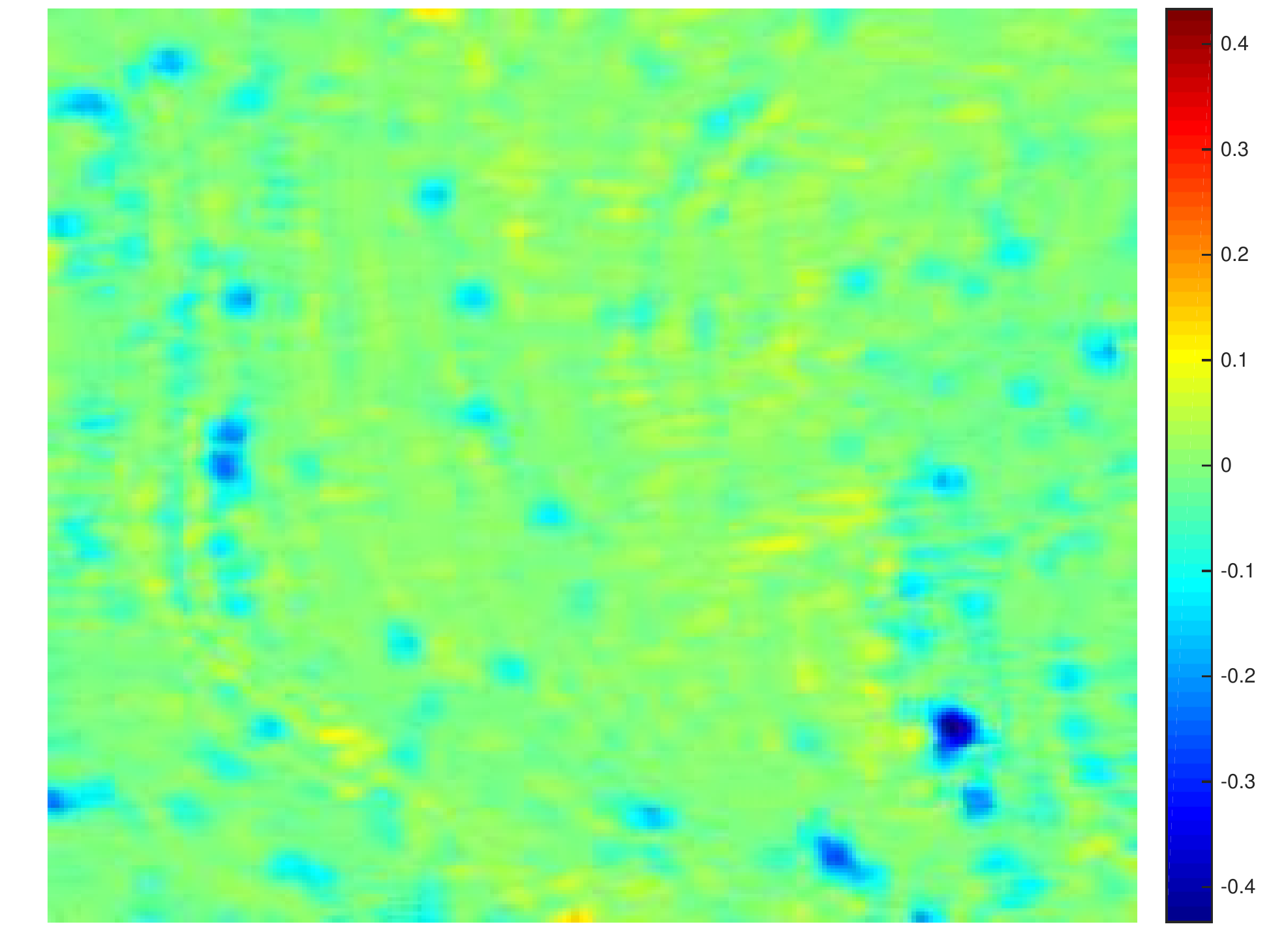}&
    \includegraphics[width = 2.5cm]{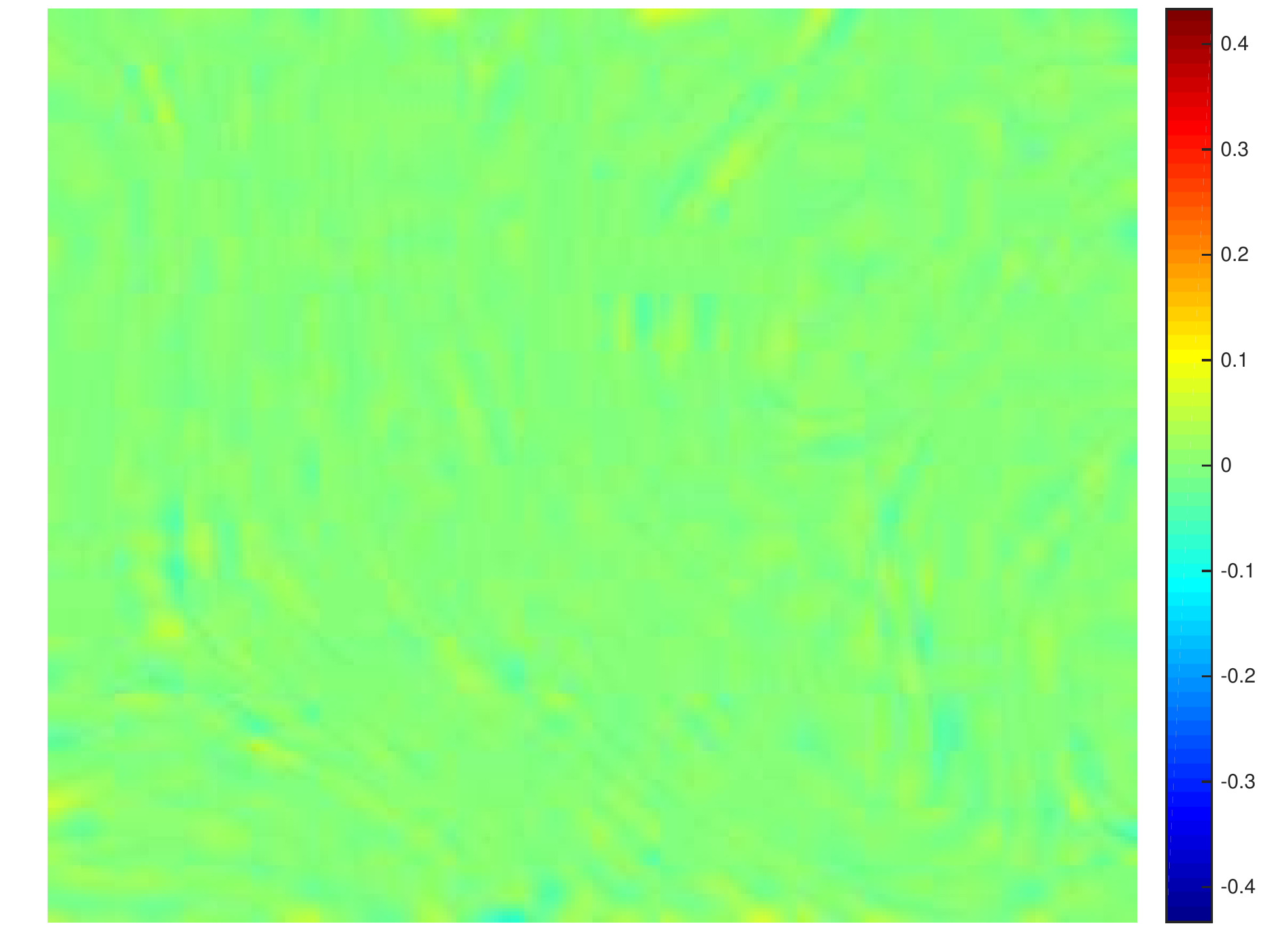}\\
    Original& EBI (33.93dB)& PLE (25.80dB)  & LDMM (\textbf{40.07dB})\\
    \includegraphics[width = 2.5cm]{plasma_3d_original_band_29}&
    \includegraphics[width = 2.5cm]{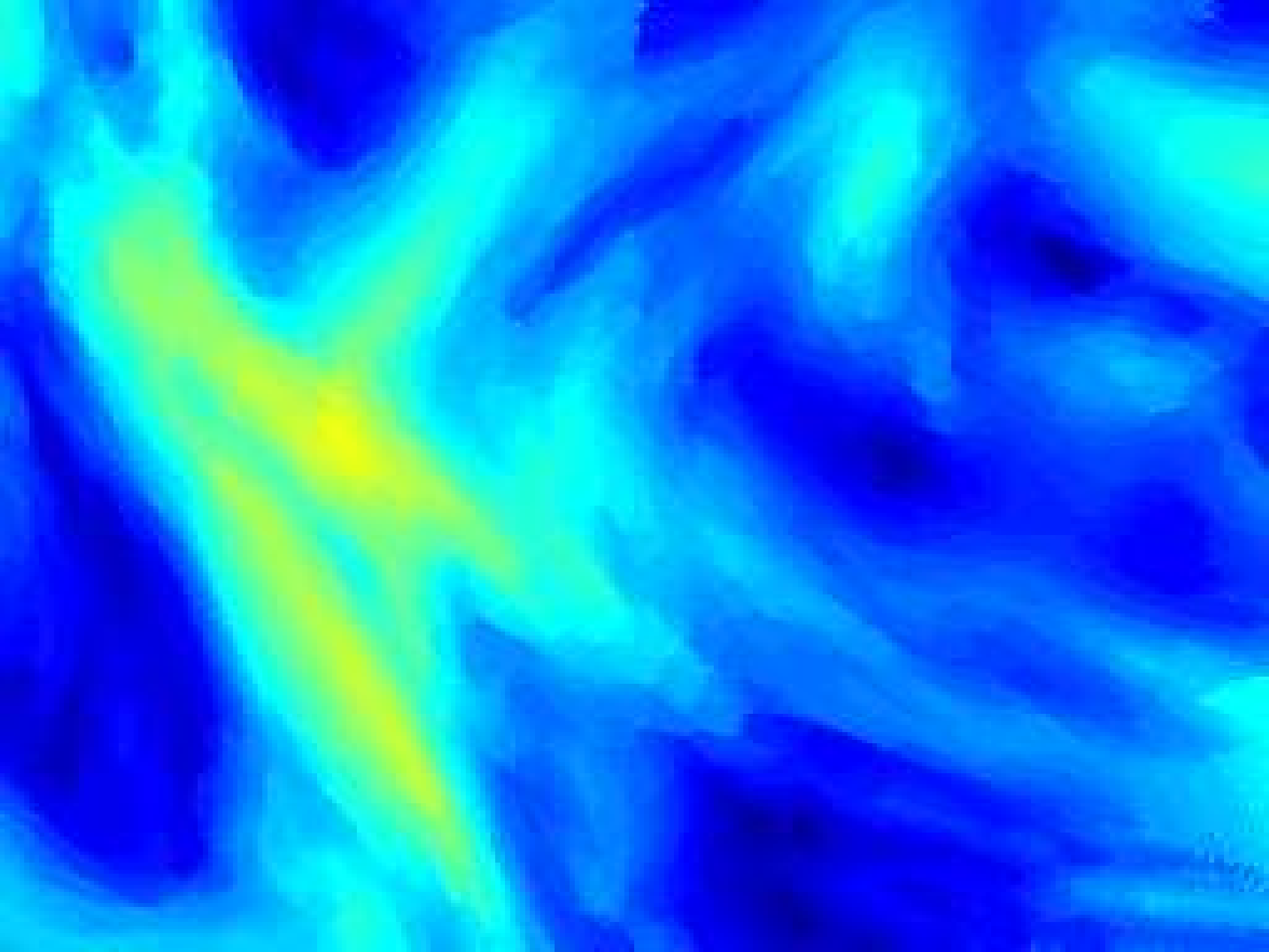}&
    \includegraphics[width = 2.5cm]{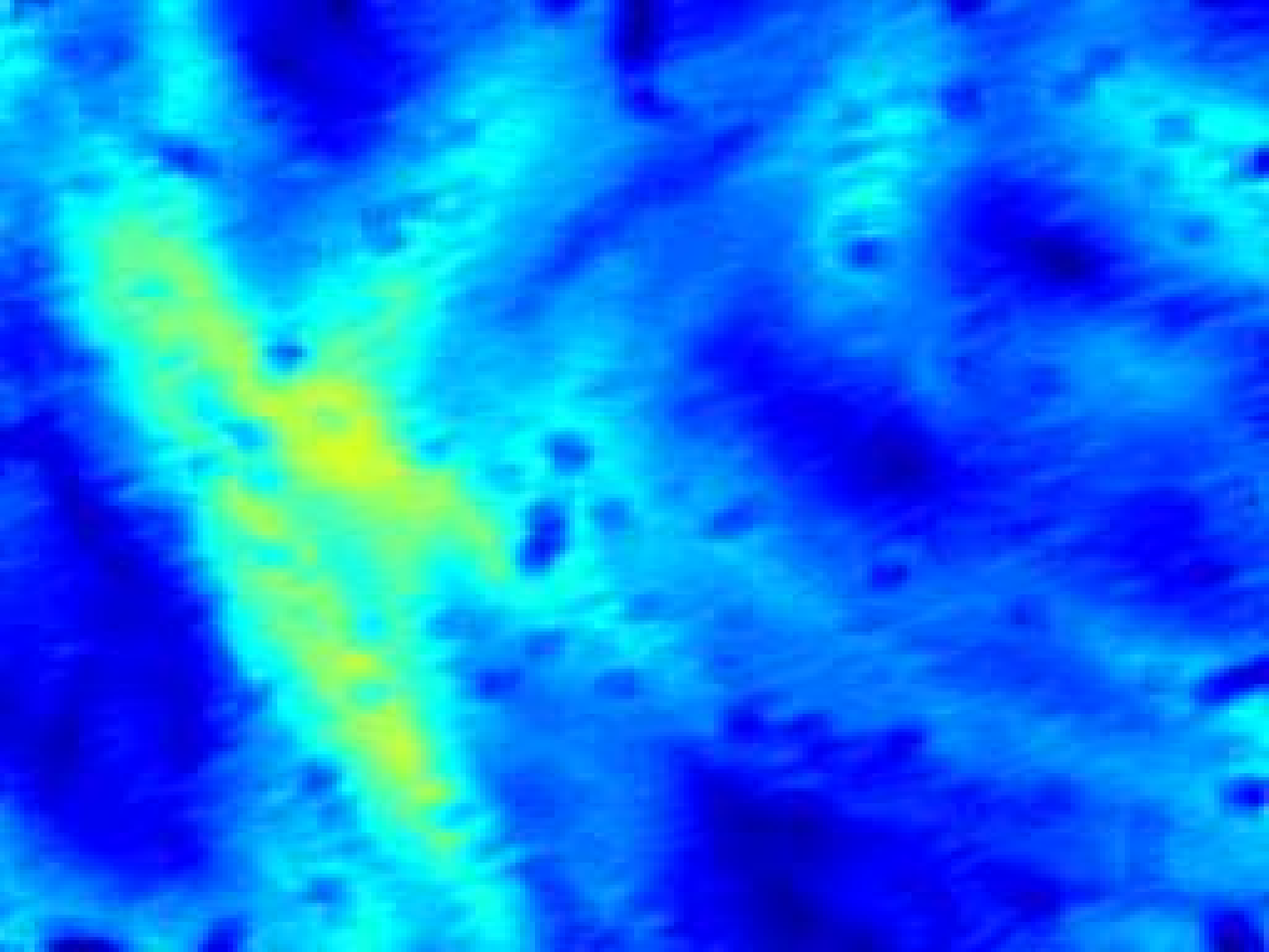}&
    \includegraphics[width = 2.5cm]{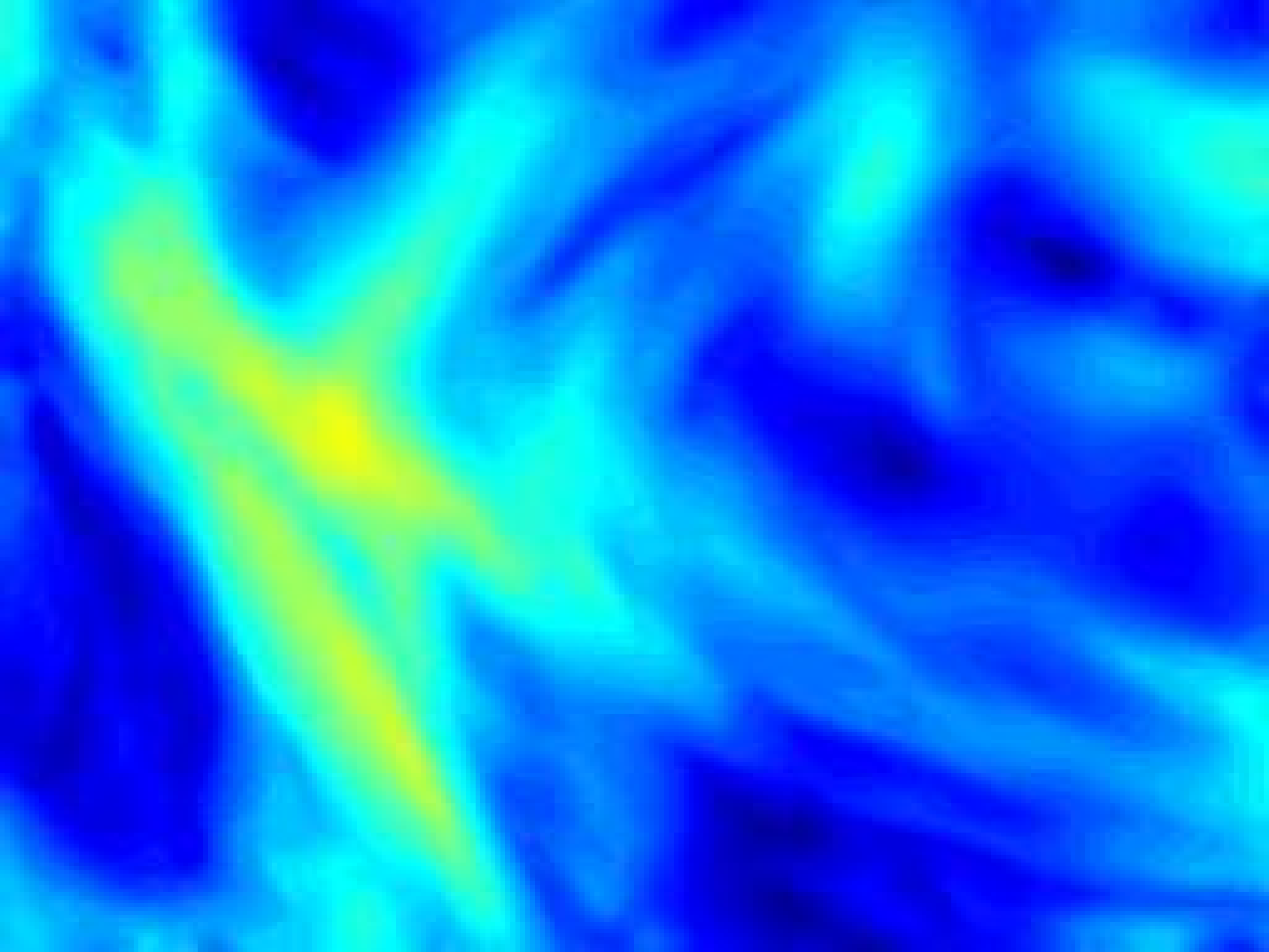}\\
    Subsample & Error & Error & Error\\
    \includegraphics[width = 2.5cm]{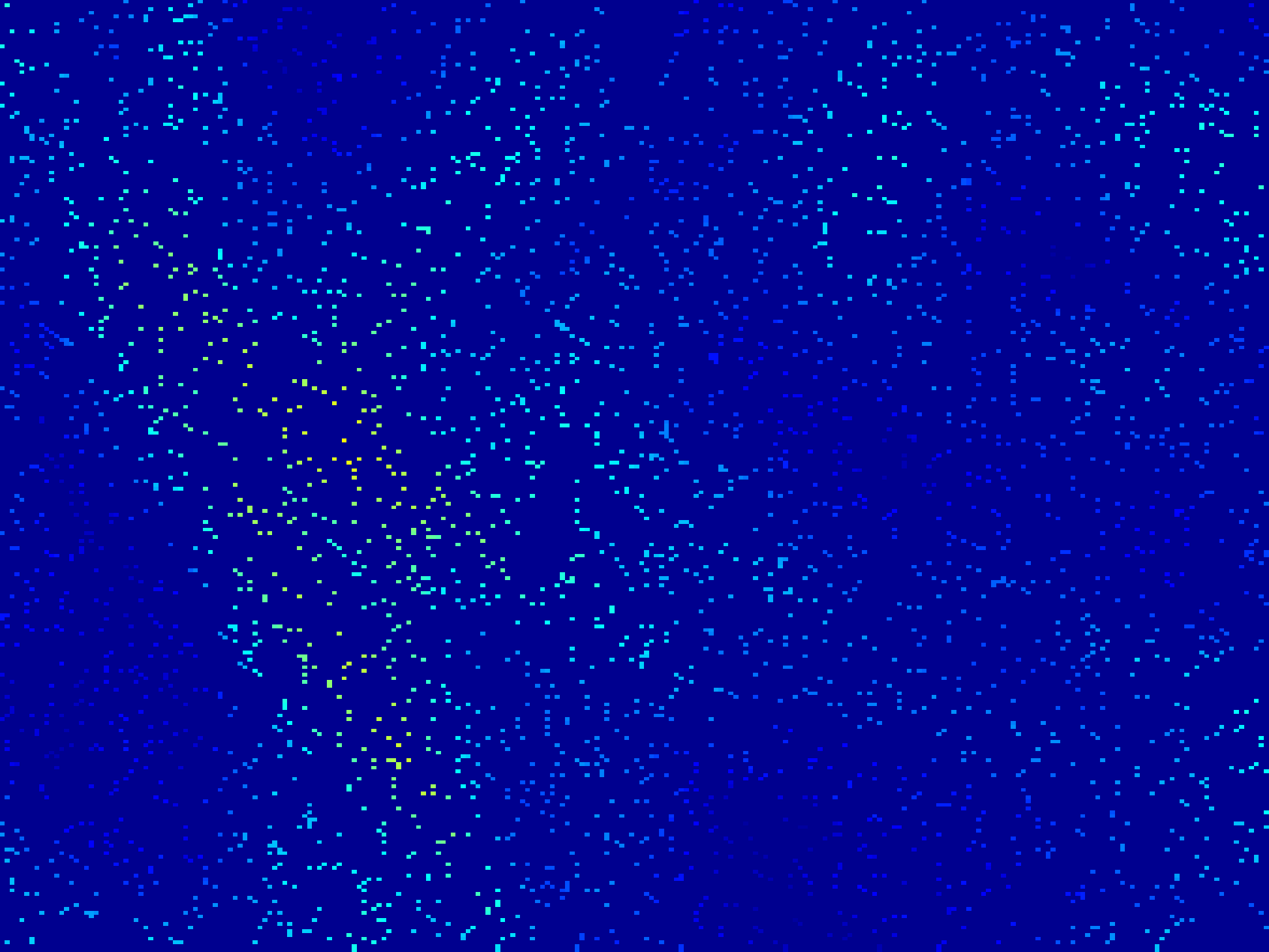}&    
    \includegraphics[width = 2.5cm]{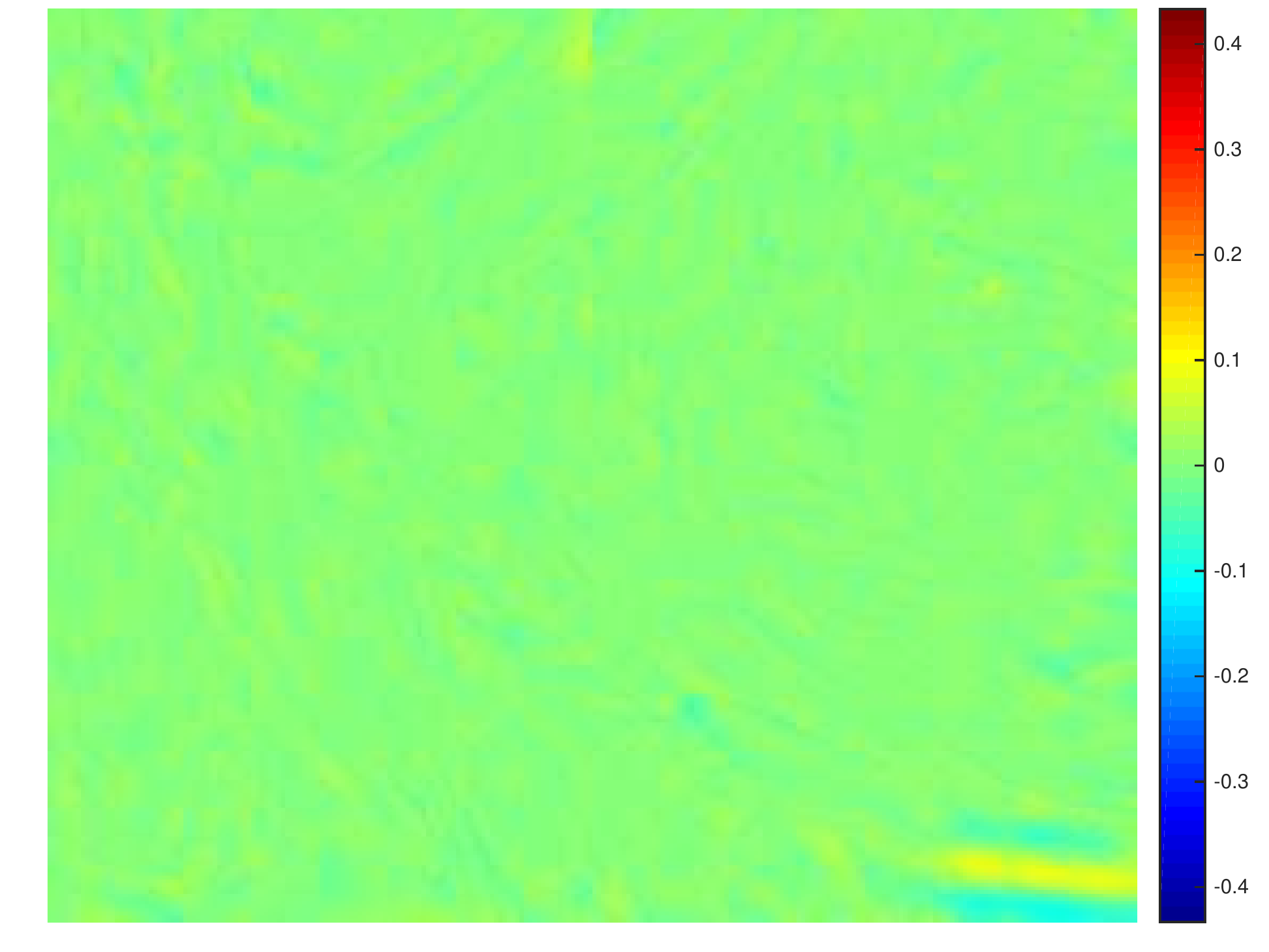}&    
    \includegraphics[width = 2.5cm]{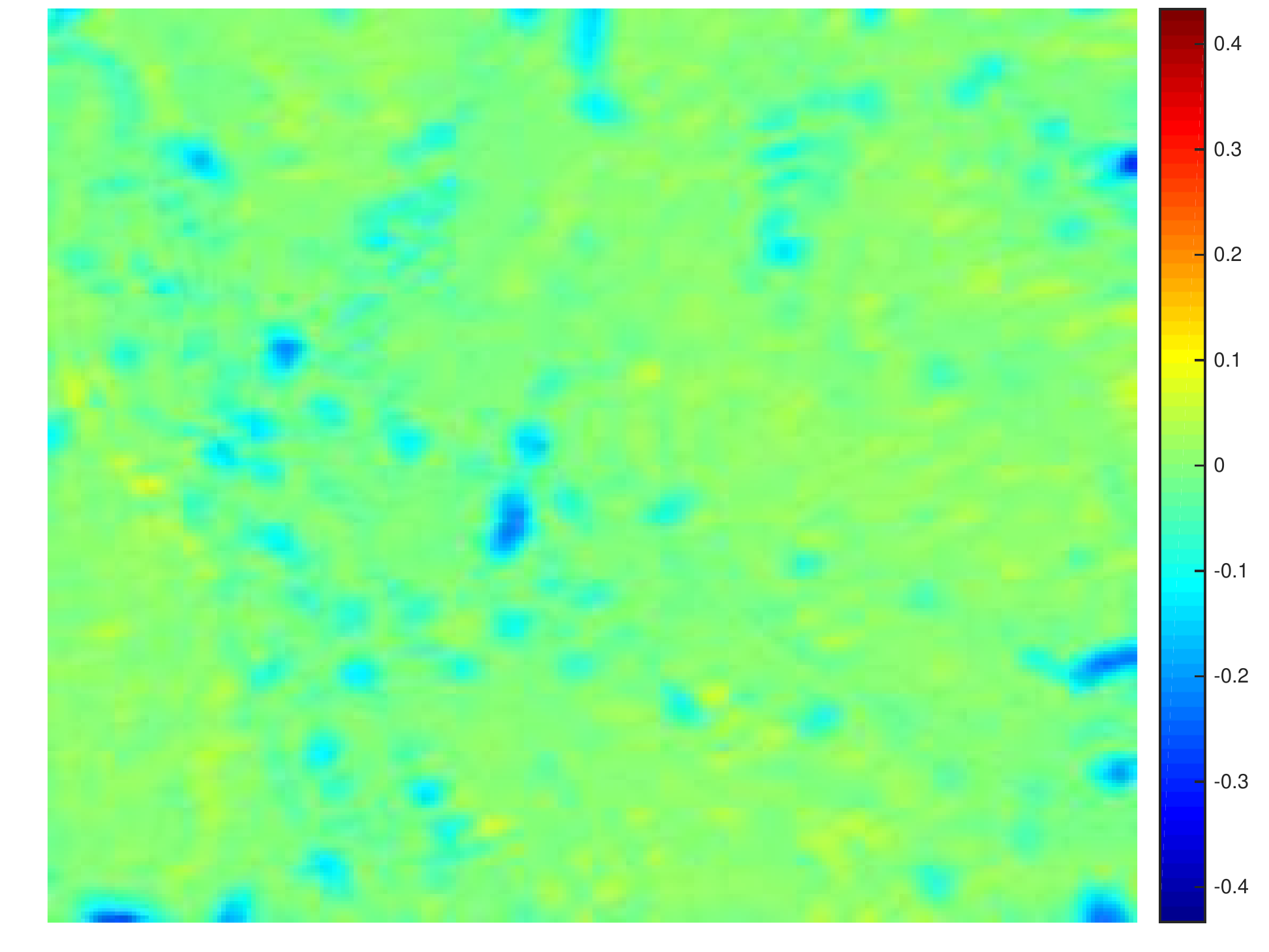}&
    \includegraphics[width = 2.5cm]{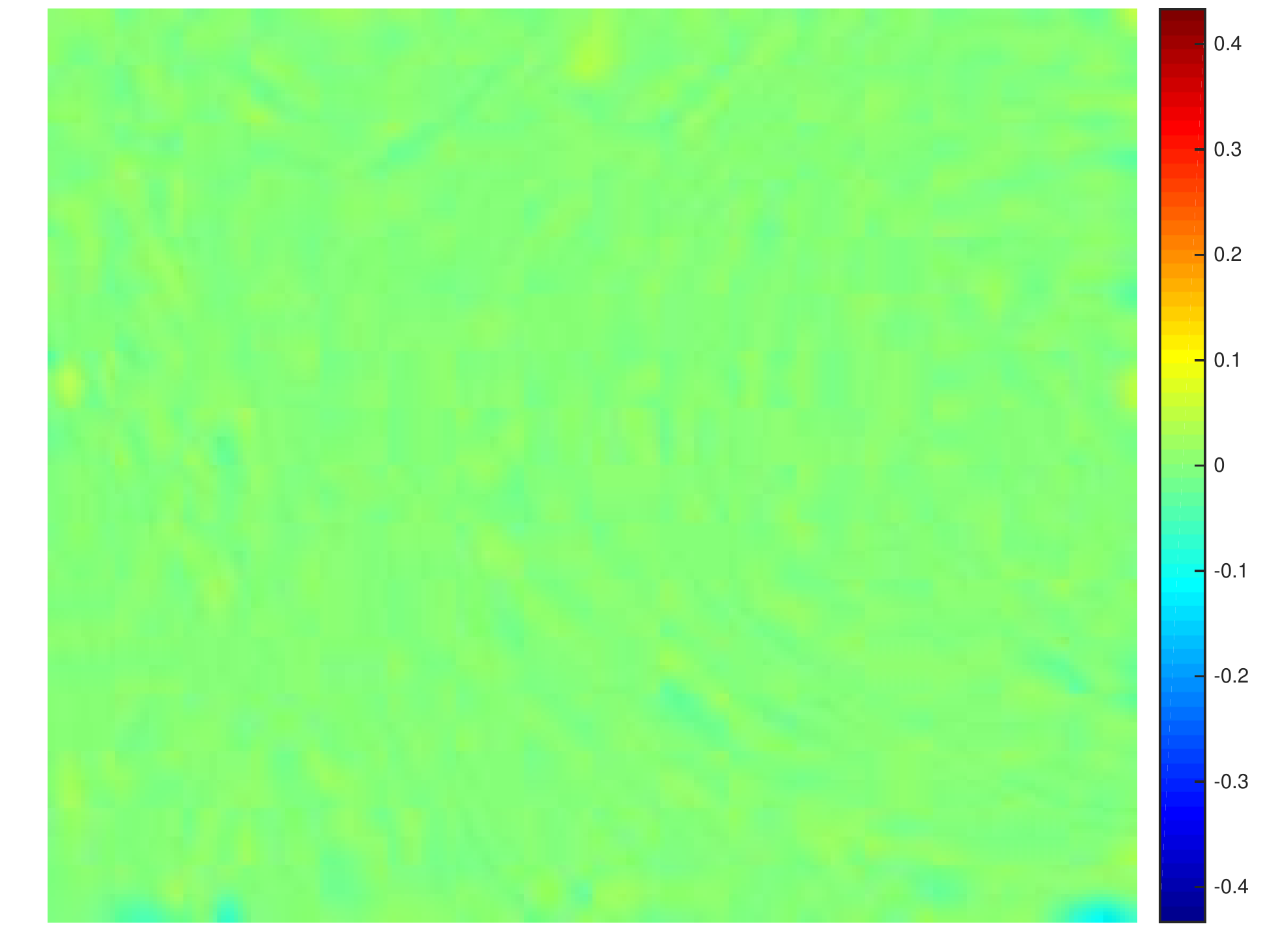}\\
  \end{tabular}
  \caption{Interpolation of the 3D plasma (magnetic field) data set from $5\%$ random sampling. The figures in the first column are two spatial cross sections of the original and subsampled data. The figures in the other three columns are the results and errors of the competing algorithms.}
  \label{fig:result_random_plasma_3d_5p}
\end{figure}

\begin{figure}[H]
  \centering
  \begin{tabular}{cccc}
    Original& EBI (30.24dB)& PLE (35.60dB)  & LDMM (\textbf{48.43dB})\\
    \includegraphics[width = 2.5cm]{lattice_3d_original_band_31}&
    \includegraphics[width = 2.5cm]{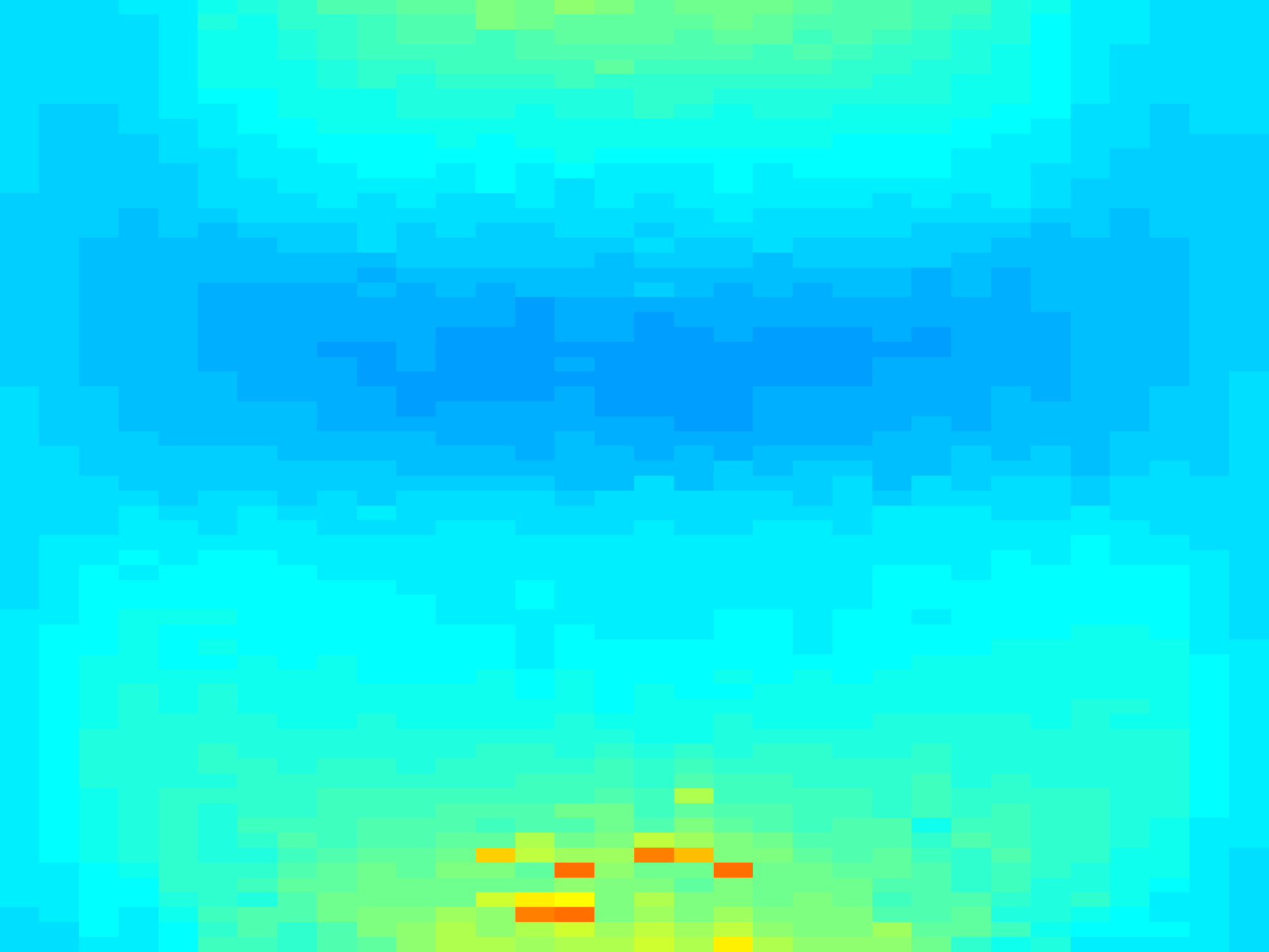}&
    \includegraphics[width = 2.5cm]{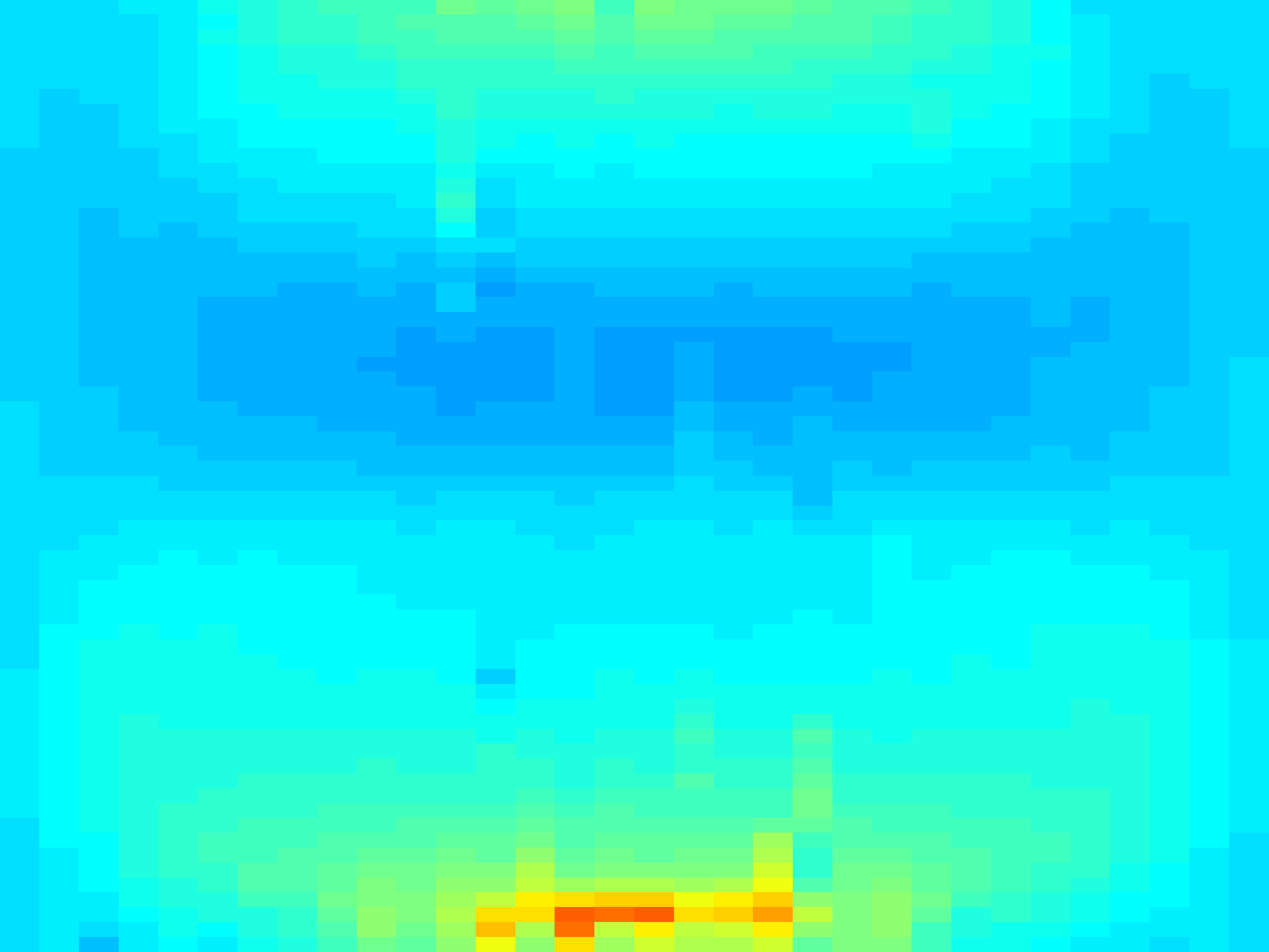}&
    \includegraphics[width = 2.5cm]{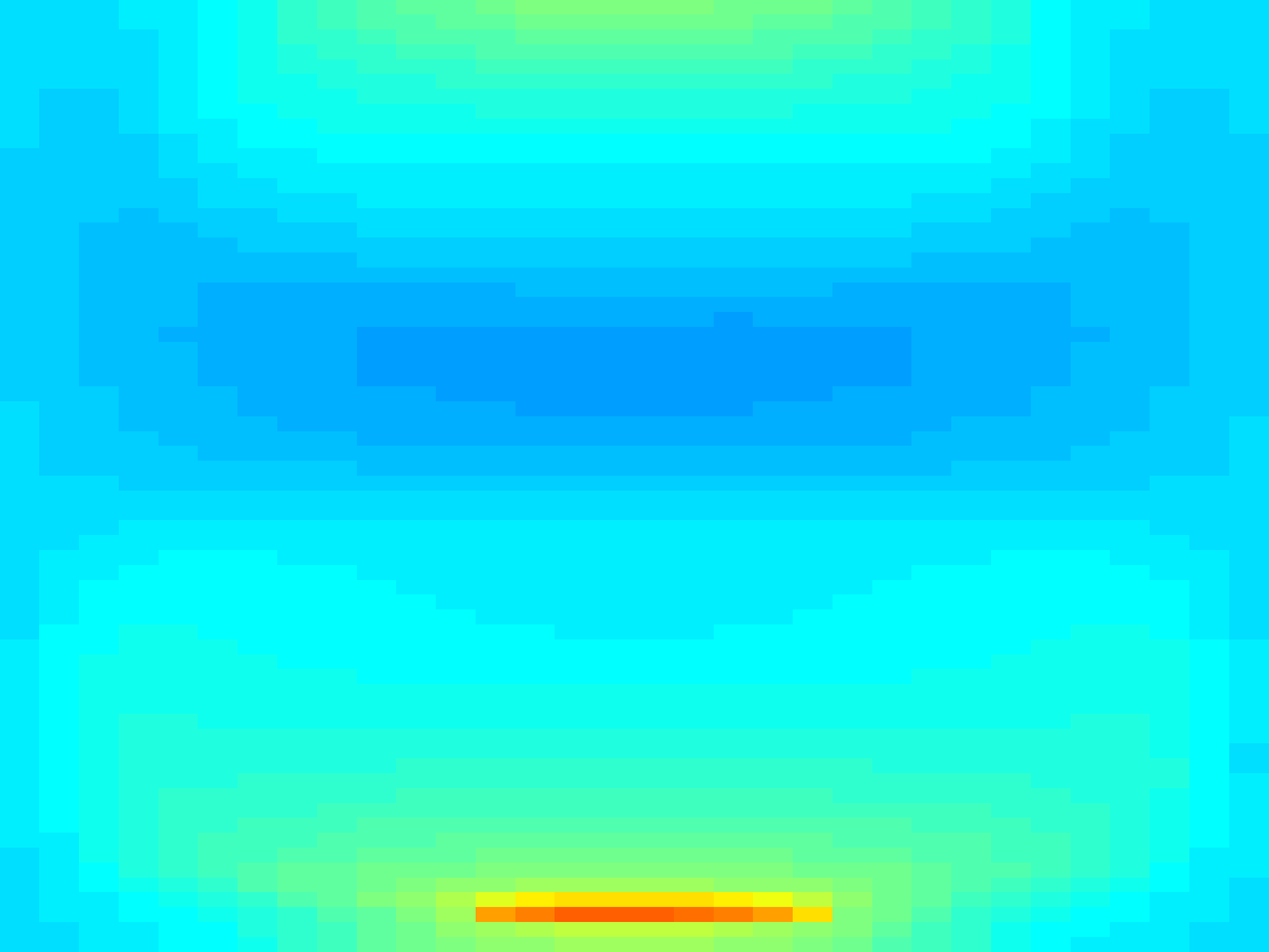}\\
    Subsample & Error & Error & Error\\
    \includegraphics[width = 2.5cm]{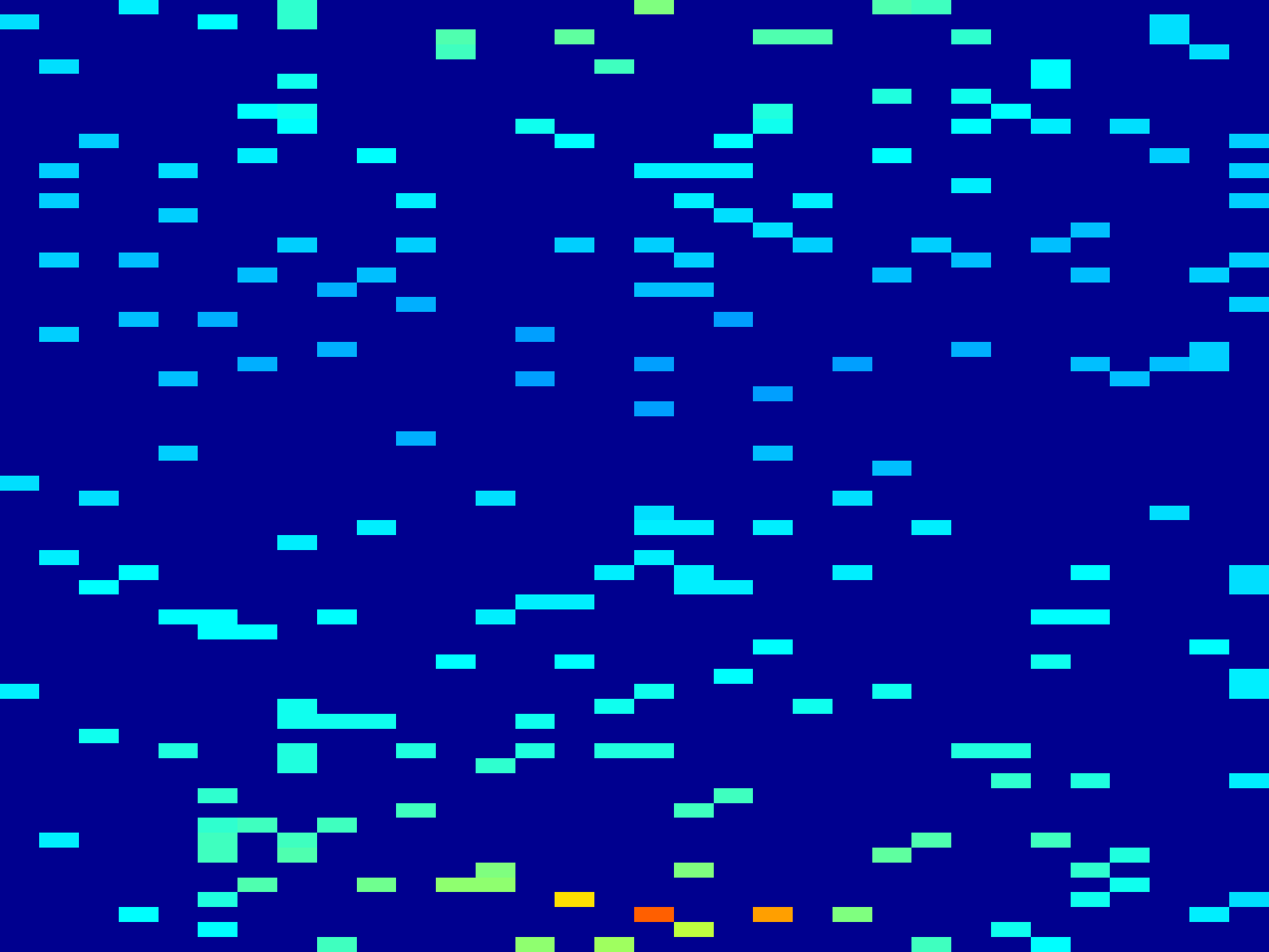}&    
    \includegraphics[width = 2.5cm]{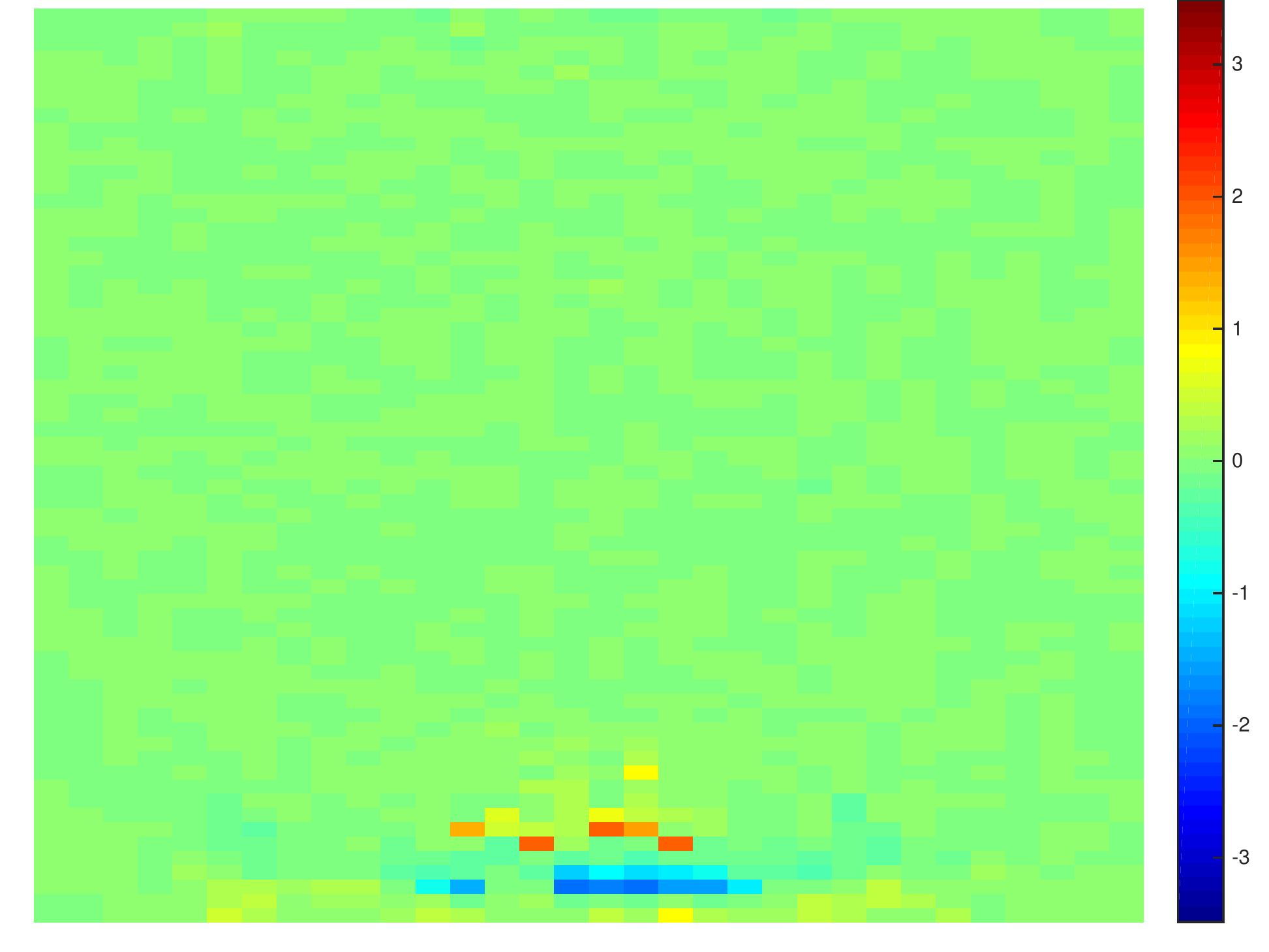}&    
    \includegraphics[width = 2.5cm]{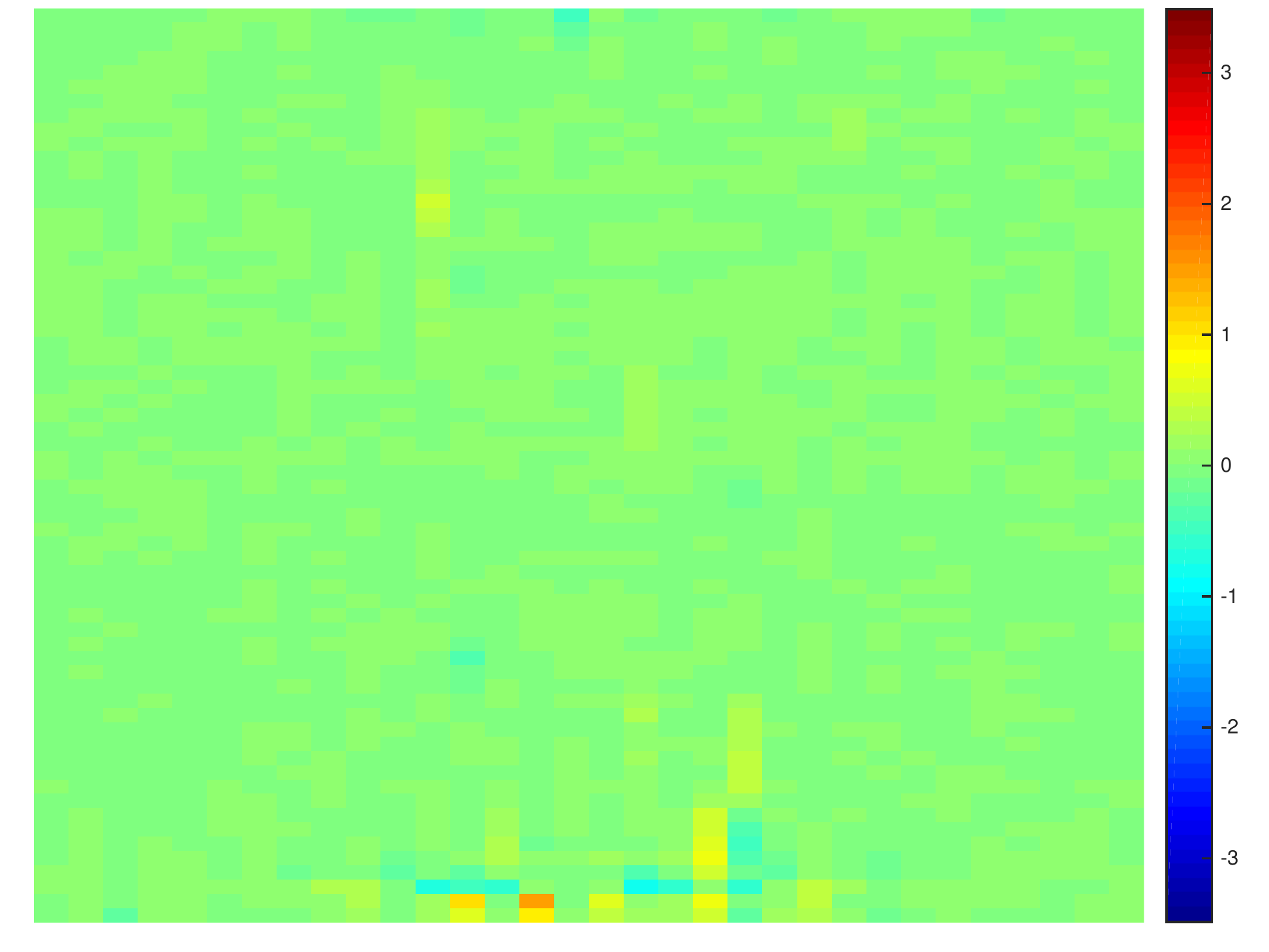}&
    \includegraphics[width = 2.5cm]{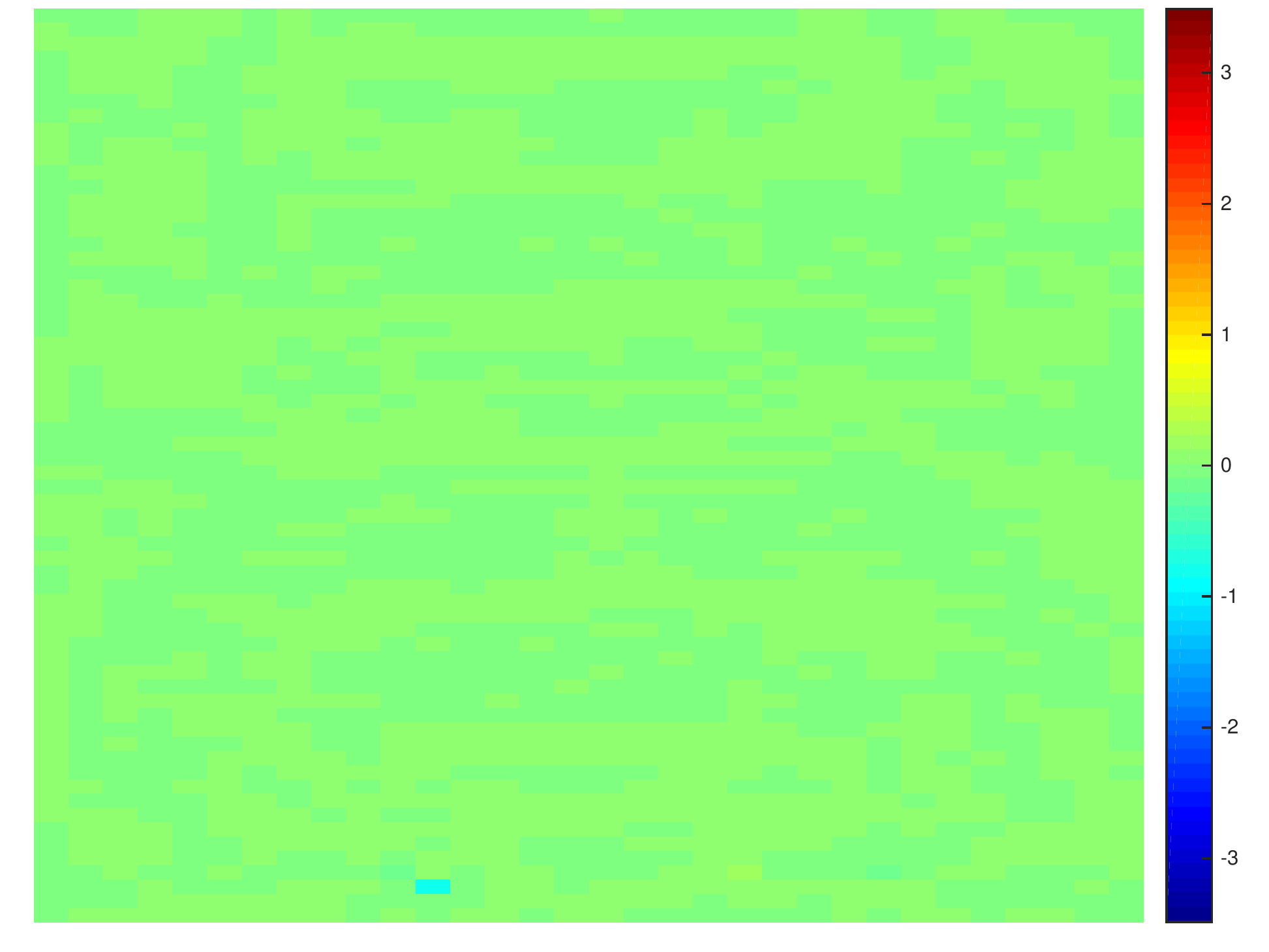}\\
    Original& EBI (30.24dB)& PLE (35.60dB)  & LDMM (\textbf{48.43dB})\\
    \includegraphics[width = 2.5cm]{lattice_3d_original_band_151}&
    \includegraphics[width = 2.5cm]{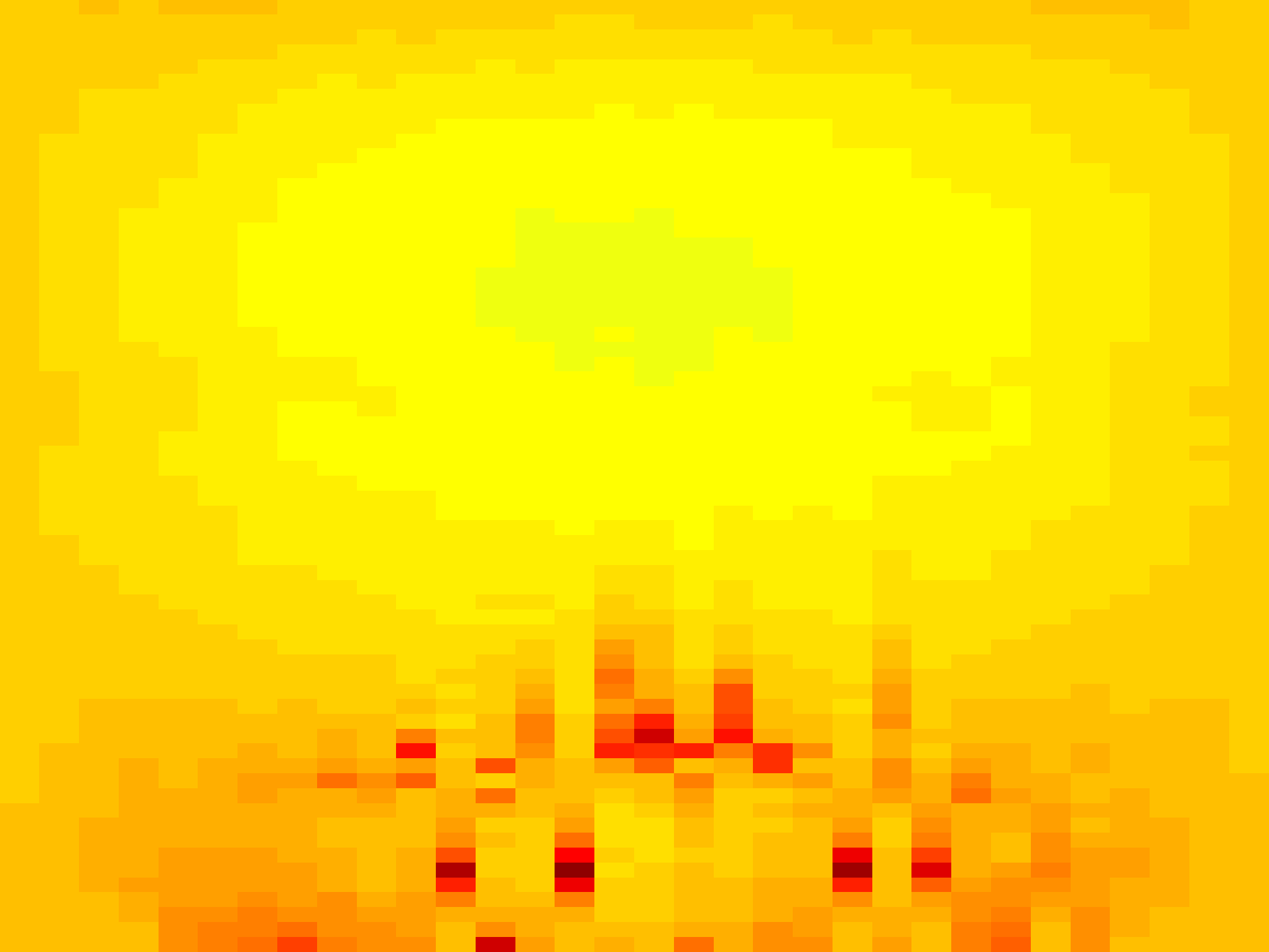}&
    \includegraphics[width = 2.5cm]{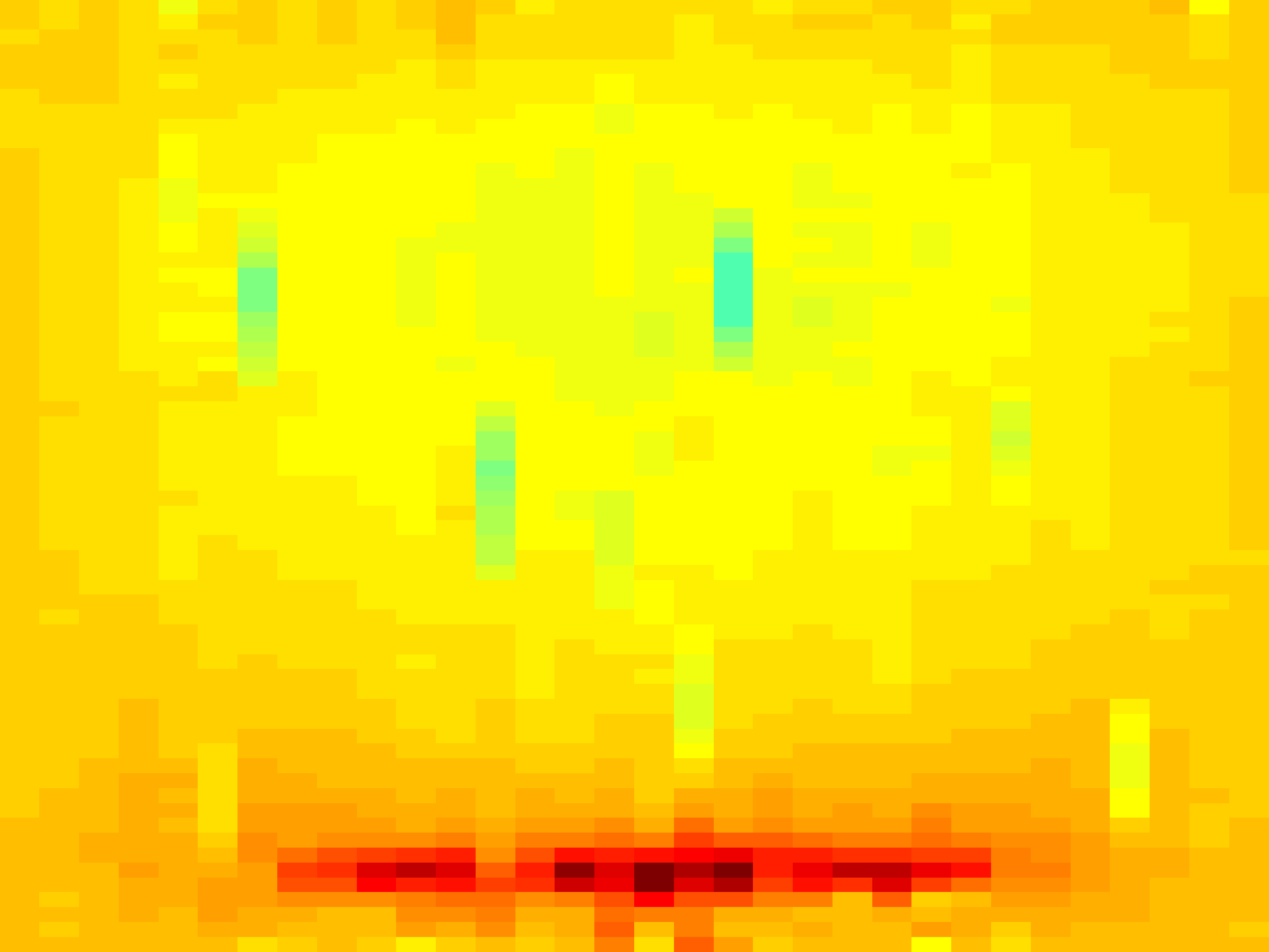}&
    \includegraphics[width = 2.5cm]{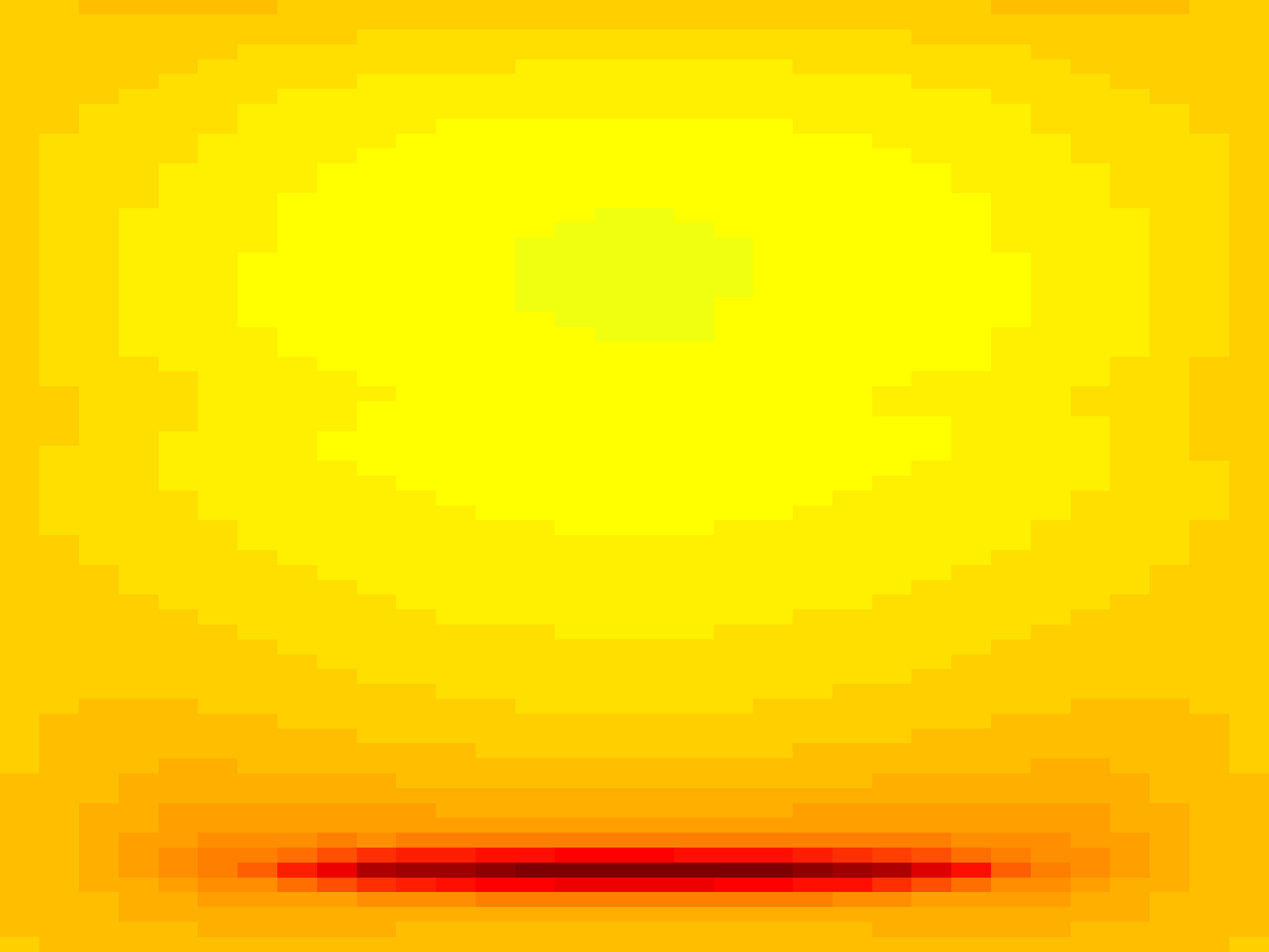}\\
    Subsample & Error & Error & Error\\
    \includegraphics[width = 2.5cm]{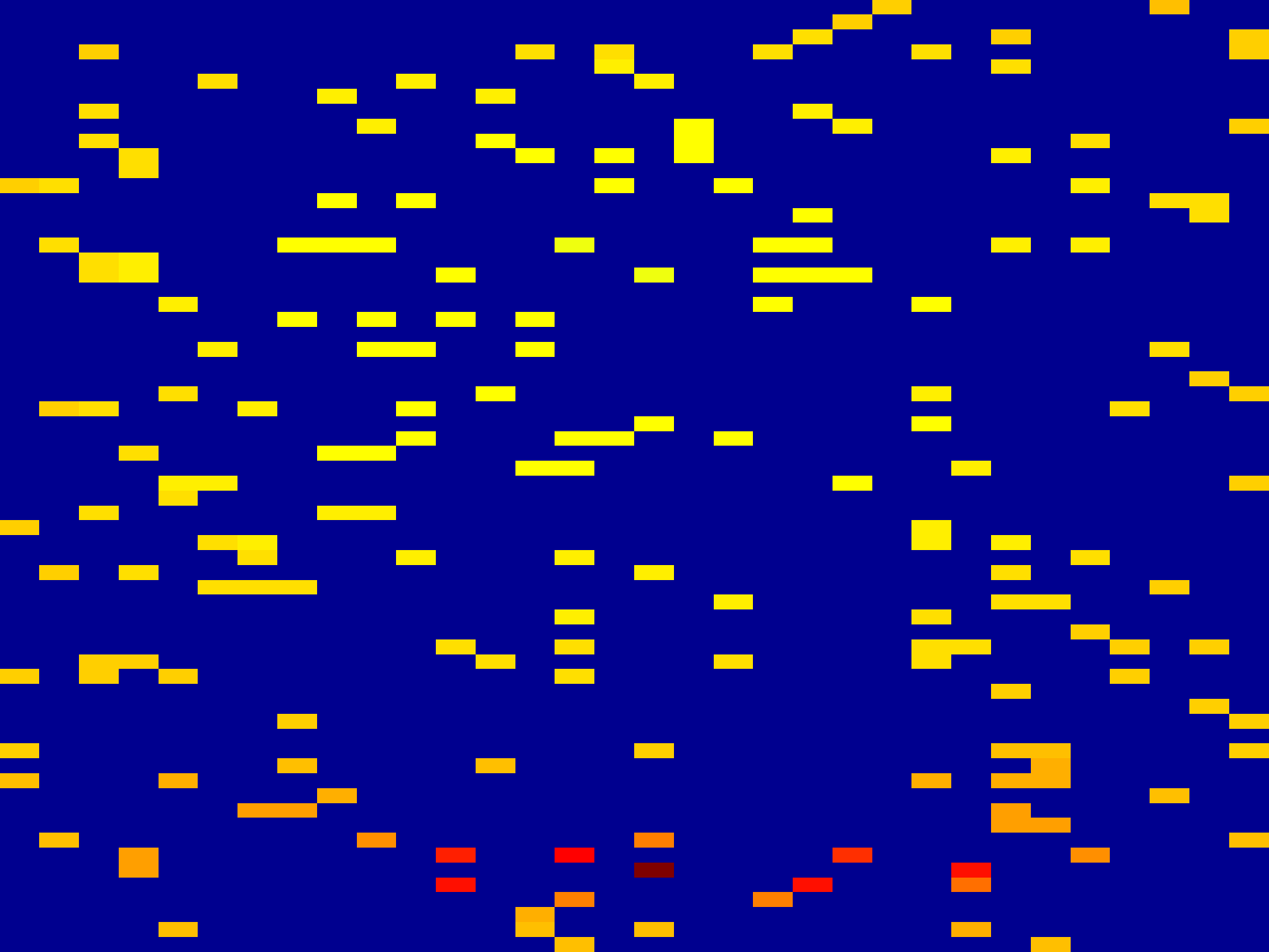}&    
    \includegraphics[width = 2.5cm]{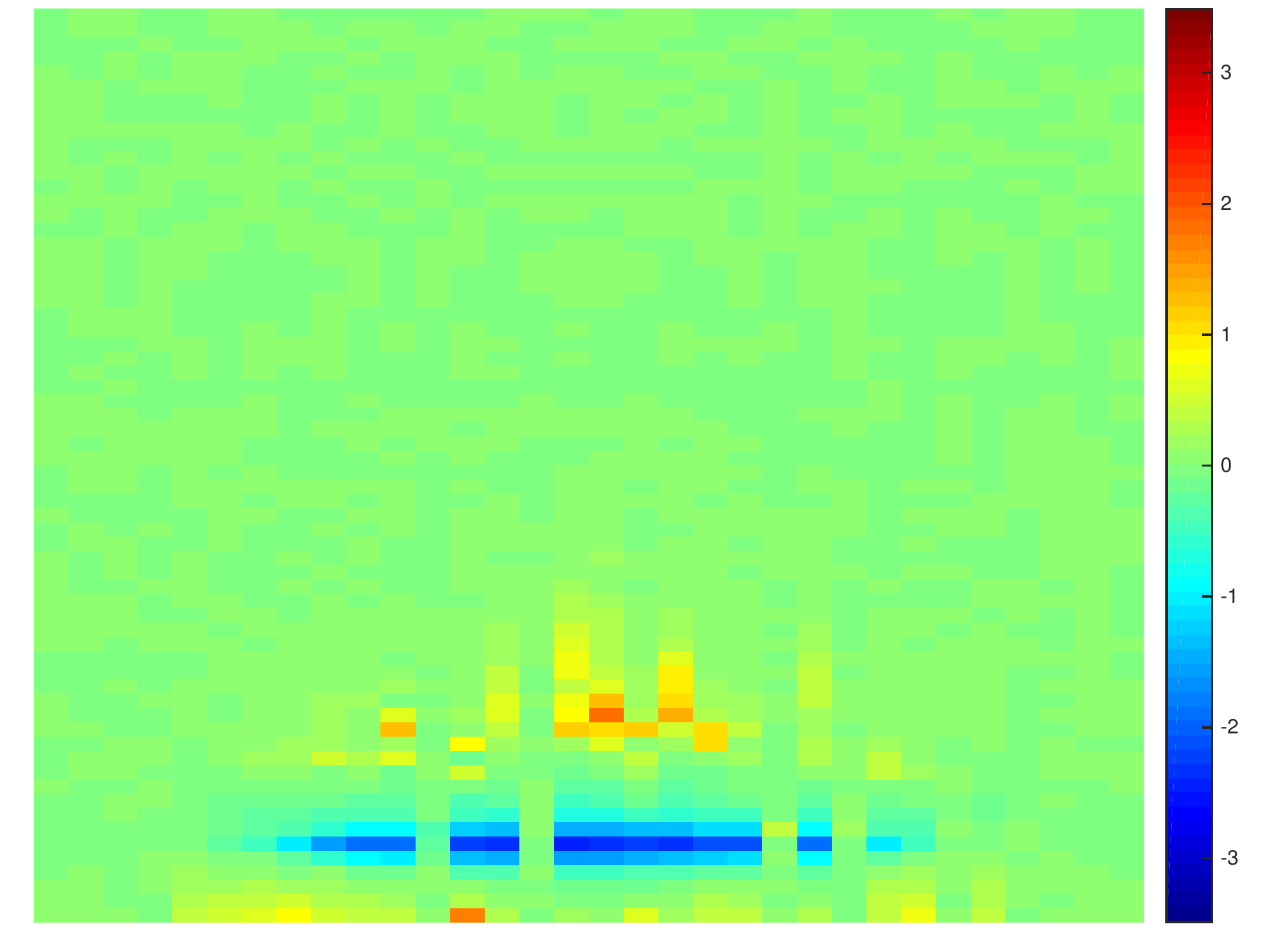}&    
    \includegraphics[width = 2.5cm]{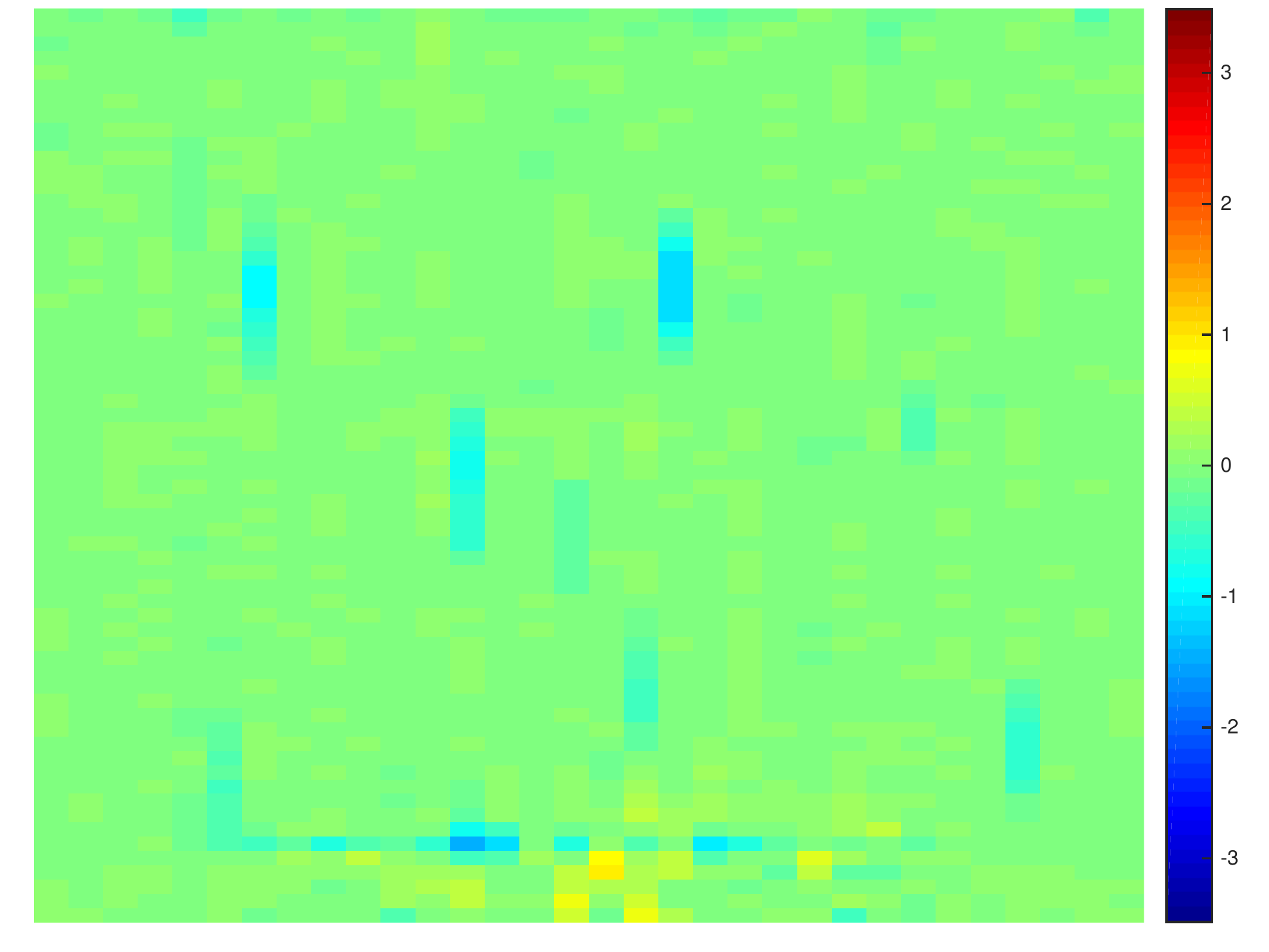}&
    \includegraphics[width = 2.5cm]{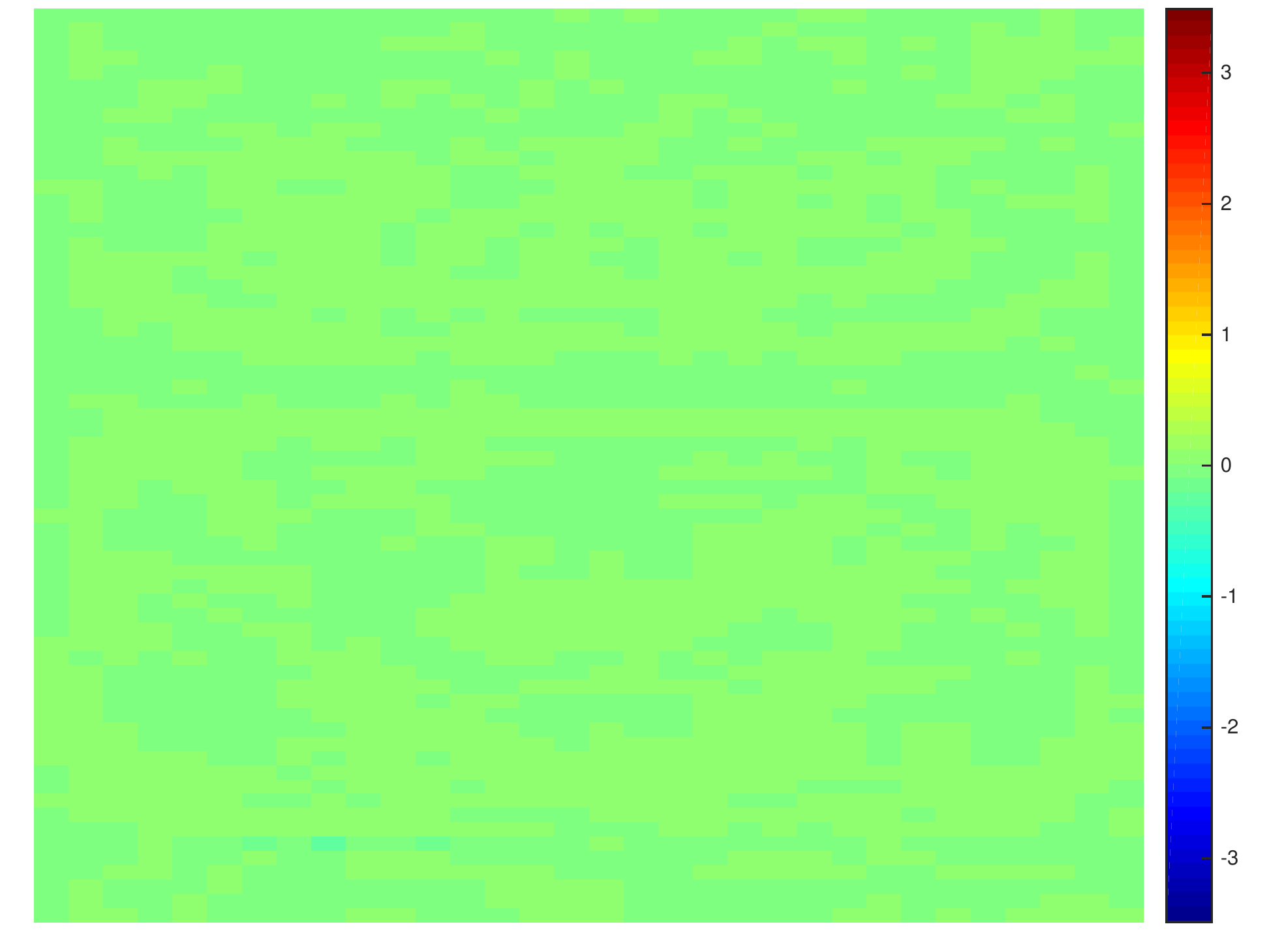}\\
  \end{tabular}
  \caption{Interpolation of the 3D lattice data set from $10\%$ random sampling. The figures in the first column are the original and subsampled angular flux at $x = 0.24$ and $x = 1.18$ . The figures in the other three columns are the results and errors of the competing algorithms.}
  \label{fig:result_random_lattice_3d_10p}
\end{figure}

\begin{figure}[H]
  \centering
  \begin{tabular}{cccc}
    Original& EBI (29.48dB)& PLE (20.93B)  & LDMM (\textbf{45.82dB})\\
    \includegraphics[width = 2.5cm]{lattice_3d_original_band_31}&
    \includegraphics[width = 2.5cm]{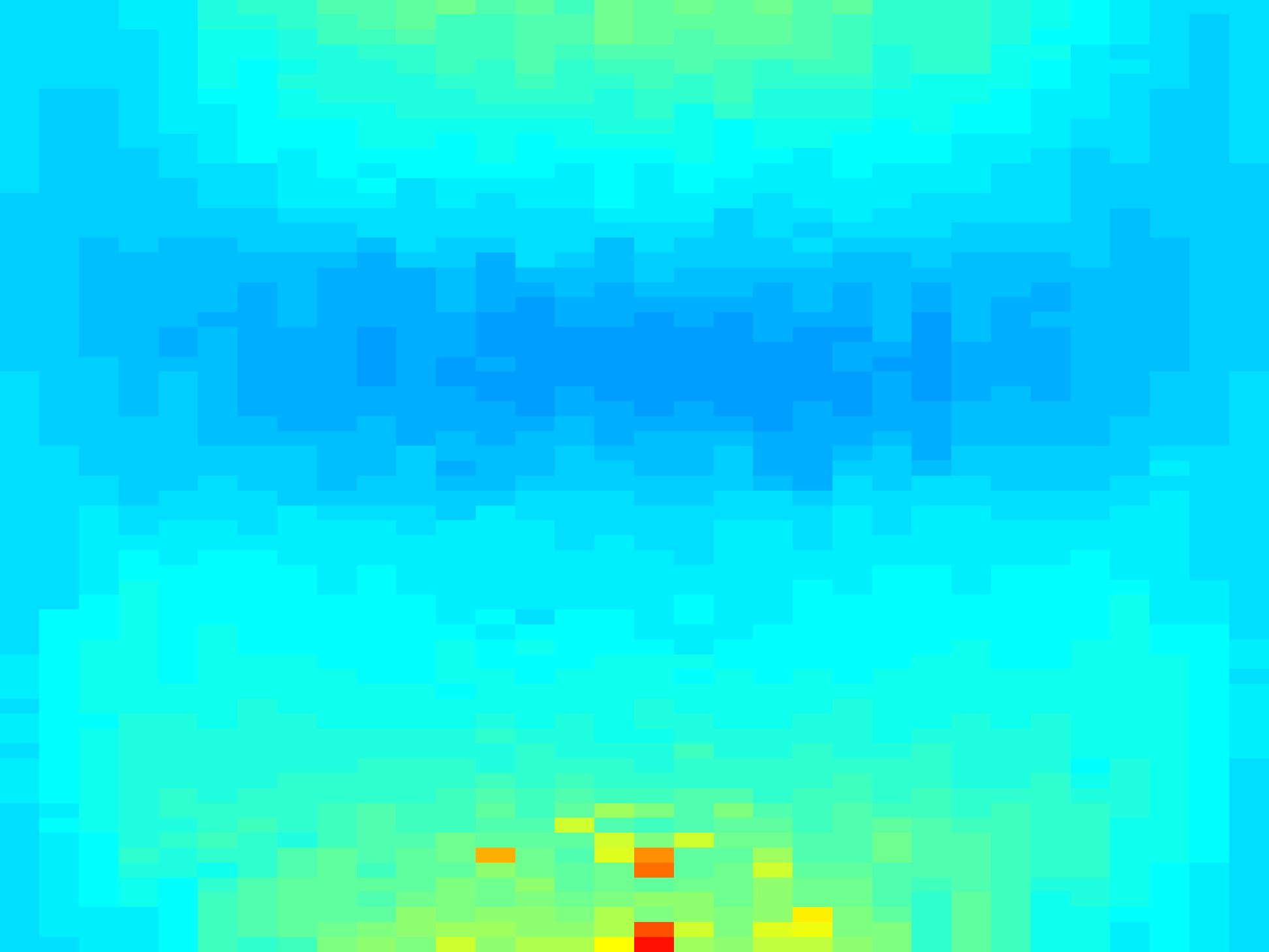}&
    \includegraphics[width = 2.5cm]{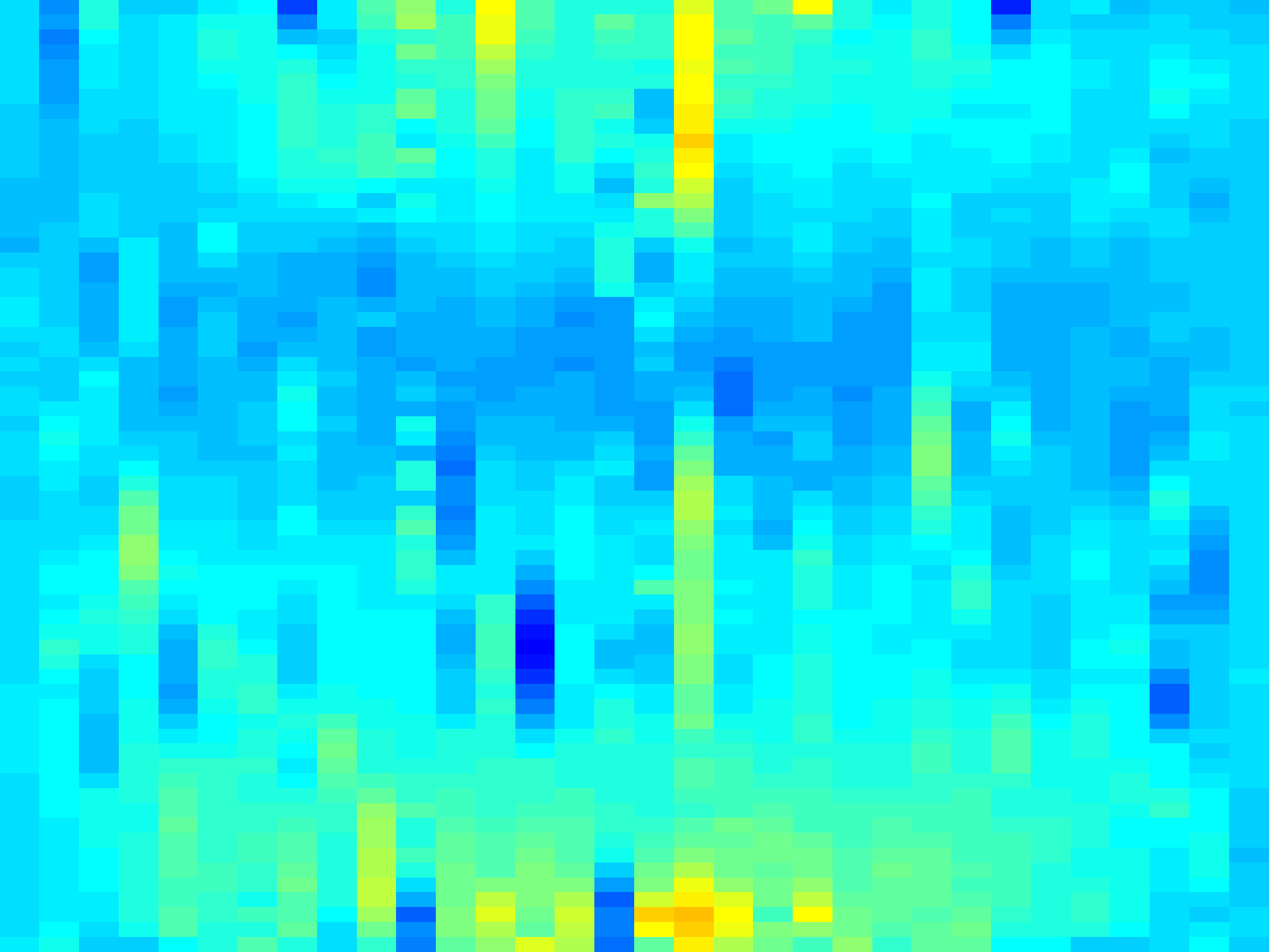}&
    \includegraphics[width = 2.5cm]{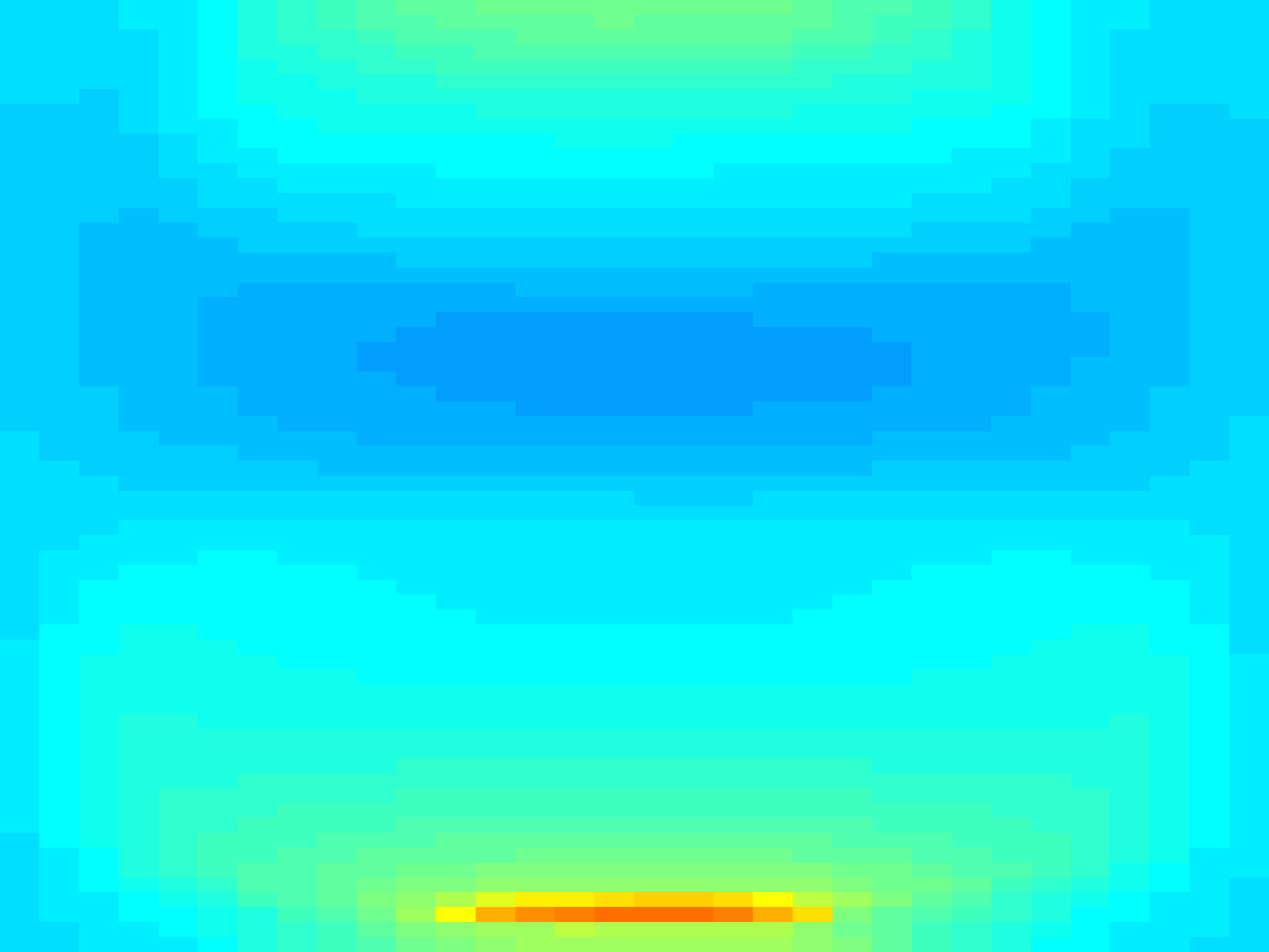}\\
    Subsample & Error & Error & Error\\
    \includegraphics[width = 2.5cm]{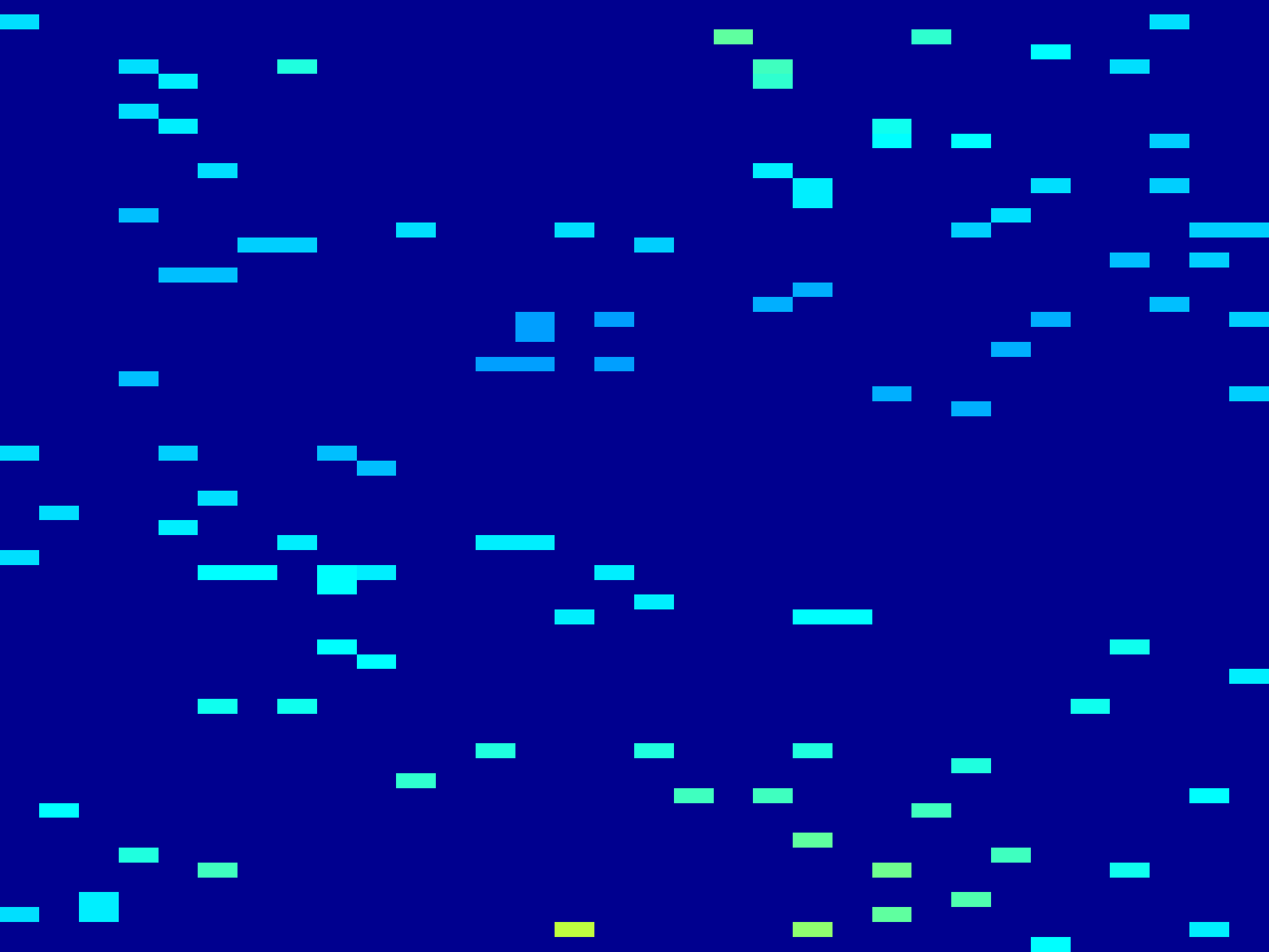}&    
    \includegraphics[width = 2.5cm]{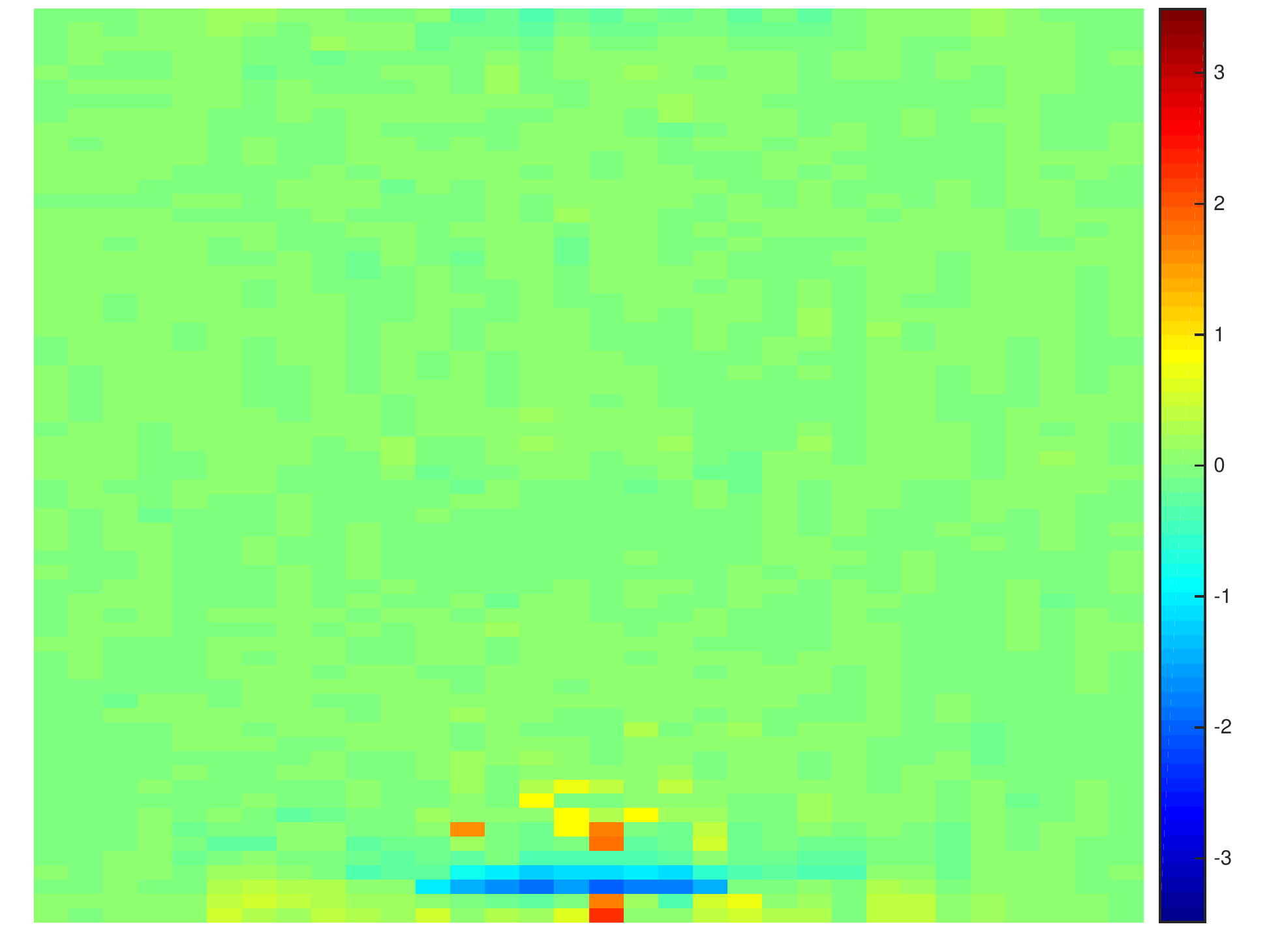}&    
    \includegraphics[width = 2.5cm]{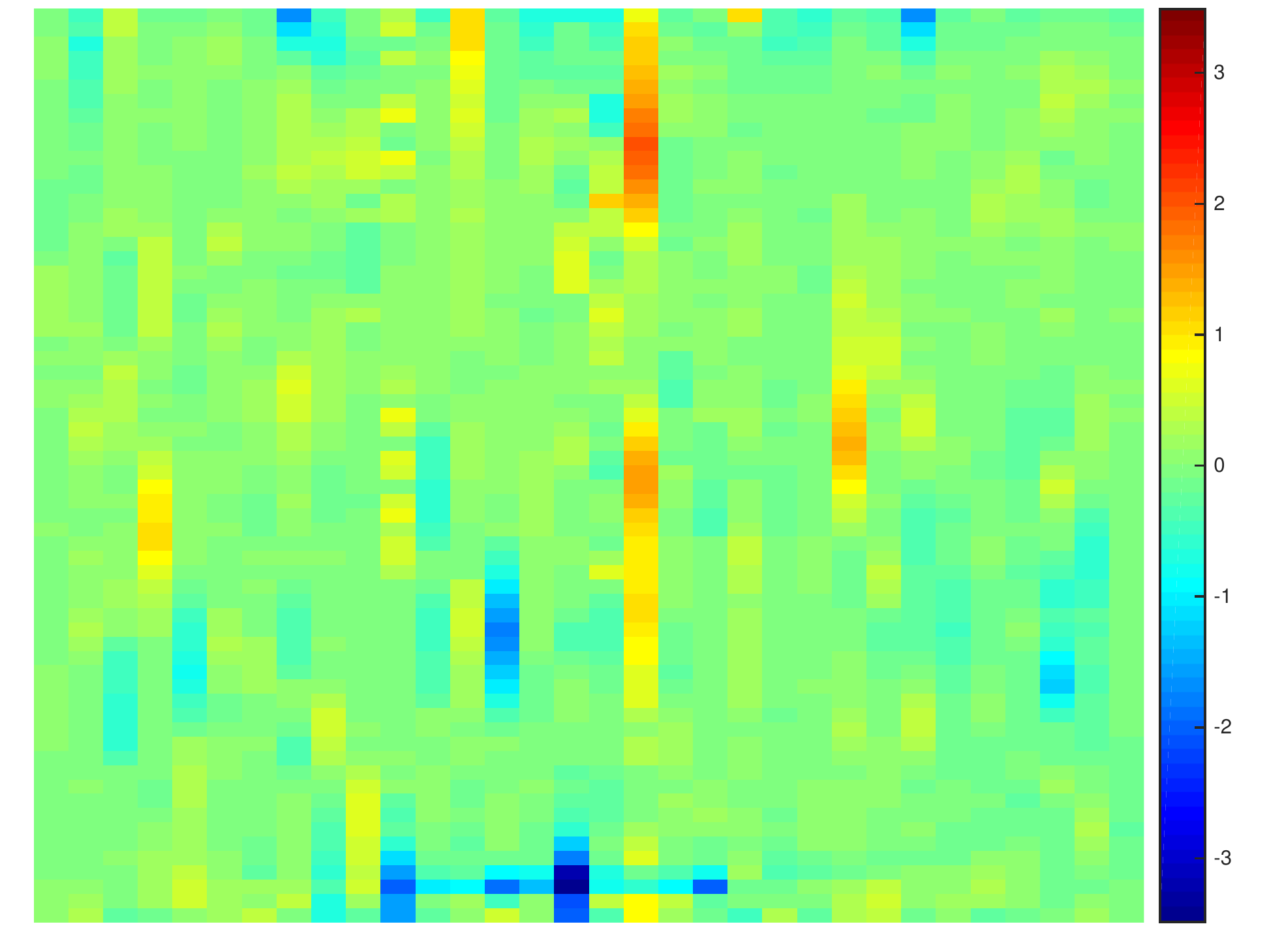}&
    \includegraphics[width = 2.5cm]{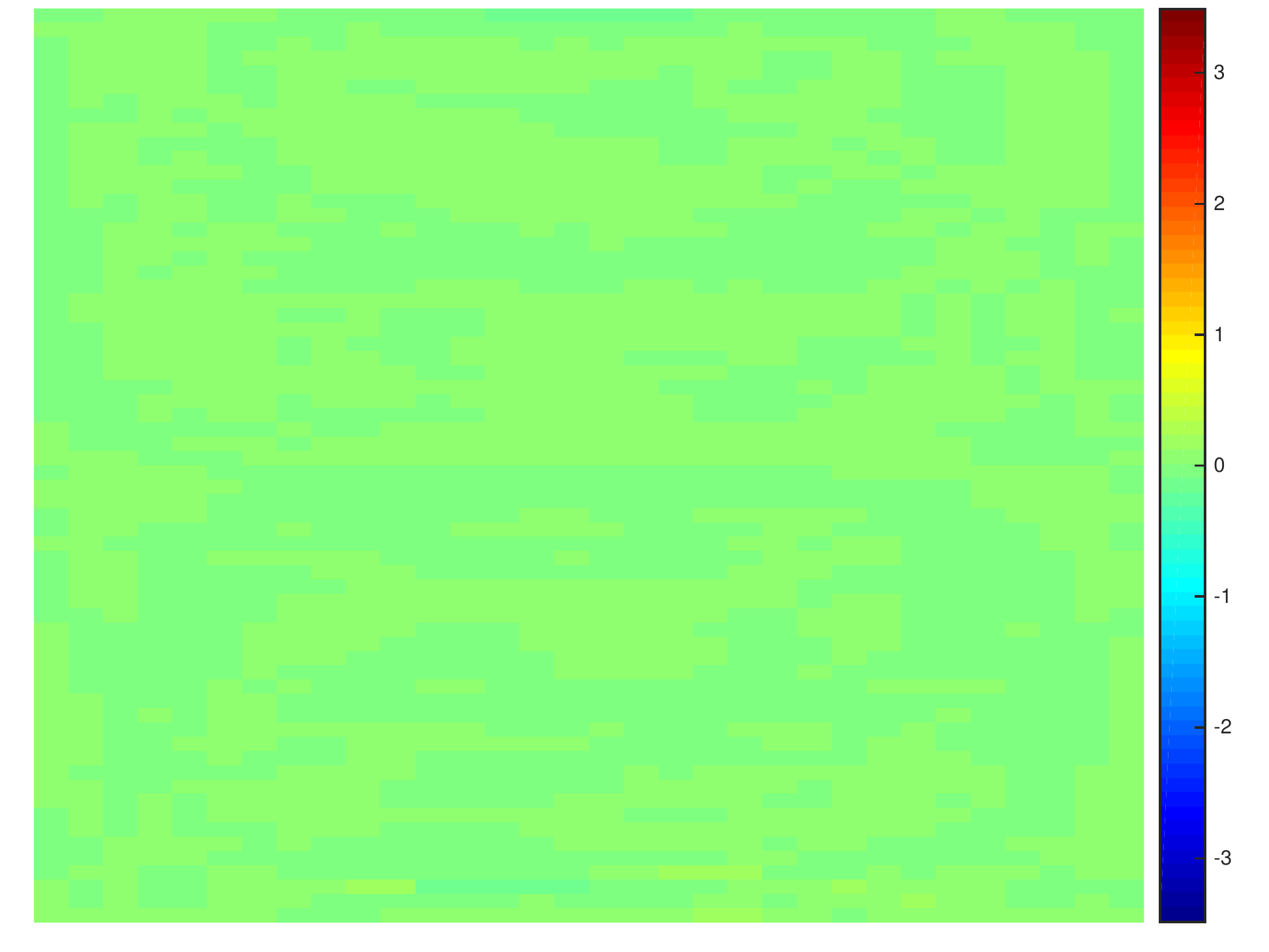}\\
    Original& EBI (29.48dB)& PLE (20.93B)  & LDMM (\textbf{45.82dB})\\
    \includegraphics[width = 2.5cm]{lattice_3d_original_band_151}&
    \includegraphics[width = 2.5cm]{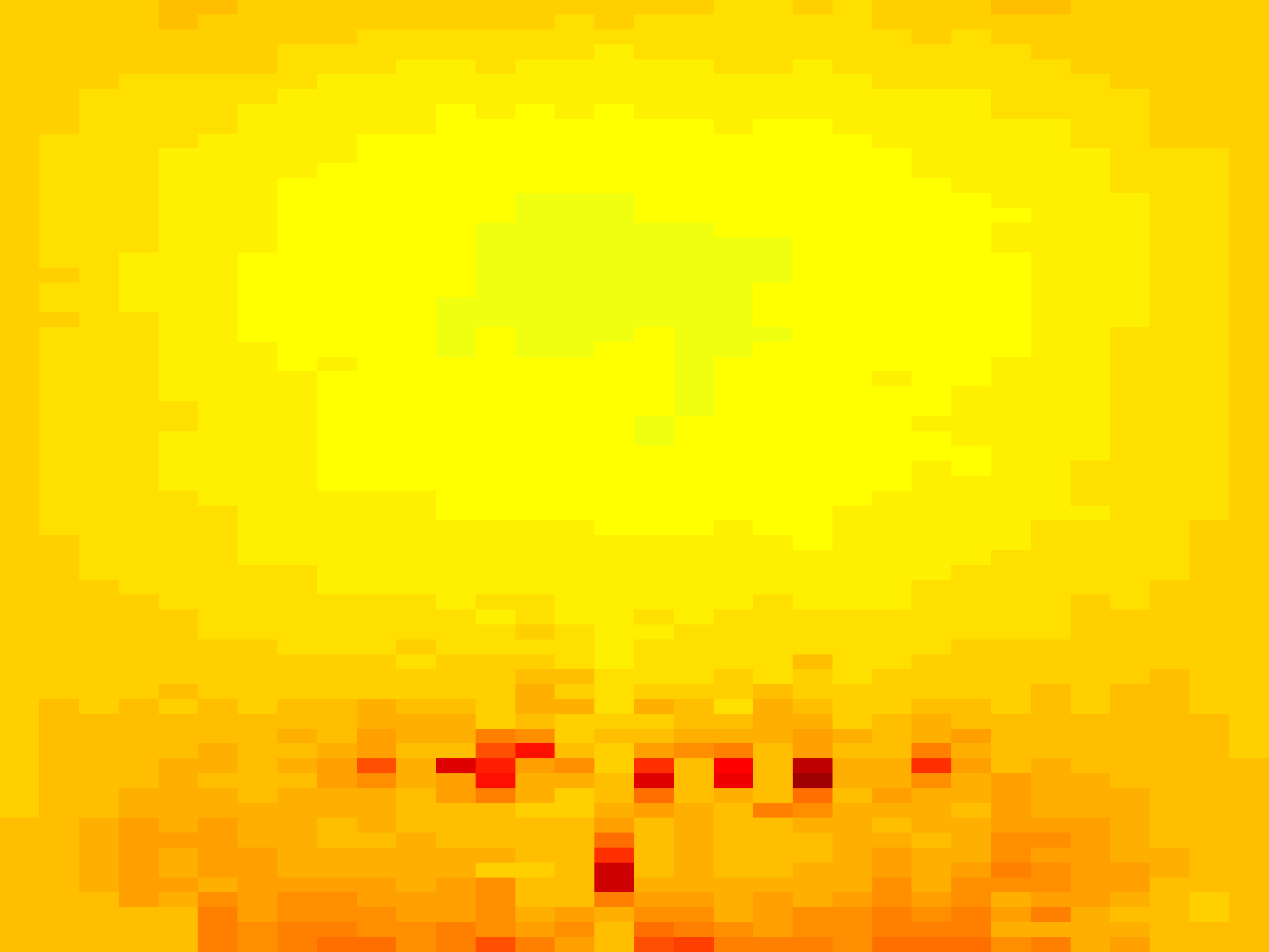}&
    \includegraphics[width = 2.5cm]{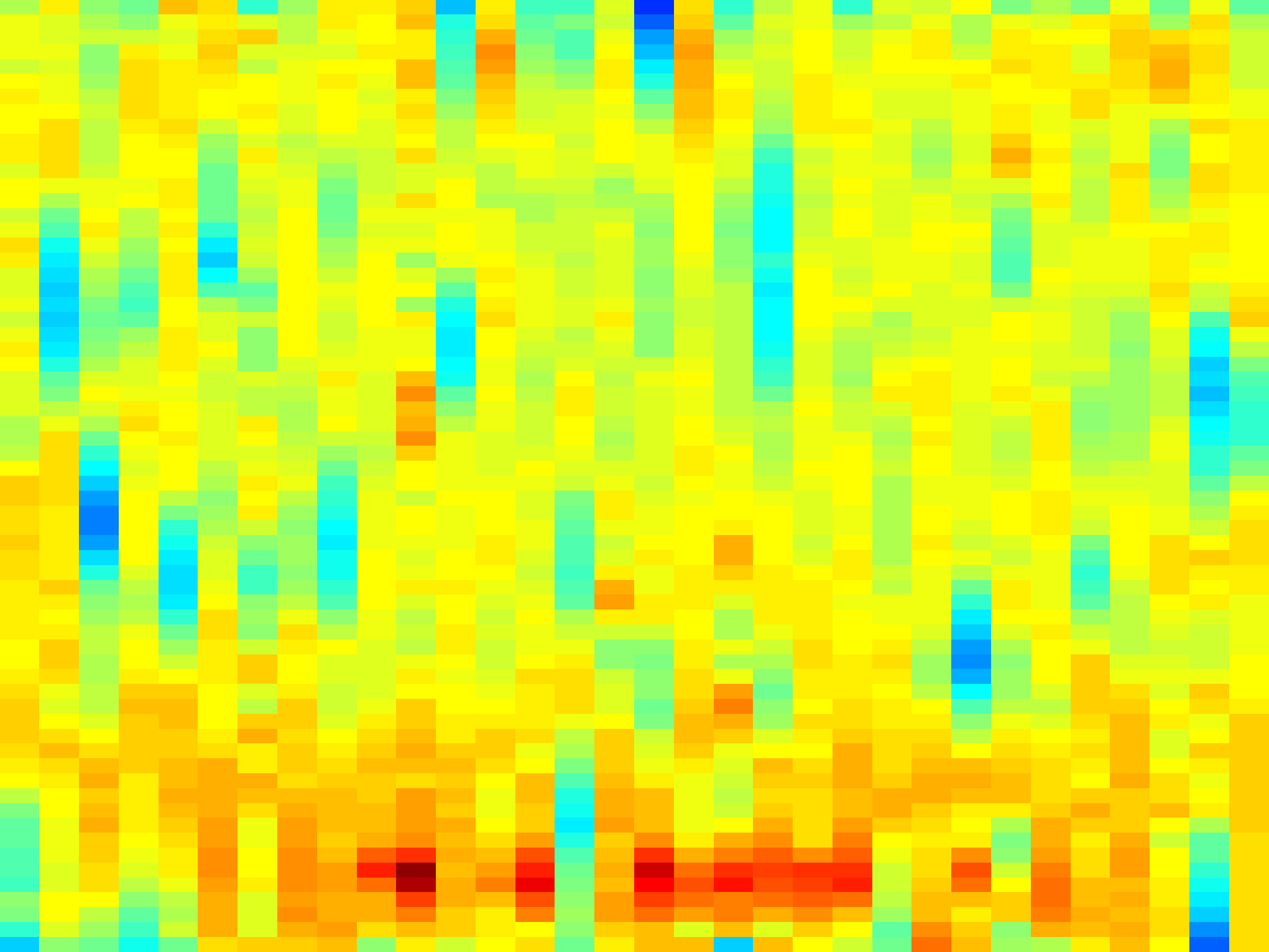}&
    \includegraphics[width = 2.5cm]{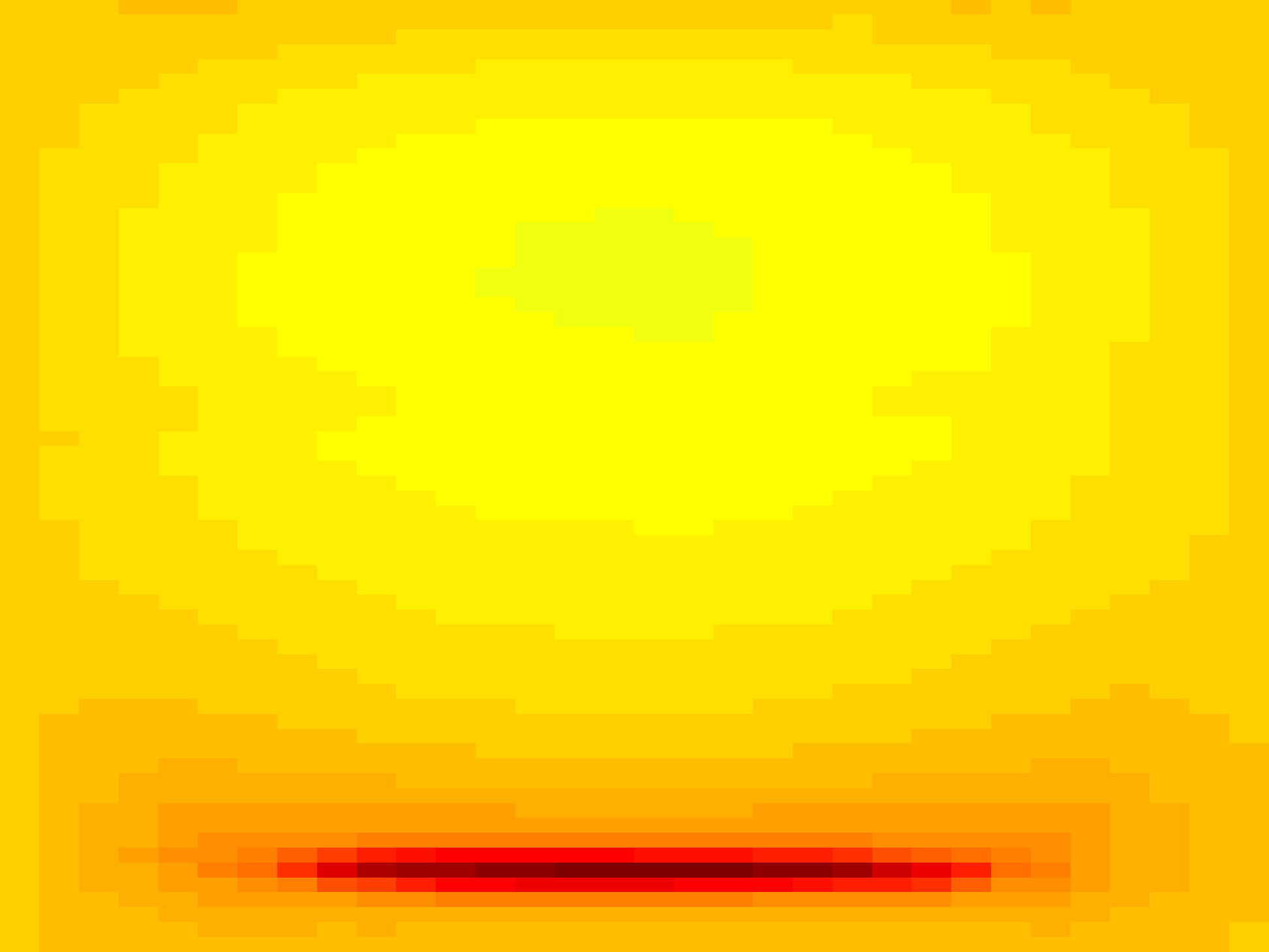}\\
    Subsample & Error & Error & Error\\
    \includegraphics[width = 2.5cm]{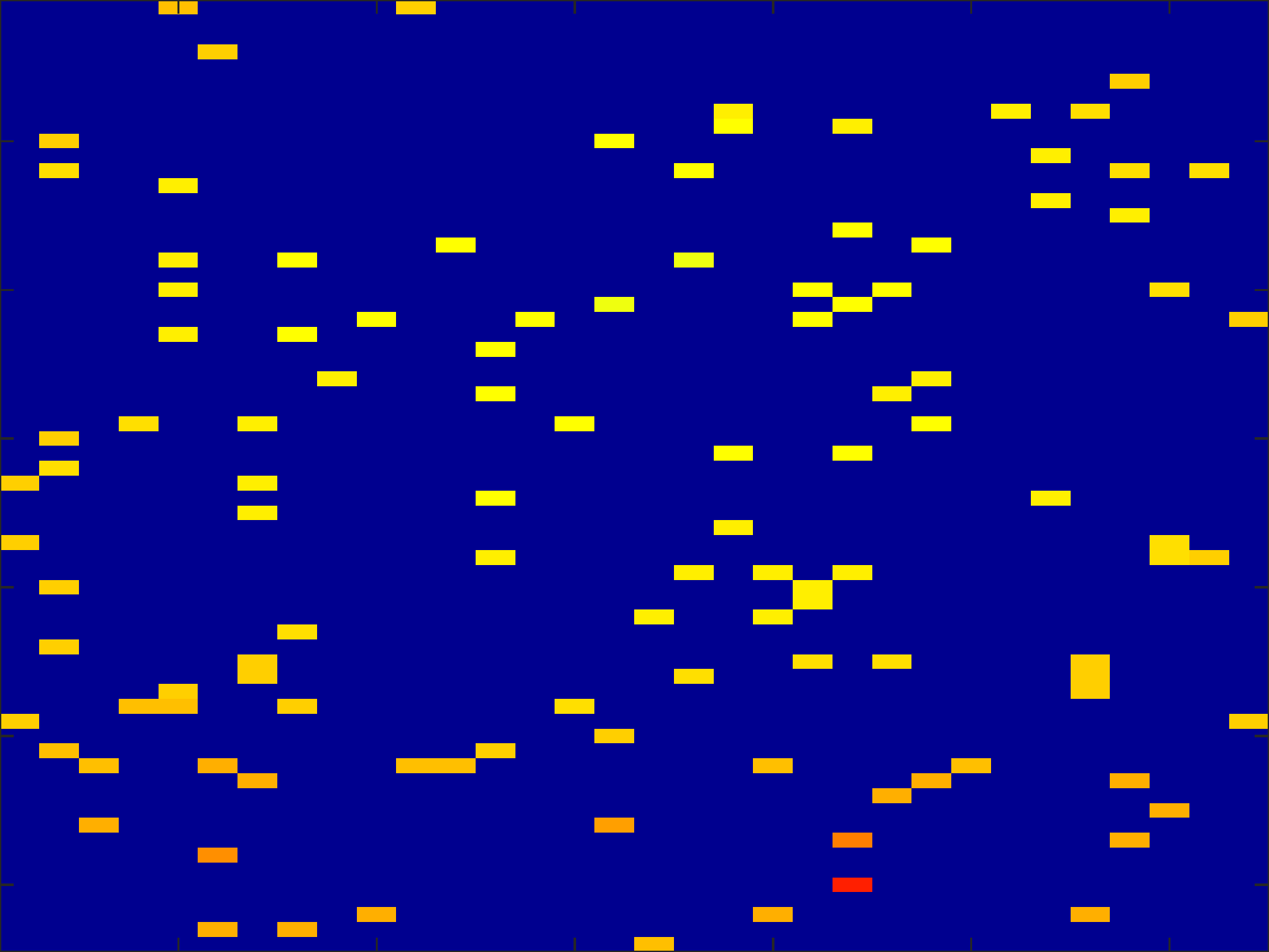}&    
    \includegraphics[width = 2.5cm]{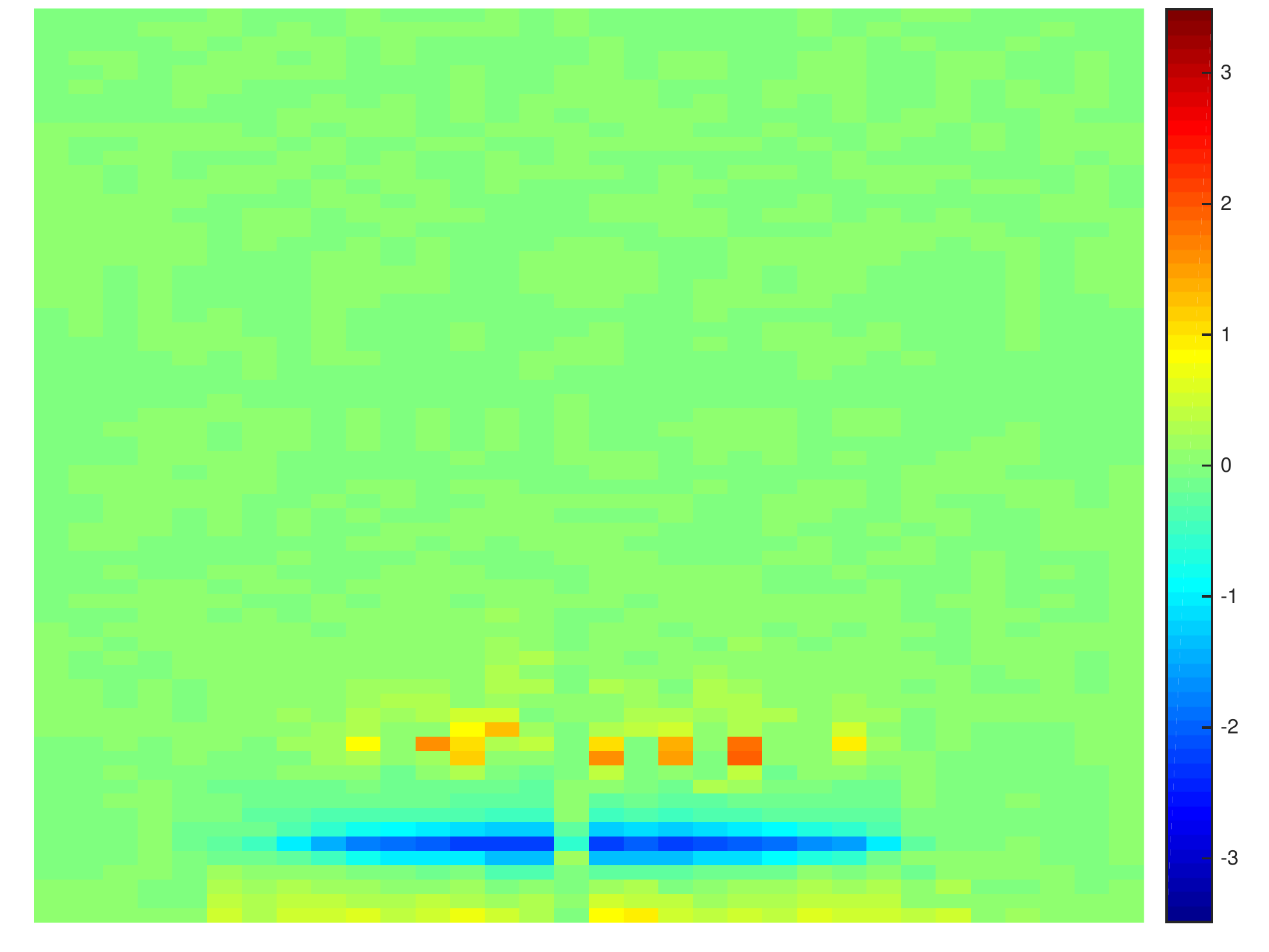}&    
    \includegraphics[width = 2.5cm]{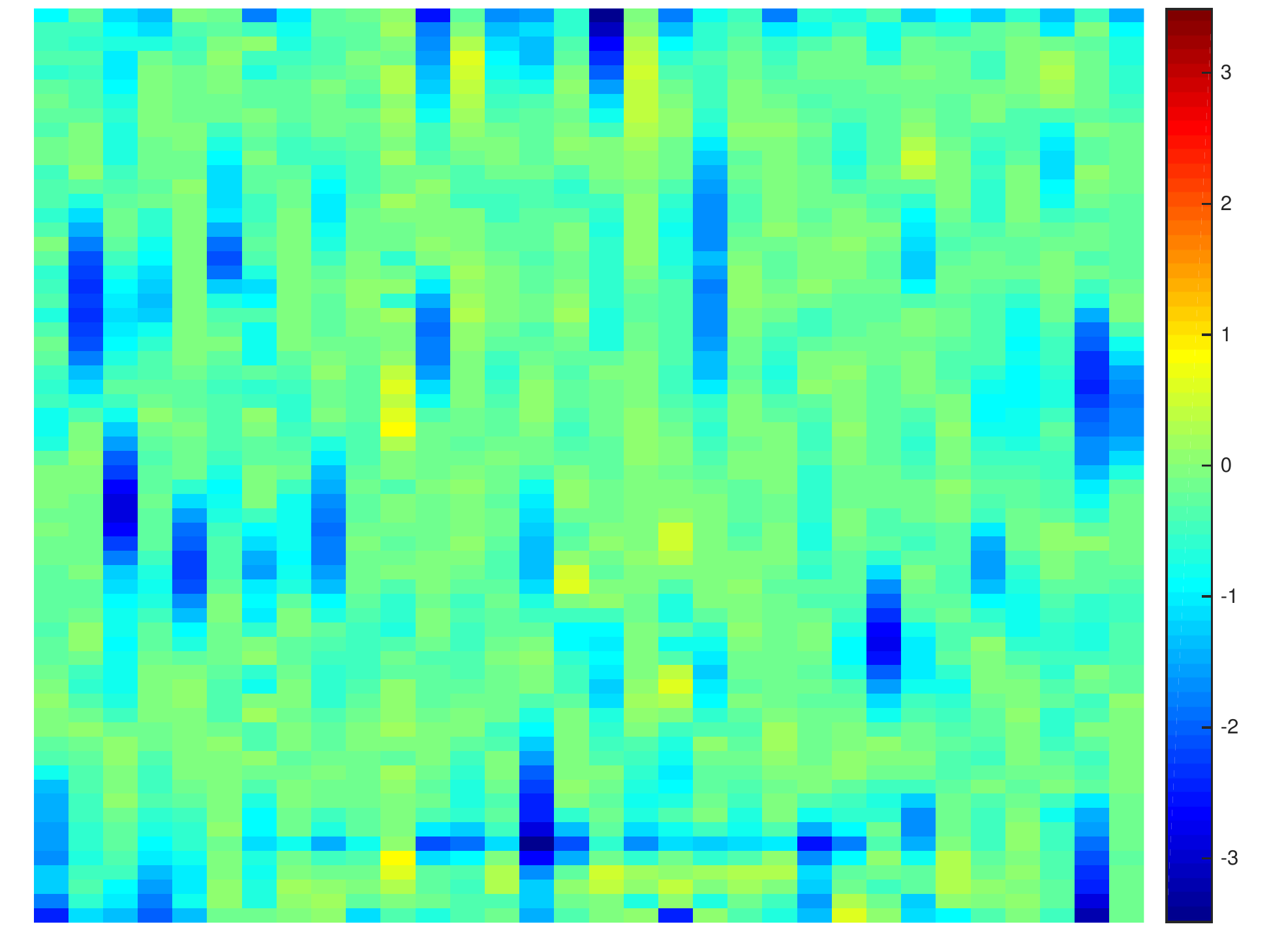}&
    \includegraphics[width = 2.5cm]{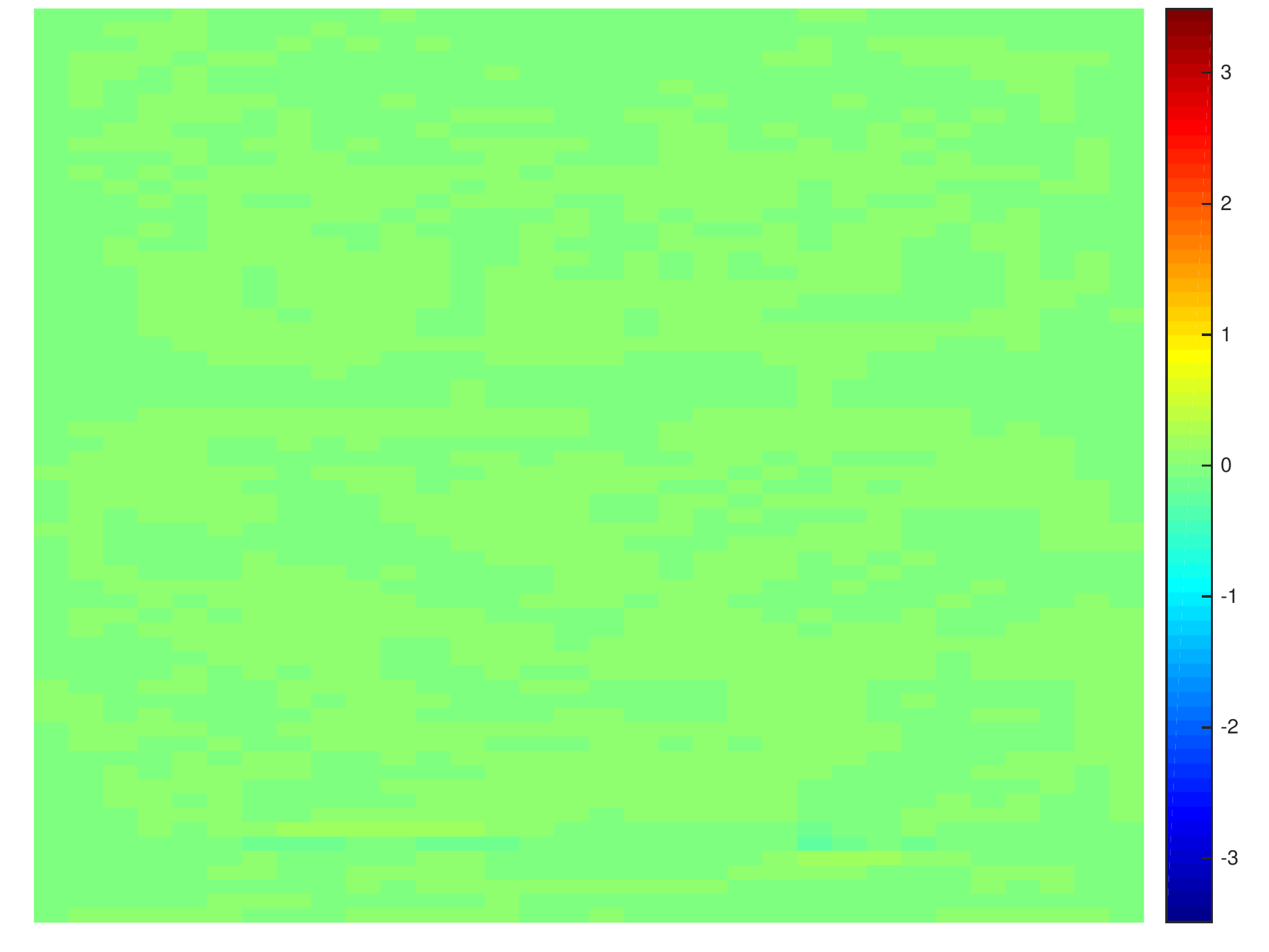}\\
  \end{tabular}
  \caption{Interpolation of the 3D lattice data set from $5\%$ random sampling. The figures in the first column are the original and subsampled angular flux at $x = 0.24$ and $x = 1.18$ . The figures in the other three columns are the results and errors of the competing algorithms.}
  \label{fig:result_random_lattice_3d_5p}
\end{figure}

\begin{figure}[H]
  \centering
  \begin{tabular}{cccc}
    Original& EBI (35.54dB)& PLE (37.20dB)  & LDMM (\textbf{39.54dB})\\
    \includegraphics[width = 2.5cm]{shock_3d_original_band_19}&
    \includegraphics[width = 2.5cm]{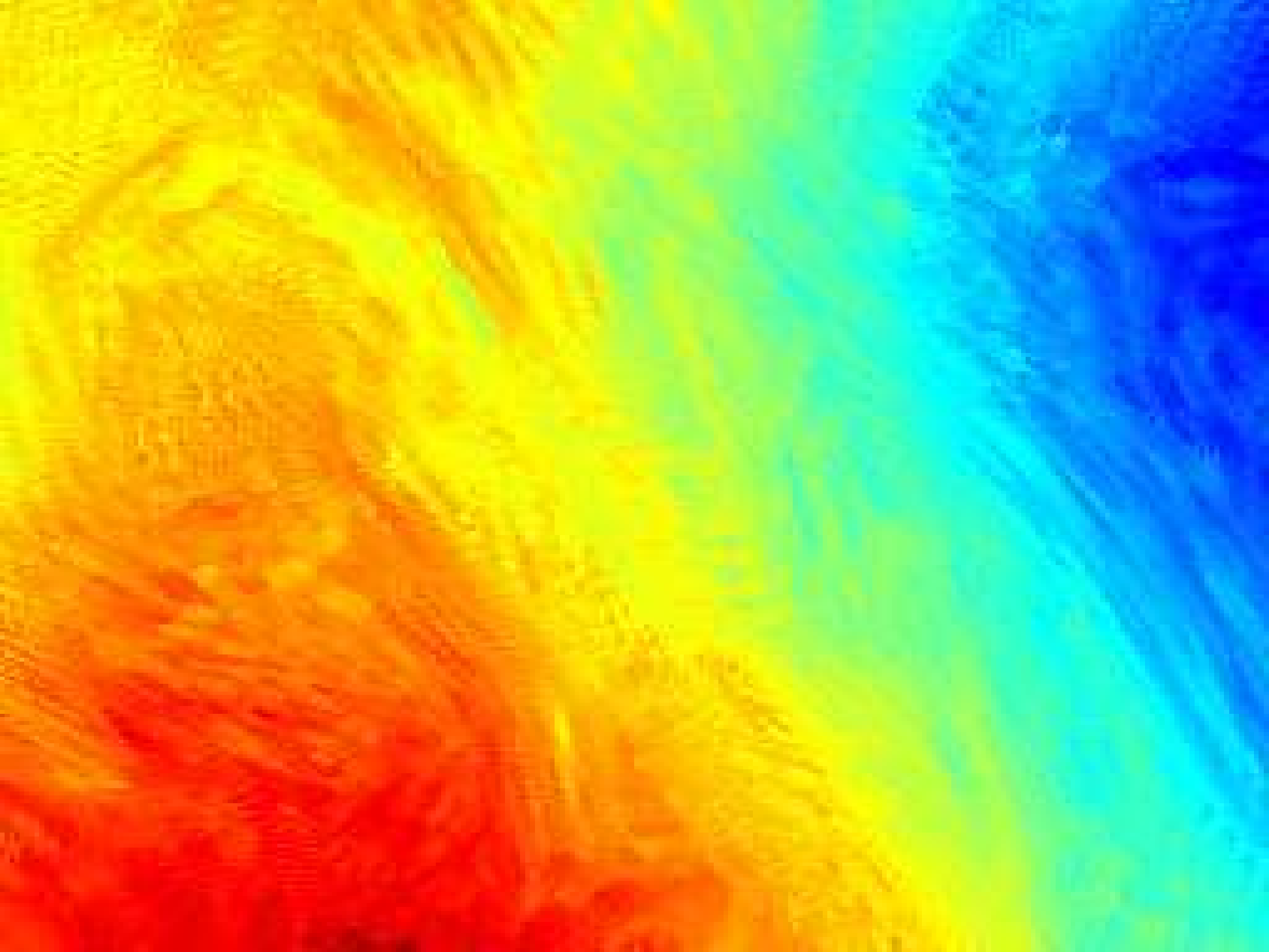}&
    \includegraphics[width = 2.5cm]{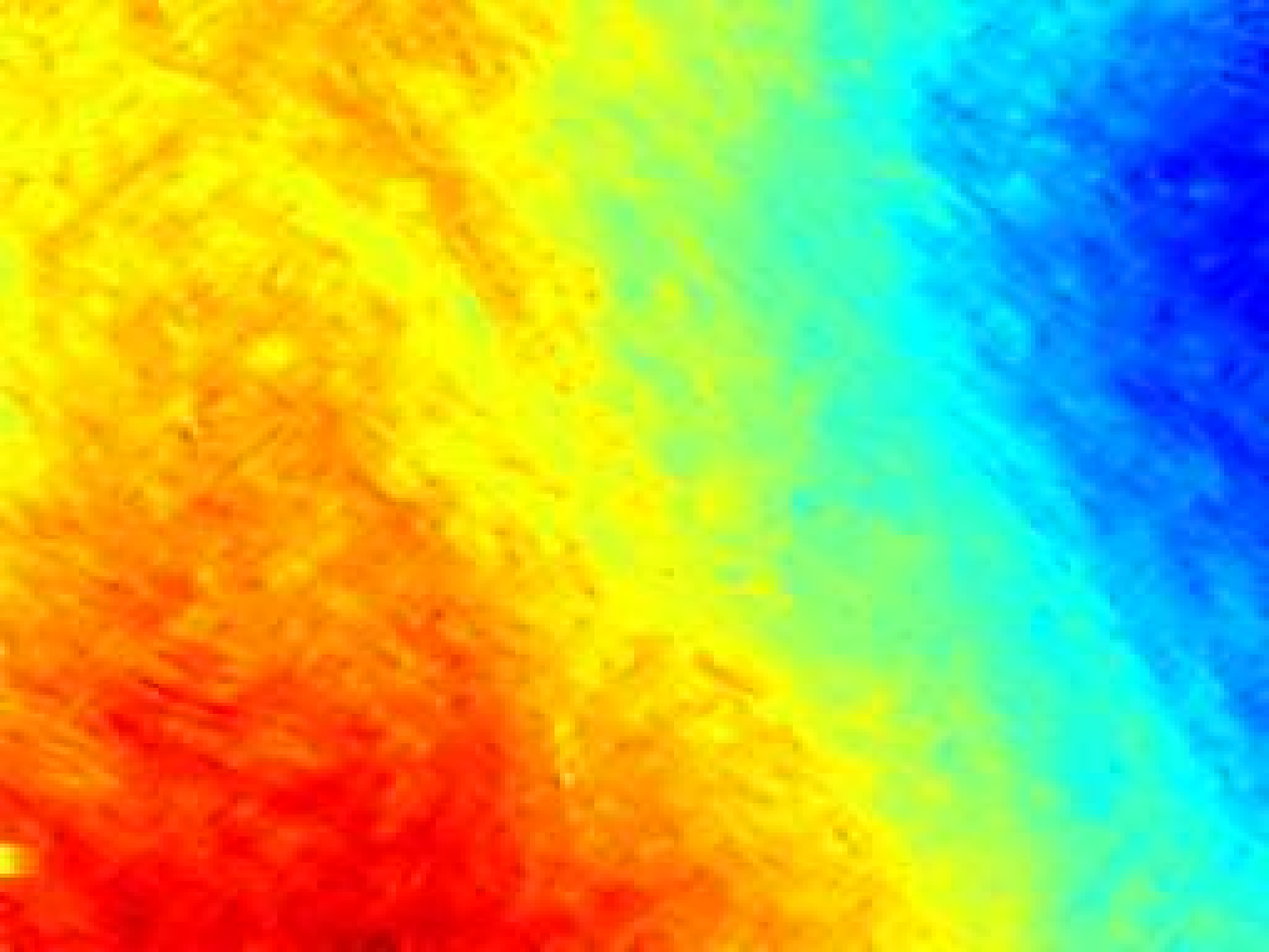}&
    \includegraphics[width = 2.5cm]{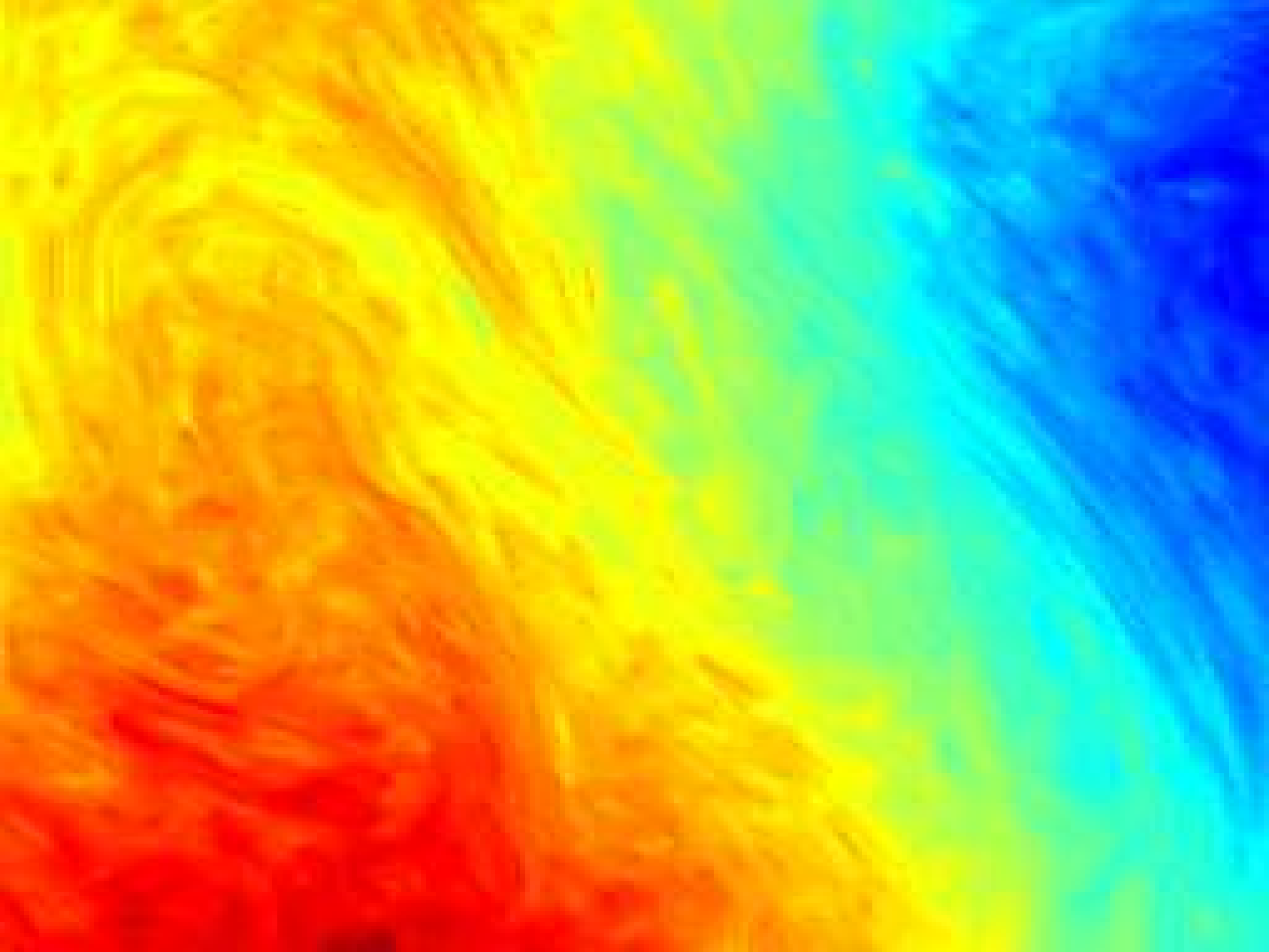}\\
    Subsample & Error & Error & Error\\
    \includegraphics[width = 2.5cm]{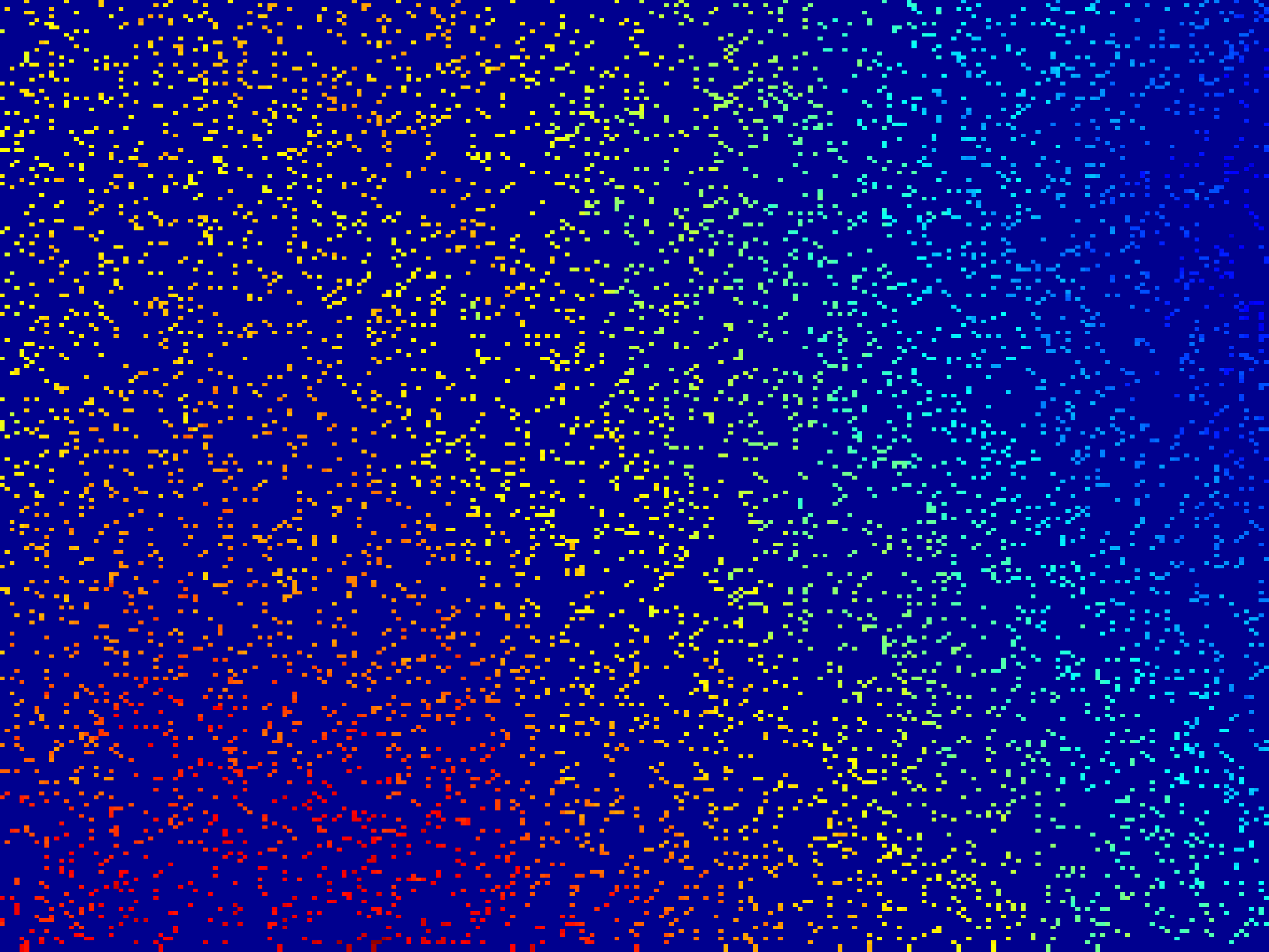}&    
    \includegraphics[width = 2.5cm]{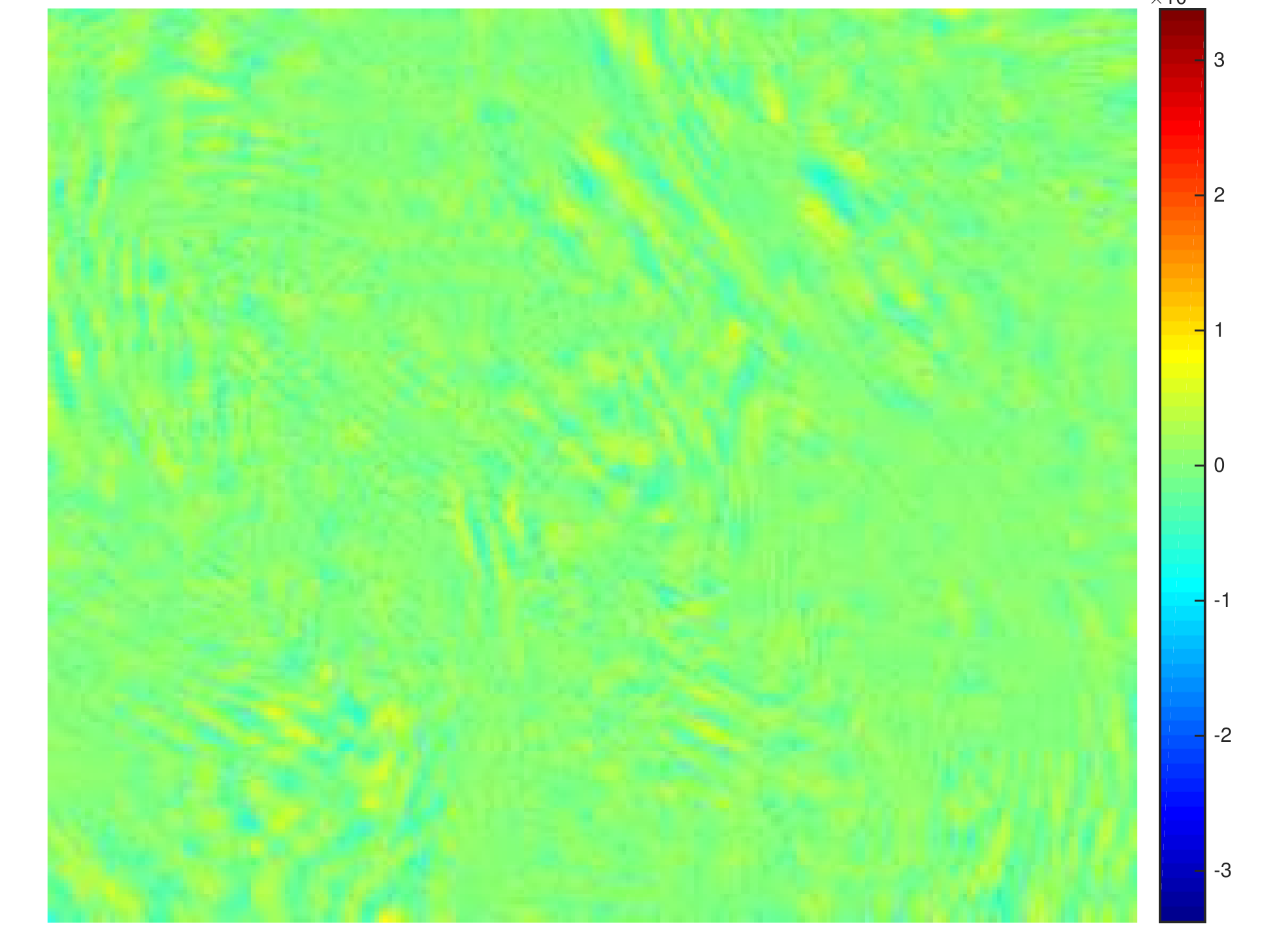}&    
    \includegraphics[width = 2.5cm]{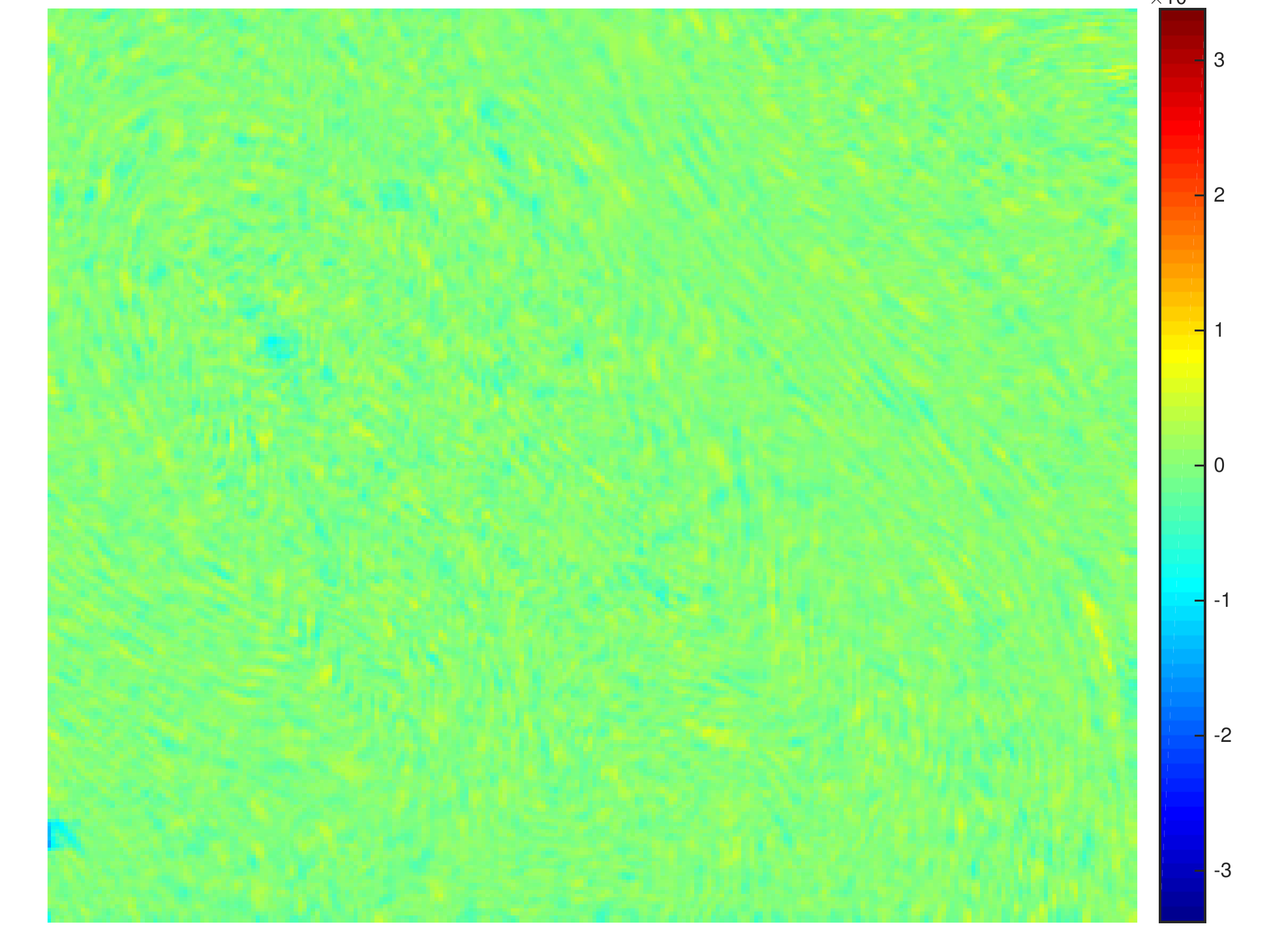}&
    \includegraphics[width = 2.5cm]{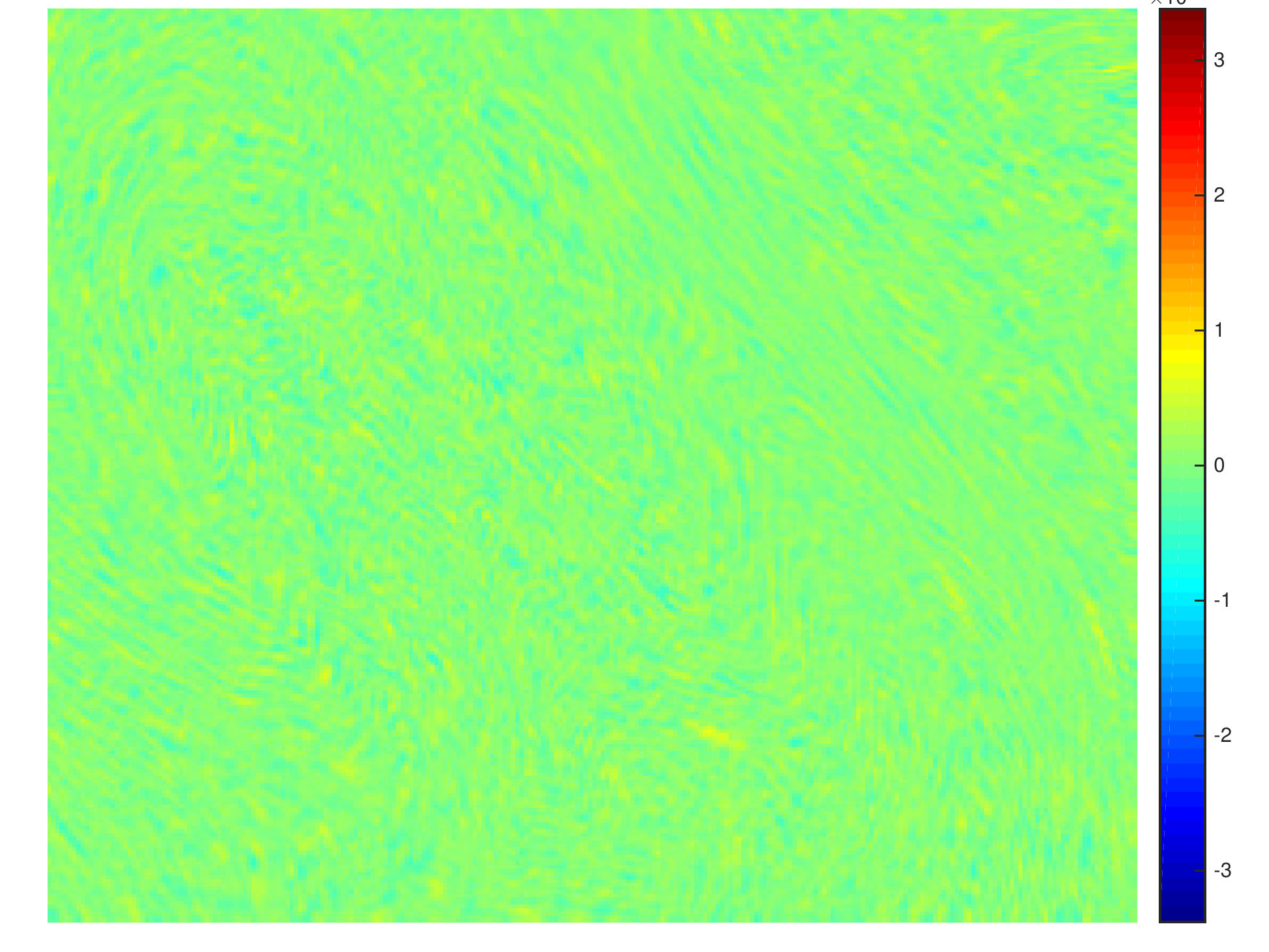}\\
    Original& EBI (35.54dB)& PLE (37.20dB)  & LDMM (\textbf{39.54dB})\\
    \includegraphics[width = 2.5cm]{shock_3d_original_band_29}&
    \includegraphics[width = 2.5cm]{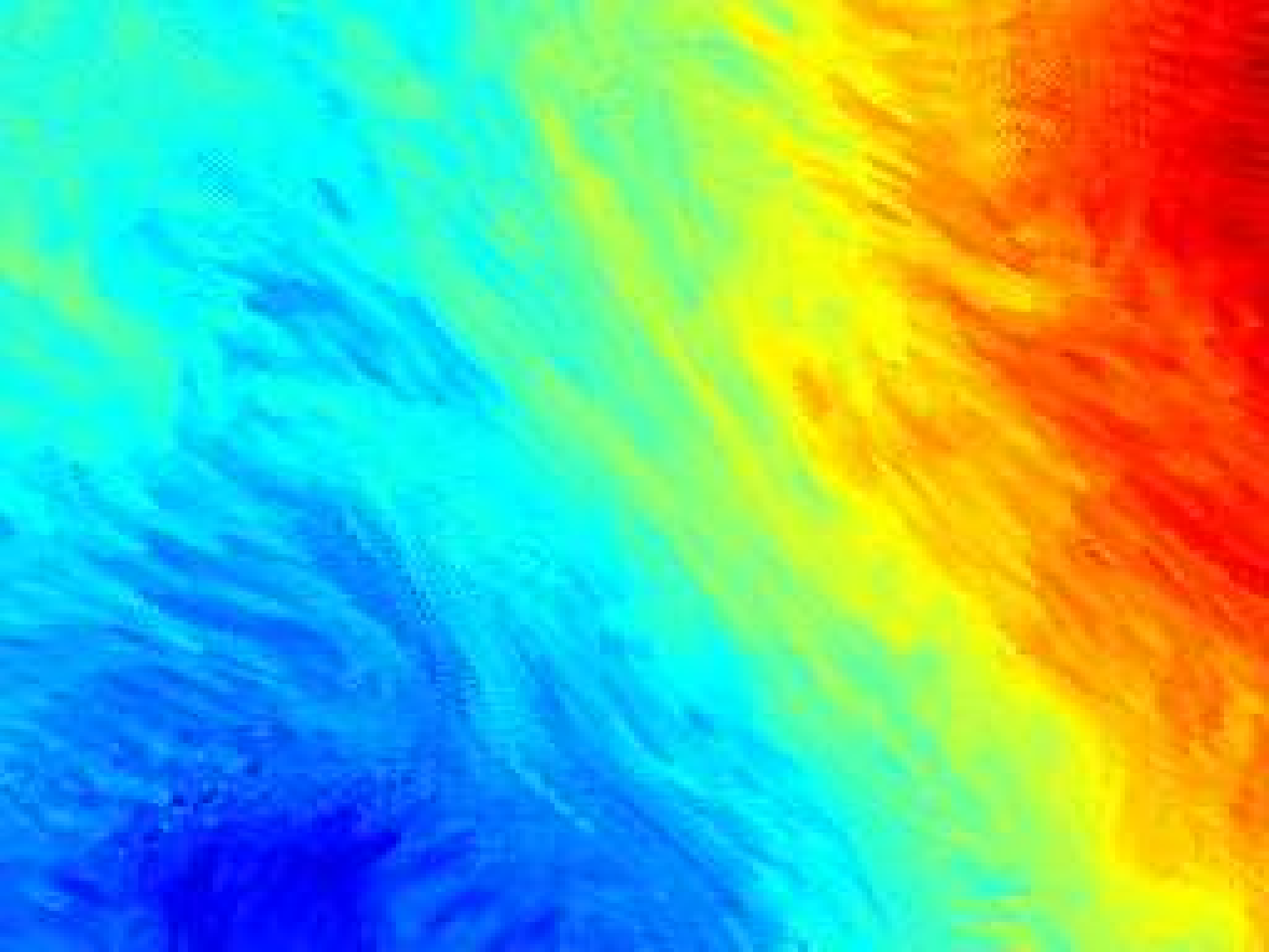}&
    \includegraphics[width = 2.5cm]{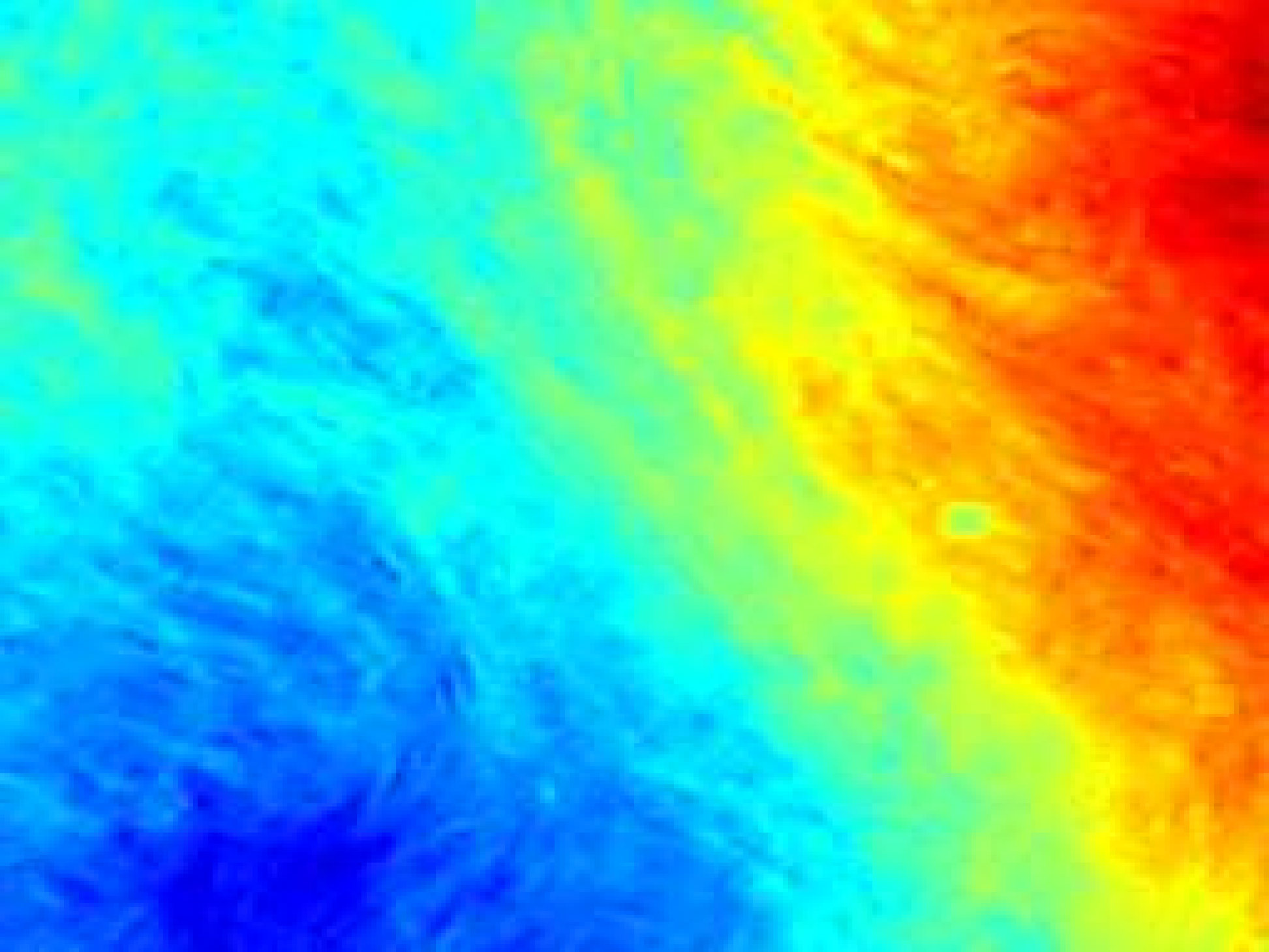}&
    \includegraphics[width = 2.5cm]{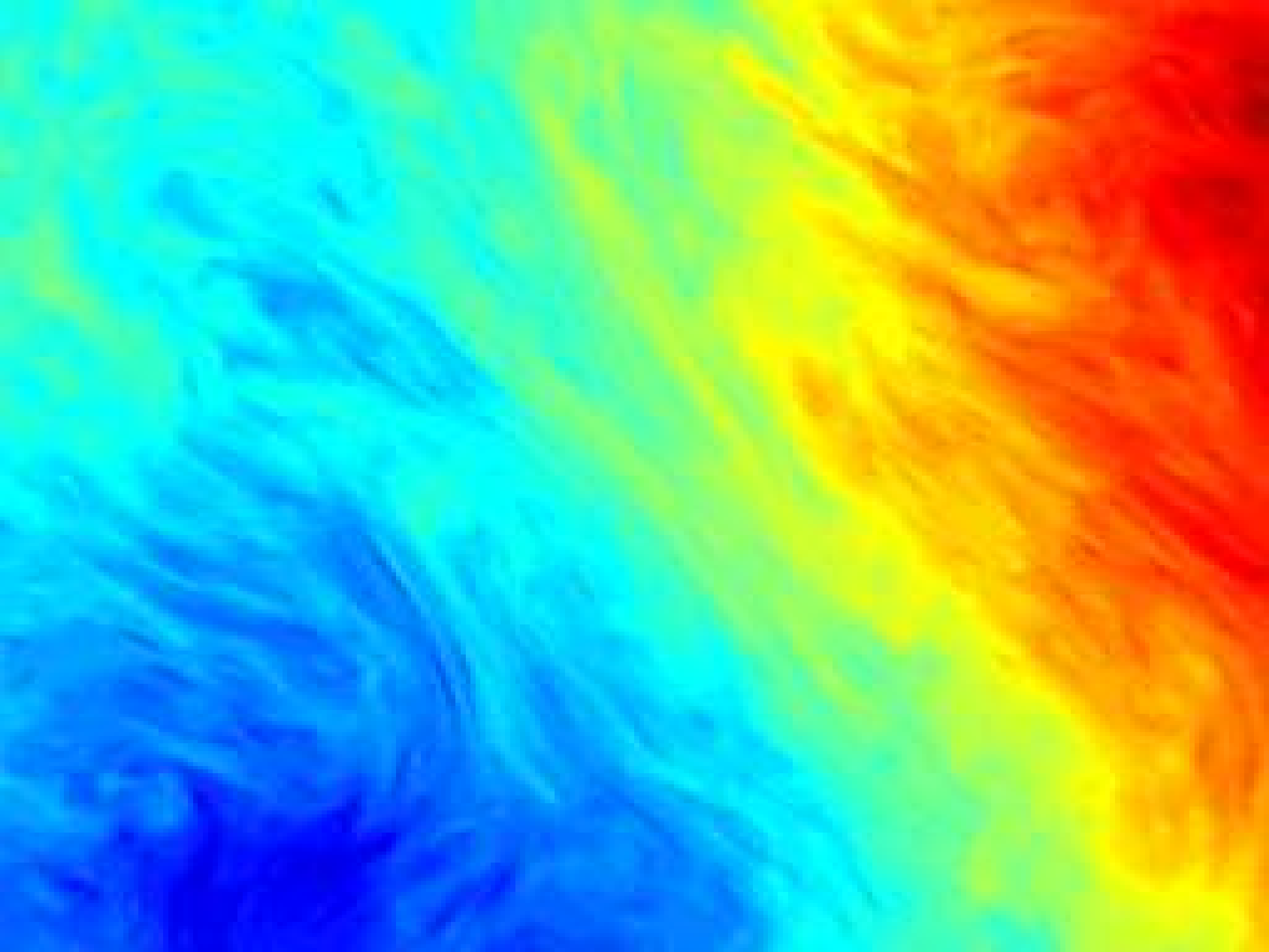}\\
    Subsample & Error & Error & Error\\
    \includegraphics[width = 2.5cm]{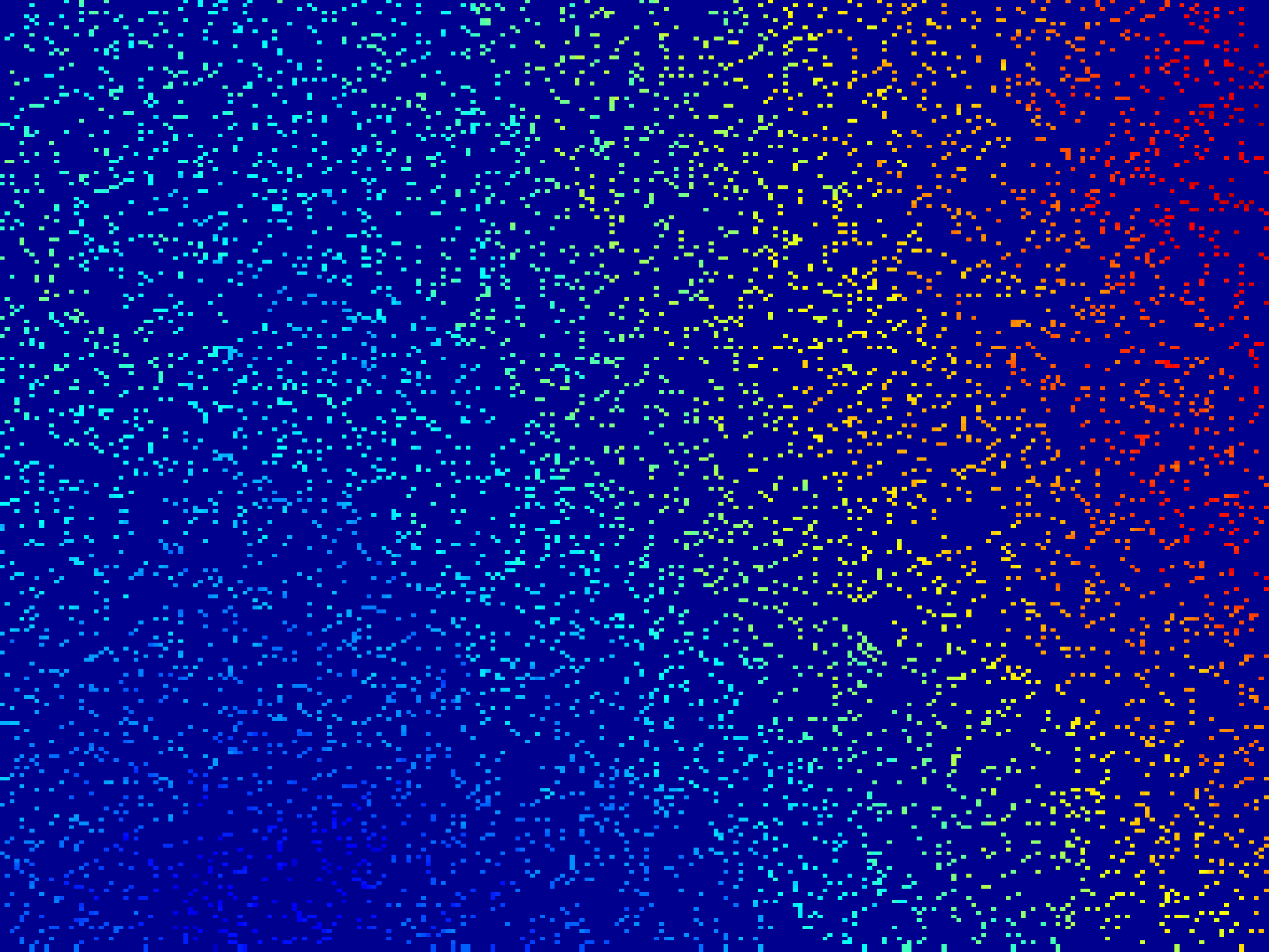}&    
    \includegraphics[width = 2.5cm]{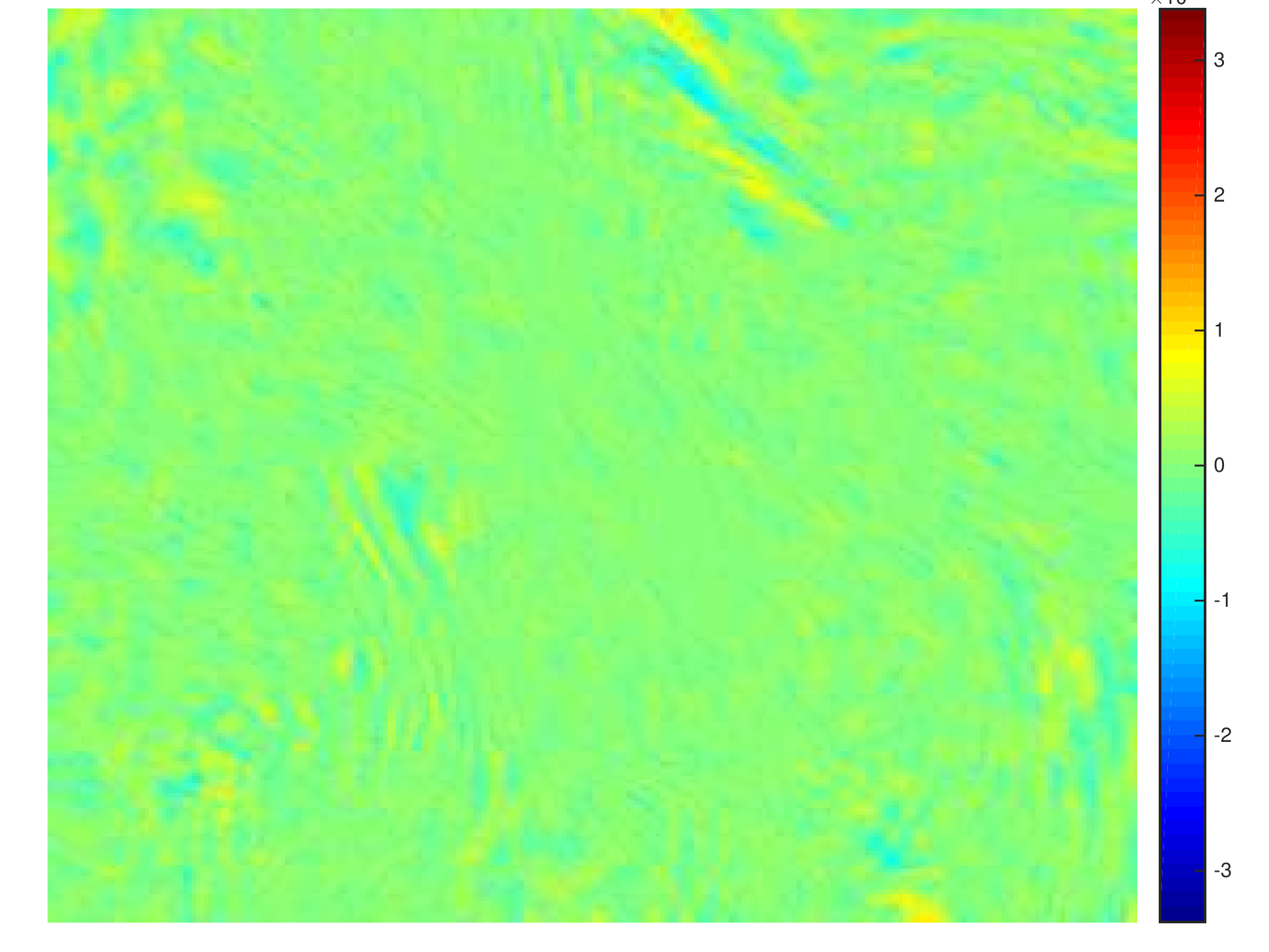}&    
    \includegraphics[width = 2.5cm]{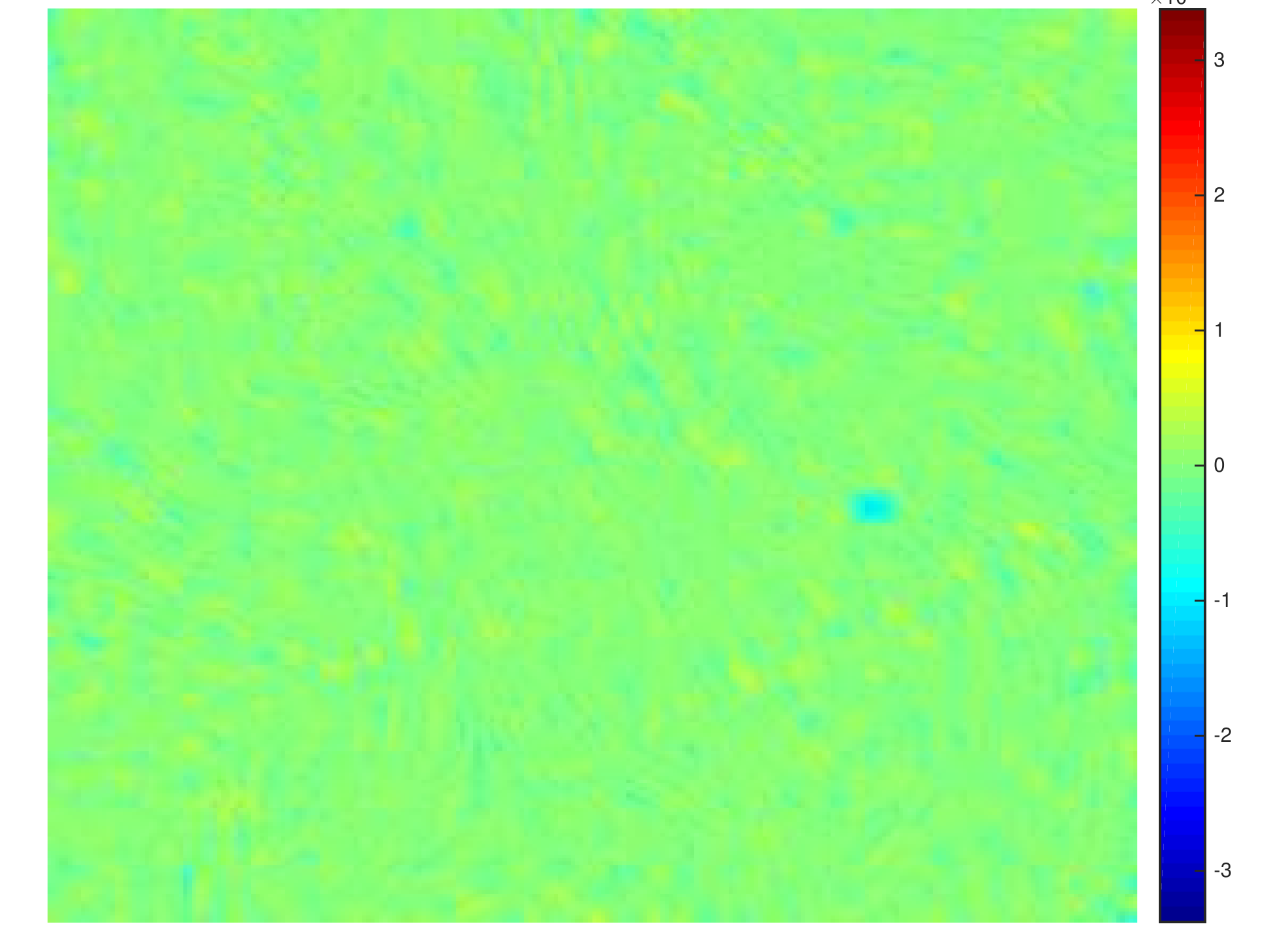}&
    \includegraphics[width = 2.5cm]{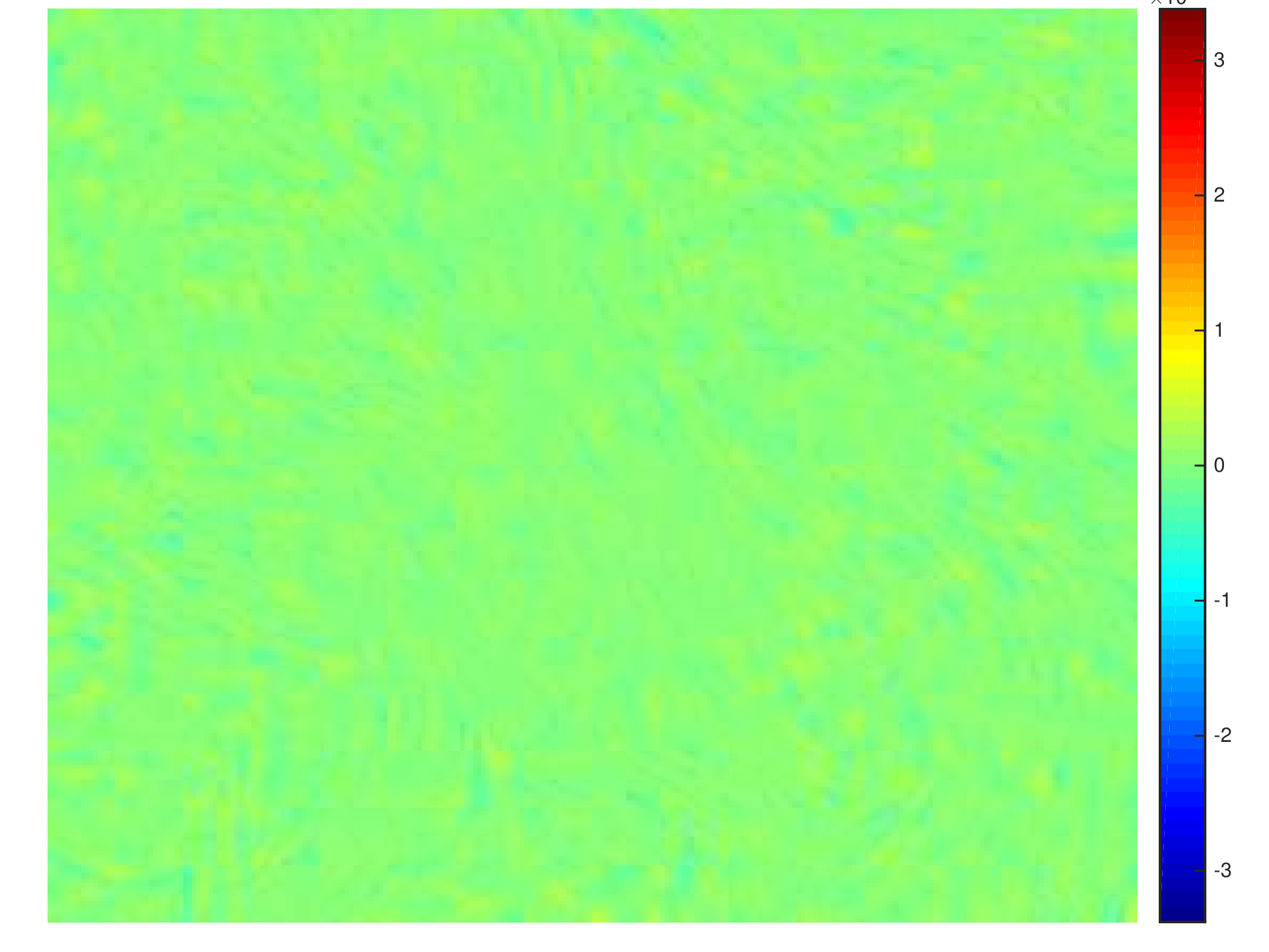}
  \end{tabular}
  \caption{Interpolation of the 3D plasma (distribution function) data set from $10\%$ random sampling. The figures in the first column are two spatial cross sections of the original and subsampled data. The figures in the other three columns are the results and errors of the competing algorithms.}
  \label{fig:result_random_shock_3d_10p}
\end{figure}

\begin{figure}[H]
  \centering
  \begin{tabular}{cccc}
    Original& EBI (34.87dB)& PLE (20.96dB)  & LDMM (\textbf{37.72dB})\\
    \includegraphics[width = 2.5cm]{shock_3d_original_band_19}&
    \includegraphics[width = 2.5cm]{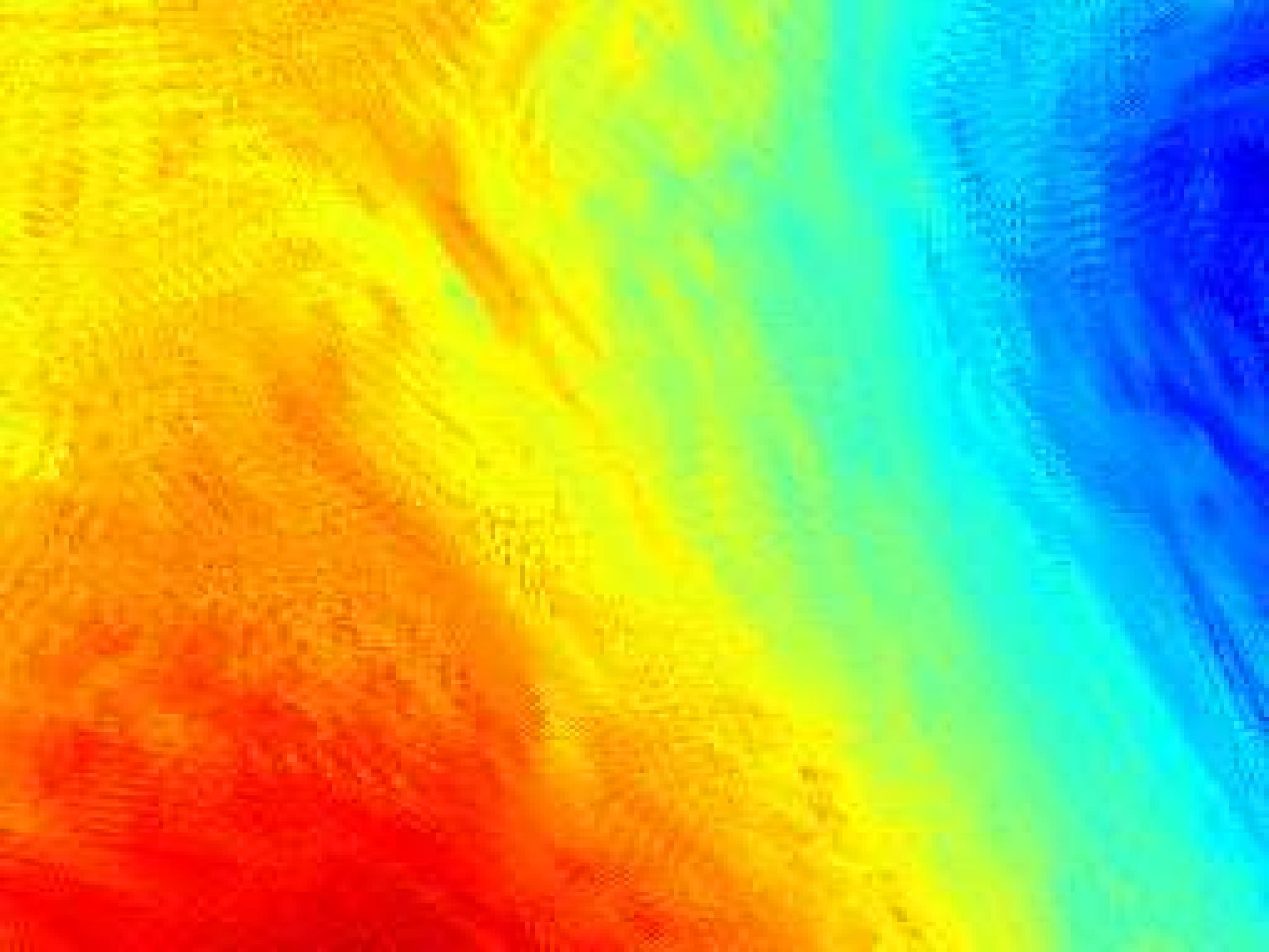}&
    \includegraphics[width = 2.5cm]{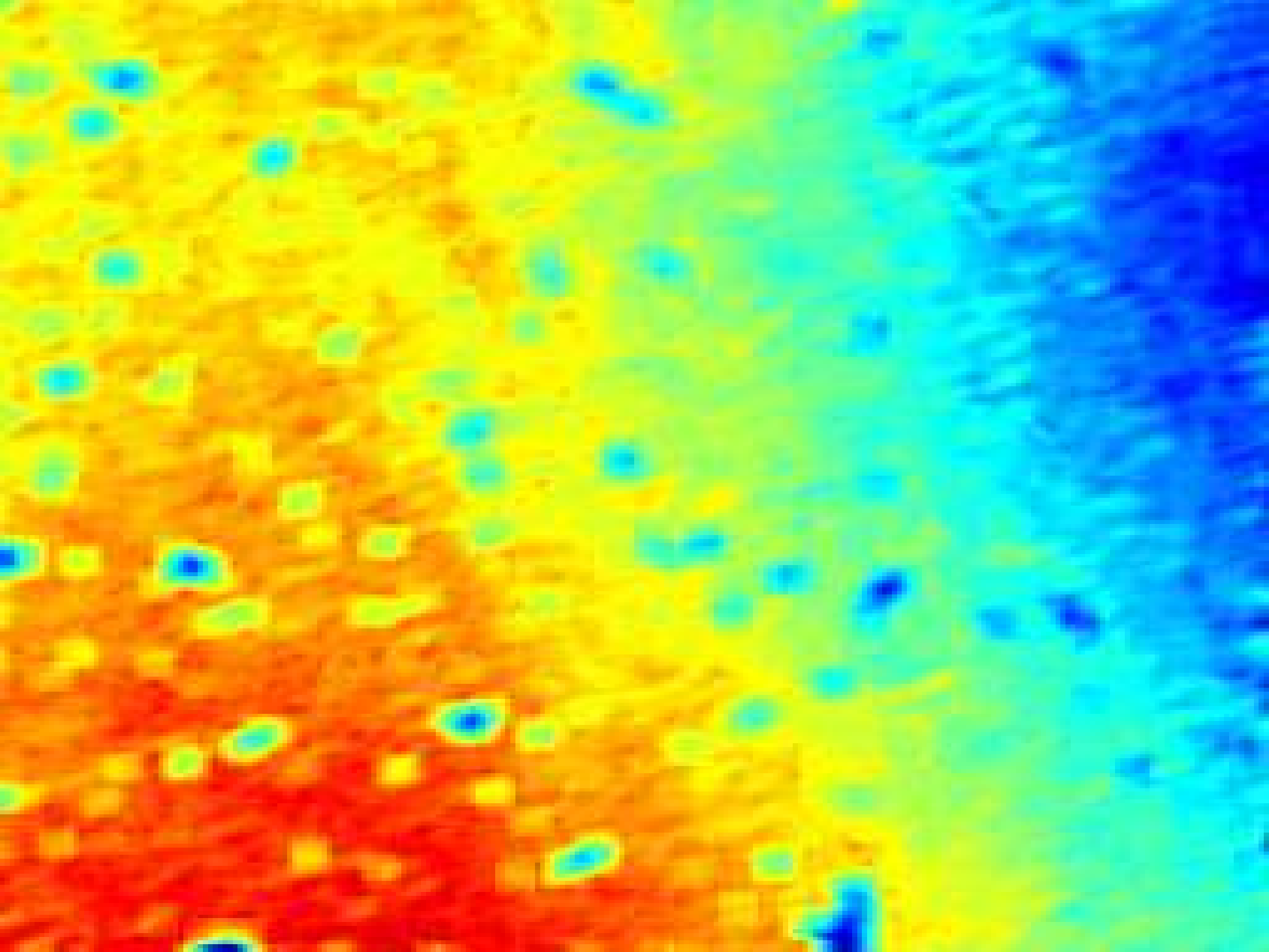}&
    \includegraphics[width = 2.5cm]{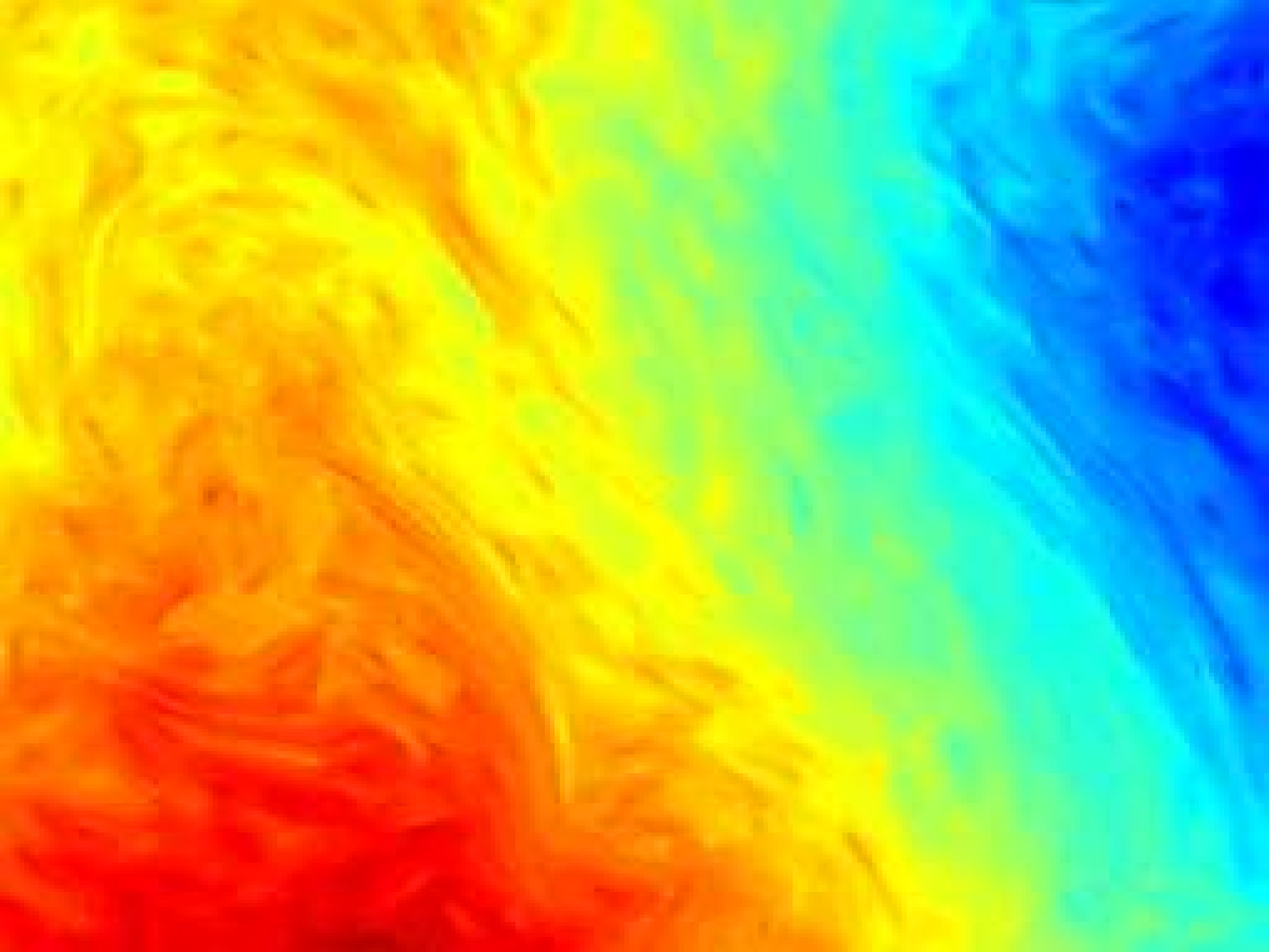}\\
    Subsample & Error & Error & Error\\
    \includegraphics[width = 2.5cm]{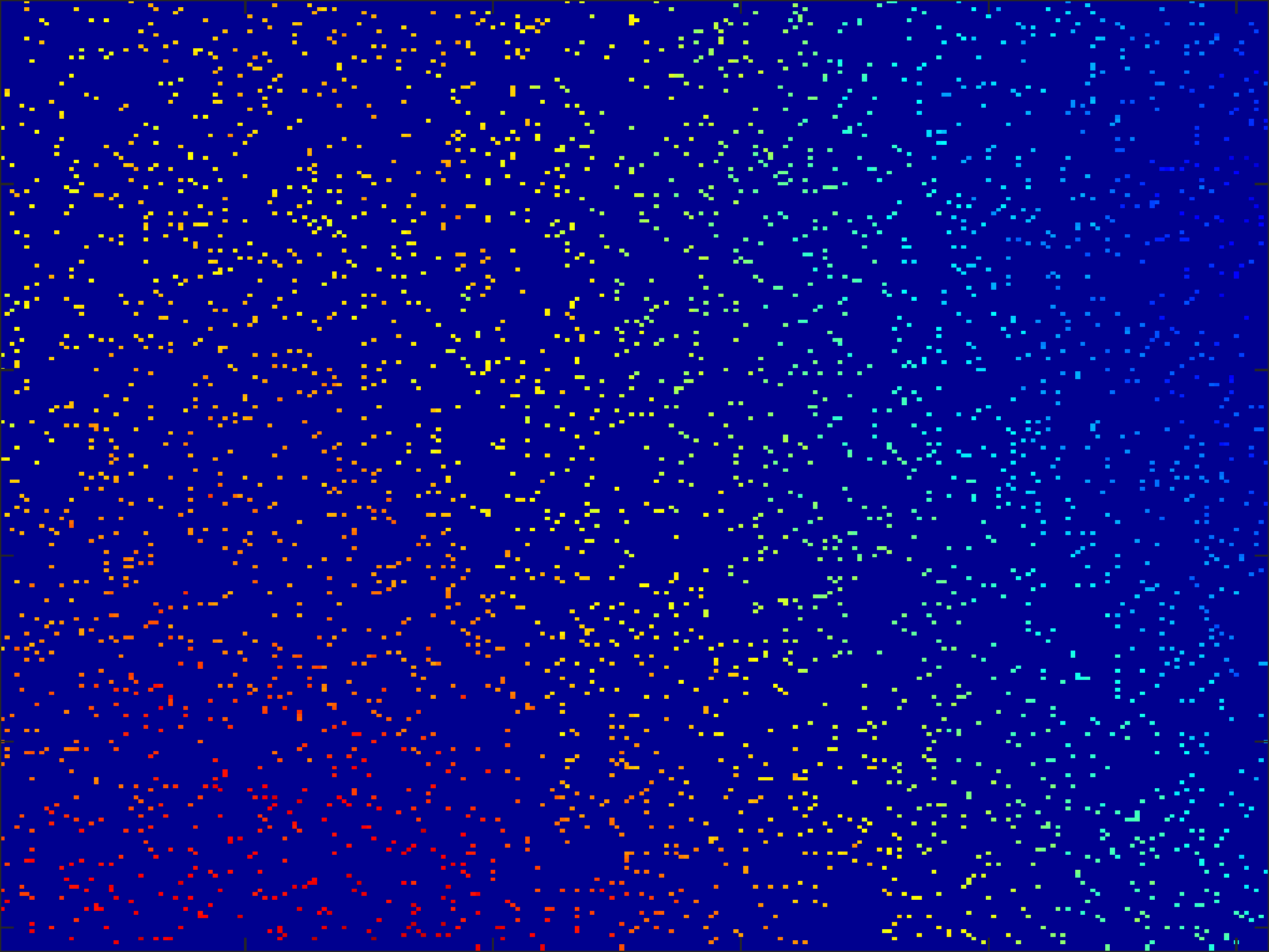}&    
    \includegraphics[width = 2.5cm]{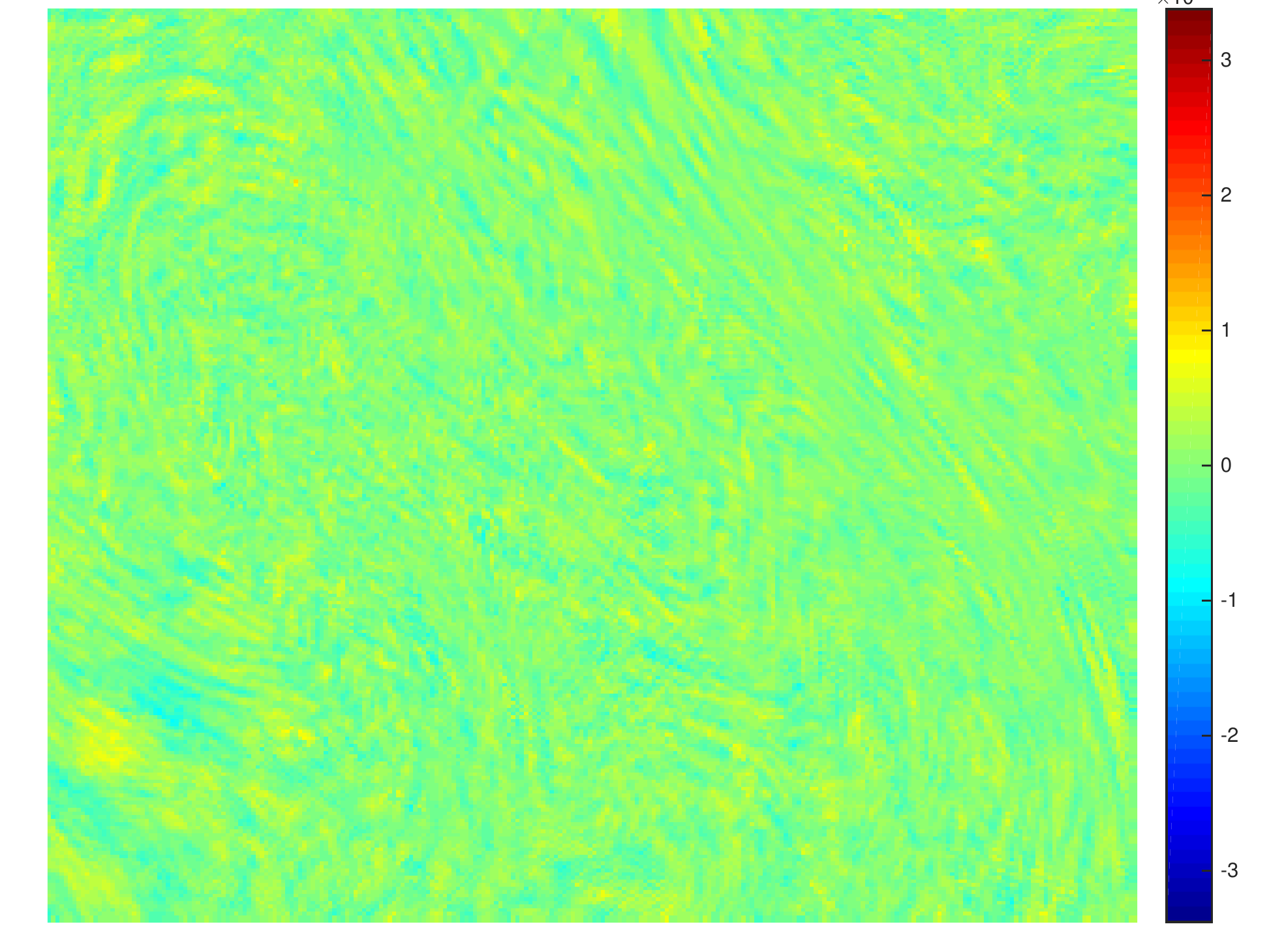}&    
    \includegraphics[width = 2.5cm]{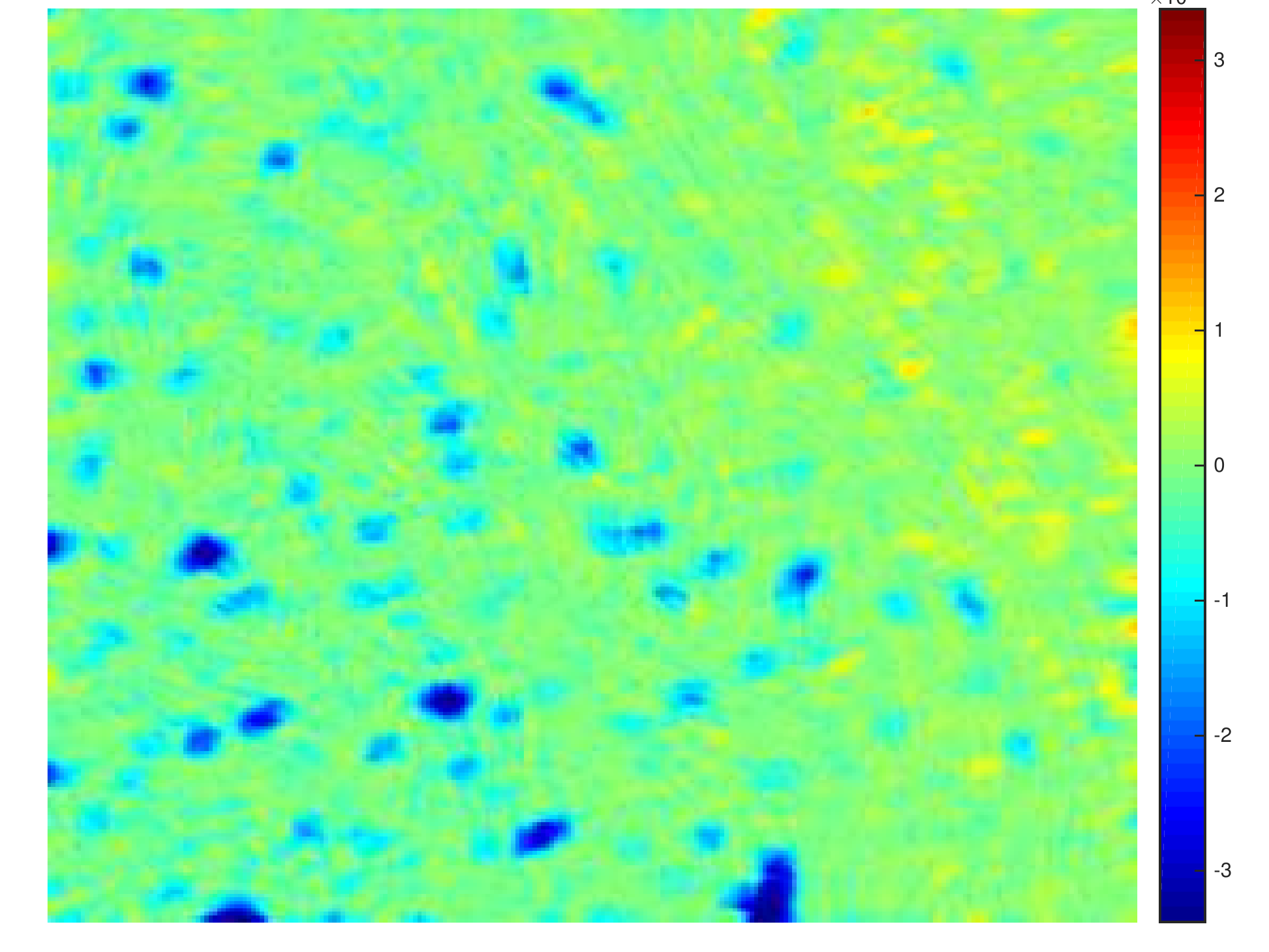}&
    \includegraphics[width = 2.5cm]{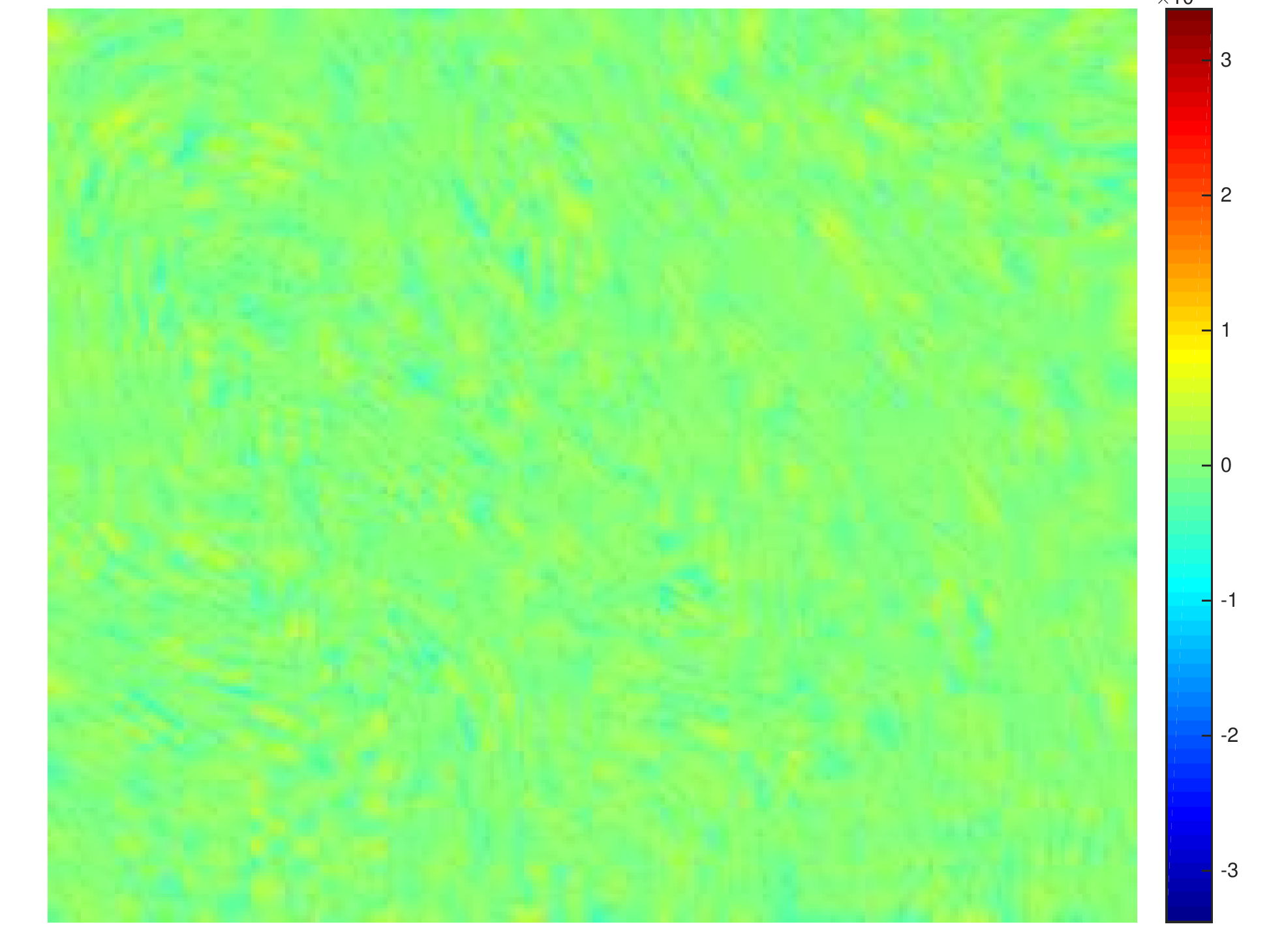}\\
    Original& EBI (34.87dB)& PLE (20.96dB)  & LDMM (\textbf{37.72dB})\\
    \includegraphics[width = 2.5cm]{shock_3d_original_band_29}&
    \includegraphics[width = 2.5cm]{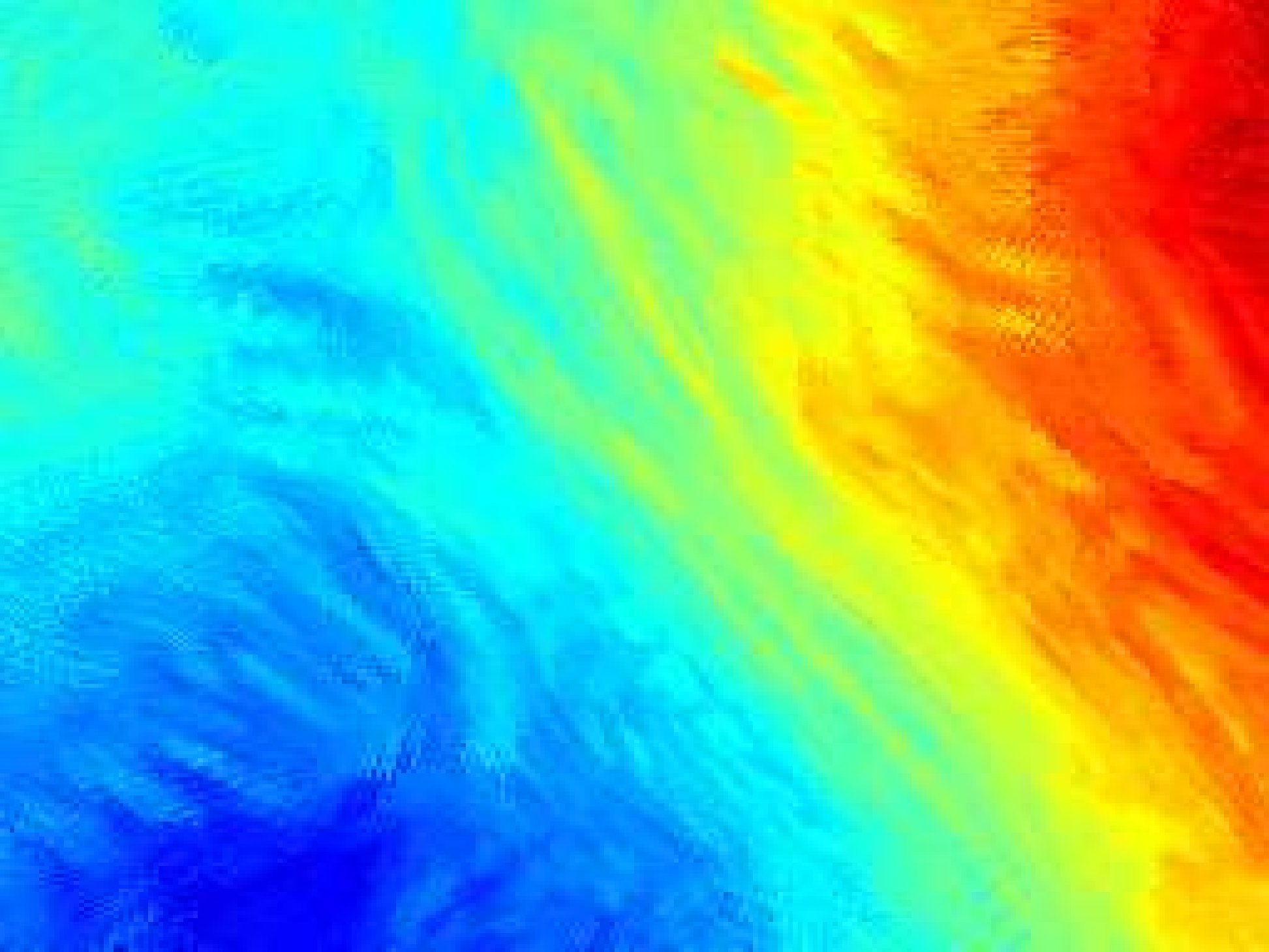}&
    \includegraphics[width = 2.5cm]{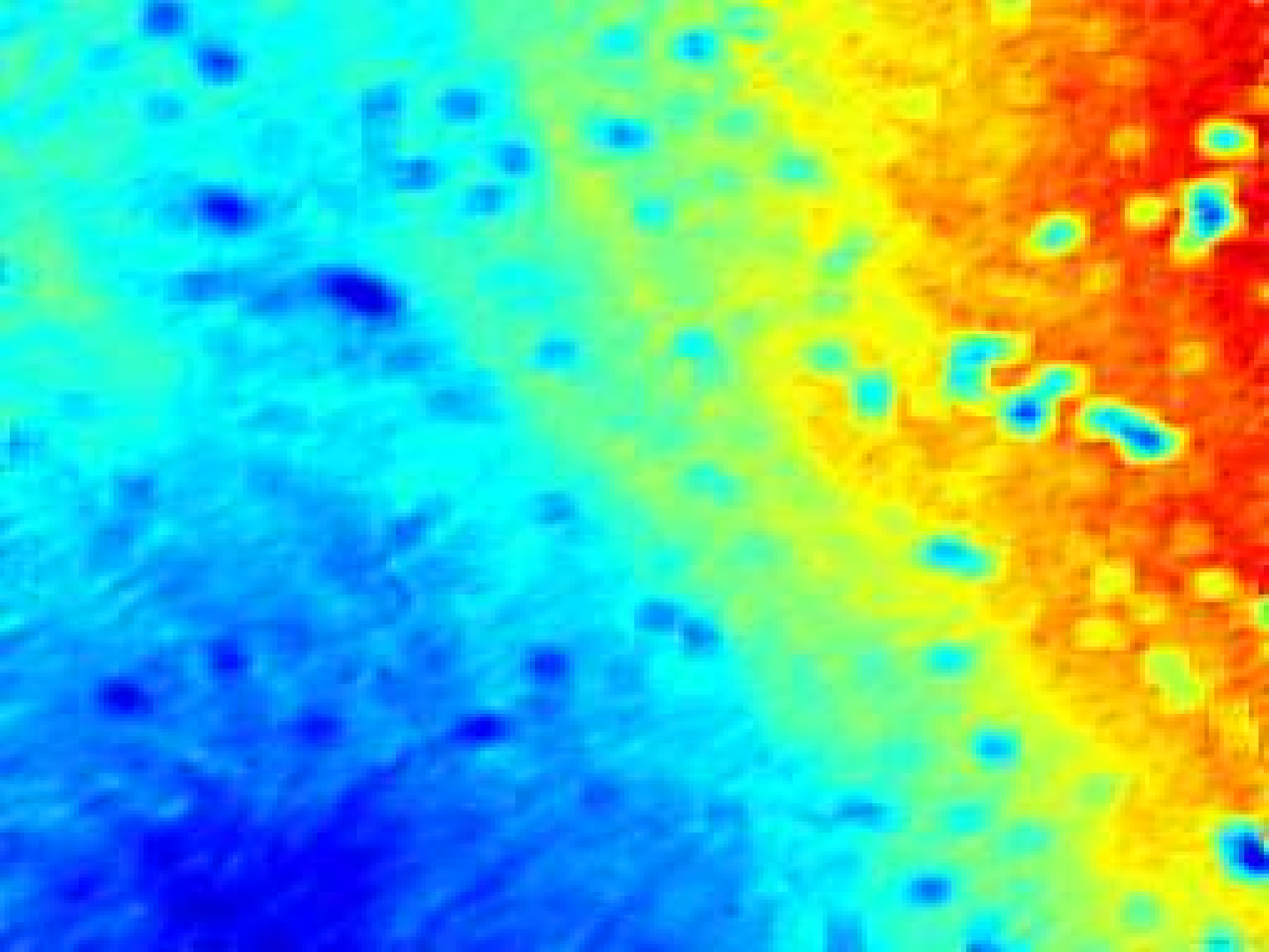}&
    \includegraphics[width = 2.5cm]{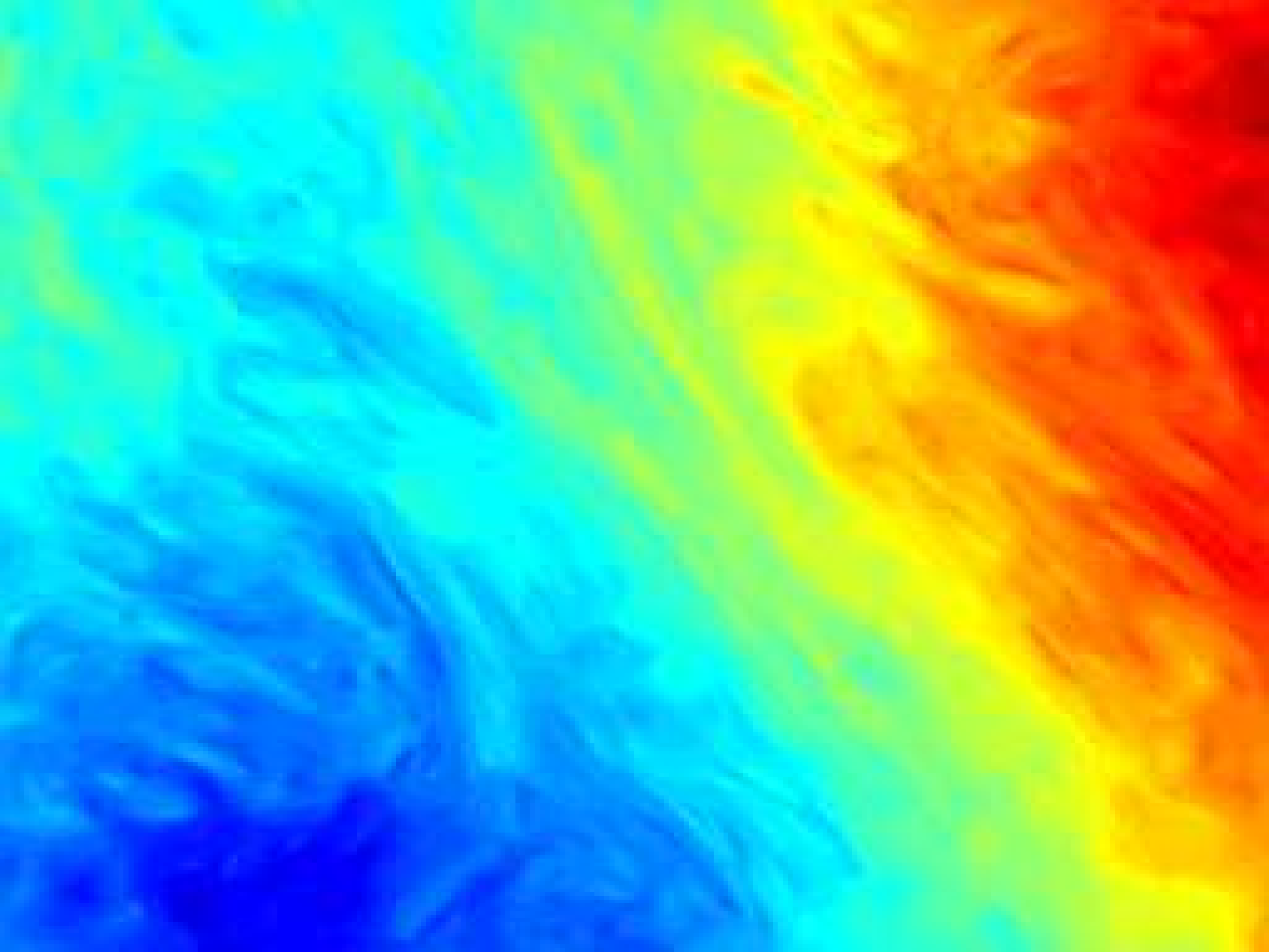}\\
    Subsample & Error & Error & Error\\
    \includegraphics[width = 2.5cm]{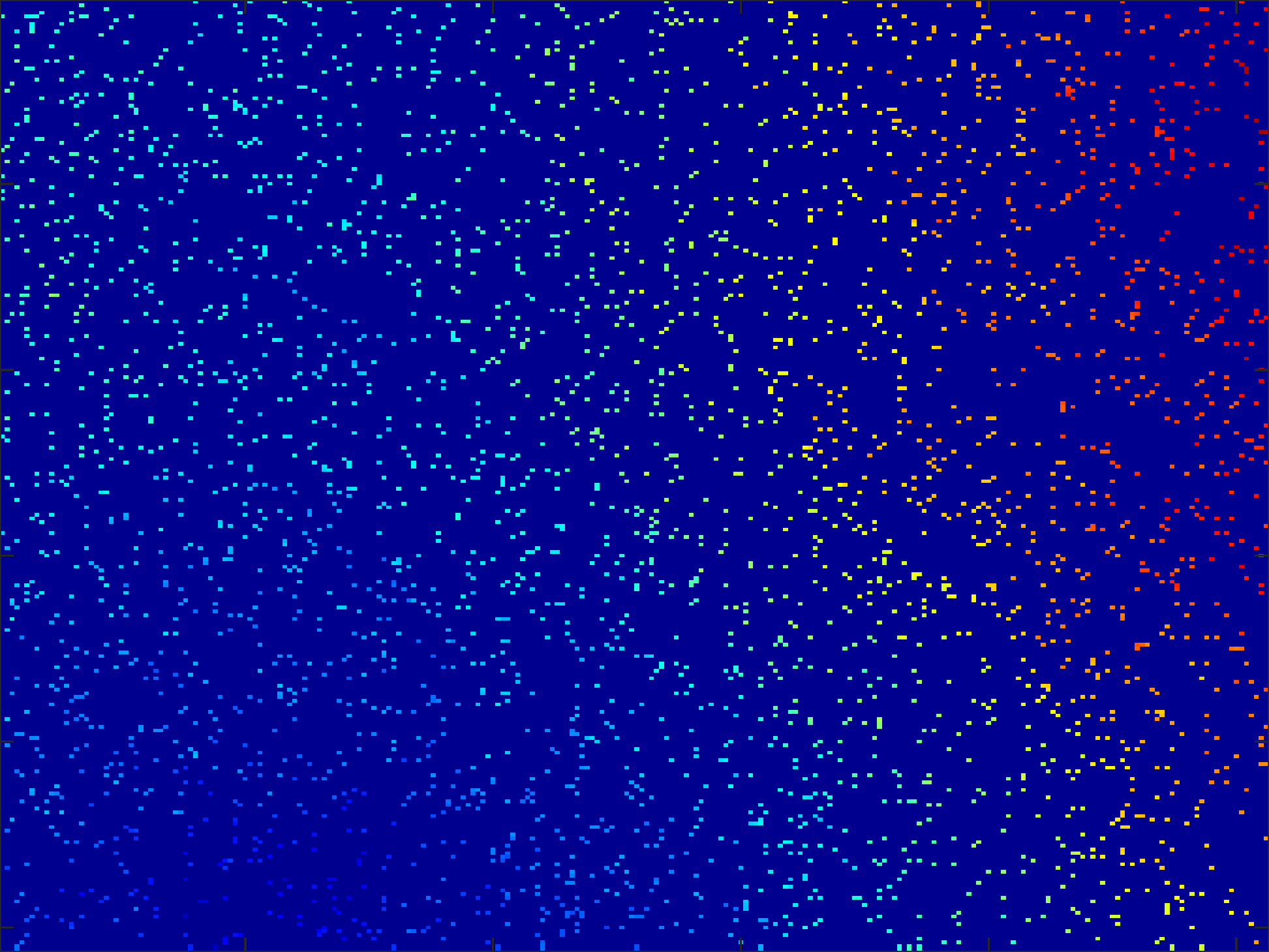}&    
    \includegraphics[width = 2.5cm]{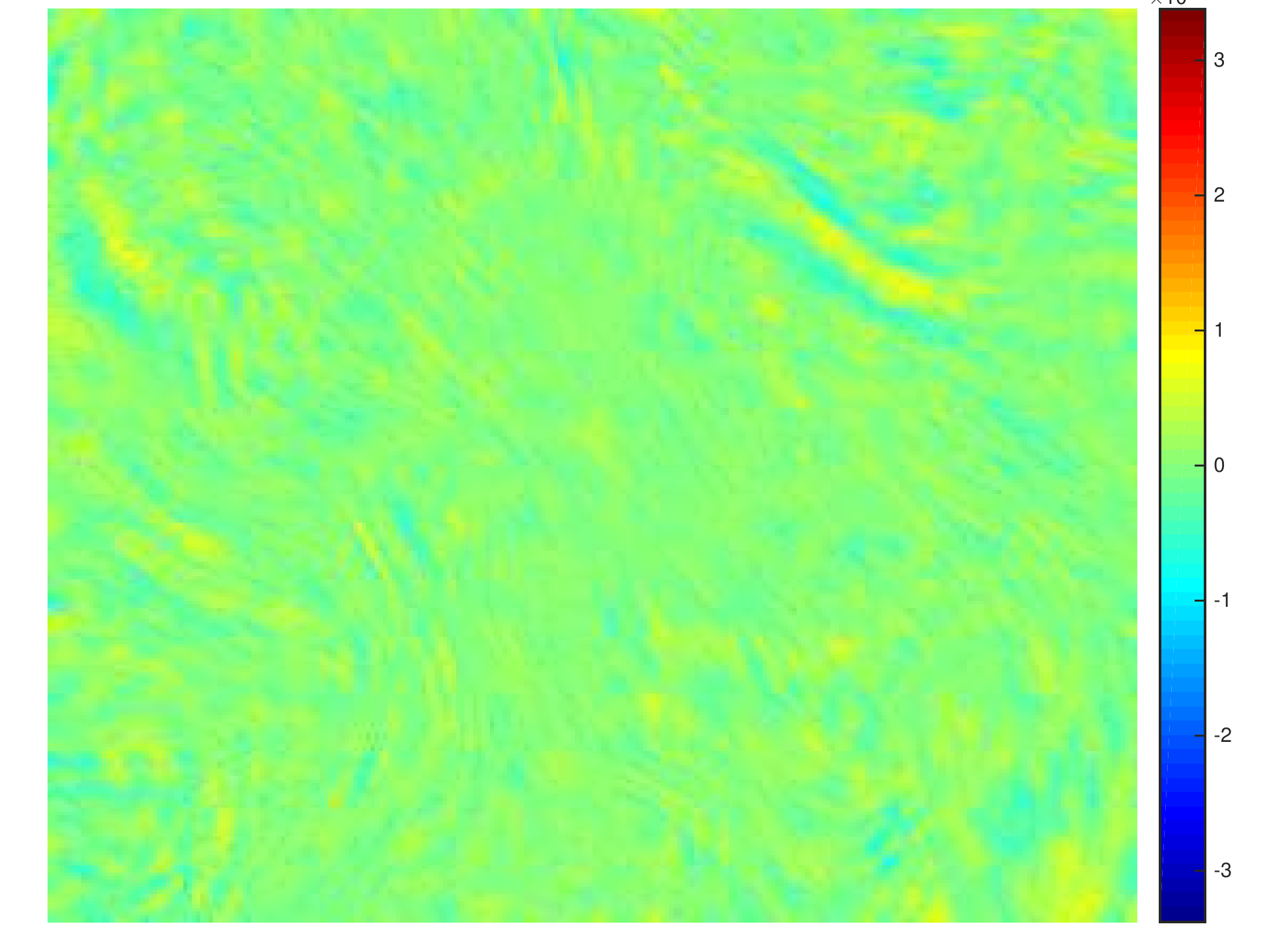}&    
    \includegraphics[width = 2.5cm]{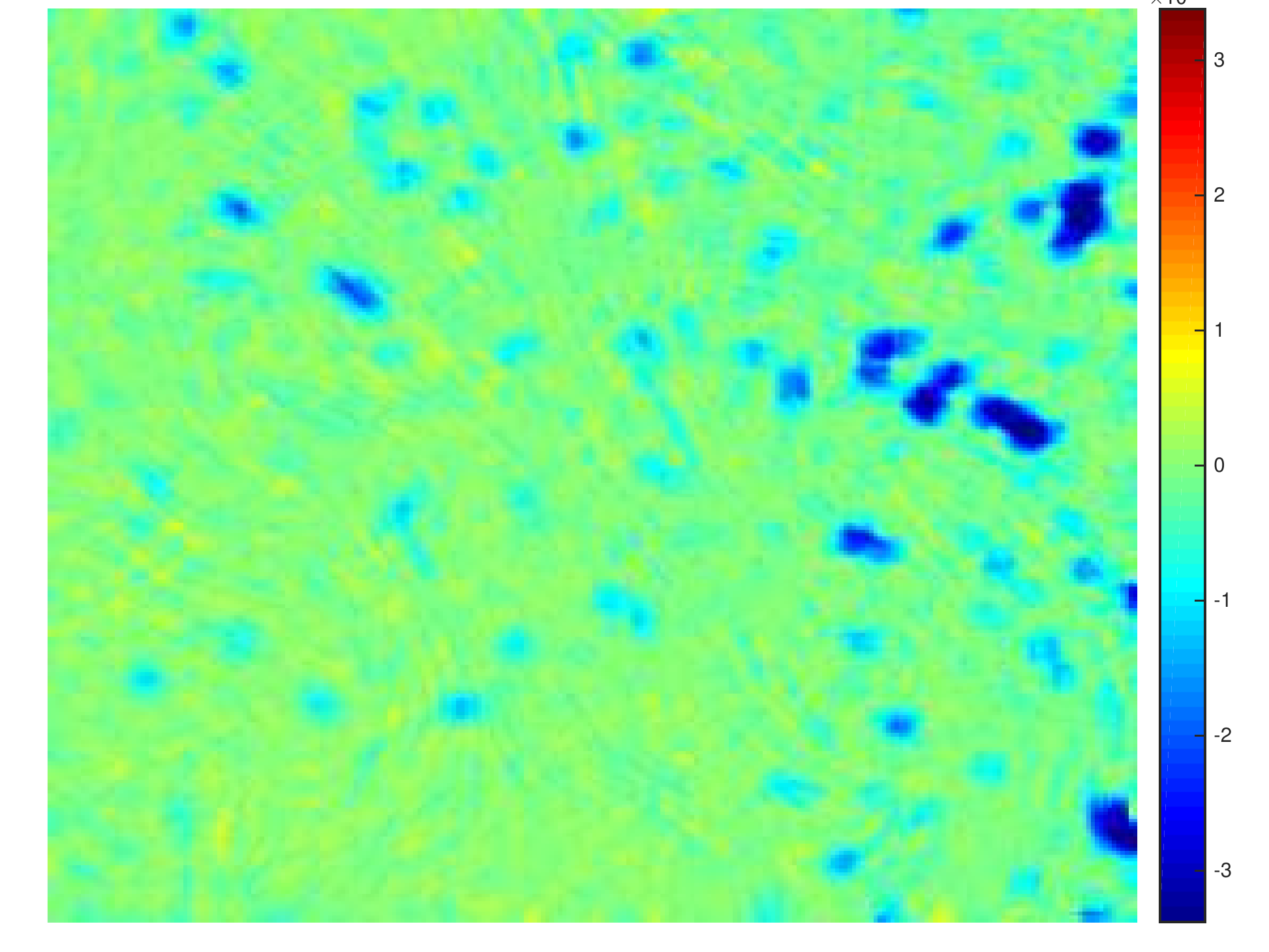}&
    \includegraphics[width = 2.5cm]{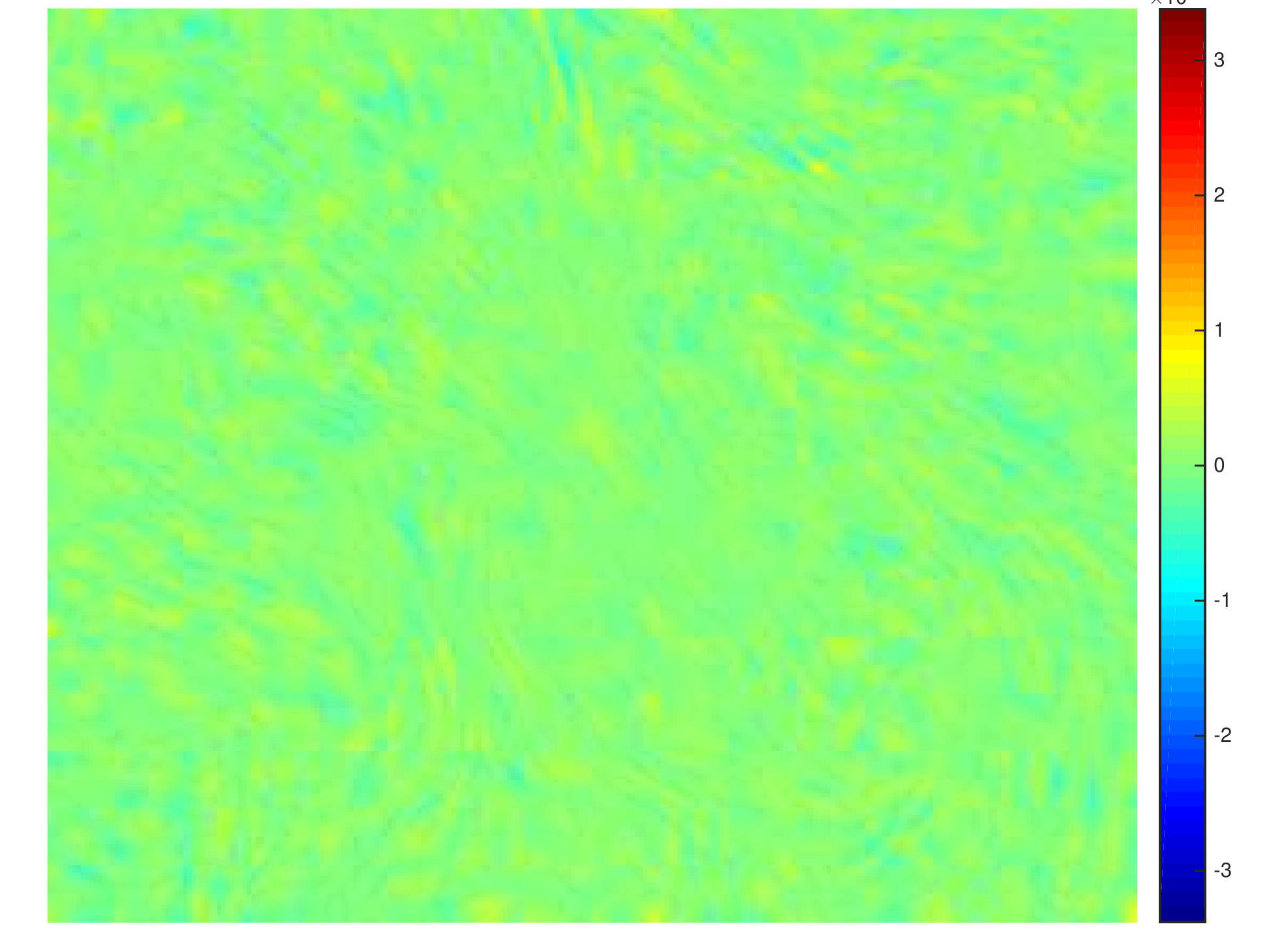}
  \end{tabular}
  \caption{Interpolation of the 3D plasma (distribution function) data set from $5\%$ random sampling. The figures in the first column are two spatial cross sections of the original and subsampled data. The figures in the other three columns are the results and errors of the competing algorithms.}
  \label{fig:result_random_shock_3d_5p}
\end{figure}

The visual of the interpolation with $10\%$ and $5\%$ are shown in Figure \ref{fig:result_random_2d_10p}-\ref{fig:result_random_shock_3d_5p}. The errors of the reconstruction in different norms are displayed in Table \ref{tab:error_random_antonio_2d}-\ref{tab:error_random_shock_3d}. It can be observed that LDMM consistently performs at a higher accuracy than EBI and PLE either visually or numerically. The superiority of LDMM is more dramatic when the sample rate is very low (5\%), in which case PLE fails to achieve reasonable results. LDMM also manages to yield smoother results, whereas EBI tends to create artificial patchy patterns. We point out that the reconstruction of the 3D data sets with PLE and EBI are obtained by applying the algorithms to 2D cross sections because of a lack of 3D implementations of both algorithms. Therefore it is not entirely fair to compare LDMM to PLE and EBI on the 3D data sets. This is especially clear on the 3D lattice data set, where values change smoothly on each direction. Nonetheless, the vast superiority of LDMM on 2D examples illustrates its advantage over the competing algorithms.

\begin{figure}[H]
  \centering
  \includegraphics[width=12cm]{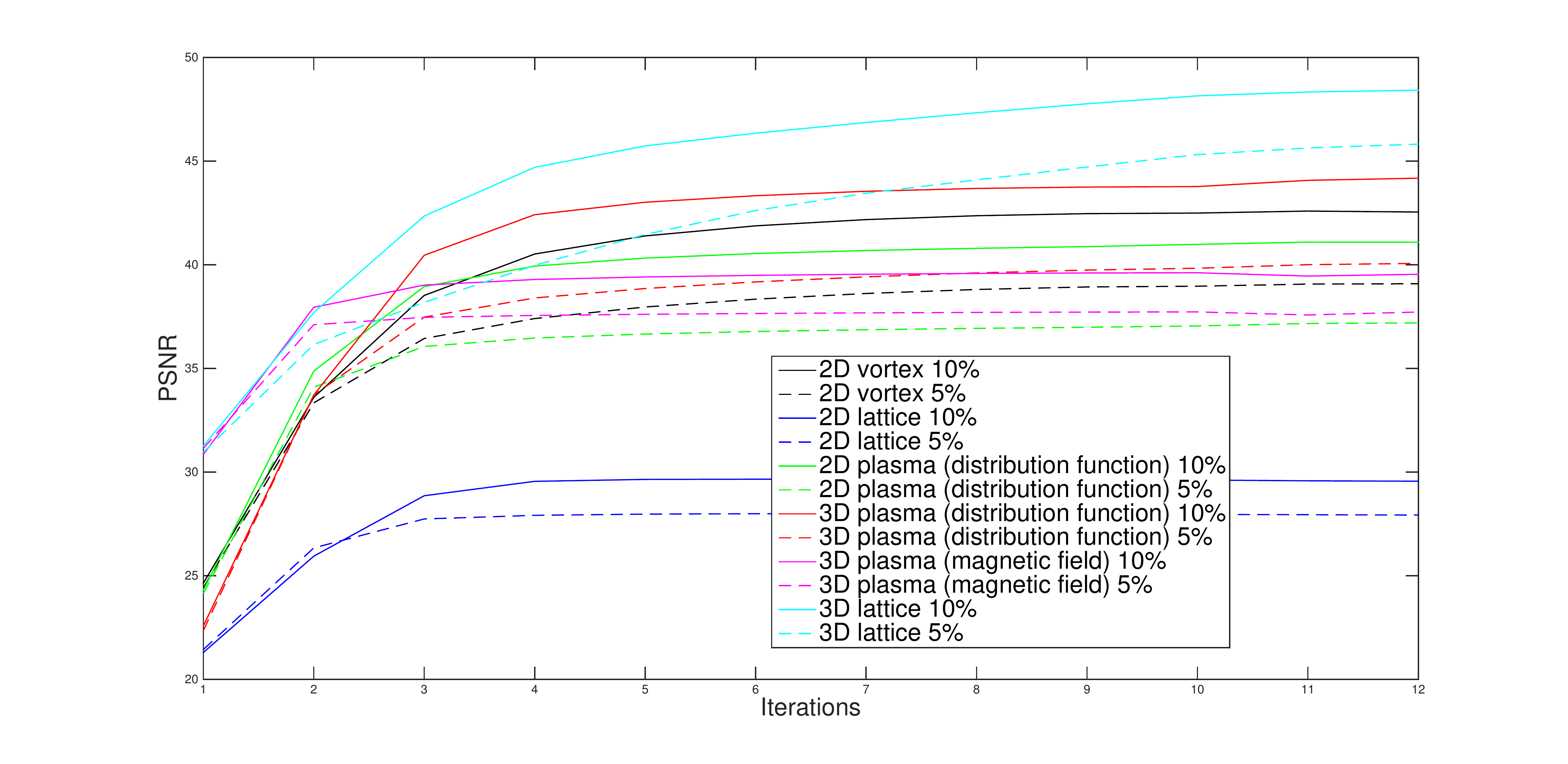}
  \caption{Numerical convergence in PSNR of LDMM on various data sets.}
  \label{fig:psnr}
\end{figure}

The numerical convergence of LDMM in PSNR is shown in Figure \ref{fig:psnr}. It can be observed that the algorithm converges fairly fast, usually within 10 iterations, and the result does not deteriorate as the iteration goes on.

\begin{table}[H]
  \centering
  \begin{tabular}{||c| c  c c||c|  c c c||}
    \hline
    $10\%$ & EBI & PLE& LDMM & $5\%$& EBI & PLE & LDMM\\
    \hline
    $L_1$ & 0.0082 & 0.0053 & \textbf{0.0034} & $L_1$ & 0.0148 & 0.0296 & \textbf{0.0056}\\
    \hline
    $L_2$ & 0.0153 & 0.0100 & \textbf{0.0075} & $L_2$ & 0.0270 & 0.0573 & \textbf{0.0111}\\
    \hline
    $L_\infty$ & 0.2280 & \textbf{0.1232} & 0.1376 & $L_\infty$ & 0.3327 & 0.7872 & \textbf{0.1102}\\
    \hline
    PSNR & 36.32 & 40.01 & \textbf{42.55} & PSNR & 31.36 & 24.84 & \textbf{39.09}\\
    \hline
  \end{tabular}
  \caption{Errors of the interpolation of the 2D vortex data set from $10\%$ and $5\%$ random sampling.}
  \label{tab:error_random_antonio_2d}
\end{table}

\begin{table}[H]
  \centering
  \begin{tabular}{||c| c  c c||c|  c c c||}
    \hline
    $10\%$ & EBI & PLE& LDMM & $5\%$& EBI & PLE & LDMM\\
    \hline
    $L_1$ & 0.0335 & 0.0272 & \textbf{0.0243} & $L_1$ & 0.0393 & 0.0535 & \textbf{0.0303}\\
    \hline
    $L_2$ & 0.0459 & 0.0377 & \textbf{0.0333} & $L_2$ & 0.0522 & 0.0805 & \textbf{0.0401}\\
    \hline
    $L_\infty$ & 0.3782 & 0.2158 & \textbf{0.1882} & $L_\infty$ & 0.2588 & 0.7148 & \textbf{0.2063}\\
    \hline
    PSNR & 26.77 & 28.48 & \textbf{29.56} & PSNR & 25.65 & 21.88 & \textbf{27.93}\\
    \hline
  \end{tabular}
  \caption{Errors of the interpolation of the 2D plasma (distribution function) data set from $10\%$ and $5\%$ random sampling.}
  \label{tab:error_random_shock_2d}
\end{table}

\begin{table}[H]
  \centering
  \begin{tabular}{||c| c  c c||c|  c c c||}
    \hline
    $10\%$ & EBI & PLE& LDMM & $5\%$& EBI & PLE & LDMM\\
    \hline
    $L_1$ & 0.0033 & 0.0030 & \textbf{0.0013} & $L_1$ & 0.0048 & 0.0187 & \textbf{0.0022}\\
    \hline
    $L_2$ & 0.0066 & 0.0077 & \textbf{0.0040} & $L_2$ & 0.0124 & 0.0442 & \textbf{0.0062}\\
    \hline
    $L_\infty$ & 0.2172 & 0.2889 & \textbf{0.1979} & $L_\infty$ & 0.8758 & 0.6156 & \textbf{0.2097}\\
    \hline
    PSNR & 43.62 & 42.32 & \textbf{47.98} & PSNR & 38.16 & 27.08 & \textbf{44.15}\\
    \hline
  \end{tabular}
  \caption{Errors of the interpolation of the 2D lattice data set from $10\%$ and $5\%$ random sampling.}
  \label{tab:error_random_lattice_2d}
\end{table}

\begin{table}[H]
  \centering
  \begin{tabular}{||c| c  c c||c|  c c c||}
    \hline
    $10\%$ & EBI & PLE& LDMM & $5\%$& EBI & PLE & LDMM\\
    \hline
    $L_1$ & 0.0075 & 0.0053 & \textbf{0.0038} & $L_1$ & 0.0115 & 0.0285 & \textbf{0.0062}\\
    \hline
    $L_2$ & 0.0128 & 0.0126 & \textbf{0.0062} & $L_2$ & 0.0201 & 0.0513 & \textbf{0.0099}\\
    \hline
    $L_\infty$ & 0.3510 & 0.9432 & \textbf{0.1330} & $L_\infty$ & 0.3740 & 0.7531 & \textbf{0.2012}\\
    \hline
    PSNR & 37.88 & 37.96 & \textbf{44.18} & PSNR & 33.93 & 25.80 & \textbf{40.07}\\
    \hline
  \end{tabular}
  \caption{Errors of the interpolation of the 3D plasma (magnetic field) data set from $10\%$ and $5\%$ random sampling.}
  \label{tab:error_random_plasma_3d}
\end{table}

\begin{table}[H]
  \centering
  \begin{tabular}{||c| c  c c||c|  c c c||}
    \hline
    $10\%$ & EBI & PLE& LDMM & $5\%$& EBI & PLE & LDMM\\
    \hline
    $L_1$ & 0.0094 & 0.0062 & \textbf{0.0008} & $L_1$ & 0.0112 & 0.0545 & \textbf{0.0013}\\
    \hline
    $L_2$ & 0.0308 & 0.0166 & \textbf{0.0038} & $L_2$ & 0.0336 & 0.0899 & \textbf{0.0051}\\
    \hline
    $L_\infty$ & 0.5291 & 0.6635 & \textbf{0.4262} & $L_\infty$ & 0.4768 & 0.0.8595 & \textbf{0.4530}\\
    \hline
    PSNR & 30.24 & 35.60 & \textbf{48.43} & PSNR & 29.48 & 20.93 & \textbf{45.82}\\
    \hline
  \end{tabular}
  \caption{Errors of the interpolation of the 3D lattice data set from $10\%$ and $5\%$ random sampling.}
  \label{tab:error_random_lattice_3d}
\end{table}

\begin{table}[H]
  \centering
  \begin{tabular}{||c| c  c c||c|  c c c||}
    \hline
    $10\%$ & EBI & PLE& LDMM & $5\%$& EBI & PLE & LDMM\\
    \hline
    $L_1$ & 0.0098 & 0.0085 & \textbf{0.0060} & $L_1$ & 0.0108 & 0.0593 & \textbf{0.0075}\\
    \hline
    $L_2$ & 0.0167 & 0.0138 & \textbf{0.0105} & $L_2$ & 0.0181 & 0.0895 & \textbf{0.0130}\\
    \hline
    $L_\infty$ & 0.2005 & 0.2912 & \textbf{0.1181} & $L_\infty$ & 0.1865 & 0.9093 & \textbf{0.1793}\\
    \hline
    PSNR & 35.54 & 37.20 & \textbf{39.54} & PSNR & 34.87 & 20.96 & \textbf{37.72}\\
    \hline
  \end{tabular}
  \caption{Errors of the interpolation of the 3D plasma (distribution function) data set from $10\%$ and $5\%$ random sampling.}
  \label{tab:error_random_shock_3d}
\end{table}

\subsection{Interpolation with Regular Sampling}

Unlike the random sampling interpolation in the previous section, reasonable initializations of LDMM can be obtained from other standard algorithms for regular sampling interpolation. In the numerical experiments on all the data sets, the results of DCT and cubic spline have been used as the initial iterates for LDMM, and the final results of LDMM initialized with DCT (LDMM (D)) and cubic spline (LDMM (C))  are obtained after three iterations of manifold updates.

%marker

\begin{figure}[H]
  \centering
  \begin{tabular}{ccc}
    Original& Cubic Spline (42.98dB)& DCT (42.88dB)\\
    \includegraphics[width=3cm]{antonio_2d_original}&
    \includegraphics[width=3cm]{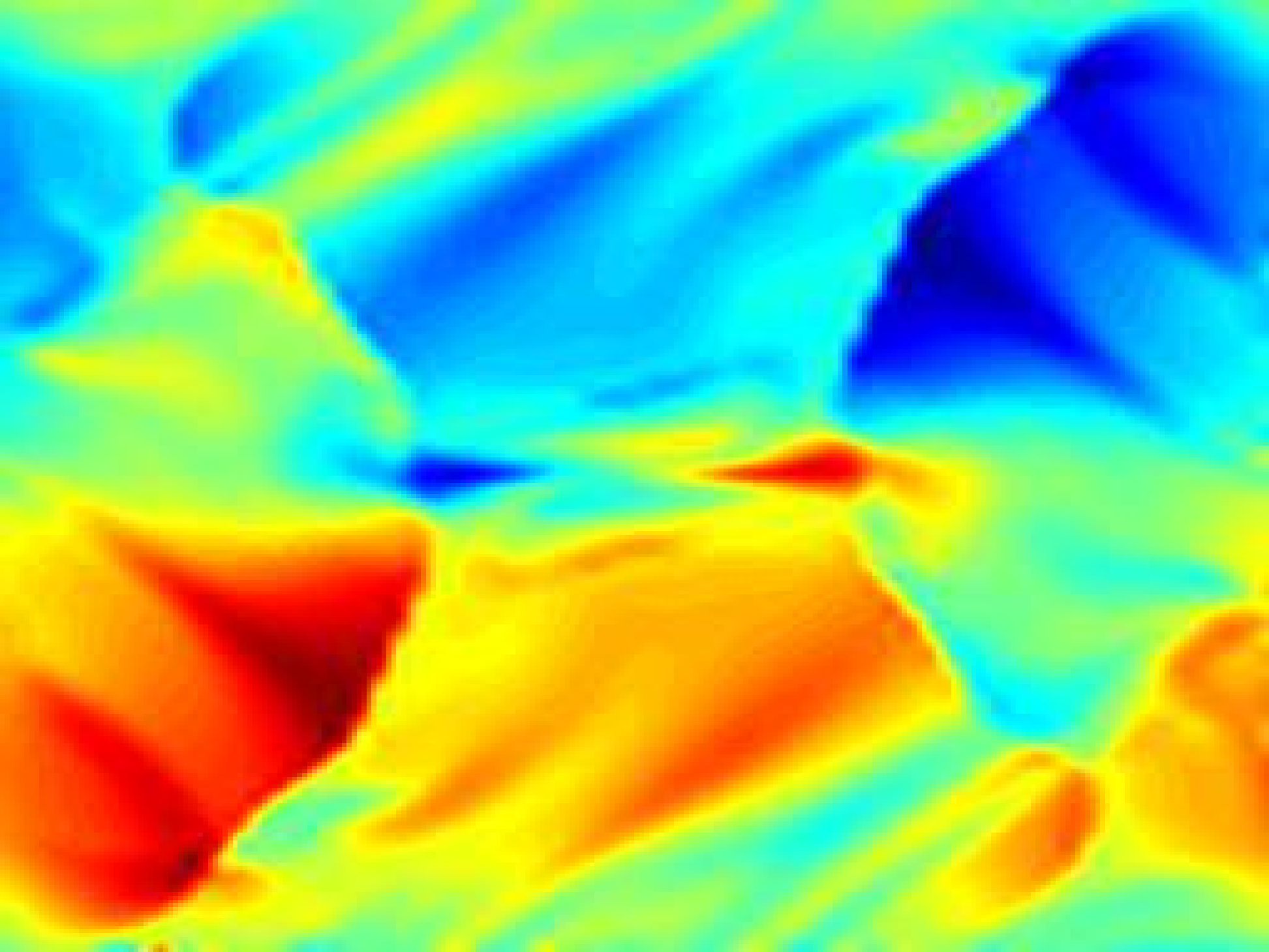}&
    \includegraphics[width=3cm]{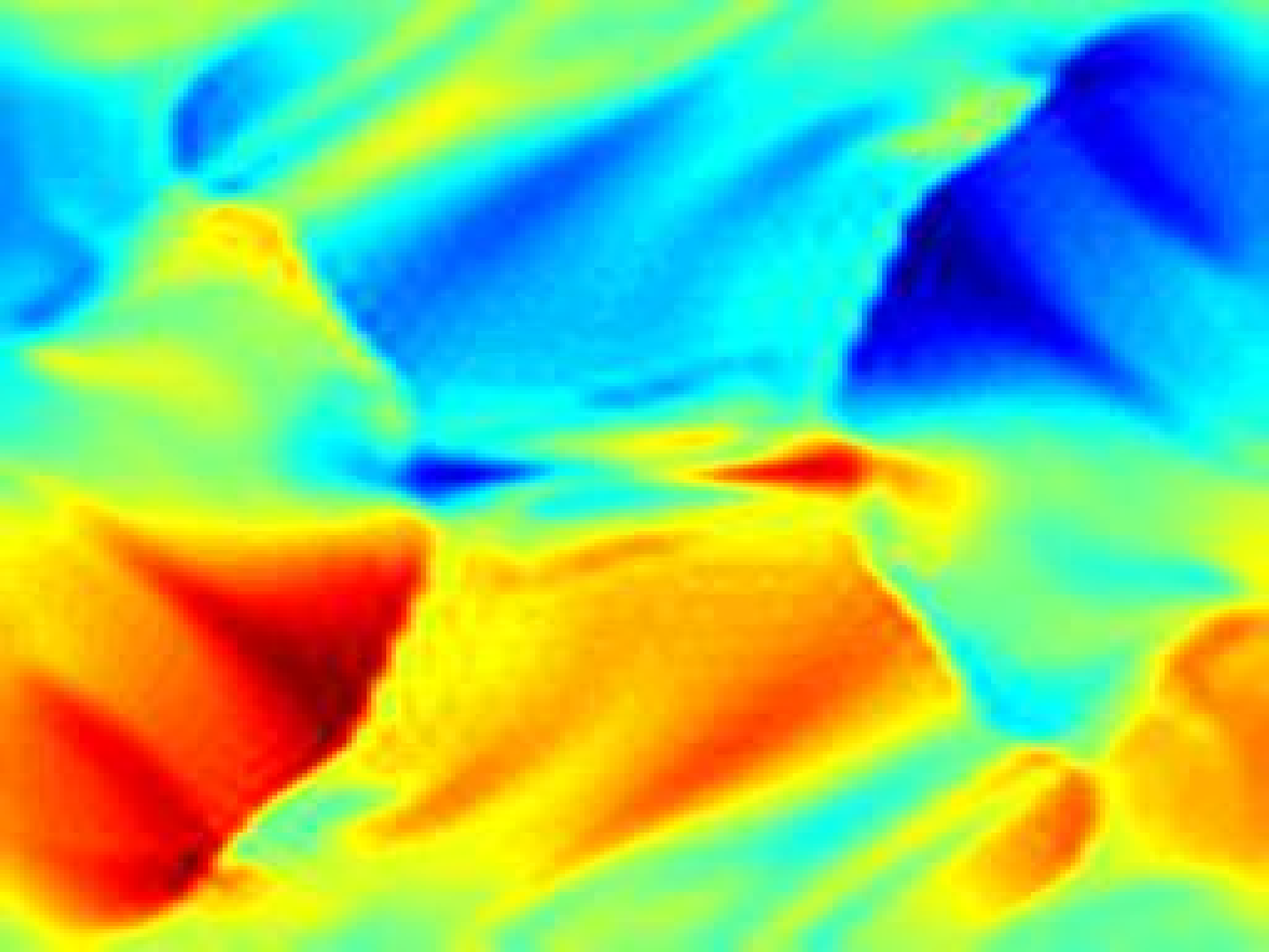}\\
    DFT (43.19dB)& Wavelet (40.48dB)& LDMM (\textbf{44.40dB})\\    \vspace{.5cm}
    \includegraphics[width=3cm]{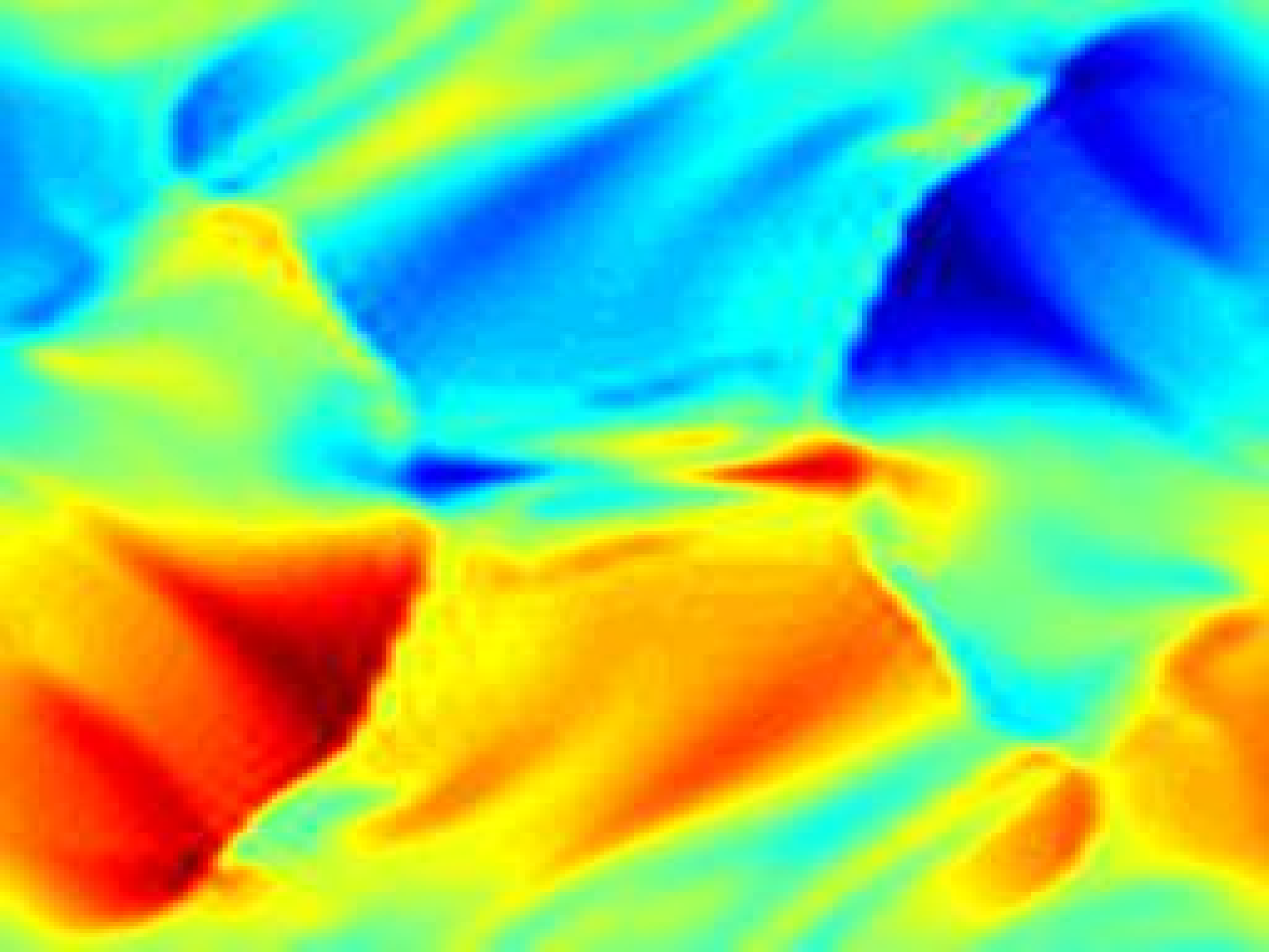}&
    \includegraphics[width=3cm]{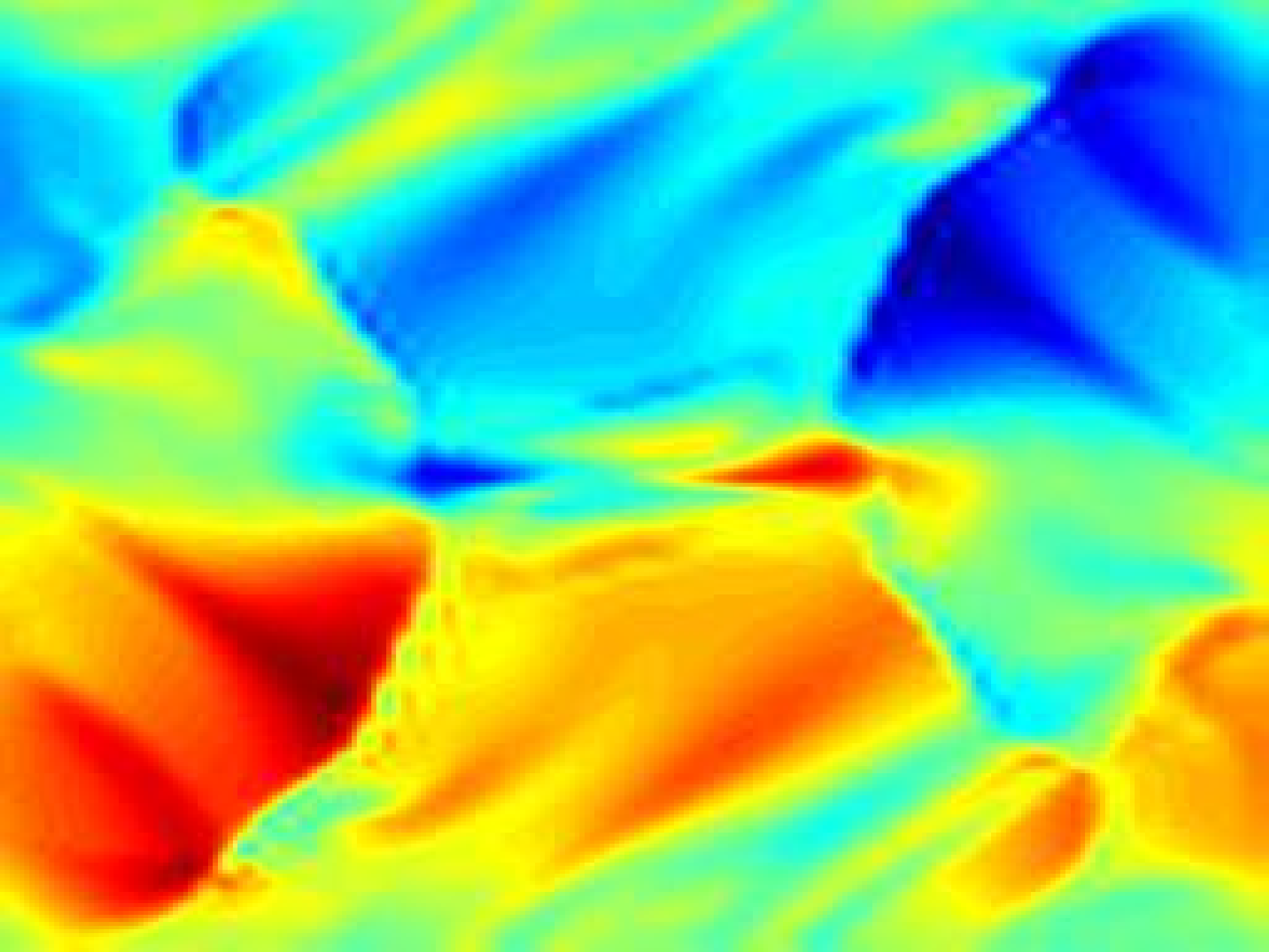}&
    \includegraphics[width=3cm]{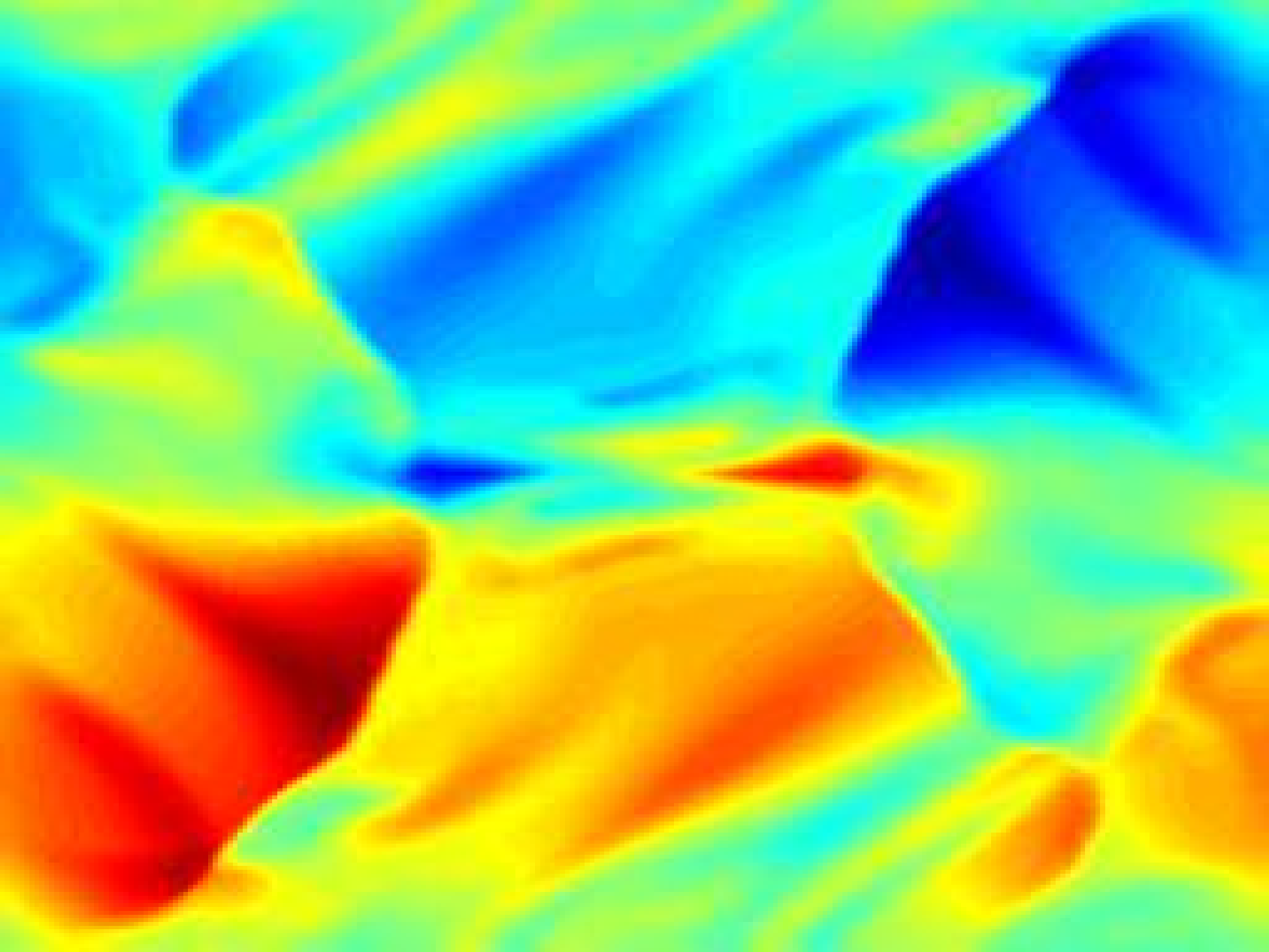}\\
    Original& Cubic Spline (26.81dB)& DCT (27.68dB)\\
    \includegraphics[width=3cm]{shock_2d_original}&
    \includegraphics[width=3cm]{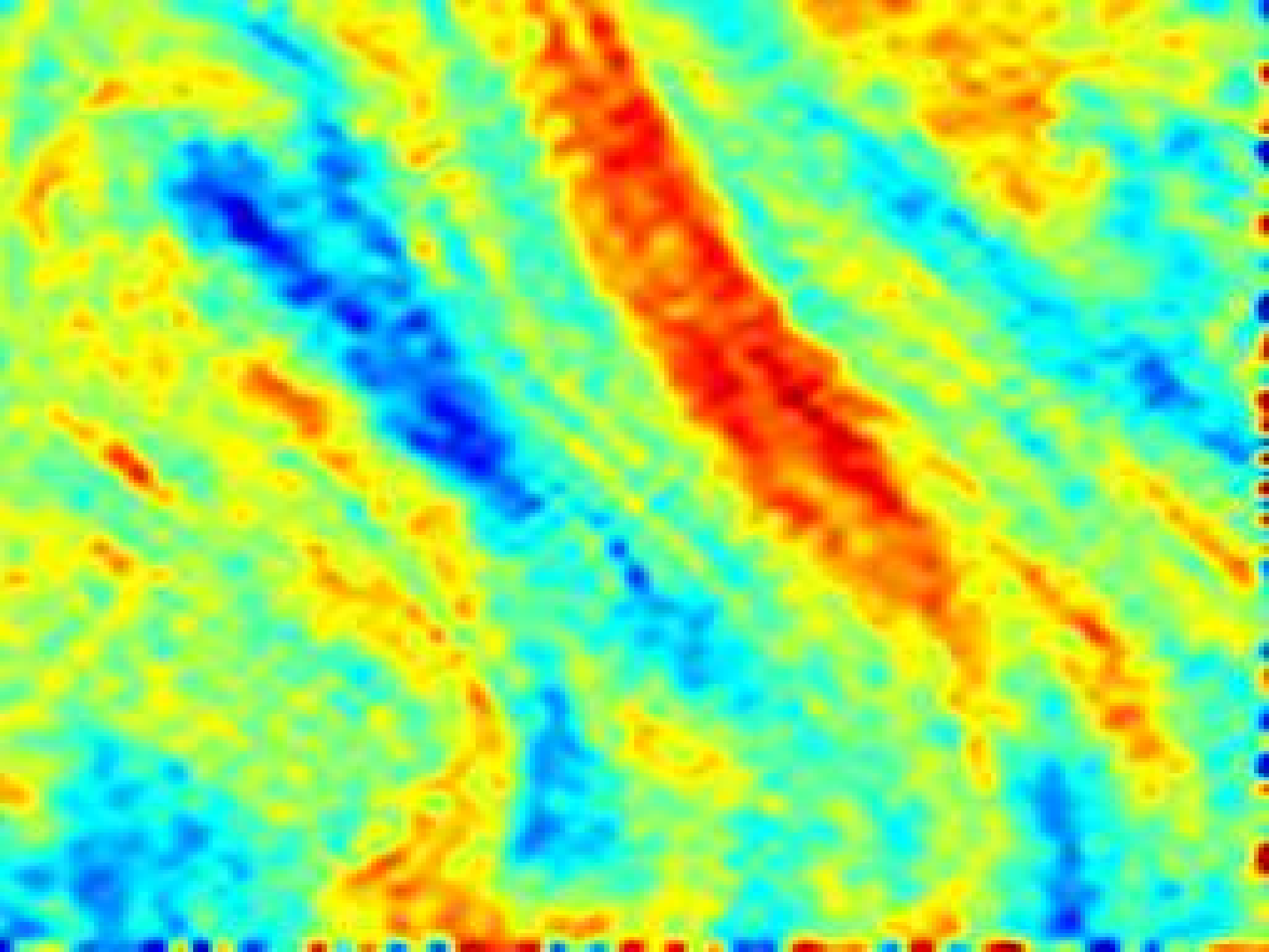}&
    \includegraphics[width=3cm]{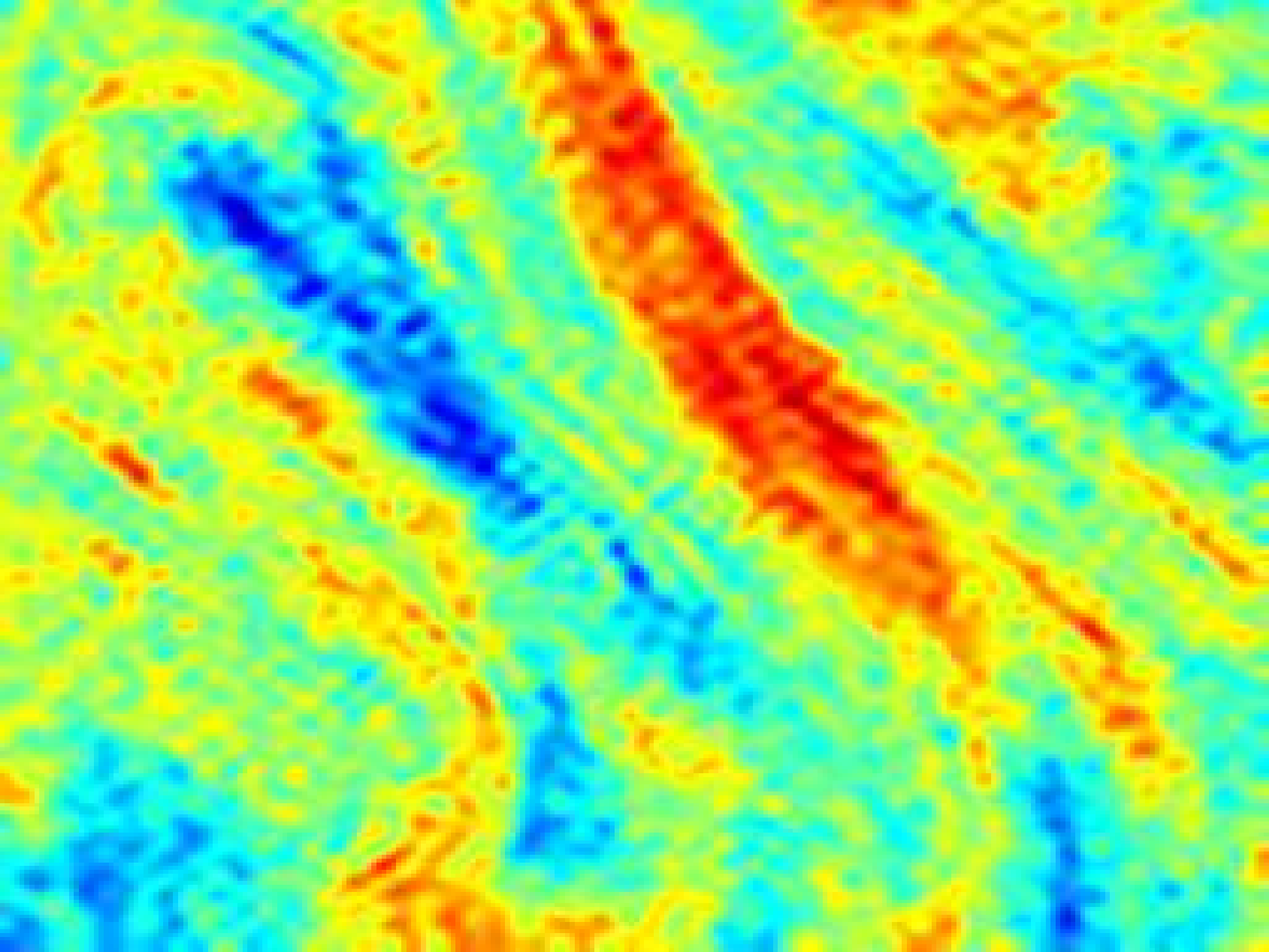}\\
    DFT (27.43dB)& Wavelet (27.34dB)& LDMM (\textbf{29.66dB})\\    \vspace{.5cm}
    \includegraphics[width=3cm]{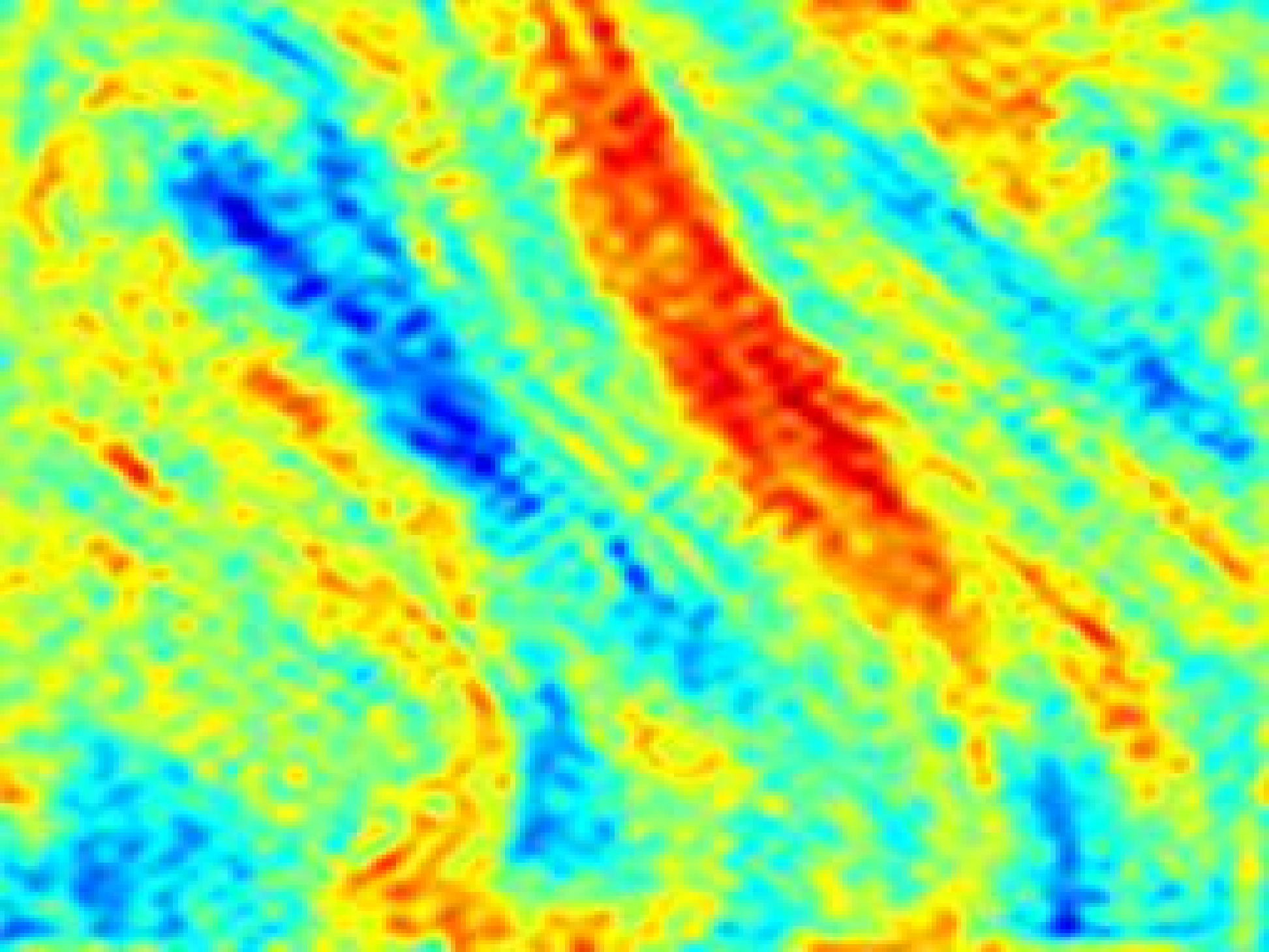}&
    \includegraphics[width=3cm]{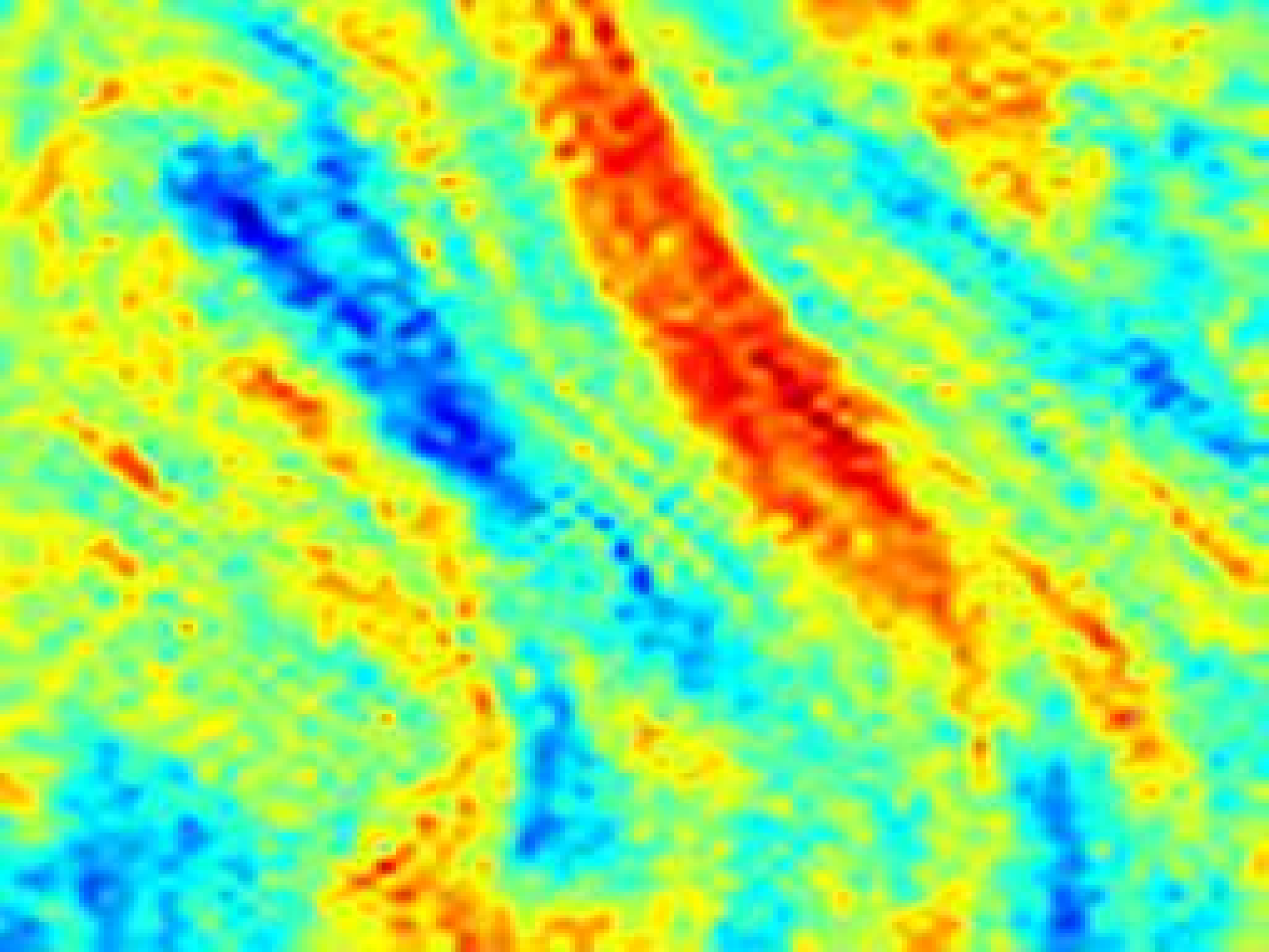}&
    \includegraphics[width=3cm]{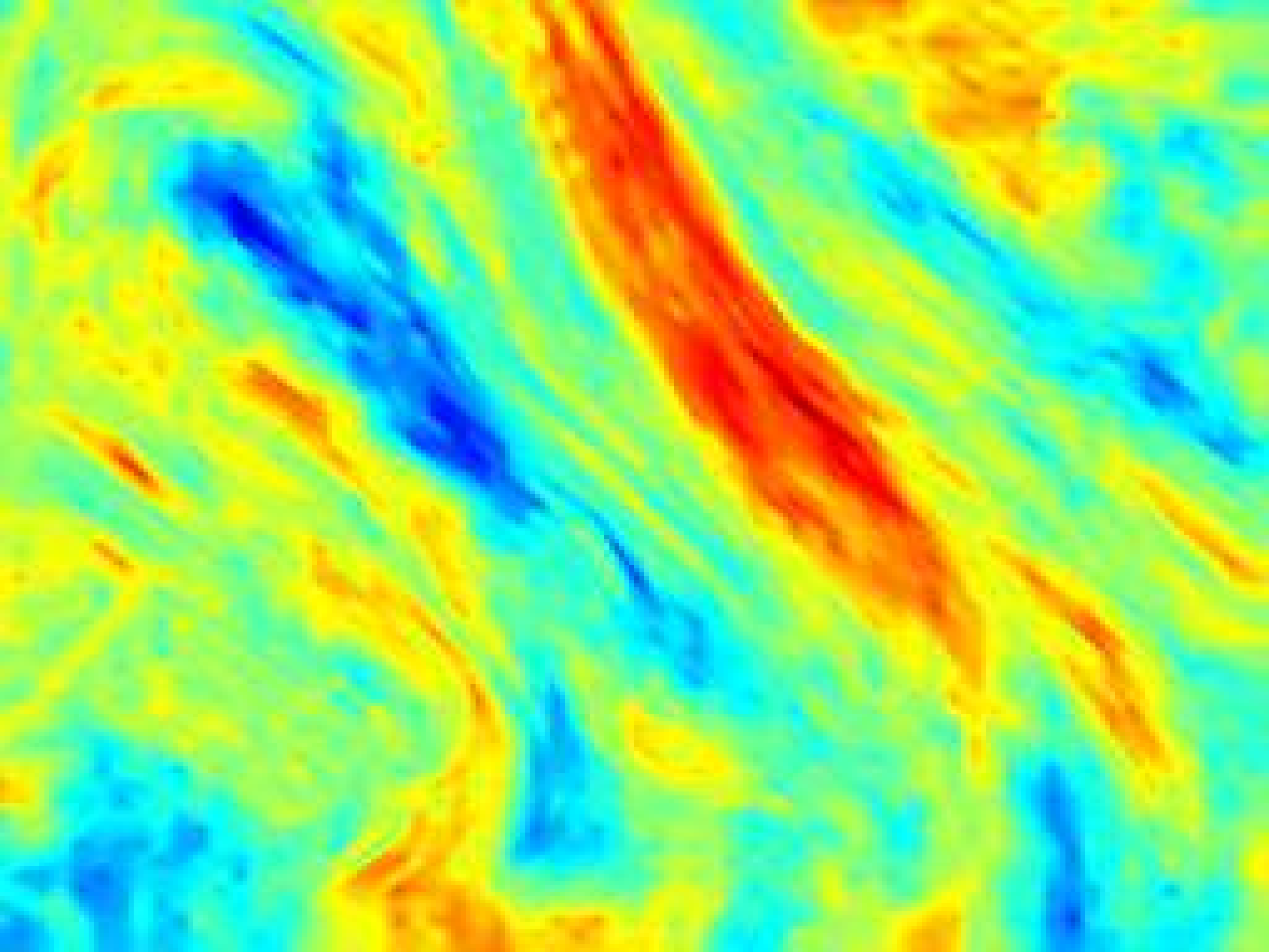}\\
    Original& Cubic Spline (46.97dB)& DCT (45.77dB)\\
    \includegraphics[width=3cm]{latticebig_2d_original}&
    \includegraphics[width=3cm]{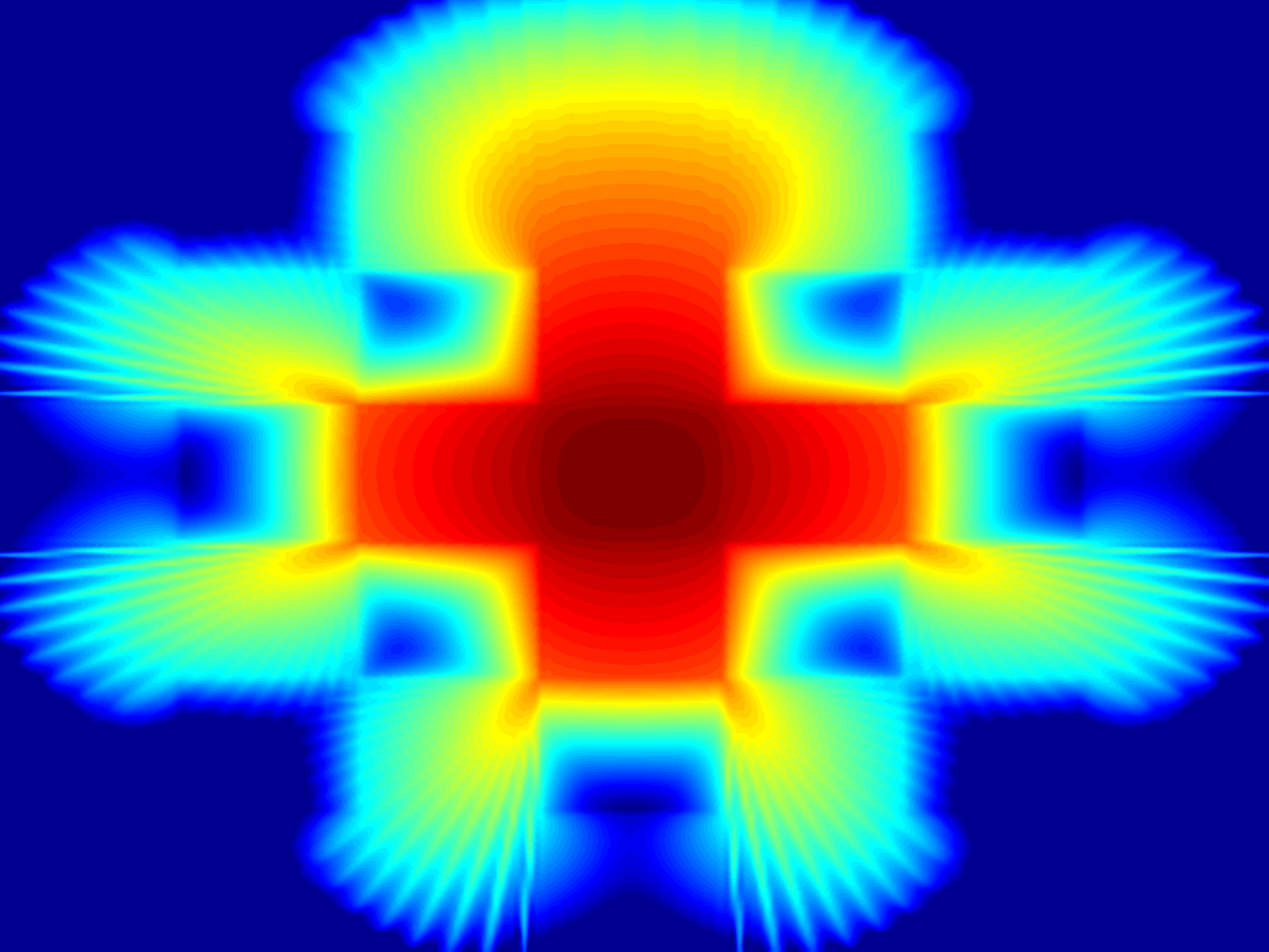}&
    \includegraphics[width=3cm]{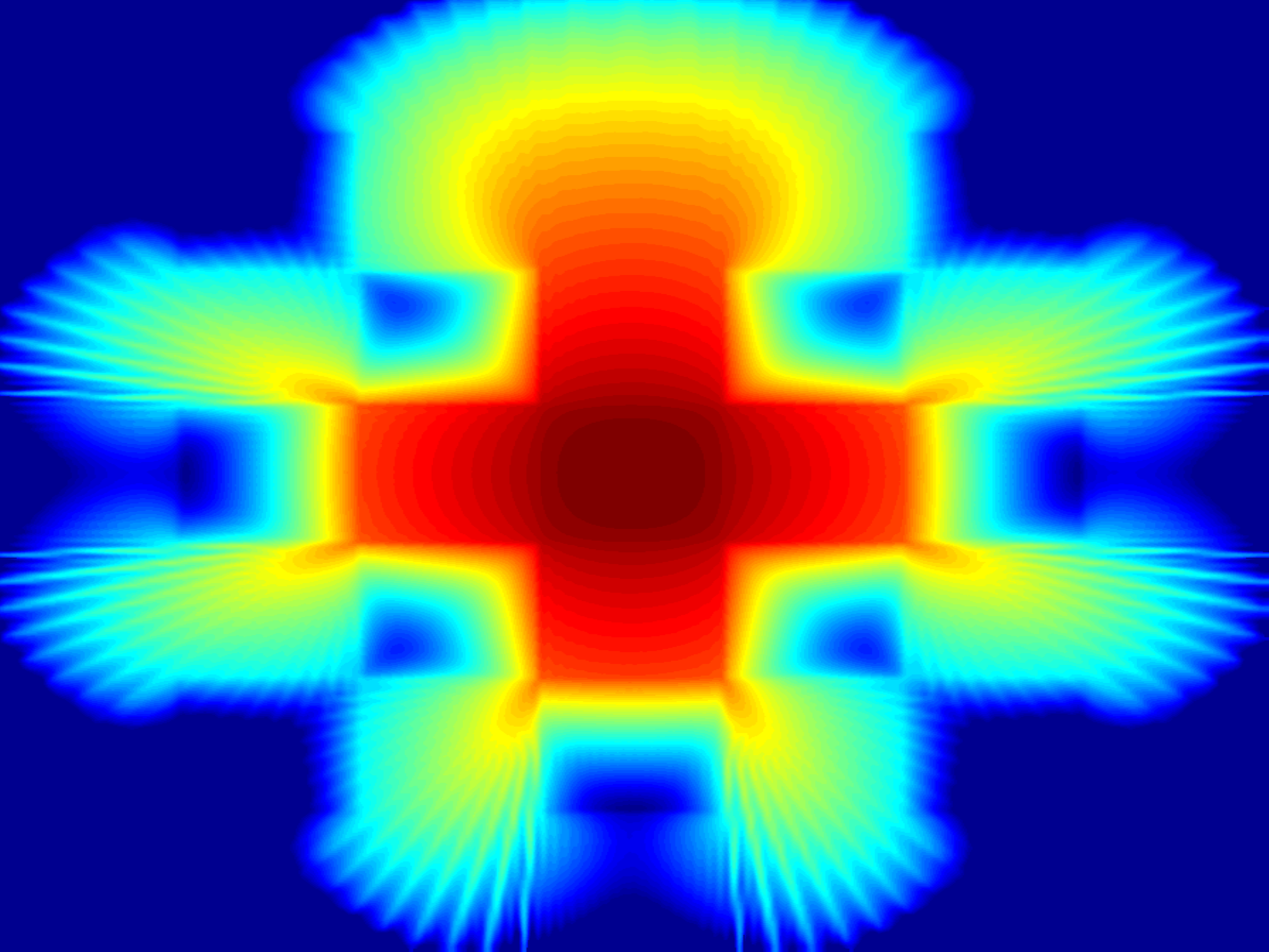}\\
    DFT (45.20dB)& Wavelet (44.31dB)& LDMM (\textbf{47.43dB})\\
    \includegraphics[width=3cm]{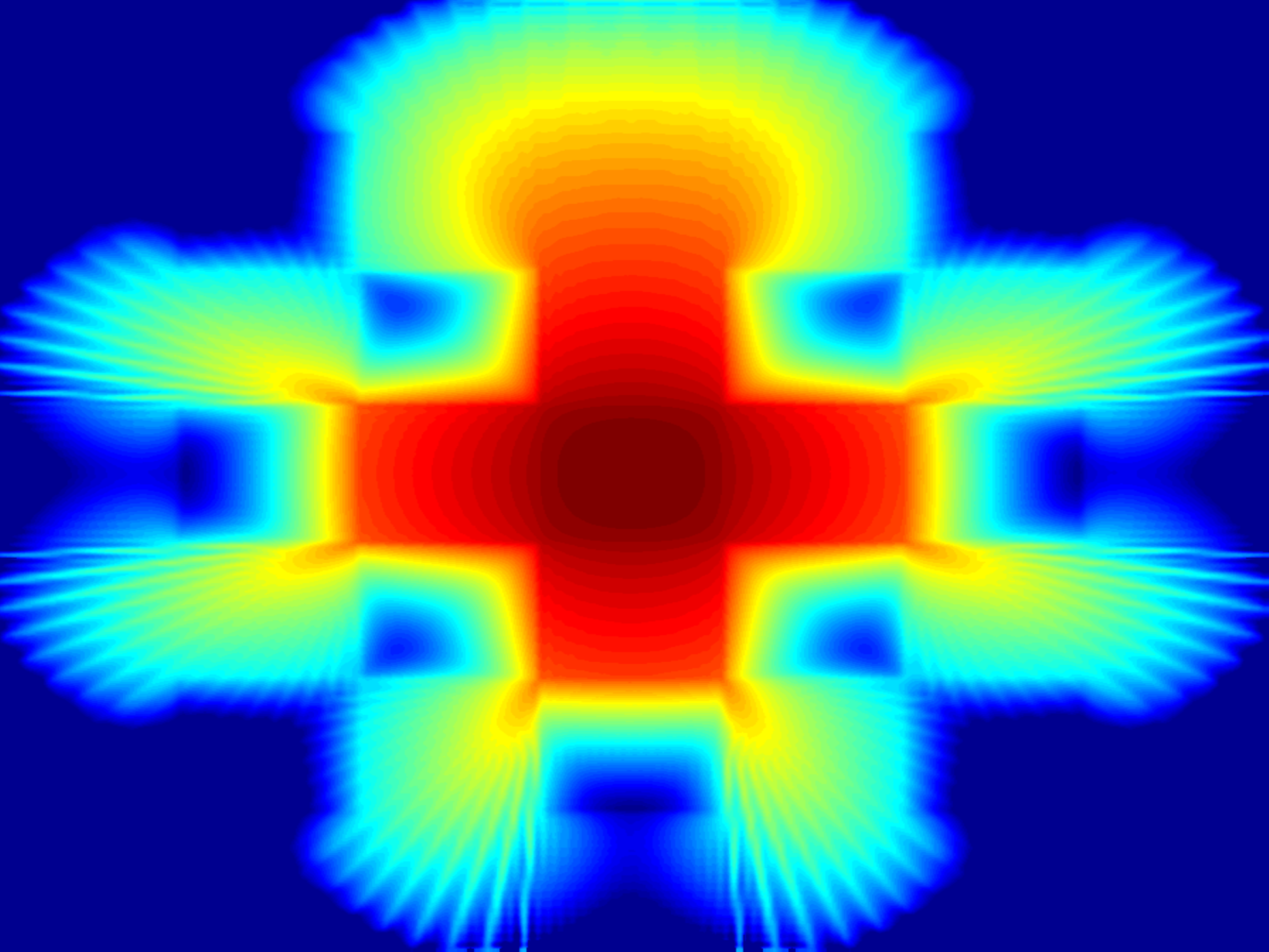}&
    \includegraphics[width=3cm]{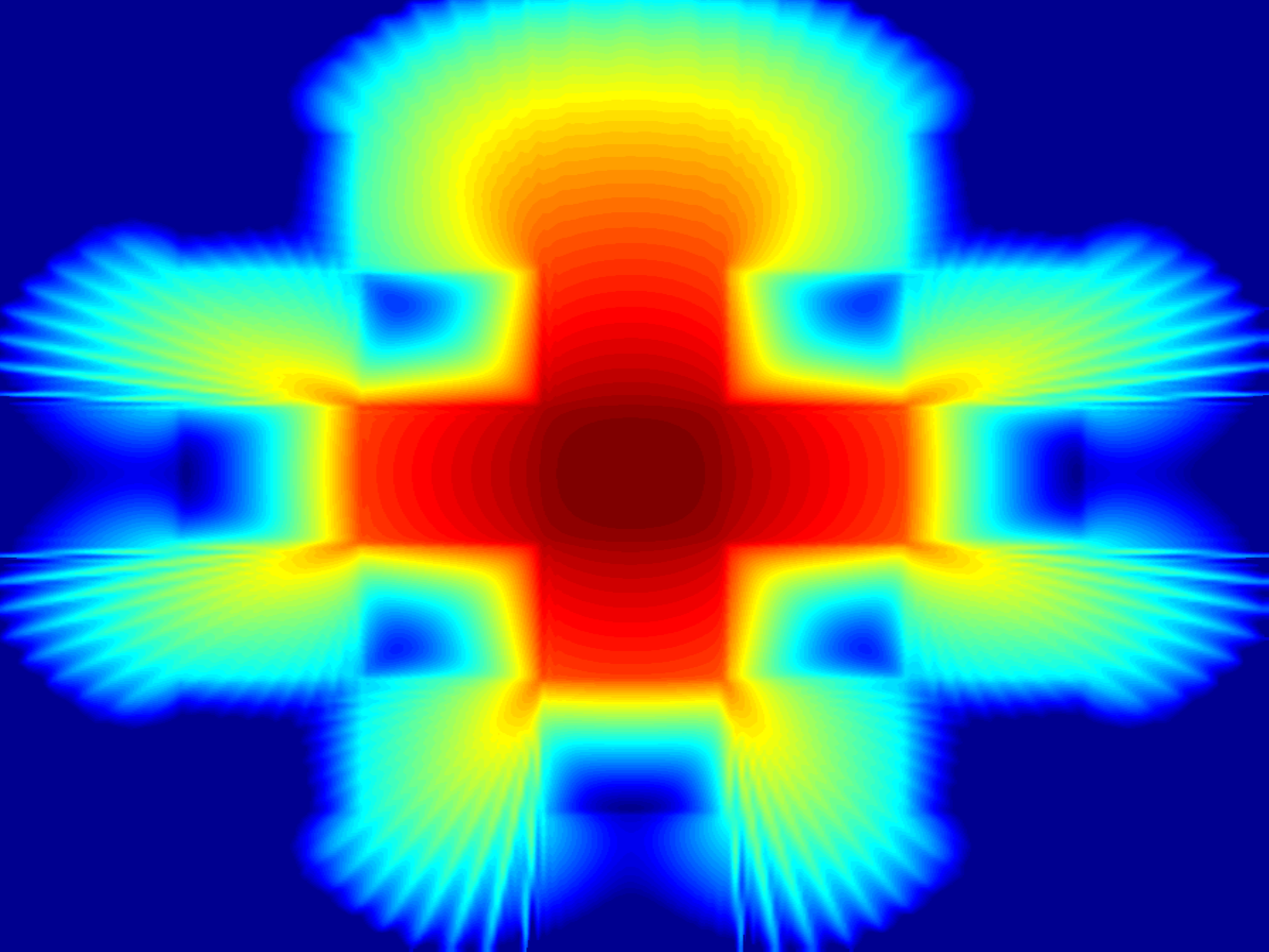}&
    \includegraphics[width=3cm]{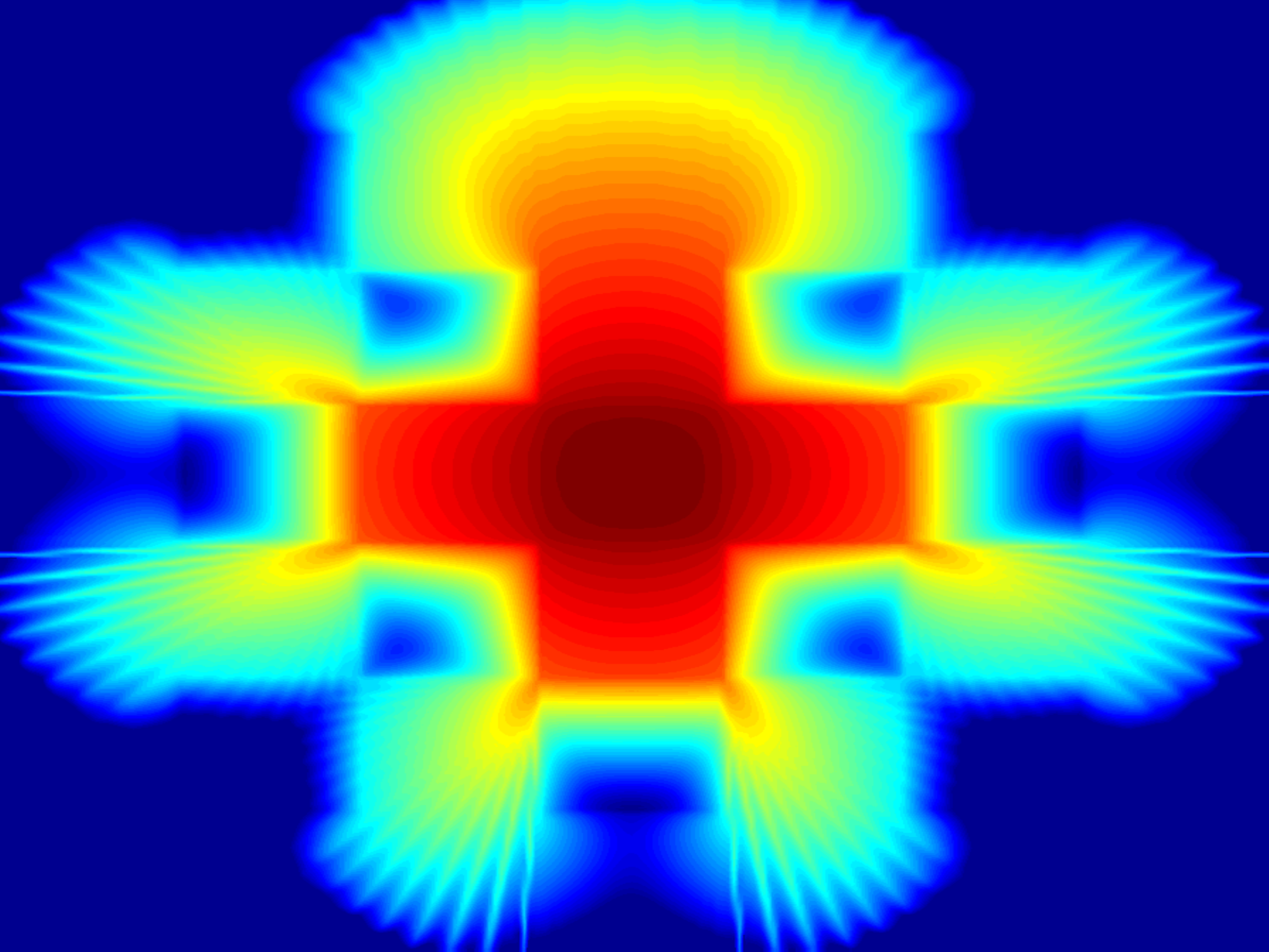}
  \end{tabular}
  \caption{Interpolation of 2D scientific data sets from regular sampling with spacing $4\times 4$. The original data are shown on the upper left corners for each data set. The results of cubic spline, DCT, DFT, wavelet, and LDMM are shown in the remaining five figures.}
  \label{fig:down_2d}
\end{figure}

\begin{figure}[H]
  \centering
  \begin{tabular}{ccc}
    Original& Cubic Spline (41.38dB)& DCT (43.76dB)\\
    \includegraphics[width=3cm]{plasma_3d_original_band_19}&
    \includegraphics[width=3cm]{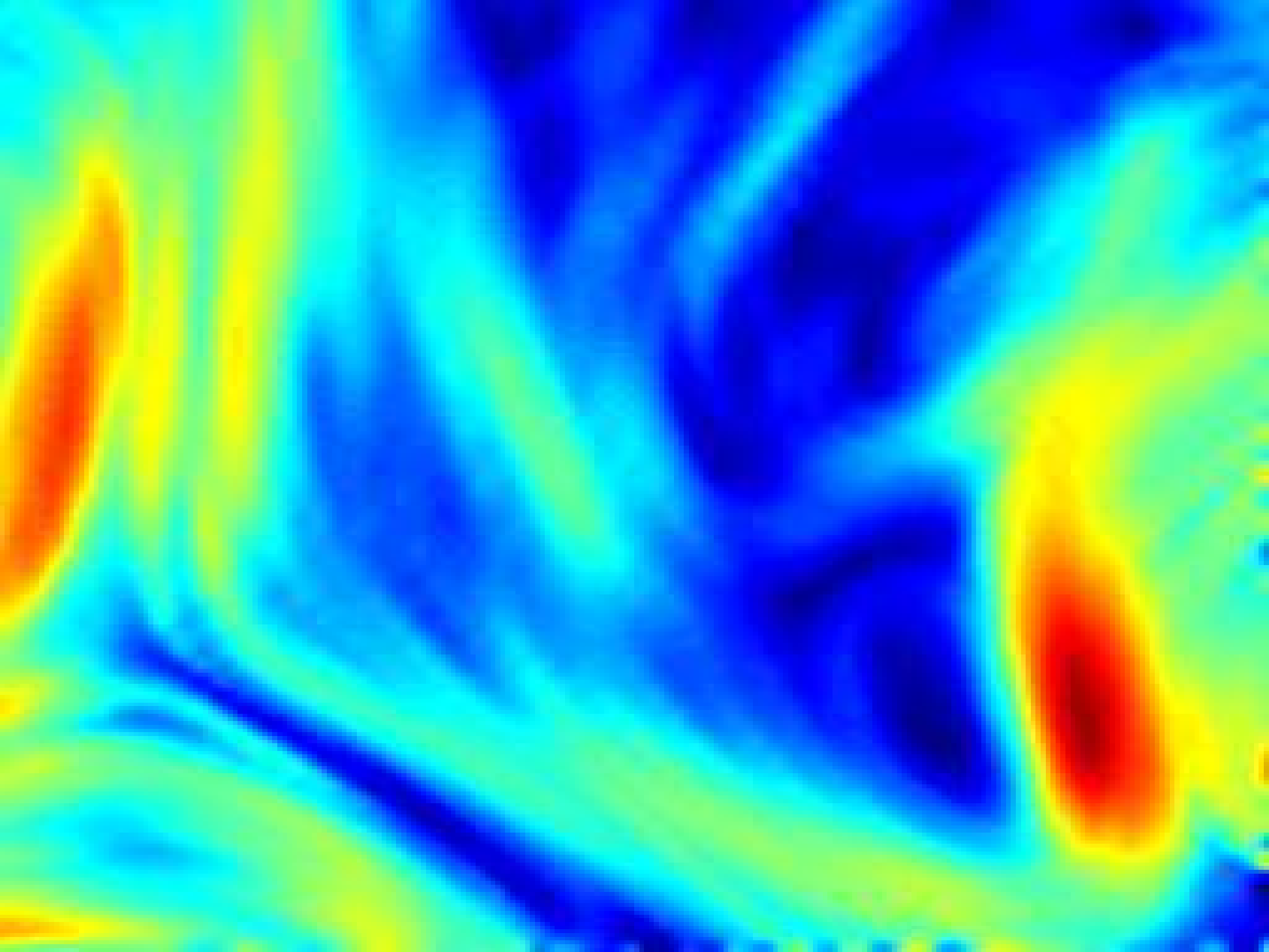}&
    \includegraphics[width=3cm]{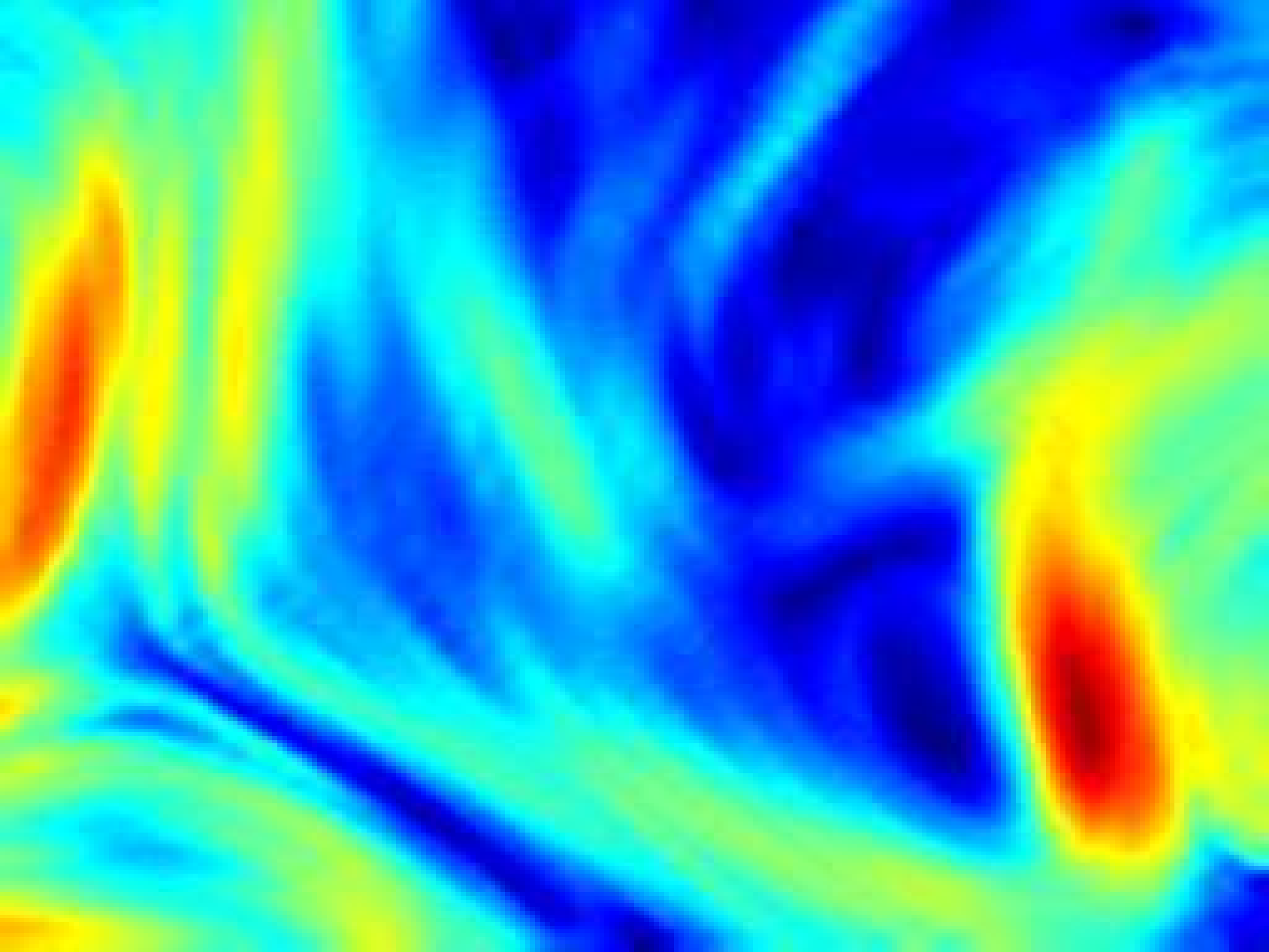}\\
    DFT (33.99dB)& Wavelet (42.15dB)& LDMM (\textbf{44.53dB})\\
    \includegraphics[width=3cm]{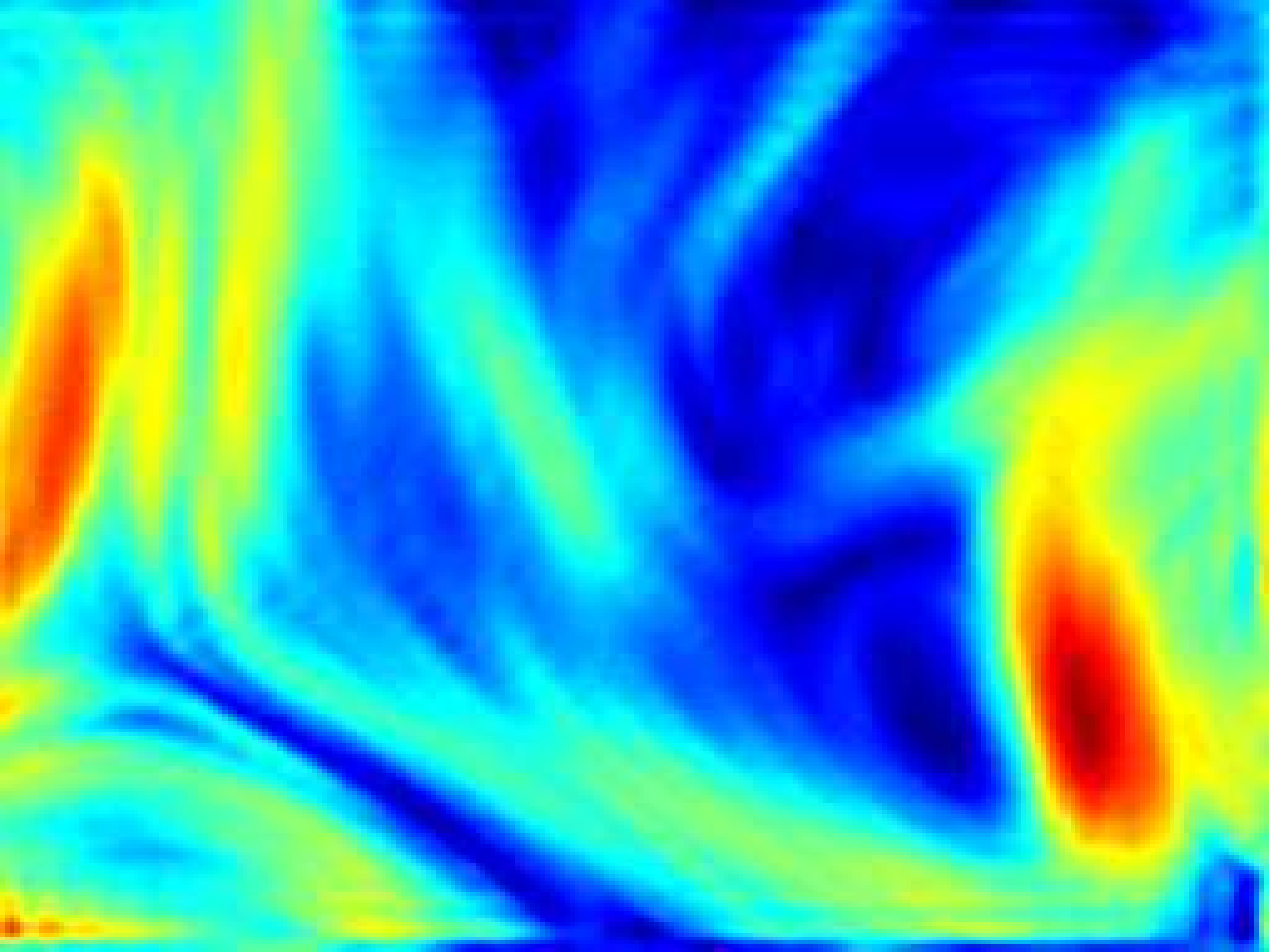}&
    \includegraphics[width=3cm]{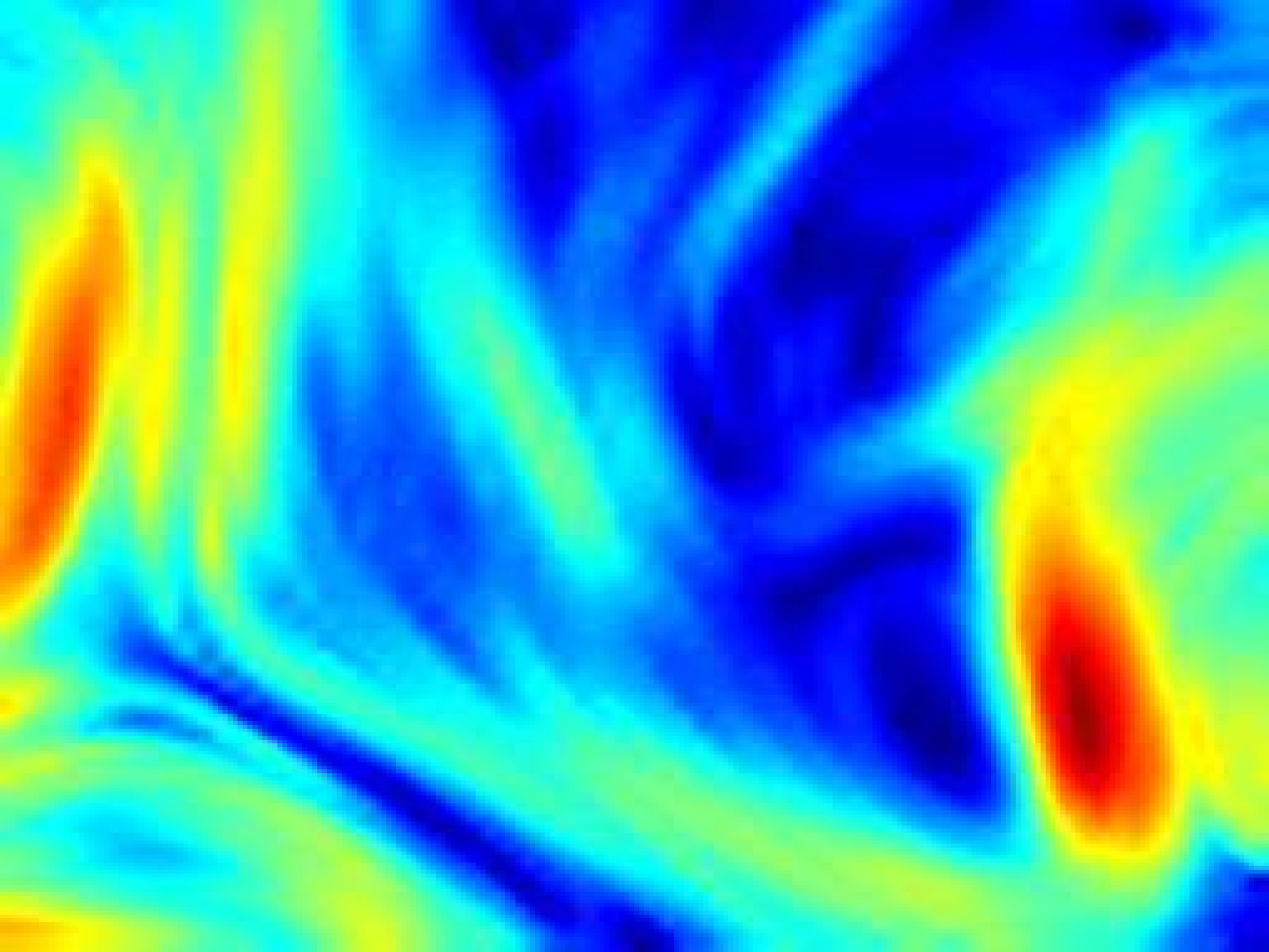}&
    \includegraphics[width=3cm]{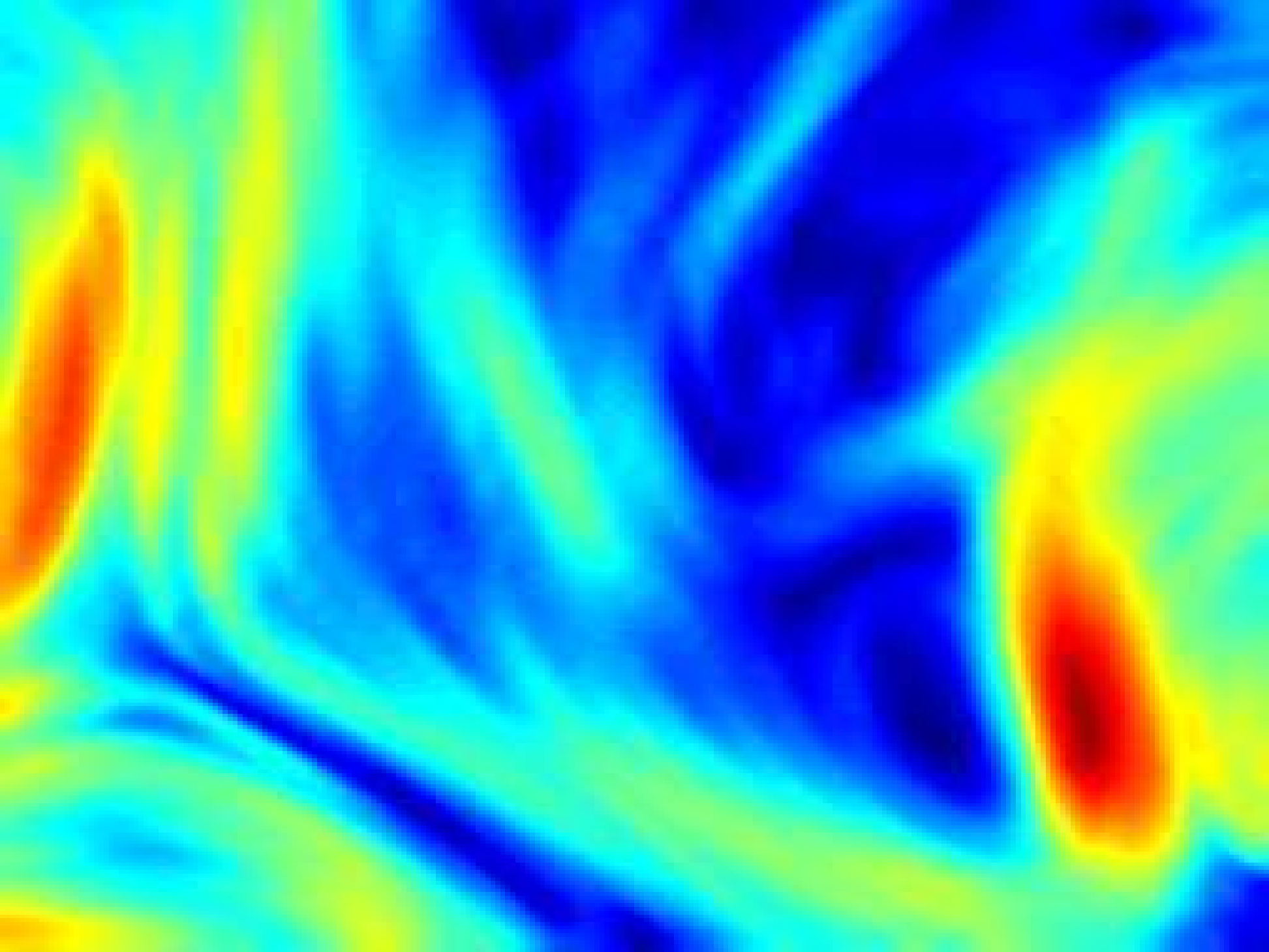}\\
    Original & Cubic Spline (41.38dB)& DCT (43.76dB)\\
    \includegraphics[width=3cm]{plasma_3d_original_band_29}&
    \includegraphics[width=3cm]{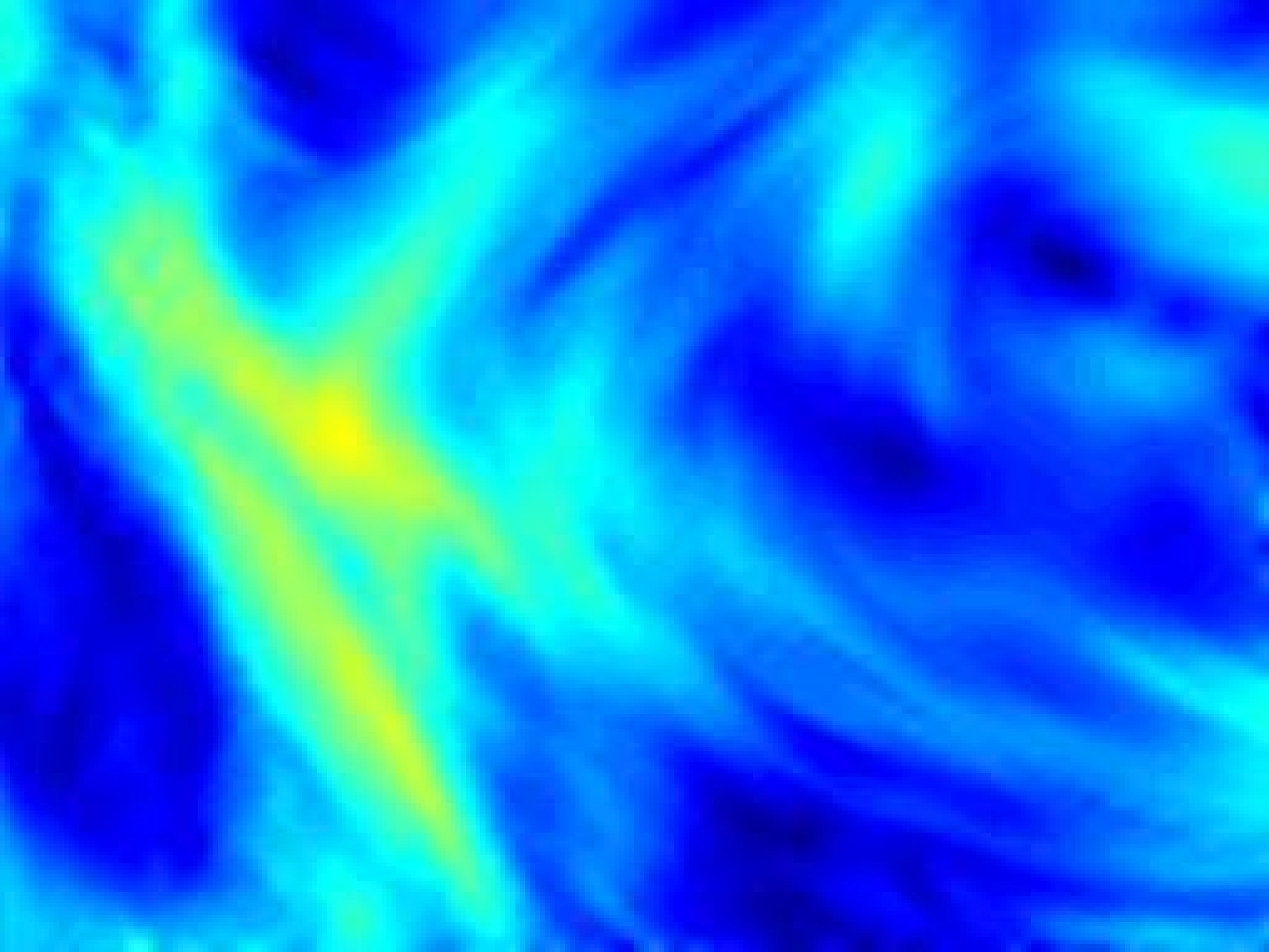}&
    \includegraphics[width=3cm]{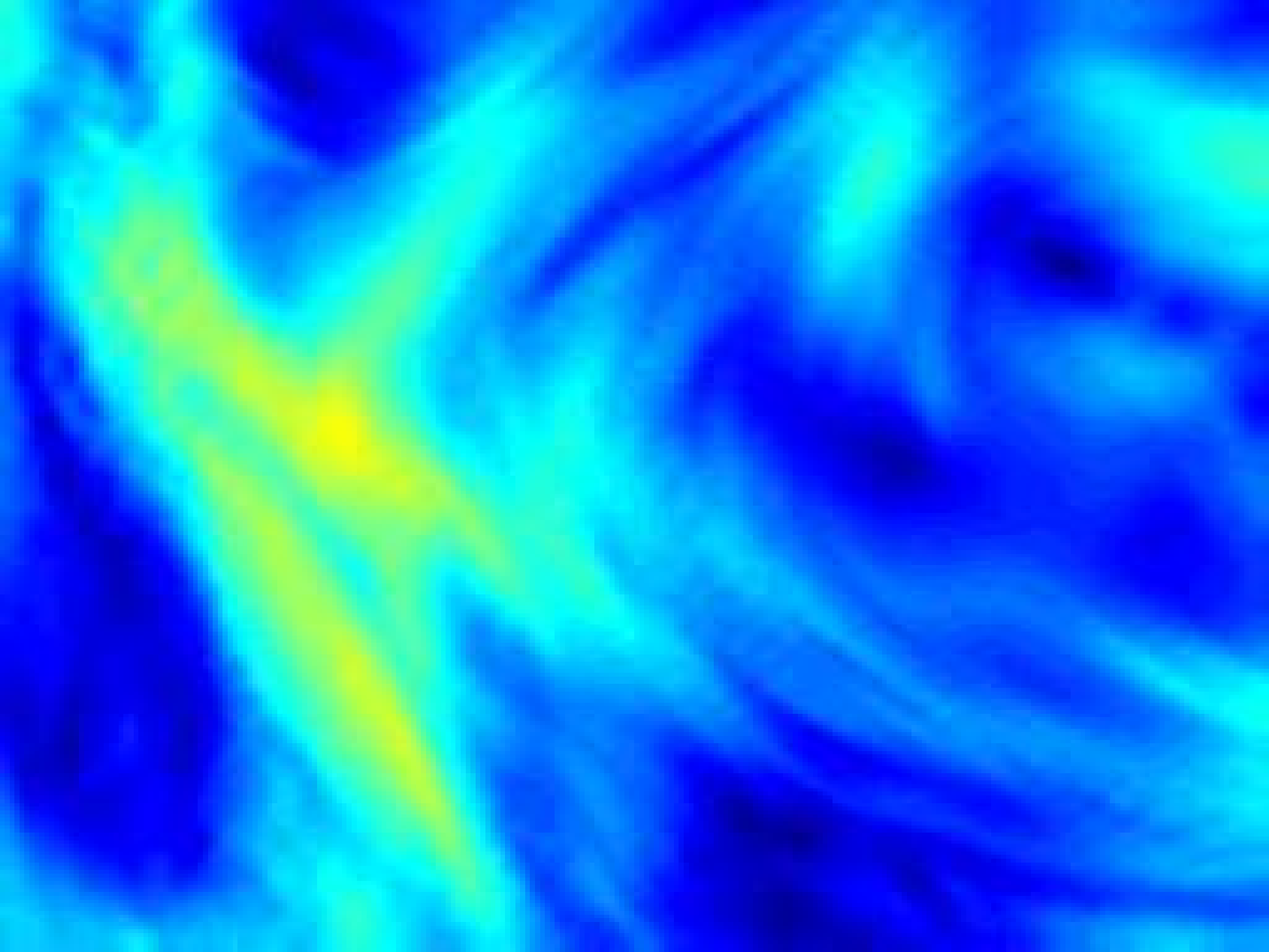}\\
    DFT (33.99dB)& Wavelet (42.15dB)& LDMM (\textbf{44.53dB})\\
    \includegraphics[width=3cm]{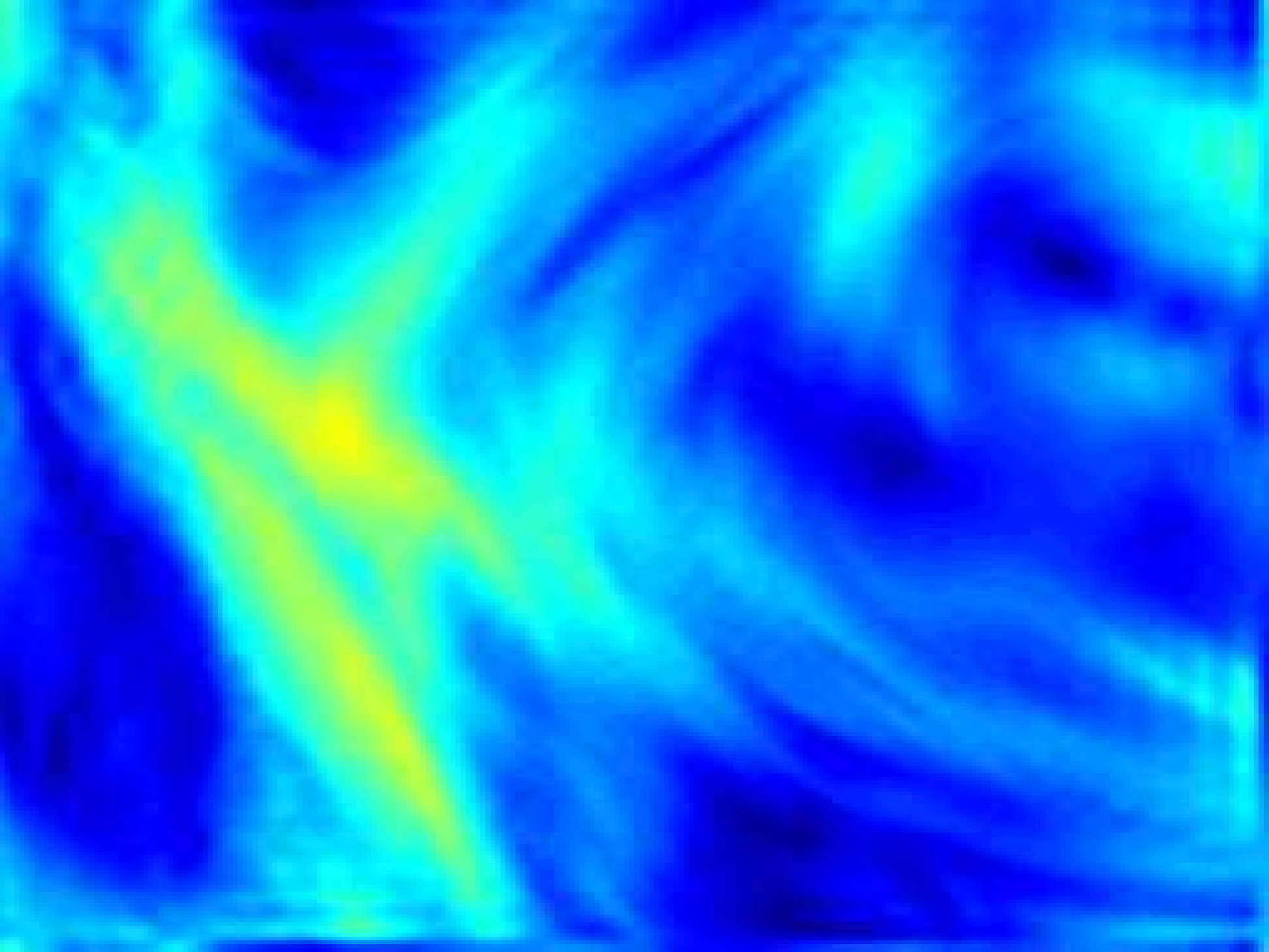}&
    \includegraphics[width=3cm]{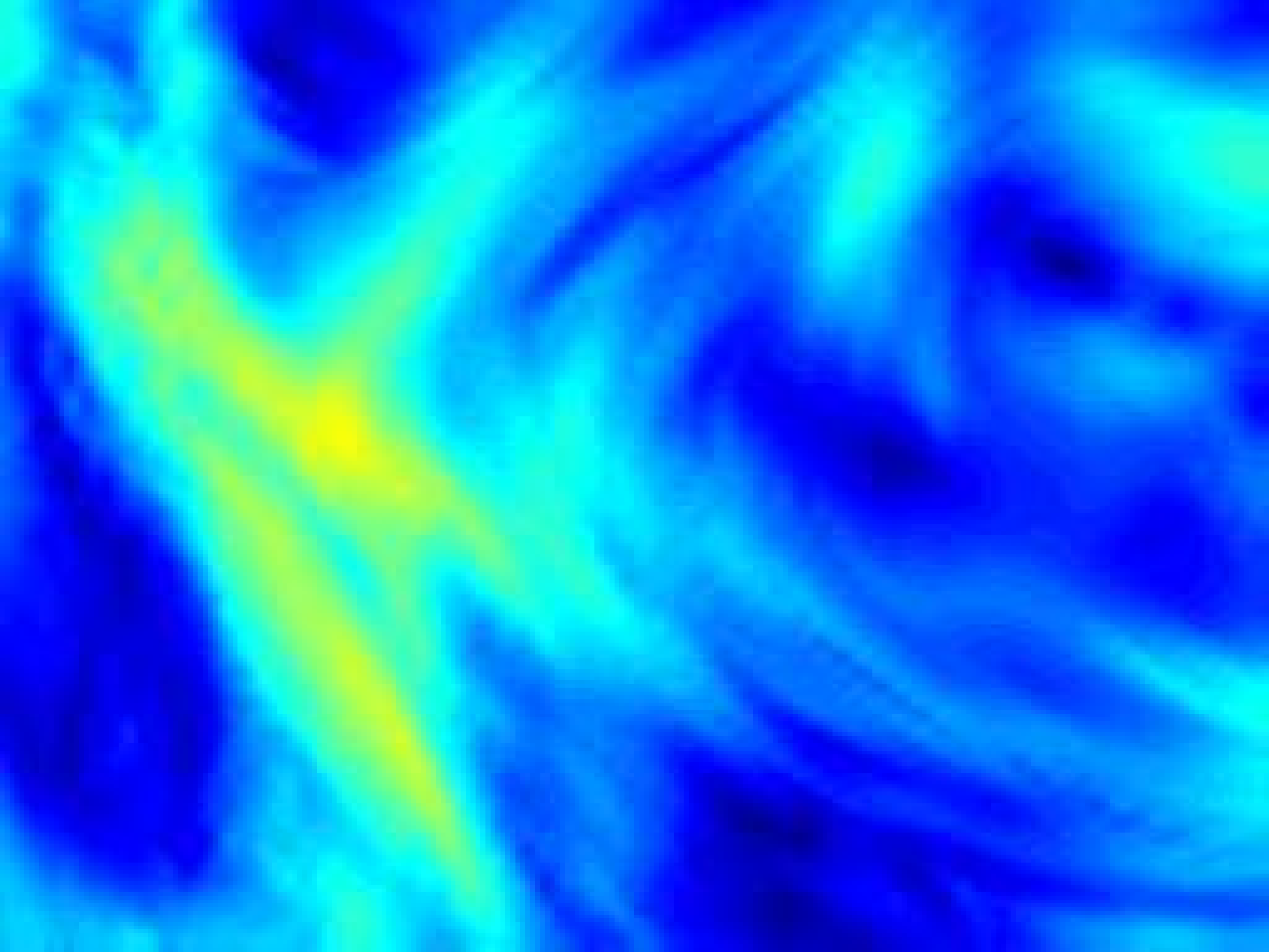}&
    \includegraphics[width=3cm]{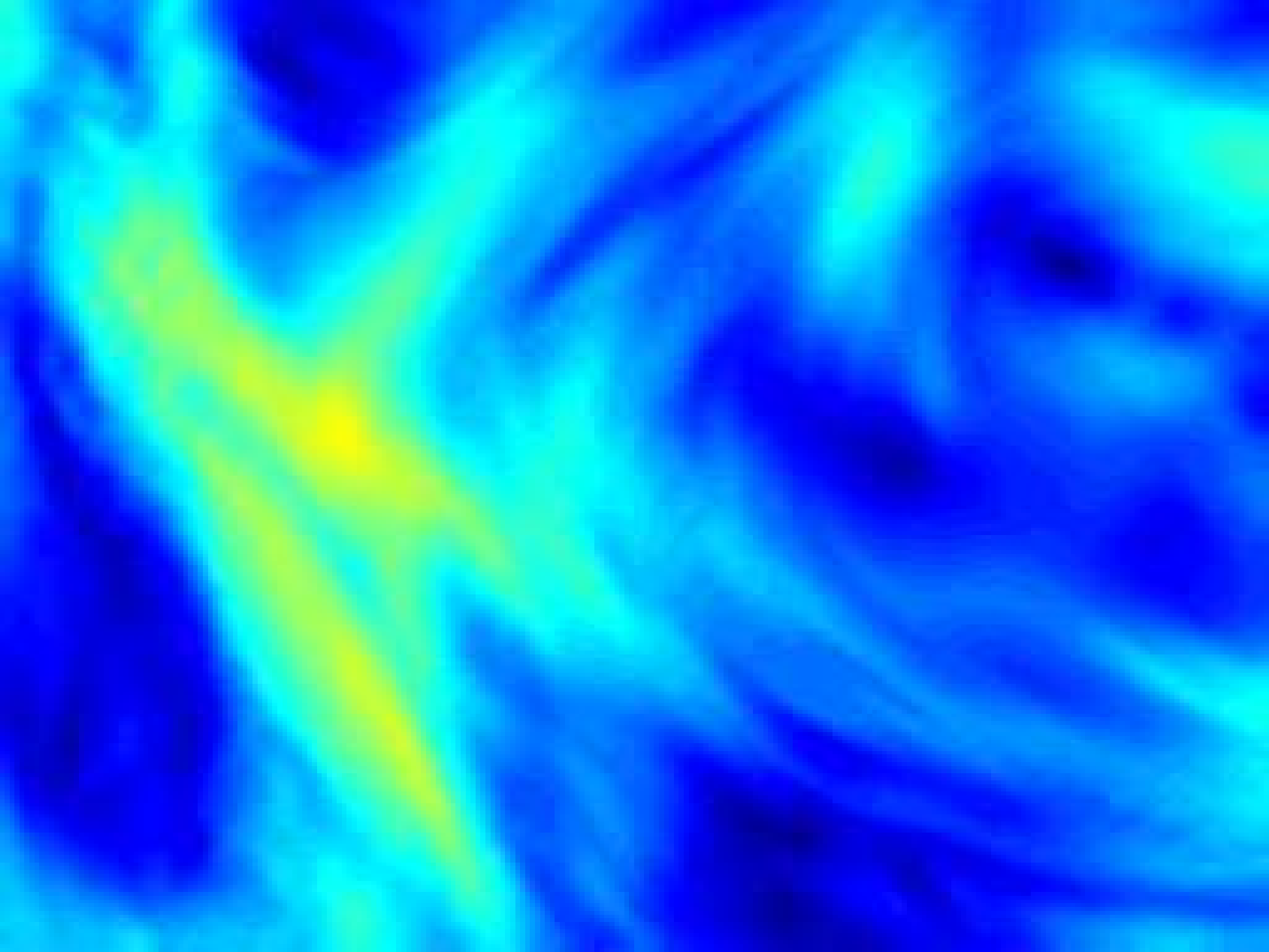}
  \end{tabular}
  \caption{Interpolation of the 3D plasma (magnetic field) data set from regular sampling with spacing $4\times 4\times 1$. Two  spatial cross  sections of the original data are shown in the first figures on the first and third row. The results of cubic spline, DCT, DFT, wavelet, and LDMM are shown in the remaining five figures.}
  \label{fig:down_plasma_3d_441}
\end{figure}

\begin{figure}[H]
  \centering
  \begin{tabular}{ccc}
    Original & Cubic Spline (22.93dB)& DCT (24.54dB)\\
    \includegraphics[width=3cm]{plasma_3d_original_band_19}&
    \includegraphics[width=3cm]{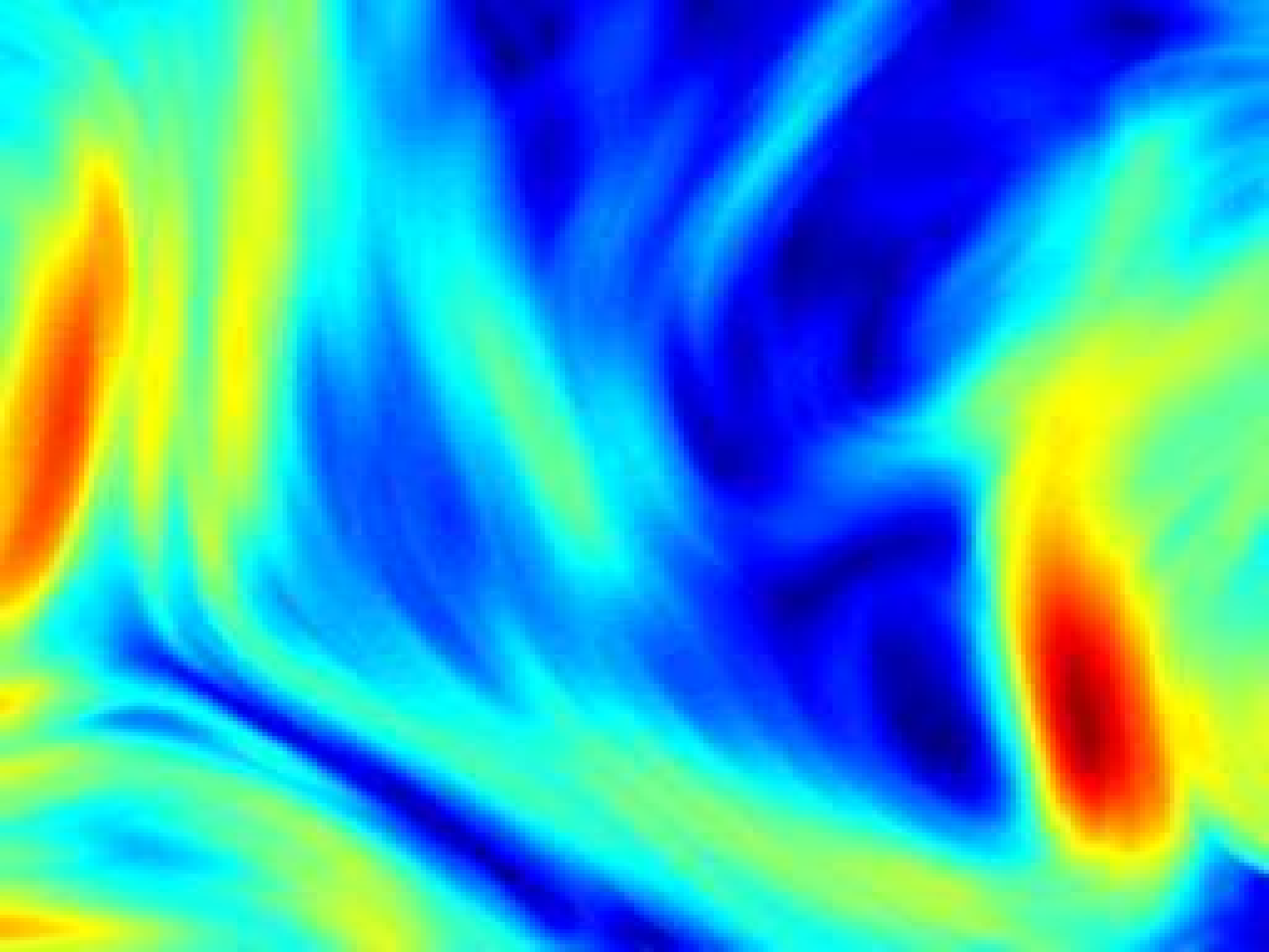}&
    \includegraphics[width=3cm]{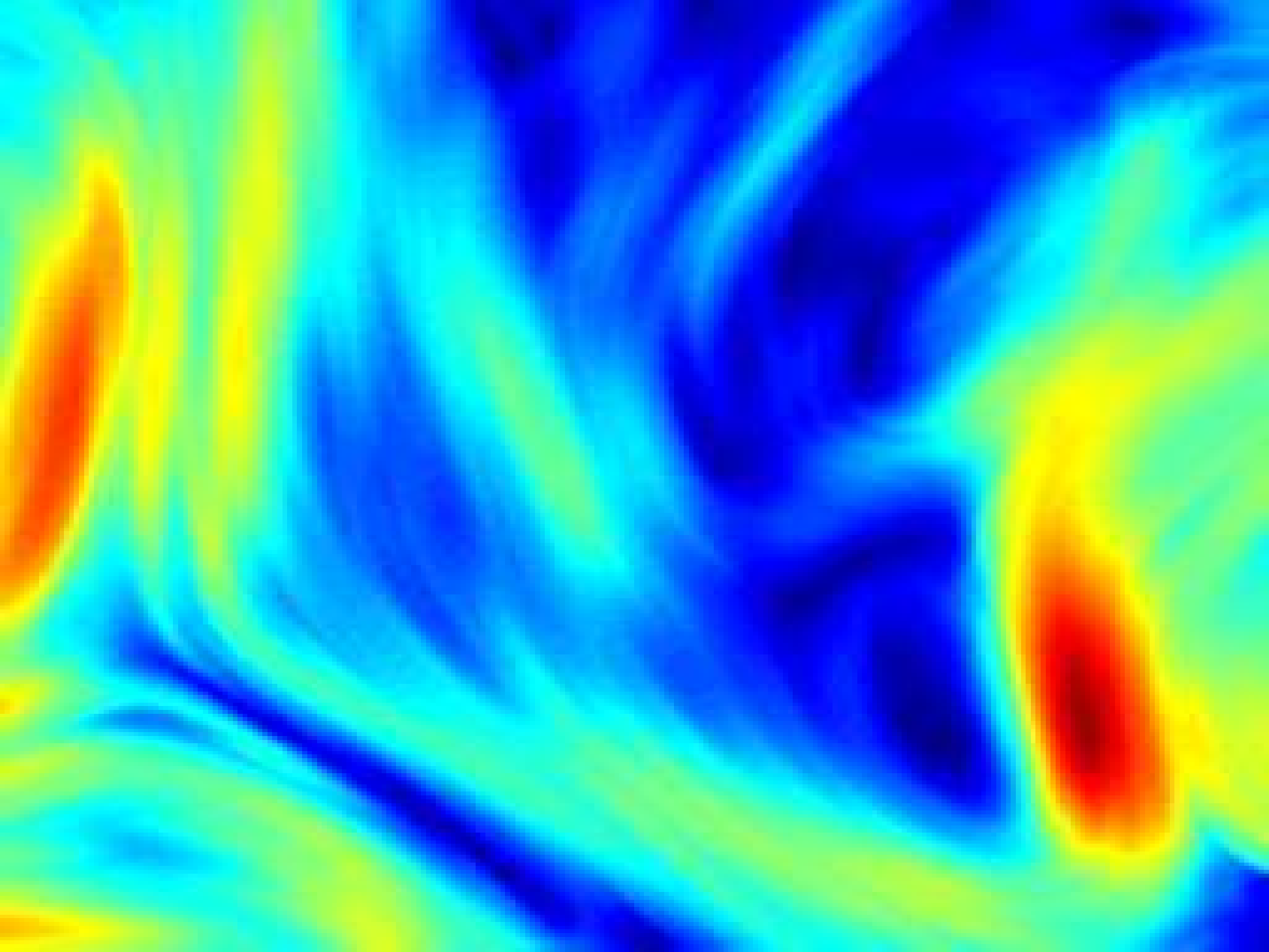}\\
    DFT (23.99dB)& Wavelet (24.25dB)& LDMM (\textbf{25.43dB})\\
    \includegraphics[width=3cm]{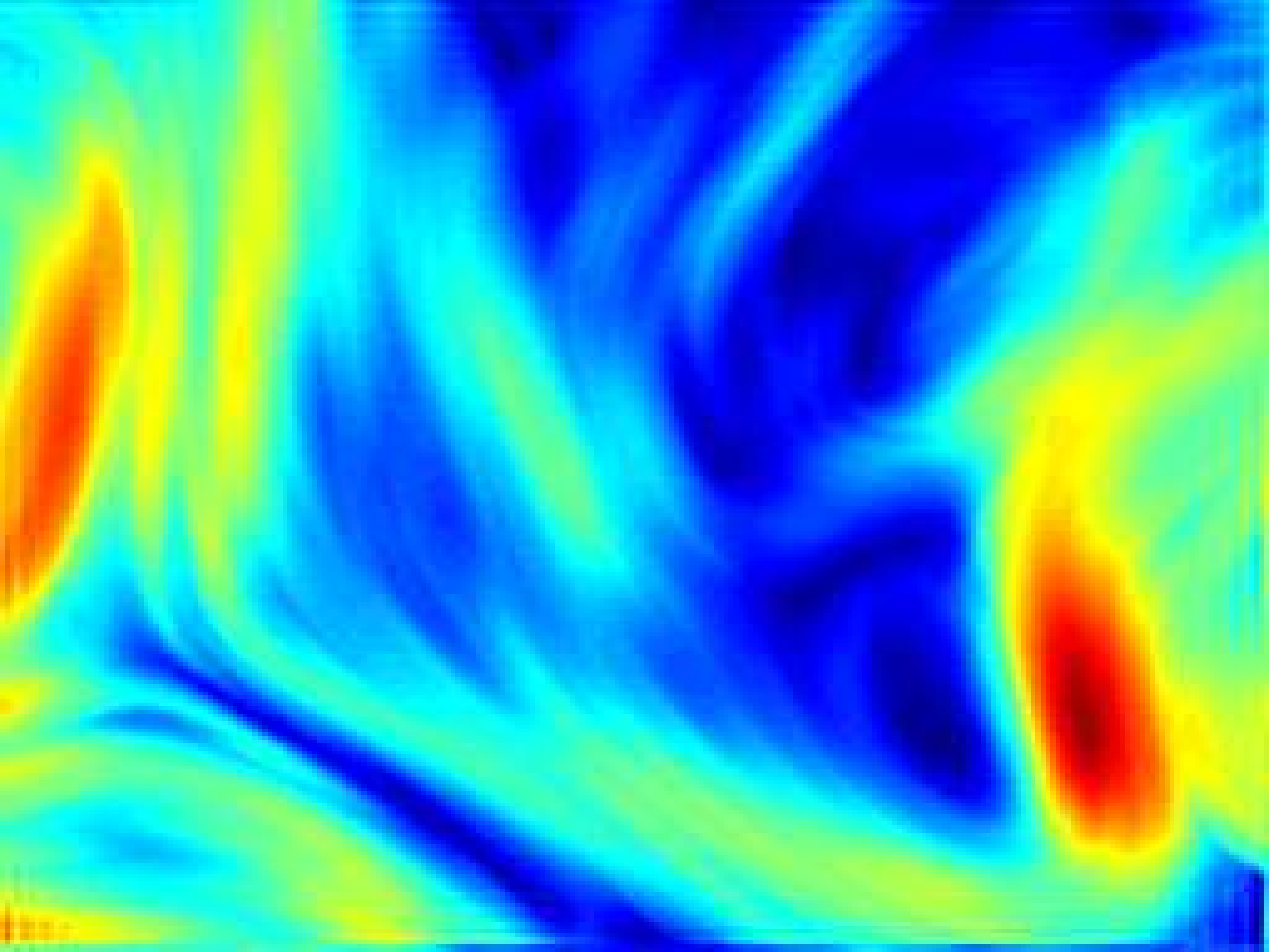}&
    \includegraphics[width=3cm]{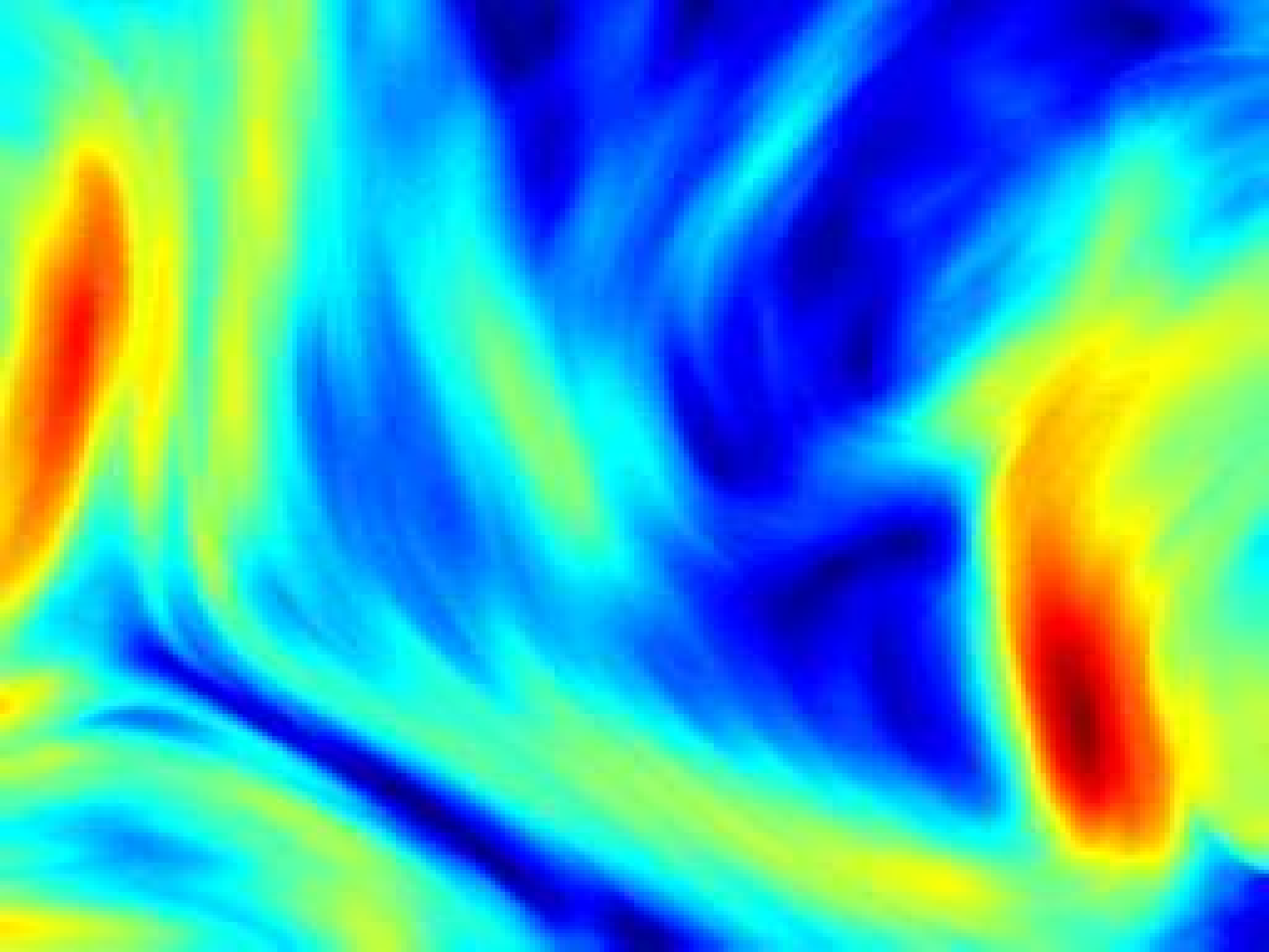}&
    \includegraphics[width=3cm]{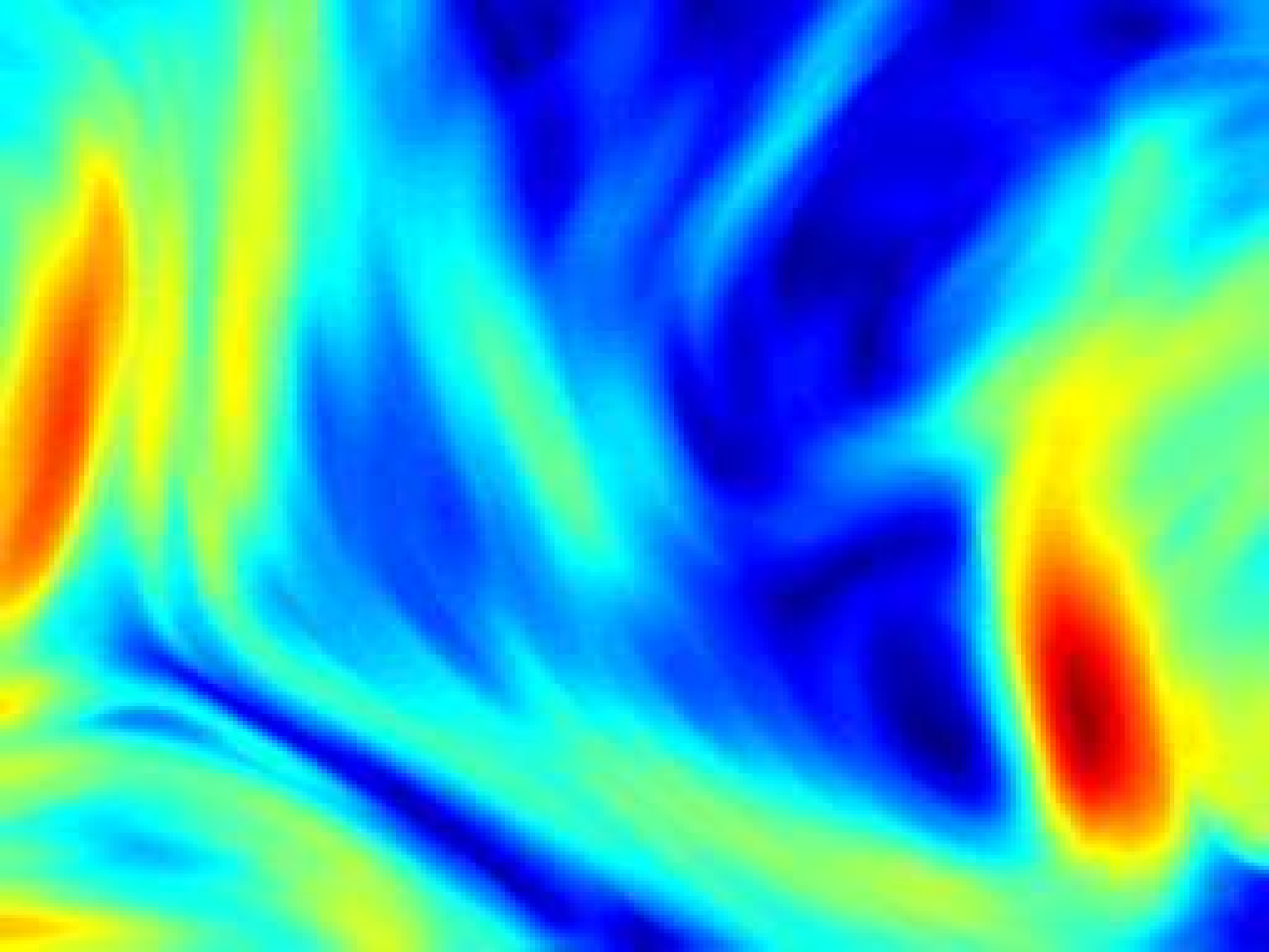}\\
    Original & Cubic Spline (22.93dB)& DCT (24.54dB)\\
    \includegraphics[width=3cm]{plasma_3d_original_band_29}&
    \includegraphics[width=3cm]{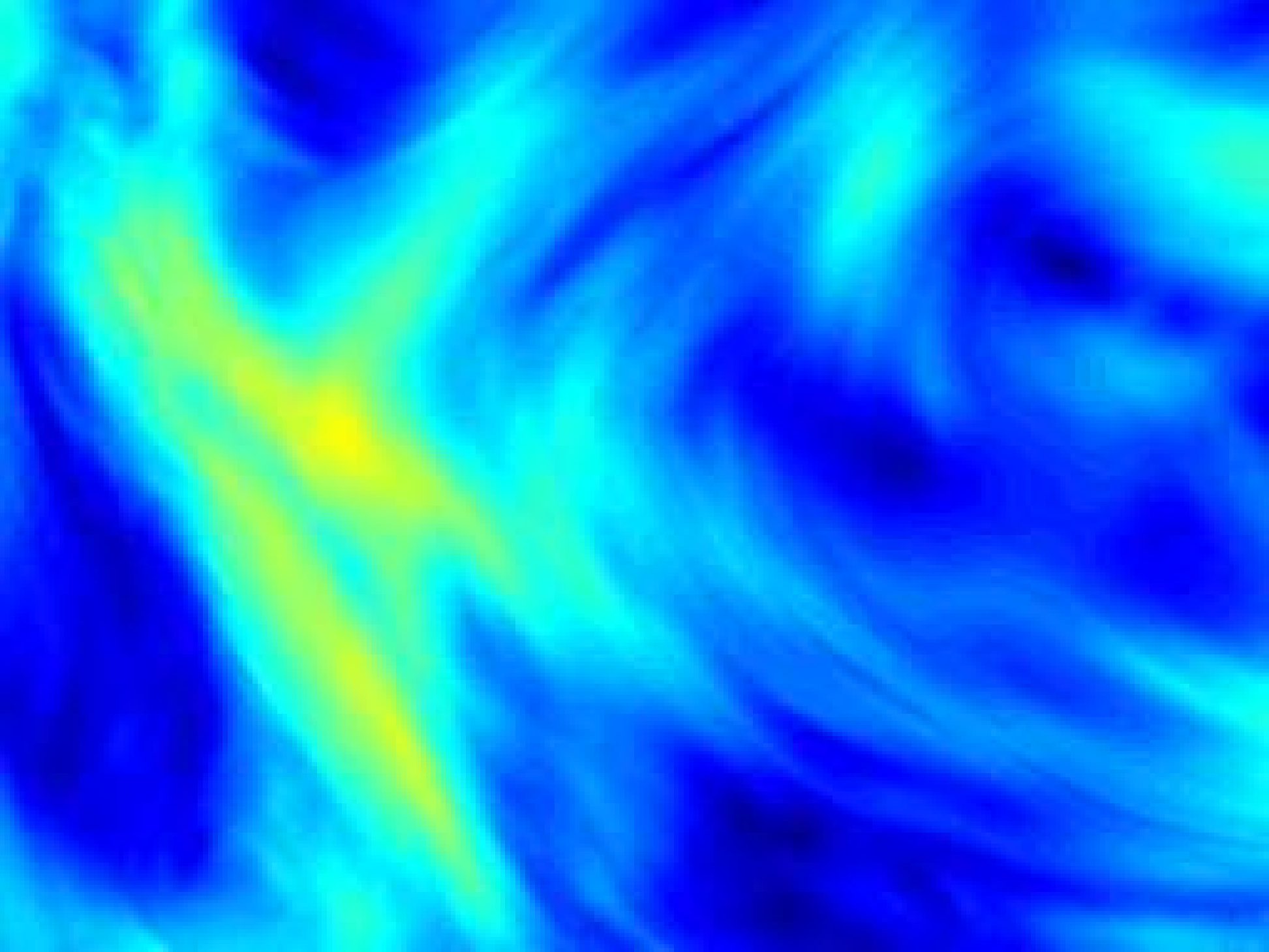}&
    \includegraphics[width=3cm]{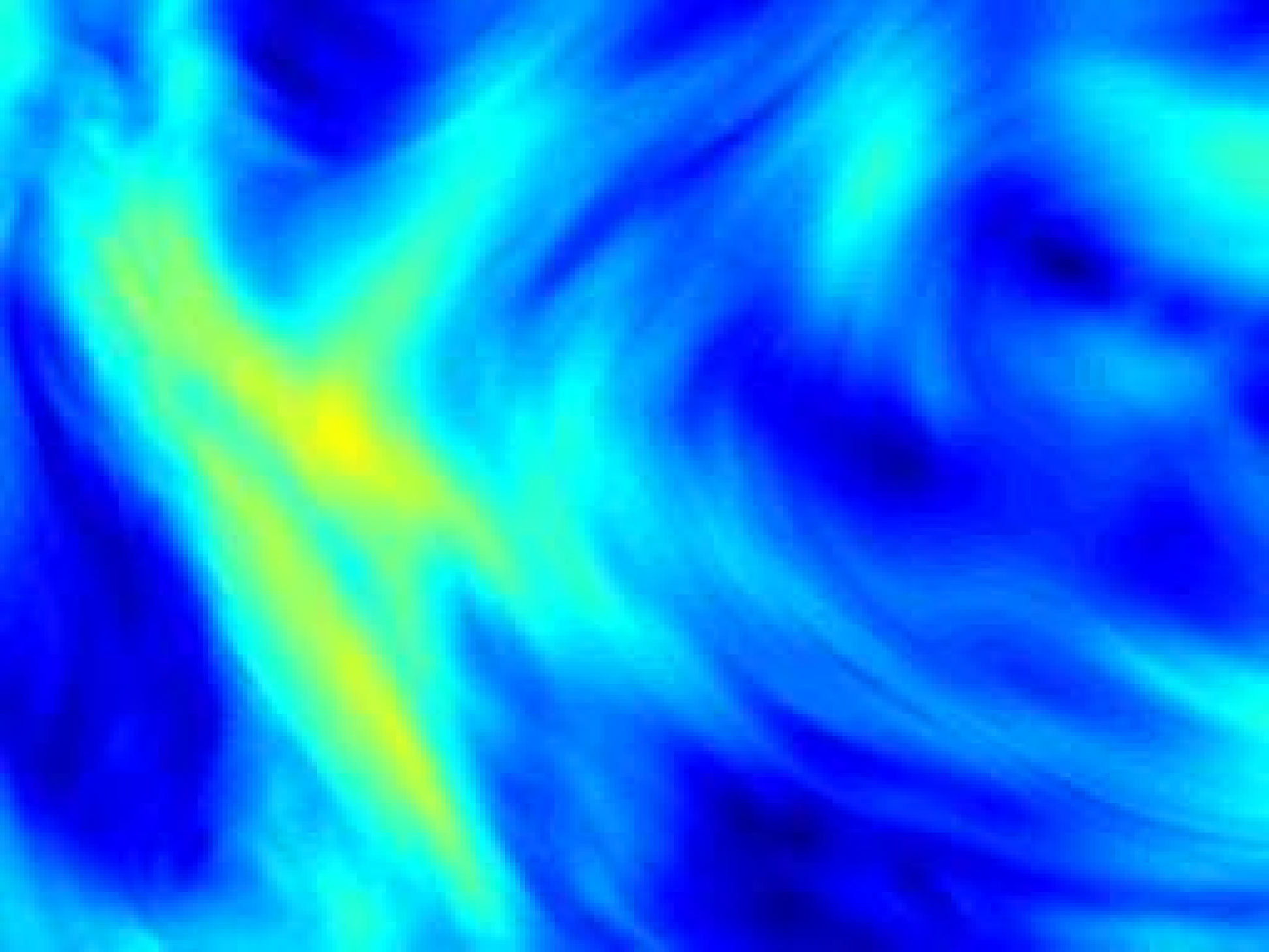}\\
    DFT (23.99dB)& Wavelet (24.25dB)& LDMM (\textbf{25.43dB})\\
    \includegraphics[width=3cm]{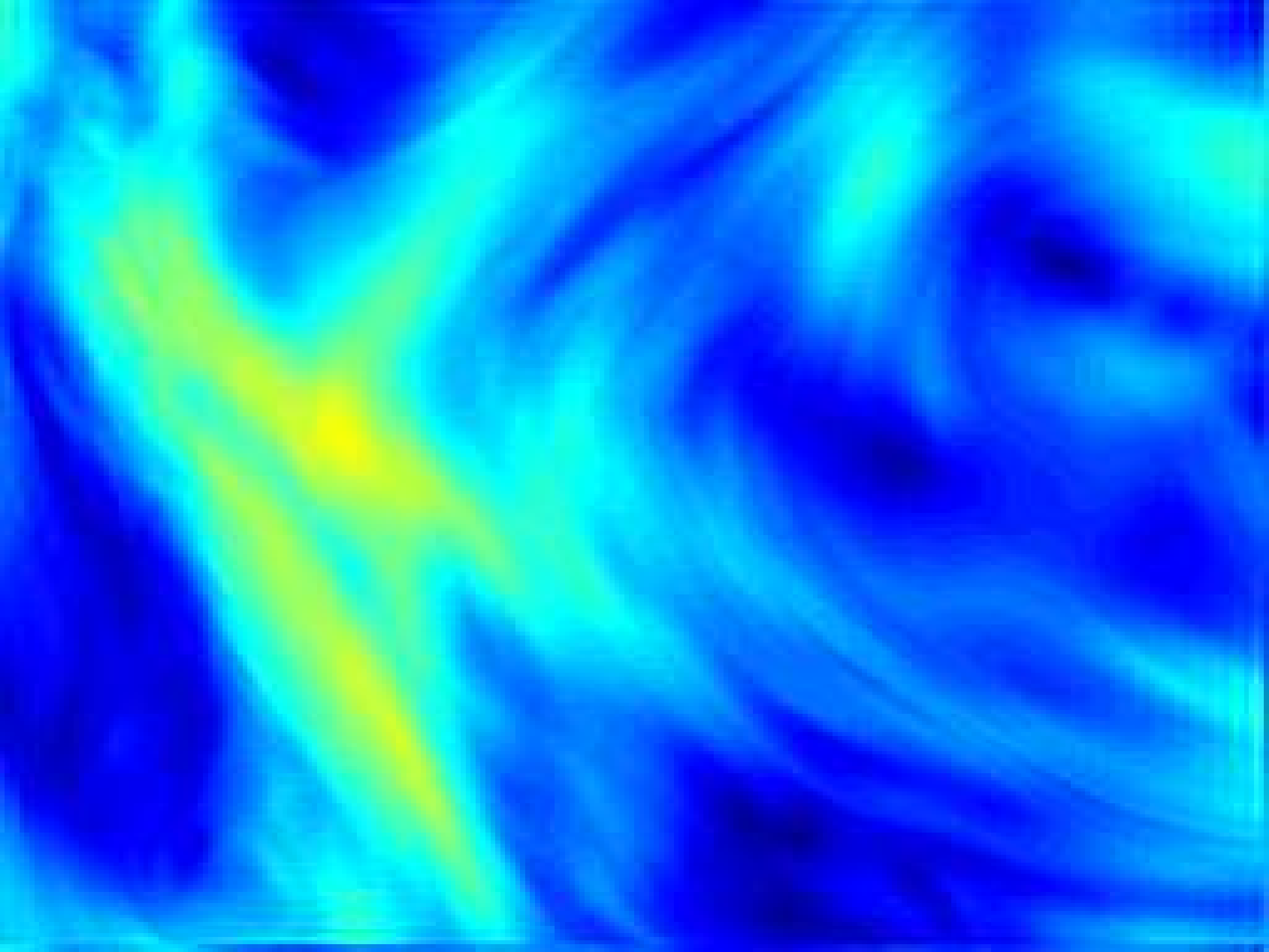}&
    \includegraphics[width=3cm]{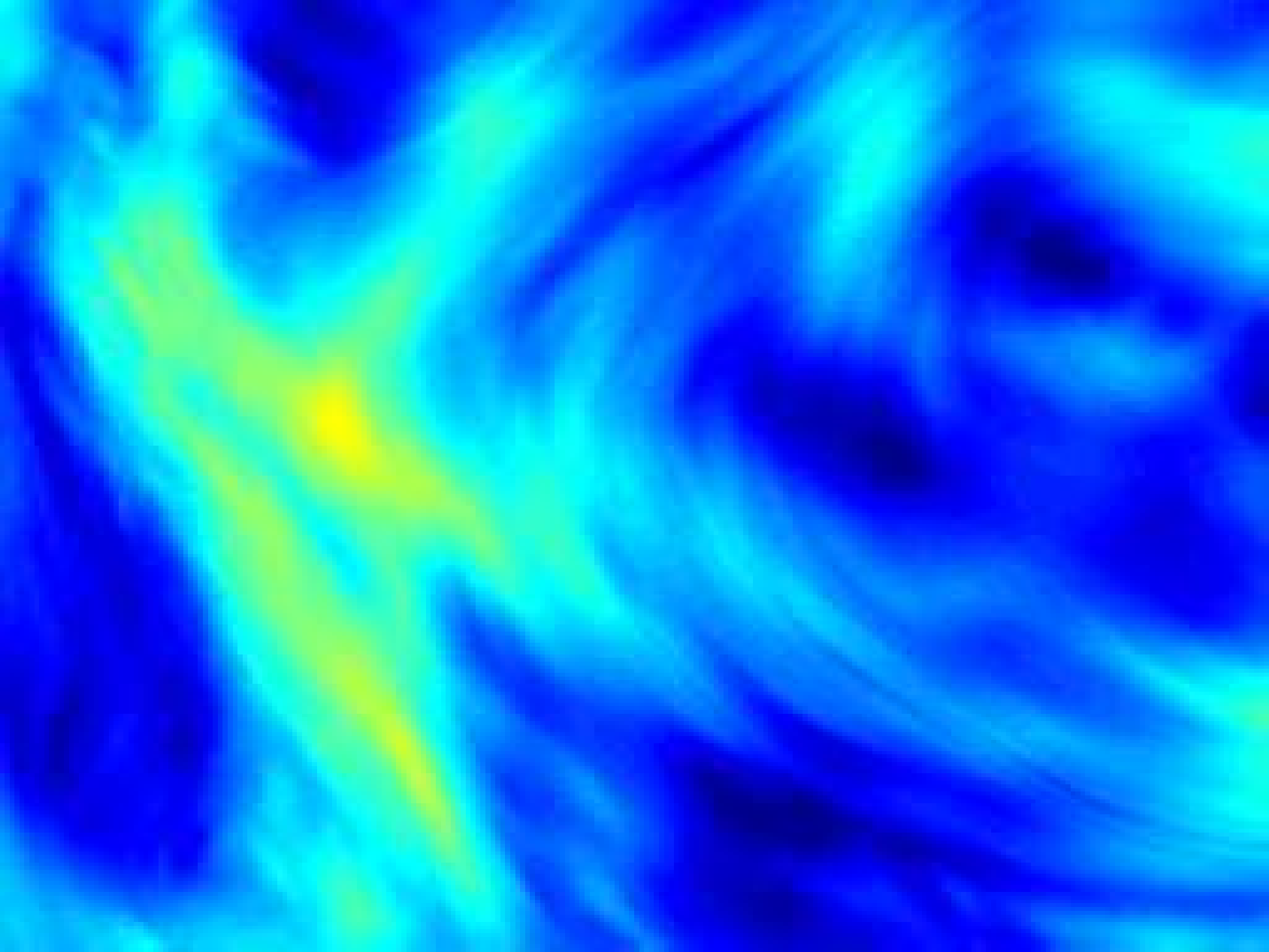}&
    \includegraphics[width=3cm]{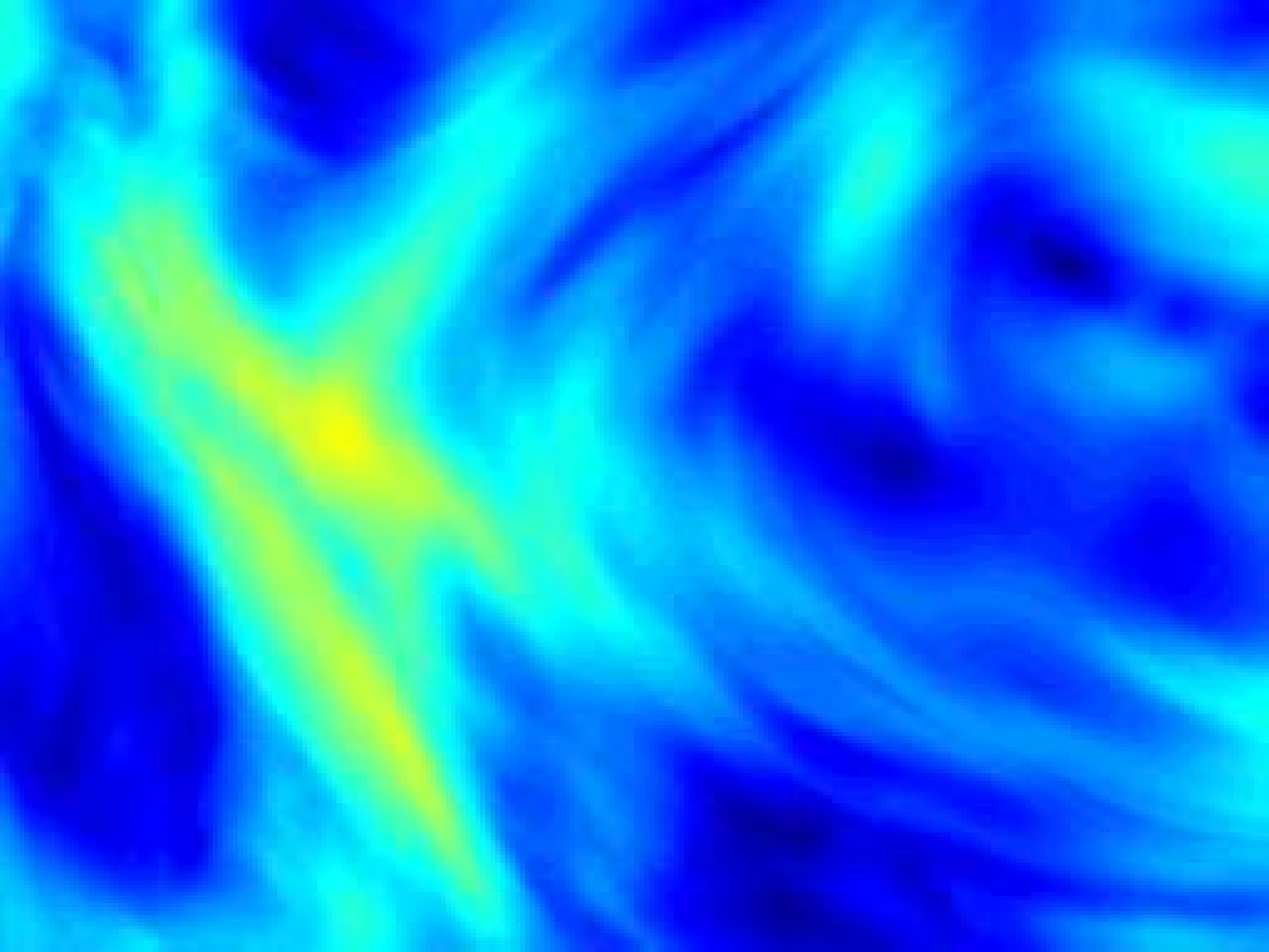}
  \end{tabular}
  \caption{Interpolation of the 3D plasma (magnetic field) data set from regular sampling with spacing $2\times 2\times 2$. Two spatial cross sections of the original data are shown in the first figures on the first and third row. The results of cubic spline, DCT, DFT, wavelet, and LDMM are shown in the remaining five figures.}
  \label{fig:down_plasma_3d_222}
\end{figure}

\begin{figure}[H]
  \centering
  \begin{tabular}{ccc}
    Original & Cubic Spline (24.54dB)& DCT (30.69dB)\\
    \includegraphics[width=3cm]{lattice_3d_original_band_31}&
    \includegraphics[width=3cm]{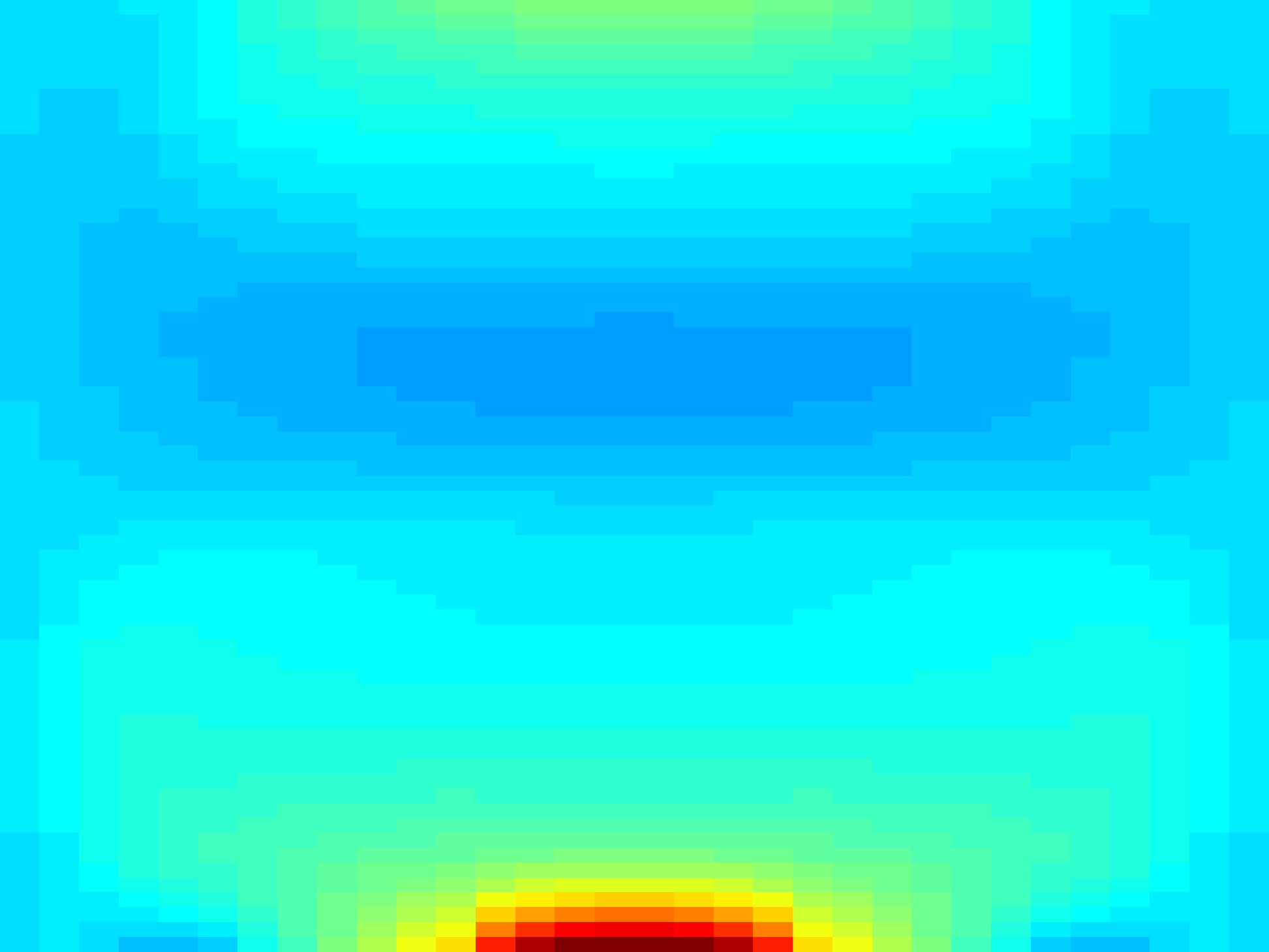}&
    \includegraphics[width=3cm]{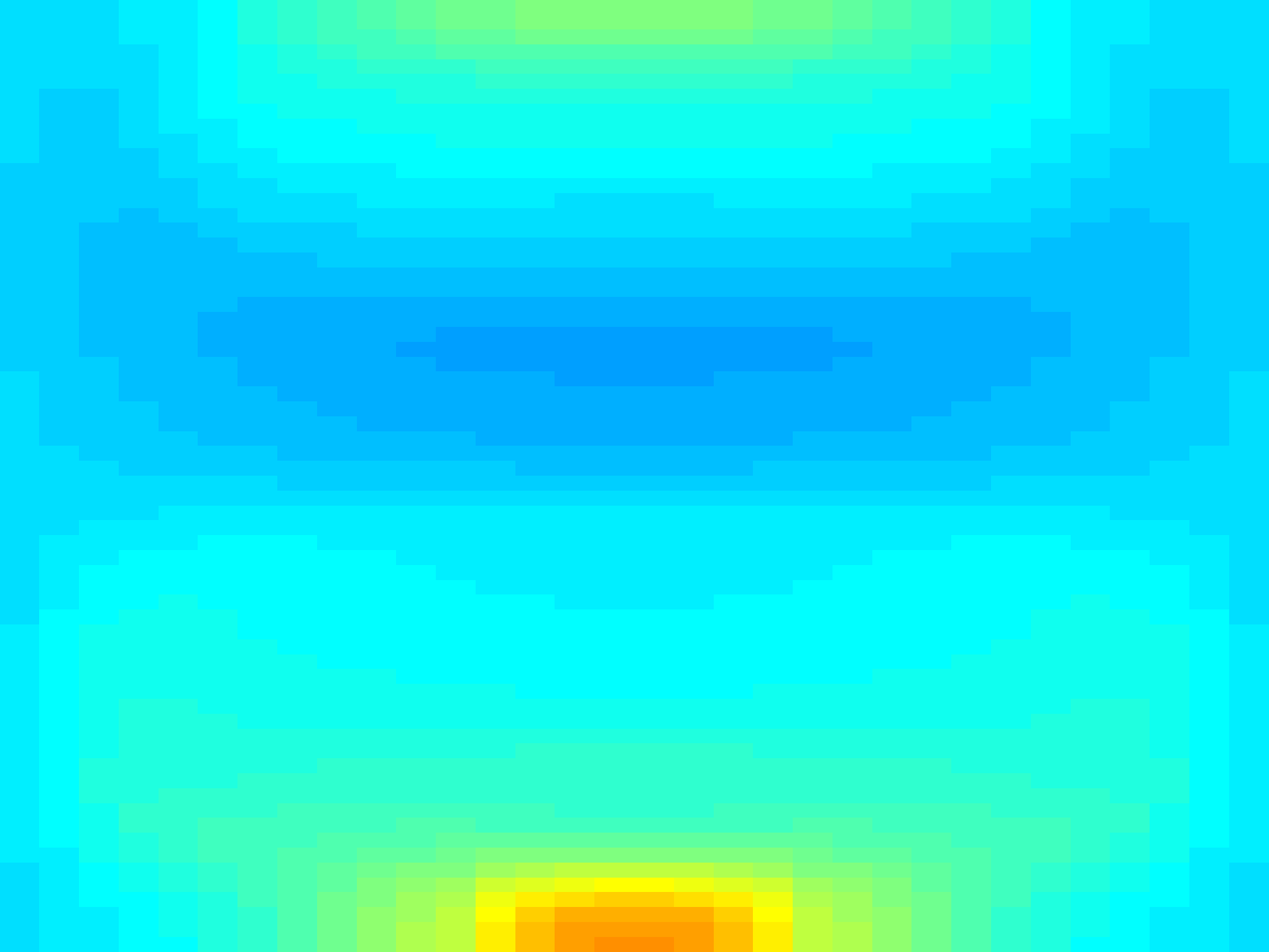}\\
    DFT (27.25dB)& Wavelet (31.03dB)& LDMM (\textbf{32.64dB})\\
    \includegraphics[width=3cm]{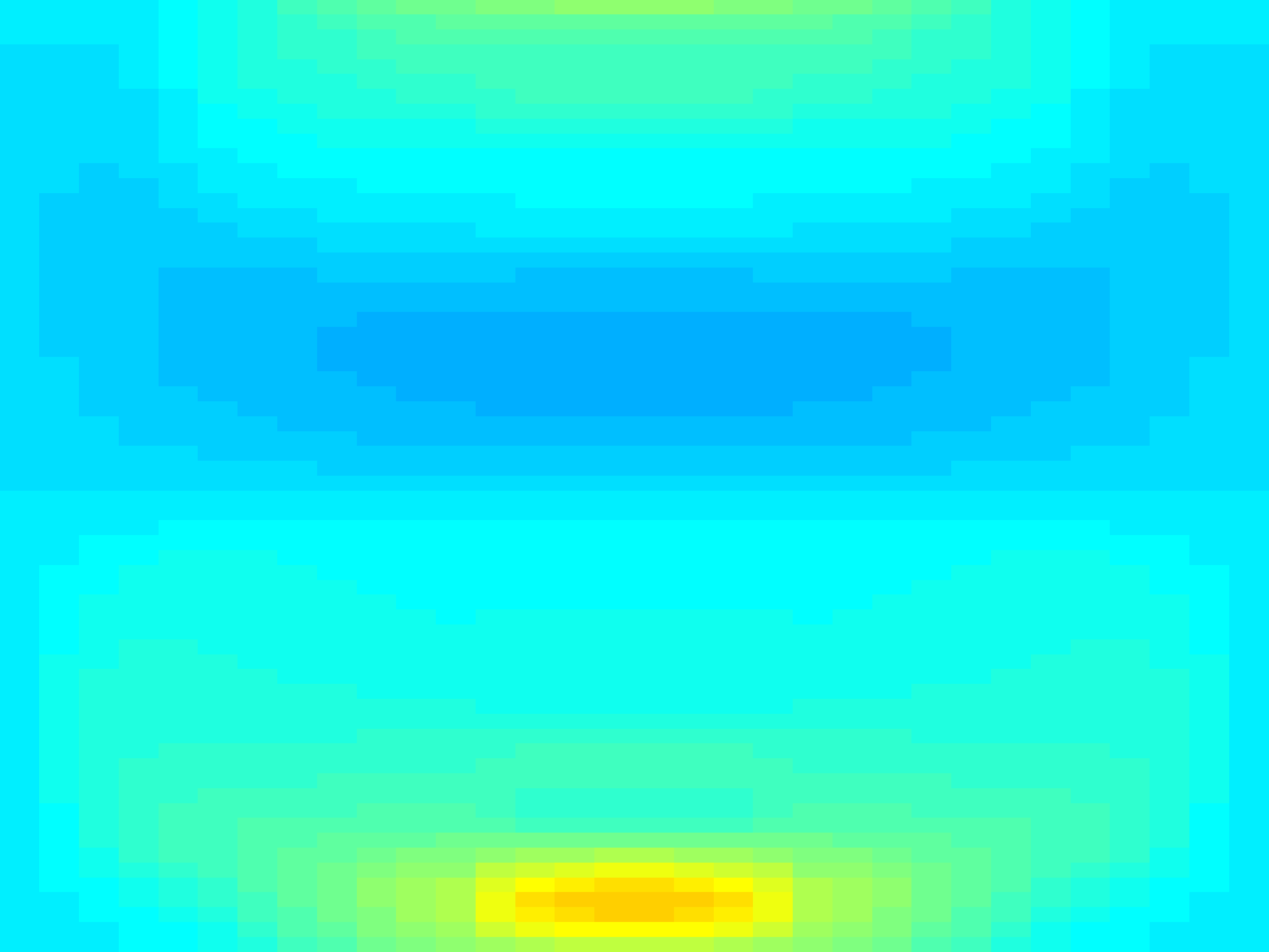}&
    \includegraphics[width=3cm]{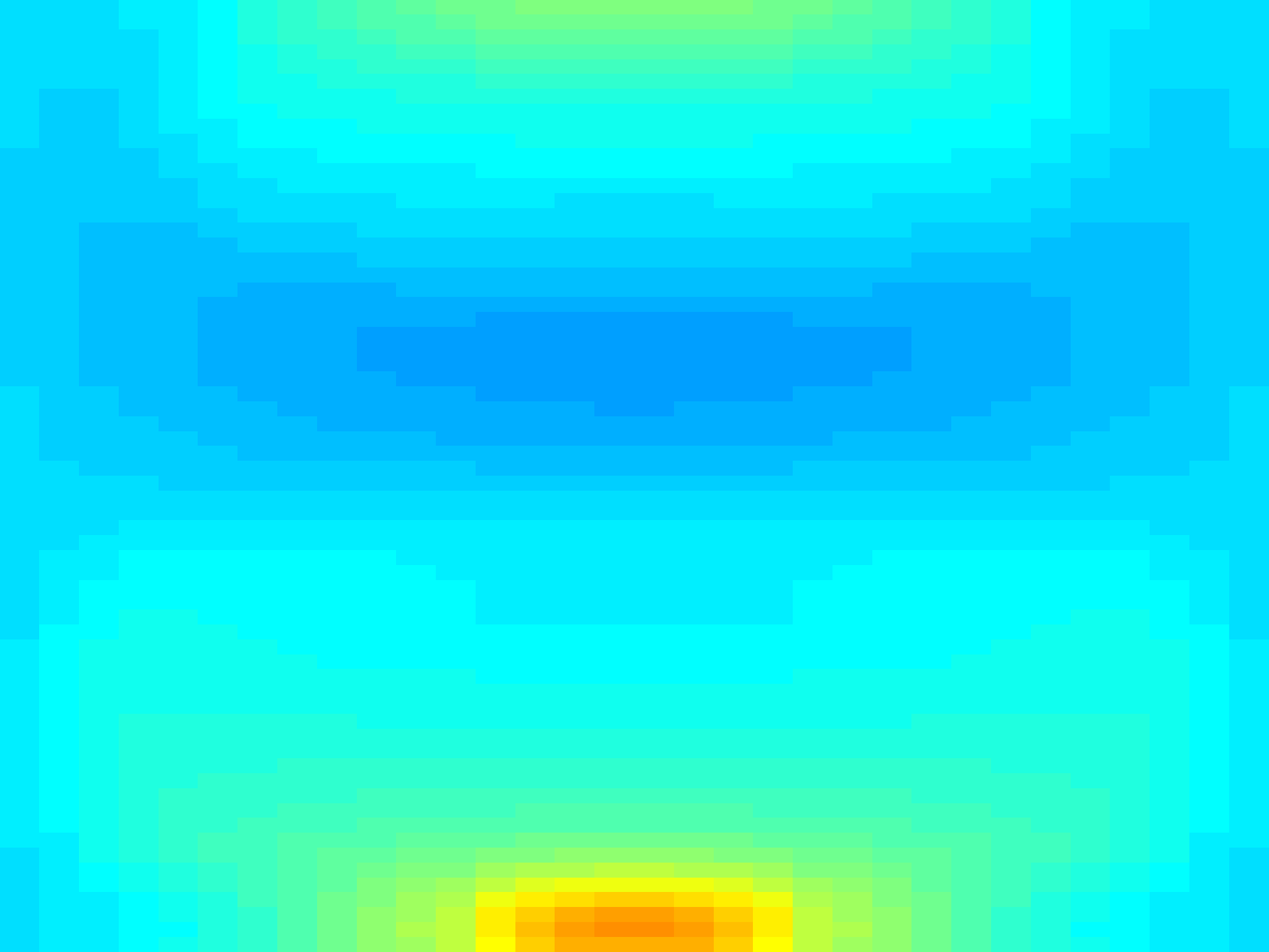}&
    \includegraphics[width=3cm]{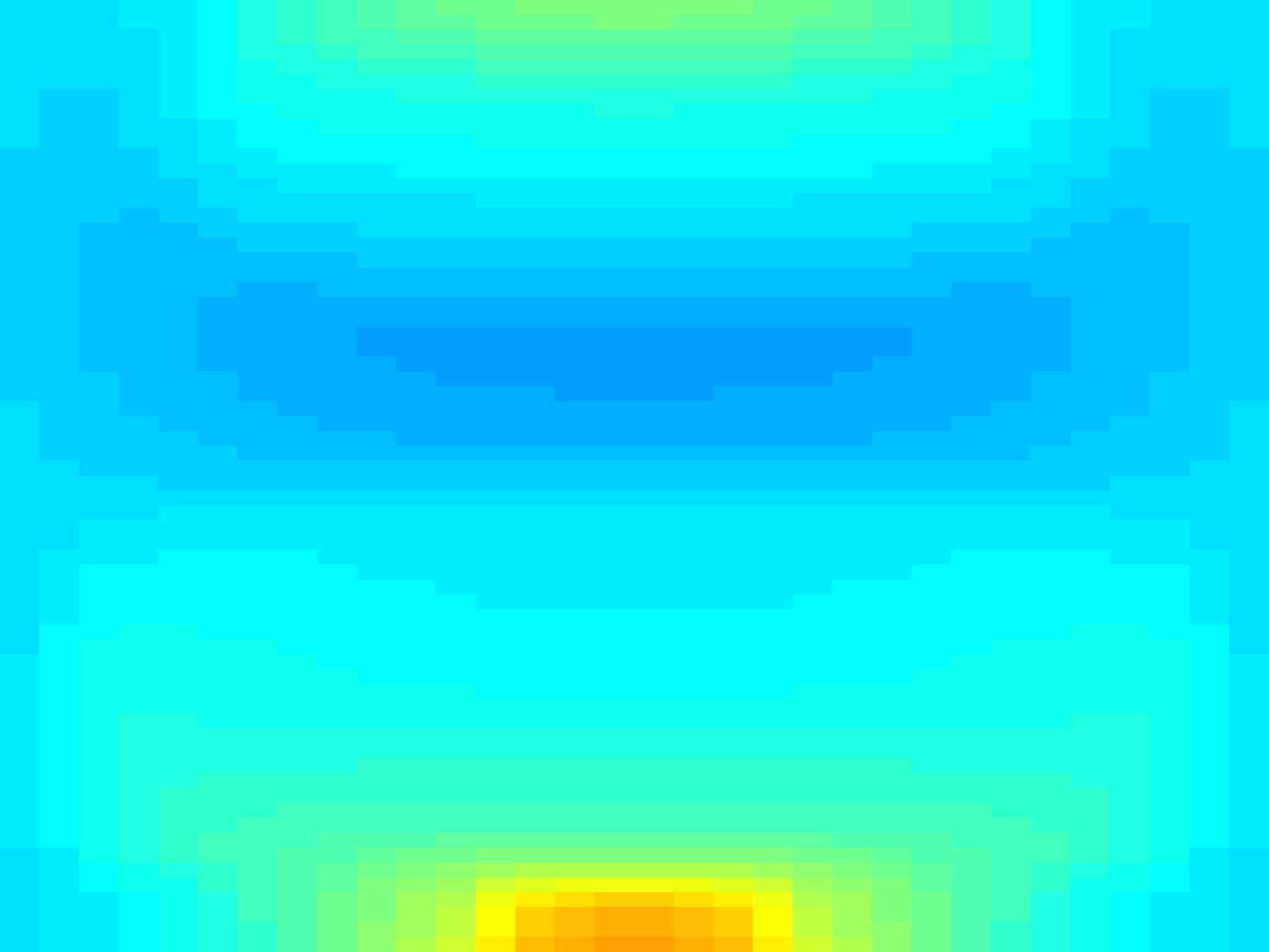}\\
    Original & Cubic Spline (24.54dB)& DCT (30.69dB)\\
    \includegraphics[width=3cm]{lattice_3d_original_band_151}&
    \includegraphics[width=3cm]{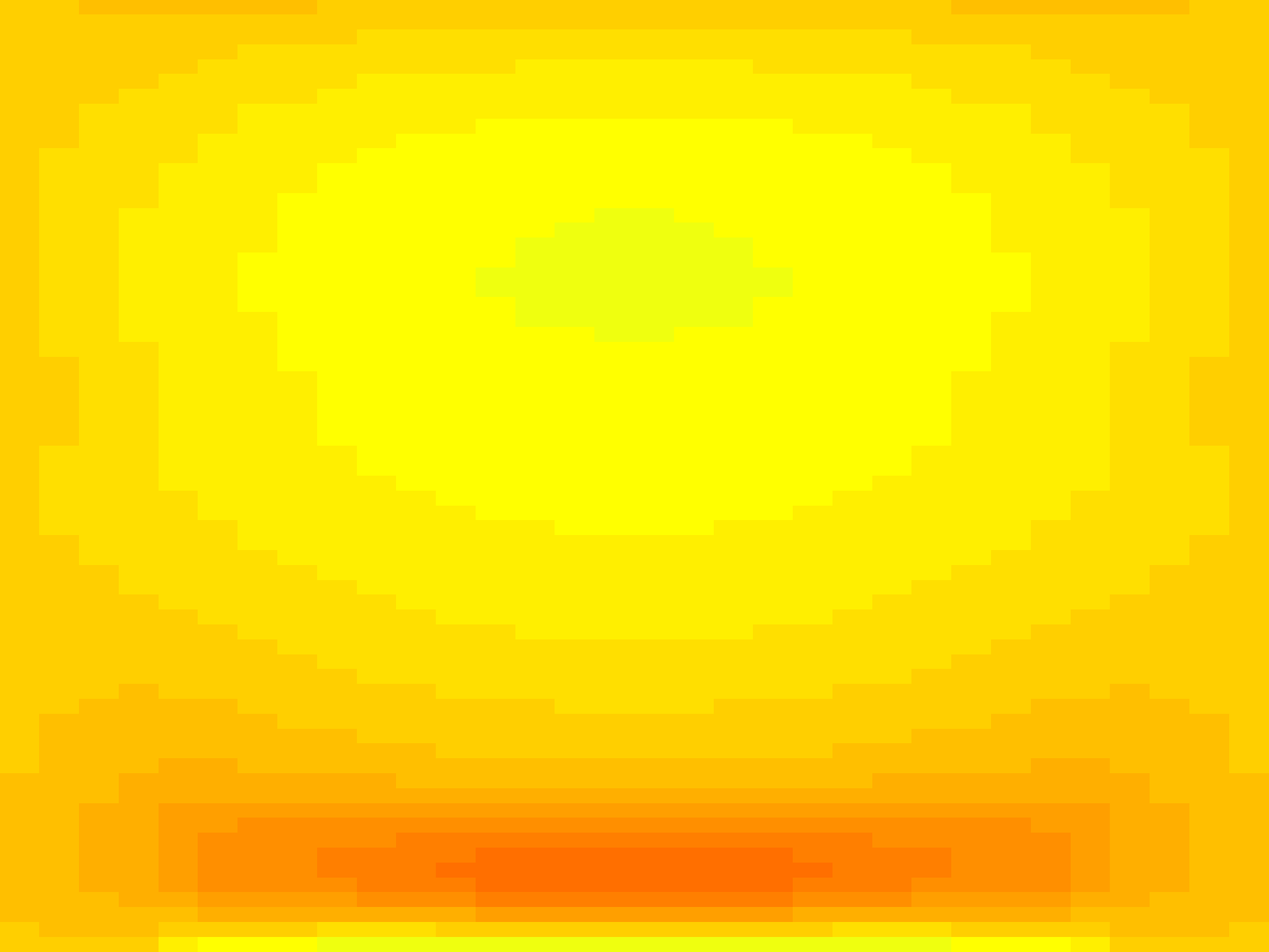}&
    \includegraphics[width=3cm]{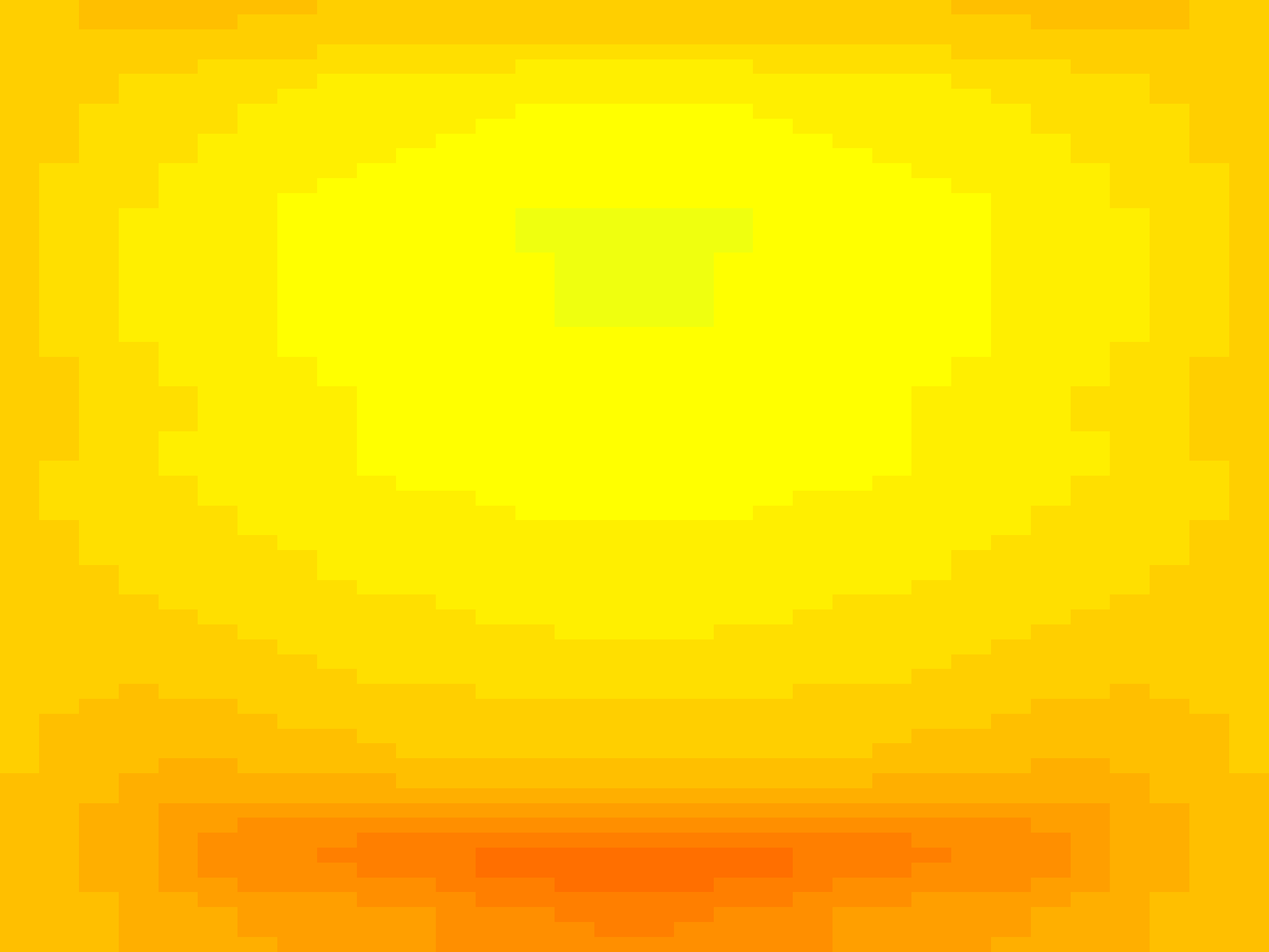}\\
    DFT (27.25dB)& Wavelet (31.03dB)& LDMM (\textbf{32.64dB})\\
    \includegraphics[width=3cm]{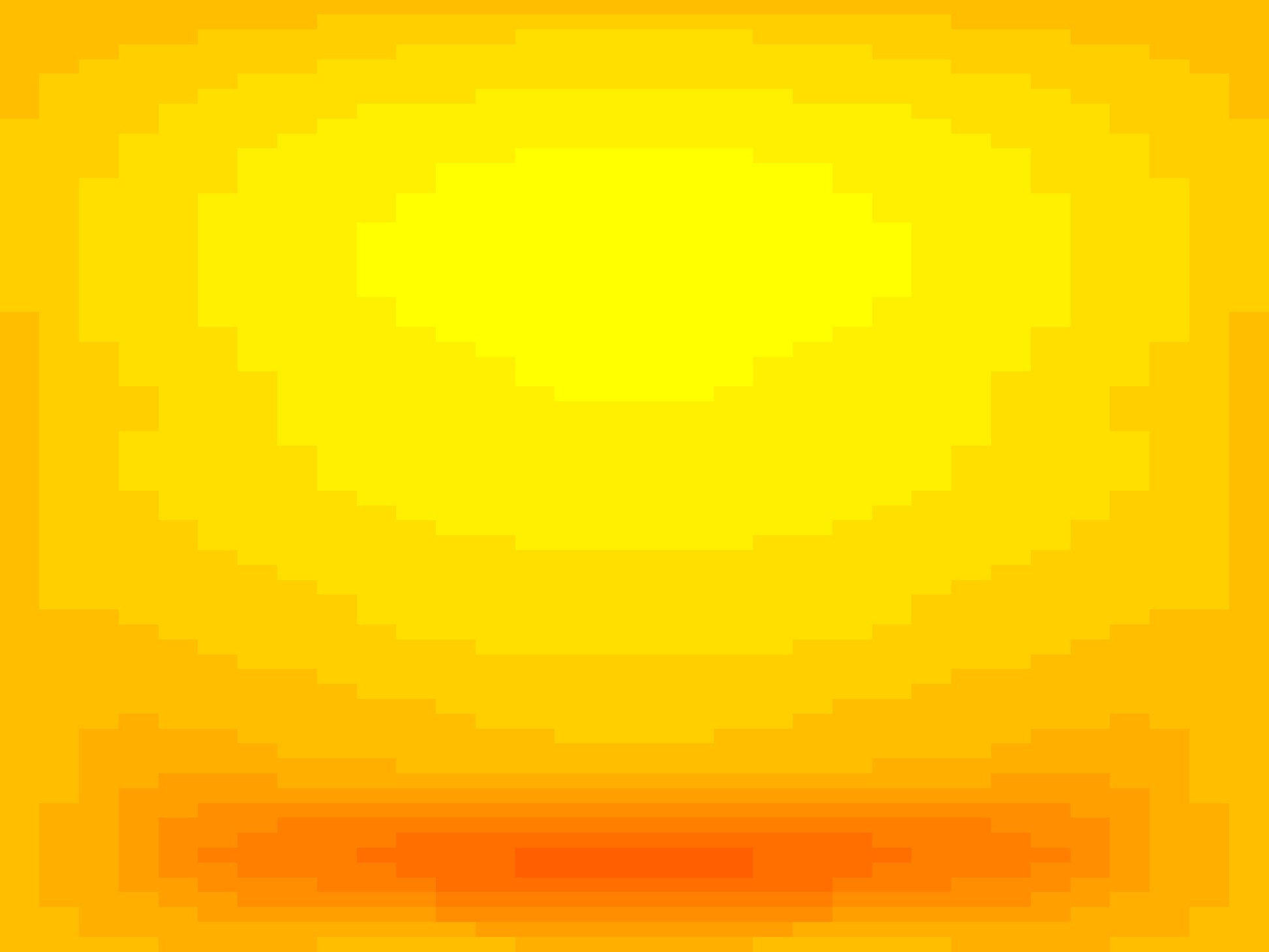}&
    \includegraphics[width=3cm]{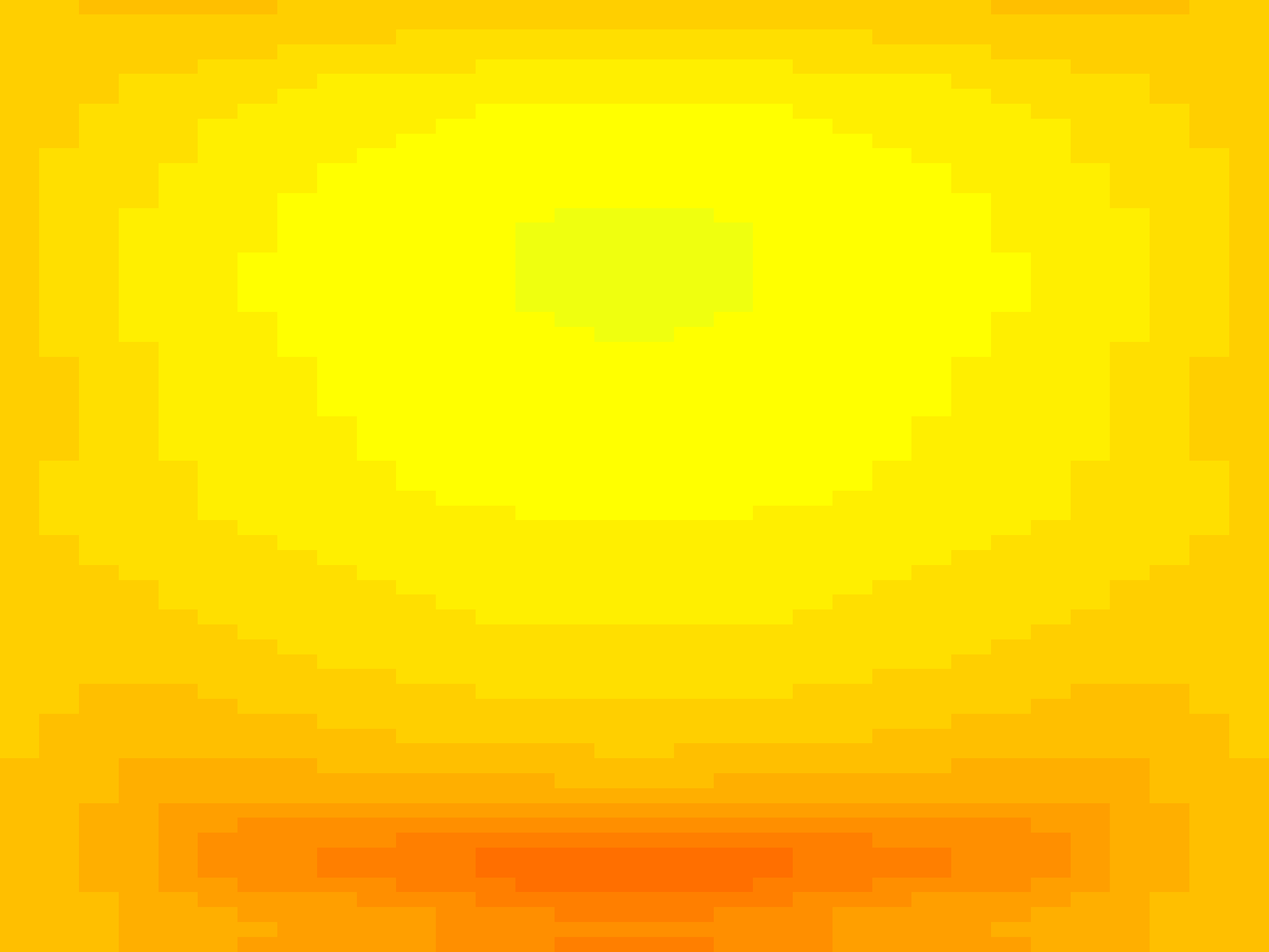}&
    \includegraphics[width=3cm]{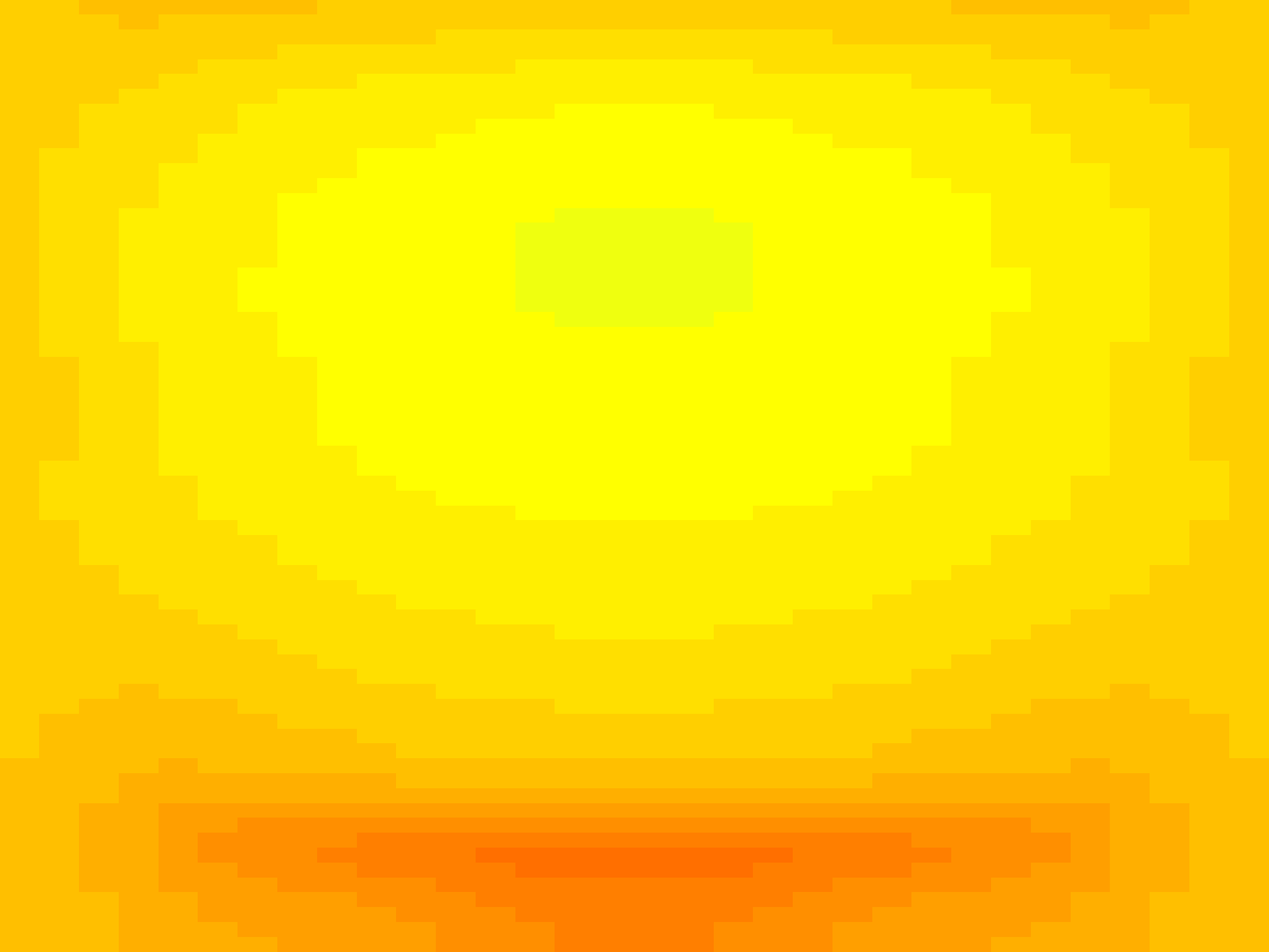}
  \end{tabular}
  \caption{Interpolation of the 3D lattice data set from regular sampling with spacing $4\times 4\times 1$. The original angular flux at $x = 0.24$ and $x = 1.18$ are shown in the first figures on the first and third row. The results of cubic spline, DCT, DFT, wavelet, and LDMM are shown in the remaining five figures.}
  \label{fig:down_lattice_3d_441}
\end{figure}

\begin{figure}[H]
  \centering
  \begin{tabular}{ccc}
    Original & Cubic Spline (30.01dB)& DCT (38.49dB)\\
    \includegraphics[width=3cm]{lattice_3d_original_band_31}&
    \includegraphics[width=3cm]{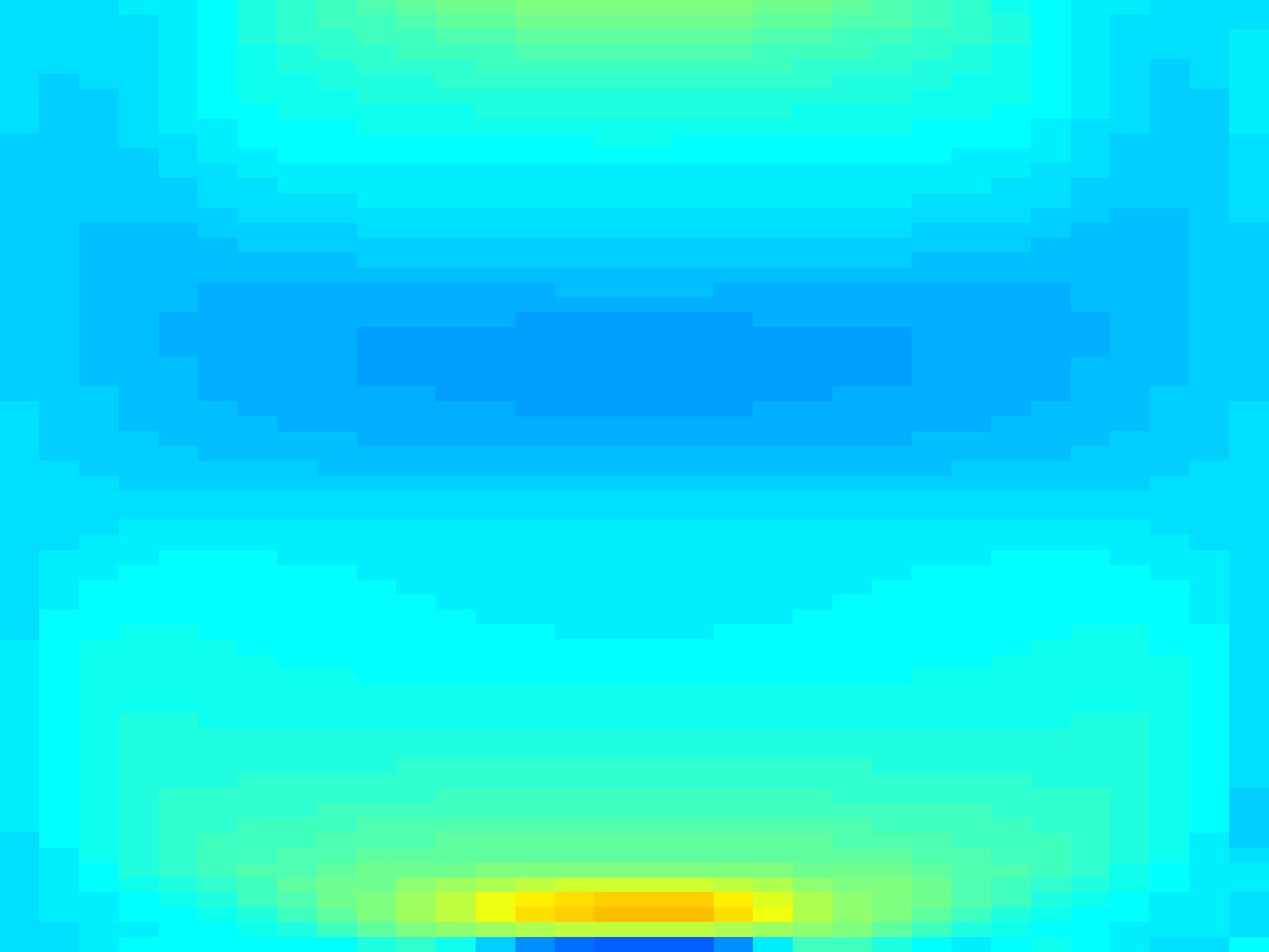}&
    \includegraphics[width=3cm]{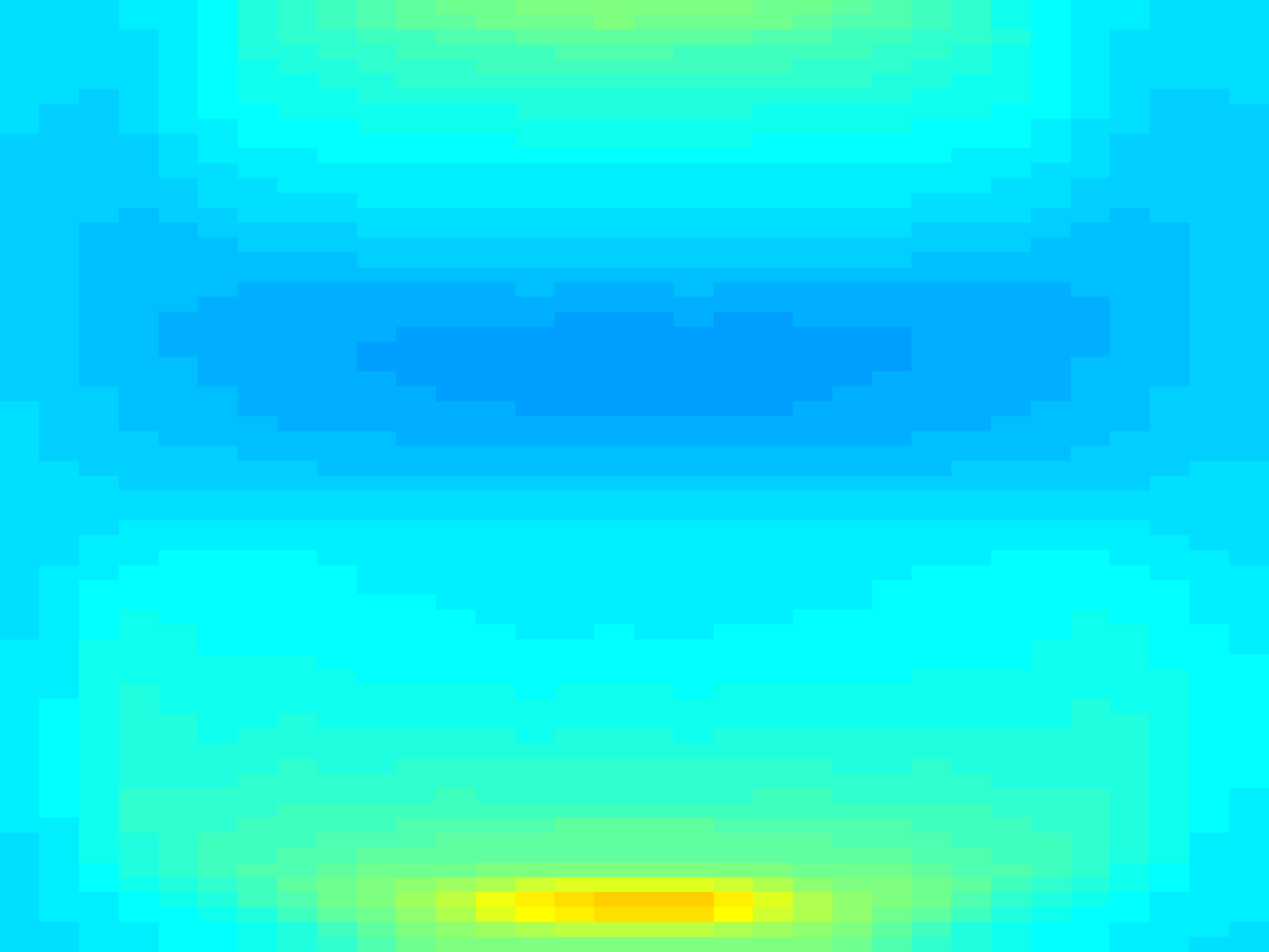}\\
    DFT (32.51dB)& Wavelet (38.15dB)& LDMM (\textbf{39.93dB})\\
    \includegraphics[width=3cm]{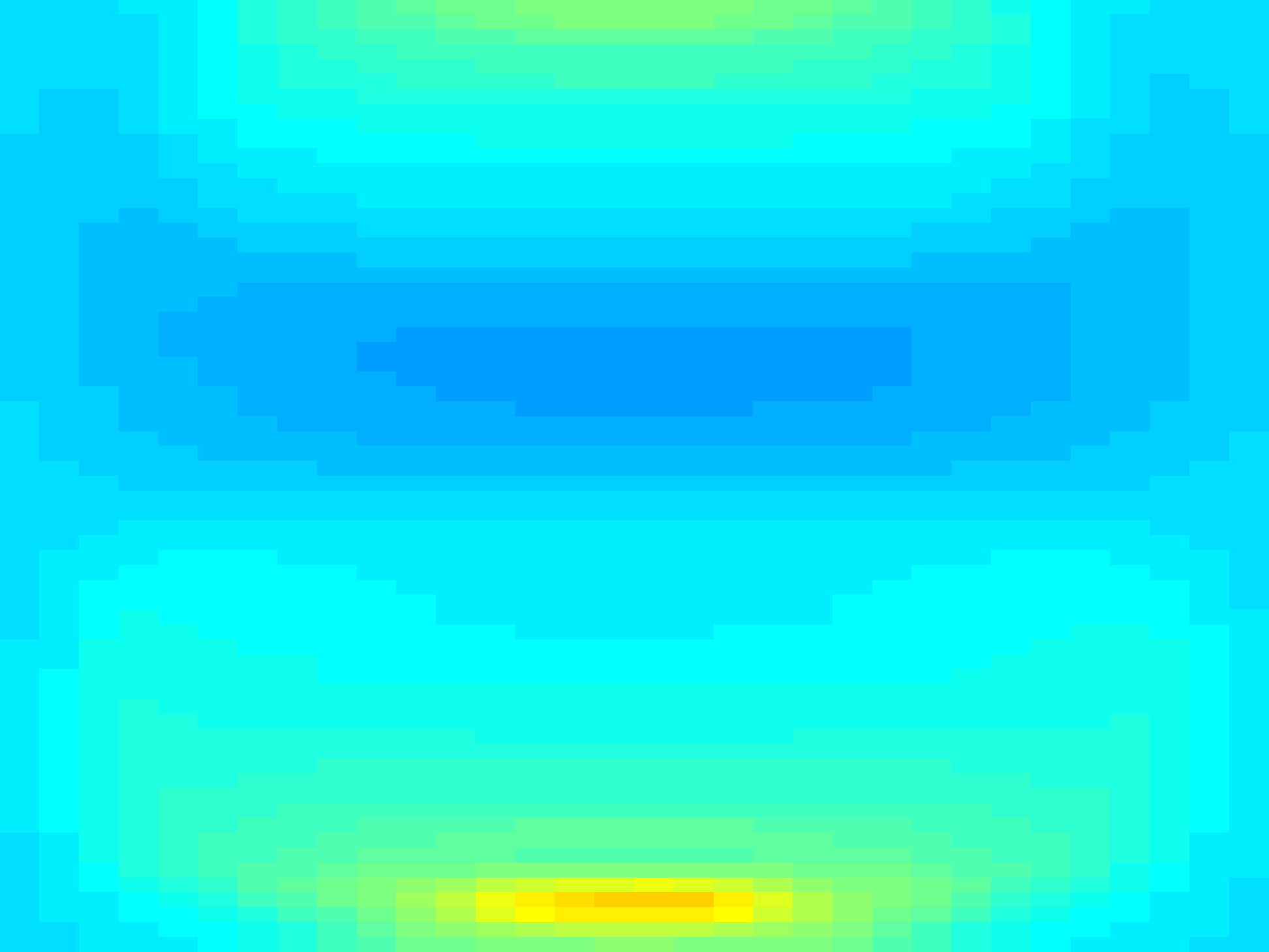}&
    \includegraphics[width=3cm]{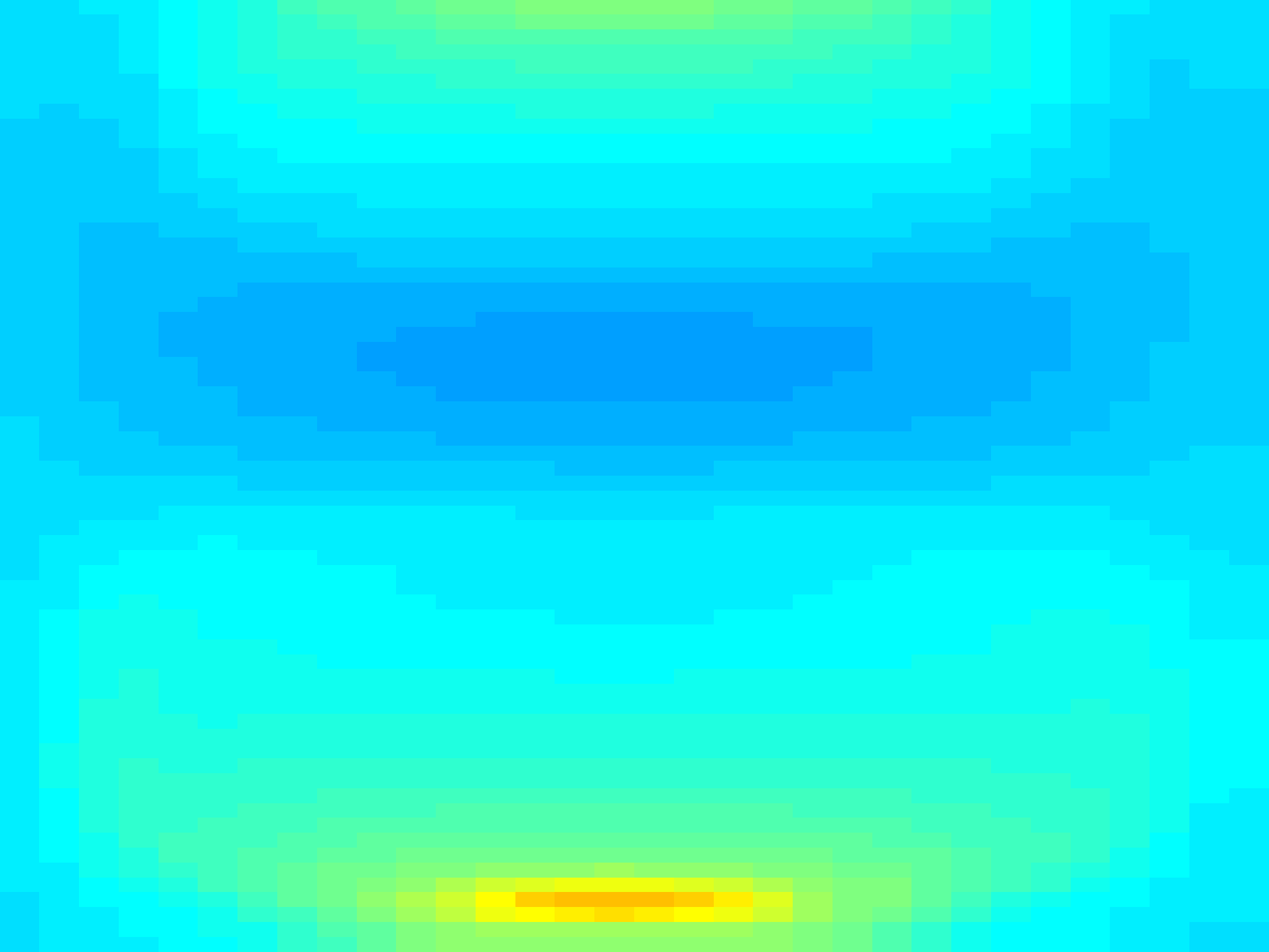}&
    \includegraphics[width=3cm]{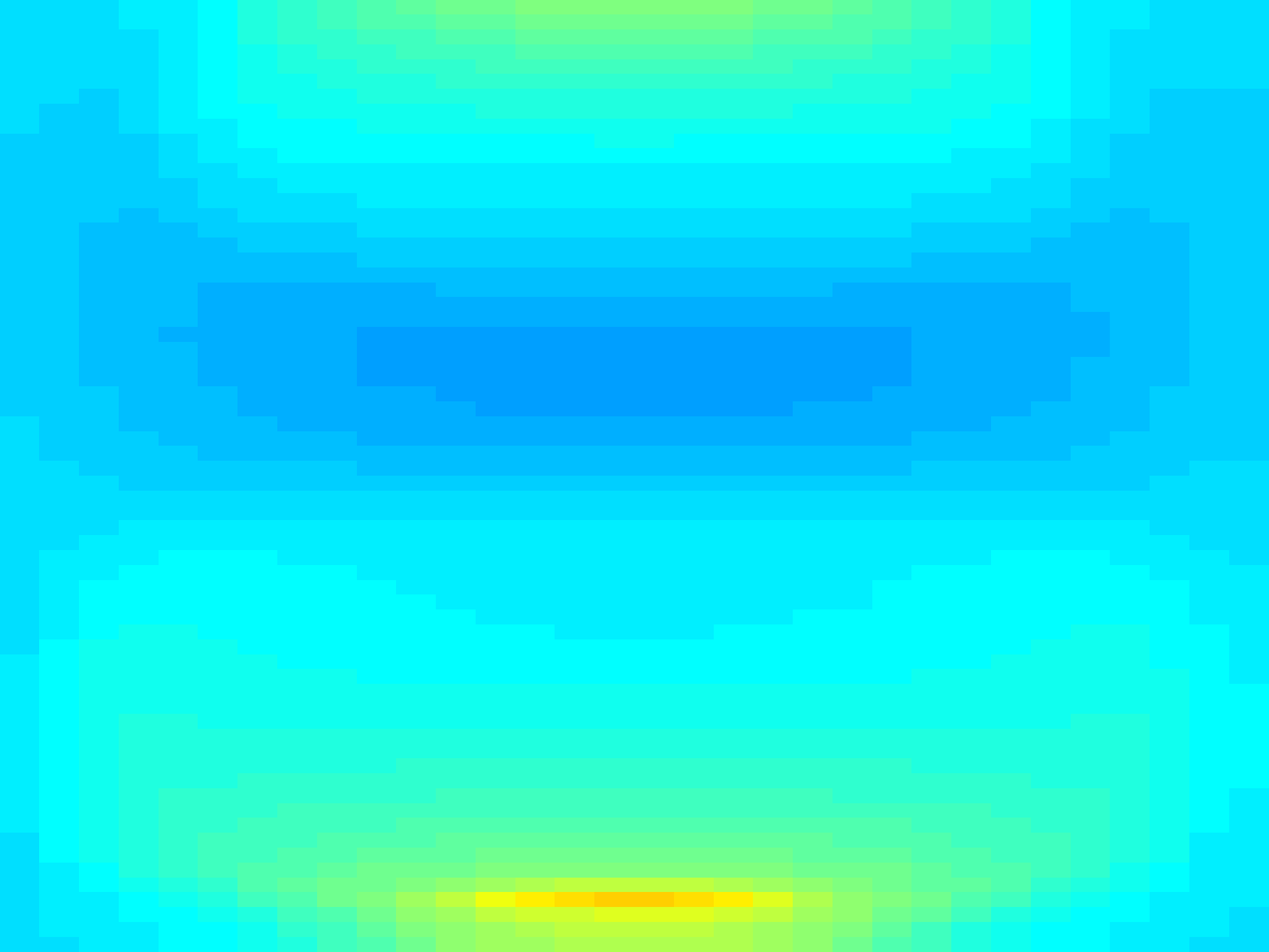}\\
    Original & Cubic Spline (30.01dB)& DCT (38.49dB)\\
    \includegraphics[width=3cm]{lattice_3d_original_band_151}&
    \includegraphics[width=3cm]{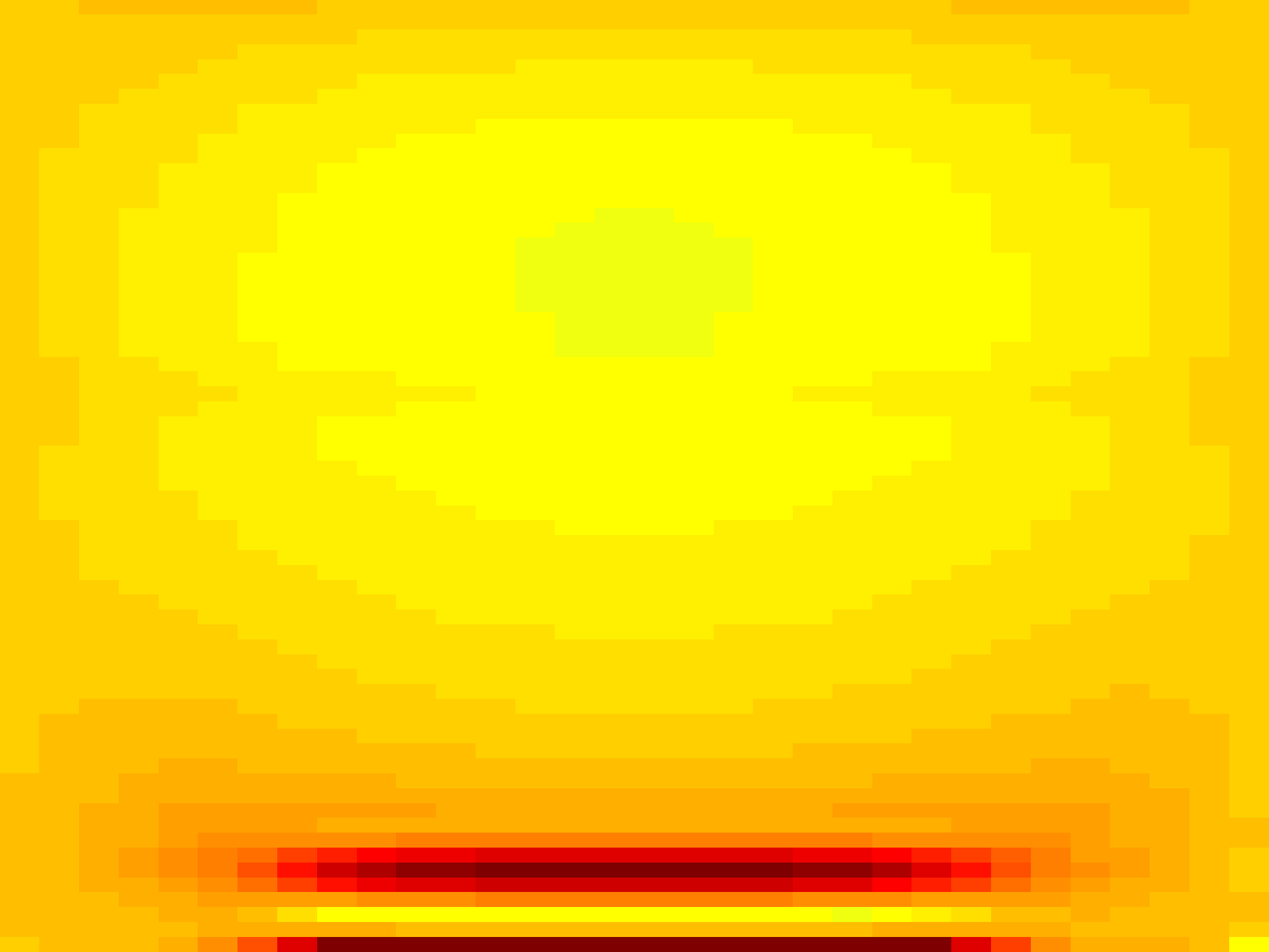}&
    \includegraphics[width=3cm]{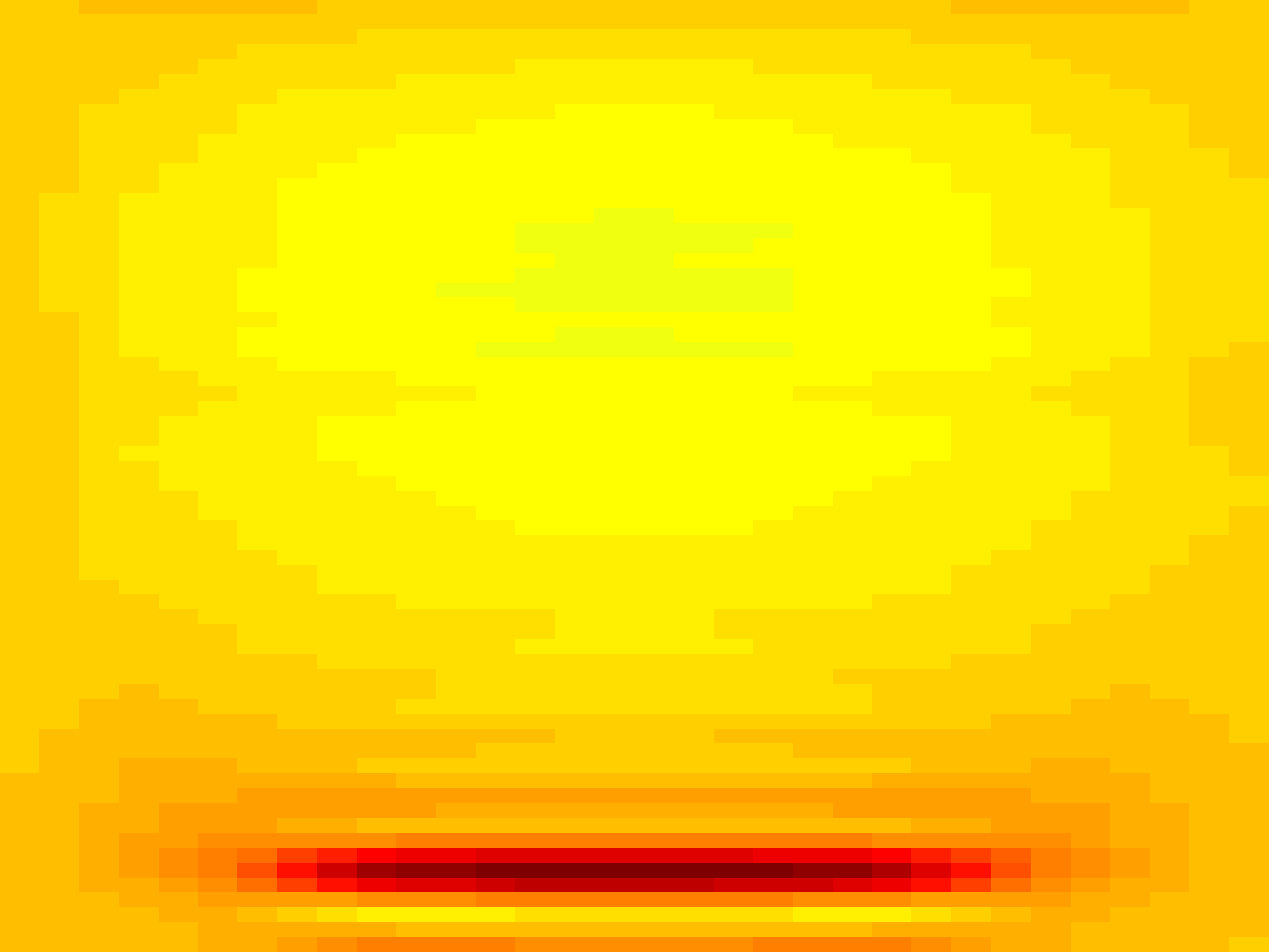}\\
    DFT (32.51dB)& Wavelet (38.15dB)& LDMM (\textbf{39.93dB})\\
    \includegraphics[width=3cm]{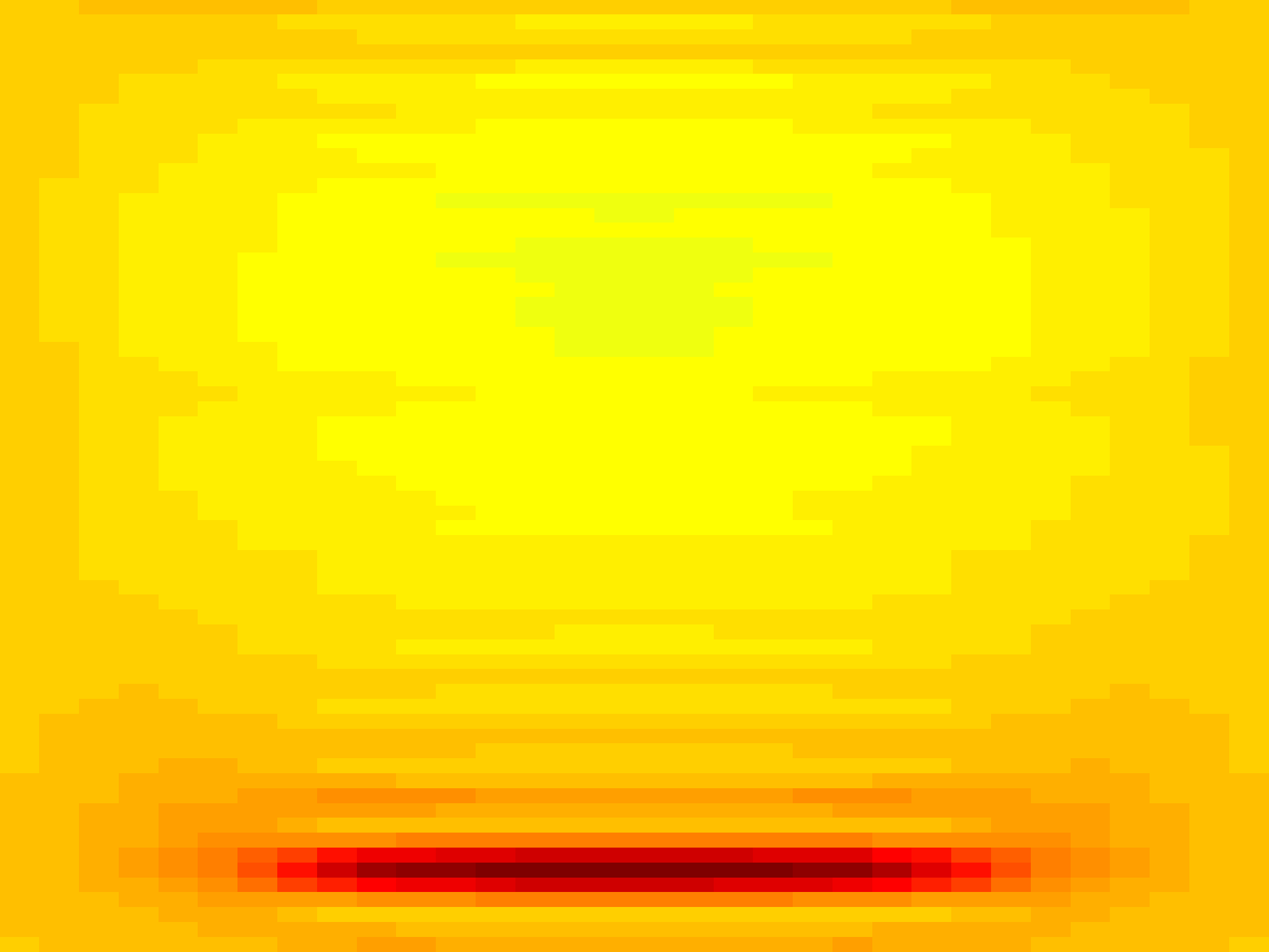}&
    \includegraphics[width=3cm]{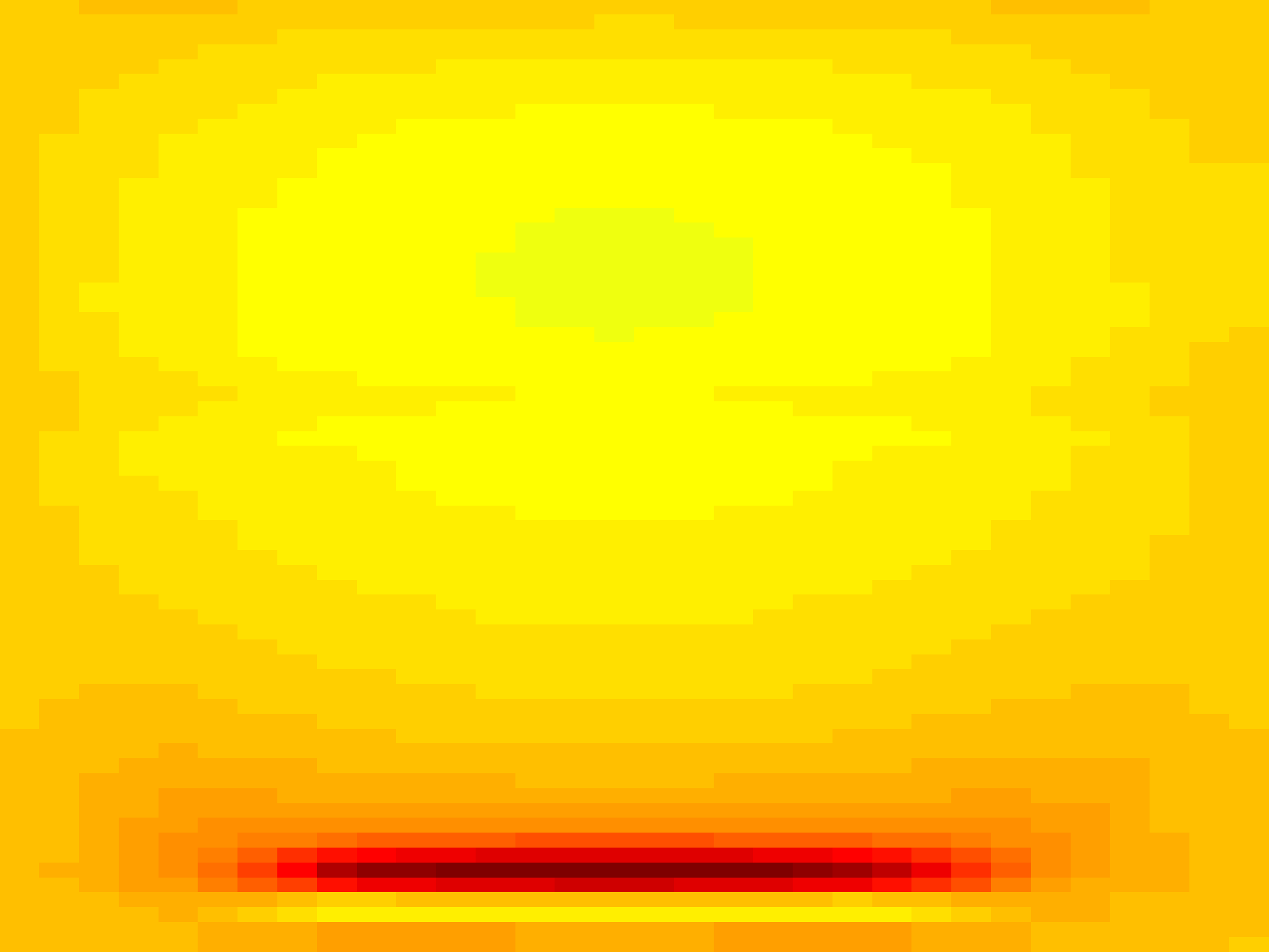}&
    \includegraphics[width=3cm]{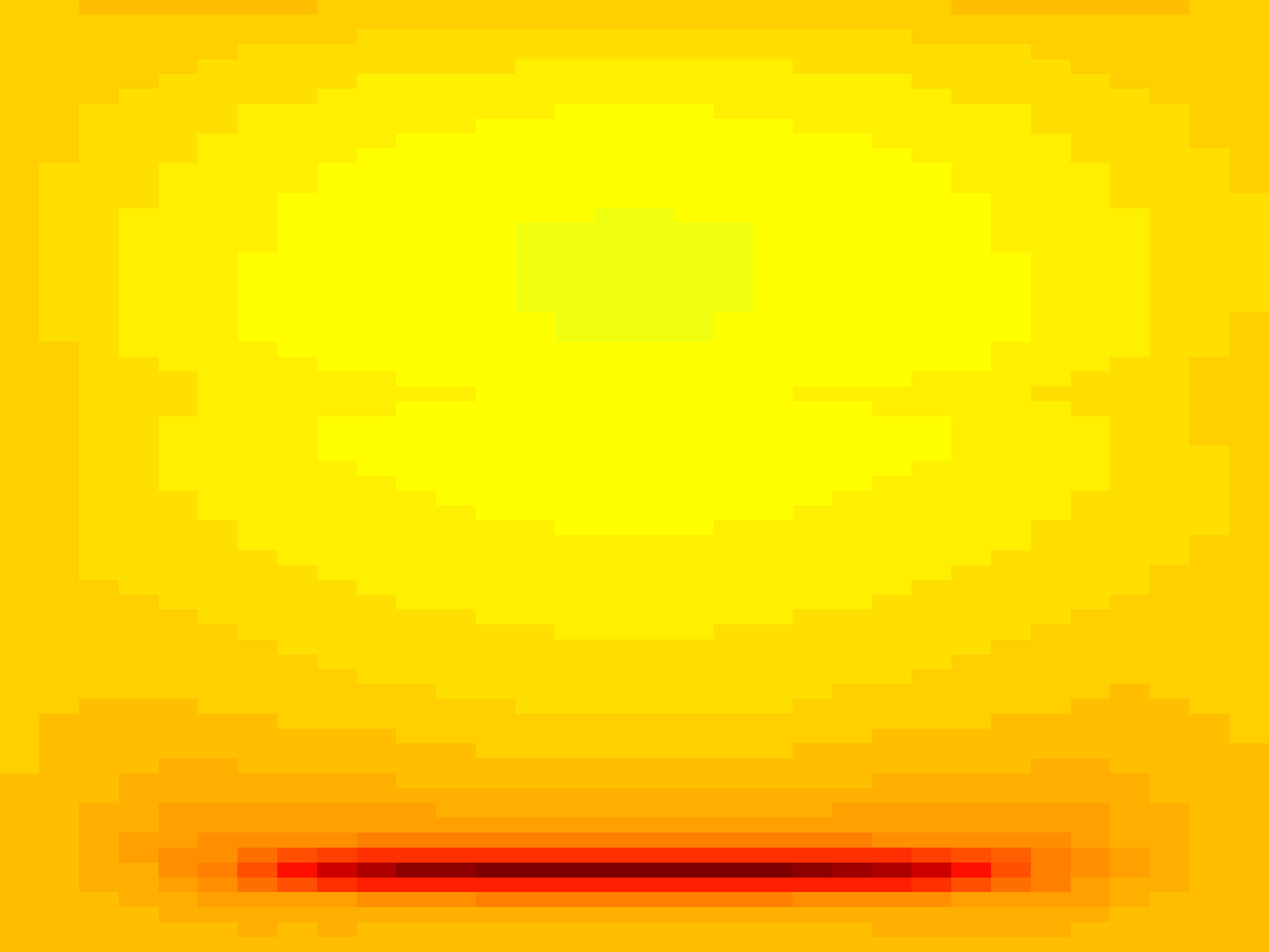}
  \end{tabular}
  \caption{Interpolation of the 3D lattice data set from regular sampling with spacing $2\times 2\times 2$. The original angular flux at $x = 0.24$ and $x = 1.18$ are shown in the first figures on the first and third row. The results of cubic spline, DCT, DFT, wavelet, and LDMM are shown in the remaining five figures.}
  \label{fig:down_lattice_3d_222}
\end{figure}

\begin{figure}[H]
  \centering
  \begin{tabular}{ccc}
    Original & Cubic Spline (36.47dB)& DCT (37.35dB)\\
    \includegraphics[width=3cm]{shock_3d_original_band_19}&
    \includegraphics[width=3cm]{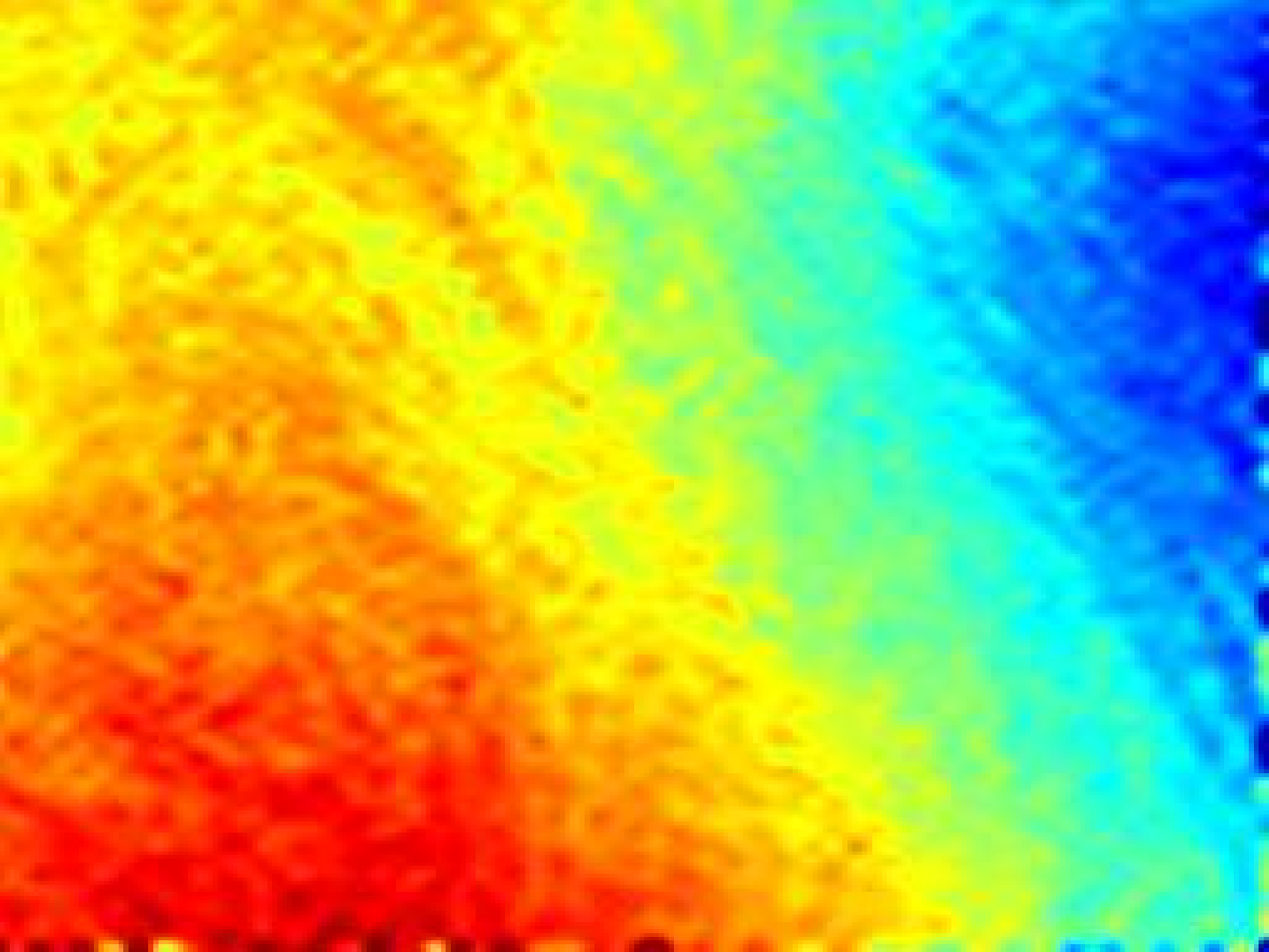}&
    \includegraphics[width=3cm]{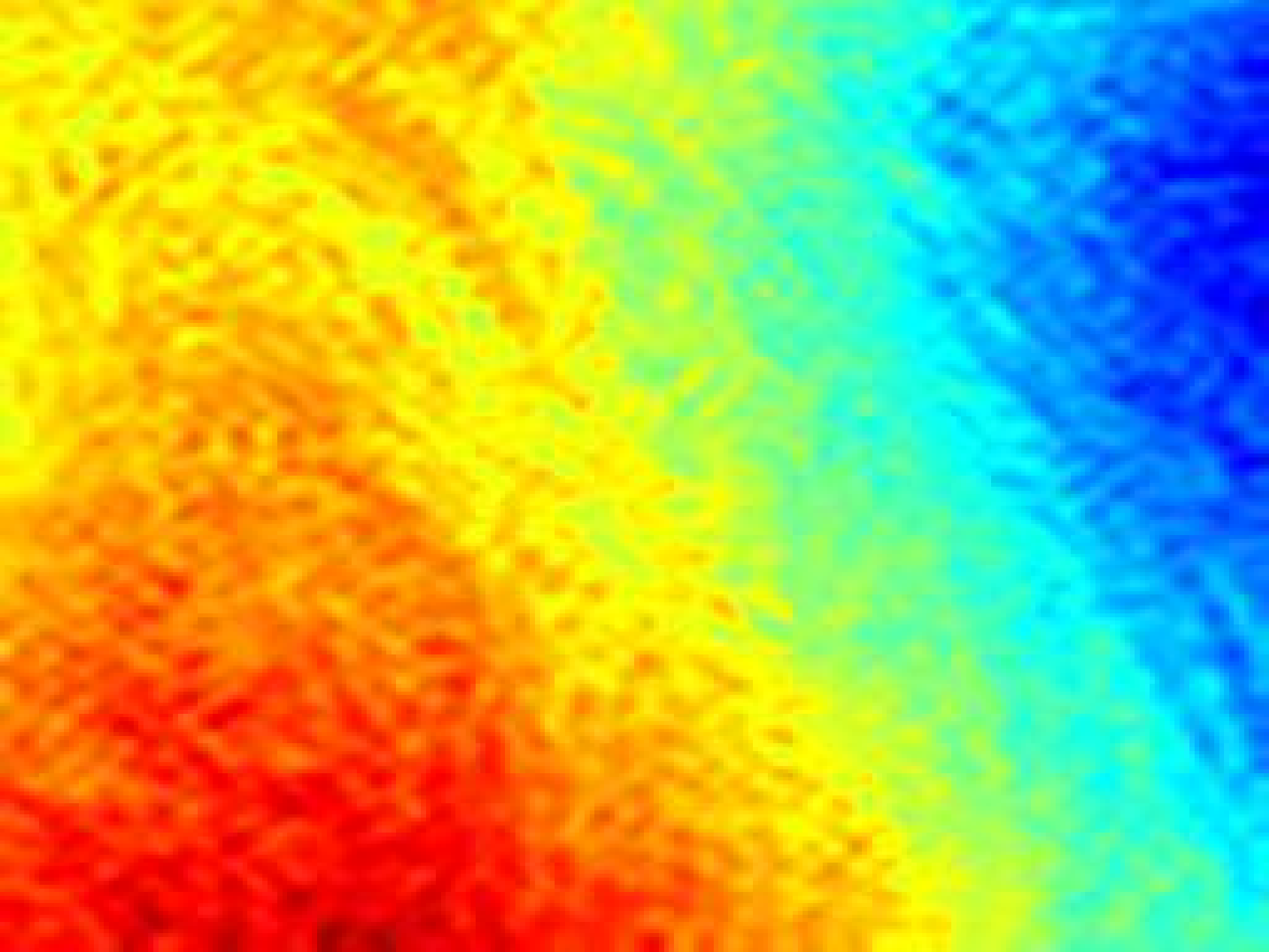}\\
    DFT (32.45dB)& Wavelet (37.02dB)& LDMM (\textbf{39.18dB})\\
    \includegraphics[width=3cm]{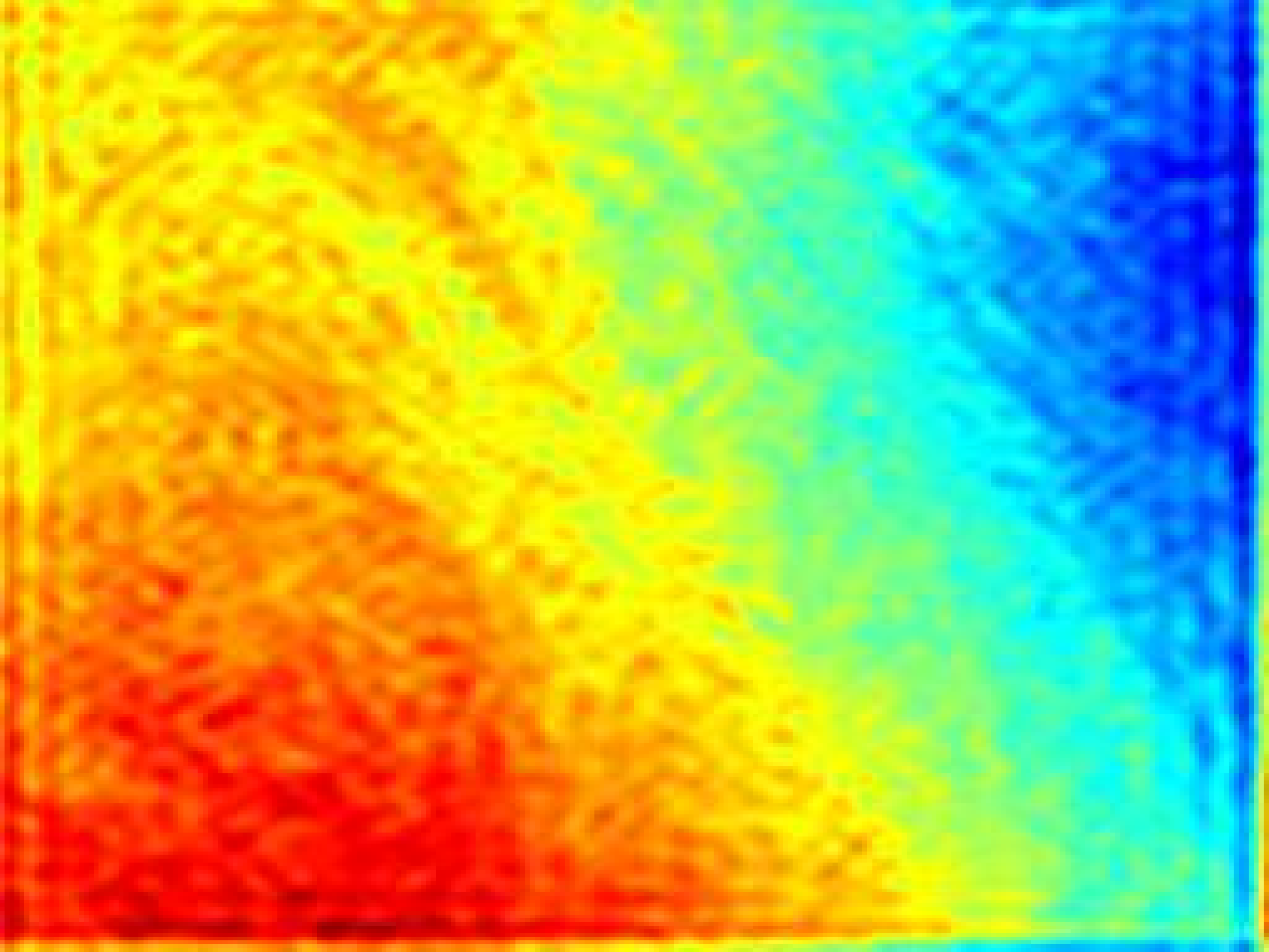}&
    \includegraphics[width=3cm]{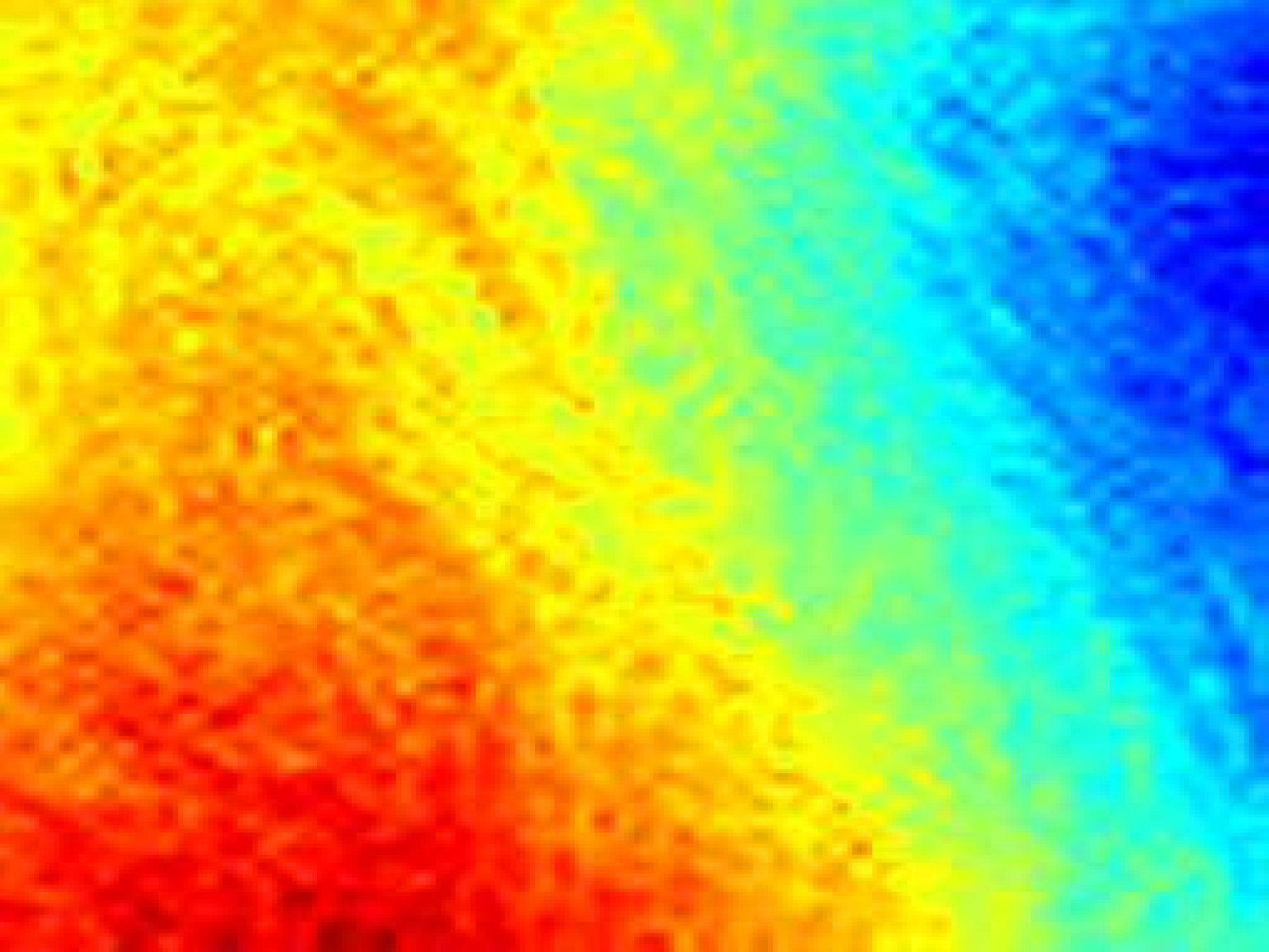}&
    \includegraphics[width=3cm]{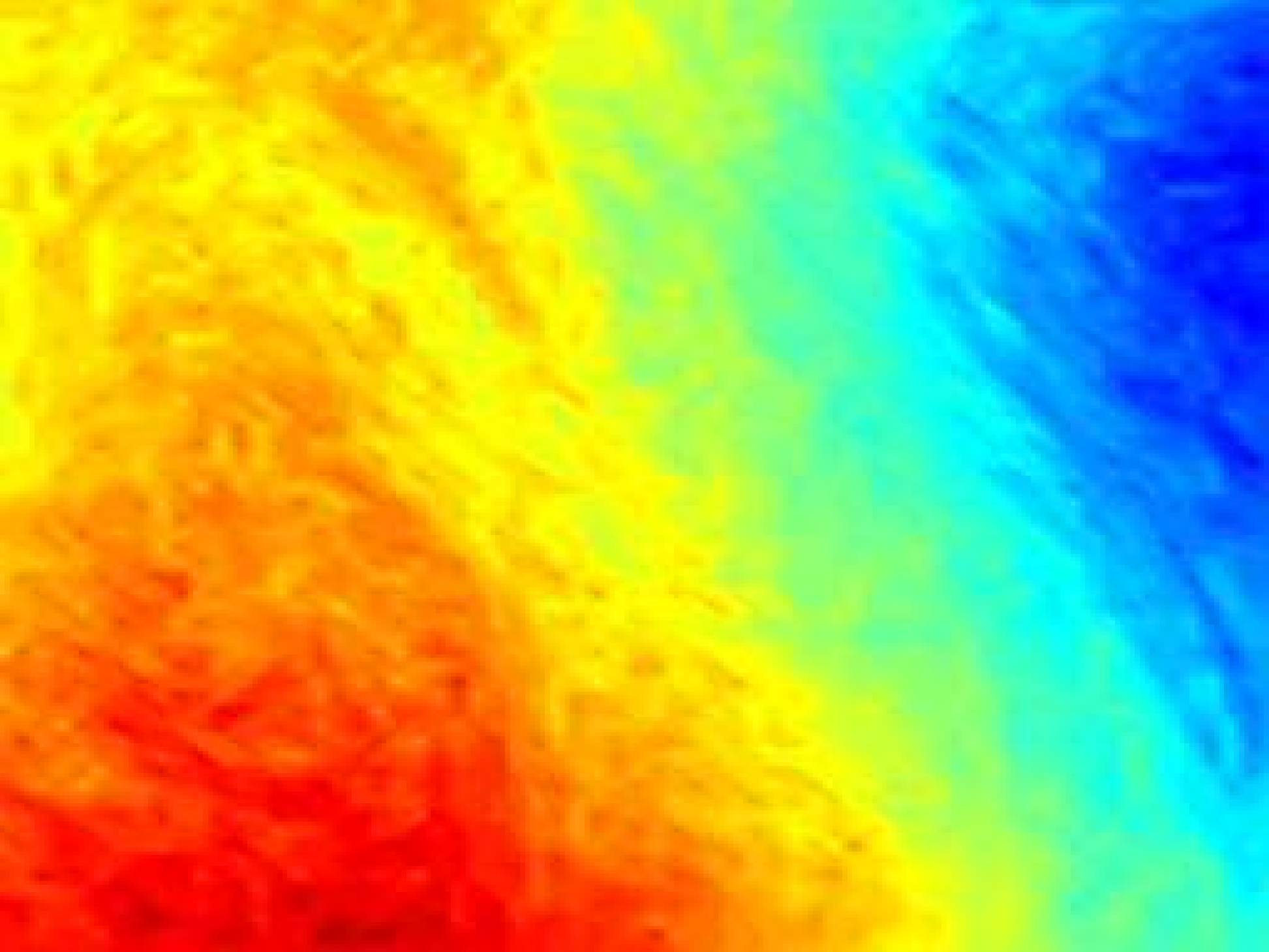}\\
    Original & Cubic Spline (36.47dB)& DCT (37.35dB)\\
    \includegraphics[width=3cm]{shock_3d_original_band_29}&
    \includegraphics[width=3cm]{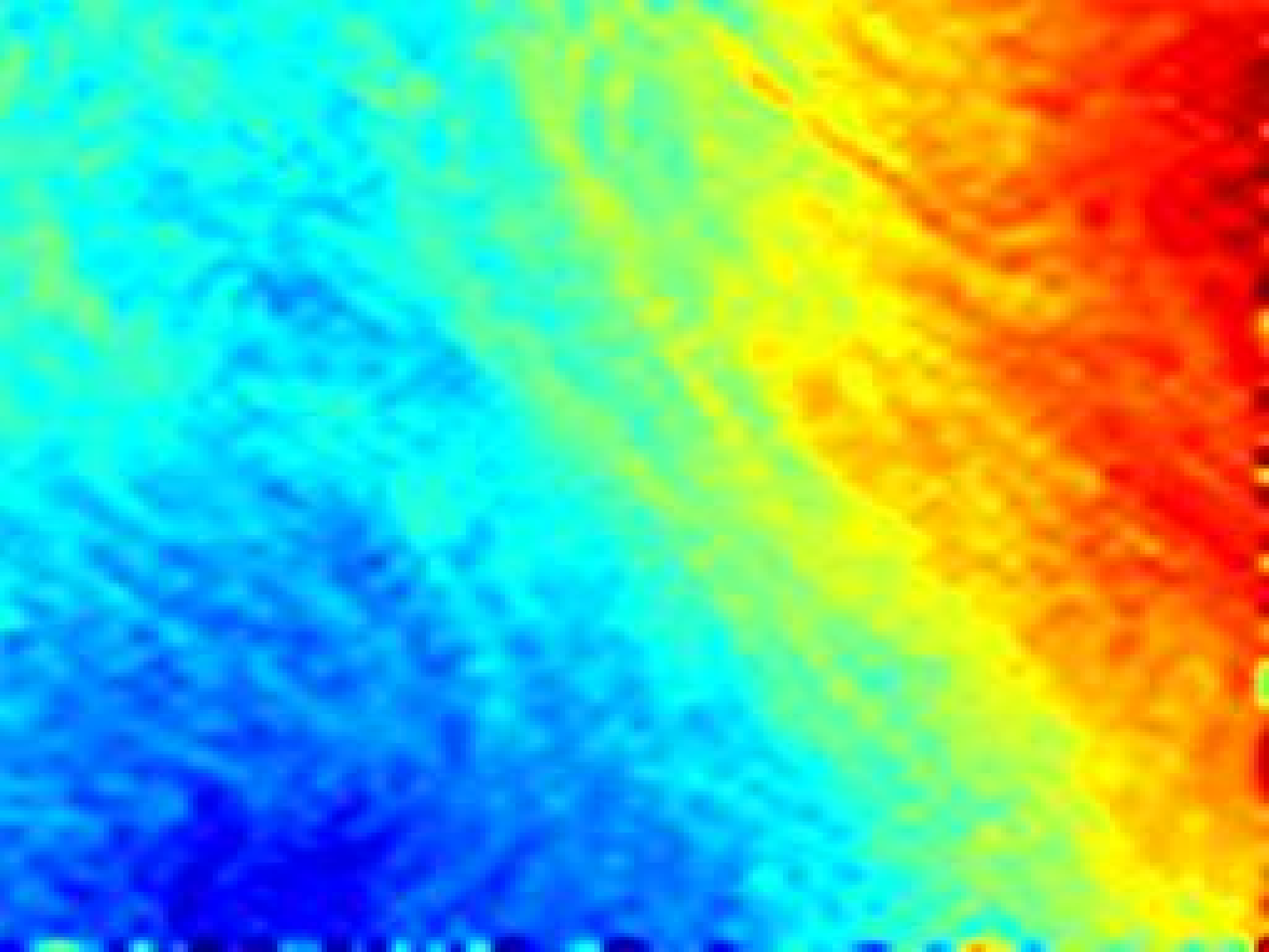}&
    \includegraphics[width=3cm]{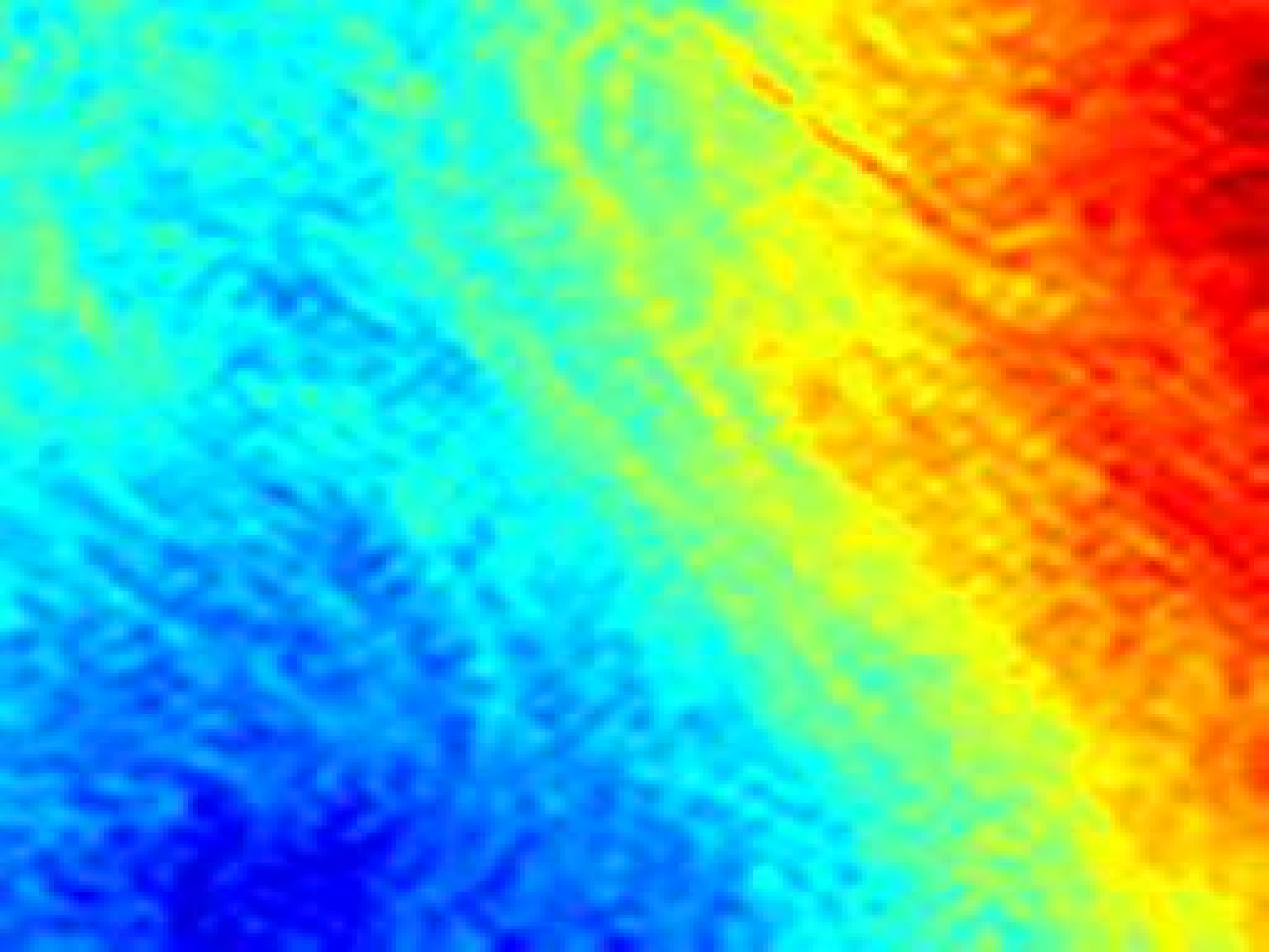}\\
    DFT (32.45dB)& Wavelet (37.02dB)& LDMM (\textbf{39.18dB})\\
    \includegraphics[width=3cm]{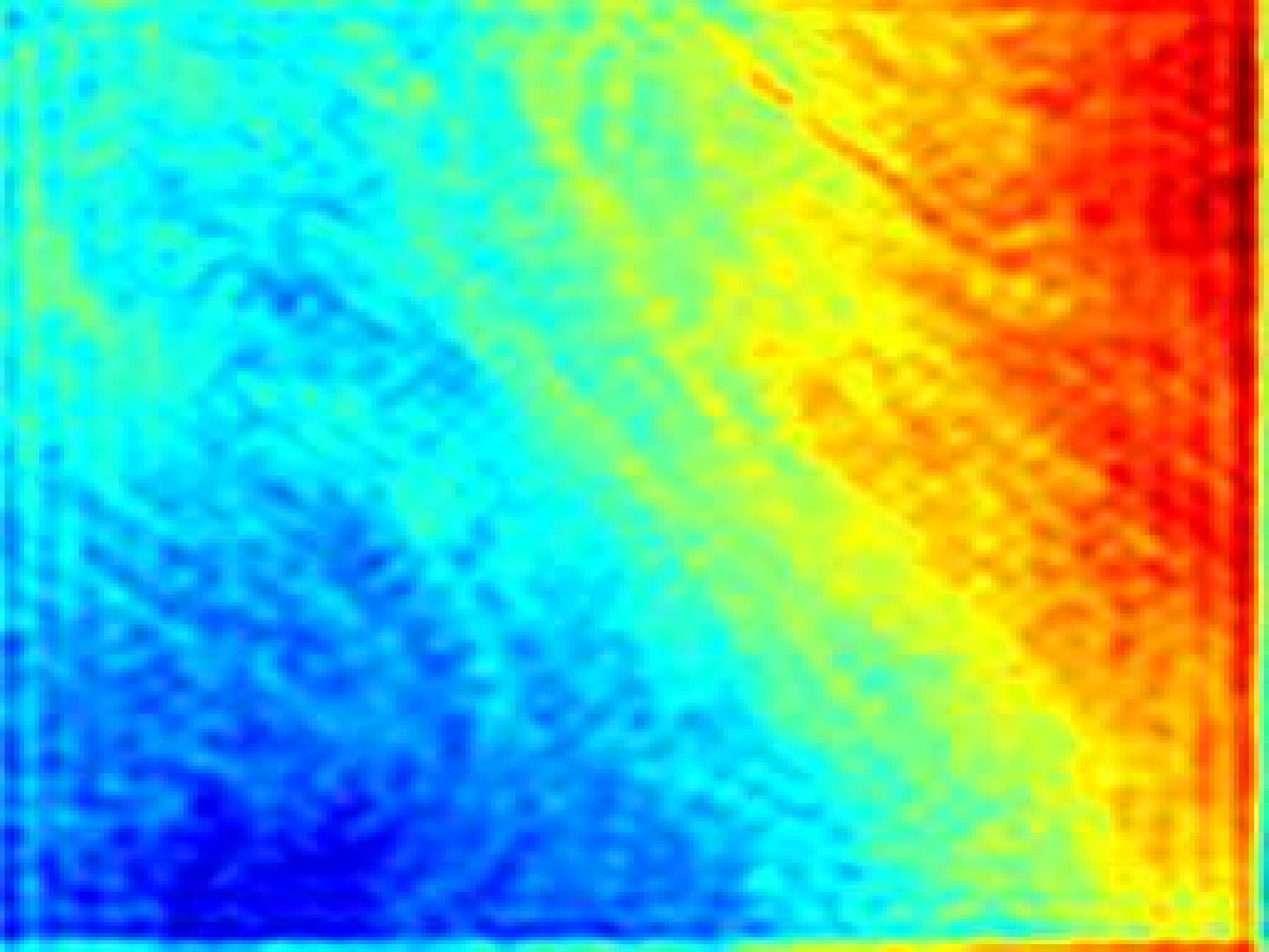}&
    \includegraphics[width=3cm]{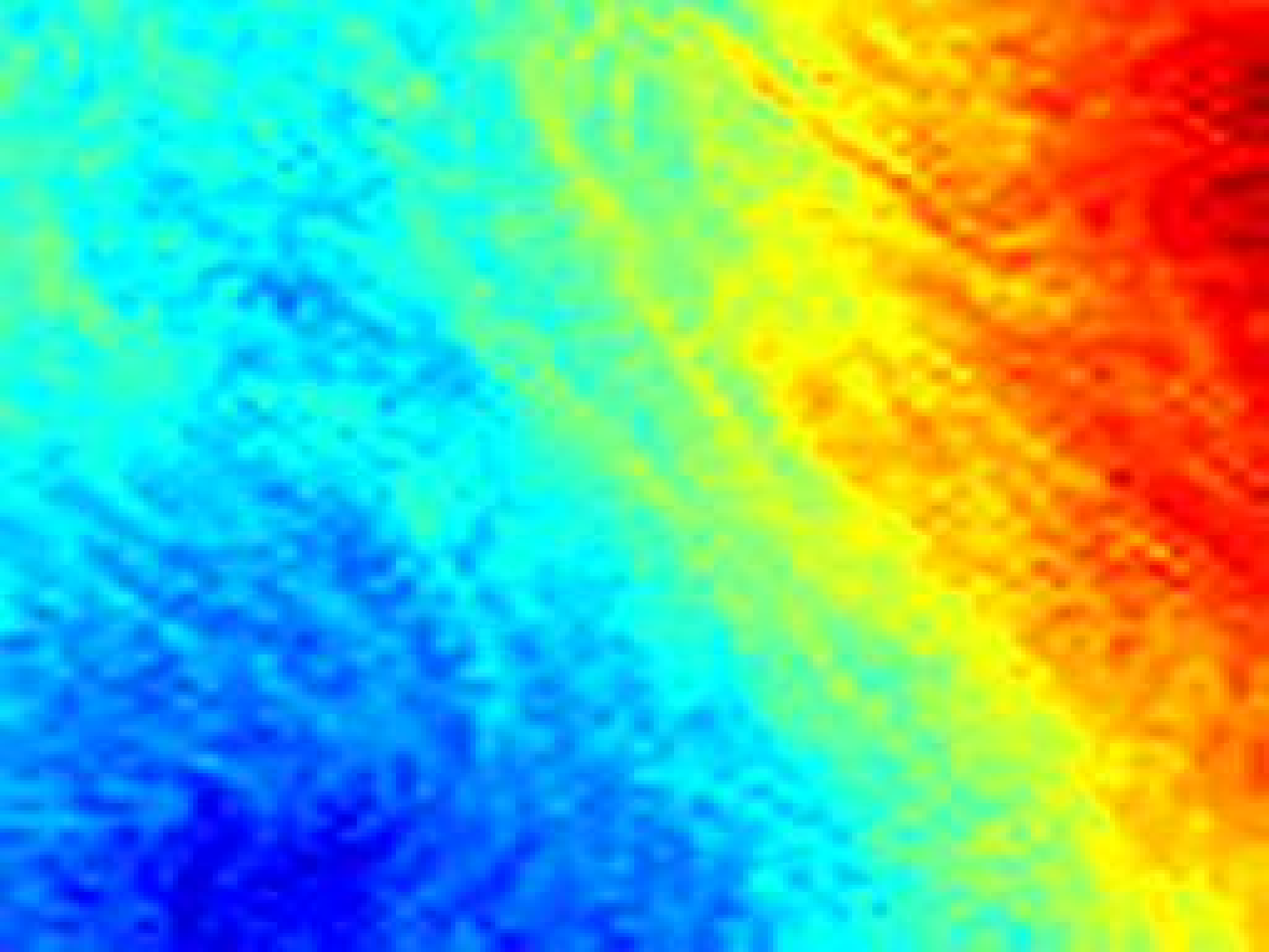}&
    \includegraphics[width=3cm]{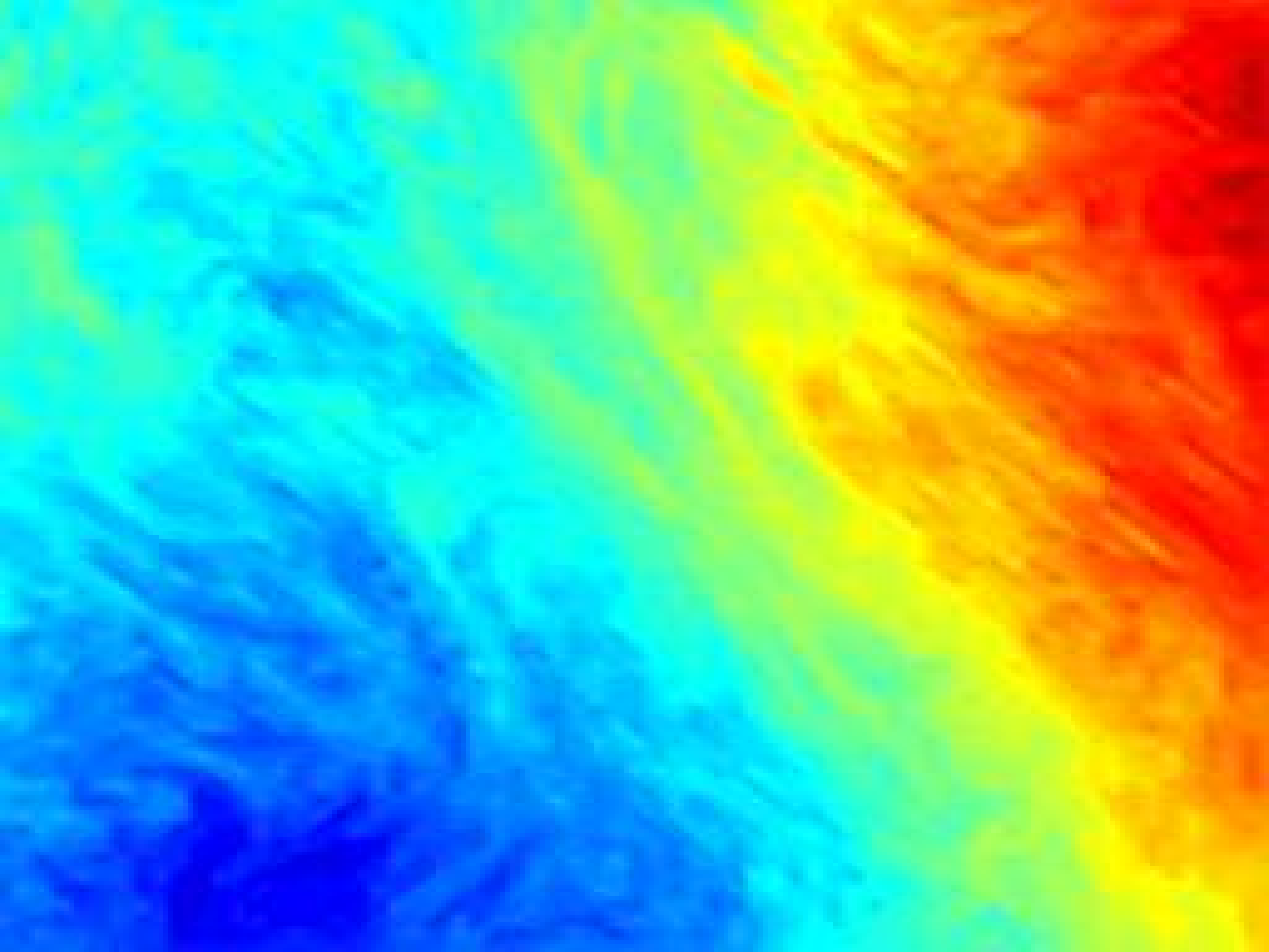}
  \end{tabular}
  \caption{Interpolation of the 3D plasma (distribution function) data set from regular sampling with spacing $4\times 4\times 1$. Two spatial cross sections of the original data are shown in the first figures on the first and third row. The results of cubic spline, DCT, DFT, wavelet, and LDMM are shown in the remaining five figures.}
  \label{fig:down_shock_3d_441}
\end{figure}

\begin{figure}[H]
  \centering
  \begin{tabular}{ccc}
    Original & Cubic Spline (30.97dB)& DCT (33.91dB)\\
    \includegraphics[width=3cm]{shock_3d_original_band_19}&
    \includegraphics[width=3cm]{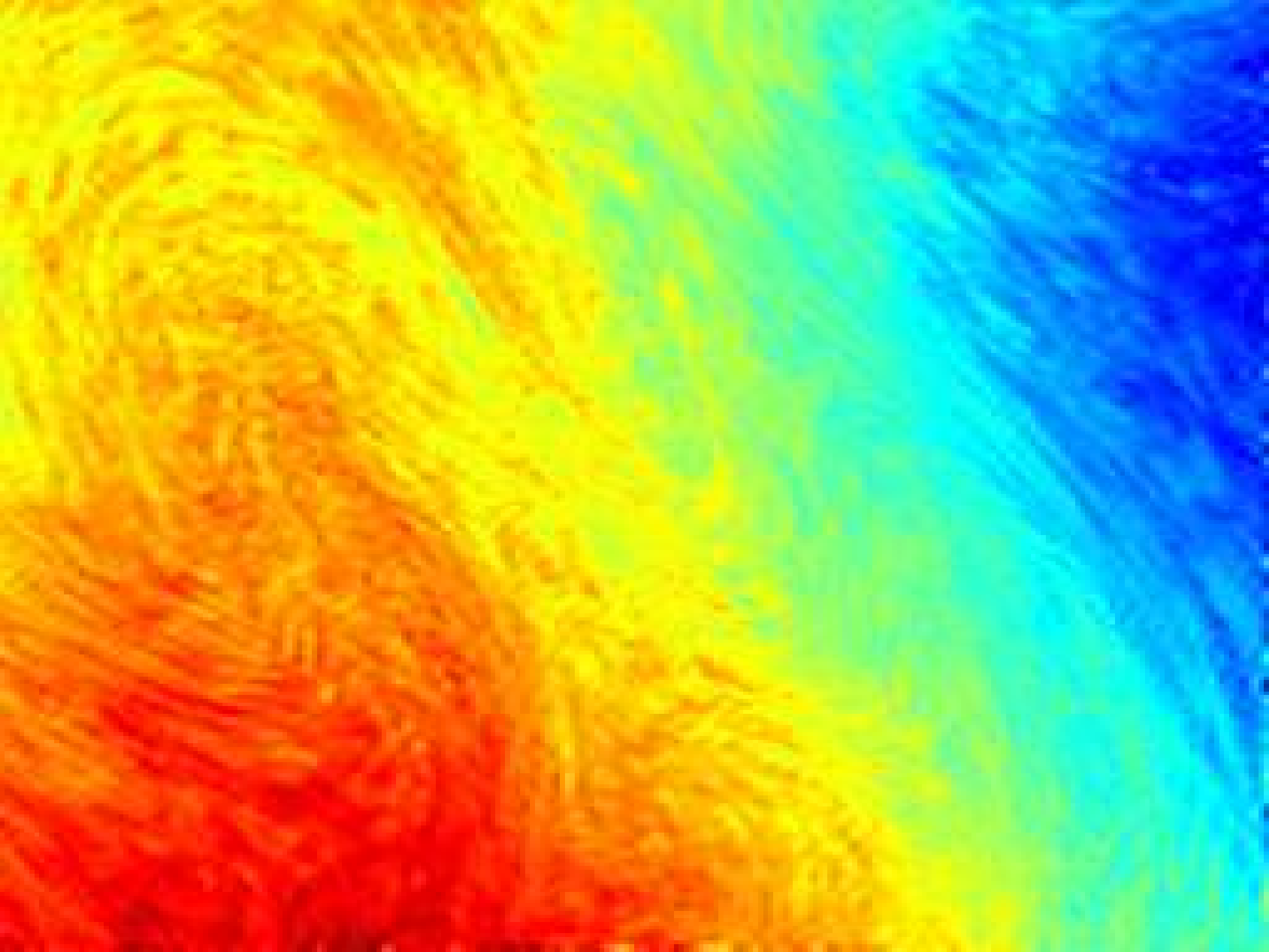}&
    \includegraphics[width=3cm]{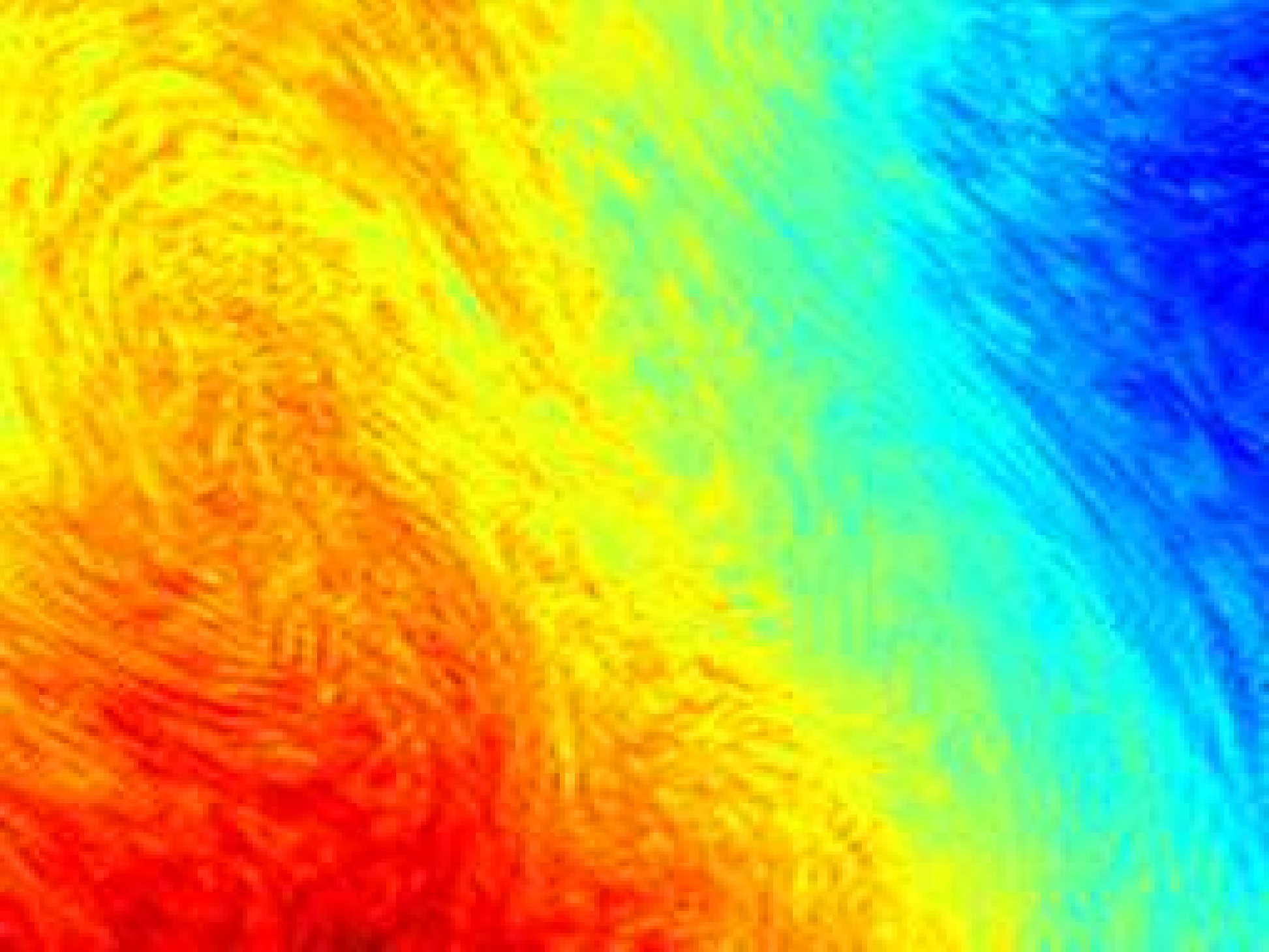}\\
    DFT (31.88dB)& Wavelet (32.81dB)& LDMM (\textbf{35.01dB})\\
    \includegraphics[width=3cm]{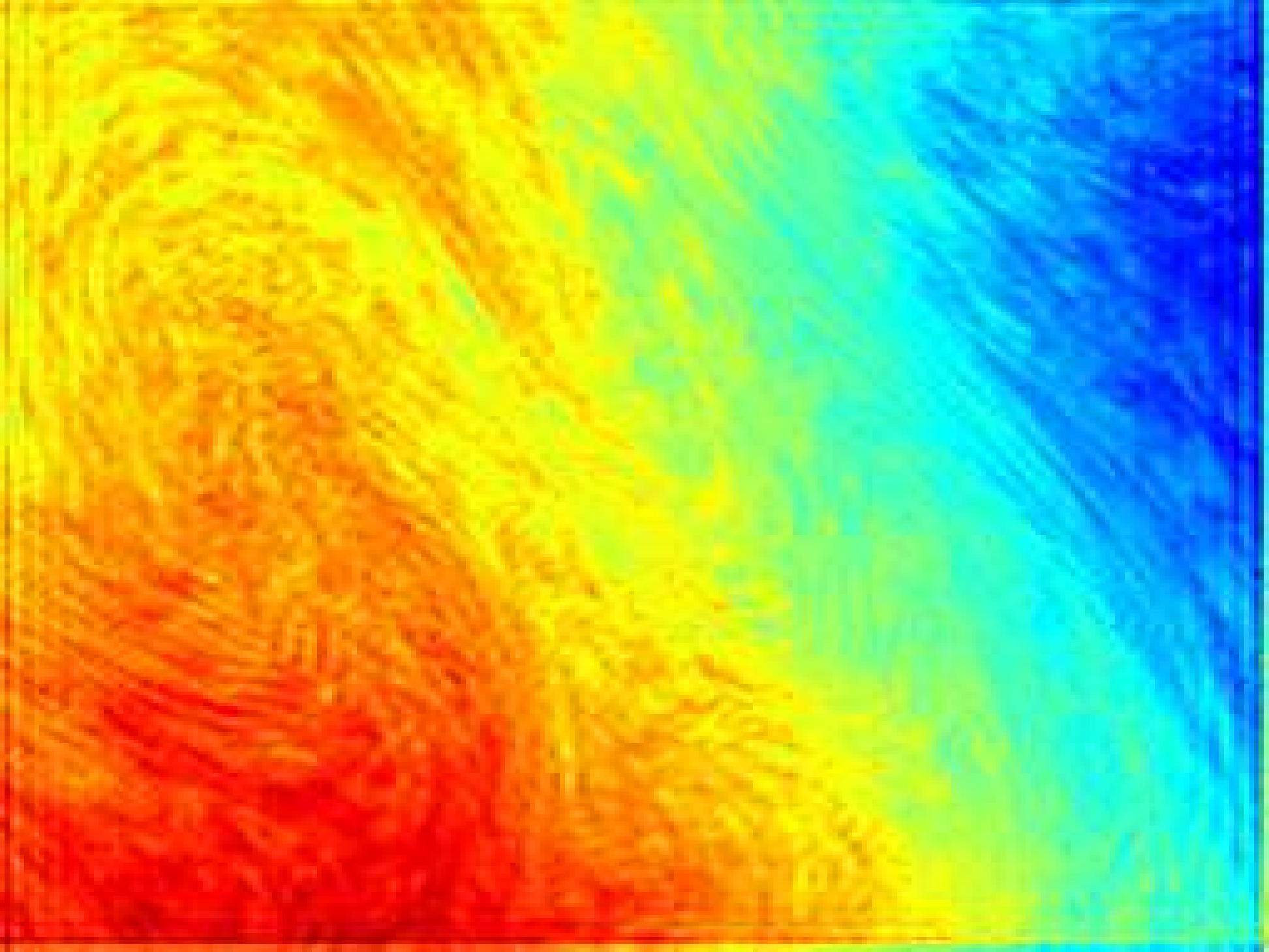}&
    \includegraphics[width=3cm]{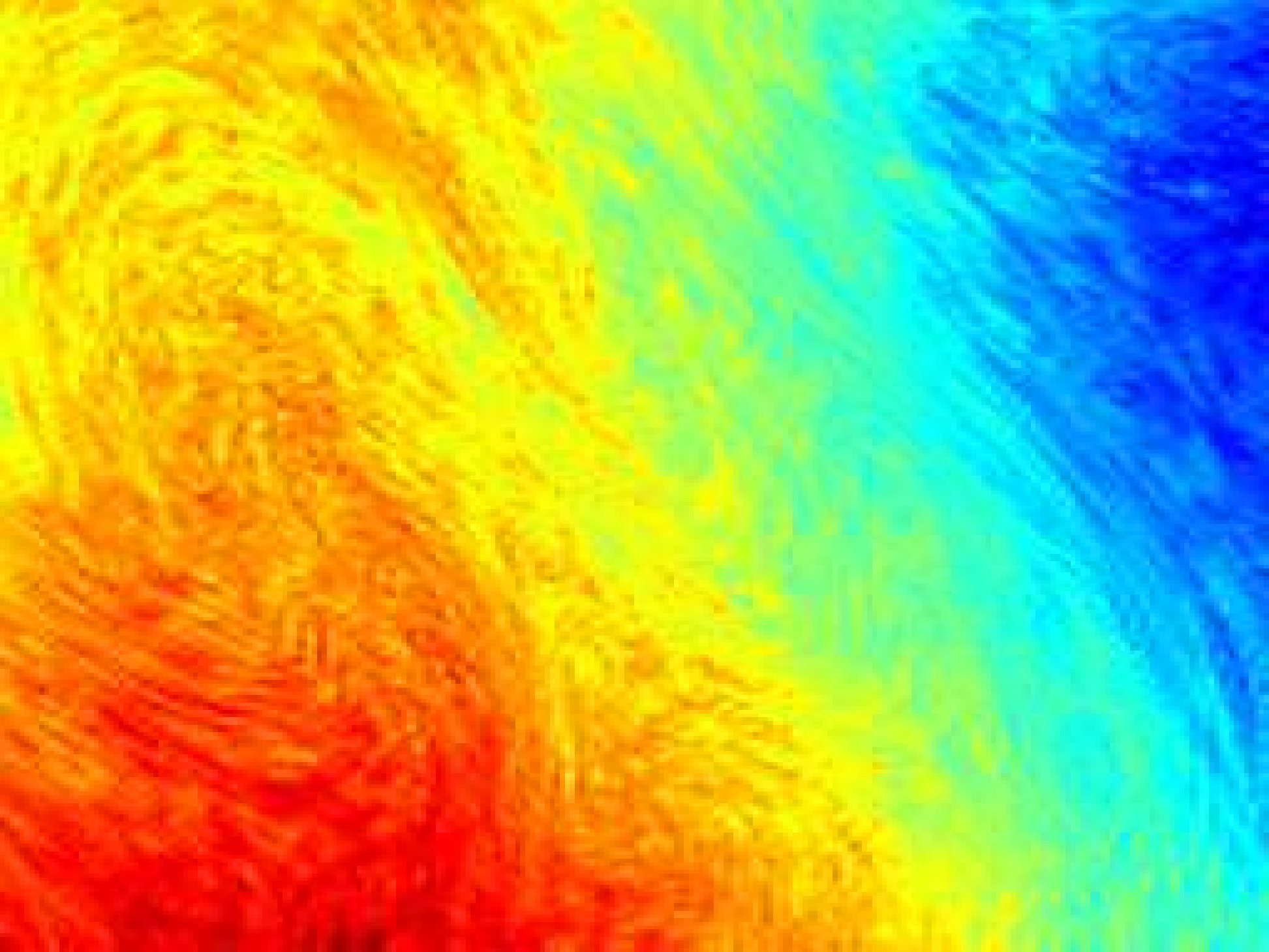}&
    \includegraphics[width=3cm]{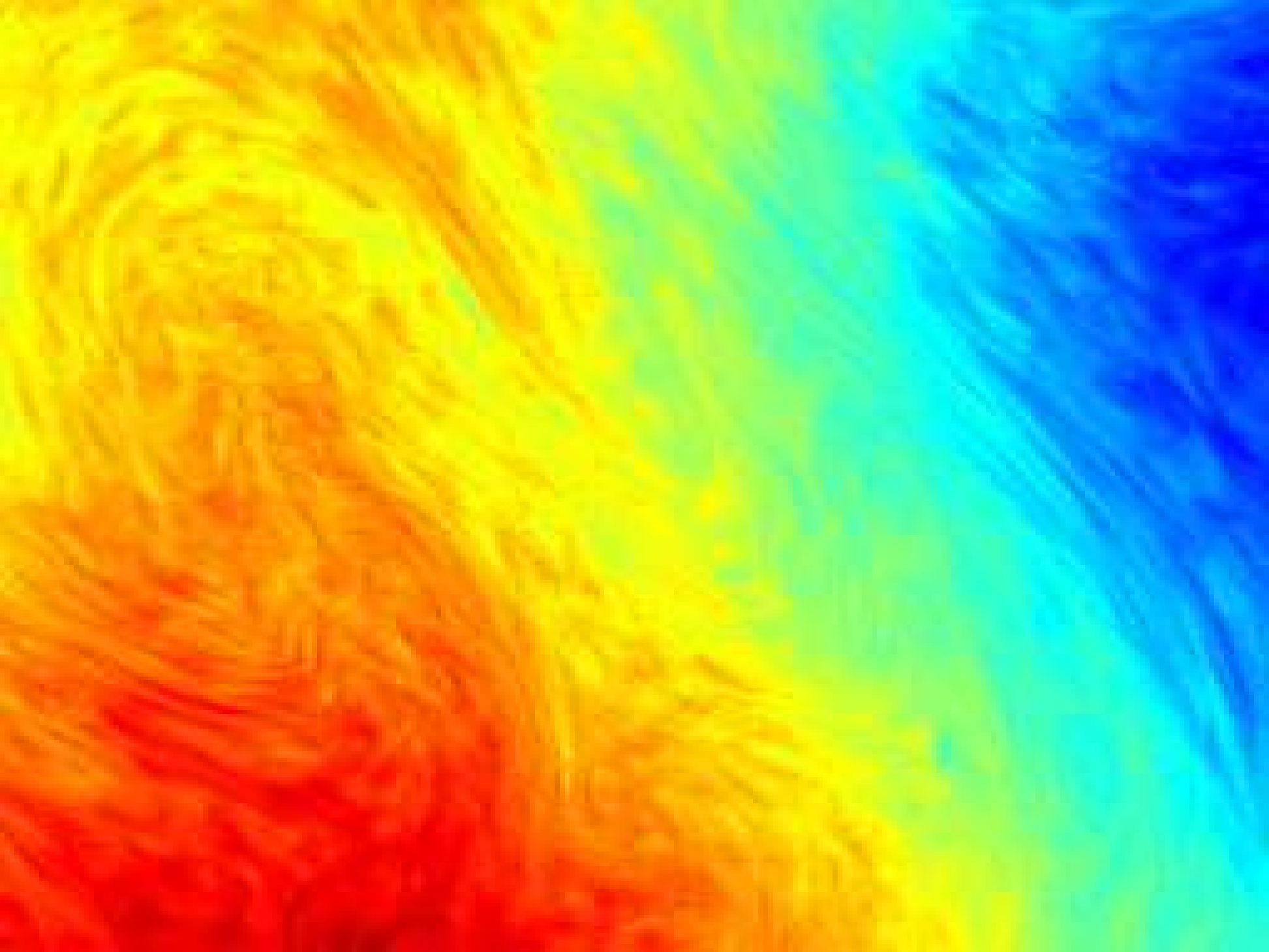}\\
    Original & Cubic Spline (30.97dB)& DCT (33.91dB)\\
    \includegraphics[width=3cm]{shock_3d_original_band_29}&
    \includegraphics[width=3cm]{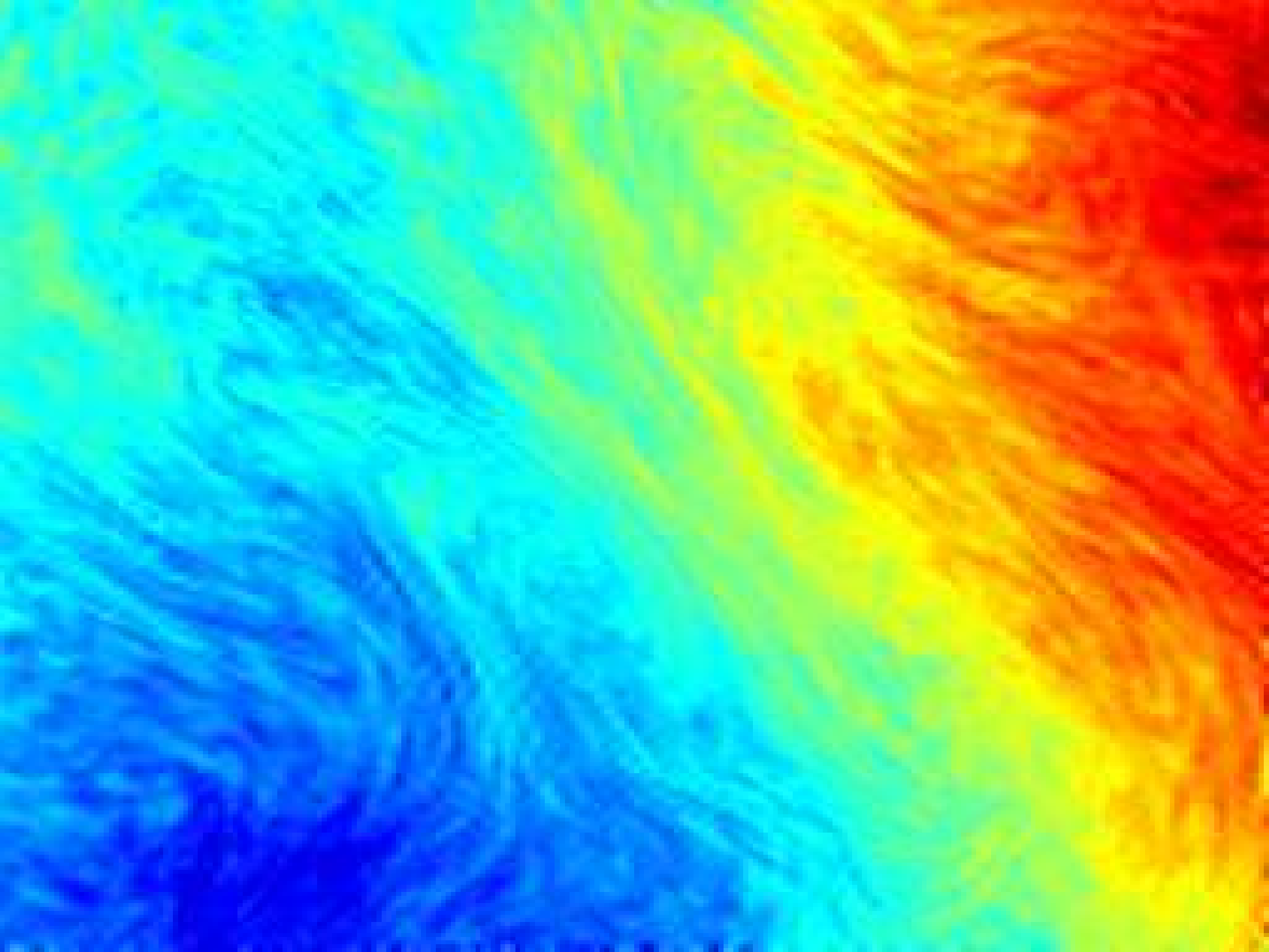}&
    \includegraphics[width=3cm]{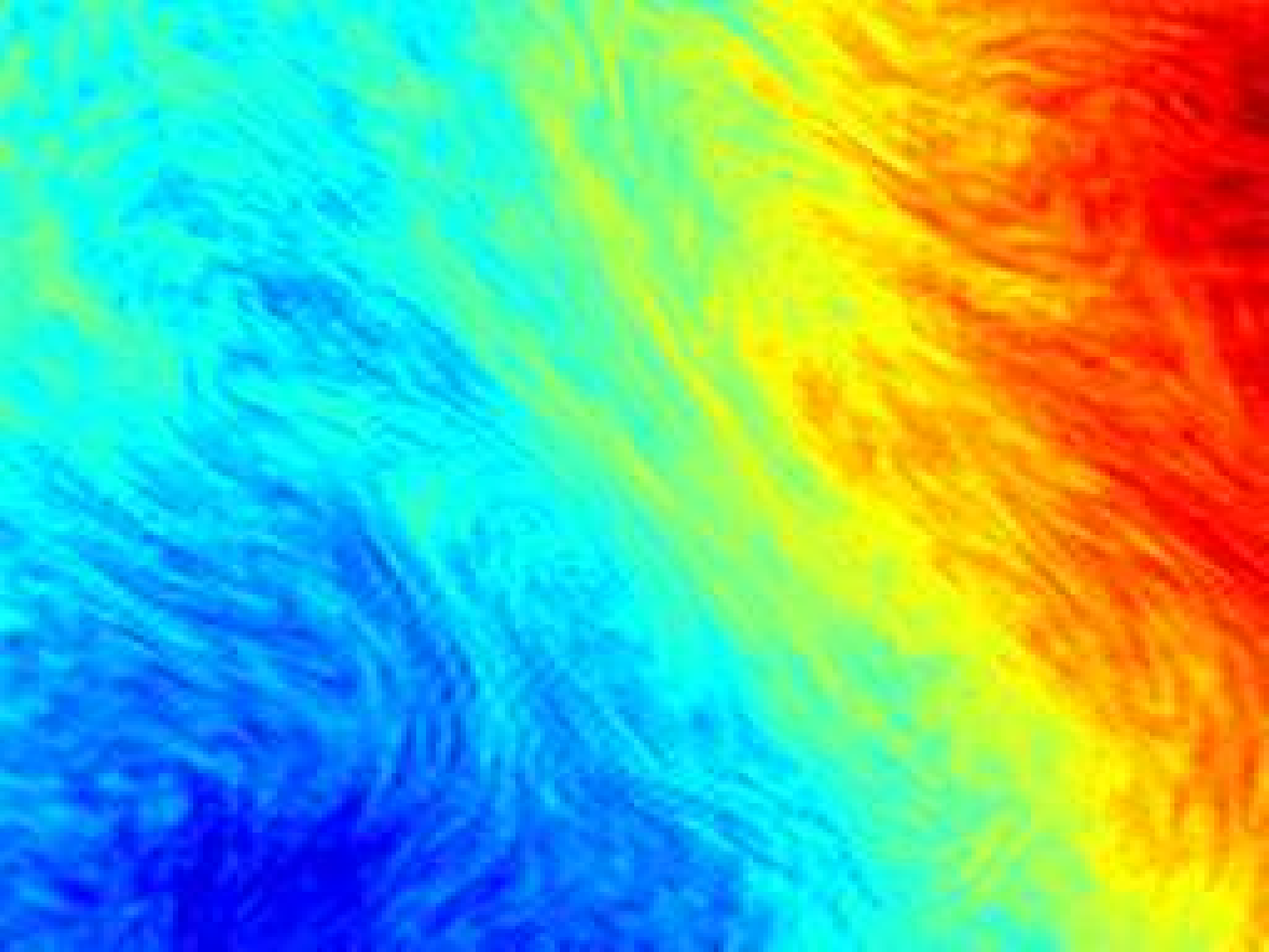}\\
    DFT (31.88dB)& Wavelet (32.81dB)& LDMM (\textbf{35.01dB})\\
    \includegraphics[width=3cm]{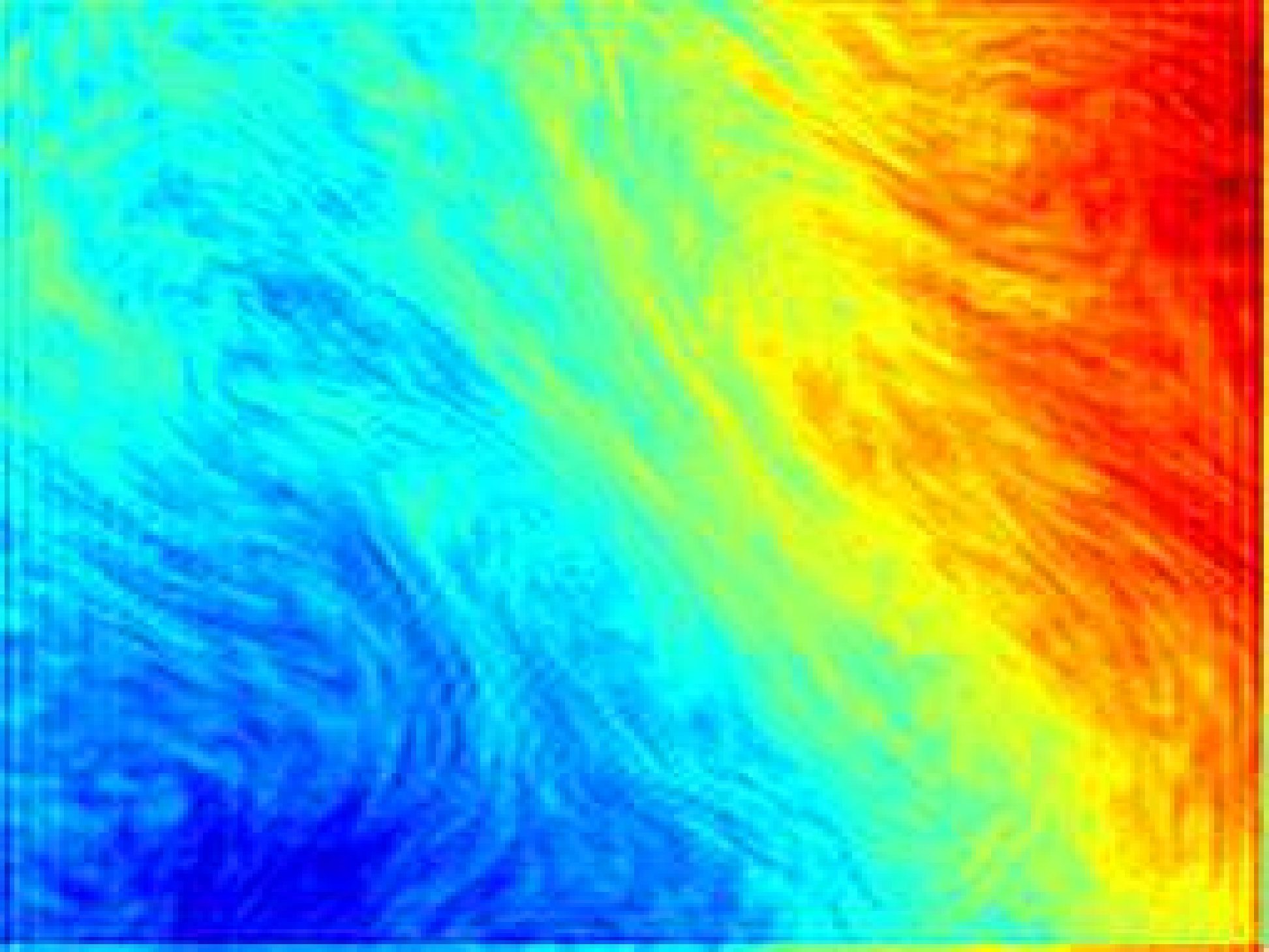}&
    \includegraphics[width=3cm]{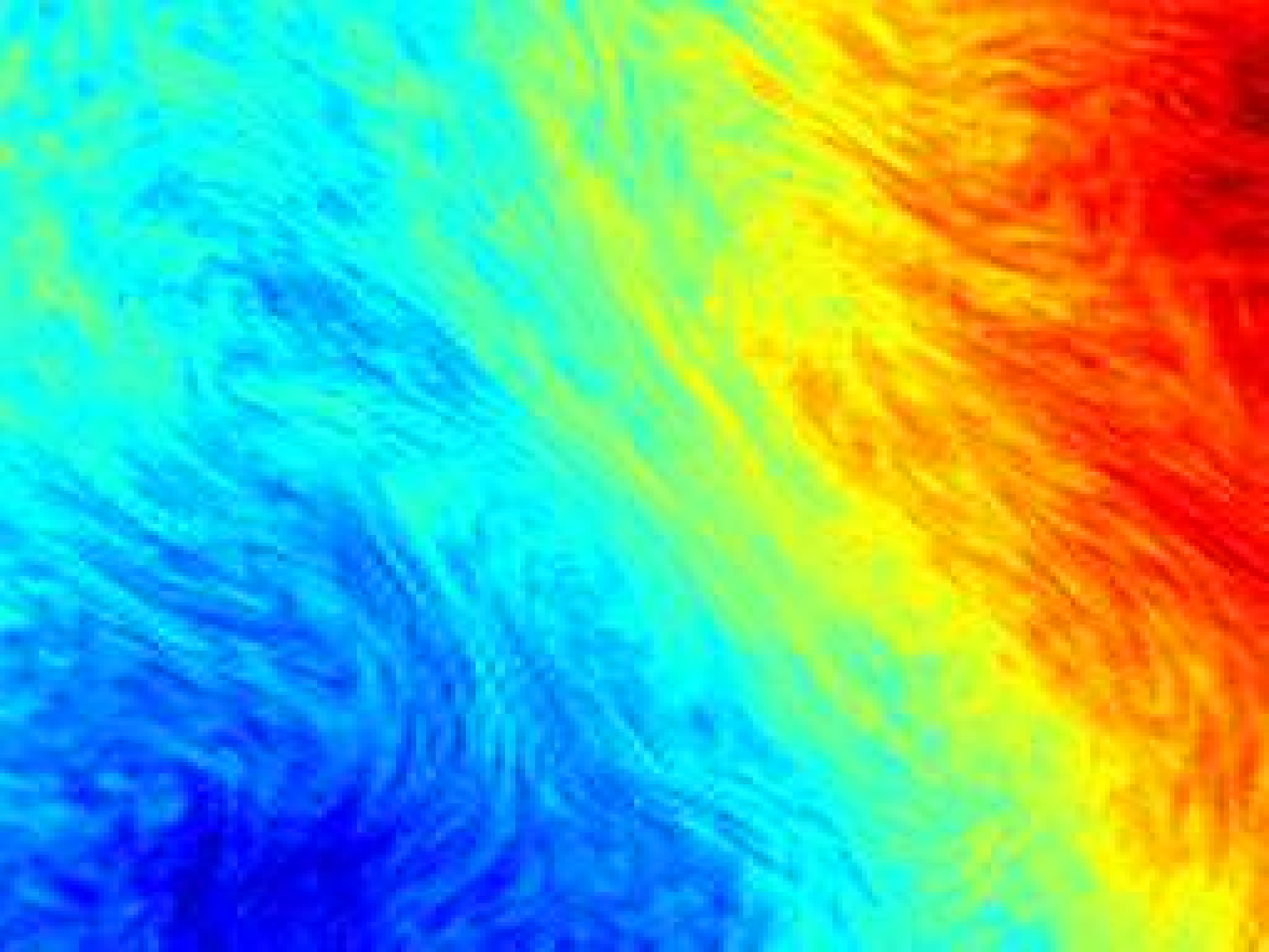}&
    \includegraphics[width=3cm]{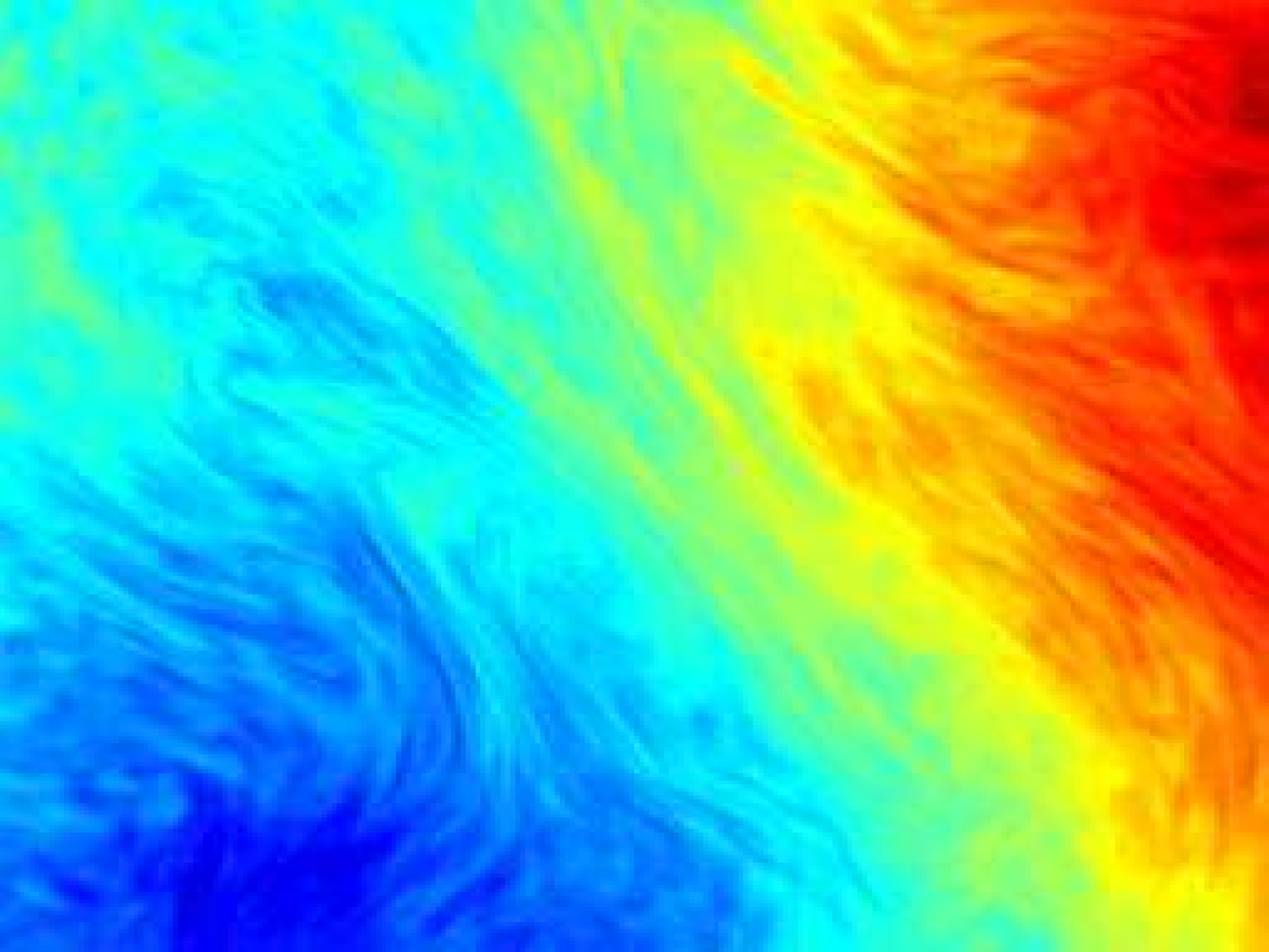}
  \end{tabular}
  \caption{Interpolation of the 3D plasma (distribution function) data set from regular sampling with spacing $2\times 2\times 2$. Two spatial cross sections of the original data are shown in the first figures on the first and third row. The results of cubic spline, DCT, DFT, wavelet, and LDMM are shown in the remaining five figures.}
  \label{fig:down_shock_3d_222}
\end{figure}

The visual of the interpolation with regular sampling ($4\times 4$ for 2D data sets, $4 \times 4 \times 1$ and $2\times 2 \times 2$ for 3D data sets) are shown in Figure \ref{fig:down_2d}-\ref{fig:down_shock_3d_222}. The errors in different norms are displayed in Table \ref{tab:error_down_antonio_2d_44}-\ref{tab:error_down_shock_3d}. It can be observed that the results of LDMM are significantly more accurate than the DCT and cubic spline initializations, and the accuracy of the result does not depend on the choice of the initialization. Moreover, LDMM consistently outperforms all the other competing algorithms on every data set, except for some rare cases where LDMM is inferior in $L_1$ or $L_{\infty}$ norms.

\begin{table}[H]
  \centering
  \begin{tabular}{||c| c  c c c c c||}
    \hline
    $4\times 4$ & Cubic & DCT& DFT& Wavelet & LDMM (D) & LDMM (C)\\
    \hline
    $L_1$       &\textbf{0.0025} &0.0038 &0.0035 &0.0049 &0.0029 & 0.0028\\
    $L_2$       &0.0071 &0.0072 &0.0069 &0.0095 &\textbf{0.0060} & 0.0061\\
    $L_\infty$   &0.1789 &\textbf{0.0937} &0.0940 &0.1122 &0.0961 & 0.1005\\
    PSNR        &42.98  &42.88  &43.19  &40.48  &\textbf{44.40} & 44.33\\
    \hline
  \end{tabular}
  \caption{Errors of the interpolation of the 2D vortex data set from regular sampling with spacing $4\times 4$.}
  \label{tab:error_down_antonio_2d_44}
\end{table}

\begin{table}[H]
  \centering
  \begin{tabular}{||c| c  c c c  c c||}
    \hline
    $4\times 4$ & Cubic & DCT& DFT& Wavelet & LDMM (D) & LDMM (C)\\
    \hline
    $L_1$       &0.0302 &0.0310 &0.0314 &0.0326 &0.0249 & \textbf{0.0248}\\
    $L_2$       &0.0456 &0.0413 &0.0425 &0.0430 &0.0329 & \textbf{0.0329}\\
    $L_\infty$   &0.7629 &0.2411 &0.3776 &0.2514 &0.1779 &\textbf{0.1741}\\
    PSNR        &26.81  &27.68  &27.43  &27.34  &29.64 & \textbf{29.66}\\
    \hline
  \end{tabular}
  \caption{Errors of the interpolation of the 2D plasma (distribution function) data set from regular sampling with spacing $4\times 4$.}
  \label{tab:error_down_shock_2d_44}
\end{table}

\begin{table}[H]
  \centering
  \begin{tabular}{||c| c  c c c c c||}
    \hline
    $4\times 4$ & Cubic & DCT& DFT& Wavelet & LDMM (D) & LDMM (C)\\
    \hline
    $L_1$       &\textbf{0.0009} &0.0015 &0.0016 &0.0020 &0.0013 & 0.0012\\
    $L_2$       &0.0045 &0.0051 &0.0055 &0.0061 &0.0044 & \textbf{0.0041}\\
    $L_\infty$   &0.1461 &0.1547 &0.2202 &0.1892 &0.1393 & \textbf{0.1278}\\
    PSNR        &46.97  &45.77  &45.20  &44.31  &47.18 & \textbf{47.43}\\
    \hline
  \end{tabular}
  \caption{Errors of the interpolation of the 2D lattice data set from regular sampling with spacing $4\times 4$.}
  \label{tab:error_down_lattice_2d_44}
\end{table}

\begin{table}[H]
  \centering
  \begin{tabular}{||c| c  c c c c c||}
    \hline
    $4\times 4\times 1$ & Cubic & DCT& DFT& Wavelet & LDMM (D) & LDMM (C)\\
    \hline
    $L_1$       &0.0038 &0.0040 &0.0071 &0.0052 &0.0037 & \textbf{0.0036}\\
    $L_2$       &0.0085 &0.0065 &0.0200 &0.0078 &\textbf{0.0059} & 0.0065\\
    $L_\infty$   &0.9649 &0.1366 &0.6449 &0.1357 &\textbf{0.1259} & 0.1911\\
    PSNR        &41.38  &43.76  &33.99  &42.15  &\textbf{44.53}  & 43.73\\
    \hline
    $2\times 2\times 2$ & Cubic & DCT& DFT& Wavelet & LDMM (D) & LDMM (C)\\
    \hline
    $L_1$       &0.0356 &0.0334 &0.0352 &0.0439 &\textbf{0.0305} & 0.0313\\
    $L_2$       &0.0714 &0.0593 &0.0632 &0.0613 &\textbf{0.0535} & 0.0559\\
    $L_\infty$   &0.8770 &0.4073 &0.5203 &0.4283 &\textbf{0.3711} & 0.4060\\
    PSNR        &22.93  &24.54  &23.99  &24.25  &\textbf{25.43}  & 25.05\\
    \hline
  \end{tabular}
  \caption{Errors of the interpolation of the 3D plasma (magnetic field) data set from regular sampling with spacing $4\times 4 \times 1$ and $2 \times 2 \times 2$.}
  \label{tab:error_down_plasma_3d}
\end{table}

\begin{table}[H]
  \centering
  \begin{tabular}{||c| c  c c c c c||}
    \hline
    $4\times 4\times 1$ & Cubic & DCT& DFT& Wavelet & LDMM (D) & LDMM (C)\\
    \hline
    $L_1$       &0.0094 &0.0072 &0.0168 &0.0066 &0.0058 & \textbf{0.0056}\\
    $L_2$       &0.0593 &0.0292 &0.0434 &0.0281 &\textbf{0.0233} & 0.0254\\
    $L_\infty$   &1.1890 &0.4223 &0.5405 &0.4245 &\textbf{0.4164} & 0.4362\\
    PSNR        &24.54  &30.69  &27.25  &31.03  &\textbf{32.64} & 31.90\\
    \hline
    $2\times 2\times 2$ & Cubic & DCT& DFT& Wavelet & LDMM (D) & LDMM (C)\\
    \hline
    $L_1$       &0.0039 &0.0027 &0.0069 &0.0045 &0.0017 & \textbf{0.0015}\\
    $L_2$       &0.0316 &0.0119 &0.0237 &0.0124 &\textbf{0.0101} & 0.0101\\
    $L_\infty$   &0.7459 &0.4109 &0.4282 &0.4233 &\textbf{0.4078} & 0.4096\\
    PSNR        &30.01 &38.49  &32.51  &38.15  &\textbf{39.93} & 39.92\\
    \hline
  \end{tabular}
  \caption{Errors of the interpolation of the 3D lattice data set from regular sampling with spacing $4\times 4 \times 1$ and $2 \times 2 \times 2$.}
  \label{tab:error_down_lattice_3d}
\end{table}

\begin{table}[H]
  \centering
  \begin{tabular}{||c| c  c c c c c||}
    \hline
    $4\times 4\times 1$ & Cubic & DCT& DFT& Wavelet & LDMM (D) & LDMM (C)\\
    \hline
    $L_1$       &0.0076 &0.0078 &0.0103 &0.0083 &\textbf{0.0064} & 0.0064\\
    $L_2$       &0.0150 &0.0136 &0.0238 &0.0141 &\textbf{0.0110} & 0.0111\\
    $L_\infty$   &0.8851 &0.1551 &0.4805 &0.1469 &\textbf{0.1093} & 0.1417\\
    PSNR        &36.47  &37.35  &32.45  &37.02  &\textbf{39.18} & 39.13\\
    \hline
    $2\times 2\times 2$ & Cubic & DCT& DFT& Wavelet & LDMM (D) & LDMM (C)\\
    \hline
    $L_1$       &0.0109 &0.0098 &0.0127 &0.0139 &\textbf{0.0089} & 0.0092\\
    $L_2$       &0.0283 &0.0202 &0.0255 &0.0229 &\textbf{0.0178} & 0.0181\\
    $L_\infty$   &0.7388 &0.2976 &0.3438 &0.2993 &\textbf{0.2088} & 0.2097\\
    PSNR        &30.97  &33.91  &31.88  &32.81  &\textbf{35.01}  & 34.85\\
    \hline
  \end{tabular}
  \caption{Errors of the interpolation of the 3D plasma (distribution function) data set from regular sampling with spacing $4\times 4 \times 1$ and $2 \times 2 \times 2$.}
  \label{tab:error_down_shock_3d}
\end{table}

\subsection{Data Compression}
Finally, we compare the performance of LDMM as a sampling-based data compression technique to other standard compression methods including singular value/ Tucker Decomposition, DFT, DCT, and the wavelet transformations. We point out that, unlike the other testing methods which usually involve hard thresholding of the expansion coefficients with respect to a particular basis, LDMM does not require access to the original full data set. Therefore we do not expect LDMM to perform equally well compared to other data compression methods. However, using the sampling-based method as a data compression technique has its own advantages:
\begin{itemize}
\item During the data compression step, sampling-based algorithms like LDMM are very easy to implement compared to other standard compression methods. Moreover, in a parallel setting, sampling based methods can be implemented independently on each node without communication, while other methods involving global transforms cannot.
\item It is also faster for sampling-based methods to reconstruct a small portion of the data set if only that part of the data set is required.
\end{itemize}

\begin{figure}[H]
  \centering
  \begin{tabular}{ccc}
    Original& SVD (33.65dB)& DCT (\textbf{66.87dB})\\
    \includegraphics[width=3cm]{antonio_2d_original}&
    \includegraphics[width=3cm]{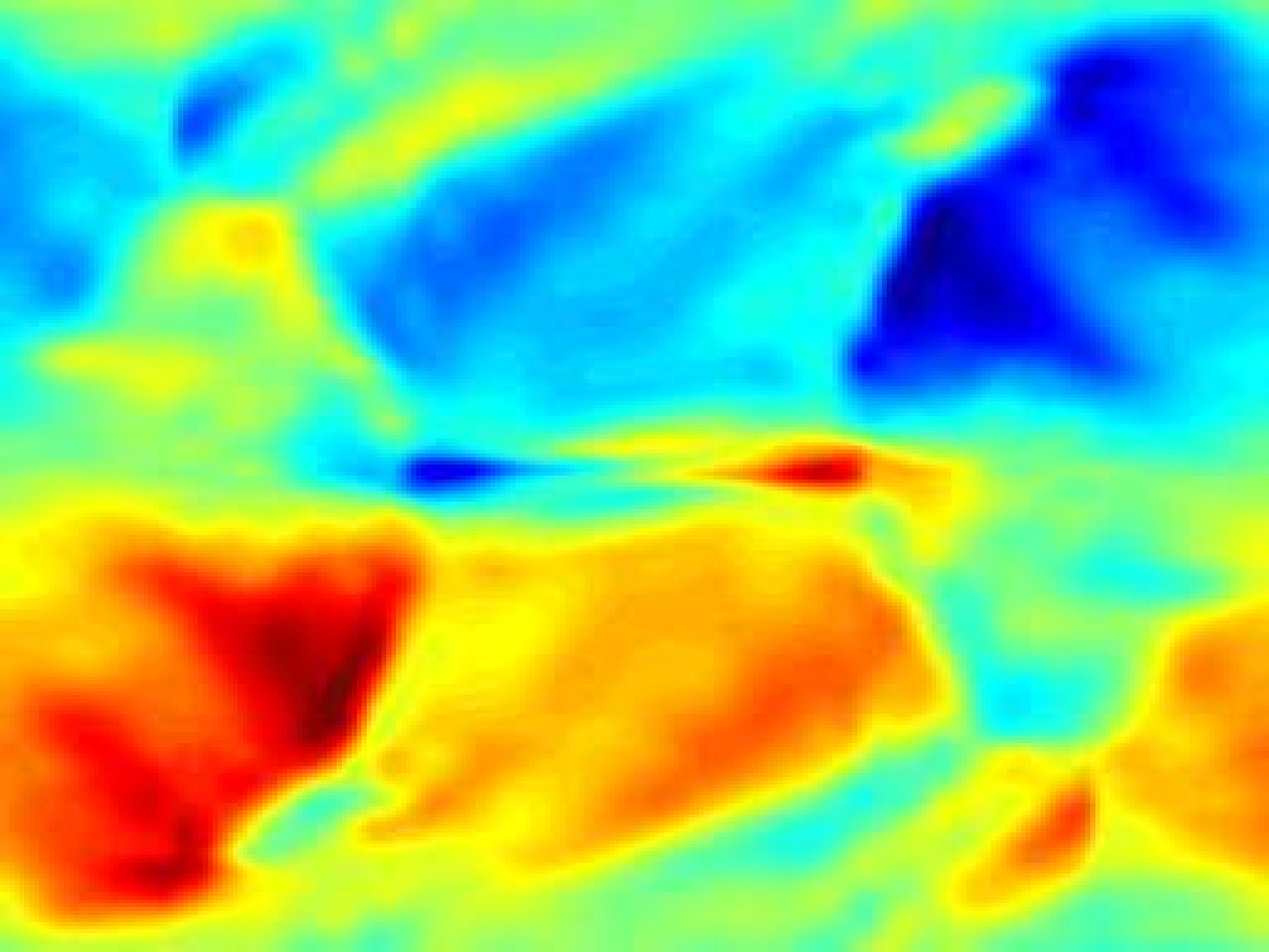}&
    \includegraphics[width=3cm]{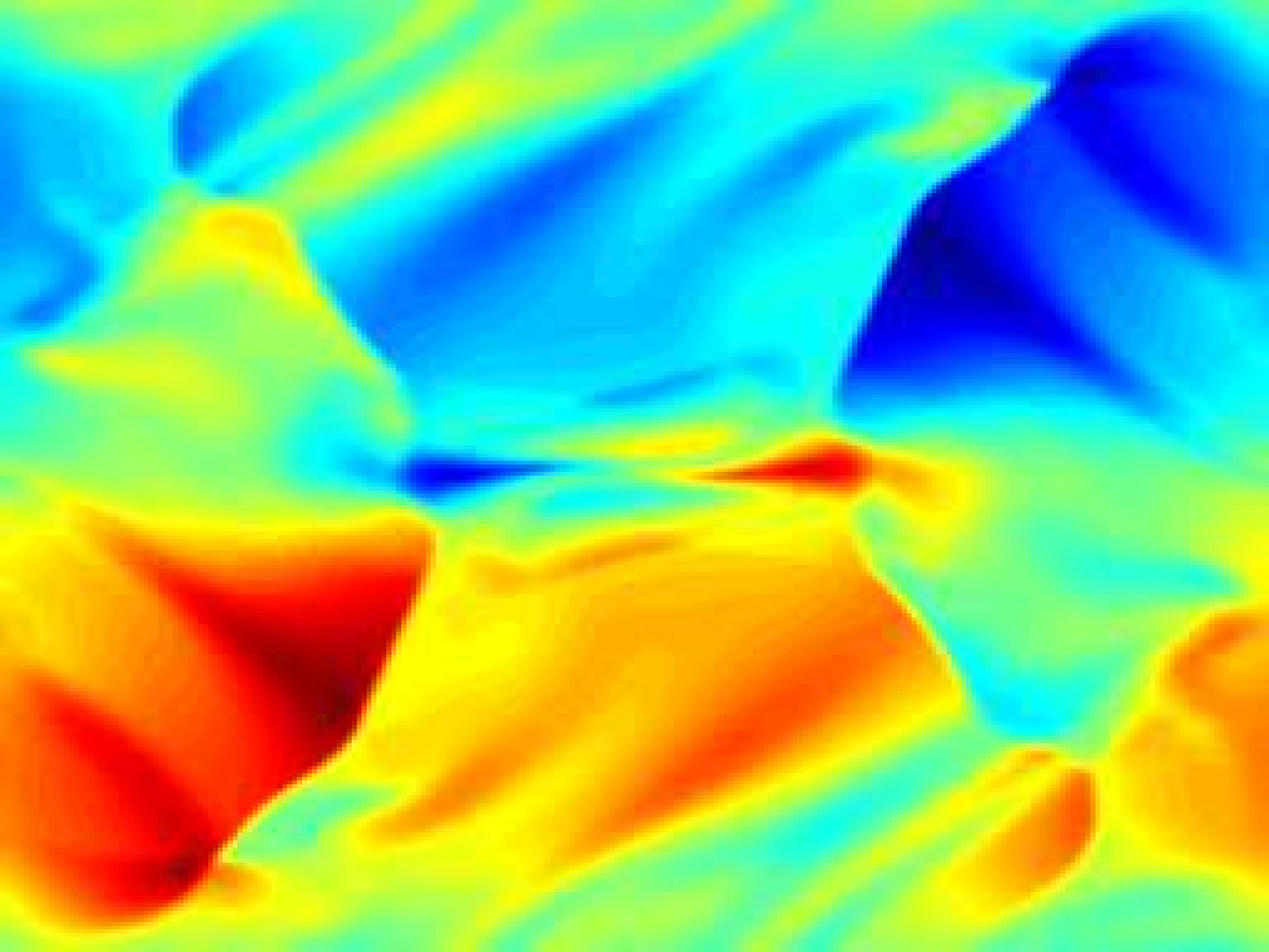}\\
    DFT (57.03dB)& Wavelet (63.01dB)& LDMM (42.55dB)\\    \vspace{.5cm}
    \includegraphics[width=3cm]{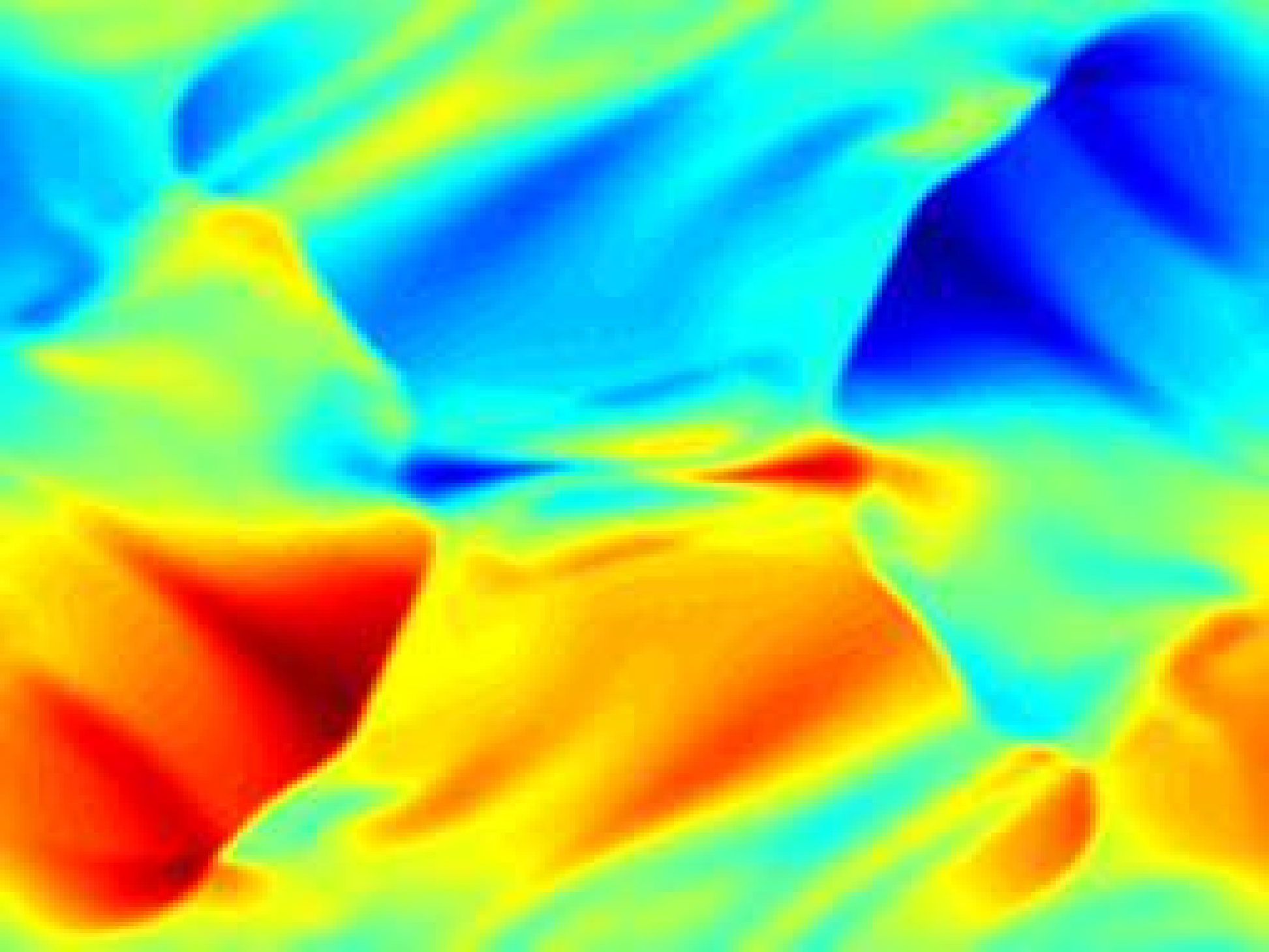}&
    \includegraphics[width=3cm]{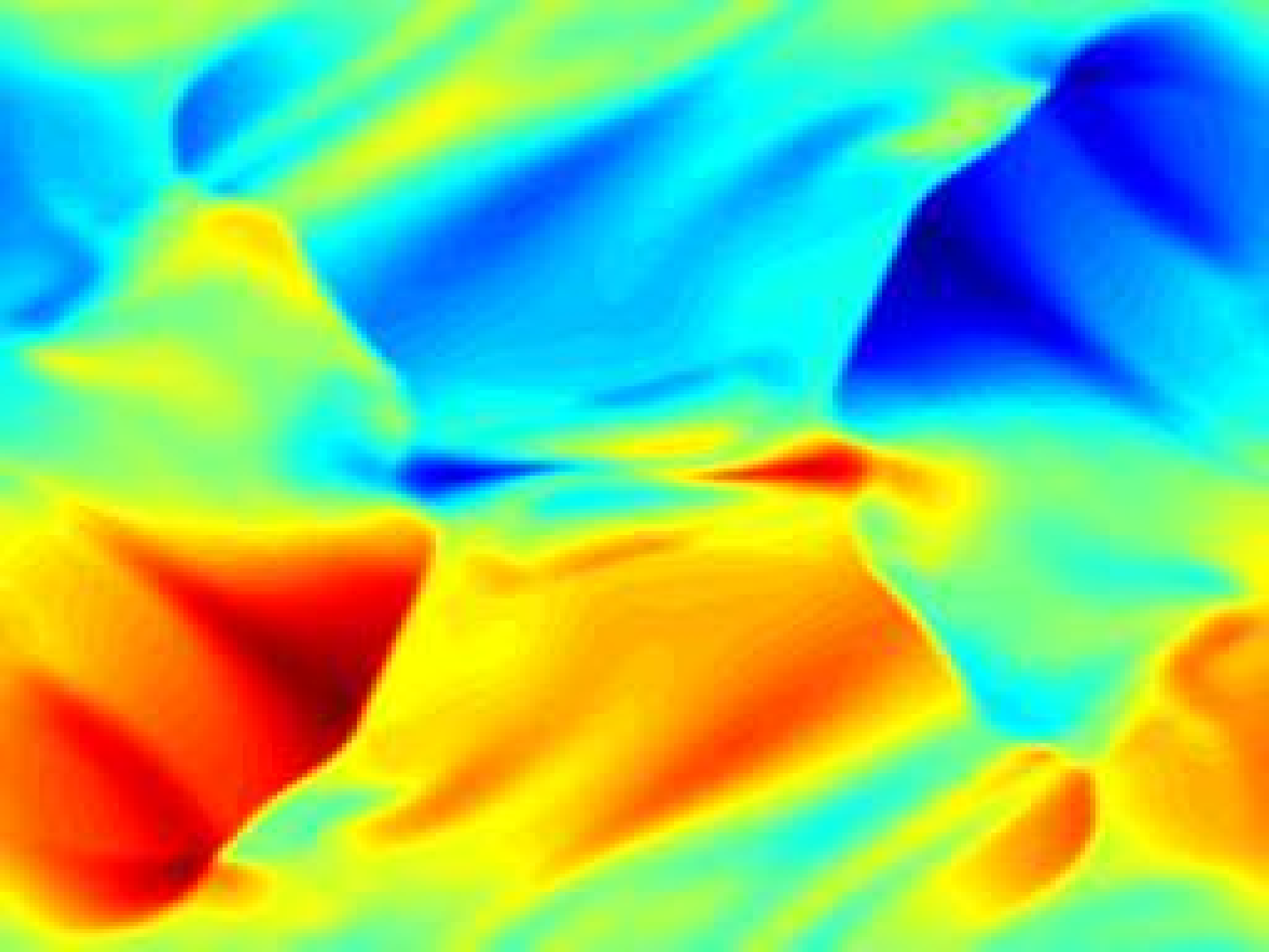}&
    \includegraphics[width=3cm]{antonio_2d_ldmm_10p}\\
    Original& SVD (27.19dB)& DCT (\textbf{36.63dB})\\
    \includegraphics[width=3cm]{shock_2d_original}&
    \includegraphics[width=3cm]{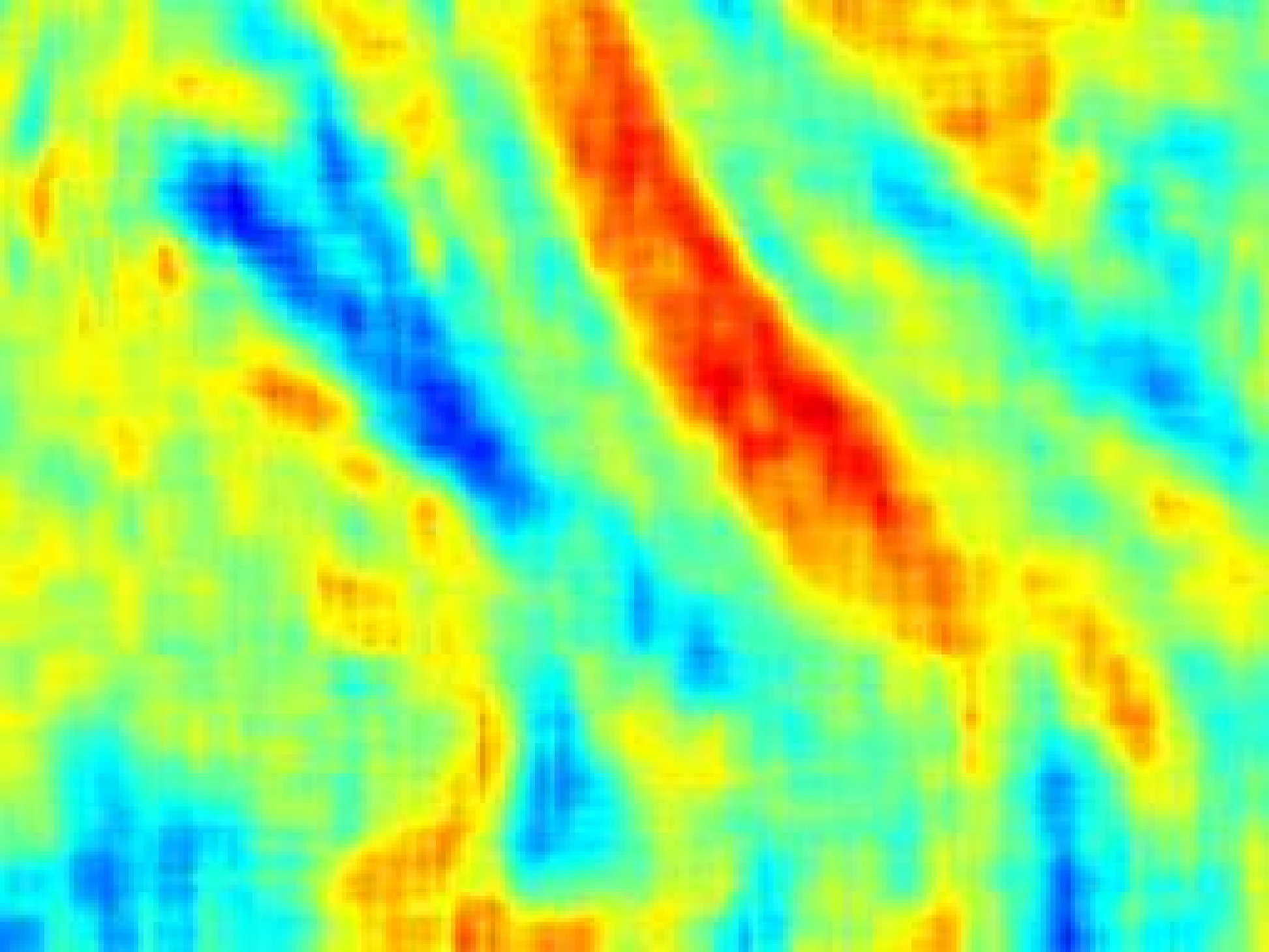}&
    \includegraphics[width=3cm]{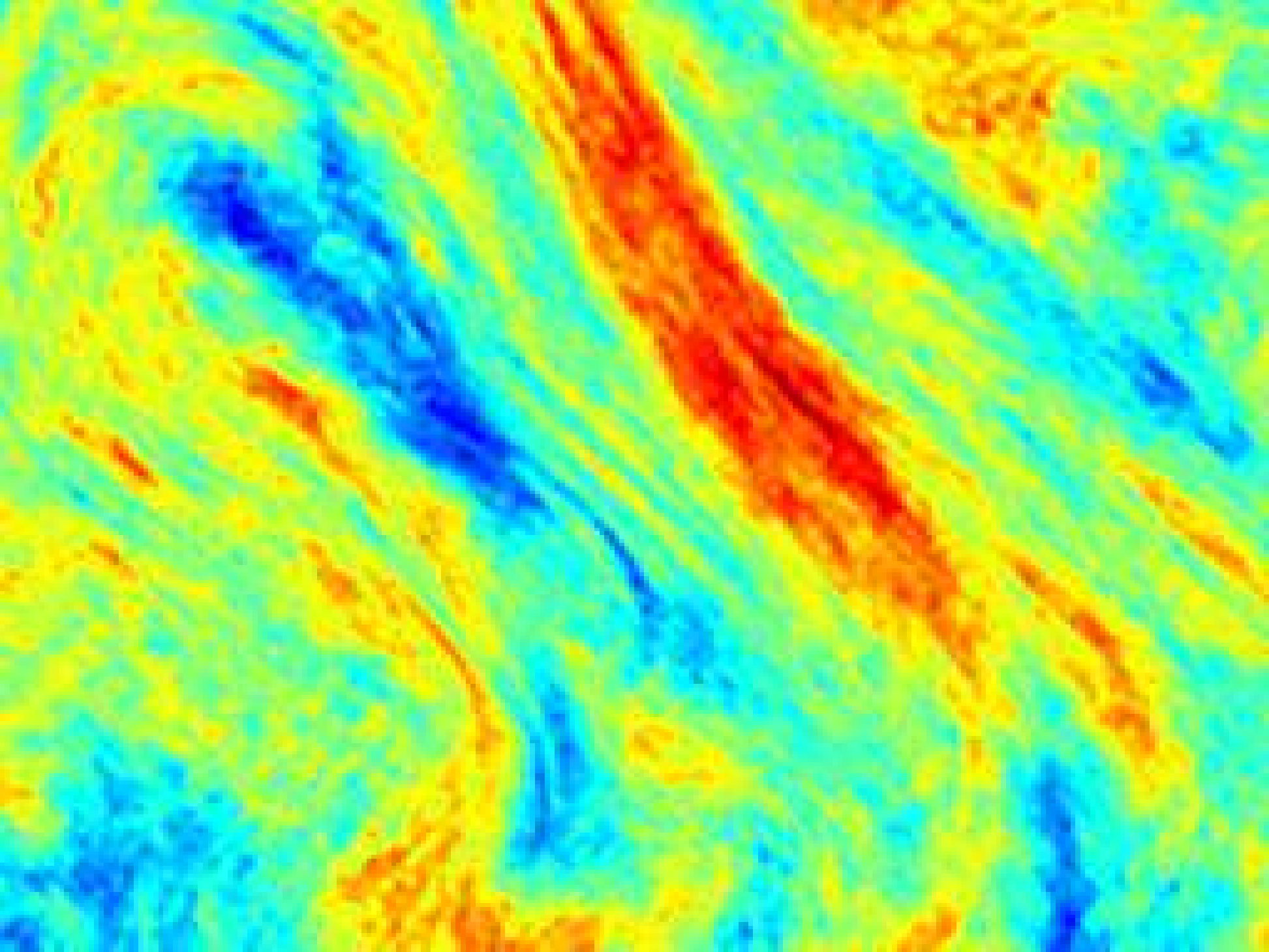}\\
    DFT (34.44dB)& Wavelet (34.78dB)& LDMM (29.56dB)\\    \vspace{.5cm}
    \includegraphics[width=3cm]{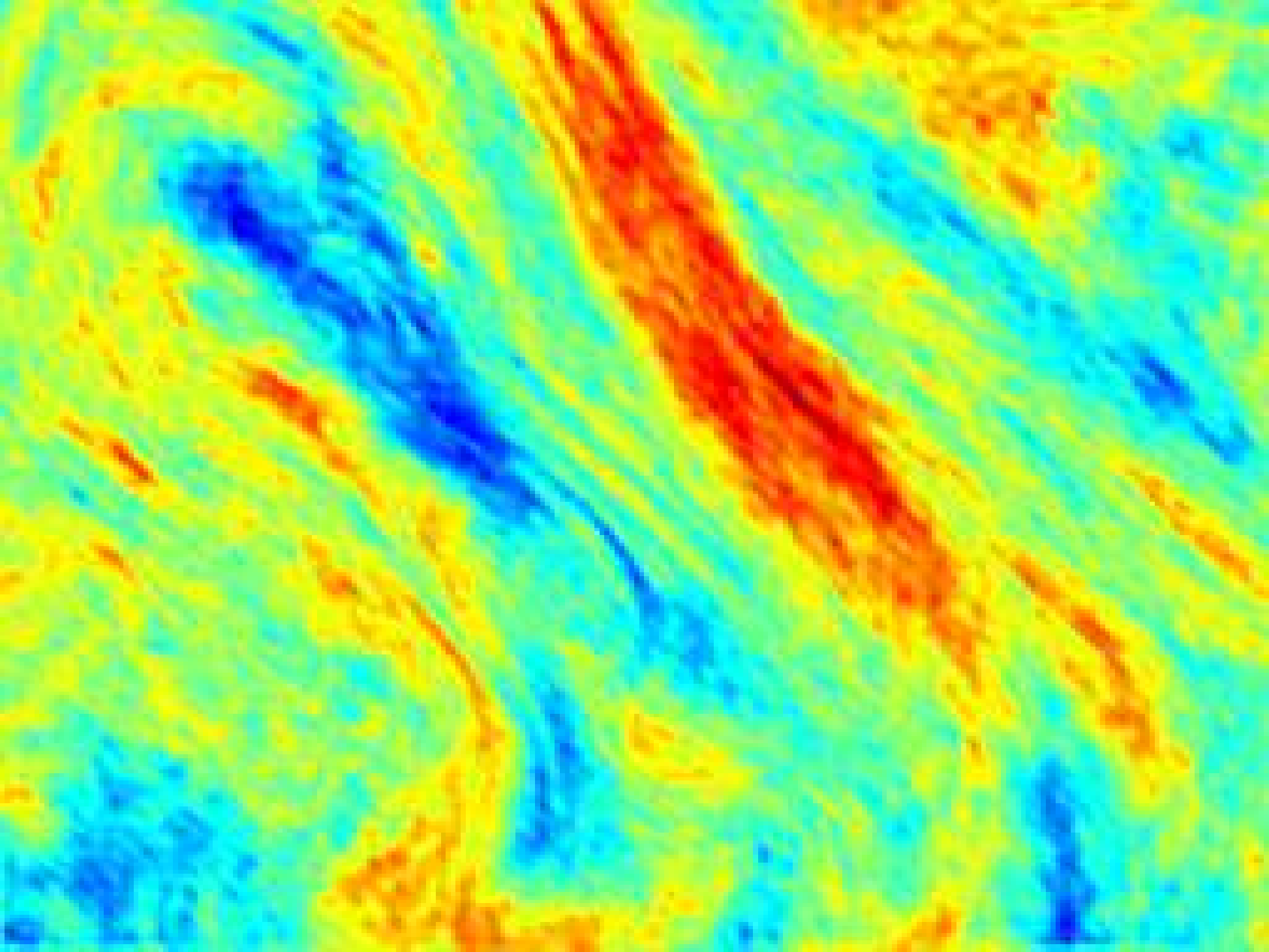}&
    \includegraphics[width=3cm]{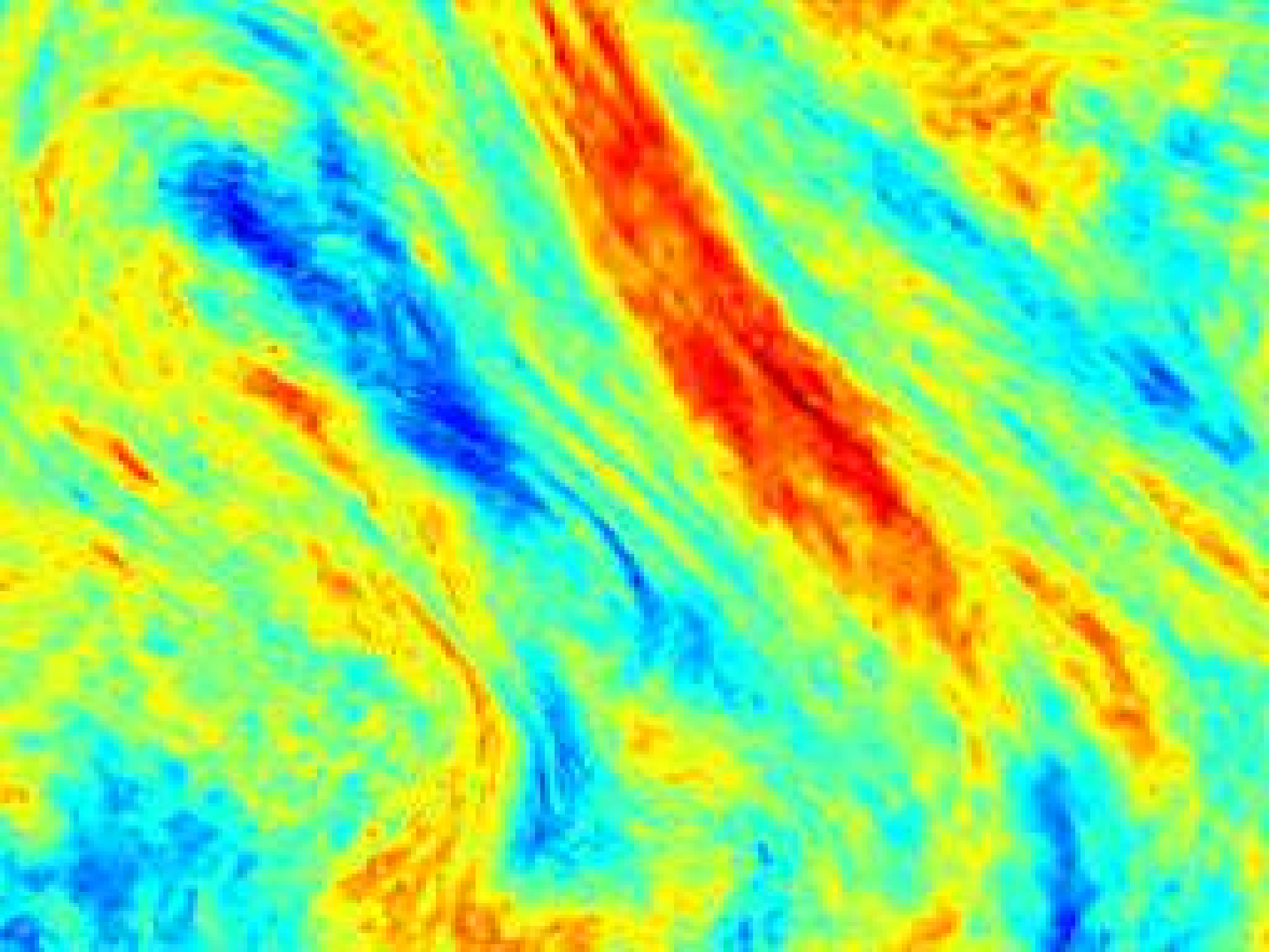}&
    \includegraphics[width=3cm]{shock_2d_ldmm_10p}\\
    Original& SVD (57.24dB)& DCT (75.73dB)\\
    \includegraphics[width=3cm]{latticebig_2d_original}&
    \includegraphics[width=3cm]{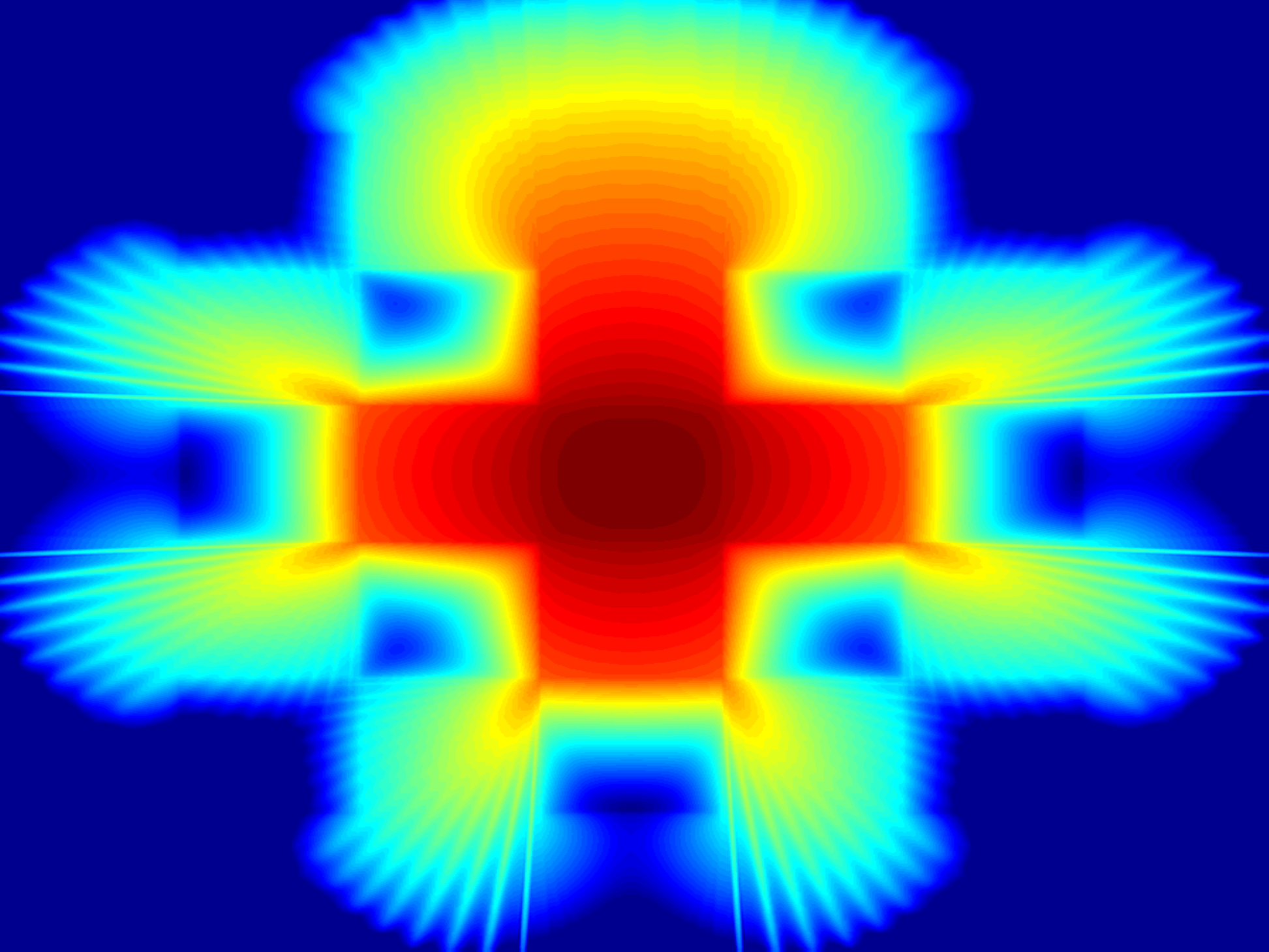}&
    \includegraphics[width=3cm]{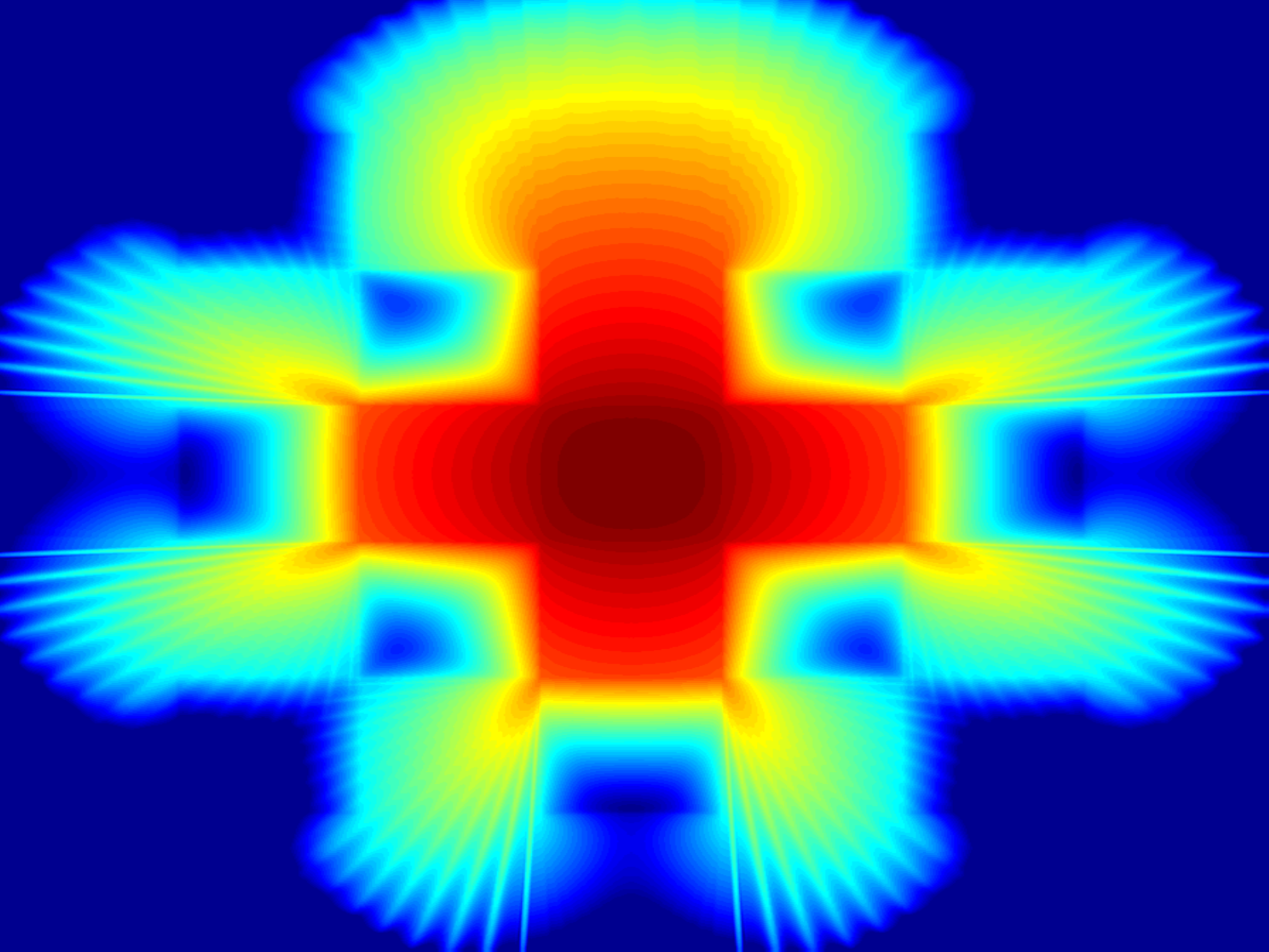}\\
    DFT (61.94dB)& Wavelet (\textbf{80.33dB})& LDMM (47.98dB)\\
    \includegraphics[width=3cm]{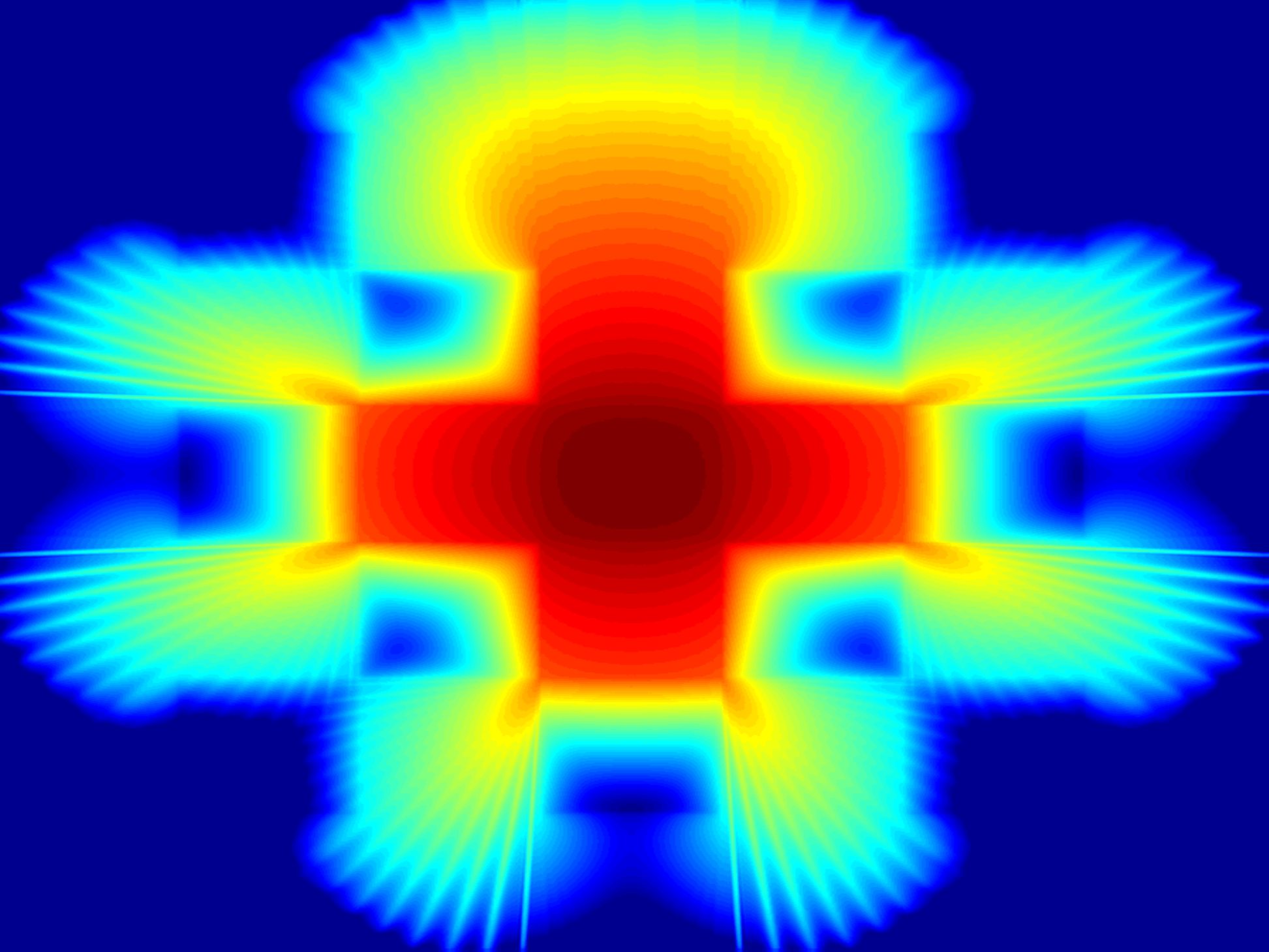}&
    \includegraphics[width=3cm]{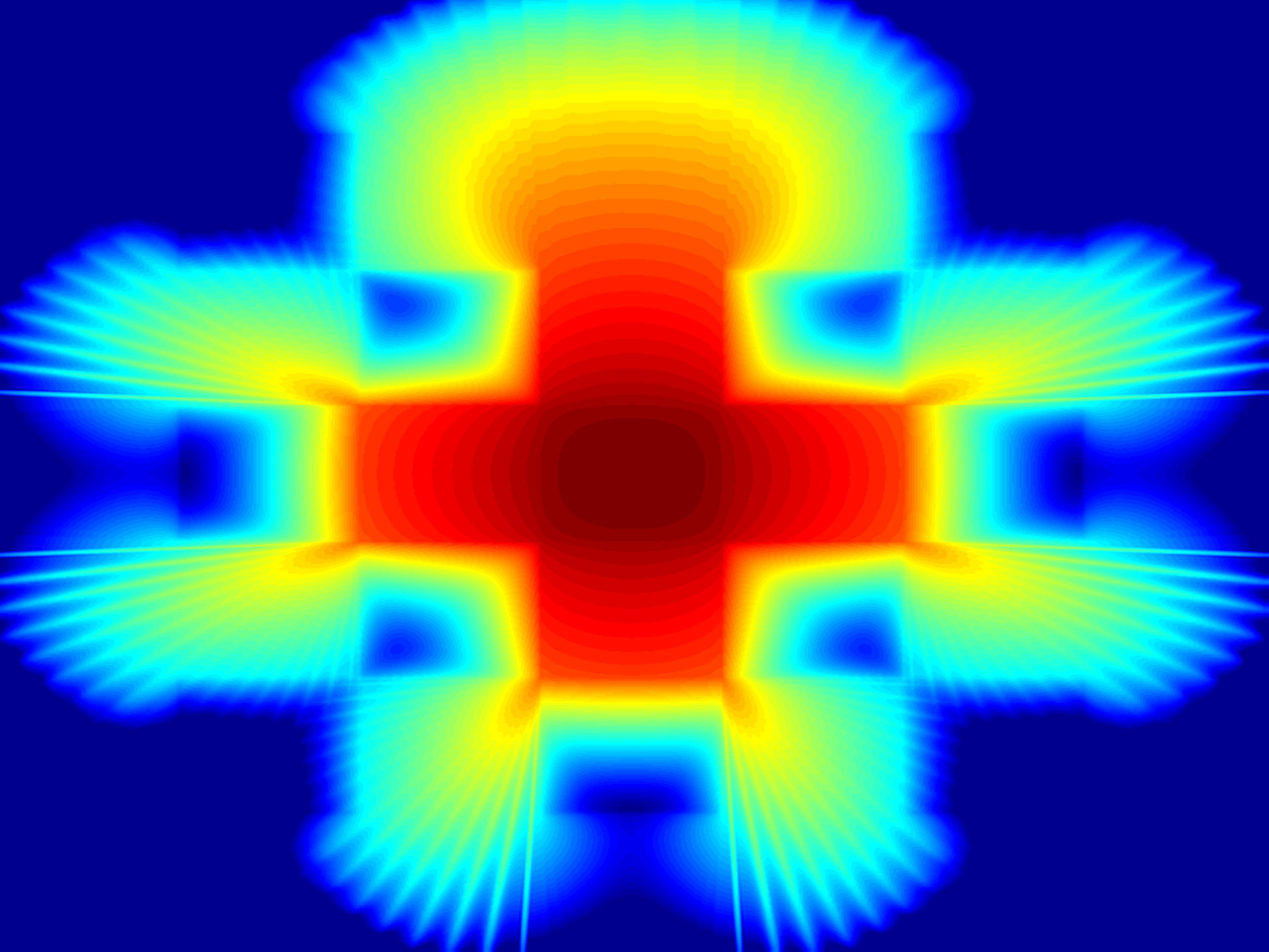}&
    \includegraphics[width=3cm]{latticebig_2d_ldmm_10p}
  \end{tabular}
  \caption{Compression of 2D scientific data sets with a 10\% data compression rate. The original data are shown on the upper left corners for each data set. The results of SVD, DCT, DFT, wavelet, and LDMM are shown in the remaining five figures.}
  \label{fig:compre_2d_10p}
\end{figure}

\begin{figure}[H]
  \centering
  \begin{tabular}{ccc}
    Original& SVD (28.29dB)& DCT (56.36dB)\\
    \includegraphics[width=3cm]{antonio_2d_original}&
    \includegraphics[width=3cm]{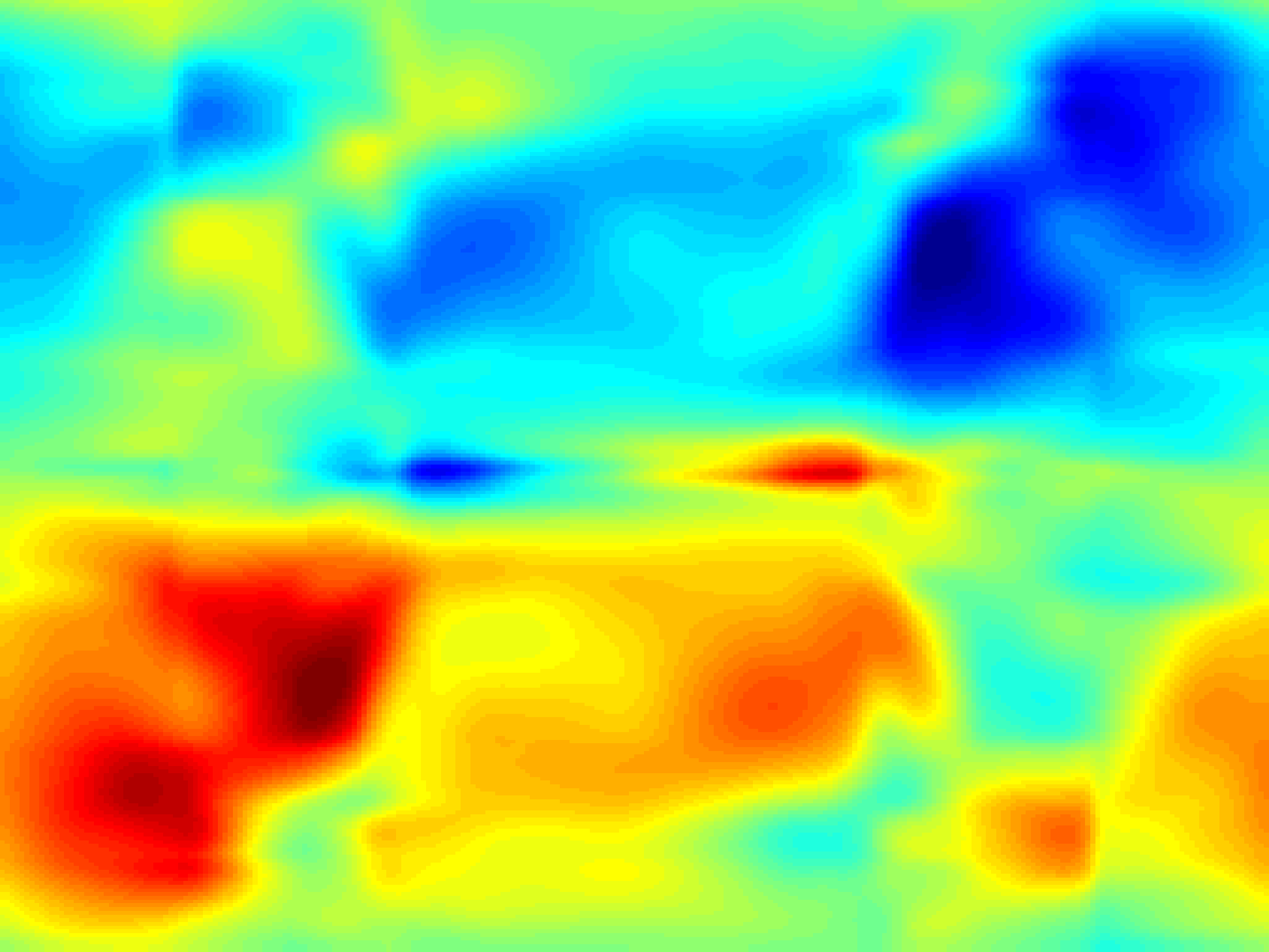}&
    \includegraphics[width=3cm]{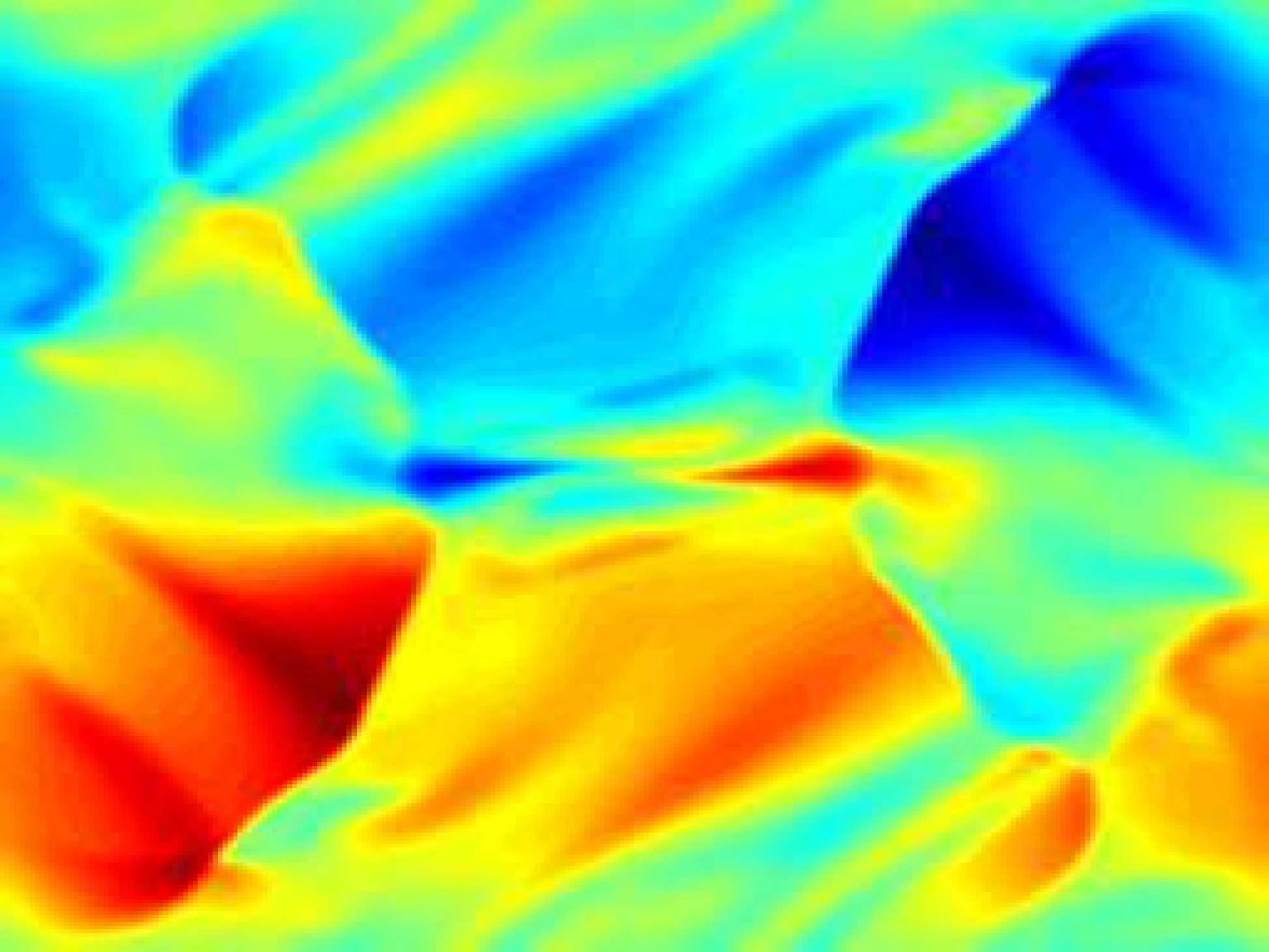}\\
    DFT (49.12dB)& Wavelet (\textbf{53.36dB})& LDMM (39.09dB)\\    \vspace{.5cm}
    \includegraphics[width=3cm]{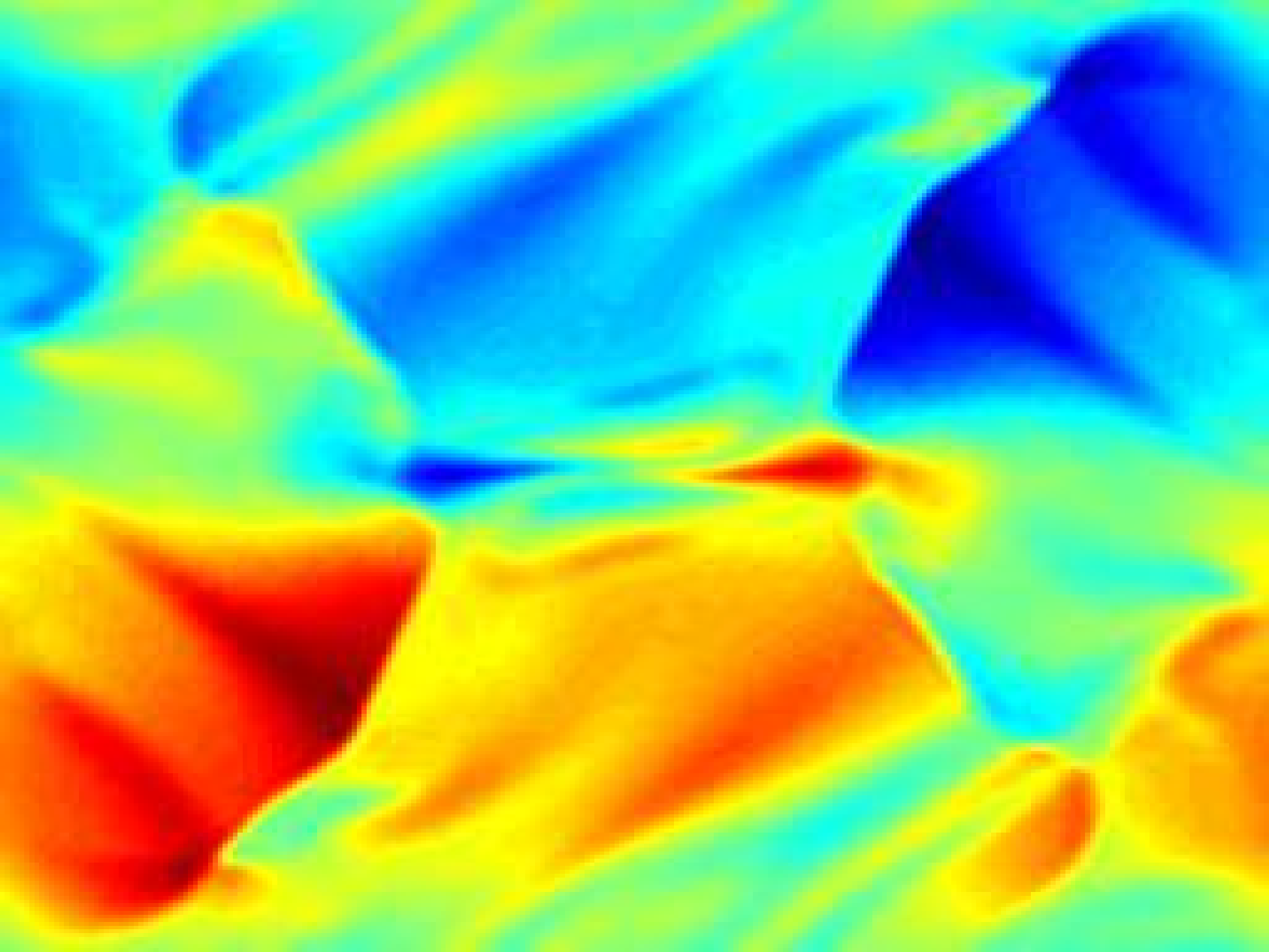}&
    \includegraphics[width=3cm]{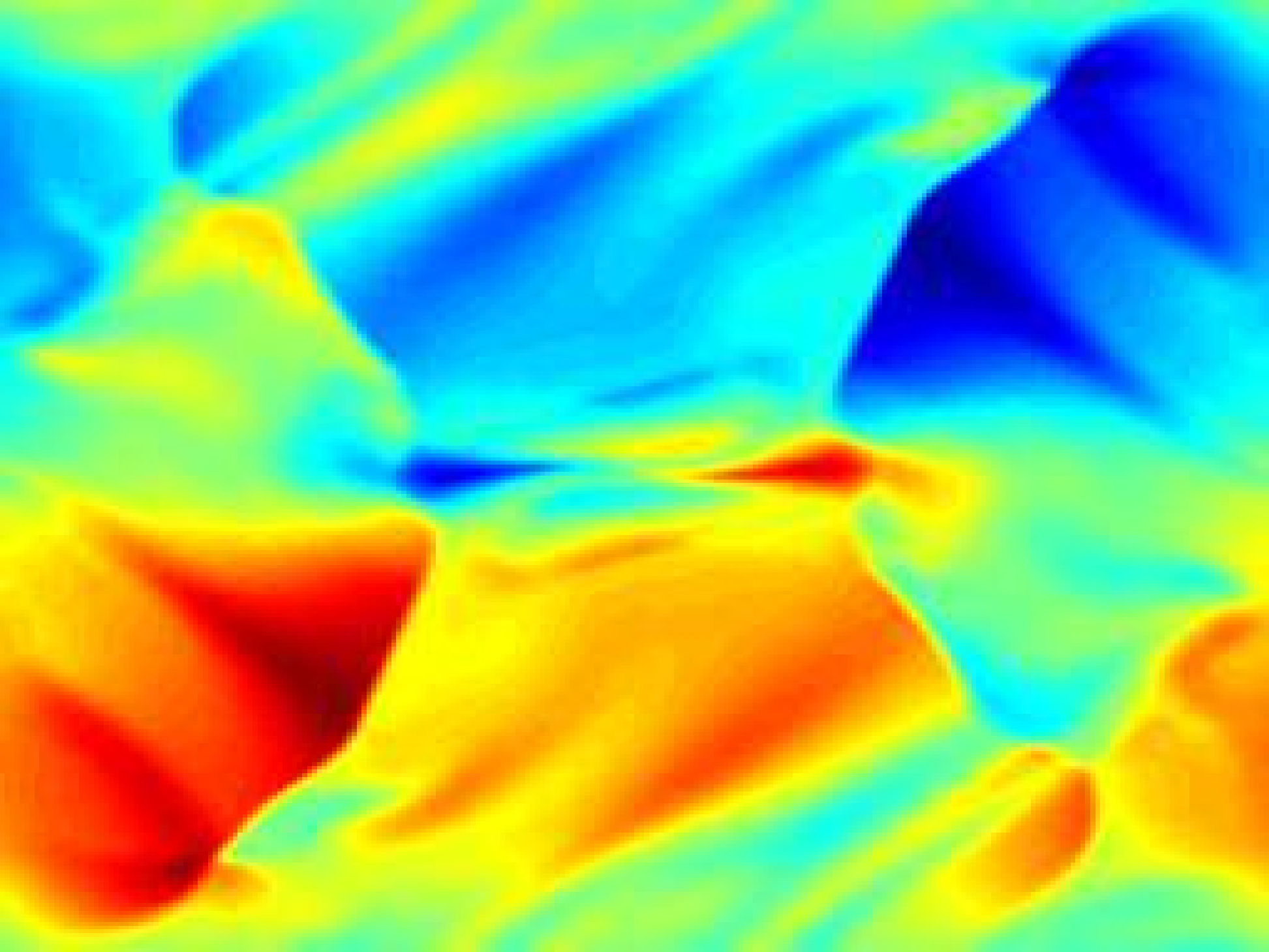}&
    \includegraphics[width=3cm]{antonio_2d_ldmm_5p}\\
    Original& SVD (24.10dB)& DCT (\textbf{32.47dB})\\
    \includegraphics[width=3cm]{shock_2d_original}&
    \includegraphics[width=3cm]{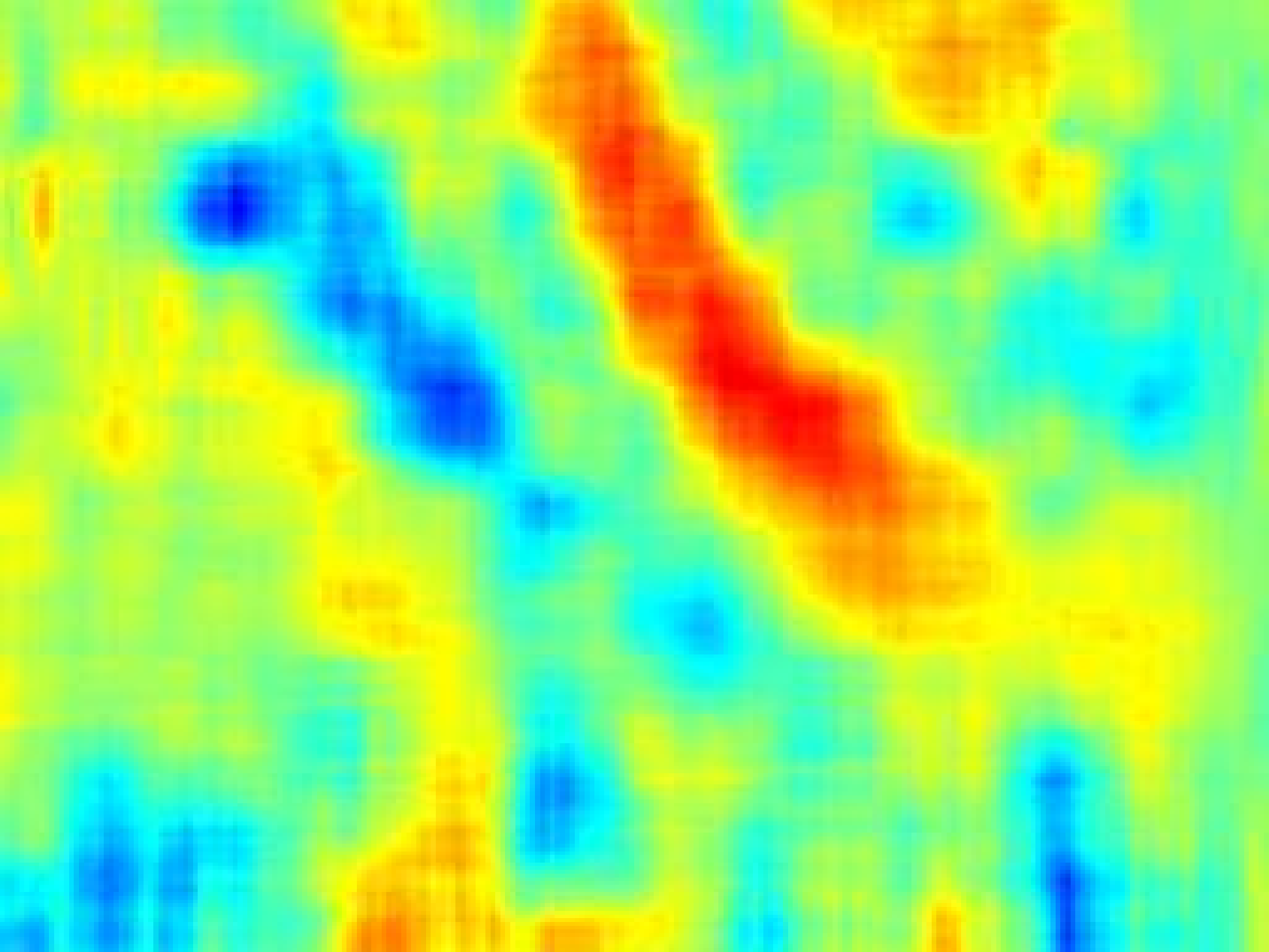}&
    \includegraphics[width=3cm]{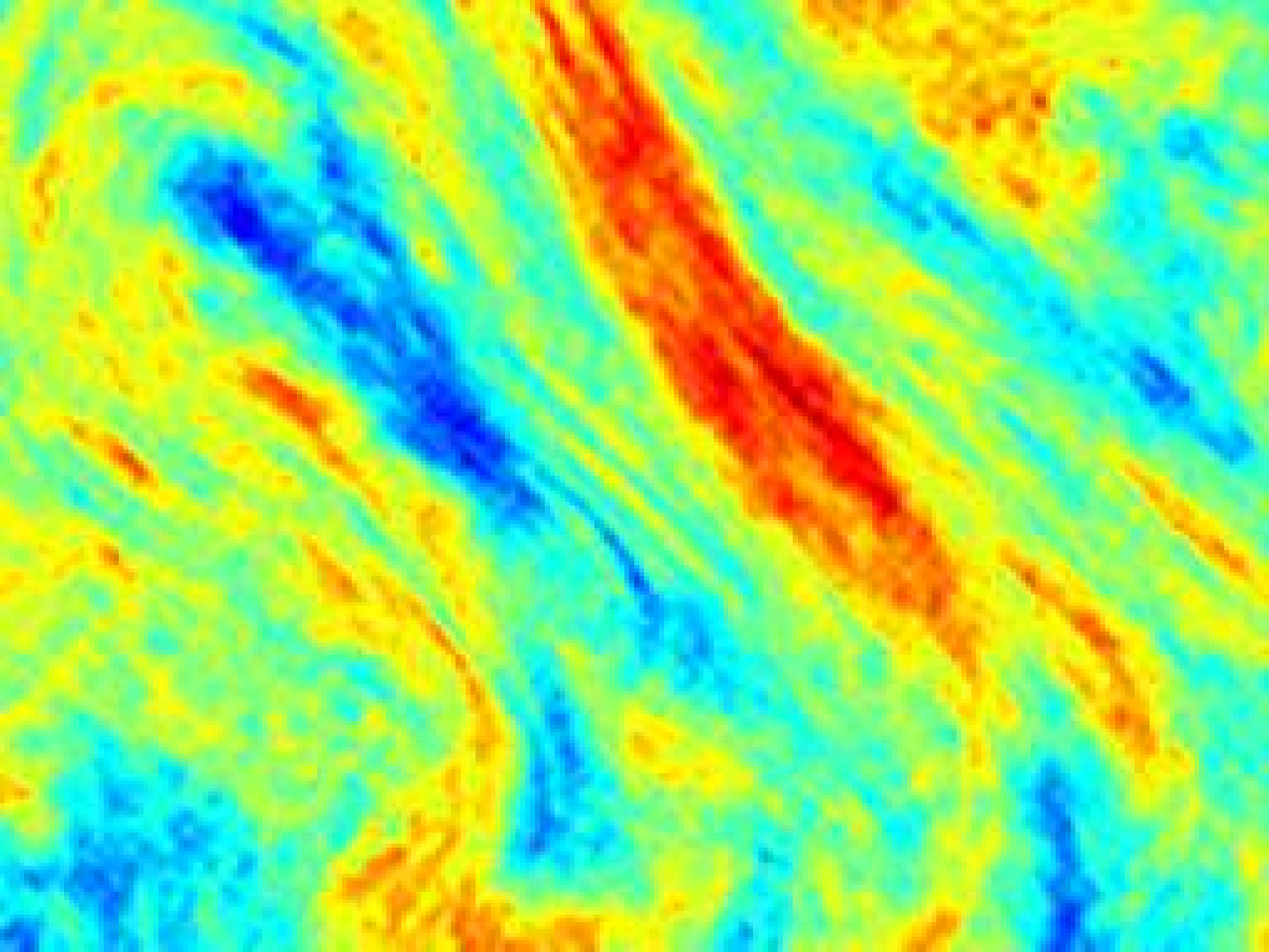}\\
    DFT (31.59dB)& Wavelet (31.94dB)& LDMM (27.93dB)\\    \vspace{.5cm}
    \includegraphics[width=3cm]{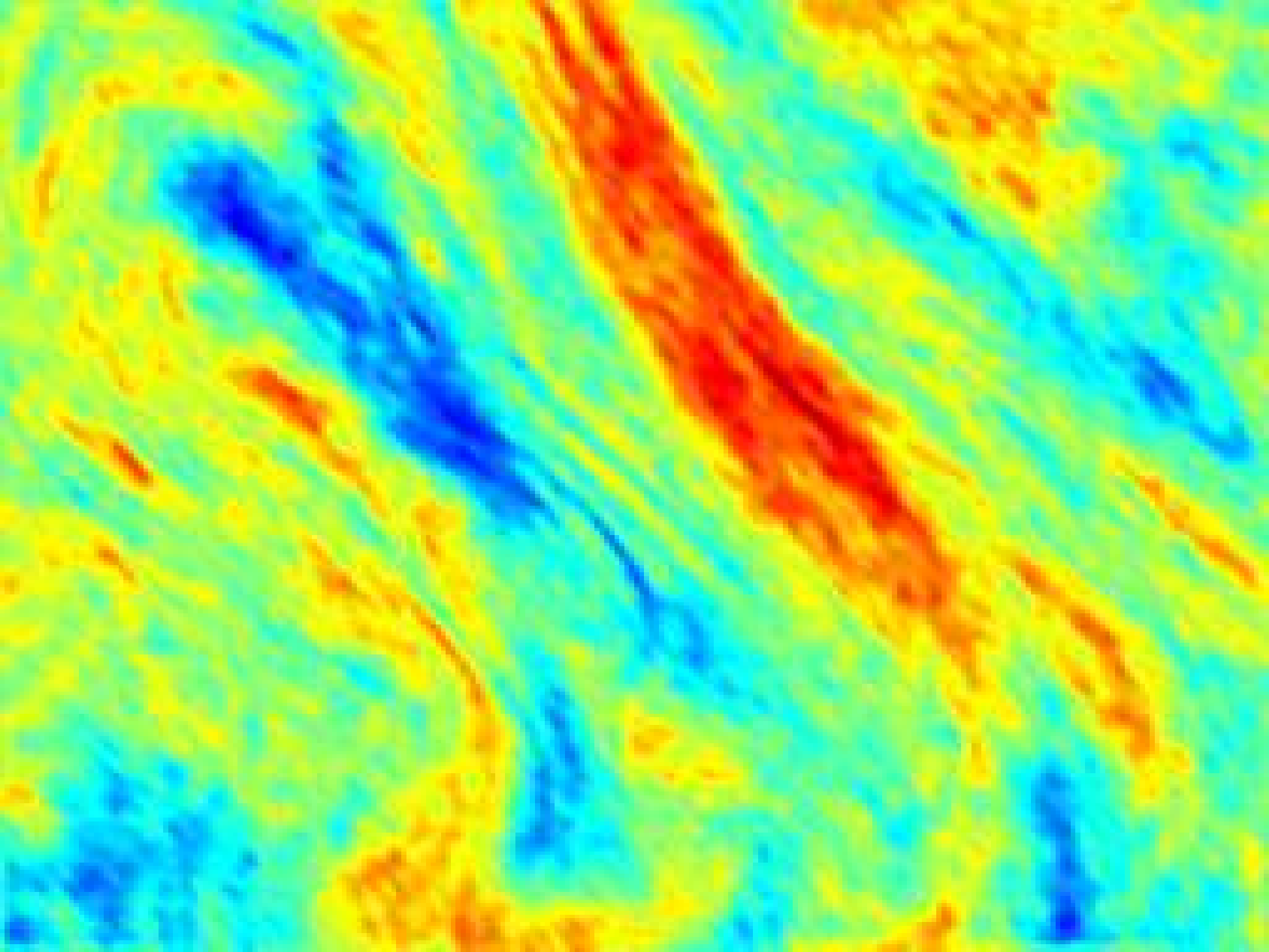}&
    \includegraphics[width=3cm]{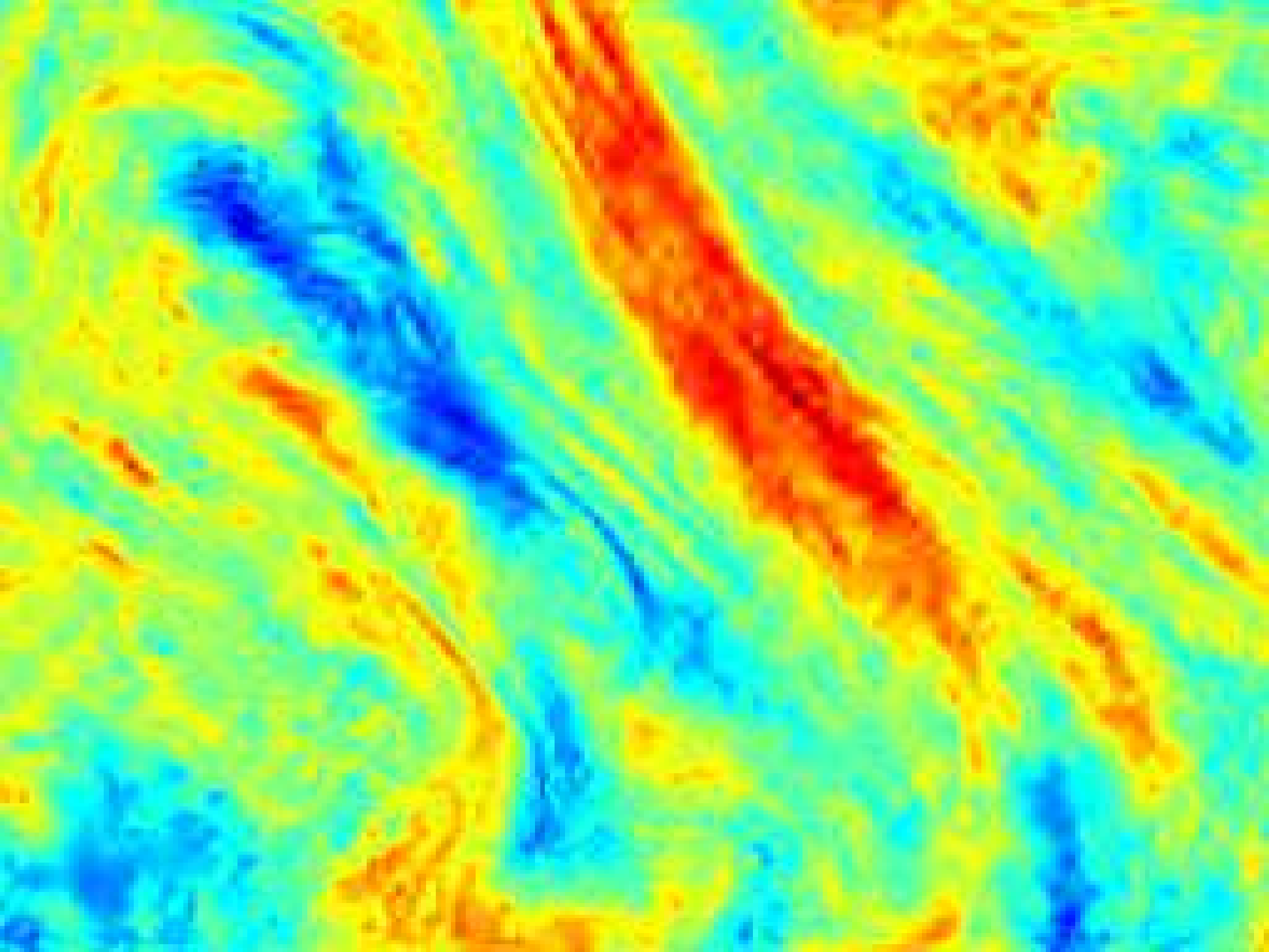}&
    \includegraphics[width=3cm]{shock_2d_ldmm_5p}\\
    Original& SVD (47.17dB)& DCT (67.59dB)\\
    \includegraphics[width=3cm]{latticebig_2d_original}&
    \includegraphics[width=3cm]{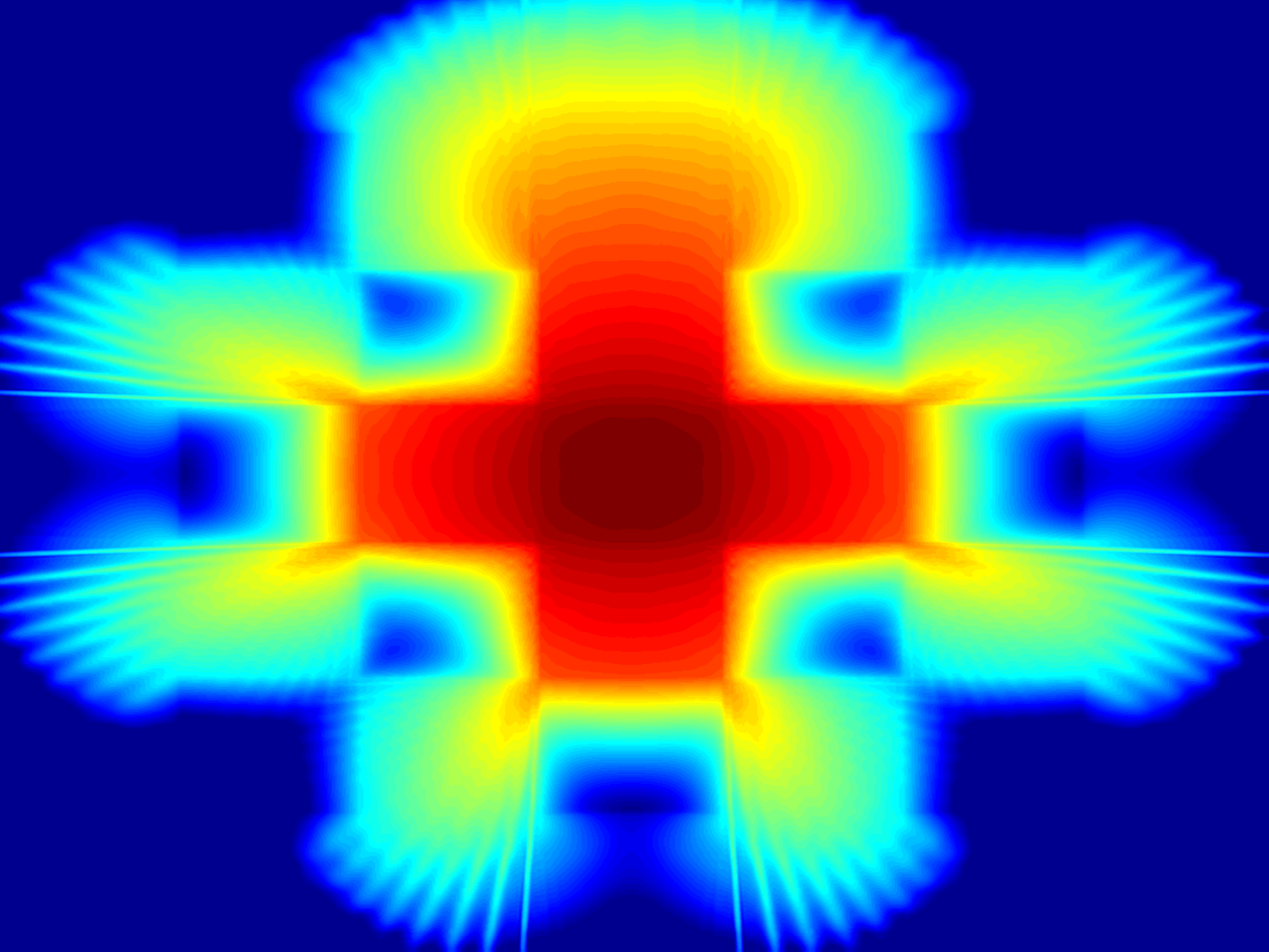}&
    \includegraphics[width=3cm]{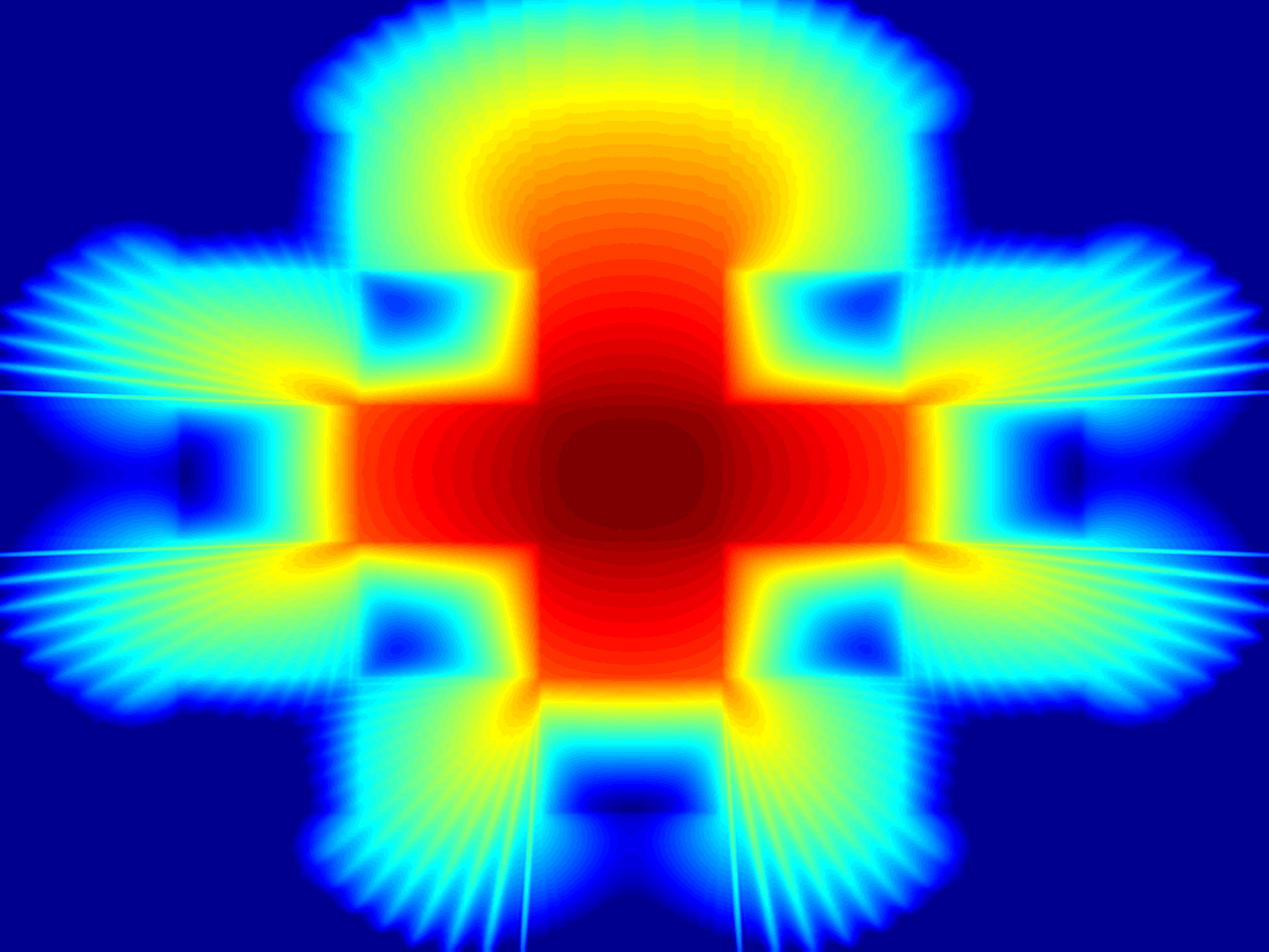}\\
    DFT (55.49dB)& Wavelet (\textbf{68.34dB})& LDMM (44.15dB)\\
    \includegraphics[width=3cm]{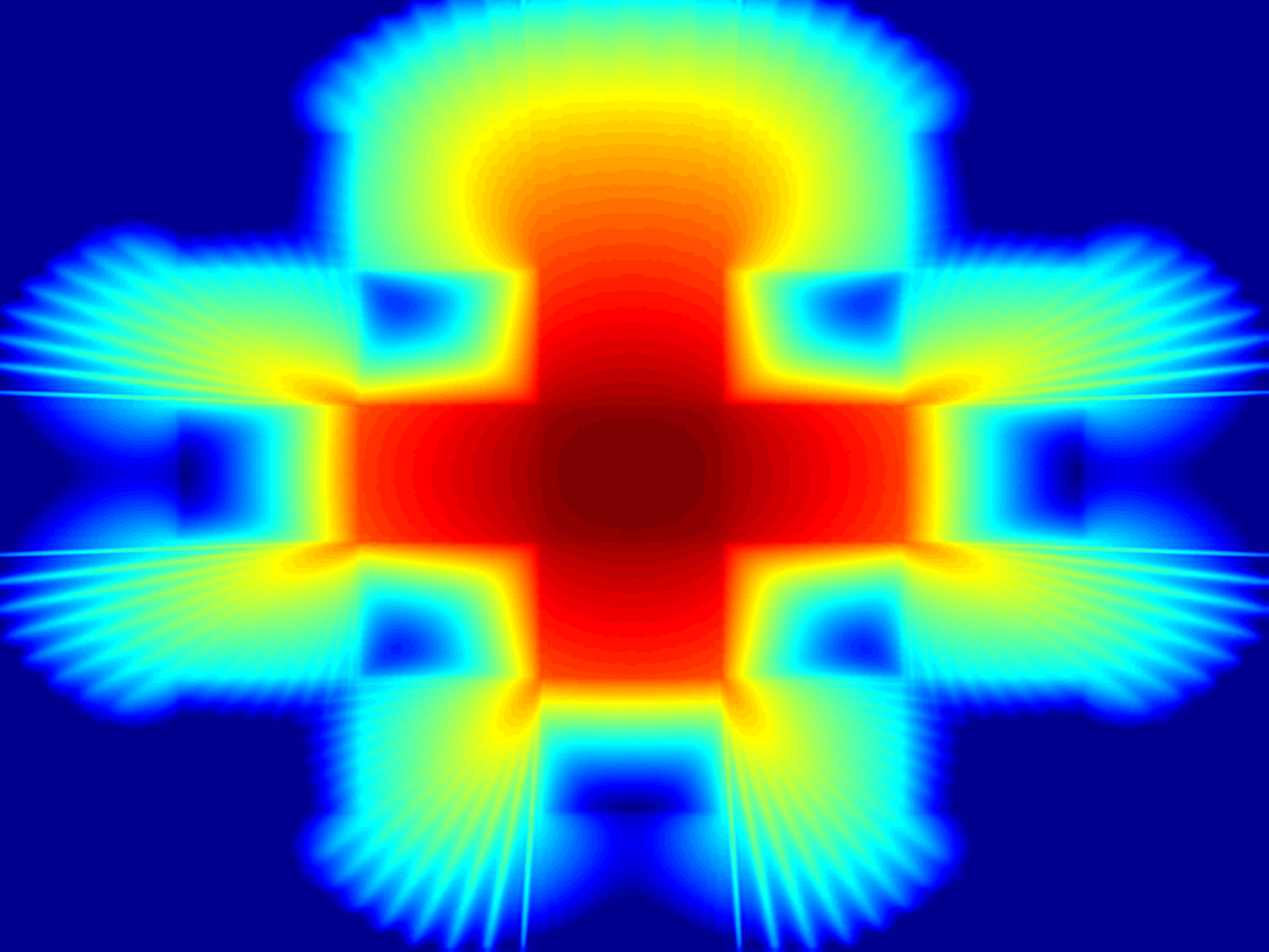}&
    \includegraphics[width=3cm]{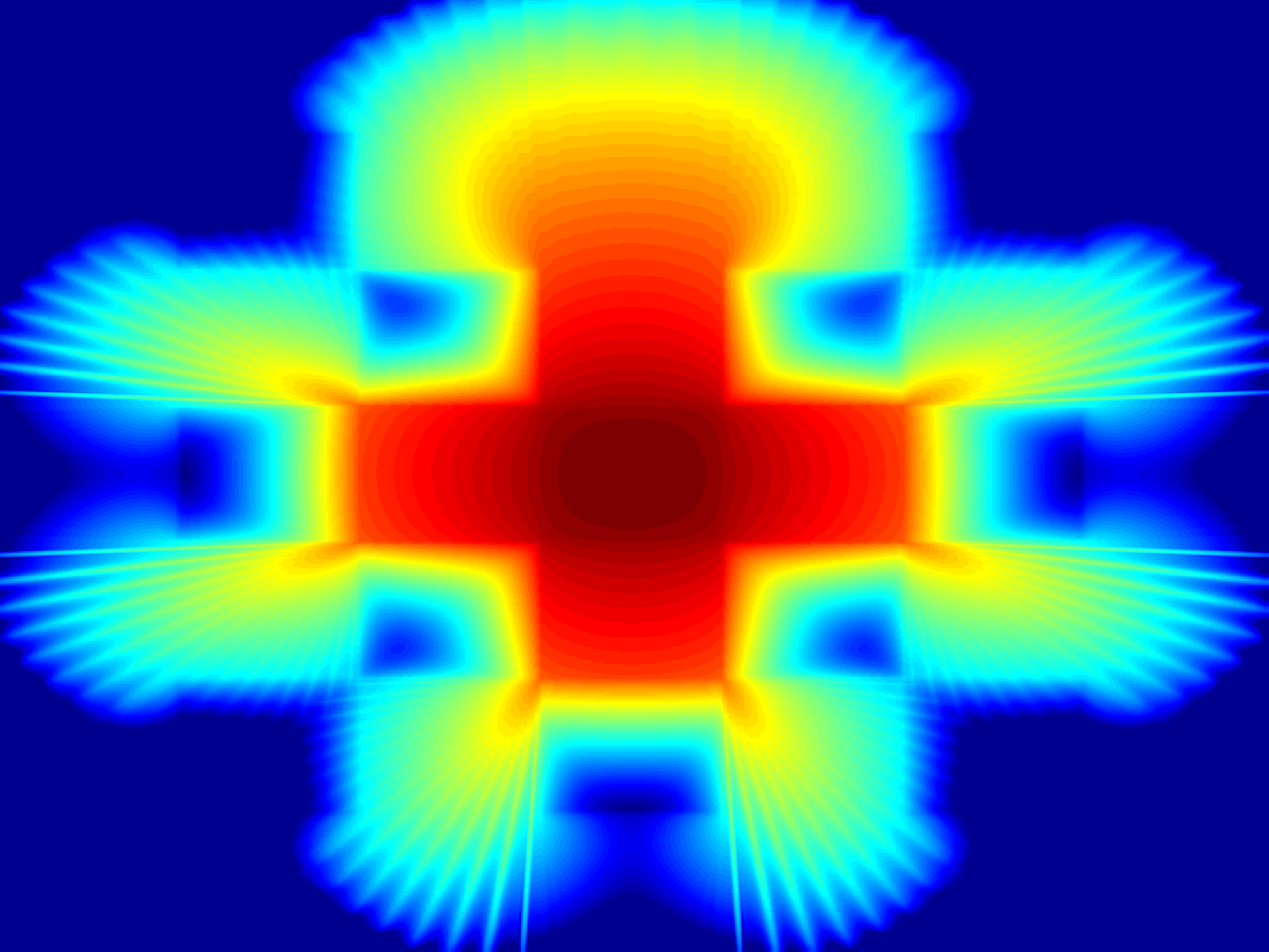}&
    \includegraphics[width=3cm]{latticebig_2d_ldmm_5p}
  \end{tabular}
  \caption{Compression of 2D scientific data sets with a 5\% data compression rate. The original data are shown on the upper left corners for each data set. The results of SVD, DCT, DFT, wavelet, and LDMM are shown in the remaining five figures.}
  \label{fig:compre_2d_5p}
\end{figure}

\begin{figure}[H]
  \centering
  \begin{tabular}{ccc}
    Original & Tucker (50.91dB)& DCT (\textbf{54.90dB})\\
    \includegraphics[width=3cm]{plasma_3d_original_band_19}&
    \includegraphics[width=3cm]{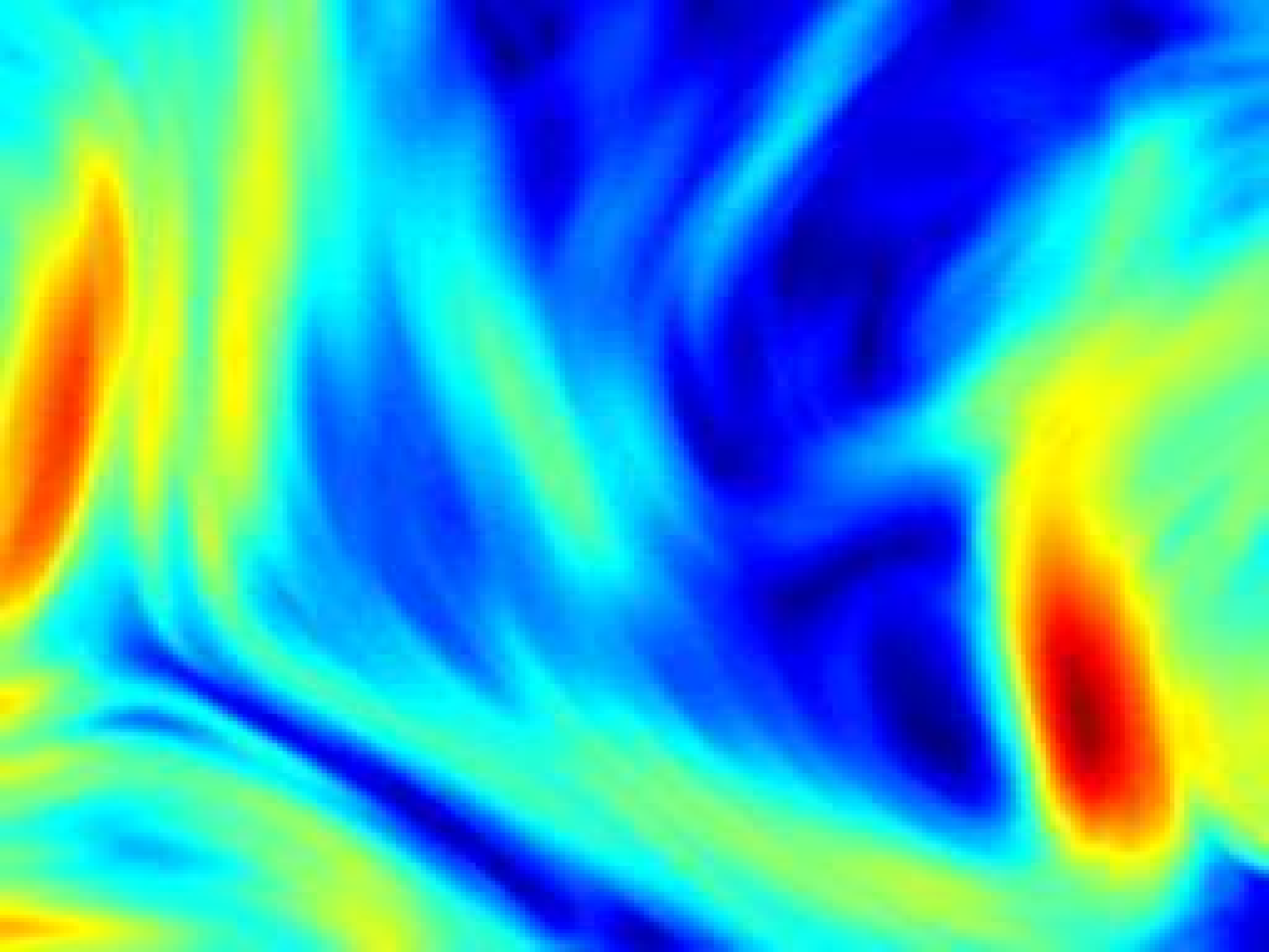}&
    \includegraphics[width=3cm]{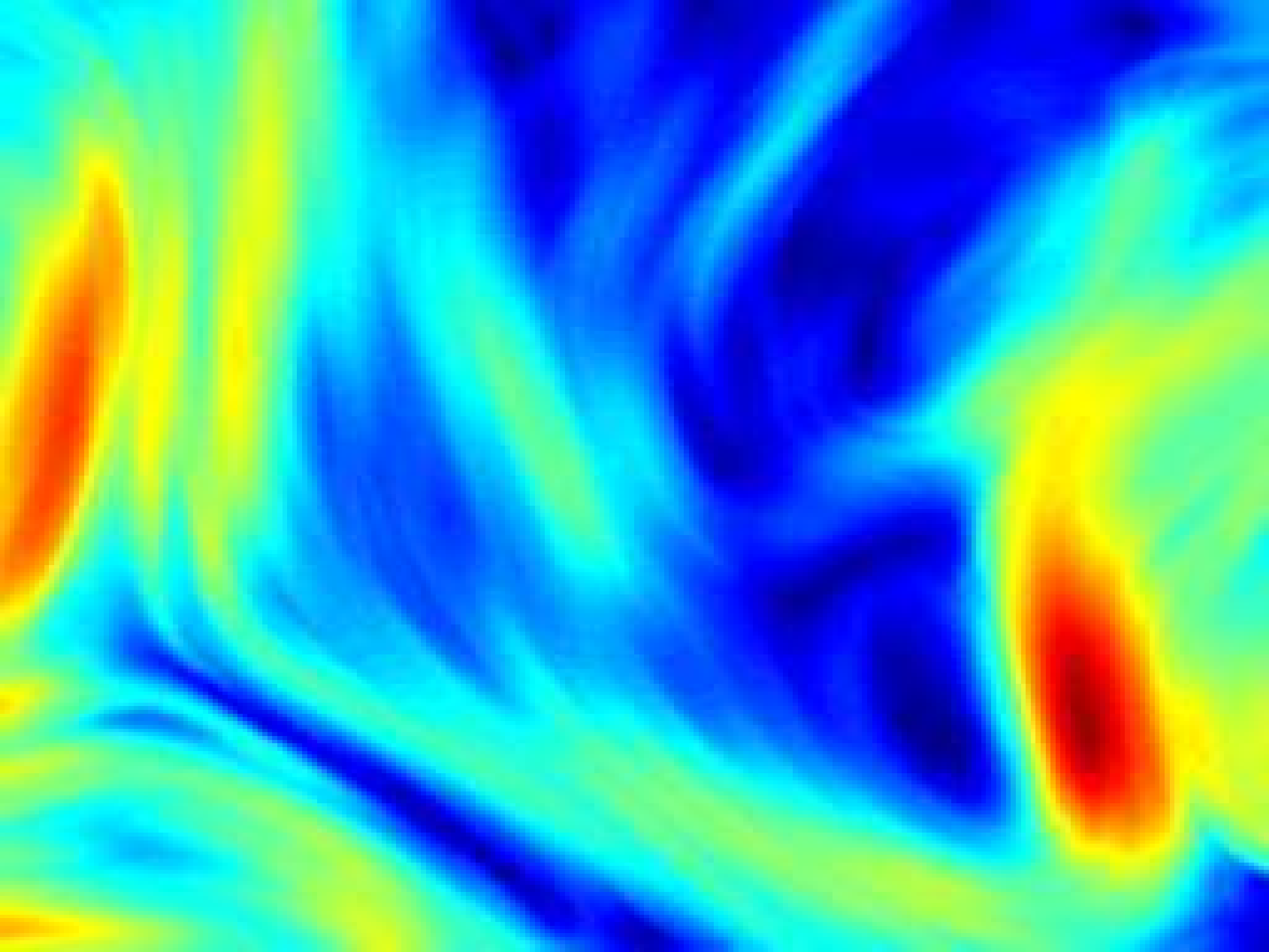}\\
    DFT (48.42dB)& Wavelet (41.01dB)& LDMM (44.18dB)\\
    \includegraphics[width=3cm]{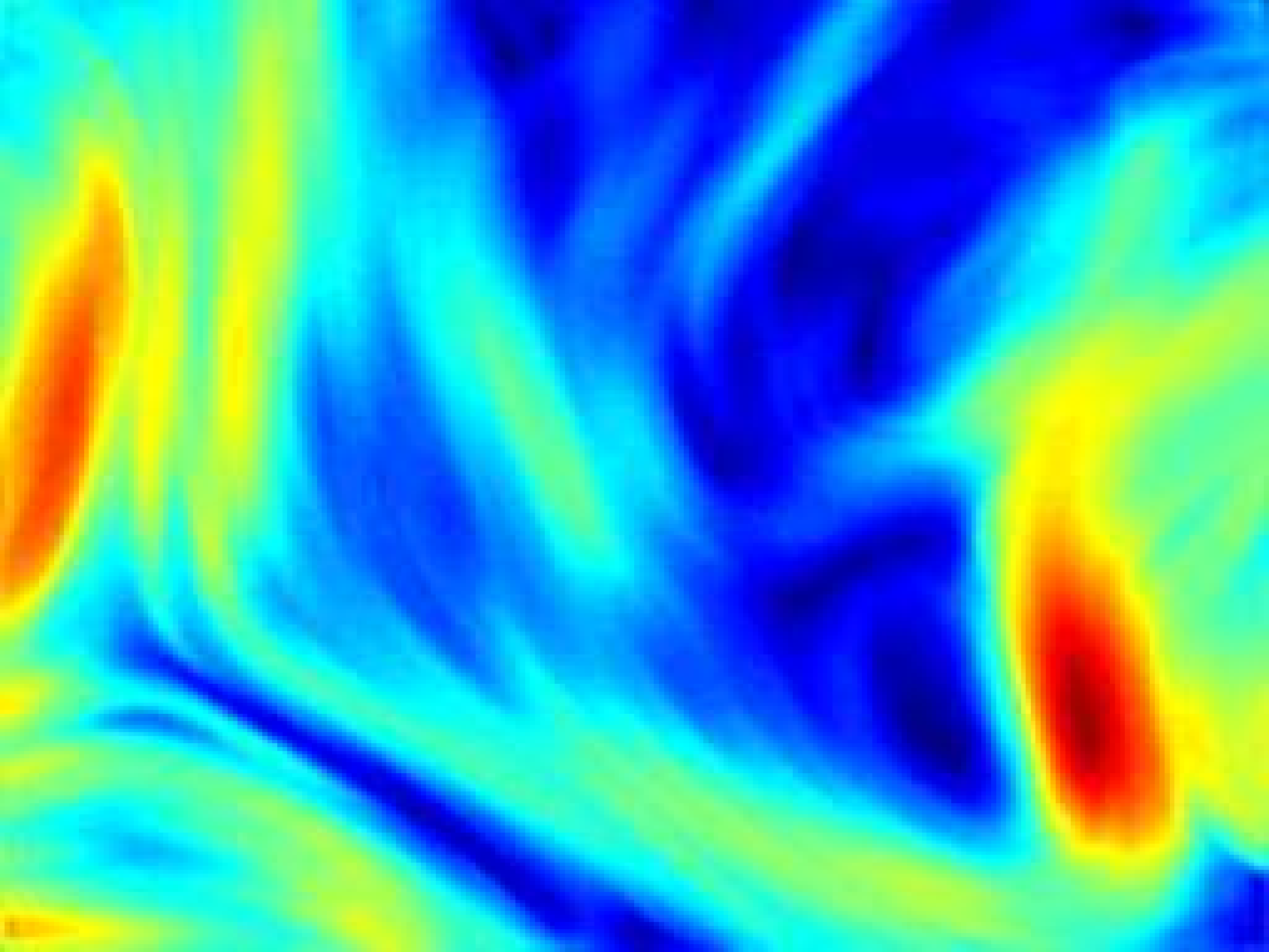}&
    \includegraphics[width=3cm]{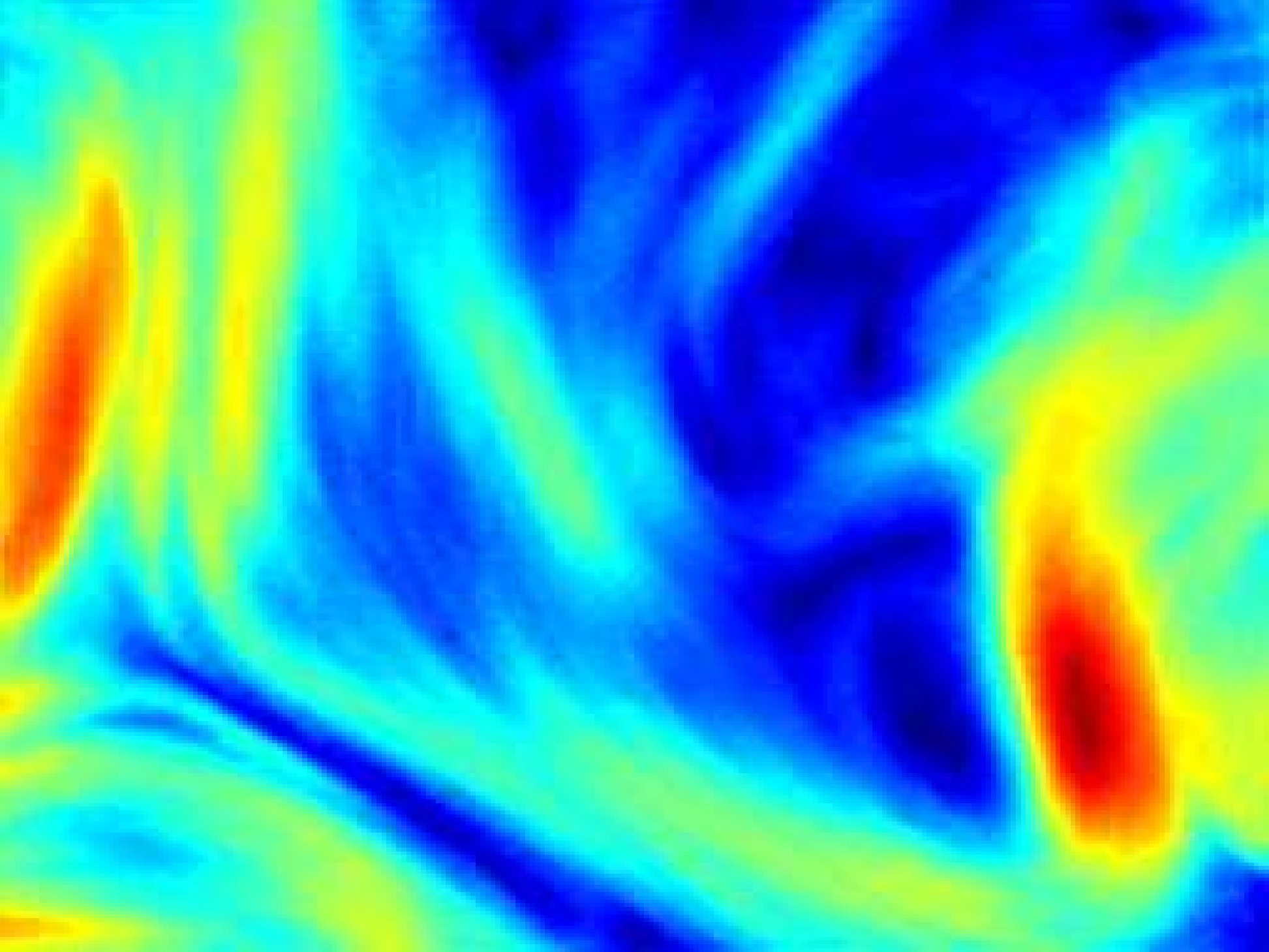}&
    \includegraphics[width=3cm]{plasma_3d_ldmm_10p_band_19}\\
    Original & Tucker (50.91dB)& DCT (\textbf{54.90dB})\\
    \includegraphics[width=3cm]{plasma_3d_original_band_29}&
    \includegraphics[width=3cm]{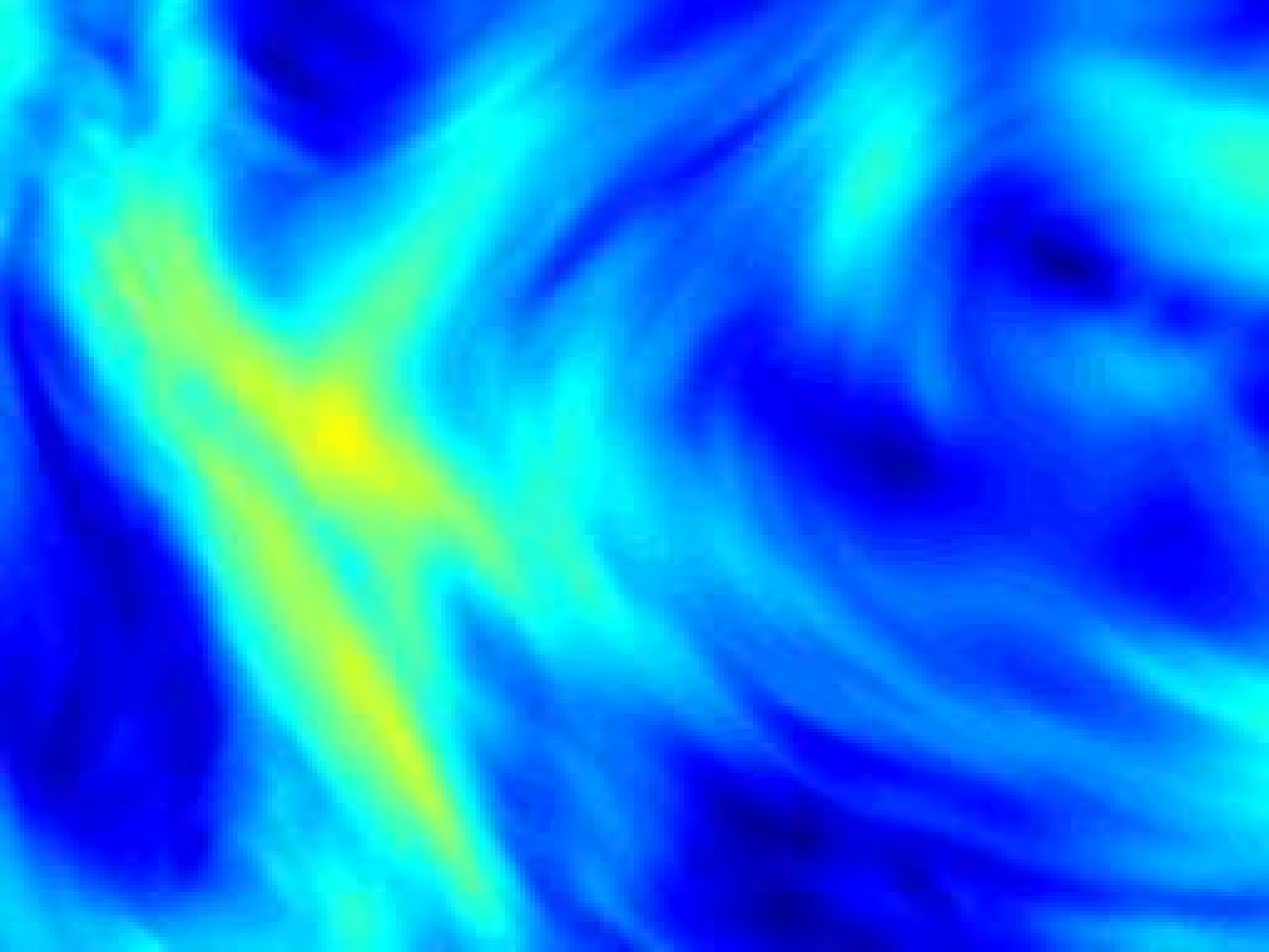}&
    \includegraphics[width=3cm]{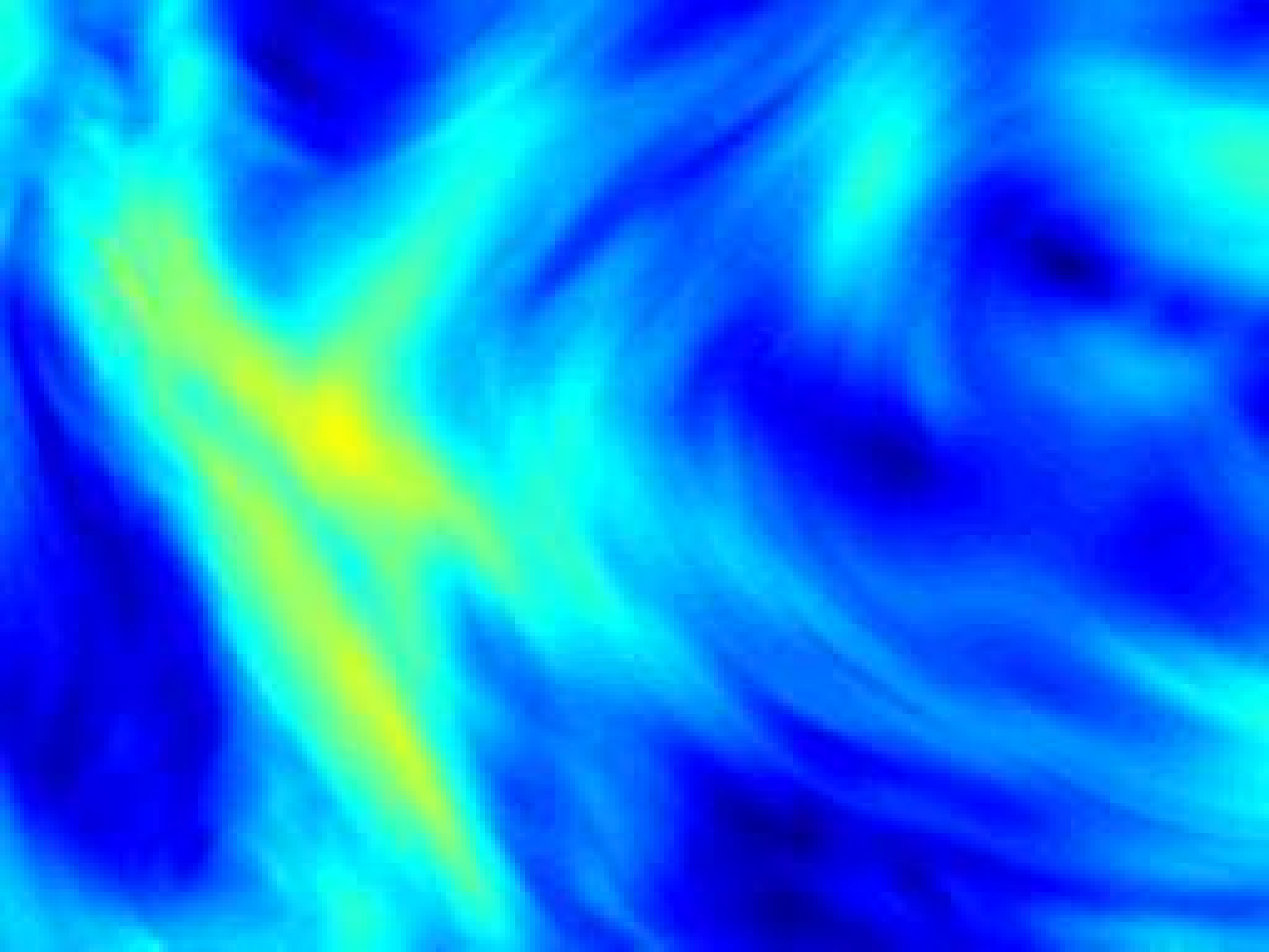}\\
    DFT (48.42dB)& Wavelet (41.01dB)& LDMM (44.18dB)\\
    \includegraphics[width=3cm]{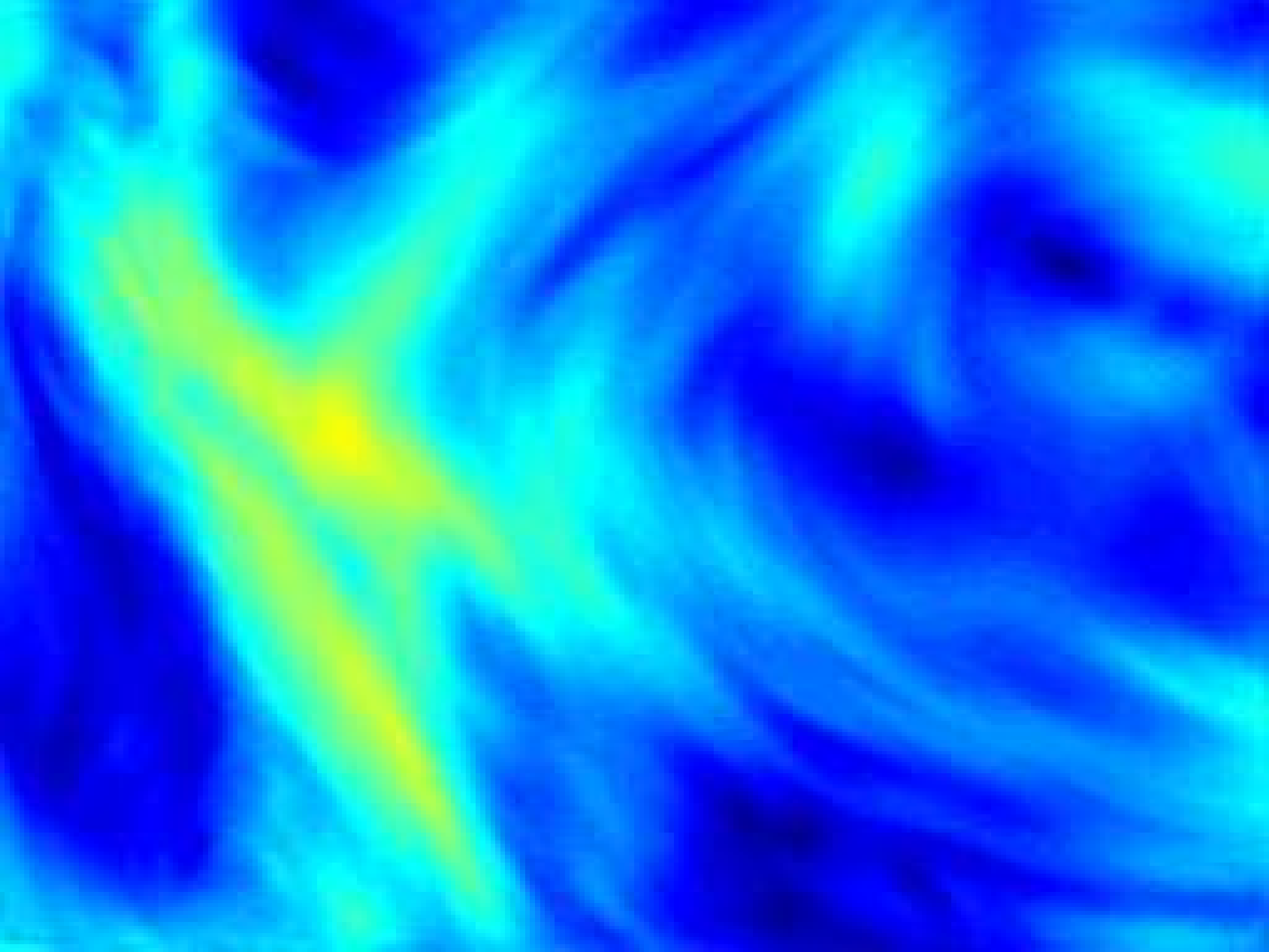}&
    \includegraphics[width=3cm]{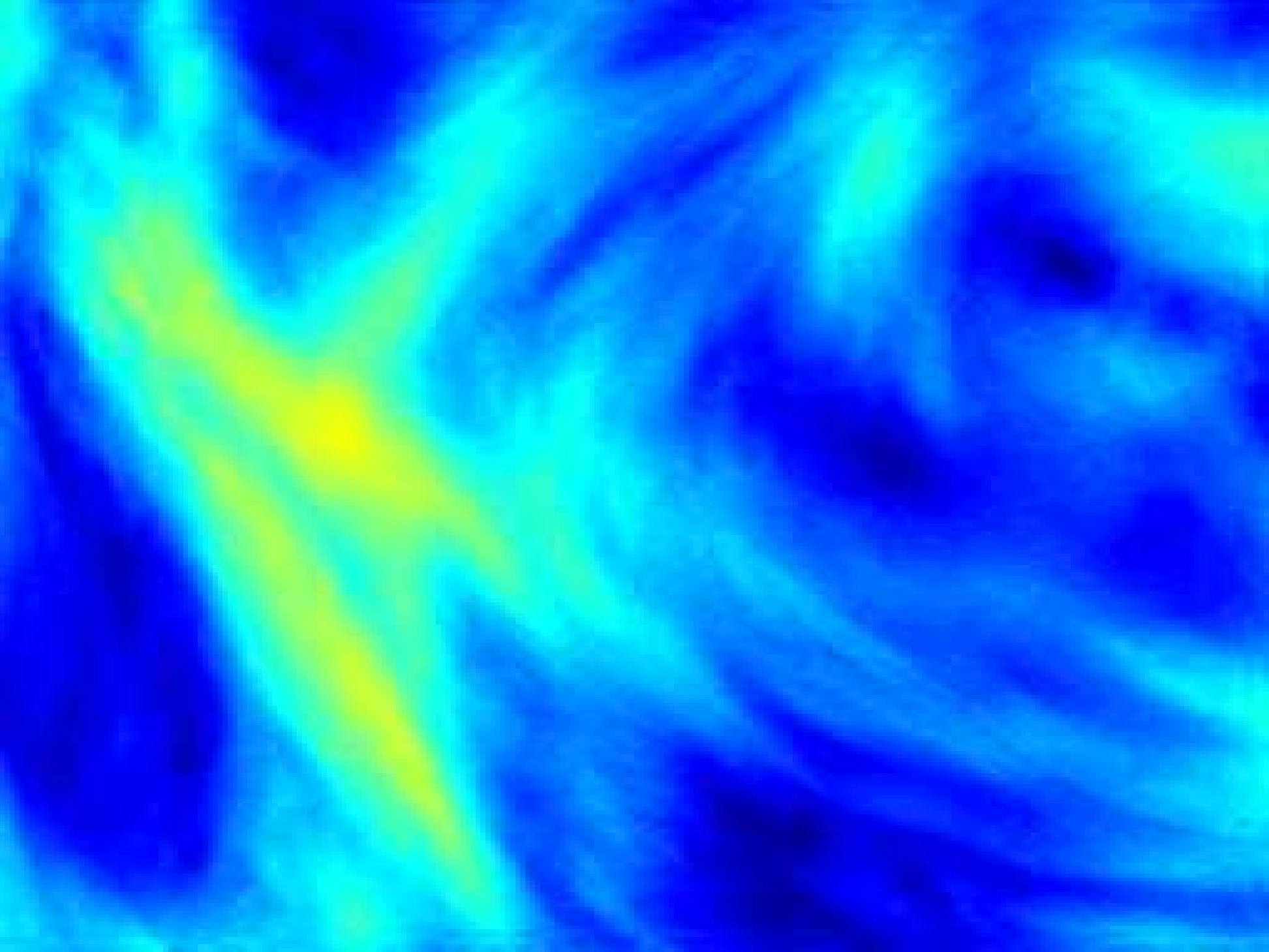}&
    \includegraphics[width=3cm]{plasma_3d_ldmm_10p_band_29}
  \end{tabular}
  \caption{Compression of the 3D plasma (magnetic field) data set with a 10\% data compression rate. Two  spatial cross  sections of the original data set are shown in the first figures on the first and third row. The results of Tucker decomposition, DCT, DFT, wavelet, and LDMM are shown in the remaining five figures.}
  \label{fig:compre_plasma_3d_10p}
\end{figure}

\begin{figure}[H]
  \centering
  \begin{tabular}{ccc}
    Original & Tucker (45.36dB)& DCT (\textbf{49.70dB})\\
    \includegraphics[width=3cm]{plasma_3d_original_band_19}&
    \includegraphics[width=3cm]{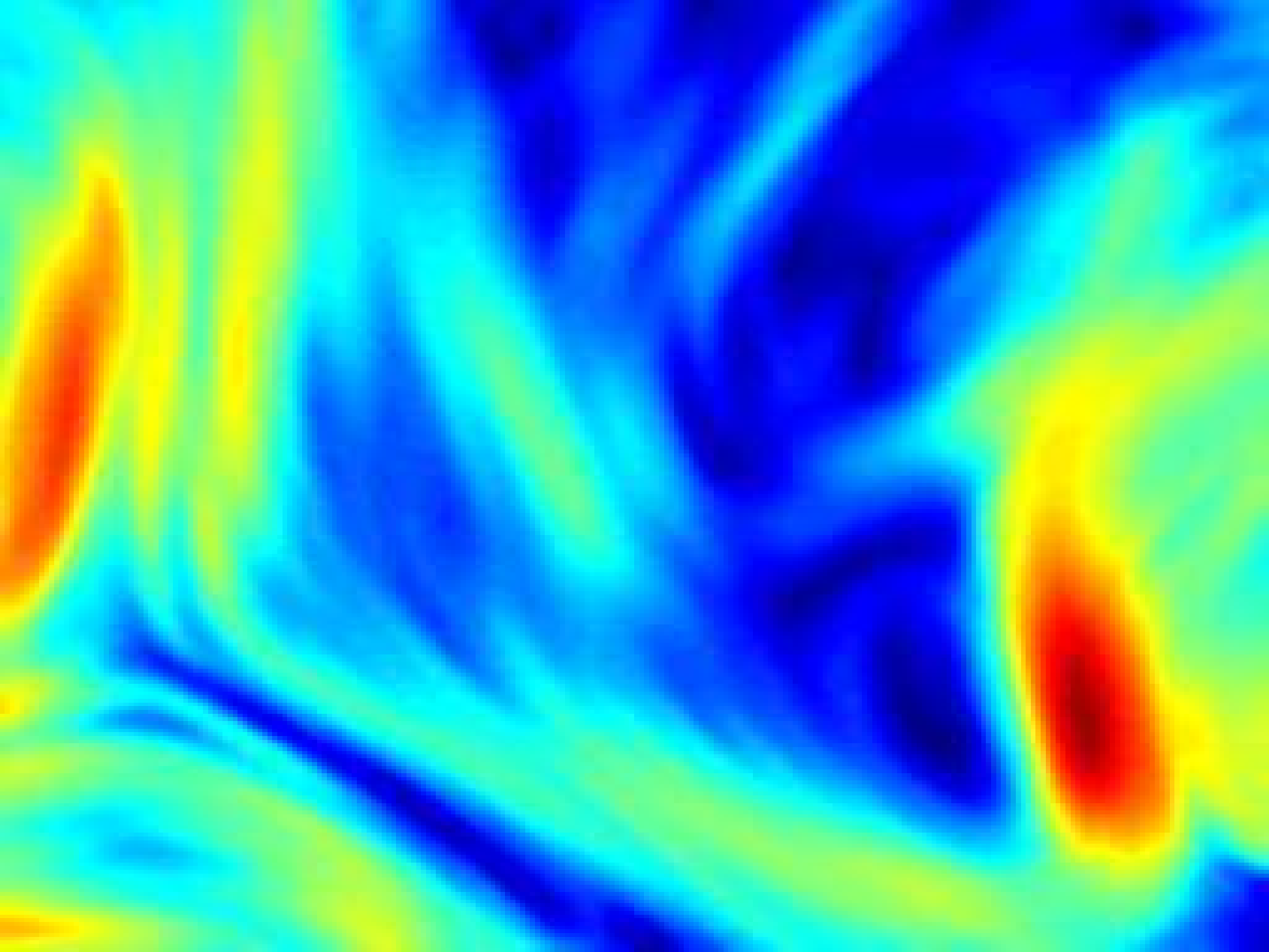}&
    \includegraphics[width=3cm]{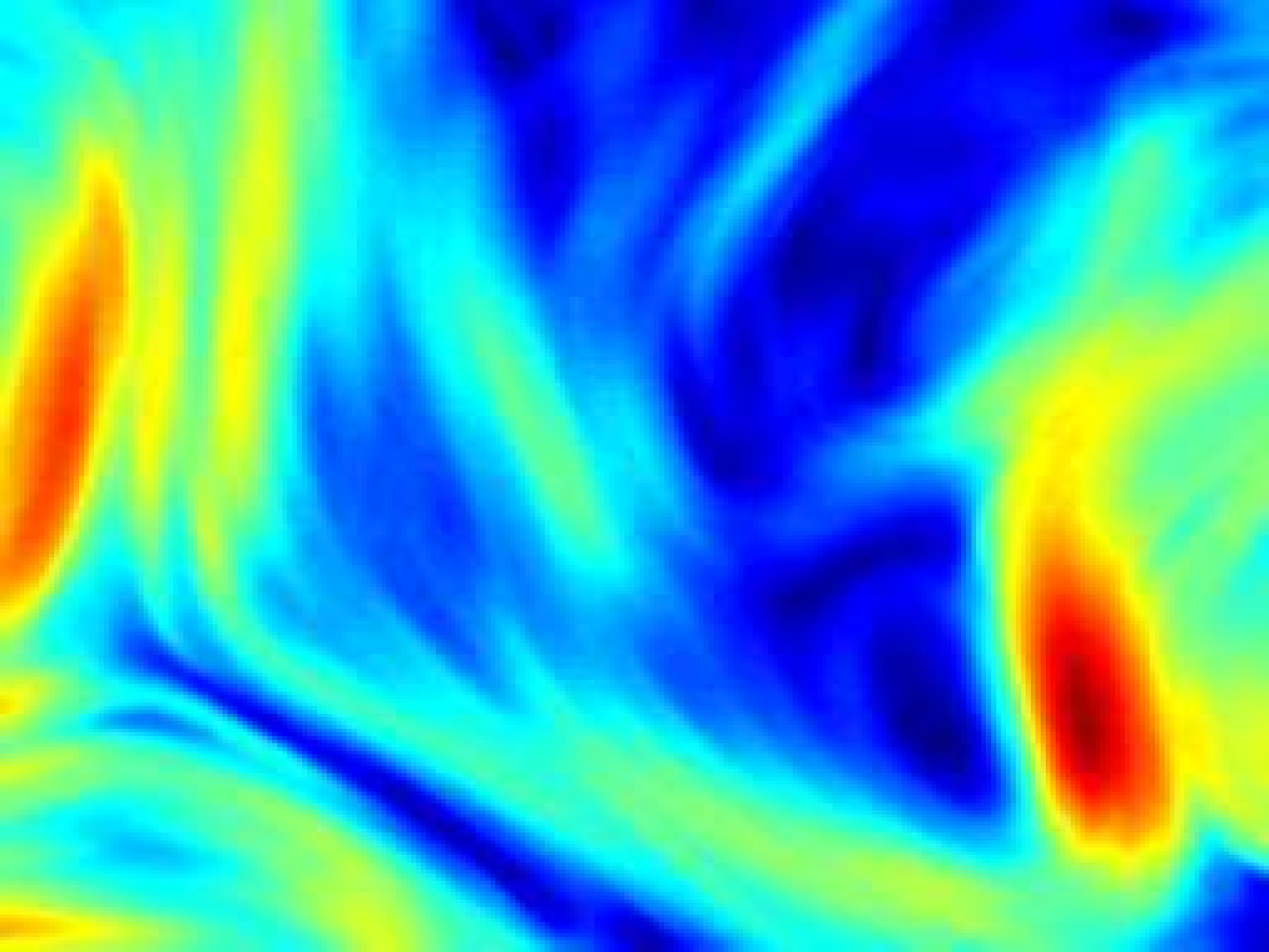}\\
    DFT (43.56dB)& Wavelet (32.74dB)& LDMM (40.07dB)\\
    \includegraphics[width=3cm]{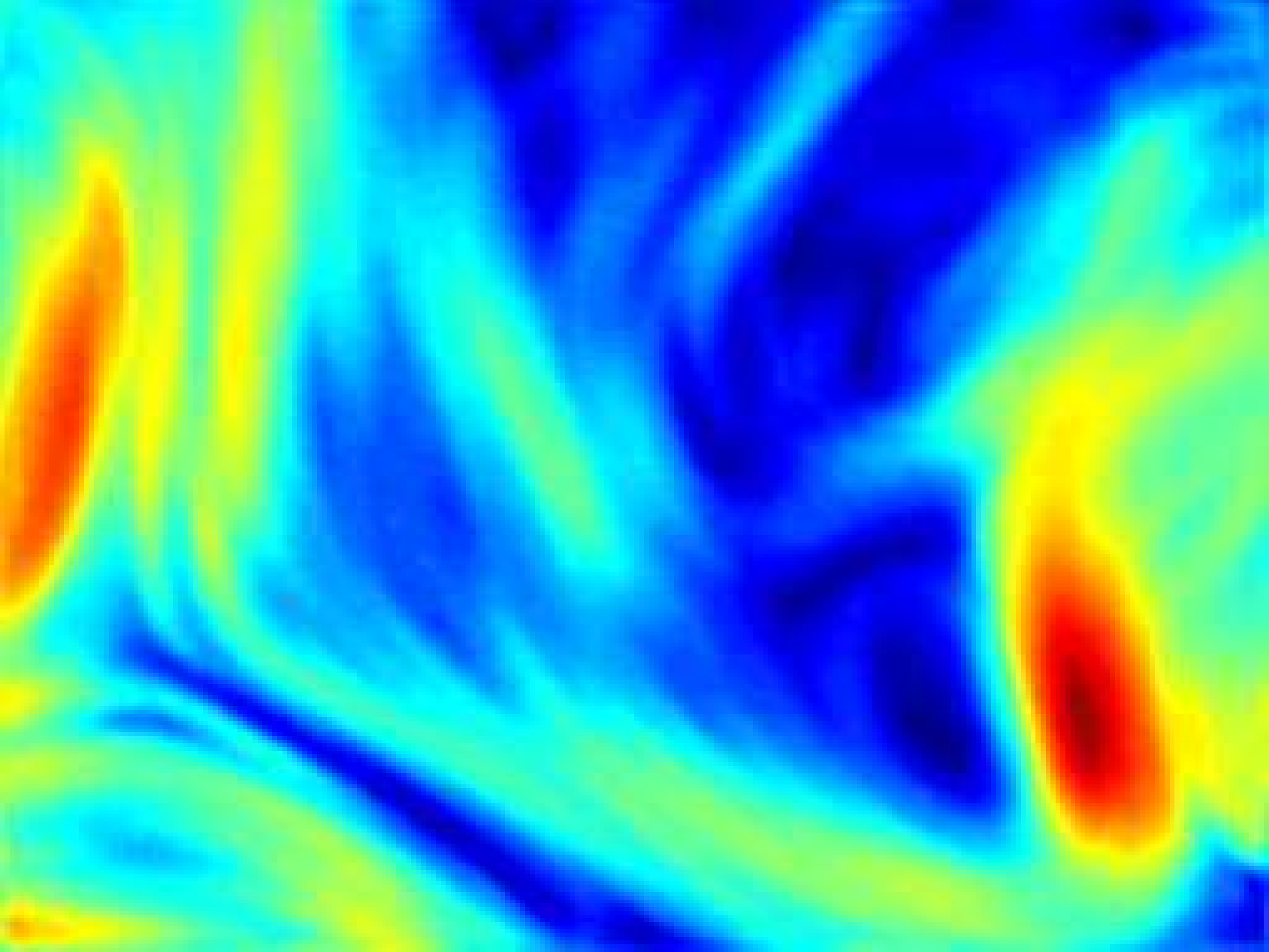}&
    \includegraphics[width=3cm]{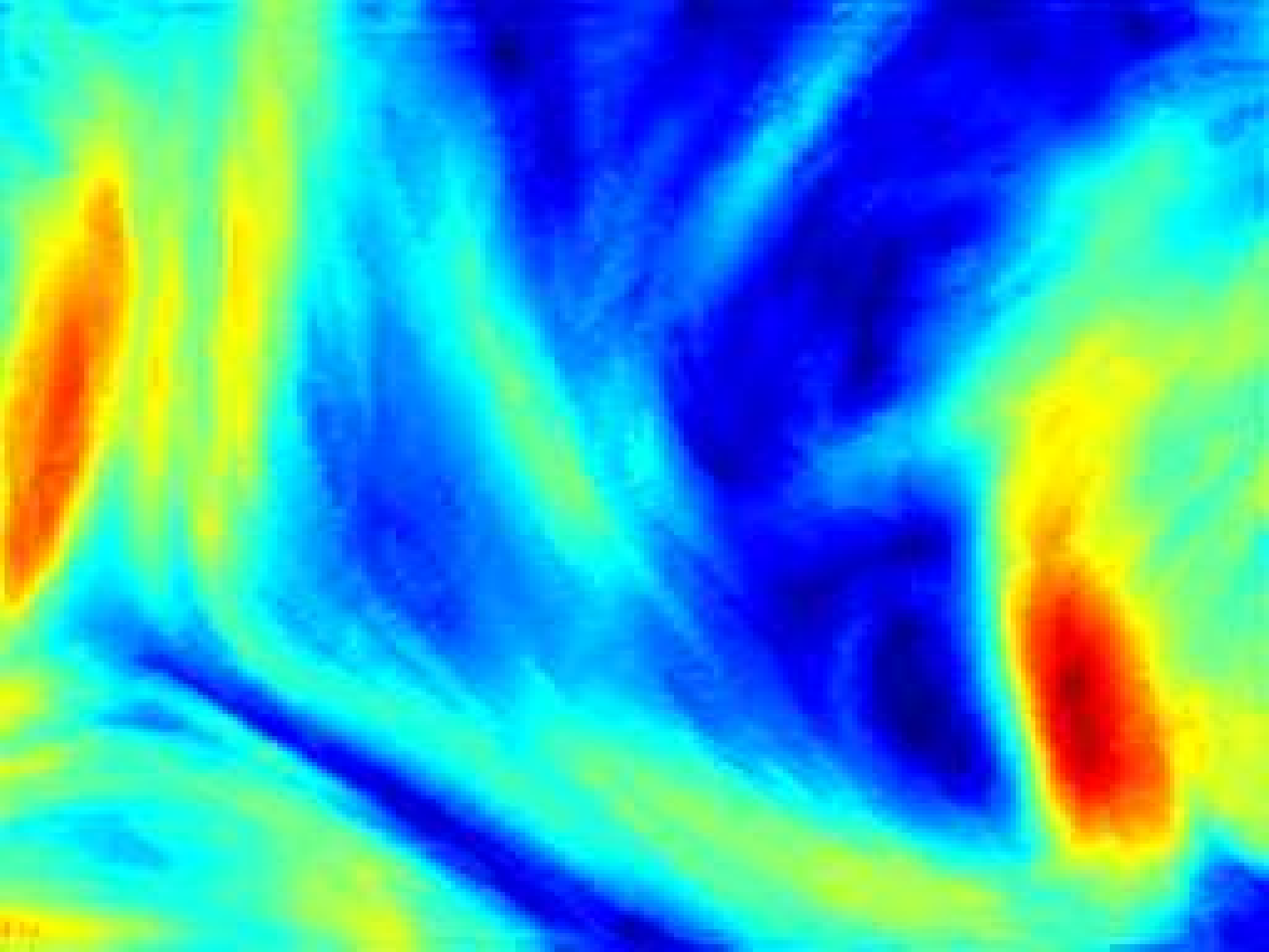}&
    \includegraphics[width=3cm]{plasma_3d_ldmm_5p_band_19}\\
    Original & Tucker (45.36dB)& DCT (\textbf{49.70dB})\\
    \includegraphics[width=3cm]{plasma_3d_original_band_29}&
    \includegraphics[width=3cm]{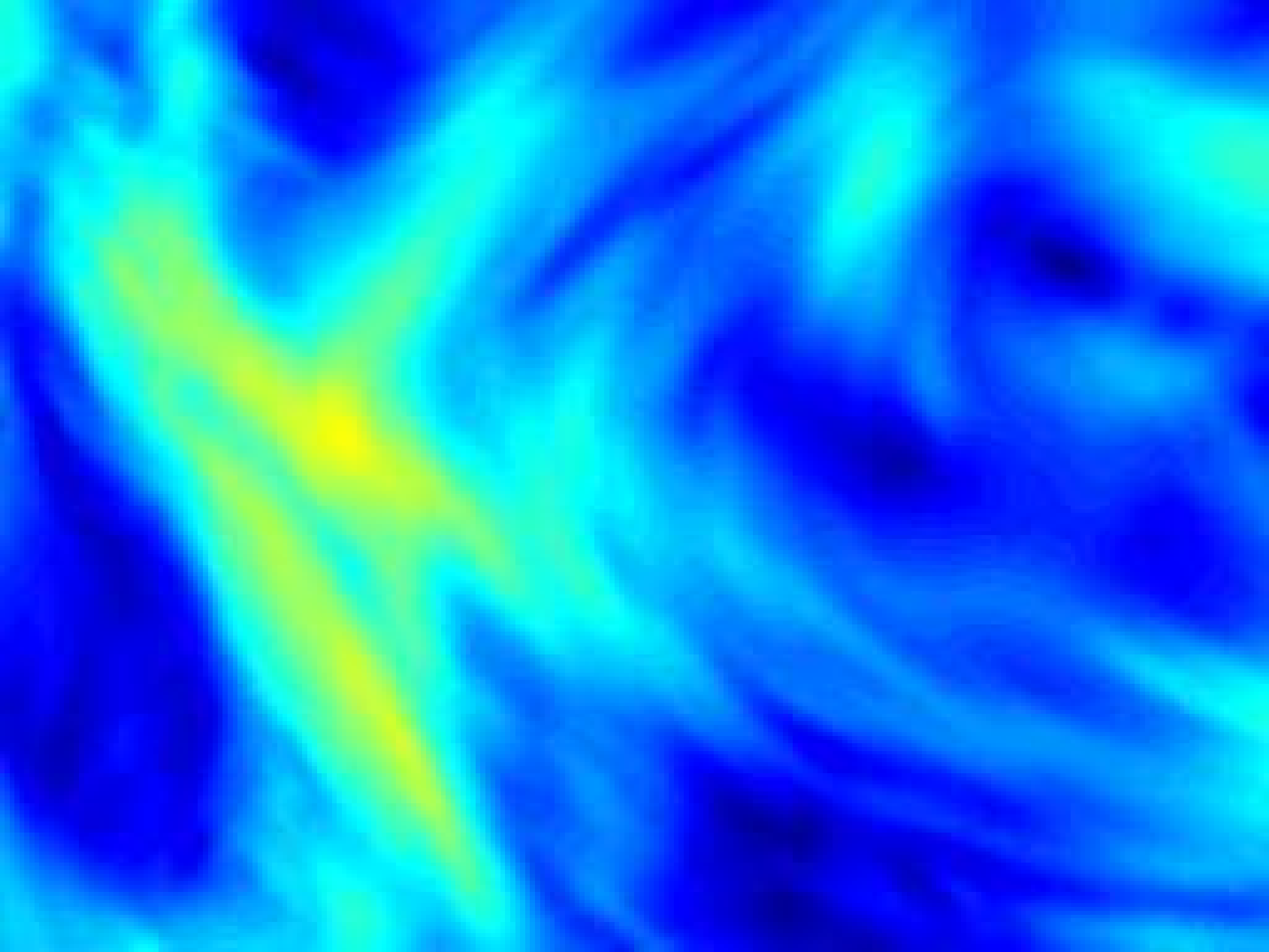}&
    \includegraphics[width=3cm]{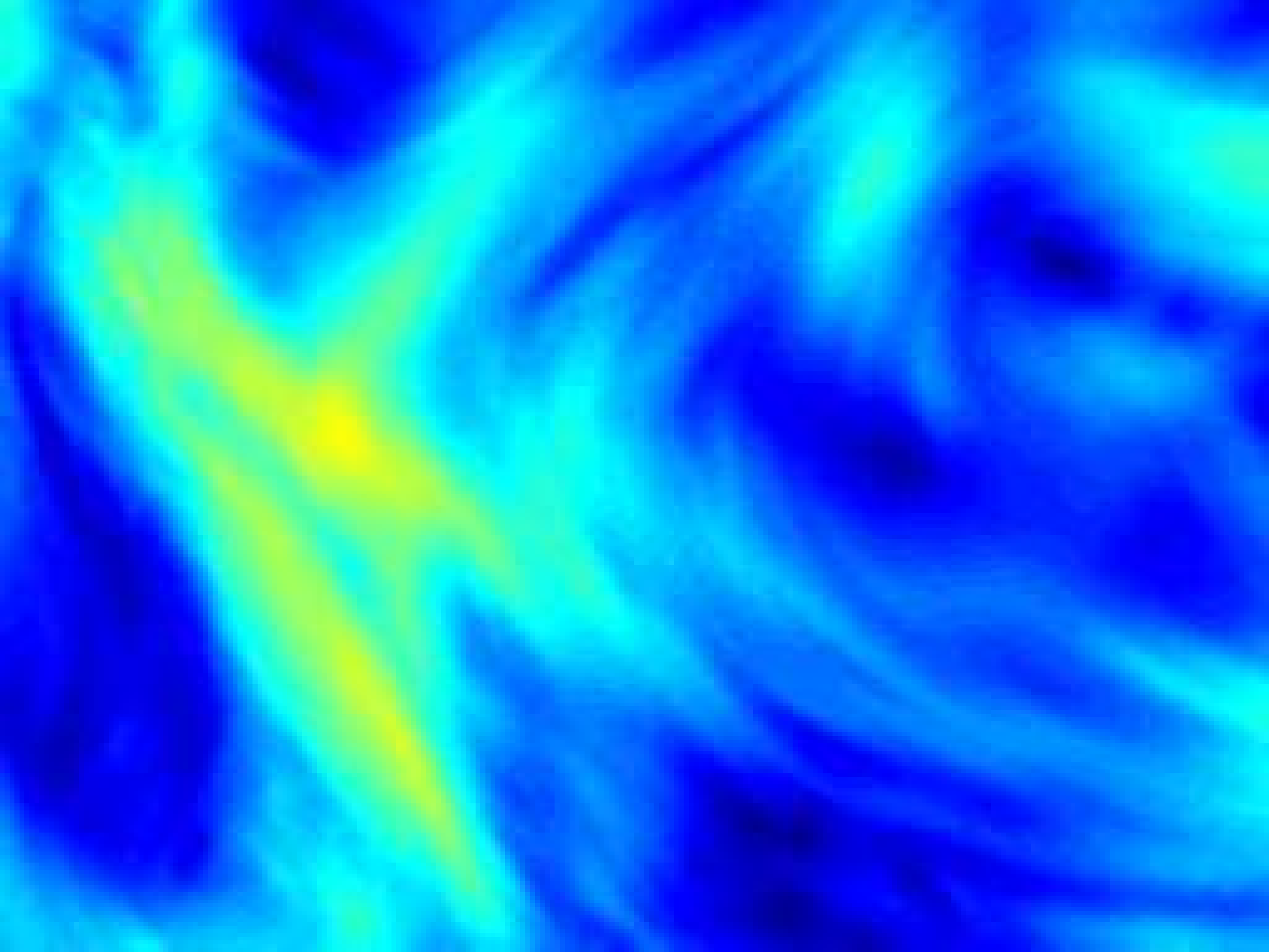}\\
    DFT (43.56dB)& Wavelet (32.74dB)& LDMM (40.07dB)\\
    \includegraphics[width=3cm]{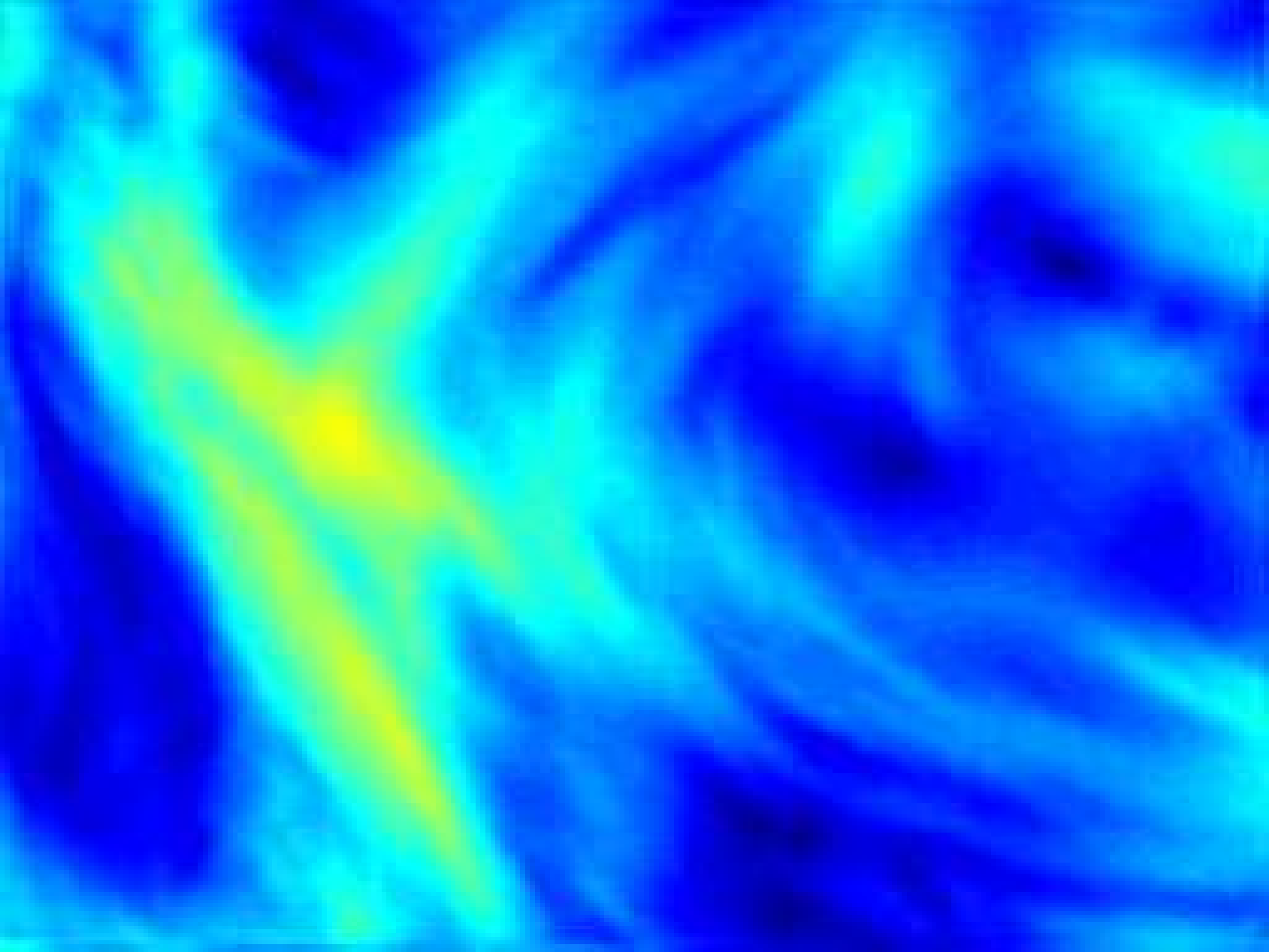}&
    \includegraphics[width=3cm]{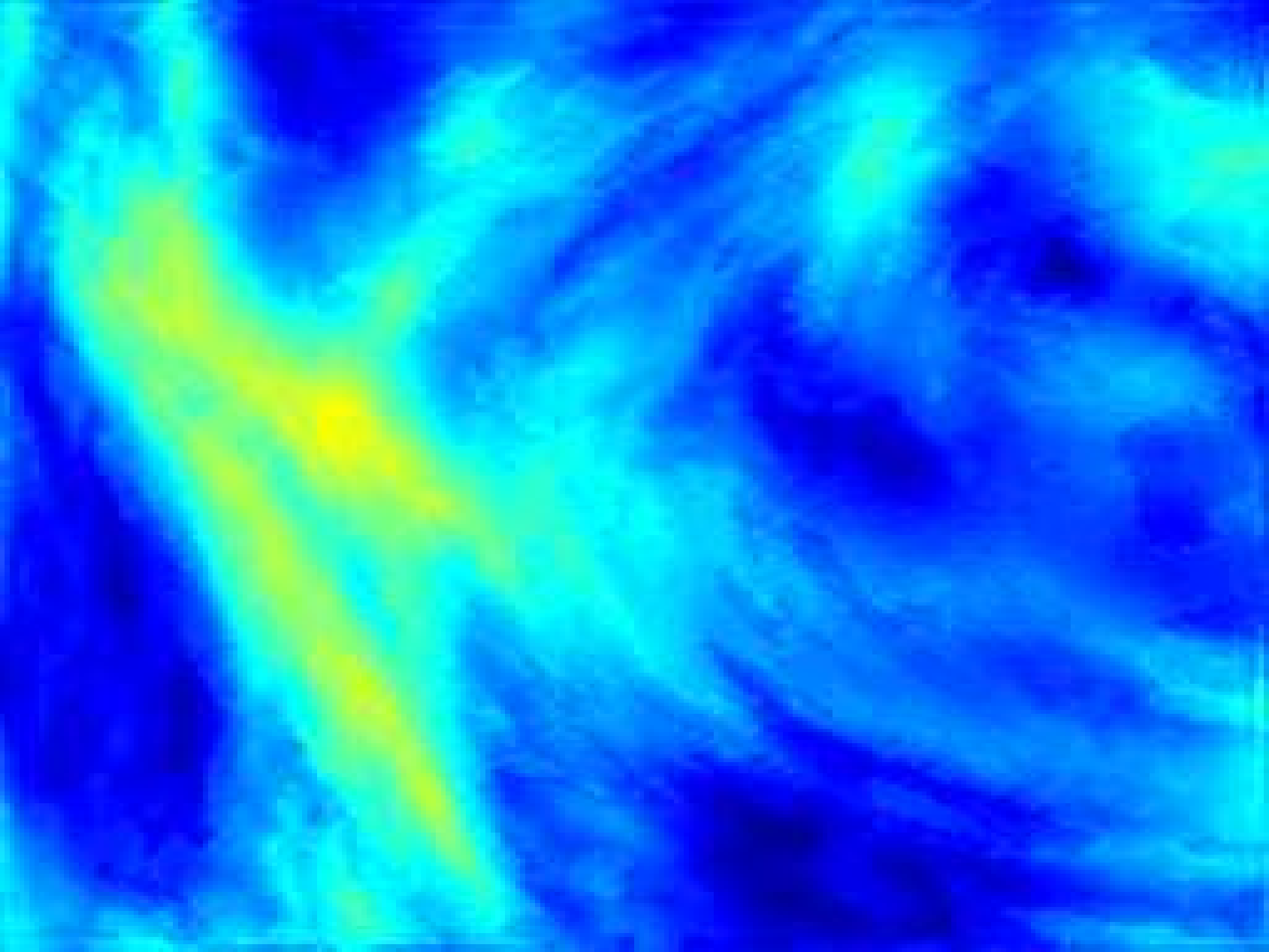}&
    \includegraphics[width=3cm]{plasma_3d_ldmm_5p_band_29}
  \end{tabular}
  \caption{Compression of the 3D plasma (magnetic field) data set with a 5\% data compression rate. Two  spatial cross  sections of the original data set are shown in the first figures on the first and third row. The results of Tucker decomposition, DCT, DFT, wavelet, and LDMM are shown in the remaining five figures.}
  \label{fig:compre_plasma_3d_5p}
\end{figure}

\begin{figure}[H]
  \centering
  \begin{tabular}{ccc}
    Original & Tucker (\textbf{97.43dB})& DCT (65.44dB)\\
    \includegraphics[width=3cm]{lattice_3d_original_band_31}&
    \includegraphics[width=3cm]{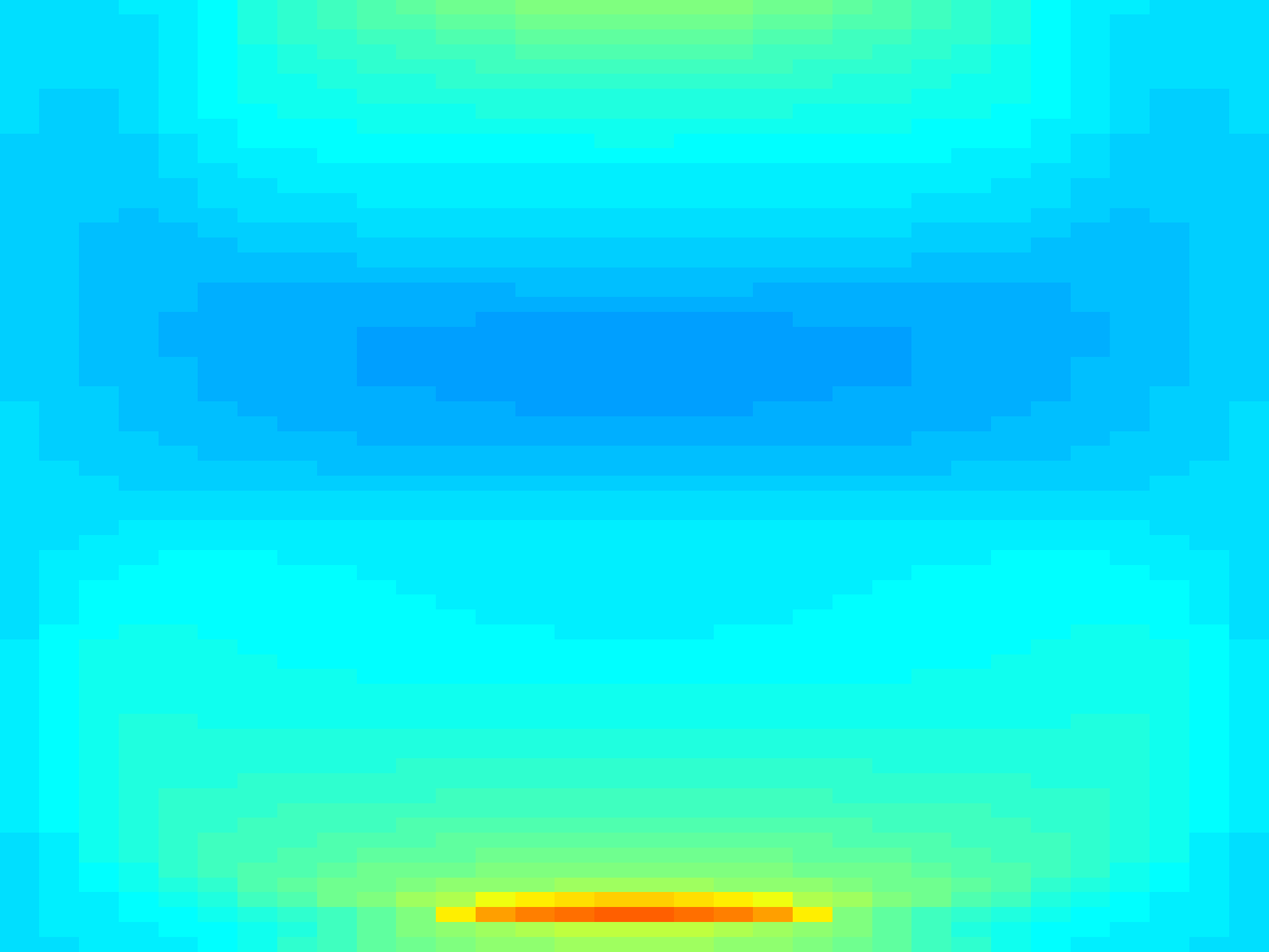}&
    \includegraphics[width=3cm]{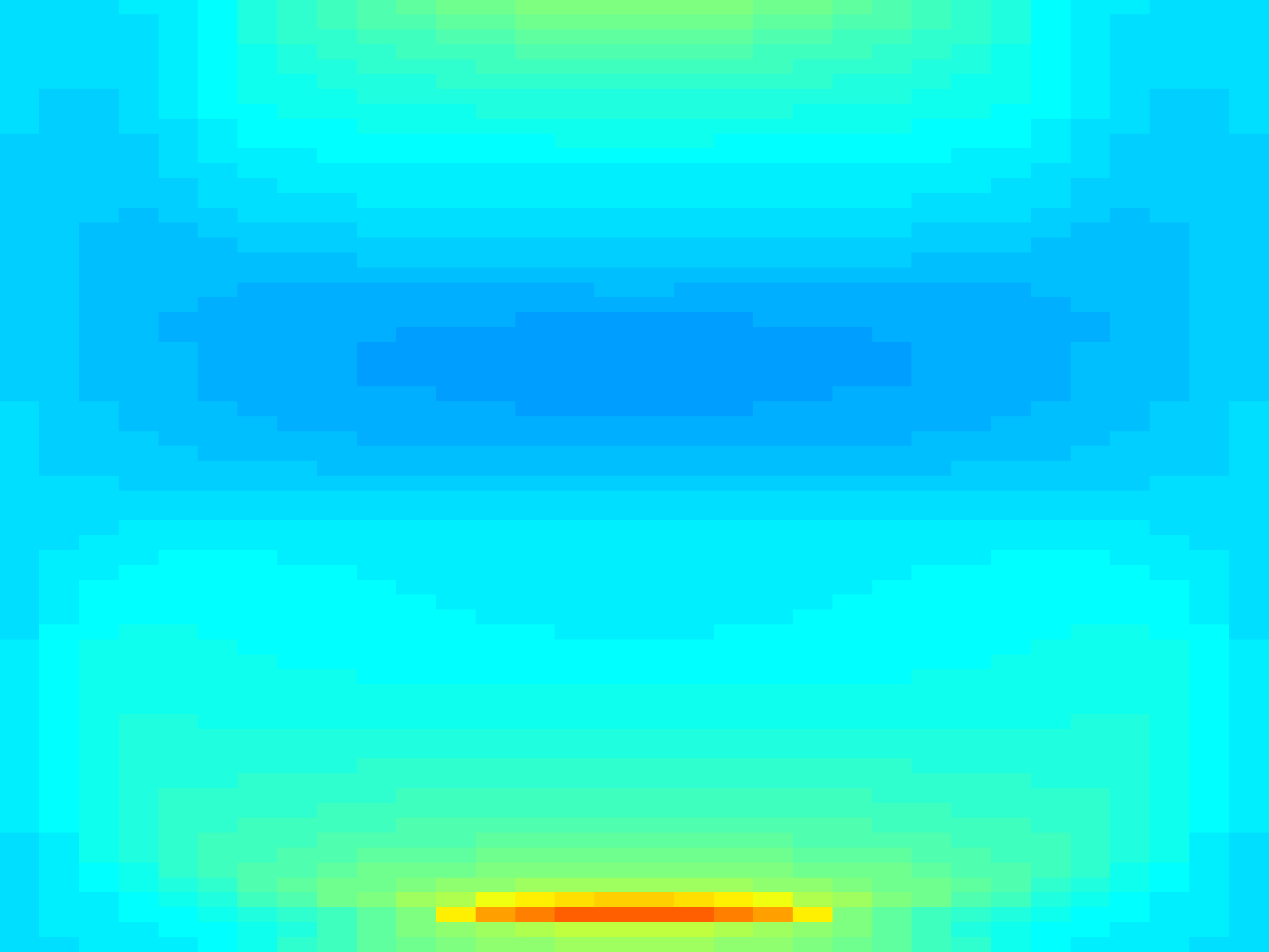}\\
    DFT (52.96dB)& Wavelet (72.61dB)& LDMM (48.43dB)\\
    \includegraphics[width=3cm]{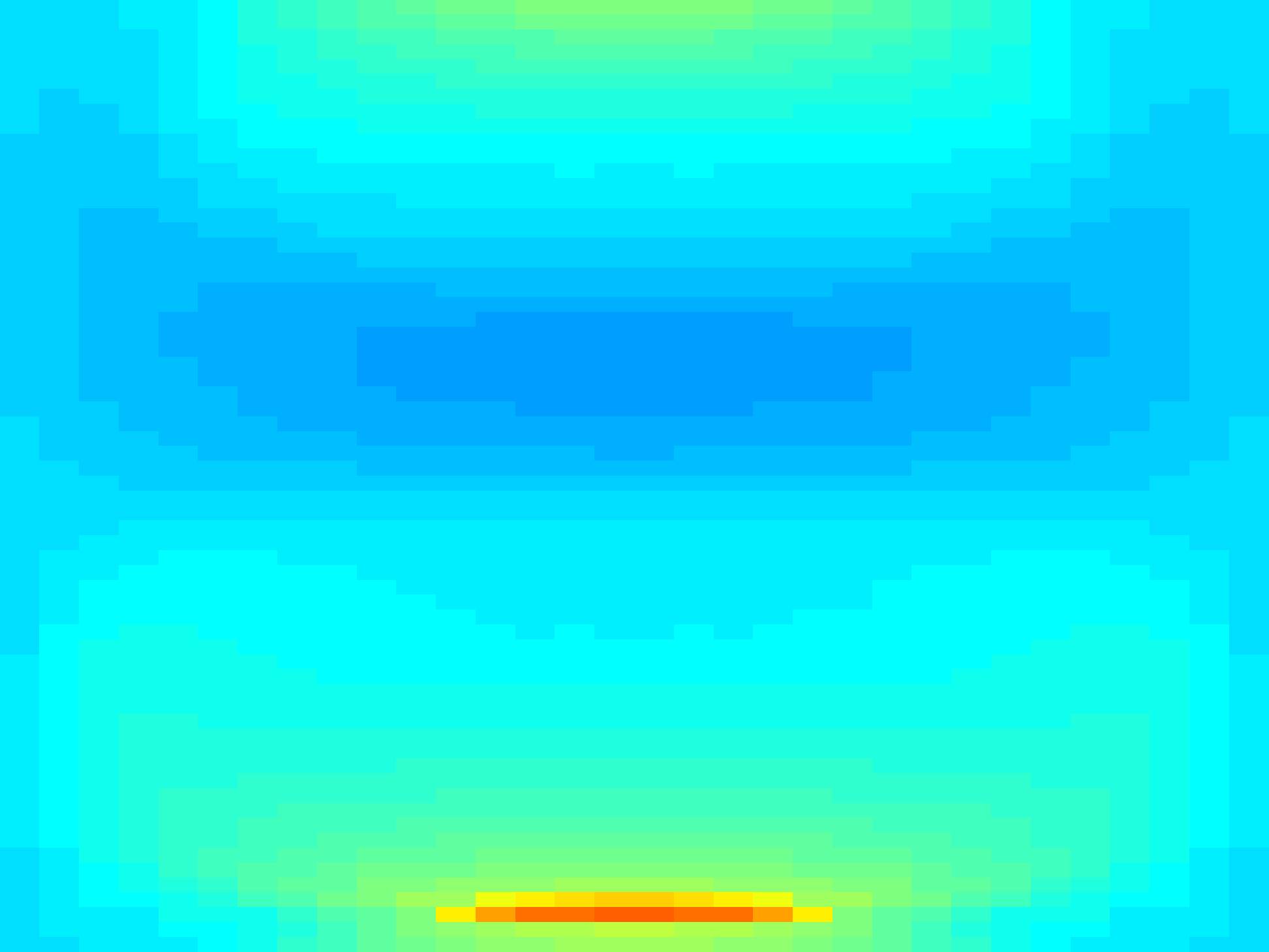}&
    \includegraphics[width=3cm]{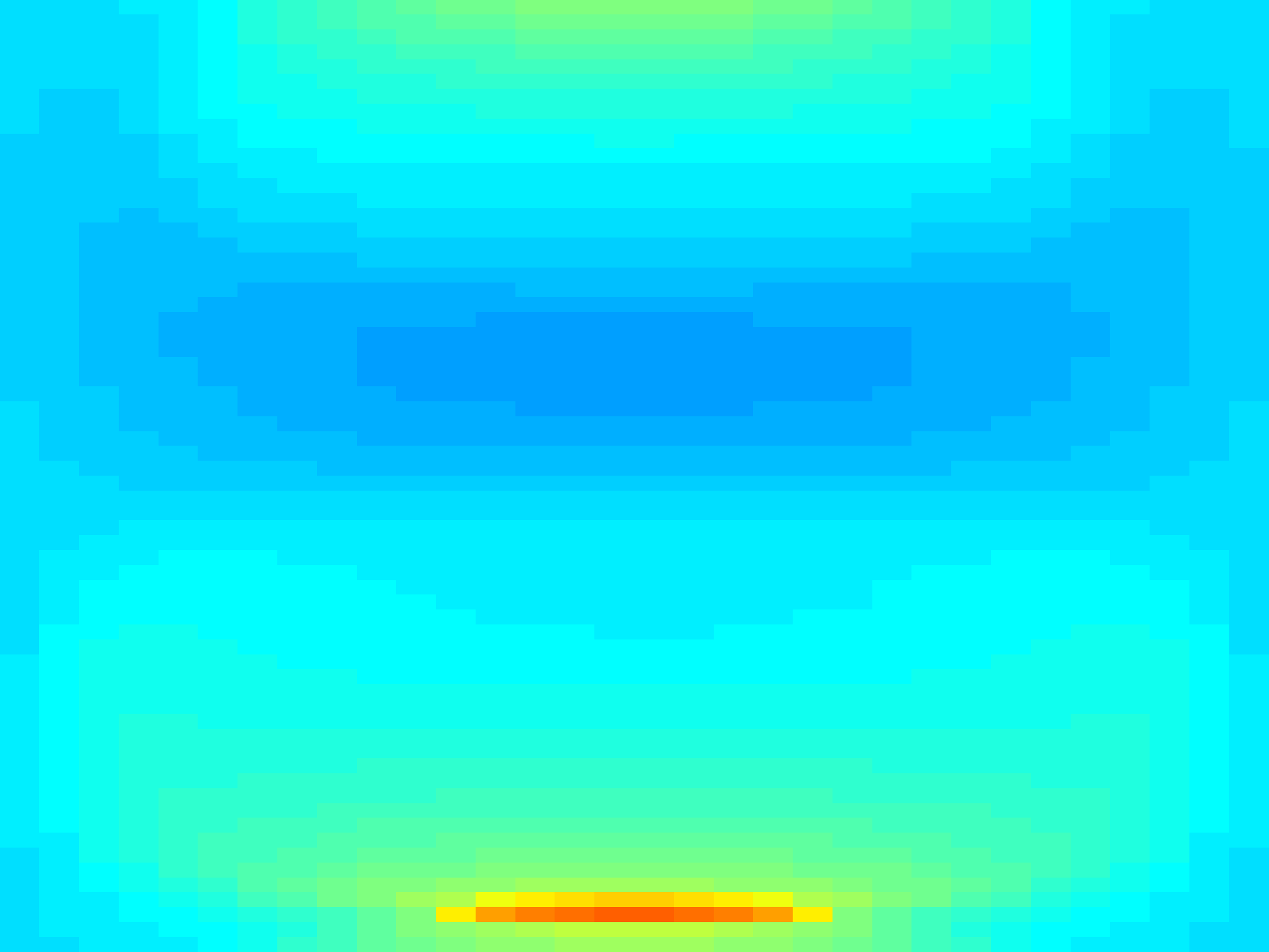}&
    \includegraphics[width=3cm]{lattice_3d_ldmm_10p_band_31}\\
    Original & Tucker (\textbf{97.43dB})& DCT (65.44dB)\\
    \includegraphics[width=3cm]{lattice_3d_original_band_151}&
    \includegraphics[width=3cm]{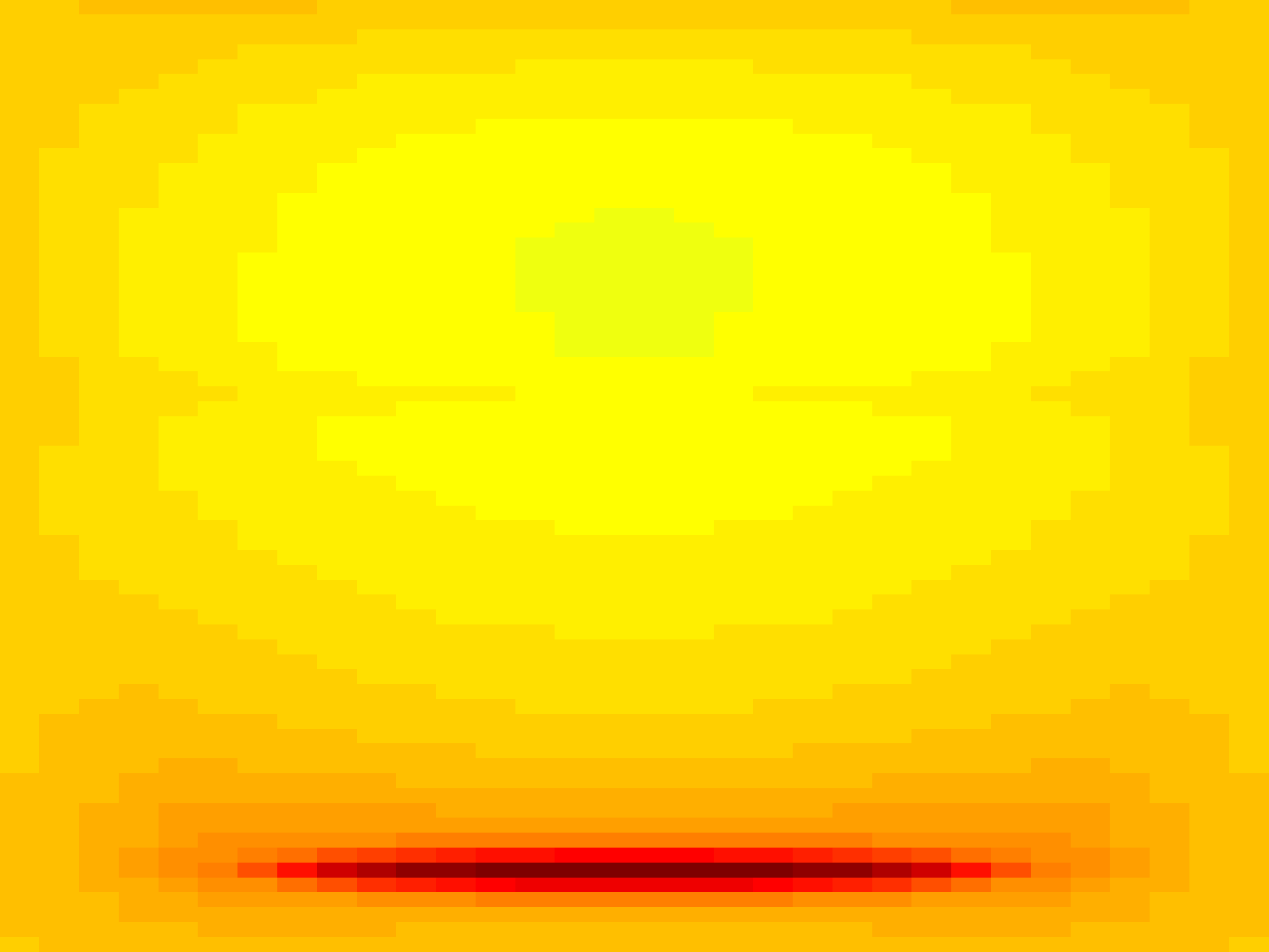}&
    \includegraphics[width=3cm]{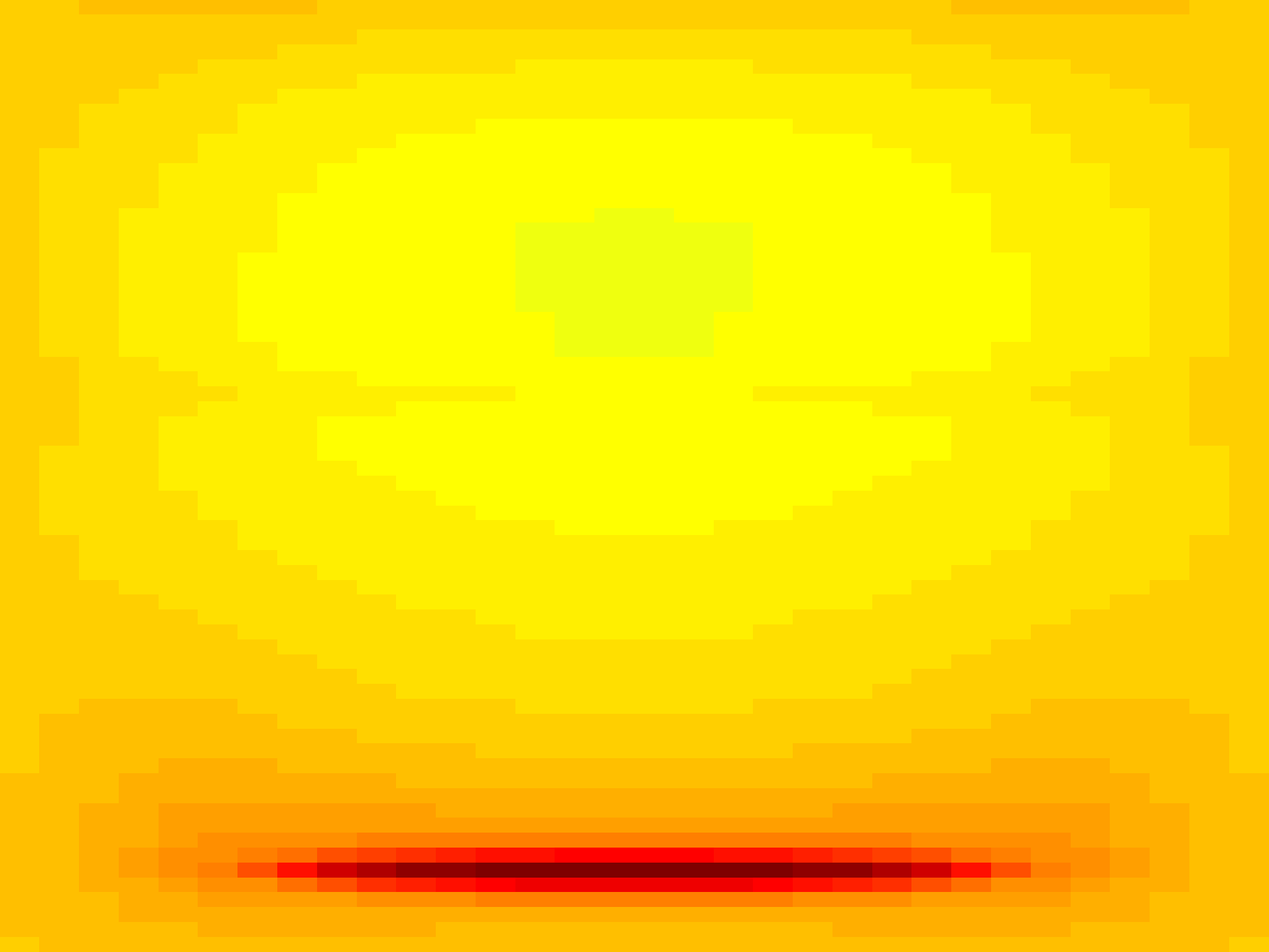}\\
    DFT (52.96dB)& Wavelet (72.61dB)& LDMM (48.43dB)\\
    \includegraphics[width=3cm]{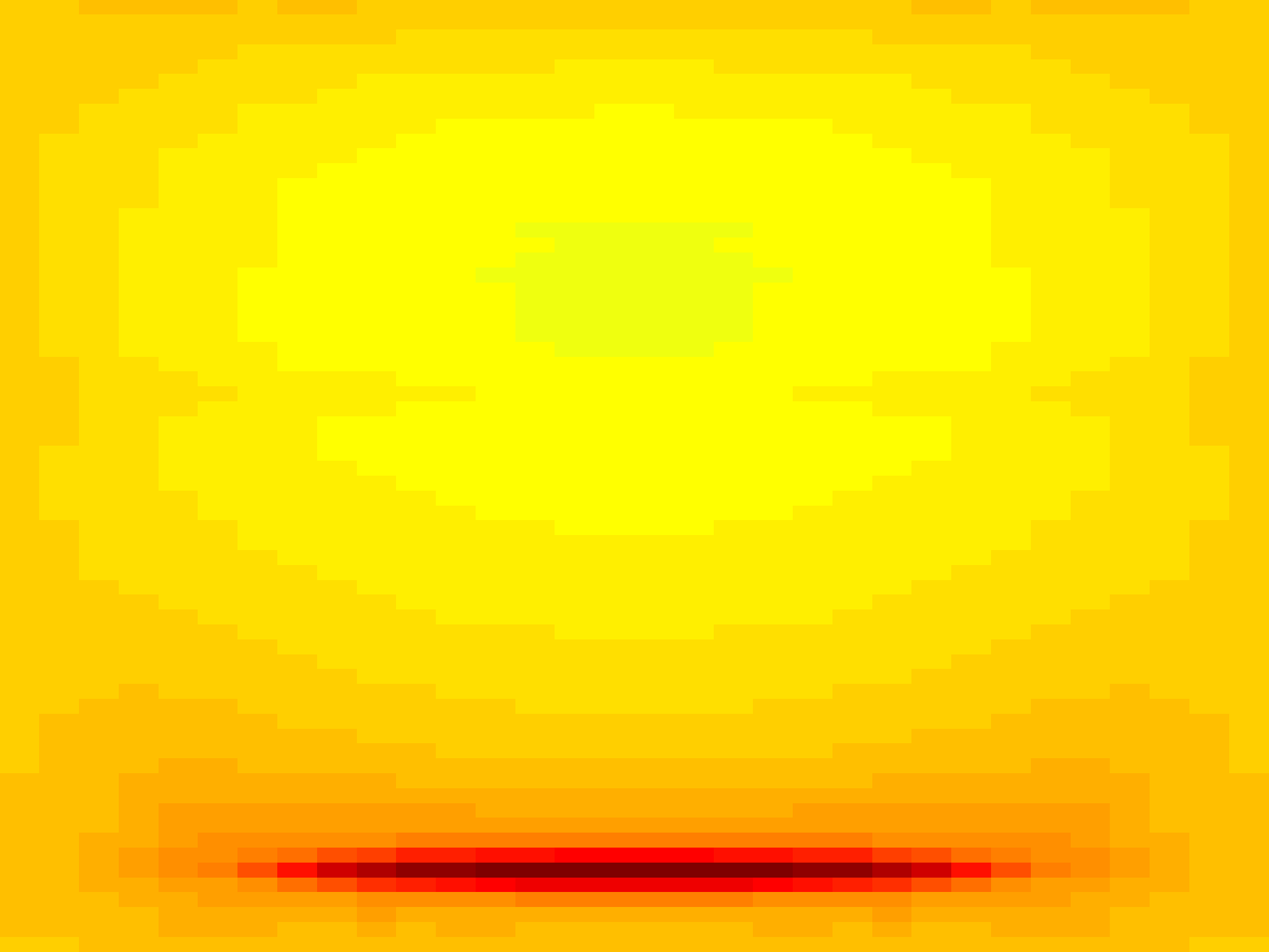}&
    \includegraphics[width=3cm]{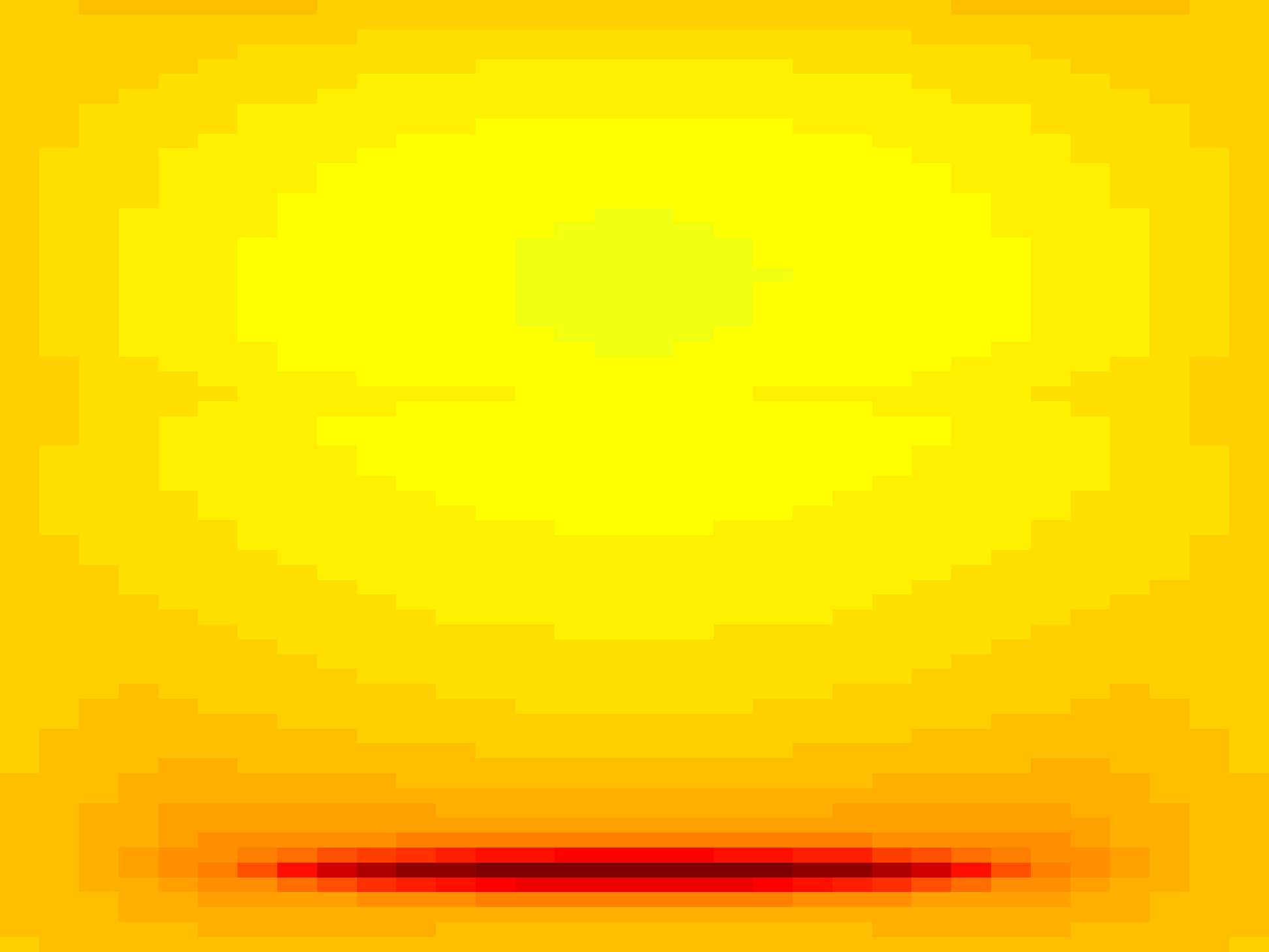}&
    \includegraphics[width=3cm]{lattice_3d_ldmm_10p_band_151}
  \end{tabular}
  \caption{Compression of the 3D lattice  data set with a 10\% data compression rate. The original angular flux at $x = 0.24$ and $x = 1.18$ are shown in the first figures on the first and third row. The results of Tucker decomposition, DCT, DFT, wavelet, and LDMM are shown in the remaining five figures.}
  \label{fig:compre_lattice_3d_10p}
\end{figure}

\begin{figure}[H]
  \centering
  \begin{tabular}{ccc}
    Original & Tucker (\textbf{78.28dB})& DCT (60.52dB)\\
    \includegraphics[width=3cm]{lattice_3d_original_band_31}&
    \includegraphics[width=3cm]{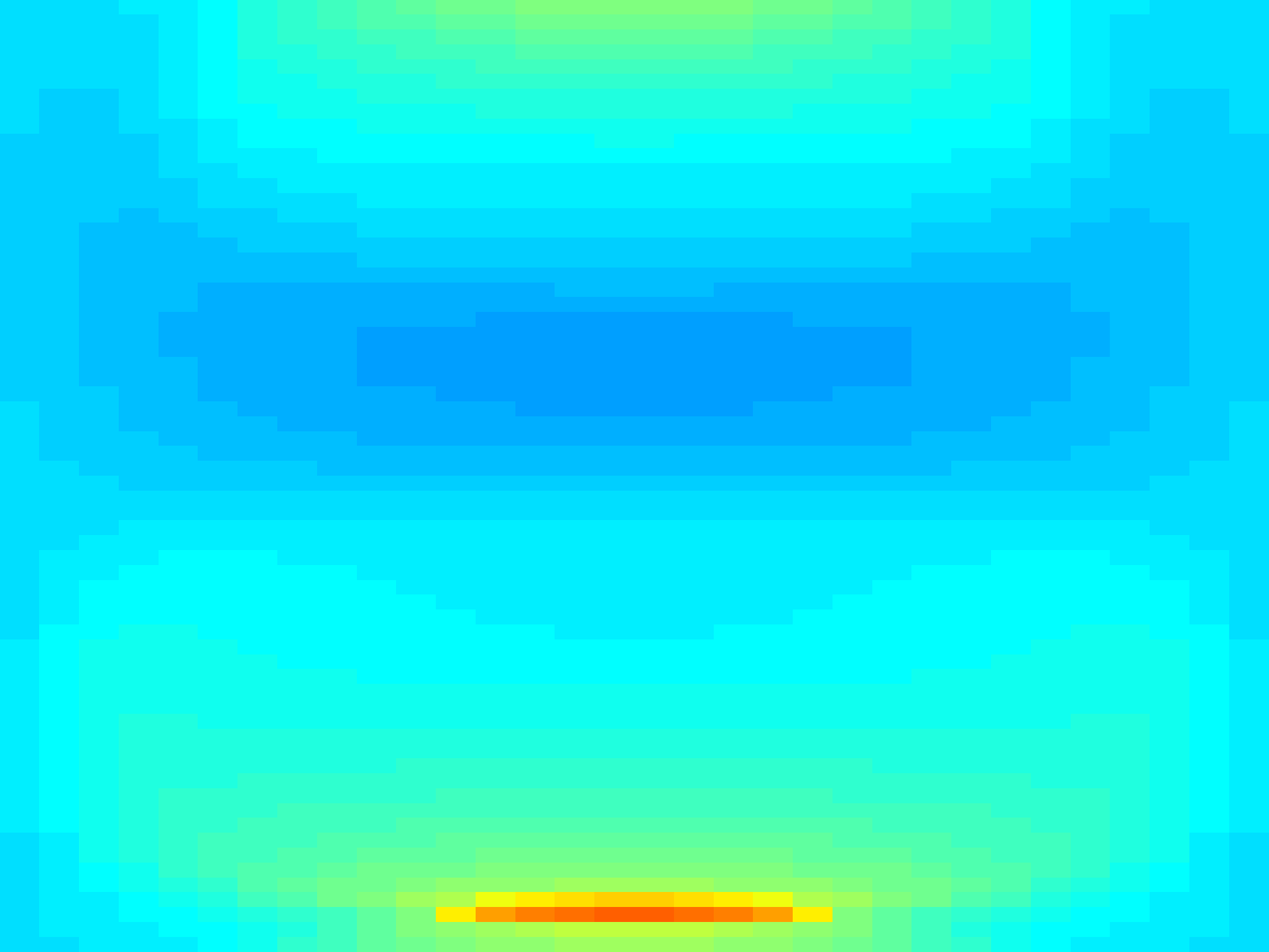}&
    \includegraphics[width=3cm]{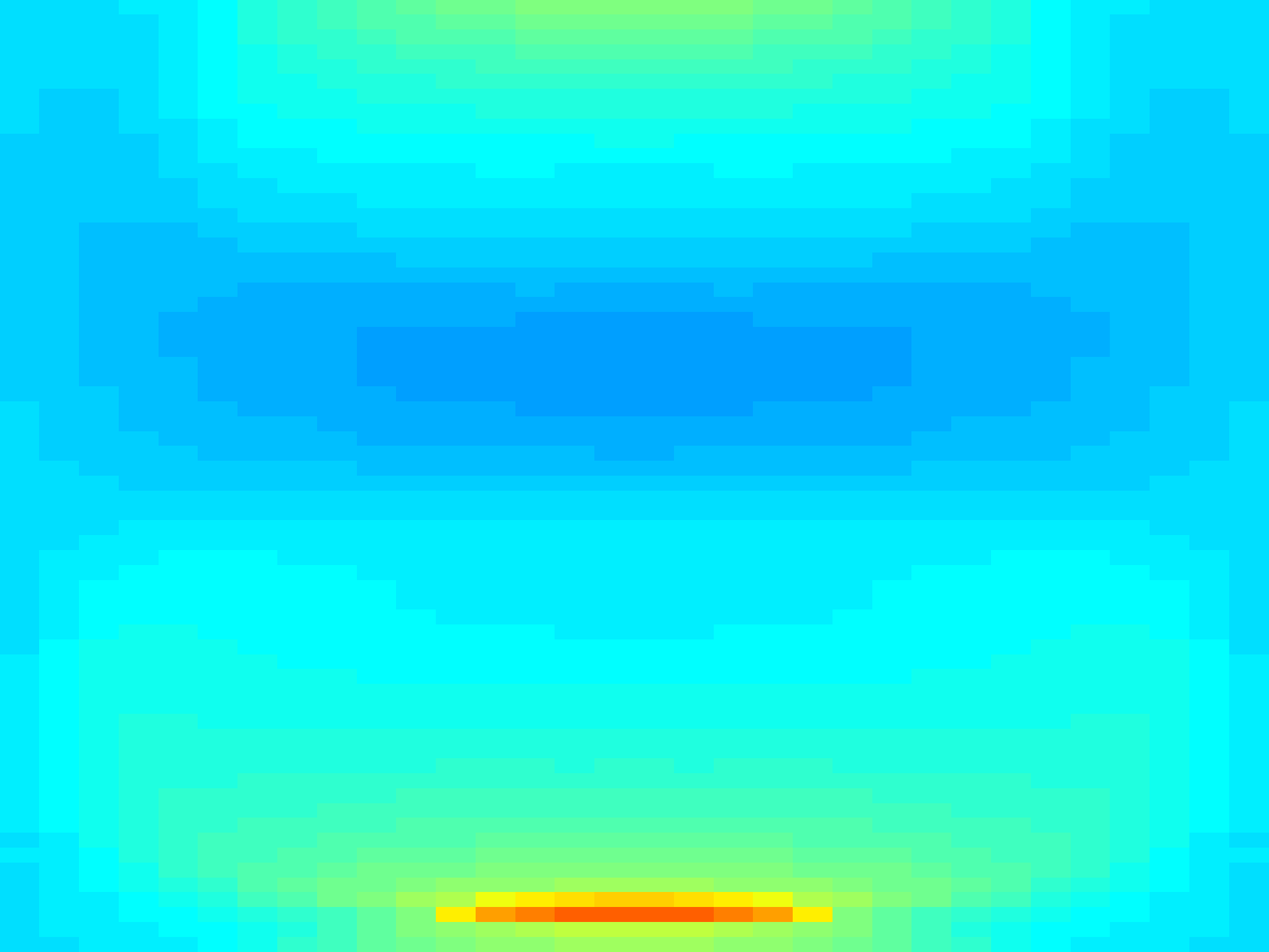}\\
    DFT (50.22dB)& Wavelet (61.25dB)& LDMM (45.82dB)\\
    \includegraphics[width=3cm]{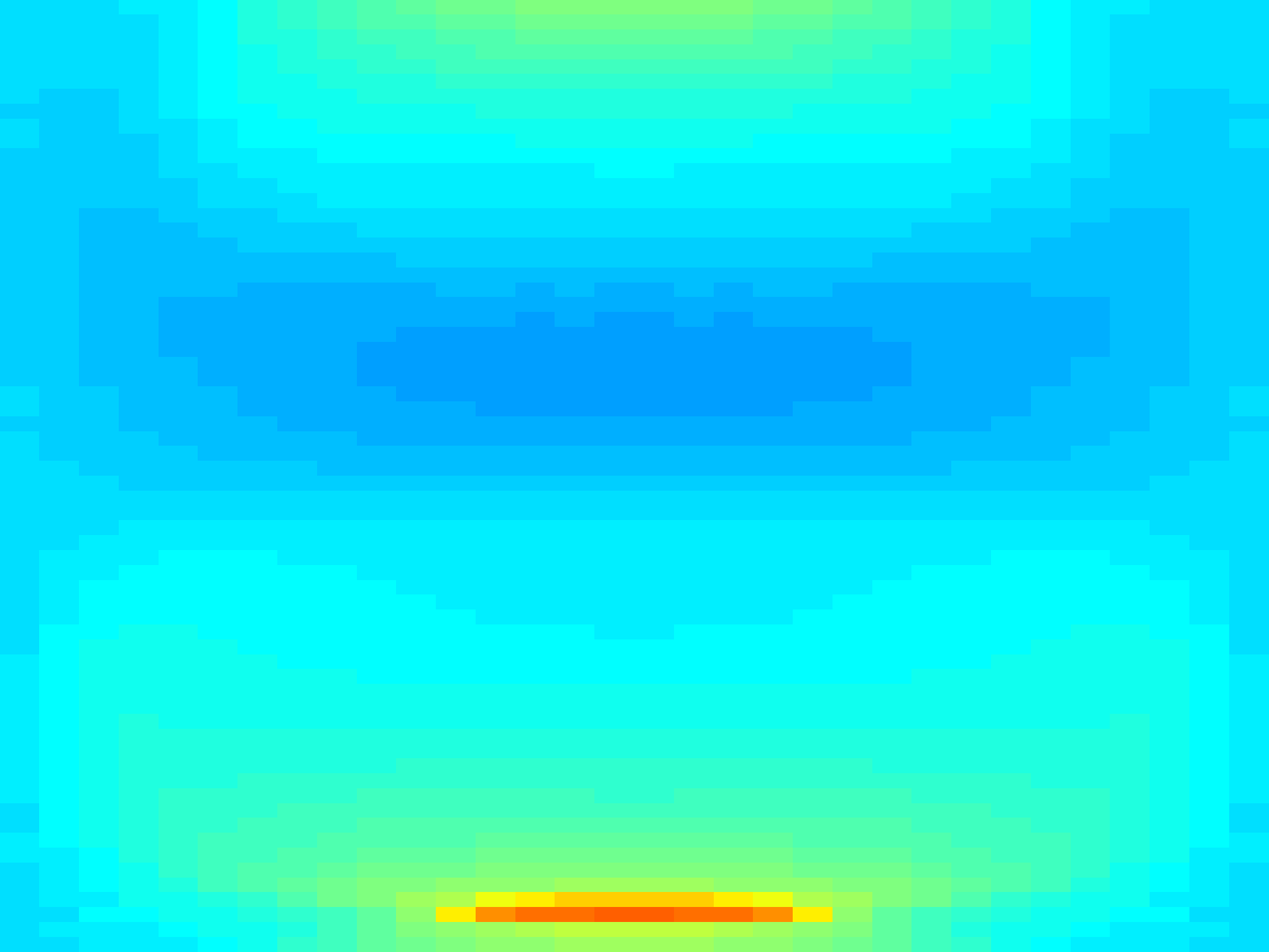}&
    \includegraphics[width=3cm]{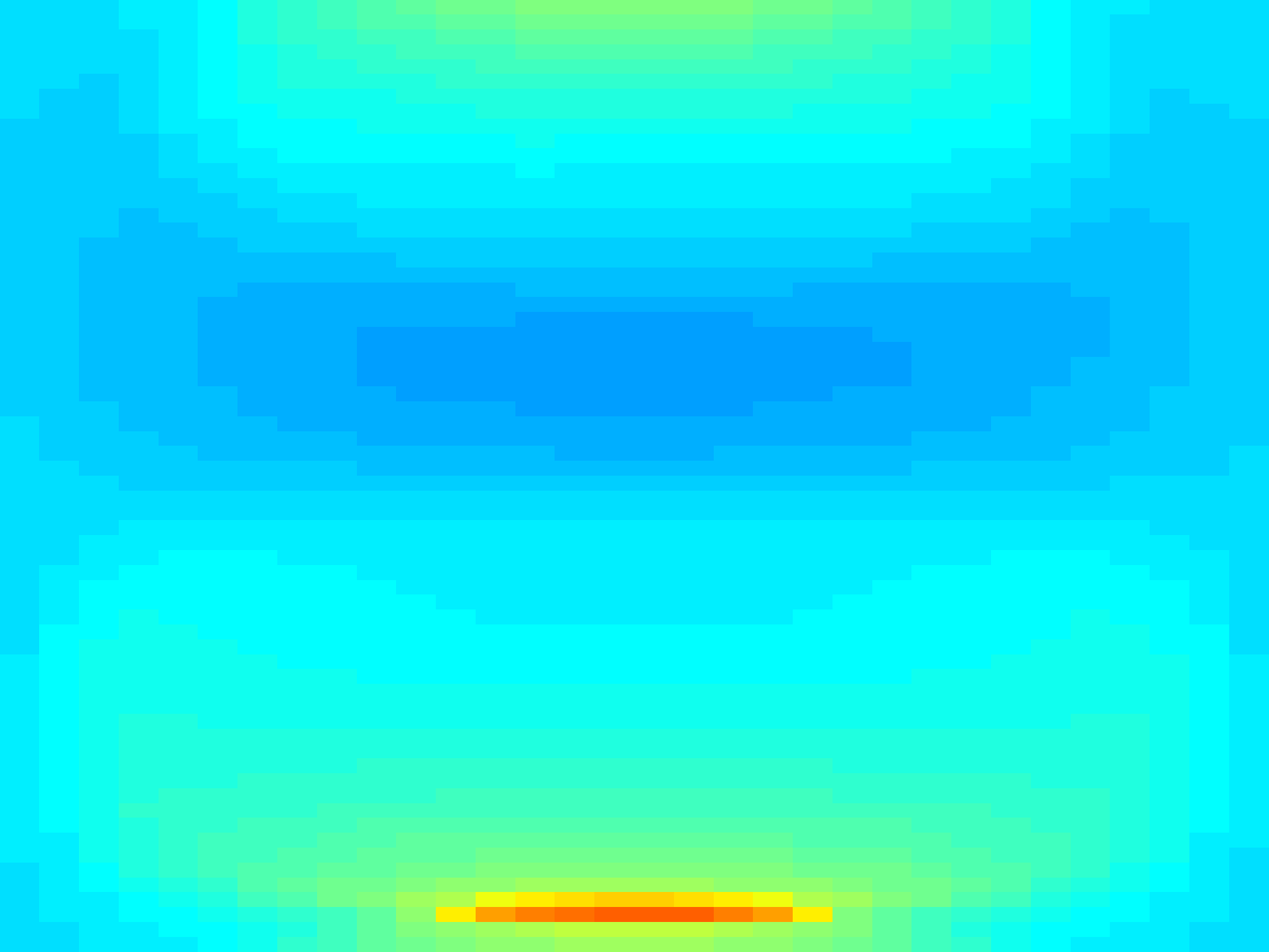}&
    \includegraphics[width=3cm]{lattice_3d_ldmm_5p_band_31}\\
    Original & Tucker (\textbf{78.28dB})& DCT (60.52dB)\\
    \includegraphics[width=3cm]{lattice_3d_original_band_151}&
    \includegraphics[width=3cm]{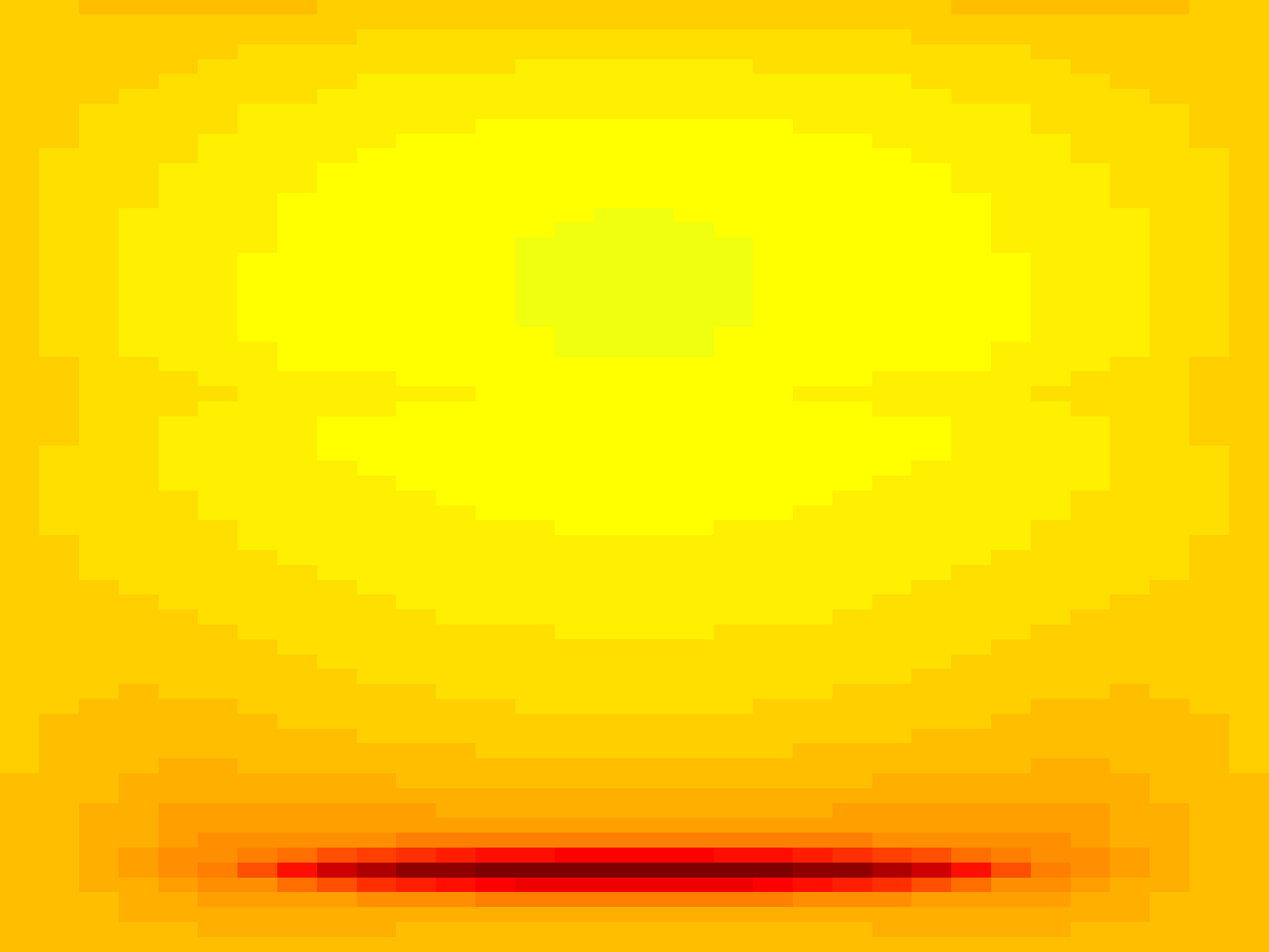}&
    \includegraphics[width=3cm]{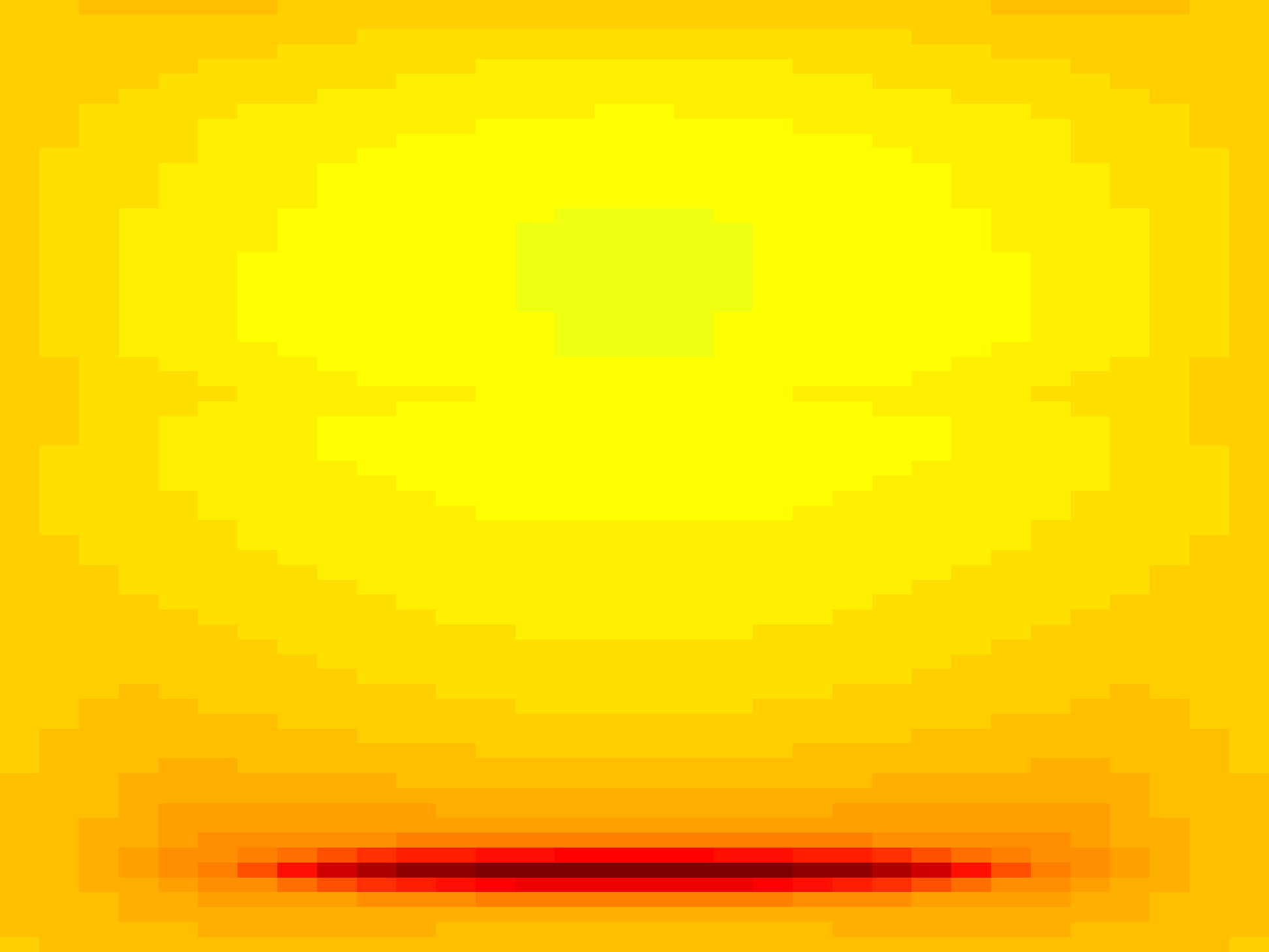}\\
    DFT (50.22dB)& Wavelet (61.25dB)& LDMM (45.82dB)\\
    \includegraphics[width=3cm]{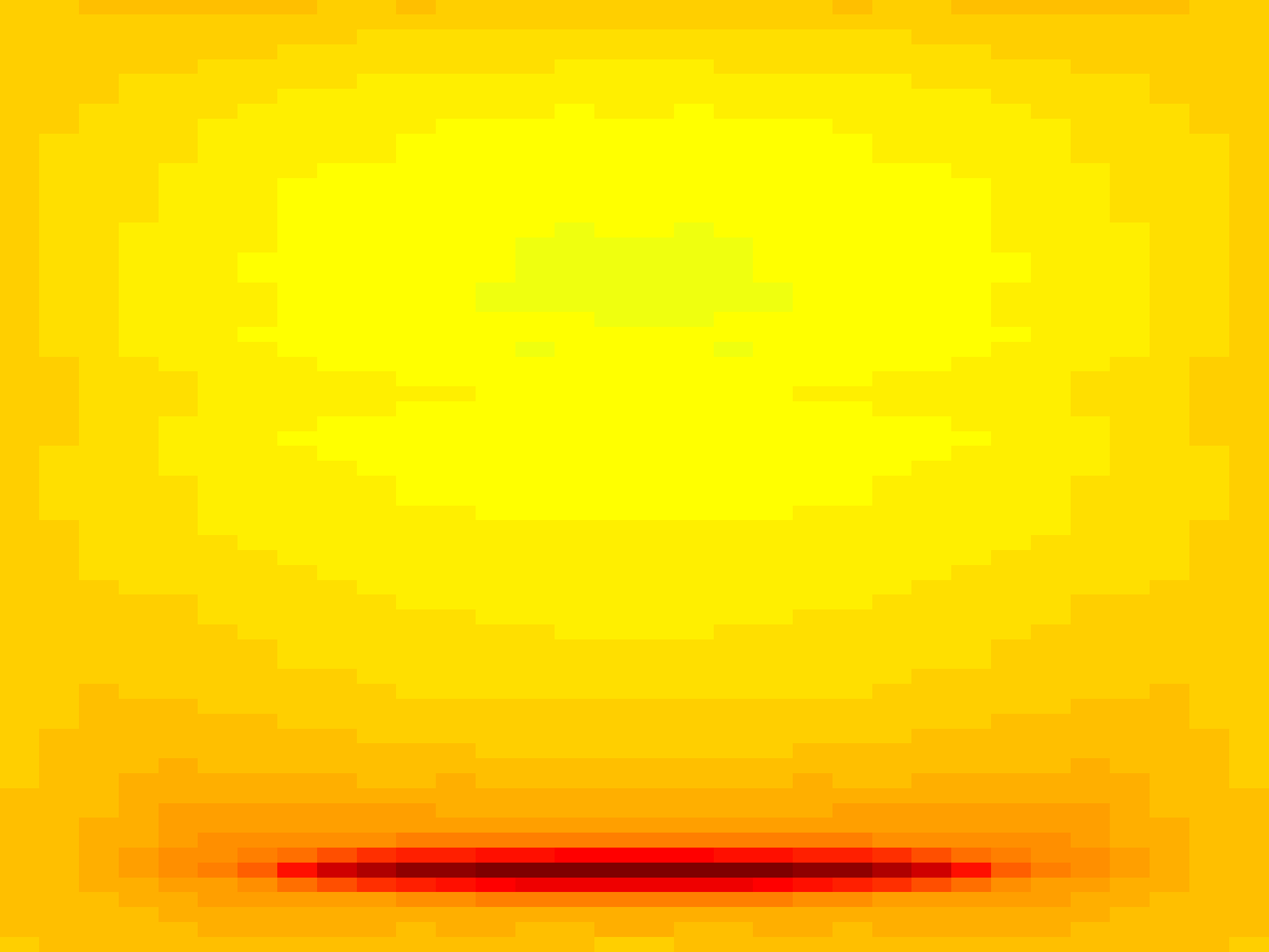}&
    \includegraphics[width=3cm]{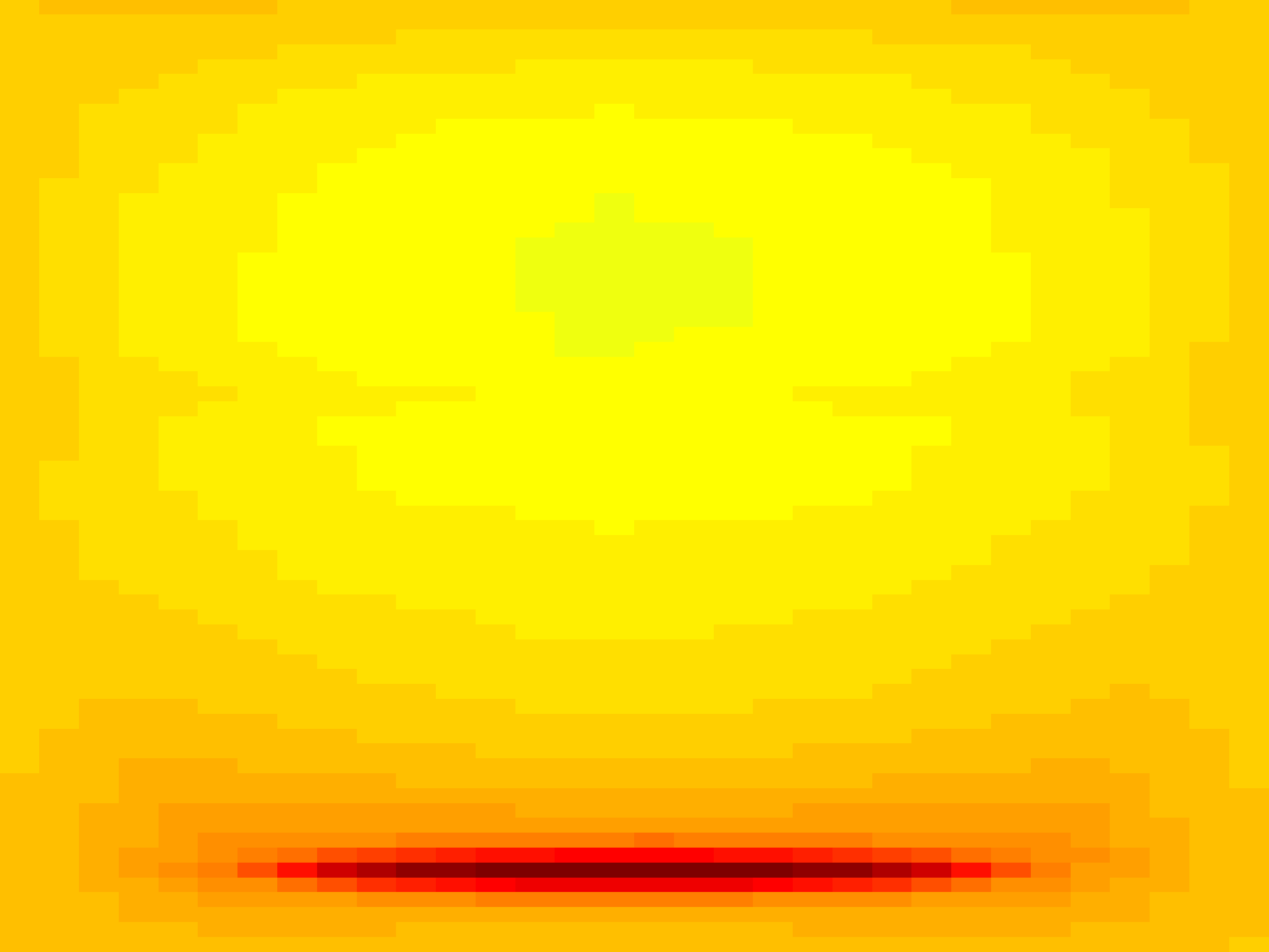}&
    \includegraphics[width=3cm]{lattice_3d_ldmm_5p_band_151}
  \end{tabular}
  \caption{Compression of the 3D lattice  data set with a 5\% data compression rate. The original angular flux at $x = 0.24$ and $x = 1.18$ are shown in the first figures on the first and third row. The results of Tucker decomposition, DCT, DFT, wavelet, and LDMM are shown in the remaining five figures.}
  \label{fig:compre_lattice_3d_5p}
\end{figure}

\begin{figure}[H]
  \centering
  \begin{tabular}{ccc}
    Original & Tucker (43.89dB)& DCT (\textbf{45.65dB})\\
    \includegraphics[width=3cm]{shock_3d_original_band_19}&
    \includegraphics[width=3cm]{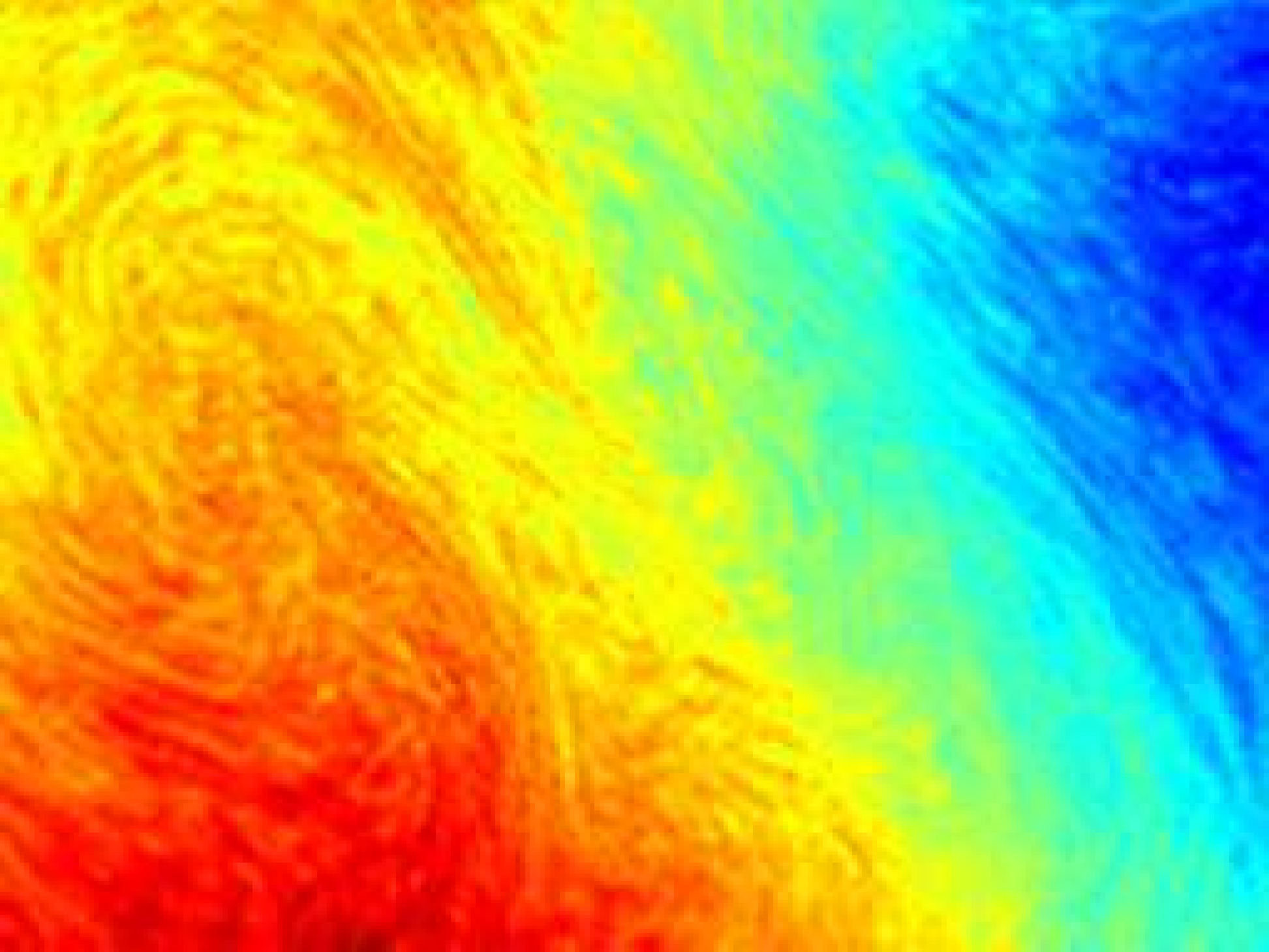}&
    \includegraphics[width=3cm]{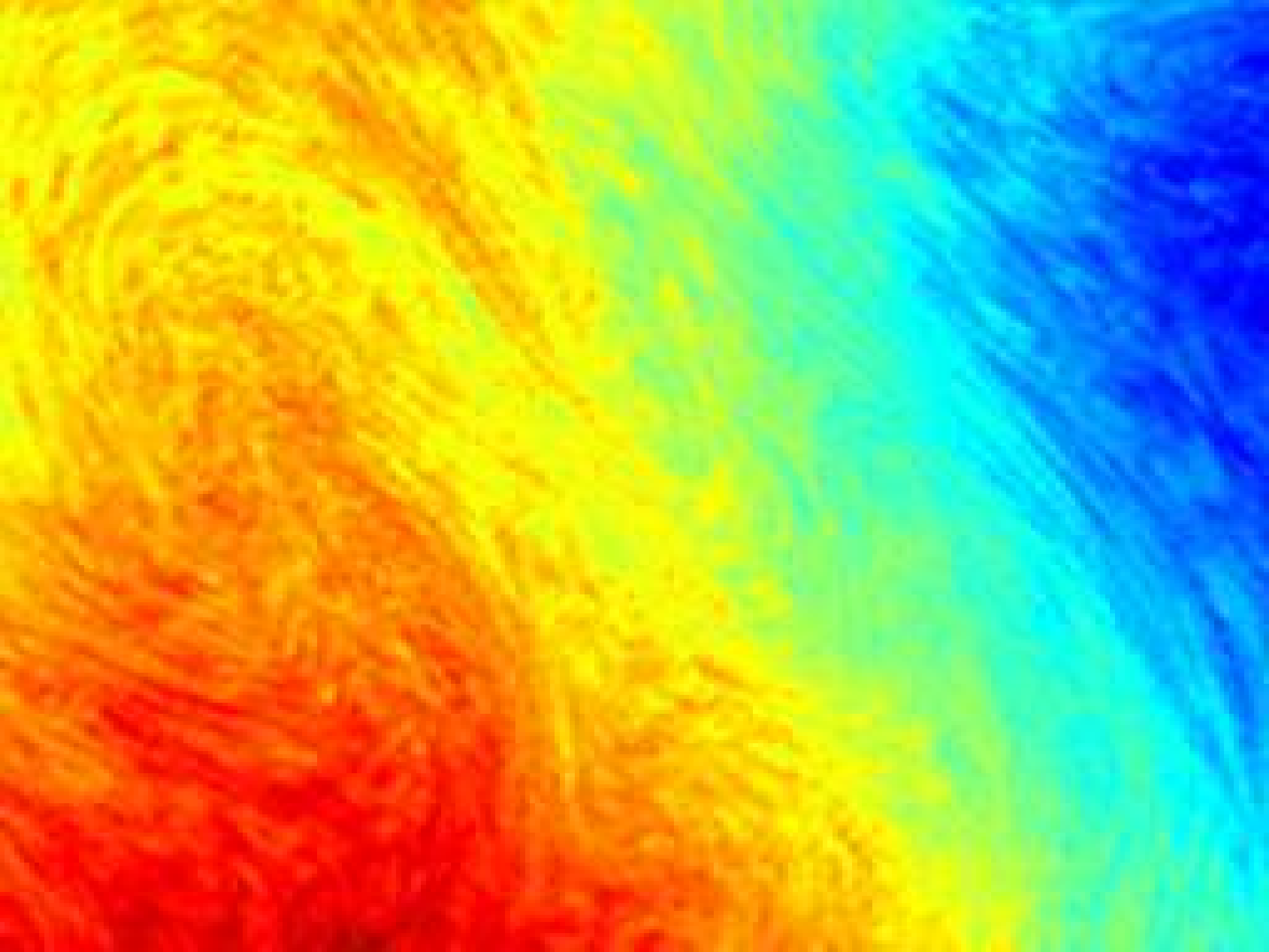}\\
    DFT (44.26dB)& Wavelet (45.17dB)& LDMM (39.54)\\
    \includegraphics[width=3cm]{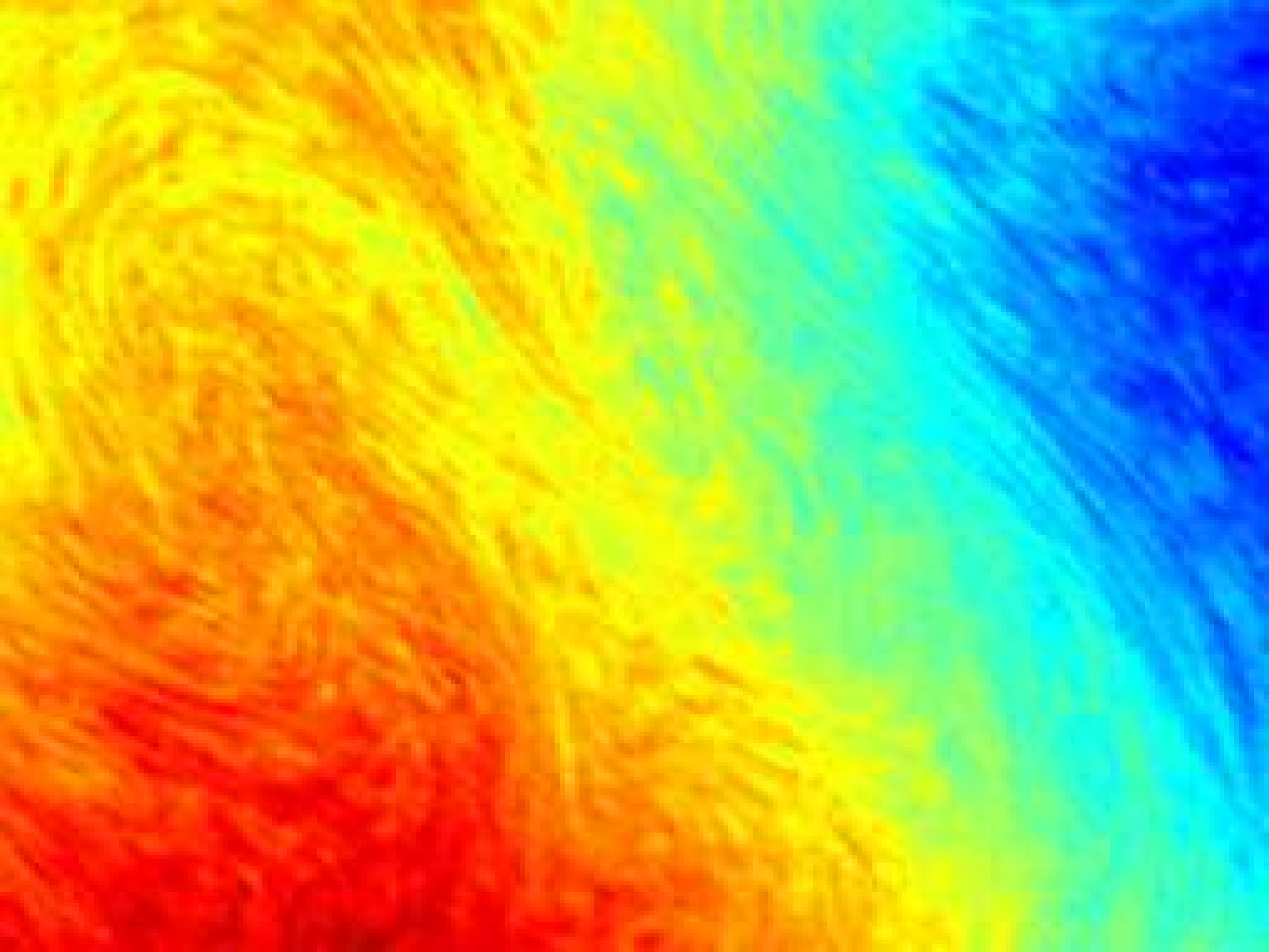}&
    \includegraphics[width=3cm]{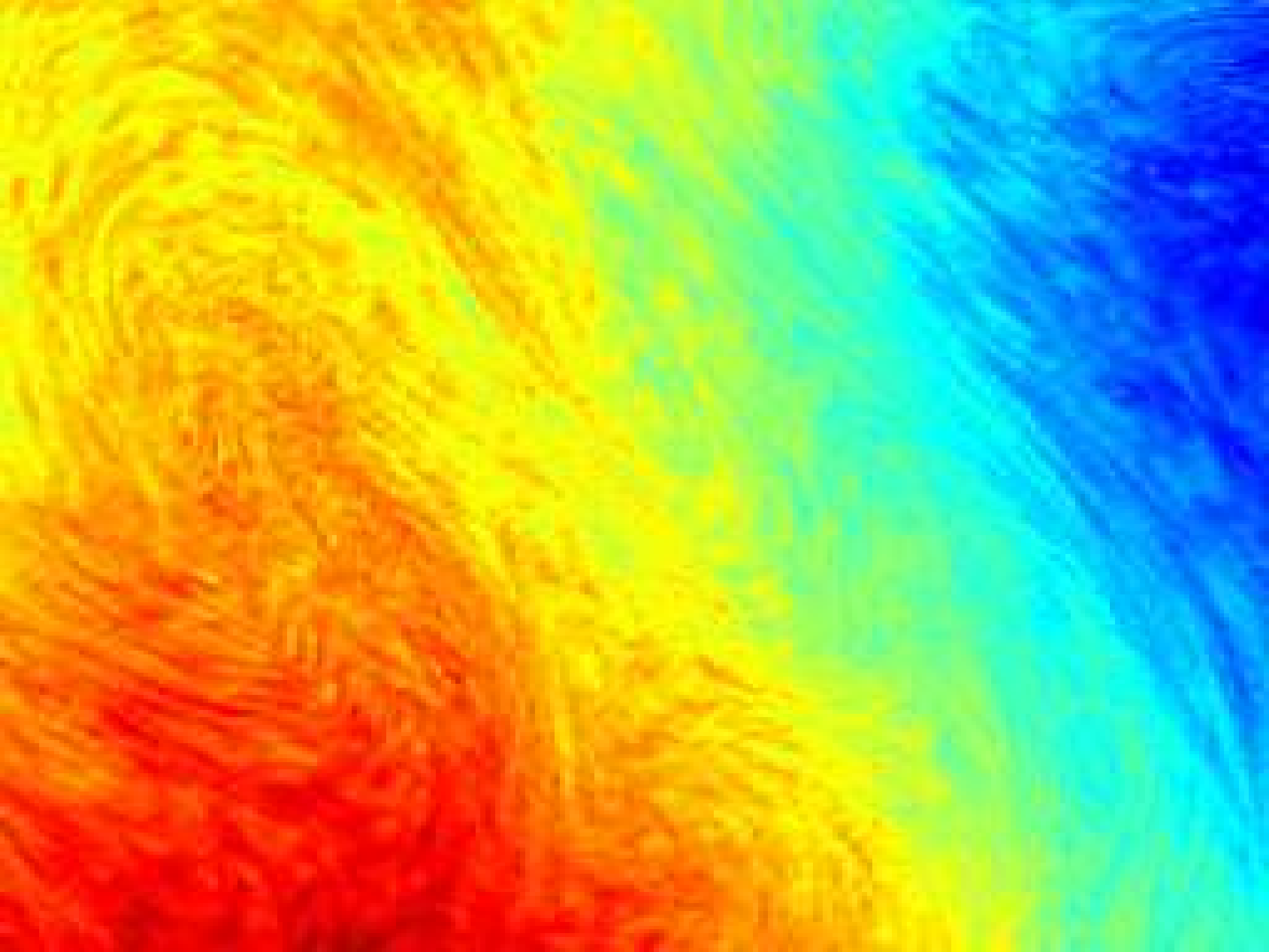}&
    \includegraphics[width=3cm]{shock_3d_ldmm_10p_band_19}\\
    Original & Tucker (43.89dB)& DCT (\textbf{45.65dB})\\
    \includegraphics[width=3cm]{shock_3d_original_band_29}&
    \includegraphics[width=3cm]{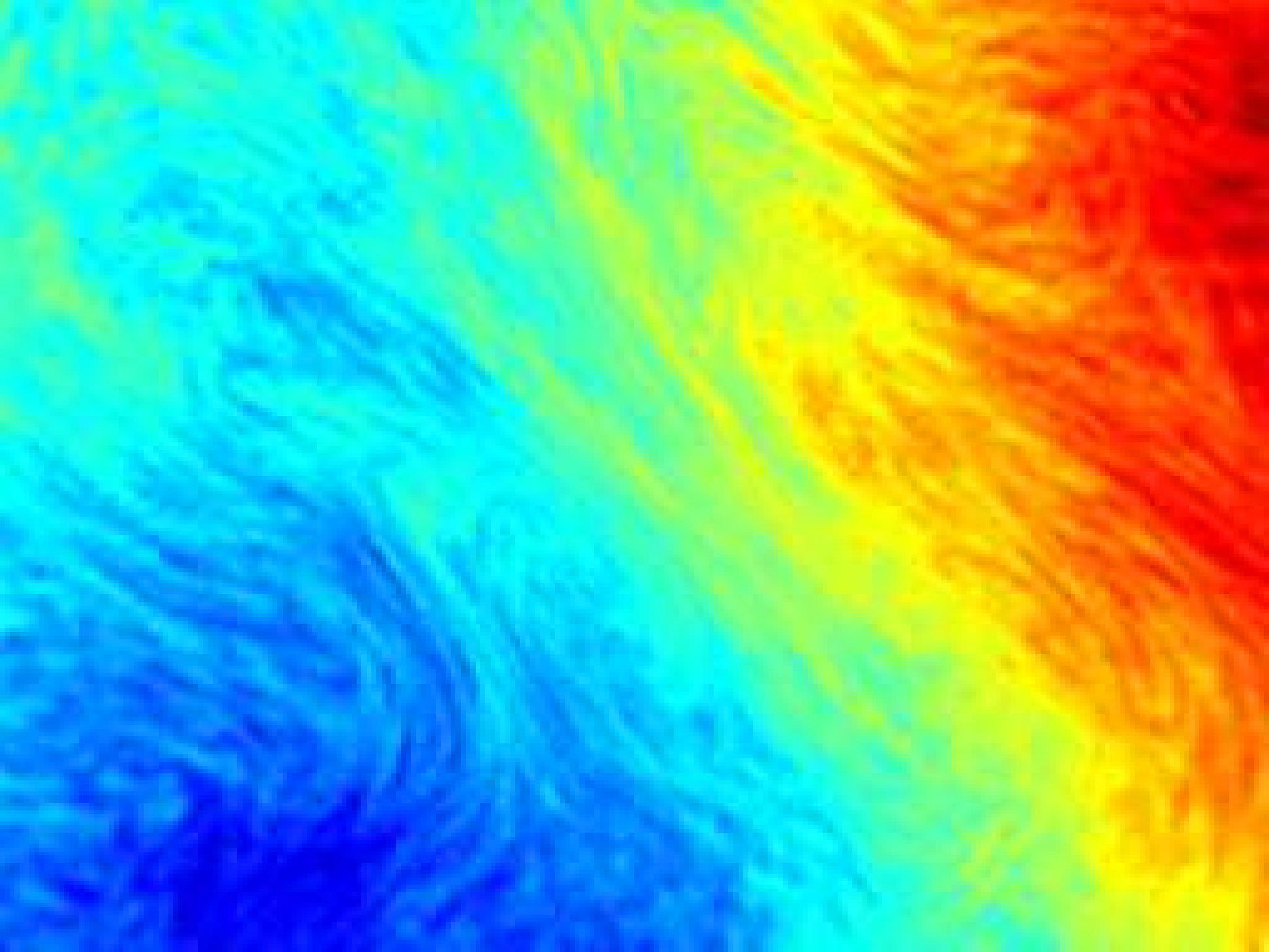}&
    \includegraphics[width=3cm]{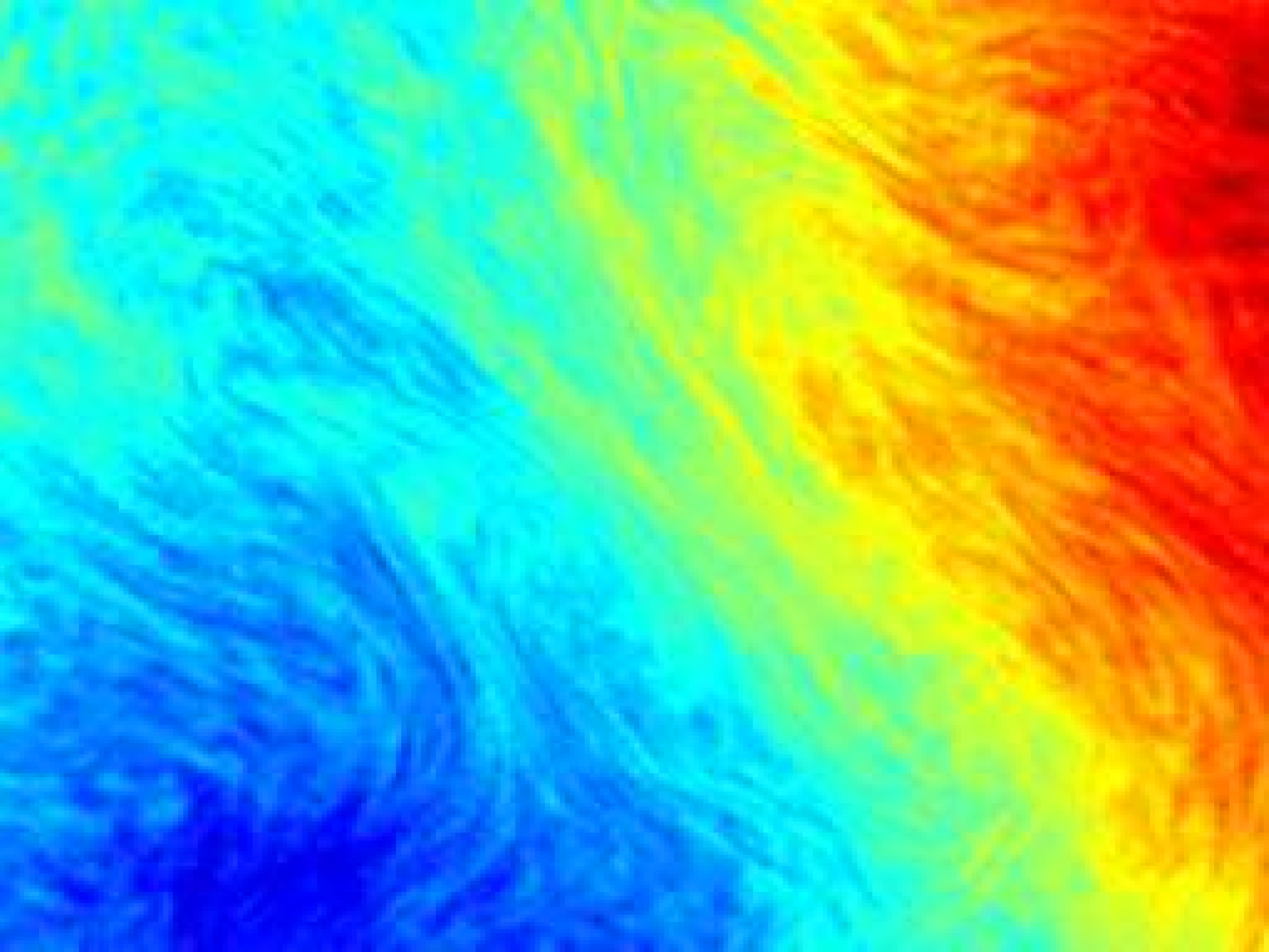}\\
    DFT (44.26dB)& Wavelet (45.17dB)& LDMM (39.54)\\
    \includegraphics[width=3cm]{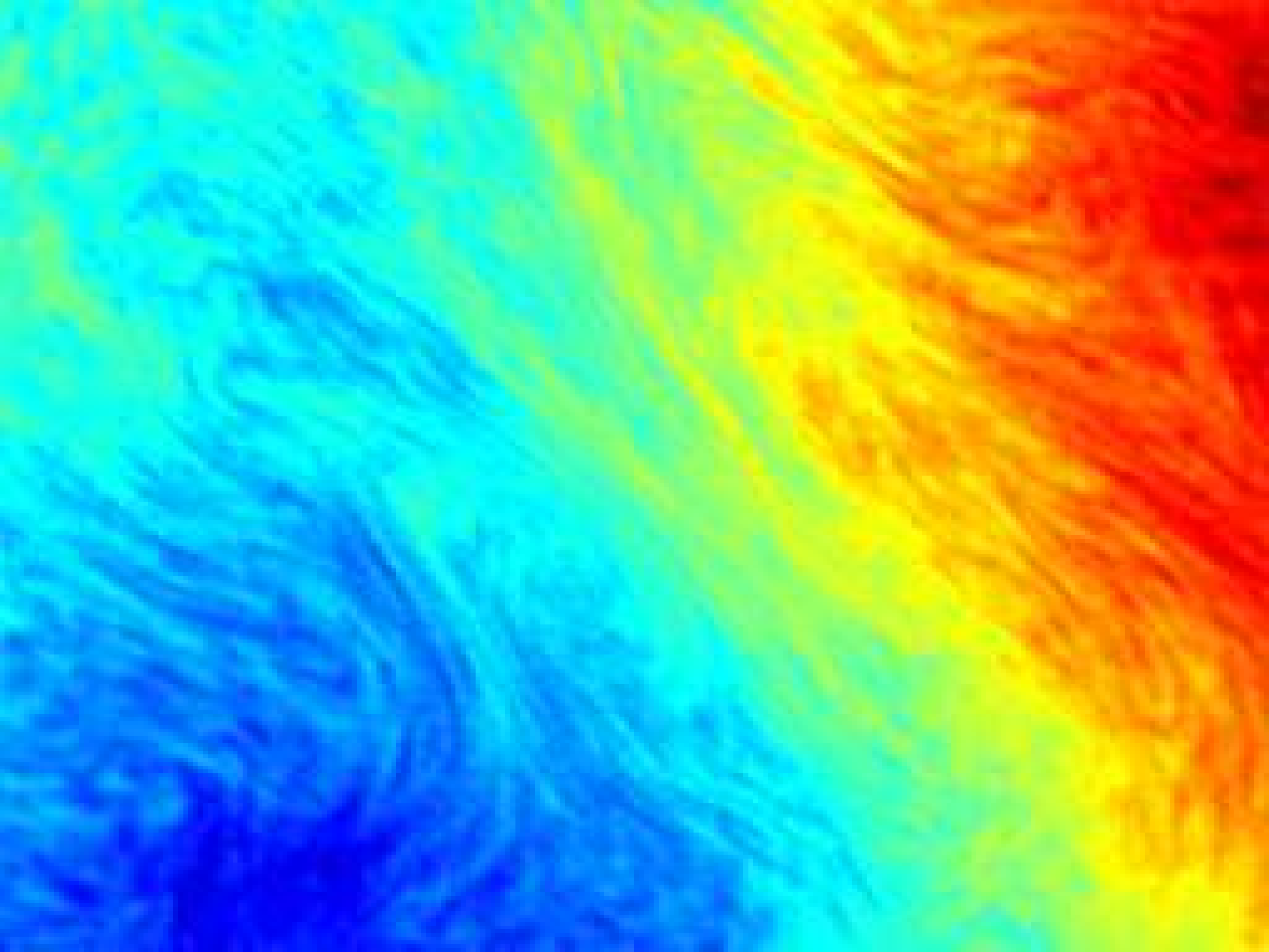}&
    \includegraphics[width=3cm]{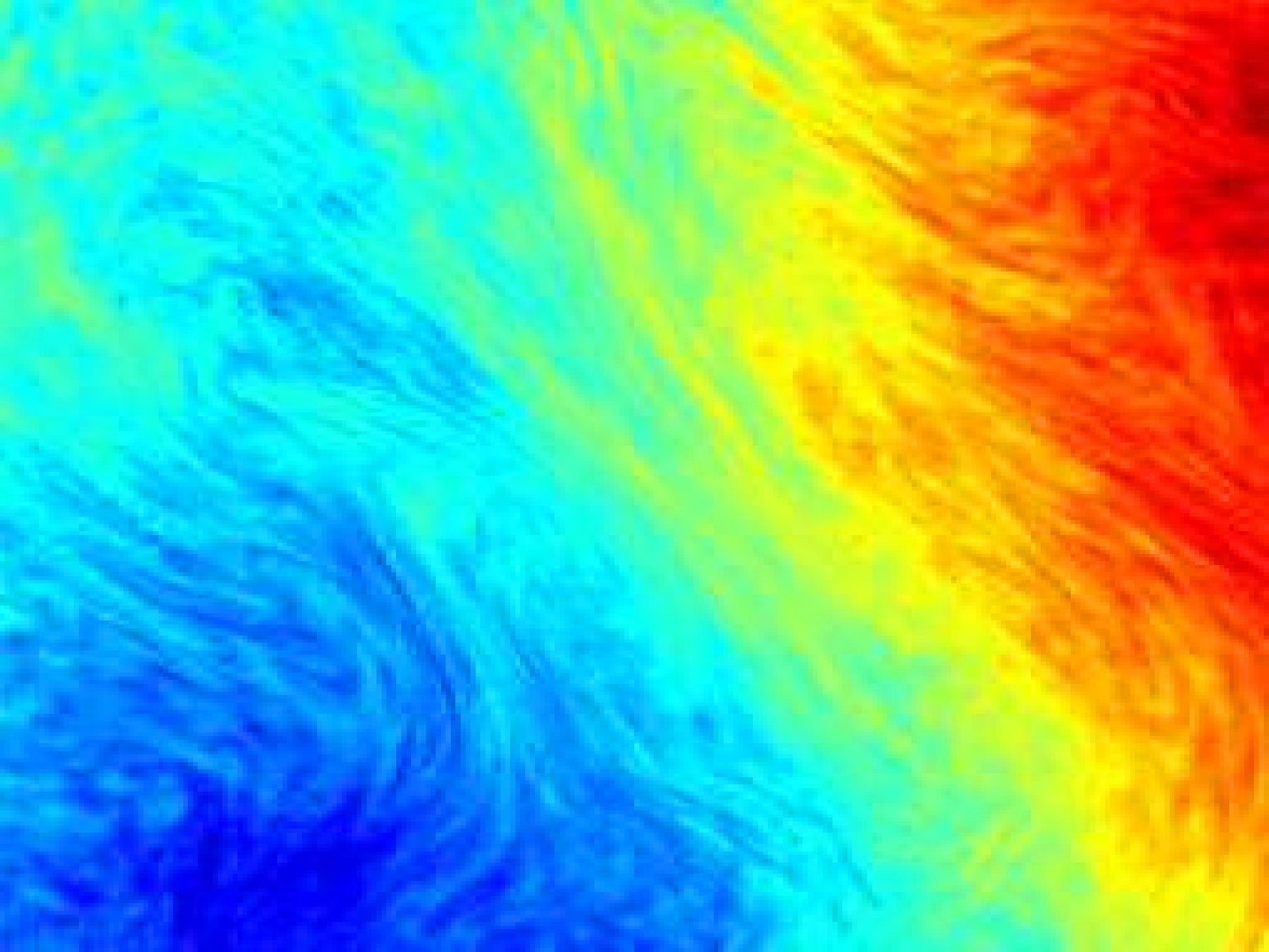}&
    \includegraphics[width=3cm]{shock_3d_ldmm_10p_band_29}
  \end{tabular}
  \caption{Compression of the 3D plasma (distribution function) data set with a 10\% data compression rate. Two  spatial cross  sections of the original data set are shown in the first figures on the first and third row. The results of Tucker decomposition, DCT, DFT, wavelet, and LDMM are shown in the remaining five figures.}
  \label{fig:compre_shock_3d_10p}
\end{figure}

\begin{figure}[H]
  \centering
  \begin{tabular}{ccc}
    Original & Tucker (40.75dB)& DCT (42.29dB)\\
    \includegraphics[width=3cm]{shock_3d_original_band_19}&
    \includegraphics[width=3cm]{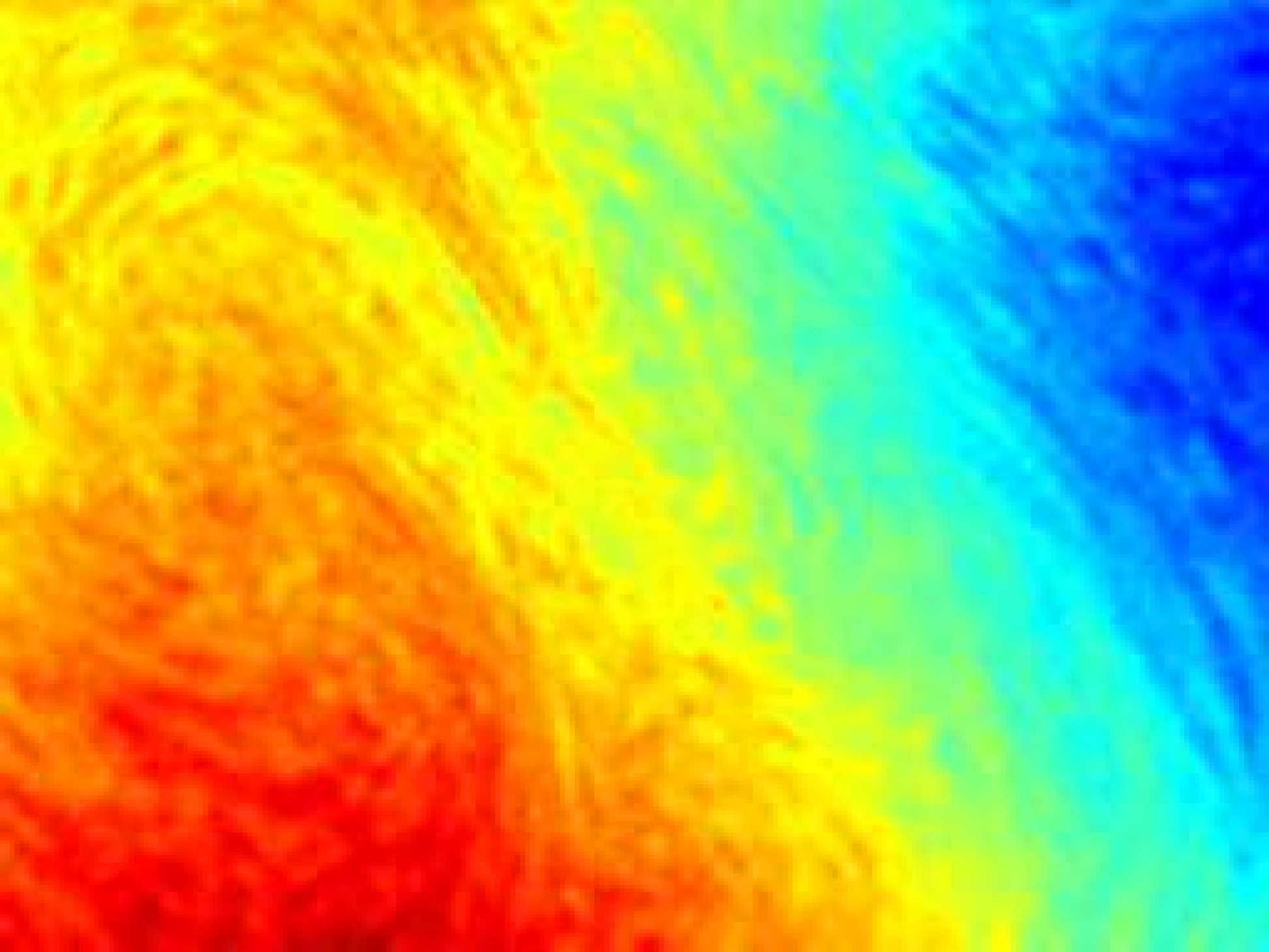}&
    \includegraphics[width=3cm]{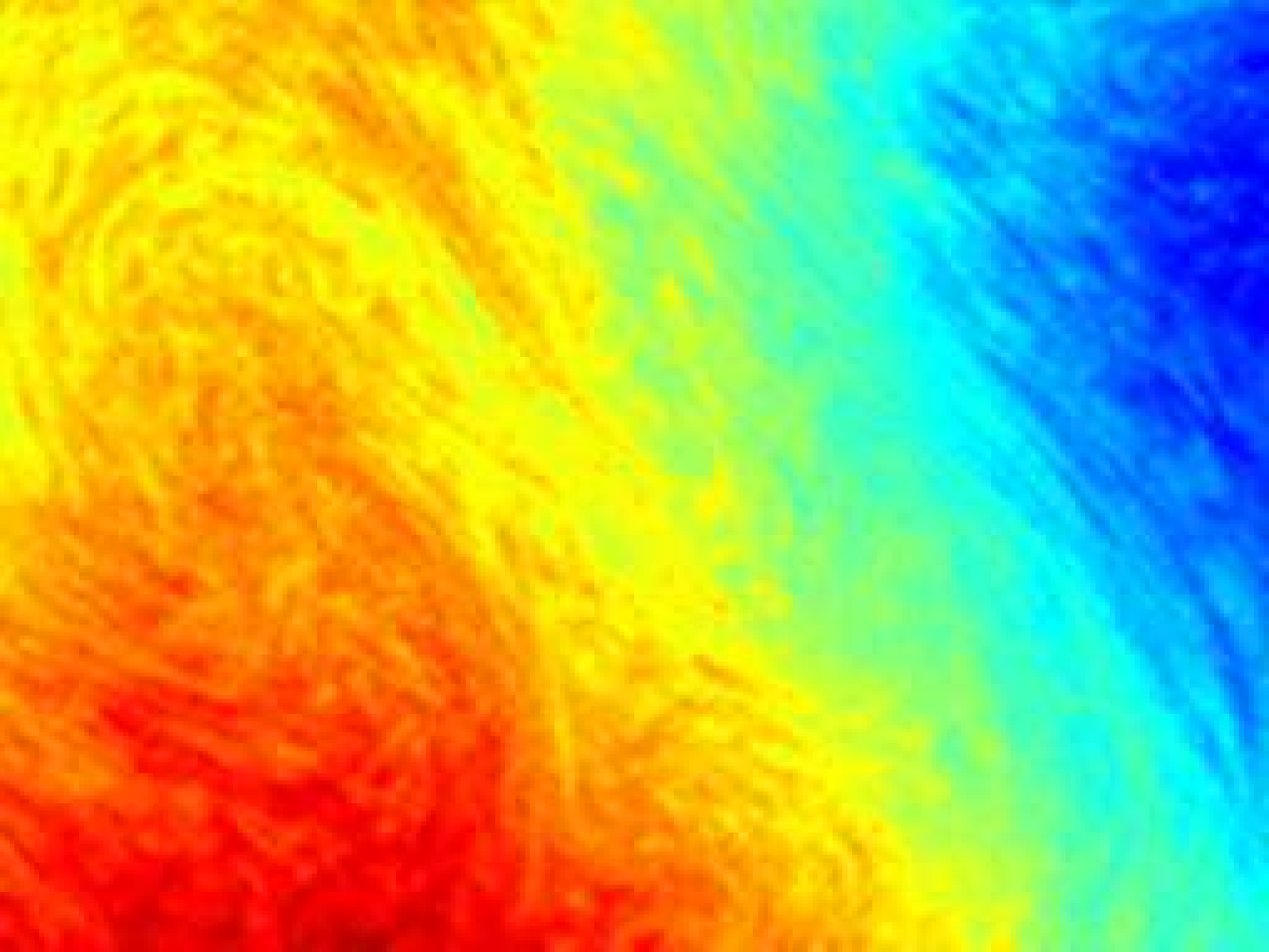}\\
    DFT (41.28dB)& Wavelet (40.97dB)& LDMM (37.72)\\
    \includegraphics[width=3cm]{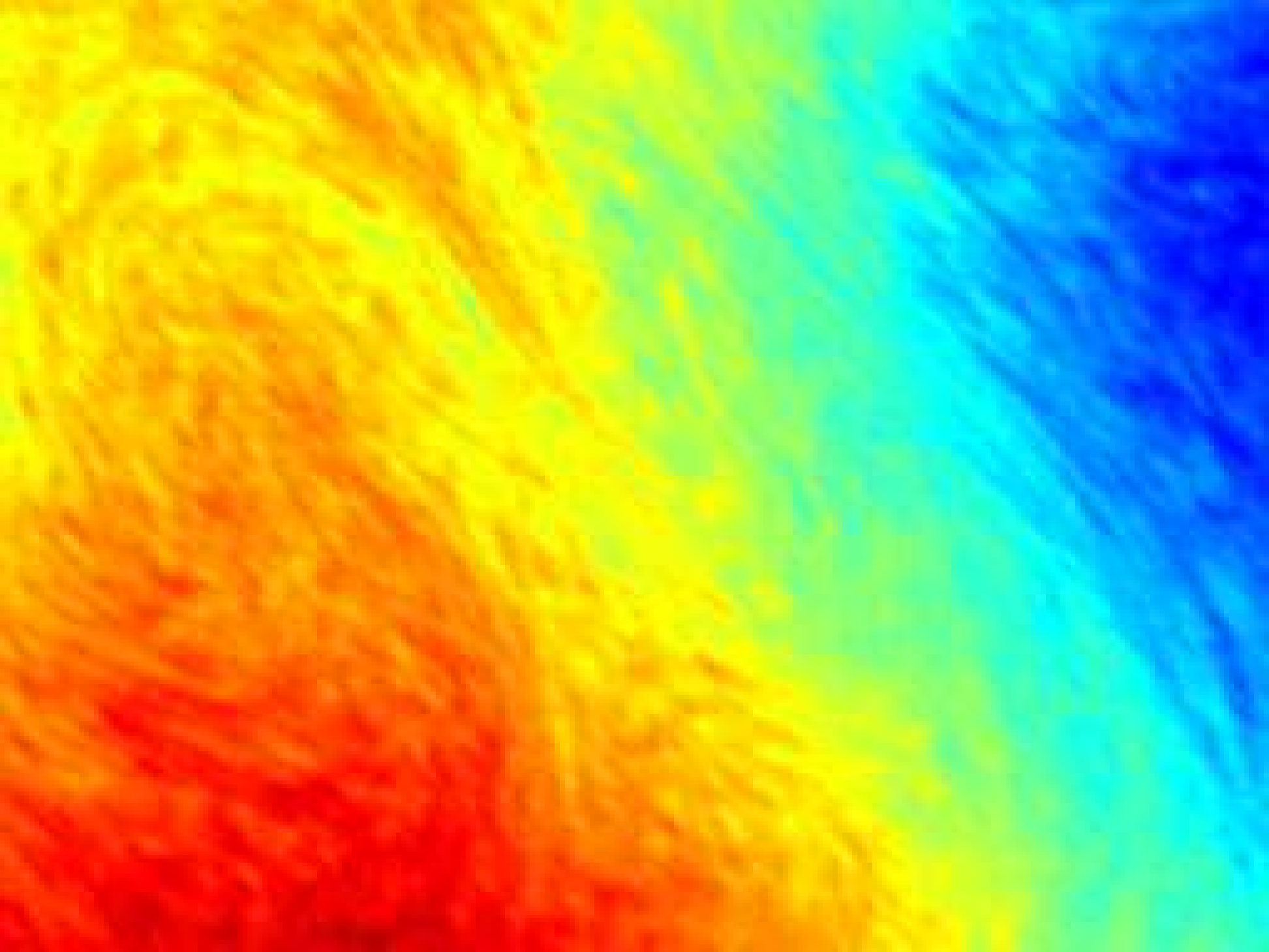}&
    \includegraphics[width=3cm]{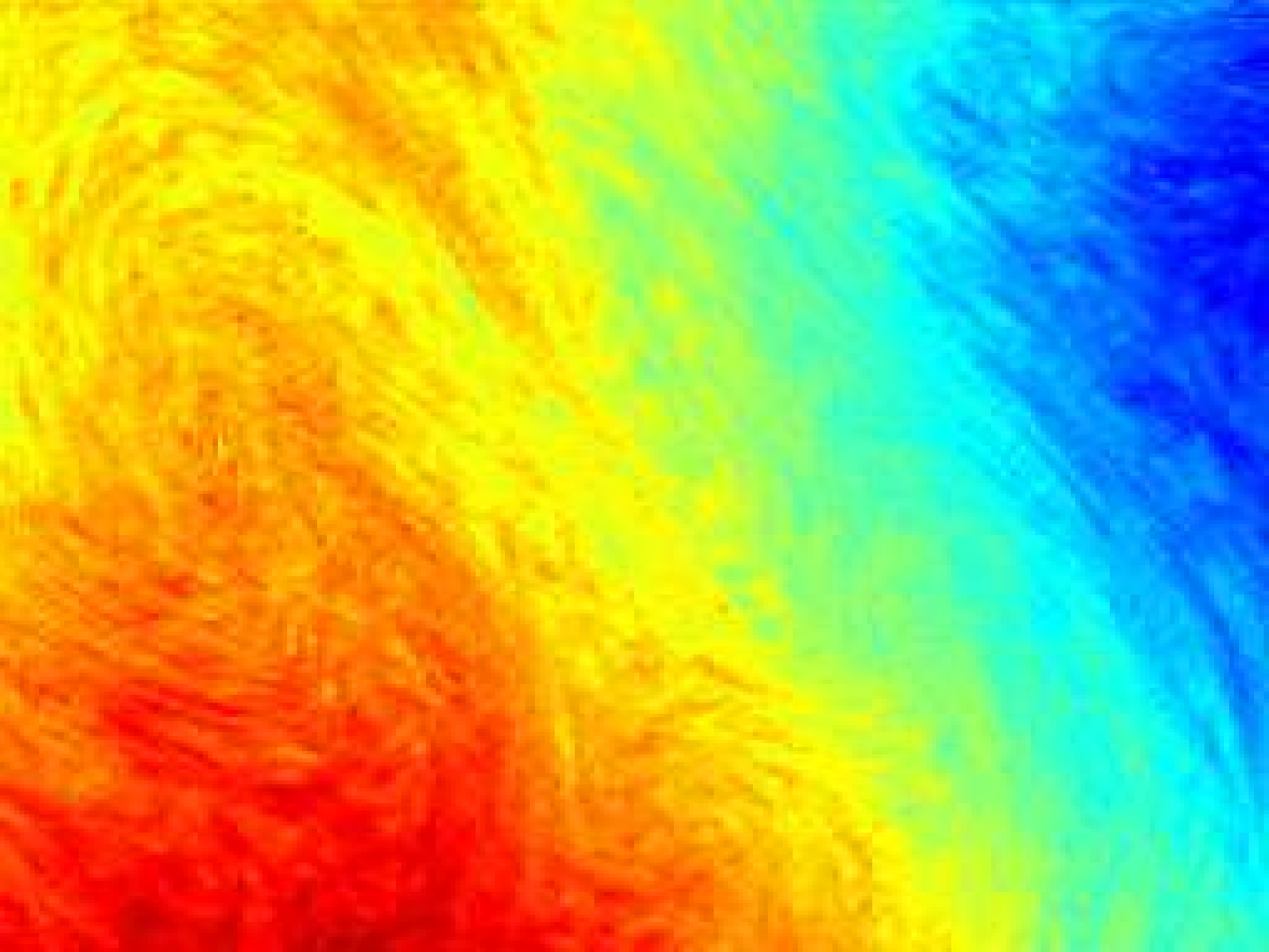}&
    \includegraphics[width=3cm]{shock_3d_ldmm_5p_band_19}\\
    Original & Tucker (40.75dB)& DCT (42.29dB)\\
    \includegraphics[width=3cm]{shock_3d_original_band_29}&
    \includegraphics[width=3cm]{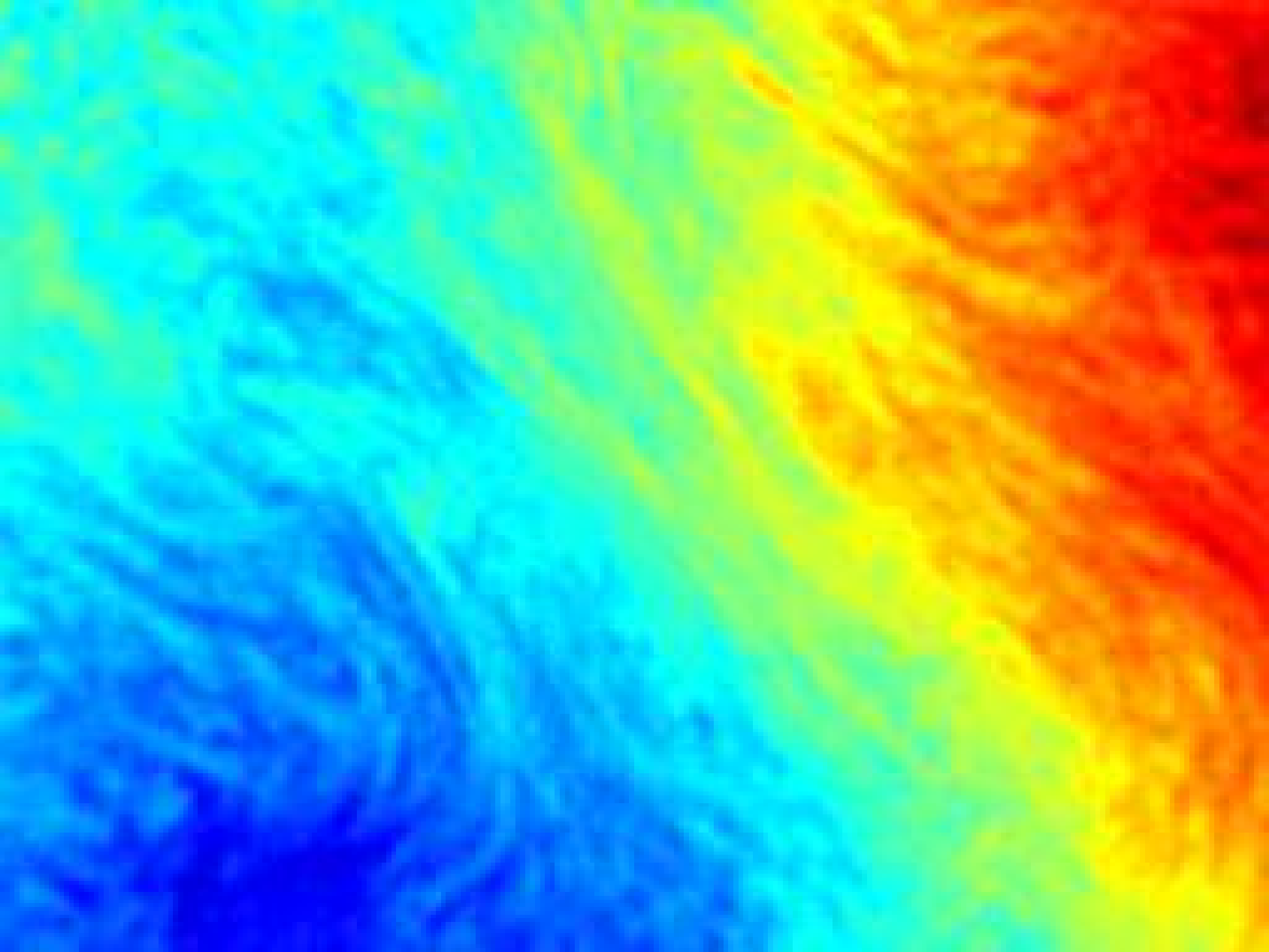}&
    \includegraphics[width=3cm]{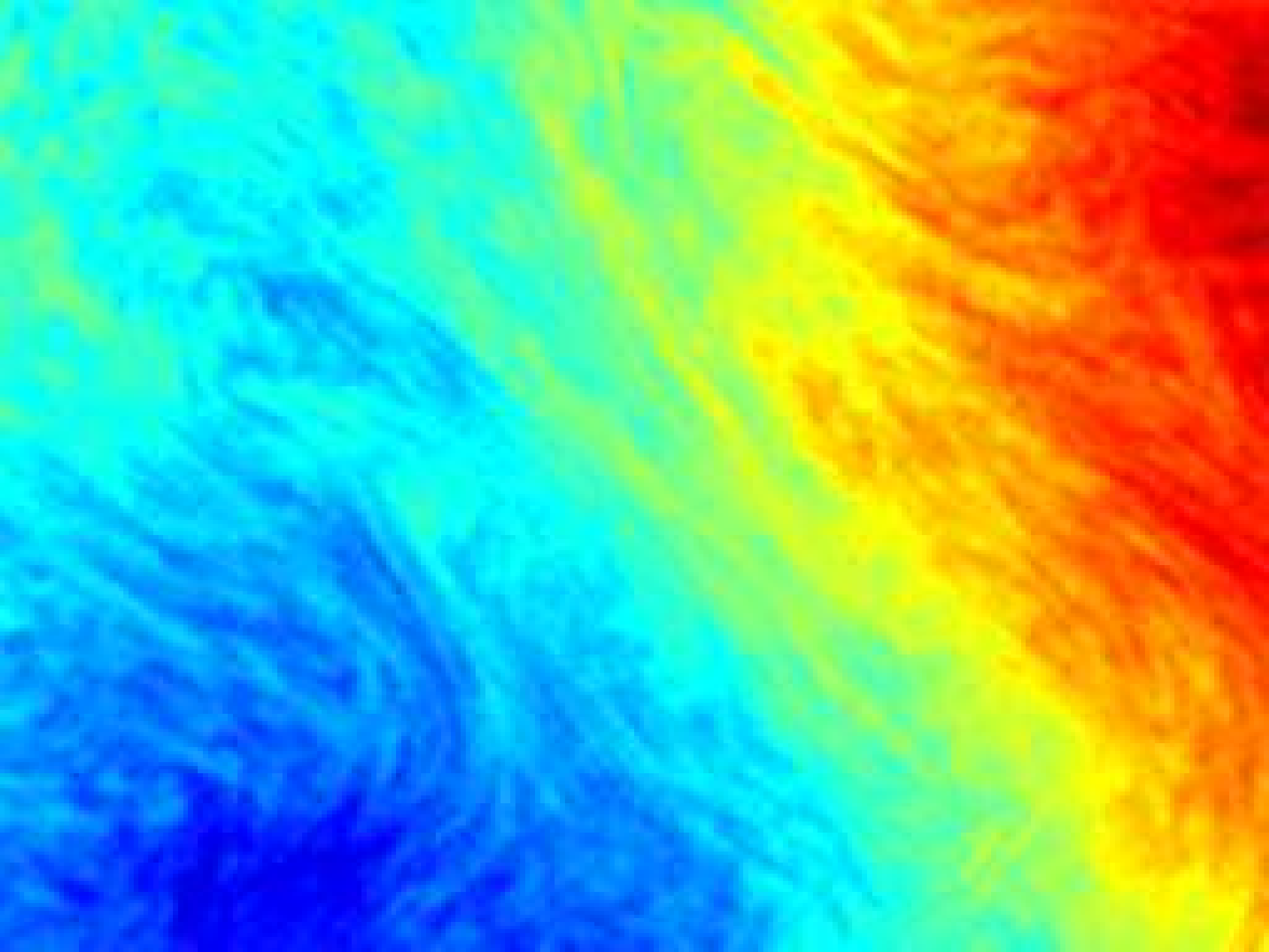}\\
    DFT (41.28dB)& Wavelet (40.97dB)& LDMM (37.72)\\
    \includegraphics[width=3cm]{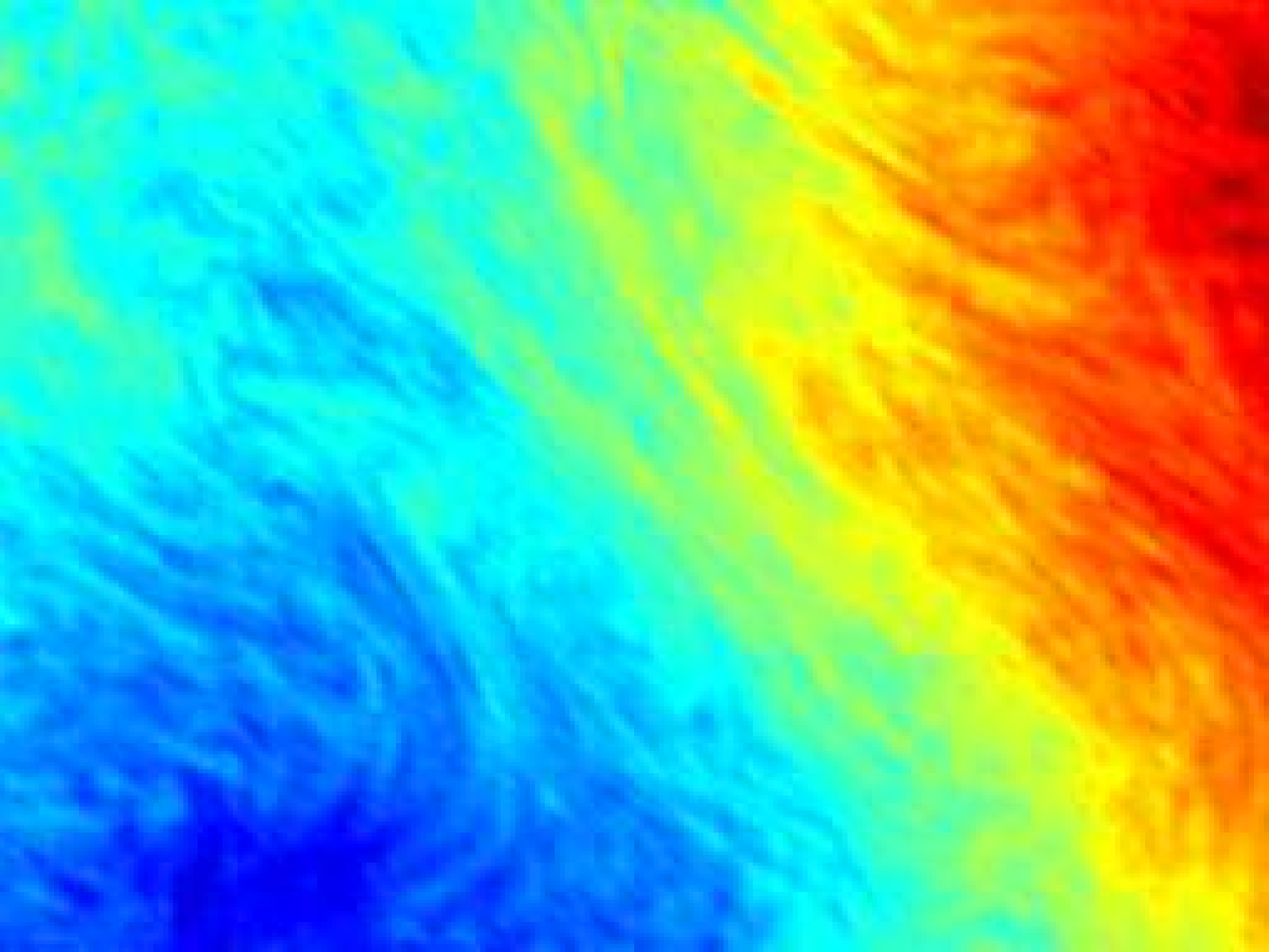}&
    \includegraphics[width=3cm]{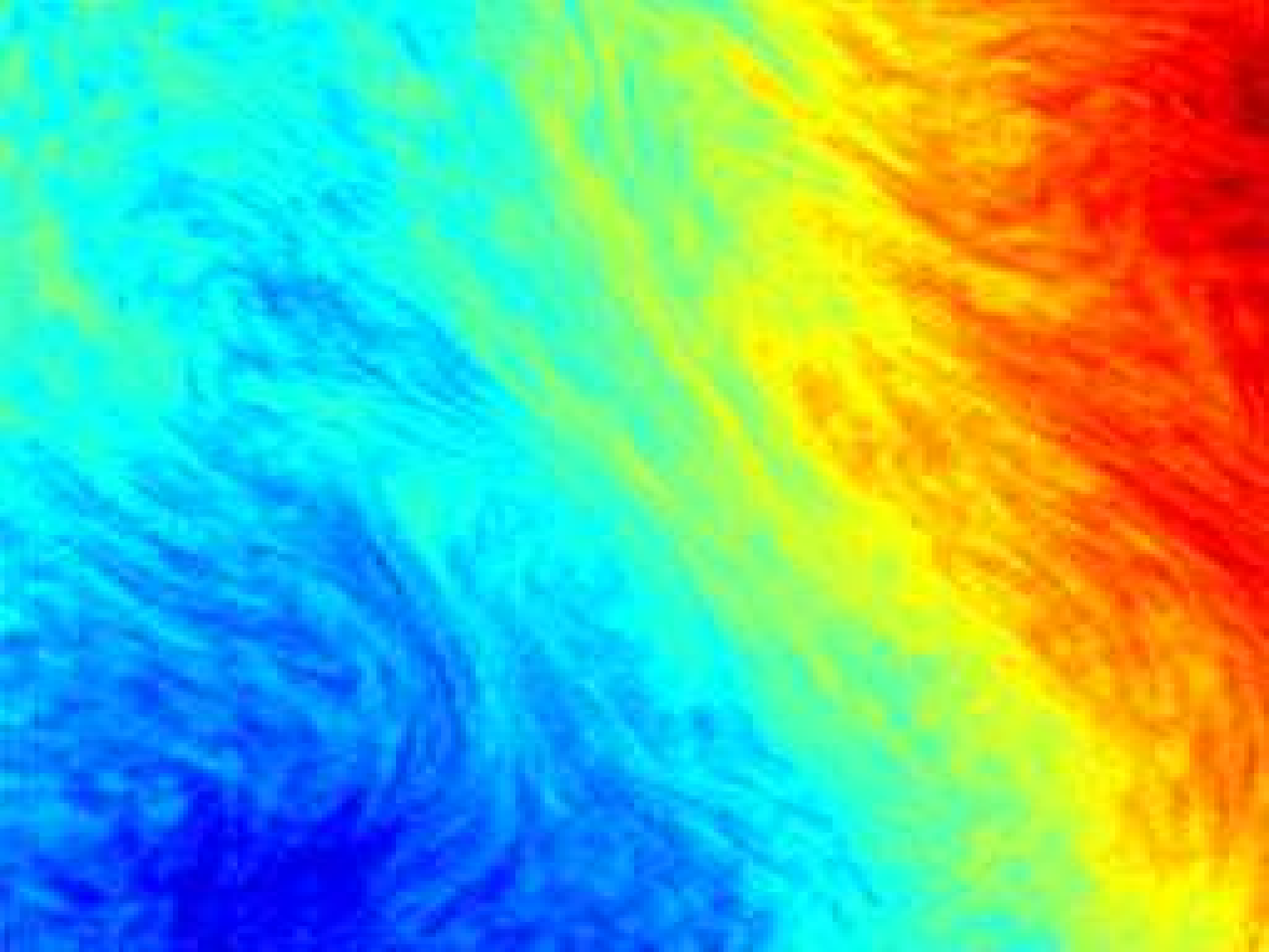}&
    \includegraphics[width=3cm]{shock_3d_ldmm_5p_band_29}
  \end{tabular}
  \caption{Compression of the 3D plasma (distribution function) data set with a 5\% data compression rate. Two  spatial cross  sections of the original data set are shown in the first figures on the first and third row. The results of Tucker decomposition, DCT, DFT, wavelet, and LDMM are shown in the remaining five figures.}
  \label{fig:compre_shock_3d_5p}
\end{figure}

In the numerical experiments, LDMM with random sampling has been used for each data set. The storage of SVD involves thresholded singular values along with the correponding singular vectors, and the storage of Tucker Decomposition involves a 3D core tensor with reduced size and three matrices for three different modes. For the other methods using global transforms, we store the coefficients with the largest magnitudes with constraint to the given budget. We mention that the results of Tucker Decomposition on 3D data sets are quite sensitive to the dimension of the core tensor along each direction. In our experiments, we choose the best result among all the possible decompositions satisfying the budget. This typically causes Tucker Decomposition to run for about two days on the 3D data sets reported in this paper. The visual and numerical results of the competing methods are reported in Figure \ref{fig:compre_2d_10p}-\ref{fig:compre_shock_3d_5p} and Table \ref{tab:error_compre_antonio_2d}-\ref{tab:error_compre_shock_3d}. As expected, the performance of LDMM in data compression is usually inferior compared to the other competing methods. However, it does outperform SVD in two of the more complicated 2D data sets (2D vortex and 2D plasma (distribution)) and the wavelet transform in the 3D plasma (magnetic field) data set. DCT almost consistently yields the best result among all the methods, and it can also be observed that tensor decomposition methods tend to achieve better results when the dimension of the data set becomes larger. Therefore, we can conclude that, at least at current stage, LDMM is a viable choice for data compression if the data set is complicated to begin with, and the user is willing to sacrifice accuracy for easy implementation in the compression step.

We point out that although LDMM does not perform equally well in data compression when compared to other methods that assume full access to the entire data set, there is still much room for improvement for LDMM. For instance, instead of randomly sampling the data set in the physical domain, we may strategically choosing pixels to sample if certain prior information is available. Moreover, if the original data set is known to the user, we can also modify the LDMM algorithm by sampling gradient values or certain entries in the weight matrices. Modifying LDMM for it to work as a data compression method will be the focus of our future work.

\begin{table}[H]
  \centering
  \begin{tabular}{||c| c  c c c  c||}
    \hline
    $10\%$ & SVD & DCT& DFT& Wavelet & LDMM\\
    \hline
    $L_1$       &0.0152 &\textbf{0.0003} &0.0010 &0.0005 &0.0034\\
    $L_2$       &0.0208 &\textbf{0.0005} &0.0014 &0.0007 &0.0075\\
    $L_\infty$   &0.1357 &\textbf{0.0056} &0.0132 &0.0067 &0.1376\\
    PSNR        &33.65  &\textbf{66.87}  &57.03  &63.01  &42.55\\
    \hline
    $5\%$ & SVD & DCT& DFT& Wavelet & LDMM\\
    \hline
    $L_1$       &0.0295 &\textbf{0.0011} &0.0024 &0.0016 &0.0056\\
    $L_2$       &0.0385 &\textbf{0.0015} &0.0035 &0.0021 &0.0111\\
    $L_\infty$   &0.1964 &0.0154 &0.0314 &\textbf{0.0149} &0.1102\\
    PSNR        &28.29  &\textbf{56.36}  &49.12  &53.36  &39.09\\
    \hline
  \end{tabular}
  \caption{Errors of the compression of the 2D vortex data set.}
  \label{tab:error_compre_antonio_2d}
\end{table}

\begin{table}[H]
  \centering
  \begin{tabular}{||c| c  c c c  c||}
    \hline
    $10\%$ & SVD & DCT& DFT& Wavelet & LDMM\\
    \hline
    $L_1$       &0.0345 &\textbf{0.0131} &0.0147 &0.0145 &0.0243\\
    $L_2$       &0.0437 &\textbf{0.0165} &0.0190 &0.0182 &0.0333\\
    $L_\infty$   &0.2597 &\textbf{0.0844} &0.1499 &0.0861 &0.1882\\
    PSNR        &27.19  &\textbf{35.63}  &34.44  &34.78  &29.56\\
    \hline
    $5\%$ & SVD & DCT& DFT& Wavelet & LDMM\\
    \hline
    $L_1$       &0.0494 &\textbf{0.0189} &0.0206 &0.0202 &0.0303\\
    $L_2$       &0.0624 &\textbf{0.0238} &0.0263 &0.0253 &0.0401\\
    $L_\infty$   &0.2794 &0.1121 &0.1920 &\textbf{0.1057} &0.2063\\
    PSNR        &24.10  &\textbf{32.47}  &31.59  &31.94  &27.93\\
    \hline
  \end{tabular}
  \caption{Errors of the compression of the 2D plasma (distribution) data set.}
  \label{tab:error_compre_shock_2d}
\end{table}

\begin{table}[H]
  \centering
  \begin{tabular}{||c| c  c c c  c||}
    \hline
    $10\%$ & SVD & DCT& FFT& Wavelet & LDMM\\
    \hline
    $L_1$       &0.0009 &\textbf{0.0001} &0.0004 &0.0006 &0.0013\\
    $L_2$       &0.0014 &0.0002 &0.0008 &\textbf{0.0001} &0.0040\\
    $L_\infty$   &0.0186 &0.0101 &0.0603 &\textbf{0.0011} &0.1979\\
    PSNR        &57.24  &75.73  &61.94  &\textbf{80.33}  &47.98\\
    \hline
    $5\%$ & SVD & DCT& DFT& Wavelet & LDMM\\
    \hline
    $L_1$       &0.0029 &0.0003 &0.0010 &\textbf{0.0002} &0.0022\\
    $L_2$       &0.0044 &0.0004 &0.0017 &\textbf{0.0004} &0.0062\\
    $L_\infty$   &0.0539 &0.0244 &0.0743 &\textbf{0.0049} &0.2097\\
    PSNR        &47.17  &67.59  &55.49  &\textbf{68.34}  &44.15\\
    \hline
  \end{tabular}
  \caption{Errors of the compression of the 2D lattice data set.}
  \label{tab:error_compre_lattice_2d}
\end{table}

\begin{table}[H]
  \centering
  \begin{tabular}{||c| c  c c c  c||}
    \hline
    $10\%$ & Tucker & DCT& DFT& Wavelet & LDMM\\
    \hline
    $L_1$       &0.0021 &\textbf{0.0014} &0.0024 &0.0068 &0.0038\\
    $L_2$       &0.0028 &\textbf{0.0018} &0.0038 &0.0089 &0.0062\\
    $L_\infty$   &0.0613 &\textbf{0.0433} &0.1757 &0.0739 &0.1330\\
    PSNR        &50.91  &\textbf{54.90}  &48.42  &41.01  &44.18\\
    \hline
    $5\%$ & Tucker & DCT& DFT& Wavelet & LDMM\\
    \hline
    $L_1$       &0.0040 &\textbf{0.0025} &0.0043 &0.0183 &0.0062\\
    $L_2$       &0.0054 &\textbf{0.0033} &0.0066 &0.0231 &0.0099\\
    $L_\infty$   &0.0911 &\textbf{0.0698} &0.2141 &0.1558 &0.2012\\
    PSNR        &45.36  &\textbf{49.70}  &43.56  &32.74  &40.07\\
    \hline
  \end{tabular}
  \caption{Errors of the compression of the 3D plasma (magnetic field) data set.}
  \label{tab:error_compre_plasma_3d}
\end{table}

\begin{table}[H]
  \centering
  \begin{tabular}{||c| c  c c c  c||}
    \hline
    $10\%$ & Tucker & DCT& DFT& Wavelet & LDMM\\
    \hline
    $L_1$       &\textbf{$\bm{9 \times 10^{-6}}$} &0.0002 &0.0007 &0.0002 &0.0008\\
    $L_2$       &\textbf{$\bm{1 \times 10^{-5}}$} &0.0005 &0.0022 &0.0002 &0.0038\\
    $L_\infty$   &\textbf{0.0002}            &0.1338 &0.2843 &0.0020 &0.4262\\
    PSNR        &\textbf{97.43}             &65.44  &52.96  &72.61  &48.43\\
    \hline
    $5\%$ & Tucker & DCT& DFT& Wavelet & LDMM\\
    \hline
    $L_1$       &\textbf{0.0001} &0.0004 &0.0010 &0.0006 &0.0013\\
    $L_2$       &\textbf{0.0001} &0.0009 &0.0031 &0.0008 &0.0051\\
    $L_\infty$   &\textbf{0.0042} &0.2053 &0.4266 &0.0095 &0.4530\\
    PSNR        &\textbf{78.28}  &60.52  &50.22  &61.25  &45.82\\
    \hline
  \end{tabular}
  \caption{Errors of the compression of the 3D lattice data set.}
  \label{tab:error_compre_lattice_3d}
\end{table}

\begin{table}[H]
  \centering
  \begin{tabular}{||c| c  c c c  c||}
    \hline
    $10\%$ & Tucker & DCT& DFT& Wavelet & LDMM\\
    \hline
    $L_1$       &0.0042 &\textbf{0.0039} &0.0045 &0.0042 &0.0060\\
    $L_2$       &0.0064 &\textbf{0.0052} &0.0061 &0.0055 &0.0105\\
    $L_\infty$   &0.0637 &0.0644 &0.0837 &\textbf{0.0373} &0.1181\\
    PSNR        &43.89  &\textbf{45.65}  &44.26  &45.17  &39.54\\
    \hline
    $5\%$ & Tucker & DCT& DFT& Wavelet & LDMM\\
    \hline
    $L_1$       &0.0060 &\textbf{0.0057} &0.0063 &0.0067 &0.0075\\
    $L_2$       &0.0092 &\textbf{0.0077} &0.0086 &0.0089 &0.0130\\
    $L_\infty$   &0.0890 &0.0766 &0.1018 &\textbf{0.0660} &0.1793\\
    PSNR        &40.75  &\textbf{42.29}  &41.28  &40.97  &37.72\\
    \hline
  \end{tabular}
  \caption{Errors of the compression of the 3D plasma (distribution function) data set.}
  \label{tab:error_compre_shock_3d}
\end{table}

\section{Conclusion}
\label{sec:conclusion}
In this paper, we propose a low dimensional manifold model for scientific data reconstruction from regular or irregular samplings. The low dimensionality of the patch manifold is used as a regularizer, and this assumption is justified through a dimension analysis of common patterns in various scientific data sets. The variational problem is solved via alternating direction of minimization, and the corresponding Laplace-Beltrami equation is discretized by weighted graph Laplacian. The proposed algorithm consistently outperforms all the competing algorithms in both regular and irregular sampling cases. The current LDMM algorithm as a data compression method does not perform as well as other standard compression algorithms that assume access to the full data set. But LDMM as a data compression method is easy to implement in the compression step, and it is also faster in the reconstruction step if only a subset of the original data set is required. Modifying LDMM for it to achieve its full potential as a data compression method will be the focus of our future work. 

\section{Acknowledgment}
The authors would like to thank Michael Crockatt and Professor Antonio Marquina for providing the neutron transport data set and the Orszag-Tang vortex data set.

\section*{References}
\bibliographystyle{elsarticle-num}
\bibliography{mybibfile}

\end{document}